\documentclass[a4paper,12pt,twoside,english]{book}
\usepackage[utf8]{inputenc}
\usepackage{caption}
\usepackage{minitoc}
\usepackage[T1]{fontenc}
\usepackage{babel}
\usepackage[a4paper,left=3cm,right=2.5cm,top=3cm,bottom=3cm, twoside]{geometry}

\usepackage{xcolor}
\usepackage{stmaryrd}
\usepackage{amssymb}
\usepackage{amsmath}
\usepackage{afterpage}
\usepackage{libertine}

\usepackage{hyperref} % for references (\ref, \label), url
\hypersetup{
    colorlinks=true,
    linkcolor=black,
    filecolor=black,      
    urlcolor=black,
    citecolor=black,
    pdfpagemode=FullScreen,
    hypertexnames=false,
    }
\usepackage[backend=bibtex, style=alphabetic, maxbibnames=6, sorting=nyt]
{biblatex} %Imports biblatex package

\usepackage{graphicx} % to include images
\graphicspath{ {figures/} }

\usepackage{subcaption}
\usepackage{float}
\usepackage{multirow}
\usepackage{array, multirow, tabularx}
\usepackage{booktabs}
\usepackage{longtable}
\usepackage{apalike}
\usepackage{minitoc}
\usepackage{tikz}
\usetikzlibrary{arrows.meta}
\usetikzlibrary{calc}
\usepackage{fancyhdr}
\usepackage[strict]{changepage}
\usepackage{siunitx}
\usepackage{enumitem}
\usepackage{epigraph}
\usepackage{lscape}
\usepackage{titletoc}
\usepackage{mdframed}

\usepackage{tabu}
\usepackage{pdfpages}

\usepackage{arydshln}
\usepackage{calc}
\definecolor{linkColor}{HTML}{32a852}

\addto\captionsfrench{%
}

\usepackage{sectsty}
\definecolor{darkseagreen}{rgb}{0.56, 0.74, 0.56}
\usepackage{lipsum}

\title{NMR global data Version2}
\author{Aleksandra Savina}
\date{October 2022}

\usepackage{emptypage}

\usepackage[acronym,xindy,toc]{glossaries} % nomain, if you define glossaries in a file, and you use \include{INP-00-glossary}
\makeglossaries

\usepackage{amsfonts}
\usepackage{amsmath}
\usepackage{amssymb}
\usepackage{amsthm}
\usepackage{mathtools}
\usepackage{stmaryrd}
\usepackage{dsfont}
\usepackage{graphicx}
\usepackage{hyperref}
\usepackage{titlesec}
\usepackage{mathabx}
\usepackage{enumitem,tocloft}
\usepackage{xcolor}
\usepackage{fourier-orns}
\usepackage[utf8]{inputenc}
\usepackage{framed}
\usepackage{supertabular}
\hypersetup{linktocpage,breaklinks=true}
\usepackage{mathrsfs}
\usepackage{csquotes}
\usepackage{bbold}

\SetSymbolFont{stmry}{bold}{U}{stmry}{m}{n}

\DeclareMathOperator*{\mystar}{*}

\theoremstyle{plain}
\newtheorem{theorem}{Theorem}[section]
\newtheorem{corollary}[theorem]{Corollary}

\newtheorem{theoremletter}{Theorem}

\newtheorem{corollaryletter}[theoremletter]{Corollary}
\newtheorem{propositionletter}[theoremletter]{Proposition}

\makeatletter

\renewcommand{\tableofcontents}{%
  \begingroup
  \renewcommand{\contentsname}{}%
  \@starttoc{toc}%
  \endgroup
}
\makeatother

\newtheorem{théorème}[theorem]{Théorème}
\newtheorem{corollaire}[theorem]{Corollaire}
\newtheorem{définition}[theorem]{Définition}
\newtheorem{lemma}[theorem]{Lemma}
\newtheorem{proposition}[theorem]{Proposition}

\newtheorem{claim}[theorem]{Claim}
\theoremstyle{definition}
\newtheorem{definition}[theorem]{Definition}
\newtheorem{example}[theorem]{Example}
\newtheorem{question}[theorem]{Question}
\newtheorem{remark}[theorem]{Remark}
\newtheorem{fact}[theorem]{Fact}

\numberwithin{equation}{section}

\newcommand{\field}{\mathbb{f}}

\newcommand{\R}{\mathbb{R}}% Real numbers
\newcommand{\dis}{\displaystyle}

\newcommand{\la}{\langle}
\newcommand{\ra}{\rangle}
\newcommand{\ld}{\mathrm{L}}

\newcommand{\Q}{\mathbb{Q}}	% Rational numbers
\newcommand{\Z}{\mathbb{Z}}	% Integers
\newcommand{\N}{\mathbb{N}}	% Natural numbers
\newcommand{\fsym}[1]{\mathsf{FSym}(#1)}

\newcommand{\sym}[1]{\mathsf{Sym}(#1)}

\newcommand{\shuf}[1]{\mathsf{Shuffler}(#1)}

\newcommand{\shufn}[2]{\mathsf{Shuffler}^{\circ {#1}}(#2)}

\newcommand{\juggler}[2]{\mathsf{Shuffler}_{#1}(#2)}

\newcommand{\jugglern}[3]{\mathsf{Shuffler}_{#2}^{\circ {#1}}(#3)}

\newcommand{\cloner}[1]{\mathsf{Cloner}_{\field}(#1)}

\newcommand{\clonern}[2]{\mathsf{Cloner}_{\field}^{\circ {#1}}(#2)}

\newcommand{\designer}[1]{\mathsf{Designer}_{F}(#1)}

\newcommand{\act}{\curvearrowright}

\newcommand{\supp}{\mathrm{supp}\ } %support d'une fonction

\newcommand{\diam}{\mathrm{diam}\ } %diamètre d'une partie

\newcommand{\prof}[1]{j_{1,#1}}
\newcommand{\profp}[1]{j_{p,#1}}
\newcommand{\folp}[1]{\text{F\o l}_{p,#1}}

\newcommand{\halo}{\mathscr{L}}

\newcommand{\rem}{\backslash\backslash}

\makeatletter
\newcommand*{\defeq}{\mathrel{\rlap{%
                     \raisebox{0.3ex}{$\m@th\cdot$}}%
                     \raisebox{-0.3ex}{$\m@th\cdot$}}%
                     =}
\makeatother

\renewcommand{\tableofcontents}{%
  \begingroup
  \renewcommand{\contentsname}{}%
  \@starttoc{toc}%
  \endgroup
}

\usepackage[xindy]{imakeidx}
\makeindex

\mtcsettitle{minitoc}{}
\mtcsetoffset{minitoc}{-1.0em}  % To shift the minitoc to the left, if needed!
\mtcsetdepth{minitoc}{1} 
\mtcsetfont{minitoc}{*}{\small\rmfamily\upshape\mdseries} 
\mtcsetfont{minitoc}{section}{\small\rmfamily\upshape\bfseries} 
\begin{document}
\dominitoc
%\dominilof
%\dominilot

%\maketitle
%----------------------------------------------------------------------------------------
%	TITLE PAGE
%----------------------------------------------------------------------------------------
\newgeometry{left=2cm,right=2cm,top=2.5cm,bottom=2.5cm}
	\begin{titlepage}
		\begin{center}
\begin{tabular}{@{}m{0.45\textwidth}m{0.45\textwidth}@{}}
  \raggedright
  \includegraphics[height=2.5cm]{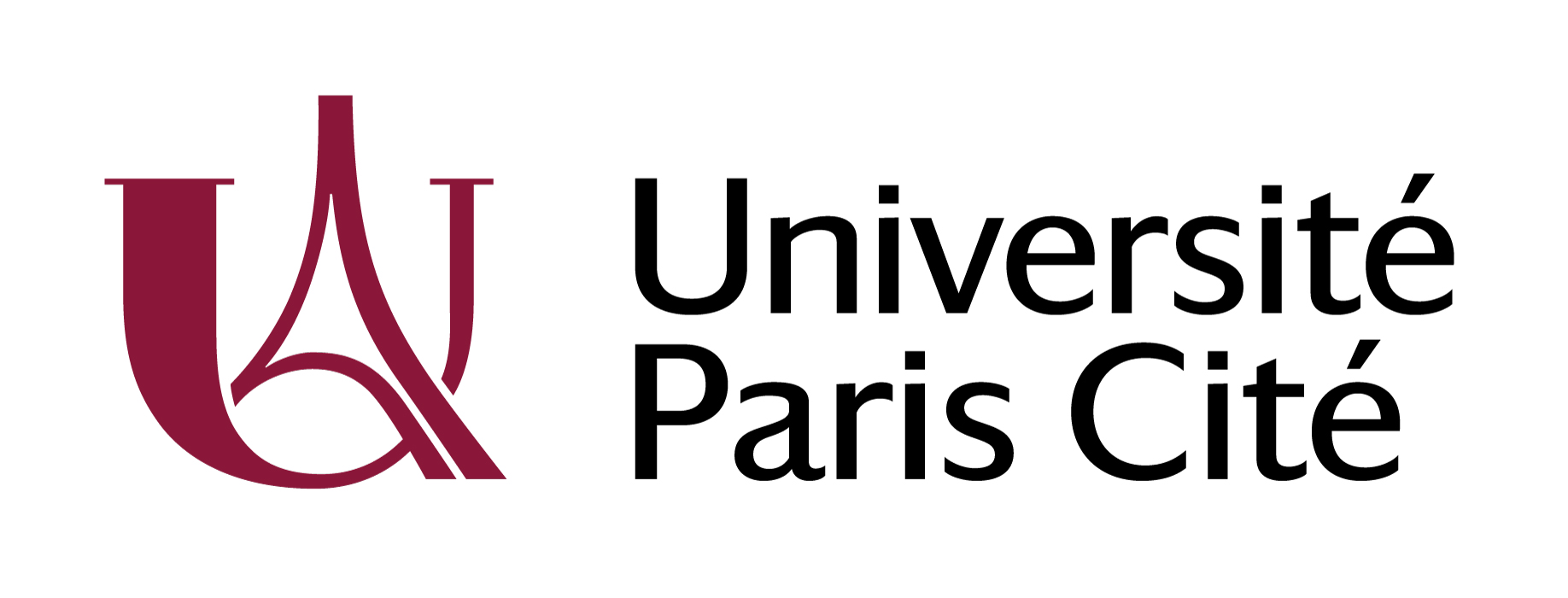}
  &
  \raggedleft
  \includegraphics[height=3cm]{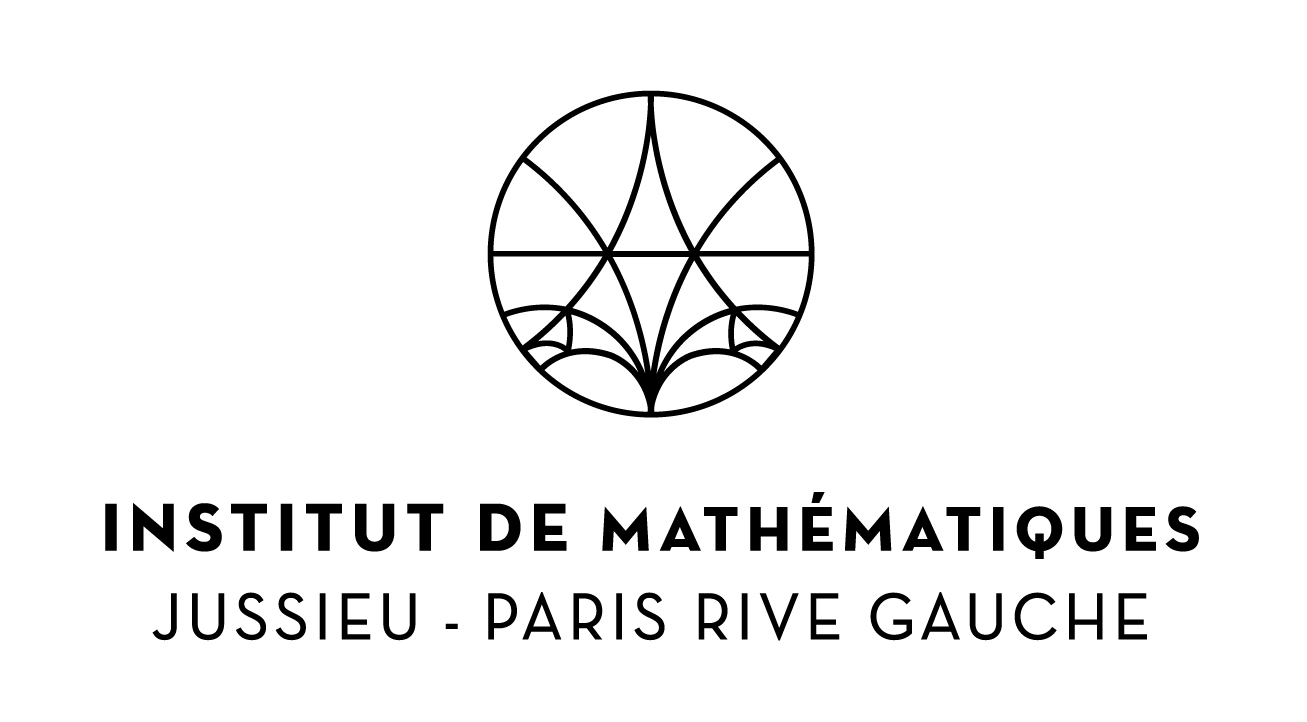}
\end{tabular}
\end{center}
	
		\begin{center}
		
\vspace*{.03\textheight}
\textsc{\LARGE Université Paris Cité}\\[0.2cm] % Univ name
		\large École doctorale de Sciences Mathématiques de Paris Centre - ED386\\
		  Institut de Mathématiques de Jussieu-Paris Rive Gauche - UMR7586\\ 
        %Équipe de recherche\\
  
  			\vfill
 
	 		\rule{\textwidth}{0.8pt} \\ % Horizontal line
	 		\vspace{10pt}
	 		 { \LARGE \bfseries Analytic, geometric and measured aspects of lamplighter-like groups} % Thesis title
	 		 \vspace{10pt}
	 		 \rule{\textwidth}{0.8pt} \\ % Horizontal line
		\end{center}
		
		\vfill
		% Author and supervisor	
		\begin{center}
			Par \textsc{\Large Vincent Dumoncel}\\[1cm] 
			\textsc{\Large Thèse de doctorat}\\[1cm]
			\textsc{\large Mathématiques Fondamentales}\\[1.2cm]
			Dirigée par \textsc{\large Anthony Genevois} et \textsc{\large Romain Tessera} \\[0.2cm] 
			Présentée et soutenue publiquement le 7 juillet 2026 %renseigner la date et choisir le type de soutenance
		\end{center}

		\vspace{1cm}
		\begin{center}
Devant un jury composé de : %(indiquer le jury complet) 
        \end{center}
		%\begin{center}
			\begin{tabular}{l@{\hskip 0.15cm}l@{\hskip 0.15cm}l@{\hskip 0.15cm}l}
				M$^{\text{me}}$ Goulnara Arzhantseva & Professor & University of Vienna & Examinatrice \\
				M$^{\text{me}}$ Amandine Escalier & Maîtresse de conférences & Université Claude Bernard Lyon 1 & Examinatrice\\
				M. David Fisher & Professor & Rice University & Rapporteur\\
                M. Anthony Genevois & Chargé de recherche  &Université de Montpellier &Membre invité\\
				M. Gabriel Pallier & Maître de conférences  &Université de Lille & Examinateur \\
				M. Christophe Pittet & Professeur des Universités  & Université Aix-Marseille & Rapporteur \\
				M. Romain Tessera  & Directeur de recherche &Université Paris Cité & Directeur \\
				M$^{\text{me}}$ Tianyi Zheng & Professor & UC San Diego & Examinatrice
			\end{tabular}\\%[0.2cm]
		%\end{center}

	\newpage %page de garde vide
	\thispagestyle{empty} %page de garde vide
	\end{titlepage}
\restoregeometry

\renewcommand{\chaptermark}[1]{\markboth{#1}{}}

\frontmatter % numérote les pages en chiffres romains jusqu'a la commande mainmatter

%----------------------------------------------------------------------------------------
%	RESUME ET ABSTRACT PAGE
%----------------------------------------------------------------------------------------
\chapter*{Résumé}
\addcontentsline{toc}{chapter}{Résumé}  
\mtcaddchapter

\vspace{1.5cm}

\textbf{Titre : Aspects analytiques, géométriques et mesurés de variantes des allumeurs de réverbères}

\vspace{0.75cm}

\noindent \textbf{Mots-clés : Quasi-isométries, produits en couronne permutationnels, produits en halo, profils isopérimétriques, équivalence orbitale quantitative.} 

\vspace{1cm}

Les produits en couronne forment une classe très étudiée en théorie des groupes. En effet, leur définition est, d'une part, suffisamment explicite pour en comprendre la structure précisément et, d'autre part, suffisamment élaborée pour refléter des comportements exotiques et inattendus, non-observés dans des classes de groupes plus « classiques ». En théorie géométrique des groupes, ils sont connus pour exhiber souvent des propriétés intéressantes et fournissent des contre-exemples à des questions cruciales. Ils permettent par exemple de démontrer que la résolubilité virtuelle n'est pas stable par quasi-isométries, ou que deux groupes moyennables peuvent être quasi-isométriques mais non bijectivement quasi-isométriques. 

\smallskip 

Dans cette thèse, on étudie, sous plusieurs angles, des groupes de type fini proches des produits en couronne, tous structurés comme des produits semi-directs. Le premier résultat majeur est une classification à quasi-isométrie près d'une certaine classe de produits en couronne permutationnels. La stratégie présentée repose sur des travaux récents, pour les produits en couronne standard, de Genevois et Tessera, mais nous mène à un phénomène de rigidité plus fort : sous certaines hypothèses, une quasi-isométrie entre deux produits en couronne permutationnels induit toujours une quasi-isométrie de paires entre les groupes de base munis des sous-groupes normaux correspondants. On montre aussi que, dans certains cas spécifiques, être quasi-isométrique est équivalent à être bijectivement quasi-isométrique. 

\smallskip 

Un chapitre est également dédié à l'étude de produits en couronne standard mais avec groupes de lampes à croissance polynomiale. On établit un résultat de rigidité pour les quasi-isométries entre de tels espaces : elles doivent toutes être « measure-scaling ». En particulier, les auto-quasi-isométries d'un tel produit en couronne sont toutes à distance bornée d'une bijection. On montre également que cette classe d'allumeurs de réverbères fournit, elle aussi, des exemples de groupes moyennables quasi-isométriques mais non bijectivement quasi-isométriques.

\smallskip 

La deuxième partie de cette thèse s'appuie sur des travaux en commun avec Corentin Correia, dans lesquels on étudie la classe des \textit{produits en halo} de Genevois et Tessera d'un point de vue mesuré. On construit des couplages d'équivalence orbitale entre deux produits en halo de même type, et, en exploitant la technologie des \textit{F\o lner tilings}, introduite récemment, on exhibe des couplages entre produits en halo et groupes abéliens libres. On établit également plusieurs estimations asymptotiques des profils $\ell^p$ de tels groupes, améliorant des résultats antérieurs d'Erschler-Zheng et Saloff-Coste-Zheng. Ces calculs nous permettent de montrer que, dans bien des cas, les couplages d'équivalence orbitale ainsi obtenus sont quantitativement optimaux.  

 %résumé en francais
\chapter*{Abstract}
\addcontentsline{toc}{chapter}{Abstract}  
\mtcaddchapter

\vspace{1.5cm}

\textbf{Title : Analytic, geometric and measured aspects of lamplighter-like groups}

\vspace{0.75cm}

\noindent \textbf{Keywords : quasi-isometries, permutational wreath products, halo products, isoperimetric profiles, quantitative orbit equivalence.} 

\vspace{1cm}

Wreath products are intensively studied in group theory. Their definition is, on the one hand, explicit enough to understand precisely their structure and, on the other hand, elaborate enough to reflect unexpected and exotic behaviours that do not occur in more “classical” classes of groups. In geometric group theory, wreath products often offer counter-examples to important questions. For instance, they allow us to prove that virtual solvability is not stable under quasi-isometries, or that two amenable groups can be quasi-isometric without being bijectively quasi-isometric. 

\smallskip 

In this thesis, we study finitely generated groups close to wreath products, all defined as semi-direct products, from various perspectives. The first main result is a quasi-isometric classification of some permutational wreath products. The strategy we follow relies on recent work of Genevois and Tessera, but yields a stronger rigidity phenomenon: under suitable assumptions, a quasi-isometry between two permutational lamplighters always induces a quasi-isometry of pairs between the base groups equipped with the corresponding normal subgroups. We also show that, in some specific cases, being quasi-isometric is equivalent to being bijectively quasi-isometric.

\smallskip 

A chapter is also dedicated to the study of standard wreath products whose lamp groups have polynomial growth. We establish a rigidity result for quasi-isometries between such wreath products: they must all be “measure-scaling”. In particular, self-quasi-isometries of such wreath products all lie at bounded distance from bijections. We also prove that this class of lamplighters provides examples of amenable quasi-isometric groups that are not bijectively quasi-isometric. 

\smallskip 

The second part of this thesis relies on joint works with Corentin Correia, in which we study \textit{halo products} of Genevois and Tessera from a measured point of view. We construct orbit equivalence couplings between two halo products of the same kind and, using the recent technology of \textit{F\o lner tilings}, we exhibit couplings between halo products and free abelian groups. We also establish several asymptotic estimates of the $\ell^p-$isoperimetric profiles of these groups, improving previous results of Erschler-Zheng and Saloff-Coste-Zheng. These computations allow us to prove that, in many cases, the orbit equivalence couplings thus obtained are quantitatively optimal.  %résumé en anglais

%----------------------------------------------------------------------------------------
%	QUOTATION PAGE ou DEDICACE
%----------------------------------------------------------------------------------------
\clearpage
\vspace*{0.2\textheight}

\begin{quote}
« Ce qui ne me tue pas me rend plus fort »
\end{quote} \bigbreak

\hfill 
Friedrich Nietzsche, \textit{Crépuscule des idoles ou Comment on philosophe avec un marteau}, 1888.

%----------------------------------------------------------------------------------------

%----------------------------------------------------------------------------------------
%	ACKNOWLEDGEMENTS
%----------------------------------------------------------------------------------------
\chapter*{Remerciements}
\addcontentsline{toc}{chapter}{Remerciements}

Contrairement à ce que beaucoup pensent, une thèse est tout sauf une expérience de solitude. On croise au fil du chemin de nombreuses personnes, la plupart vous aidant, à leur façon, à rendre cette expérience plus savoureuse. Ces quelques pages leur sont dédiées. 

\vspace{0.15cm}

Mes premiers remerciements vont d'abord et avant tout à Romain et Anthony pour avoir rendu ce projet possible et m'avoir aidé à le faire aboutir. Ils ont tous les deux apporté beaucoup à cette thèse. 

\vspace{0.15cm}

Romain, merci d'avoir ouvert ta porte quand j'y ai frappé il y a presque quatre ans. Merci de m'avoir fait confiance et d'avoir proposé l'idée d'un co-encadrement. J'en ai tiré de précieux enseignements qui m'ont beaucoup aidé sur le plan scientifique, professionnel et personnel. Merci de l'enthousiasme dont tu fais preuve quand tu expliques les maths. La façon dont tu vois et abordes la recherche me fascinera toujours. Merci pour tes encouragements lorsque je séchais ou que je doutais, ils ont été déterminants. 

\vspace{0.15cm} 

Anthony, merci d'avoir accepté de co-encadrer cette thèse et de m'avoir accueilli plusieurs fois à Montpellier. Merci pour les nombreuses heures passées devant un tableau à réfléchir, à refaire des preuves qui ne marchaient jamais et à finir par me dire « oh écoute je ne sais pas ». Merci d'avoir passé tout ce temps à m'écouter raconter ma vie et à me répondre avec tellement de conseils pertinents. Ils ont été, je crois, décisifs au moment de choisir ce que je souhaitais faire après la thèse. Merci pour tout. 

\vspace{0.15cm}

Merci à Christophe Pittet et David Fisher pour avoir accepté de rapporter cette thèse et pour leurs commentaires sur mon travail. 

\vspace{0.15cm}

Merci à Gabriel Pallier, Amandine Escalier, Goulnara Arzhantseva et Tianyi Zheng d'avoir accepté de faire partie du jury et d'examiner cette thèse. Merci aussi à Gabriel et Amandine pour leurs invitations à Lille et à Lyon, leurs accueils chaleureux et leurs nombreux retours sur mes exposés et le manuscrit. 

\vspace{0.15cm} 

Je souhaite remercier aussi toutes celles et ceux rencontrés en dehors de Paris et qui ont à chaque fois fait preuve de bienveillance et de gentillesse, au premier rang desquels Jérémie Brieussel et Thomas Haettel à Montpellier, et Damien Gaboriau, Mikaël de la Salle et Amandine Escalier un peu partout en conférence. 

\vspace{0.15cm} 

Merci à toute l'équipe d'algèbres d'opérateurs de l'IMJ pour la bonne ambiance de travail qu'elle instaure, en particulier mon tuteur parti trop tôt François Le Maître et son redoutable successeur Pierre Fima. Merci aussi à Denis pour les cours dispensés ensemble et les discussions sympathiques au détour des couloirs de l'IMJ. 

\vspace{0.15cm} 

Je salue les doctorants du labo qui y font régner une ambiance de détente plus que de travail et qui ont contribué à rendre cette thèse plus agréable. Un grand merci à Juan, Francesca, MC, Francesca, Mathieu et Laura pour leur amitié, à Mario et Fabien pour leur bienveillance et à Kostya pour l'introduction aux opérateurs of finite propagation, que je maîtrise maintenant mieux que les quasi-isométries. Merci à Werner, Ivory (qui a parlé à notre groupe de travail, un couteau sous la gorge), Bryan, Paul et Enzo de faire du 6ème étage ce qu'il est. Merci à Léo et Lorenzo de faire vivre ce bureau 752, que j'ai délaissé trop souvent. Merci à Juan pour de nombreux commentaires pertinents sur mes travaux, et pour m'avoir assuré l'année dernière avec un aplomb extraordinaire que le profil $\ell^p$ des groupes de Grigorchuk était connu... moi qui allais relayer l'information partout en séminaire, avant de me rendre compte qu'il m'avait menti ! 

\vspace{0.15cm}  

Merci à mon frère de thèse Corentin pour tous ces reels envoyés sur Insta, l'argent dépensé à Honorine et ses reprises de « Flashback », de Gazo et Favé. Je n'oublierai pas les nombreuses conneries qu'il m'a racontées lorsque nous étions dans la même chambre en conférence, ses fantastiques imitations de Romain Tessera et Didier Deschamps et ses cours improvisés de portugais. Merci aussi pour ces heures passées à travailler, qui ont abouti à de jolis résultats, certains figurant ci-dessous, de sympathiques \textit{joint talks}, dont un que j'ai planté magistralement, et à ce jour où nous avons passé huit heures à vérifier qu'une action de groupe était bien une action (mais au moins on est sûrs maintenant).

\vspace{0.15cm} 

Un grand merci également aux autres doctorants et postdocs croisés en conférence ou à d'autres occasions, qui ont tous fait en sorte de rendre ces exposés un peu moins longs et ces soirées un peu plus alcoolisées : Justin, Maximilien, Gaëtan, Patrick, Hermès, Antonio et Jorge. Merci aussi à Eduardo pour son invitation et son accueil en Allemagne, et à Cosmas (the “$\cos$ of mass”) pour son enthousiasme concernant les lampshufflers. 

\vspace{0.15cm}

La thèse ne dure que trois ans, mais sa fin est, pour moi, l'achèvement d'un chemin bien plus long, entamé il y a bientôt huit ans à l'Université de Genève, poursuivi à Lausanne et ensuite à Paris. J'ai eu la chance de rencontrer et de discuter, souvent aux bons moments, avec des personnes dont les conseils ont été décisifs. Merci à David Cimasoni, Donna Testerman, Nicolas Monod, Florian Richter, Sacha Friedli et François Genoud pour m'avoir aidé et m'avoir fait comprendre que l'ambition ne devrait pas connaître de limites. Tout particulièrement, merci à Donna Testerman pour son soutien constant depuis le début de mes études à l'EPFL et tout au long de cette thèse. 

\vspace{0.15cm} 

Durant toutes ces années, j'ai aussi pu compter sur le soutien d'amis proches, rencontrés à la fac, au lycée, ou ailleurs, et à qui je dois beaucoup. Un grand merci à Kamil, Christophe, JV, Mélissa, Baptiste, Camille, Camille, Mathilde, Rebecca et Guillaume. 

\vspace{0.15cm} 

Je suis reconnaissant à ma famille pour son soutien indéfectible, en particulier ma mère qui, même après toutes ces années, ne désespère pas de comprendre un jour la géométrie à grande échelle des lamplighters. Merci à mon frère François pour de nombreuses contributions ayant optimisé l'écriture de ce manuscrit. Merci à Christian et Eugénie, respectivement docteur et future doctorante, pour leur intérêt dans cette thèse. 

\vspace{0.15cm}

Je dédie ce manuscrit à ceux qui me regardent depuis là-haut, notamment mon grand-père qui, j'en suis sûr, aurait été fier d'assister à ma soutenance de thèse. 

\vspace{0.15cm}  

Enfin, et même si elle n'est plus là, merci à Alex. Que ce soit à Genève, Lausanne ou Paris, elle a été pendant longtemps ma principale raison de me lever chaque jour, et ma première fan dans tous mes projets. Cette thèse n'aurait pas abouti sans son soutien. Pour cela, et pour tout le reste, elle sait mon affection, ma reconnaissance et ma gratitude. 
\mtcaddchapter

%----------------------------------------------------------------------------------------
%	LIST OF CONTENTS/FIGURES/TABLES PAGES
%----------------------------------------------------------------------------------------
%\setcounter{secnumdepth}{3} % organisational level that receives a numbers
%\setcounter{tocdepth}{3}
\renewcommand\contentsname{Contents} % the empty name
\begingroup
\vspace{-10cm} % the removed space. Set as appropriate
\newpage\tableofcontents
\endgroup % Prints the main table of contents
%\listoffigures % Prints the list of figures
%\addcontentsline{toc}{chapter}{Liste des figures}
%\listoftables % Prints the list of tables
%\addcontentsline{toc}{chapter}{Liste des tableaux}
\newpage

\newcommand\blankpage{%
    \null
    \thispagestyle{empty}%
    \addtocounter{page}{-1}%
    \newpage}

%----------------------------------------------------------------------------------------
%	Liste des abréviations-
%----------------------------------------------------------------------------------------

% --------------------------------------------------------------
%                         Début du corps
% --------------------------------------------------------------

\mainmatter
%\setlength{\parskip}{.7em}

% style thèse classique) : indentation, pas d’espace entre paragraphes
\setlength{\parindent}{1.2em}
\setlength{\parskip}{0pt}

% petit espace propre
%\setlength{\parindent}{0pt}
%\setlength{\parskip}{0.25\baselineskip plus 1pt minus 1pt}

%\titlespacing*{\section}{0pt}{.9em}{.8em}
\renewcommand{\baselinestretch}{1.1}

\fancyhead[RO]{\leftmark}
\fancyhead[LE]{\textsc{\chaptername~\thechapter}}

\fancyhead{} % clear all header fields
\fancyhead[OL]{\textsc{Résumé long en français}}
\chapter*{Résumé long en français}
\addcontentsline{toc}{chapter}{Résumé long en français}  
\mtcaddchapter

\begingroup
\renewcommand{\thetheorem}{\arabic{theorem}}
\setcounter{theorem}{0}

Si $(X,d_{X})$ et $(Y,d_{Y})$ sont deux espaces métriques, une~\textit{quasi-isométrie} entre $X$ et $Y$ est une application $f\colon X\rightarrow Y$ pour laquelle il existe deux constantes $C\ge 1$, $K\ge 0$ telles que 
\begin{equation*}
    \frac{1}{C}\cdot d_{X}(x,y)-K \le d_{Y}(f(x), f(y)) \le C\cdot d_{X}(x,y)+K
\end{equation*}
pour tous $x,y\in X$, et $d_{Y}(y, f(X))\le K$ pour tout $y\in Y$. Quand une telle application existe, on dit que $X$ et $Y$ sont~\textit{quasi-isométriques}. Il est souvent difficile de déterminer si deux espaces métriques sont quasi-isométriques et, lorsque c'est le cas, de décrire explicitement l'ensemble des quasi-isométries entre ces deux espaces. On se restreint donc souvent à des espaces ayant une structure supplémentaire. 

Une telle classe est celle des~\textit{groupes de type fini}, qui ont naturellement une structure métrique. Comprendre ces espaces à quasi-isométrie près est motivé en particulier par un résultat célèbre de Gromov sur les groupes à croissance polynomiale, qui montre que la structure algébrique d'un groupe de type fini est gouvernée par sa géométrie à grande échelle :

\begin{théorème}[\cite{Gro81}]
Soit $G$ un groupe de type fini. Alors $G$ est à croissance polynomiale si et seulement s'il contient un sous-groupe nilpotent d'indice fini.
\end{théorème}

La~\textit{croissance} d'un groupe est peut-être l'exemple le plus élémentaire d'un invariant de quasi-isométrie ; cette fonction compte le nombre d'éléments dans une boule centrée en $1_{G}$ en fonction du rayon de la boule. Le fait que cette fonction ait un comportement polynomial lorsque le rayon grandit est préservé par les quasi-isométries et, avec le théorème de Gromov, il suit qu'avoir un sous-groupe nilpotent d'indice fini est aussi un invariant de quasi-isométrie. 

La thèse aborde plusieurs autres invariants de quasi-isométrie, connus sous le nom de profils isopérimétriques $\ell^p$. Ayant une définition plus élaborée, ces profils sont plus durs à calculer que les fonctions de croissance, mais renferment aussi plus d'informations sur la géométrie du groupe. On renvoie à la Section~\ref{sec:isoprof} pour plus de détails. 

Dans une autre direction, la littérature est abondante sur les groupes de type fini qui sont~\textit{quasi-isométriquement rigides}. Un groupe $G$ est~\textit{quasi-isométriquement rigide} si, pour tout groupe $G'$ quasi-isométrique à $G$, $G'$ est isomorphe à un sous-groupe d'indice fini de $G$ ou à un quotient de $G$ par un sous-groupe normal fini. Par exemple, les groupes libres abéliens~\cite{Sha04} et non-abéliens sont quasi-isométriquement rigides~\cite{Sta68, Dun85}, de même que les groupes de Baumslag-Solitar résolubles~\cite{FM98}. La question est ouverte pour beaucoup d'autres classes, comme les groupes polycycliques, les groupes d'Artin à angles droits et les groupes d'automorphismes extérieurs de groupes libres. 

Enfin, la littérature montre aussi que, même restreint à une classe particulière, des techniques poussées sont nécessaires pour espérer décrire toutes les quasi-isométries entre deux groupes de la classe. À titre d'exemple, mentionnons le résultat célèbre de Farb et Mosher sur la classification à quasi-isométrie près des groupes de Baumslag-Solitar résolubles :

\begin{théorème}[{\cite[Théorème~7.1]{FM98}}]
Soient $n,m\ge 2$. Alors $\text{BS}(1,n)$ et $\text{BS}(1,m)$ sont quasi-isométriques si et seulement si $n$ et $m$ sont puissances d'un même nombre.
\end{théorème}

Le problème de rigidité quasi-isométrique a été, jusqu'à maintenant, exploré principalement pour deux grandes familles de groupes : ceux qui exhibent des propriétés de courbure négative, et ceux qui sont moyennables (et souvent plus, par exemple résolubles). Cette thèse explore une troisième direction : la classe des~\textit{produits en couronne} et toutes leurs variantes, qu'on présente maintenant. 

\paragraph{Les principaux protagonistes de cette thèse.} Etant donnés deux groupes de type fini $A$ et $B$, leur~\textit{produit en couronne} est le groupe de type fini $A\wr B$ défini comme
\begin{equation*}
    A\wr B \defeq \left(\bigoplus_{B}A\right)\rtimes B
\end{equation*}
où $B$ agit sur la somme directe en permutant les coordonnées selon son action sur lui-même par translation à gauche. Le groupe $A$ s'appelle classiquement le groupe des « lampes », et le groupe $B$ est le groupe « de base », à cause d'une description visuelle des éléments d'un produit en couronne selon des coloriages de graphes et une flèche (ou un petit personnage) qui se déplace sur ce graphe. En effet, soit par exemple $(c,p)\in \Z/2\Z\wr\Z$. On imagine que $c\in\bigoplus_{\Z}\Z/2\Z$ est un coloriage des sommets du graphe de Cayley de $\Z$, qui est une rue infinie, avec seulement deux couleurs possibles : $0$ ou $1$. On imagine que $p\in\Z$ est la position d'un allumeur de réverbères. Un ensemble générateur de $\Z/2\Z\wr\Z$ peut être construit en deux parties : une partie correspond aux changements de coloriages, et l'autre partie correspond aux changements de position de l'allumeur. Il y a alors deux façons de se déplacer de $(c,p)$ à un voisin dans $\text{Cay}(\Z/2\Z\wr\Z,U)$ : 
\begin{itemize}
    \item soit l'allumeur se déplace dans $\Z$, allant de $p$ à un voisin de $p$, et le coloriage reste inchangé;
    \item soit l'allumeur reste à sa position, mais il modifie la couleur du sommet où il se trouve, allant à $1$ si cette couleur est $0$, et à $0$ si cette couleur est $1$. 
\end{itemize}

L'exemple~\ref{ex:descriptionoflamplighters} contient une description plus générale. 

Maintenant, dans le cadre du programme de Gromov, on peut se demander :

\begin{question}\label{INTROFRANCAISquestionGromov}
Soient $A_{1}$, $A_{2}$, $B_{1}$, $B_{2}$ des groupes de type fini. Quand est-ce que $A_{1}\wr B_{1}$ et $A_{2}\wr B_{2}$ sont quasi-isométriques ? 
\end{question}

\paragraph{Classifications à quasi-isométrie près.} En 2013, Alex Eskin, David Fisher et Kevin Whyte ont répondu complètement à cette question pour des allumeurs de réverbères au-dessus de groupes cycliques.

\begin{théorème}[{\cite[Théorème~1.2]{EFW13}}]
Soient $F_{1}$, $F_{2}$ deux groupes finis non-triviaux. Alors $F_{1}\wr\Z$ et $F_{2}\wr\Z$ sont quasi-isométriques si et seulement s'il existe des entiers $a,r,s\ge 1$ tels que $|F_{1}|=a^{r}$ et $|F_{2}|=a^{s}$. 
\end{théorème}

En particulier, ce résultat fait apparaître clairement une dépendance entre l'existence d'une quasi-isométrie entre les groupes concernés et la cardinalité des groupes de lampes qu'on choisit. 

Une deuxième contribution majeure est due à Anthony Genevois et Romain Tessera : 

\begin{théorème}[{\cite[Corollaire~1.6]{GT24b}}]\label{thm:INTROFRANCAISclassificationGT21}
Soient $F_{1}$ et $F_{2}$ deux groupes finis non-triviaux. Soient $H_{1}$ et $H_{2}$ deux groupes de présentation finie. Supposons que $H_{1}$ est à un bout.
\begin{enumerate}[label=(\roman*)]
    \item Si $H_{1}$ n'est pas moyennable, alors $F_{1}\wr H_{1}$ et $F_{2}\wr H_{2}$ sont quasi-isométriques si et seulement si $|F_{1}|$, $|F_{2}|$ ont les mêmes diviseurs premiers et il existe une quasi-isométrie $H_{1}\rightarrow H_{2}$.
    \item Si $H_{1}$ est moyennable, alors $F_{1}\wr H_{1}$ et $F_{2}\wr H_{2}$ sont quasi-isométriques si et seulement s'il existe des entiers $a,r,s\ge 1$ tels que $|F_{1}|=a^{r}$, $|F_{2}|=a^{s}$ et il existe une quasi-$\frac{s}{r}$-to-one quasi-isométrie $H_{1}\rightarrow H_{2}$.
\end{enumerate}
\end{théorème}

Deux remarques s'imposent ici : d'abord, ce résultat s'applique lorsque les groupes de base sont de présentation finie et à un bout, ce qui couvre déjà un grand nombre de cas, de telles propriétés étant courantes pour des groupes infinis de type fini. Ensuite, il rend encore plus explicite la dépendance entre la quasi-isométrie qui relie les groupes de base et les cardinalités des groupes de lampes, en recourant, dans le cas moyennable, à la notion de « quasi-(something)-to-one » quasi-isométrie. Ce raffinement quantitatif des quasi-isométries a été introduit formellement dans~\cite{GT22}, et nous occupera largement dans cette thèse (cf. en particulier la Section~\ref{sec:ScalingQI} et le Chapitre~\ref{chap:chapter4}).

Enfin, mentionnons qu'il existe une version plus générale du Théorème~\ref{thm:INTROFRANCAISclassificationGT21}, avec des hypothèses plus faibles, établie dans~\cite{GT24a}. On revient sur cet énoncé plus large dans le Chapitre~\ref{chap:chapter2}. Dans cette introduction, le Théorème~\ref{thm:INTROFRANCAISclassificationGT21} est suffisant en guise de motivation. 

\paragraph{Le match quasi-isométries vs équivalences biLipschitz.} Dans une autre direction, c'est une question importante de comprendre à quel point les notions de quasi-isométrie et de quasi-isométrie bijective (qu'on appelle aussi~\textit{équivalence biLipschitz}) sont différentes. Il suit d'un résultat de Kevin Whyte (cf. Théorème~\ref{thm:Whytethm}) que les deux notions coïncident dans le monde non-moyennable, et pour le monde moyennable, la question était restée ouverte jusqu'en 2010 et un premier contre-exemple dû à Tullia Dymarz :

\begin{théorème}[{\cite[Théorème~1]{Dym10}}]\label{thm:INTROFRANCAISBiLipvsQIDymarz}
Soit $F$ un groupe fini non-trivial, et soit $k\ge 2$ un entier. Alors $F\wr\Z$ et $F^{k}\wr\Z$ ne sont pas biLipschitz équivalents si $k$ n'est pas un produit de nombres premiers apparaissant dans la décomposition de $|F|$. 
\end{théorème}

Le Théorème~\ref{thm:INTROFRANCAISclassificationGT21} fournit aussi les premiers exemples de groupes \textit{moyennables} ayant toutes leurs auto-quasi-isométries à distance bornée d'une bijection :

\begin{corollaire}[{\cite[Corollaire~1.16]{GT24b}}]
Soit $F$ un groupe fini non-trivial, et soit $K$ un groupe moyennable de présentation finie à un bout. Alors toute quasi-isométrie $F\wr K\longrightarrow F\wr K$ est à distance bornée d'une bijection.
\end{corollaire}

\paragraph{Produits en halo.} Dans un article plus récent~\cite{GT24a}, Genevois et Tessera ont poussé leurs méthodes encore plus loin et les ont généralisées au cadre des~\textit{produits en halo} :

\begin{définition}[{cf. Définition~\ref{def:haloproducts}}]
Soit $H$ un groupe. Un~\textit{produit en halo $\halo$ au-dessus de $H$} est la donnée, pour tout $S\subset H$, d'un groupe $L(S)$ tel que :
\begin{itemize}
    \item pour tous $R,S\subset H$, si $R\subset S$ alors $L(R)\leqslant L(S)$;
    \item $L(\emptyset)=\lbrace 1\rbrace$ et $L(H)=\langle L(S) : S\subset H \;\text{fini}\rangle$;
    \item pour tous $R,S\subset H$, $L(R\cap S)=L(R)\cap L(S)$.
\end{itemize}
\end{définition}
Etant donné une action $\alpha\colon H\curvearrowright L(H)$, le~\textit{produit en halo} $\halo H$ au-dessus de $H$ est alors le produit semi-direct
\begin{equation*}
    \halo H \defeq L(H)\rtimes_{\alpha} H.
\end{equation*}

Par exemple, si $F$ est un groupe fini non-trivial, la famille $L(S)=\bigoplus_{S}F$, $S\subset H$ donne le produit en halo $\halo H=F\wr H$, et si on choisit plutôt $L(S)=\fsym{S}$ le groupe des permutations à support fini $S\rightarrow S$, le produit en halo $\halo H$ est habituellement appelé le~\textit{lampshuffler} au-dessus de $H$ :
\begin{equation*}
    \shuf{H} \defeq \fsym{H}\rtimes H.
\end{equation*}
où $\fsym{H}$ désigne le groupe des permutations à support fini $H\rightarrow H$. On fournit une liste plus détaillée de constructions similaires dans la Section~\ref{sec:defHalo}. 

La motivation de~\cite{GT24a} pour introduire un tel formalisme vient du fait que, en réalité, le phénomène de rigidité pour les quasi-isométries entre produits en couronne observé dans~\cite{GT24b} (cf. Théorème~\ref{thm:INTROFRANCAISclassificationGT21}) est rendu possible par la structure de produit semi-direct, et le fait que le groupe de base satisfait certaines hypothèses de non-séparation par certains sous-espaces. On explique cette observation et les techniques mises en place plus en détail dans le Chapitre~\ref{chap:chapter2}. 

\paragraph{Fonction de F\o lner et profil isopérimétrique.} Il n'est pas difficile de voir qu'un produit en couronne $A\wr B$ est à croissance exponentielle dès que $B$ est infini et que $A$ est non-trivial. Pour pouvoir les distinguer à quasi-isométrie près, il faut considérer d'autres invariants. L'un d'entre eux, qui se révèle très utile en pratique, est la~\textit{fonction de F\o lner} :
\begin{définition}
Soit $G$ un groupe de type fini avec $S_{G}$ un ensemble générateur fini. La~\textit{fonction de F\o lner} de $G$ par rapport à $S_{G}$ est l'application
\begin{align*}
    \text{F\o l}_{G}\colon \N &\longrightarrow \R_{+} \\
    n&\longmapsto \inf\left\lbrace |A| : \frac{|\partial_{G}A|}{|A|}\le \frac{1}{n}\right\rbrace
\end{align*}
où la~\textit{frontière} d'un ensemble fini $A\subset G$ est définie par
\begin{equation*}
    \partial_{G}A\defeq \lbrace y\in G\setminus A : \exists s\in S_{G},\exists a\in A, y=as\rbrace. 
\end{equation*}
\end{définition}

La fonction de F\o lner d'un groupe $G$ doit être vue comme un moyen de « mesurer » sa moyennabilité : plus cette fonction tend vers l'infini rapidement, moins le groupe est moyennable. 

De façon équivalente, on peut aussi travailler avec son inverse généralisé, appelé le~\textit{profil isopérimétrique} de $G$, qui est la fonction $\prof{G}\colon \N\rightarrow\R_{+}$ définie par 
\begin{equation*}
    \prof{G}(n) \defeq \sup\left\lbrace \frac{|A|}{|\partial_{G}A|} : |A|\le n\right\rbrace. 
\end{equation*}

Une formule d'Erschler est particulièrement pratique pour calculer les fonctions de F\o lner des produits en couronne :

\begin{théorème}[{\cite[Théorème~1]{Ers03}}]\label{thm:INTROFRANCAISformulaforprofilewreathproducts}
Soient $G$ et $H$ des groupes moyennables de type fini. On suppose que la fonction de F\o lner de $H$ vérifie :
\begin{equation*}
    \forall C>0,\;\exists K>0,\; C\cdot\text{F\o l}_{H}(n) \le \text{F\o l}_{H}(Kn) \;\text{pour tout $n$ assez grand}.
\end{equation*}
Alors on a
\begin{equation*}
    \text{F\o l}_{G\wr H}(n) \simeq \text{F\o l}_{G}(n)^{\text{F\o l}_{H}(n)}.
\end{equation*}
\end{théorème}

On discute l'hypothèse dans cet énoncé plus en détail dans la Section~\ref{sec:assumptionstar}, puisqu'on y fait aussi largement appel pour calculer les fonctions de F\o lner des lampshufflers et autres produits en halo. 

\paragraph{\'Equivalence orbitale quantitative...} En plus d'être un invariant de quasi-isométrie, le profil isopérimétrique est aussi étroitement relié à une méthode de comparaison~\textit{mesurée} des groupes de type fini. 

Si $\Gamma$ et $\Lambda$ sont deux groupes de type fini, un~\textit{couplage d'équivalence orbitale} entre eux est la donnée d'un espace de probabilité standard $(X,\mu)$ et de deux actions libres et p.m.p. $\Gamma,\Lambda\curvearrowright (X,\mu)$ ayant les mêmes orbites :
\begin{equation*}
    \Lambda\cdot x = \Gamma\cdot x
\end{equation*}
pour $\mu-$presque tout $x\in X$. 

Un résultat célèbre d'Ornstein et Weiss affirme en fait qu'un tel couplage existe toujours entre deux groupes moyennables infinis. 

\begin{théorème}[{\cite{OW80}}]
Deux groupes infinis moyennables sont orbitalement équivalents. 
\end{théorème}

L'équivalence orbitale est dès lors trop faible pour distinguer deux groupes moyennables, et n'encode aucune information géométrique. Il est cependant possible de renforcer la définition en imposant des contraintes~\textit{quantitatives} sur les cocycles. 

En effet, un tel couplage d'équivalence orbitale vient avec deux applications $c_{\Gamma,\Lambda}\colon \Gamma\times X\rightarrow \Lambda$, $c_{\Lambda,\Gamma}\colon \Lambda\times X\rightarrow \Gamma$, définies de façon unique par les identités 
\begin{equation*}
    \gamma\cdot x = c_{\Gamma,\Lambda}(\gamma,x)\cdot x, \; \lambda\cdot x = c_{\Lambda,\Gamma}(\lambda,x)\cdot x
\end{equation*}
pour tous $\gamma\in\Gamma$, $\lambda\in\Lambda$ et $\mu-$presque tout $x\in X$. Ces applications sont appelées des \textit{cocycles} associés au couplage, puisque $c_{\Gamma,\Lambda}$ satisfait l'identité
\begin{equation*}
    c_{\Gamma,\Lambda}(\gamma\gamma',x)=c_{\Gamma,\Lambda}(\gamma, \gamma'\cdot x)c_{\Gamma,\Lambda}(\gamma',x), \; \gamma,\gamma'\in\Gamma, \;x\in X
\end{equation*}
et une identité similaire est vérifiée par $c_{\Lambda,\Gamma}$. Pour $p\ge 0$, on dit alors que $c_{\Gamma,\Lambda}$ est $\ld^{p}-$\textit{intégrable} si 
\begin{equation*}
    \int_{X} |c_{\Gamma,\Lambda}(\gamma,x)|_{S_{\Lambda}}^{p}\;\mathrm{d}\mu(x) < +\infty
\end{equation*}
pour tout $\gamma\in S_{\Gamma}$, et on dit qu'il est $\ld^{\infty}$ si pour tout $\gamma\in\Gamma$, l'application $|c_{\Gamma,\Lambda}(\gamma,\cdot)|_{S_{\Lambda}}\colon X\rightarrow\N$ est essentiellement bornée. On peut en fait définir la $\varphi-$intégrabilité d'un cocycle pour des fonctions $\varphi$ plus générales, cf. Définition~\ref{def:integrabilityofcocycles}. 

En plus, on a des résultats de rigidité, qui nous disent que l'existence d'un couplage d'équivalence orbitale entre deux groupes de type fini avec une certaine quantification impose une borne supérieure sur cette quantification. Cette borne peut souvent être exprimée en termes des profils isopérimétriques des groupes couplés :

\begin{théorème}[{\cite[Théorème~1.1]{DKLMT22}}]\label{thm:INTROFRANCAISObstructionDKLMT}
Soient $G$ et $H$ deux groupes moyennables de type fini. Soit $\varphi\colon \R_{+}\rightarrow\R_{+}$. Supposons qu'il existe un couplage $(\varphi, \ld^0)$-intégrable d'équivalence orbitale de $G$ vers $H$. 
\begin{enumerate}[label=(\roman*)]
    \item Si $\varphi$ et $t\longmapsto \frac{t}{\varphi(t)}$ sont croissantes, alors $\varphi\circ\prof{H}(n) \preccurlyeq \prof{G}(n)$.
    \item Si $\varphi(x)=x^p$ pour un $p\ge 1$, alors $\profp{H}(n)\preccurlyeq \profp{G}(n)$.
    \end{enumerate}
\end{théorème}

\paragraph{...entre produits en couronne.} Il est établi dans~\cite{DKLMT22} que la construction du produit en couronne est compatible avec l'équivalence orbitale quantitative, dans le sens suivant : 

\begin{théorème}[{\cite[Corollaire~7.3]{DKLMT22}}]\label{thm:INTROFRANCAISstabilityforOEbetweenlamplighters}
Soit $\Lambda$ un groupe fini. Soient $H$ et $K$ deux groupes de type fini et $\varphi,\psi\colon\R_{+}\rightarrow\R_{+}$ deux applications croissantes. Supposons qu'il existe un couplage $(\varphi,\psi)-$intégrable d'équivalence orbitale de $H$ vers $K$. Alors il existe un couplage $(\varphi,\psi)-$intégrable d'équivalence orbitale de $\Lambda\wr H$ vers $\Lambda\wr K$.
\end{théorème}

Le fait de connaître le profil isopérimétrique des produits en couronne (Théorème~\ref{thm:INTROFRANCAISformulaforprofilewreathproducts}) permet alors de déduire la caractérisation suivante :

\begin{corollaire}\label{cor:INTROFRANCAISquantitativecomparisonbetweenlamplighters}
Soient $k,\ell\ge 1$ des entiers. Soit $p\ge 0$. Soit $F$ un groupe fini. Alors il existe un couplage $(\ld^p, \ld^0)-$intégrable d'équivalence orbitale de $F\wr\Z^{k+\ell}$ vers $F\wr\Z^{k}$ si et seulement si $p<\frac{k}{k+\ell}$. 
\end{corollaire}

Ce résultat fournit une borne supérieure sur la $\ld^p-$intégrabilité d'un cocycle du plus gros groupe $\Z^{k+\ell}$ vers le plus petit $\Z^{k}$. Comme les cocycles mesurent la distorsion requise pour passer d'un graphe de Cayley d'un groupe à l'autre, ce dernier énoncé peut être vu comme une comparaison précise des géométries des deux groupes.

\bigskip
\noindent \textbf{Nouveaux résultats présentés dans cette thèse} 
\bigskip

Présentons maintenant les contributions de cette thèse. 

\paragraph{Rigidité quasi-isométrique des allumeurs de réverbères permutationnels.} Dans le Chapitre~\ref{chap:chapter3}, on établit une classification quasi-isométrique de certains produits en couronne permutationnels, ayant groupes de lampes finis. L'énoncé est le suivant. 

\begin{théorème}[{cf. Théorème~\ref{thm:classificationPWPuptoQI}}]\label{thm:INTROFRANCAISclassificationPWPuptoQI}
Soient $E$, $F$ des groupes finis non-triviaux. Soient $G$ et $H$ des groupes de présentation finie, avec des sous-groupes normaux infinis de type fini $M\lhd G$, $N\lhd H$. Supposons que $M$ a indice infini dans $G$, et $G$ (resp. $H$) n'est pas grossièrement séparé par une collection de sous-espaces qui se plonge uniformément quasi-isométriquement dans $M$ (resp. $N$). Les assertions suivantes sont vraies.
\begin{itemize}
    %\vspace{0.05cm}
    \item Si $M$ n'est pas co-moyennable dans $G$, alors $E\wr_{G/M}G$ et $F\wr_{H/N}H$ sont quasi-isométriques si et seulement si \;$|E|$ et $|F|$ ont les mêmes diviseurs premiers et il existe une quasi-isométrie de paires $(G,M)\longrightarrow (H,N)$. 
    \item Si $M$ est co-moyennable dans $G$, alors $E\wr_{G/M}G$ et $F\wr_{H/N}H$ sont quasi-isométriques si et seulement si \;$|E|=n^{r}$, $|F|=n^{s}$ pour certains $n,r,s\ge 1$ et il existe une quasi-isométrie de paires $(G,M)\longrightarrow (H,N)$ induisant une quasi-$\frac{s}{r}$-to-one quasi-isométrie $G/M\rightarrow H/N$.
\end{itemize}
\end{théorème}

Comme dans le Théorème~\ref{thm:INTROFRANCAISclassificationGT21}, les cas non-moyennable et moyennable doivent être distingués, le second étant plus rigide que le premier. 

Comme application, on obtient par exemple : 

\begin{corollaire}[{cf. Corollaire~\ref{cor:classificationofPWPoverZ^duptoQI}}]
Soient $E$, $F$ des groupes finis non-triviaux, et soient $m,m',n,n'$ quatre entiers tels que $m\ge n\ge 2$, $m'\ge n'\ge 2$. Alors $E\wr_{\Z^n}\Z^m$ et $F\wr_{\Z^{n'}}\Z^{m'}$ sont quasi-isométriques si et seulement si $m=m'$, $n=n'$ et $|E|$, $|F|$ sont puissances d'un même nombre. 
\end{corollaire}

\paragraph{Plus de scaling groups.} Comme expliqué dans la Section~\ref{sec:Scalinggroupsofhaloproducts}, la description générale des quasi-isométries entre deux produits en couronne permet de déduire qu'elles sont toutes « scaling », et de calculer précisément le « scaling group » d'allumeurs de réverbères. 

Il s'avère qu'une telle description est aussi accessible pour certains allumeurs de réverbères ayant des lampes infinies, comme montré dans~\cite{BGT24}. On utilise ce fait dans le Chapitre~\ref{chap:chapter4} pour établir la propriété de rigidité suivante :

\begin{théorème}[{cf. Théorème~\ref{thm:maintheoremfortheclassMexp}}]\label{thm:INTROFRANCAISmaintheoremfortheclassMexp}
Soient $N$ et $M$ des groupes de type fini à croissance polynomiale, avec degrés de croissance $n$ et $m$ respectivement. Soient $G$ et $H$ deux groupes moyennables de présentation finie dans $\mathcal{M}_{\text{exp}}$. Alors toute quasi-isométrie $N\wr G \longrightarrow M\wr H$ est quasi-$\frac{m}{n}$-to-one. 
\end{théorème}

On renvoie à la Section~\ref{sec:introchapter4} pour des détails sur ces hypothèses et la définition de la classe $\mathcal{M}_{\text{exp}}$. Pour l'instant, il suffit de remarquer que le « scaling factor » d'une quasi-isométrie entre deux produits en couronne avec lampes à croissance polynomiale peut être calculé en fonction des degrés de croissance impliqués. Cette observation, combinée au Théorème~\ref{thm:INTROFRANCAISclassificationGT21}, permet de déduire une contrainte arithmétique qui permet de distinguer certains produits en couronne itérés à quasi-isométrie près : 

\begin{proposition}[{cf.  Proposition~\ref{prop:mixingofscalingconditions}}]\label{prop:INTROFRANCAISiteratedlamplighters}
Soient $n,m\ge 2$ deux entiers. Soient $N_{1}$ et $N_{2}$ deux groupes de type fini à croissance polynomiale, avec degrés de croissance $n_{1}$ et $n_{2}$ respectivement. Soient $G$ et $H$ deux groupes moyennables de présentation finie dans $\mathcal{M}_{\text{exp}}$. Si \;$\Z_{n}\wr(N_{1}\wr G)$ et $\Z_{m}\wr(N_{2}\wr H)$ sont quasi-isométriques, alors il existe $a,r,s\ge 1$ tels que $n=a^{r}$, $m=a^{s}$, et
\begin{equation*}
    \frac{s}{r}=\frac{n_{2}}{n_{1}}. 
\end{equation*}
\end{proposition}

Par exemple, ce résultat permet de déduire qu'il n'existe pas de quasi-isométrie
\begin{equation*}
\Z_{2}\wr(\Z^{2}\wr \text{SOL}(\Z)) \longrightarrow \Z_{4}\wr(\Z^{3}\wr \text{SOL}(\Z)).
\end{equation*}
De plus, on obtient une classification complète quand $N_{1}$ et $N_{2}$ sont virtuellement abéliens, cf. Corollaire~\ref{cor:classificationofiteratedwreathproducts}.

On présente maintenant, dans le reste de cette introduction, deux travaux en commun avec Corentin Correia~\cite{cordum25, cordum26}, dans lesquels on montre que la structure de produit en halo est aussi particulièrement adaptée pour calculer le profil isopérimétrique et construire des couplages d'équivalence orbitale. 

\paragraph{Profils isopérimétriques des lampshufflers.} Le Chapitre~\ref{chap:chapter5} contient des résultats pour beaucoup d'exemples de produits en halo, mais on se restreint ici au cas des lampshufflers, pour simplifier la présentation. 

Pour de tels groupes, Erschler et Zheng ont prouvé une borne inférieure sur leur fonction de F\o lner :

\begin{théorème}[{\cite[Corollaire~1.4]{EZ21}}]
Soit $H$ un groupe moyennable de type fini. Alors on a
\begin{equation*}
    \text{F\o l}_{\shuf{H}}(x)\succcurlyeq \beta_{H}(x)^{\beta_{H}(x)}
\end{equation*}
où $\beta_{H}$ est la fonction de croissance de $H$. De plus, cette borne est optimale lorsque $H$ est à croissance polynomiale. 
\end{théorème}

En particulier, quand $H$ est à croissance exponentielle, on obtient la borne supérieure
\begin{equation*}
    \prof{\shuf{H}}(x) \preccurlyeq \ln(\ln(x))
\end{equation*}
sur le profil de $\shuf{H}$. Cette borne ne peut pas être optimale dans beaucoup de cas, puisqu'elle ne dépend que du fait que le groupe de base est à croissance exponentielle. Par exemple, $F\wr\Z$ et $F\wr(F\wr\Z)$ sont tous deux à croissance exponentielle, mais on s'attend à ce que le profil de $\shuf{F\wr(F\wr\Z)}$ soit plus lent que celui de $\shuf{F\wr\Z}$.

Un des principaux résultats de~\cite{cordum25} fournit un encadrement plus fin du profil de $\shuf{H}$, en fonction du profil de $H$ :

\begin{théorème}[{cf. Théorème~\ref{thm:boundsForProfile intro}}]\label{thm:INTROFRANCAISboundsForProfile intro}
Soit $p\ge 1$. Soit $H$ un groupe moyennable de type fini dont le profil $\profp{H}$ satisfait l'hypothèse~$(\star)$. Alors,
le profil $\profp{\shuf{H}}$ de $\shuf{H}$ satisfait
\begin{equation*}
    \profp{H}\left(\frac{\ln(x)}{\ln(\ln(x))}\right) \preccurlyeq \profp{\shuf{H}}(x) \preccurlyeq \prof{H}(\ln(x)). 
\end{equation*}
\end{théorème}

Dans cet énoncé, l'hypothèse~$(\star)$ est celle apparaissant dans le Théorème~\ref{thm:INTROFRANCAISformulaforprofilewreathproducts}, et en effet, la stratégie de preuve consiste à trouver, dans un lampshuffler, de bonnes copies de produits en couronne et à appliquer le théorème de monotonie du profil $\ell^p$, combiné avec le Théorème~\ref{thm:INTROFRANCAISformulaforprofilewreathproducts}. 

Du Théorème~\ref{thm:INTROFRANCAISboundsForProfile intro}, il est aussi possible de dériver des estimations pour le profil $\ell^p$ des lampshufflers itérés, définis inductivement par $\shufn{0}{H}\defeq H$ et 
\begin{equation*}
    \shufn{n}{H}\defeq \shuf{\shufn{n-1}{H}}
\end{equation*}
lorsque $n\ge 1$. Par exemple, pour des groupes de base à croissance polynomiale, on a : 

\begin{proposition}[{cf. Proposition~\ref{prop:profileofshufnofpolynomialgrowthgroupsINTRO}}]\label{prop:INTROFRANCAISprofileofshufnofpolynomialgrowthgroupsINTRO}
Soit $H$ un groupe de type fini à croissance polynomiale de degré $d\ge 1$. Alors on a
\begin{equation*}
    \profp{\shufn{n}{H}}(x) \simeq \left(\frac{\ln^{\circ n}(x)}{\ln^{\circ (n+1)}(x)}\right)^{\frac{1}{d}}
\end{equation*}
pour tout entier $n\ge 1$ et tout nombre réel $p\ge 1$.
\end{proposition}

Ces estimées ont aussi des conséquences sur le programme de classification commencé dans~\cite{GT24a}, et on peut aussi prouver que, dans certains cas, être quasi-isométrique revient à être bijectivement quasi-isométrique. Par exemple : 

\begin{corollaire}[{cf. Théorème~\ref{thm:IteratedShufflersPolynomialQIBiLip intro}}]\label{cor:INTROFRANCAISiteratedlampshufflersoverZ^dQI}
Soient $d,k\ge 1$ et $n,m\ge 0$ des entiers. Alors les lampshufflers itérés $\shufn{n}{\Z^d}$ et $\shufn{m}{\Z^k}$ sont quasi-isométriques si et seulement s'ils sont biLipschitz équivalents, ce qui arrive si et seulement si $n=m$ et $d=k$. 
\end{corollaire}

\paragraph{\'Equivalence orbitale quantitative pour les lampshufflers.} Lorsque deux groupes ne sont pas quasi-isométriques, l'équivalence orbitale quantitative est appropriée pour comparer un peu plus précisément leurs géométries. Pour les allumeurs de réverbères, l'existence de couplages d'équivalence orbitale quantitativement optimaux a été établie dans~\cite{DKLMT22} et~\cite{Cor25}. Dans~\cite{cordum26}, on entreprend une étude similaire pour d'autres produits en halo. 

Par exemple, pour les lampshufflers, on établit que :

\begin{théorème}[{cf. Théorème~\ref{thm:stabilityofcouplings+quantificationLampjugglersINTRO}}]
Si deux groupes $H$ et $K$ sont orbitalement équivalents, alors $\shuf{H}$ et $\shuf{K}$ sont orbitalement équivalents. De plus, si $H$ et $K$ sont de type fini, si $\varphi,\psi\colon\R_{+}\rightarrow\R_{+}$ sont croissantes, et s'il existe un couplage $(\varphi,\psi)-$intégrable d'équivalence orbitale de $H$ vers $K$, alors il existe un tel couplage de $\shuf{H}$ vers $\shuf{K}$.
\end{théorème}

Cet énoncé doit être mis en relation avec le Théorème~\ref{thm:INTROFRANCAISstabilityforOEbetweenlamplighters}, et en fait on extrait de la preuve de ce résultat une méthode générale qui peut, en pratique, être appliquée pour tout produit en halo satisfaisant des hypothèses peu contraignantes. On renvoie à la Section~\ref{sec:generalmethodstability} pour plus de détails. 

Couplé avec nos calculs de profils isopérimétriques de lampshufflers itérés (cf. Proposition~\ref{prop:INTROFRANCAISprofileofshufnofpolynomialgrowthgroupsINTRO}), on peut déduire des couplages quantitativement optimaux entre lampshufflers itérés au-dessus de groupes abéliens libres :

\begin{corollaire}[{cf. Théorème~\ref{thm:optimalitypolynomialgrowthJugglersINTRO}}]
Soient $k,d\ge 1$ des entiers tels que $k>d$. Soit $n\ge 0$ un entier.
Alors $\shufn{n}{\Z^k}$ et $\shufn{n}{\Z^d}$ sont $\ld^{p}$ orbitalement équivalents si et seulement si $p<\frac{d}{k}$.
\end{corollaire}

\bigskip

\noindent \textbf{Structure de la thèse} 

\bigskip

Le Chapitre~\ref{chap:chapter1} contient le matériel et les prérequis nécessaires de géométrie des groupes pour la lecture des chapitres suivants. Le lecteur familier avec le domaine peut sans autre s'en exonérer et aborder le Chapitre~\ref{chap:chapter2}, dont le but est de présenter, dans les grandes lignes, le contenu de~\cite{GT24a} et~\cite{GT24b}, puisque la première partie de la thèse repose grandement sur le contenu de ces papiers. Le Chapitre~\ref{chap:chapter3} présente l'article~\cite{Dum24} et établit en particulier le Théorème~\ref{thm:INTROFRANCAISclassificationPWPuptoQI}. Le Chapitre~\ref{chap:chapter4} présente les contributions de~\cite{Dum26} concernant les allumeurs de réverbères à lampes infinies (notamment le Théorème~\ref{thm:INTROFRANCAISmaintheoremfortheclassMexp} et la Proposition~\ref{prop:INTROFRANCAISiteratedlamplighters}). On étudie les profils isopérimétriques de produits en halo dans le Chapitre~\ref{chap:chapter5}, qui contient aussi nos applications au problème d'existence de quasi-isométries et plongements réguliers entre tels produits en halo (par exemple le Corollaire~\ref{cor:INTROFRANCAISiteratedlampshufflersoverZ^dQI}). Nos résultats sur l'équivalence orbitale quantitative sont présentés dans le Chapitre~\ref{chap:chapter6} et reposent sur le contenu de~\cite{cordum26}.

%\chapter*{Introduction}
\fancyhead{} % clear all header fields
\fancyhead[OL]{\textsc{Introduction in english}}
\chapter*{Introduction in English}
\addcontentsline{toc}{chapter}{Introduction}  
\mtcaddchapter

If $(X,d_{X})$ and $(Y,d_{Y})$ are two metric spaces, a~\textit{quasi-isometry} between $X$ and $Y$ is a map $f\colon X\rightarrow Y$ for which there are two constants $C\ge 1$, $K\ge 0$ such that 
\begin{equation*}
    \frac{1}{C}\cdot d_{X}(x,y)-K \le d_{Y}(f(x), f(y)) \le C\cdot d_{X}(x,y)+K
\end{equation*}
for all $x,y\in X$, and $d_{Y}(y, f(X))\le K$ for any $y\in Y$. When such a map exists, we say that $X$ and $Y$ are~\textit{quasi-isometric}. It is in general a hard problem to determine whether two given metric spaces $X$ and $Y$ are quasi-isometric and, when they are, it is an even harder problem to reach a complete description of all quasi-isometries that exist between $X$ and $Y$. To address such questions, we therefore look at classes of metric spaces carrying an additional structure.

One such class is that of~\textit{finitely generated groups}, which are naturally endowed with their word metrics. In his seminal work in the 1980s, Gromov started the ambitious program of trying to classify finitely generated groups up to quasi-isometry. This program is motivated by his celebrated result on groups of polynomial growth, which is the most striking illustration of the fact that the large-scale geometry of a group, defined only from a chosen word metric, is in fact intimately related to its algebraic structure:

\begin{theorem}[\cite{Gro81}]
Let $G$ be a finitely generated group. Then $G$ has polynomial growth if and only if it contains a finite-index nilpotent subgroup. 
\end{theorem}

More generally, to prove that two groups are~\textit{not} quasi-isometric, it is often fruitful to seek~\textit{quasi-isometry invariants}, namely a property $\mathcal{P}$ which is invariant under quasi-isometry (i.e. if $G$ has $\mathcal{P}$ and if there is a quasi-isometry $G\rightarrow H$, then $H$ has $\mathcal{P}$). Then, if one of our groups has $\mathcal{P}$ but the other does not, then we can rule out the existence of a quasi-isometry between them.

The~\textit{growth} of a group is, perhaps, the most elementary example of a quasi-isometry invariant; it counts the number of elements in a ball of a given radius centered at $1_{G}$ when $G$ is endowed with a word-metric; see Section~\ref{subsubsectionGROWTH} for details. Since having polynomial growth is a quasi-isometry invariant, it follows from Gromov's result that being virtually nilpotent is a geometric property. 

The thesis explores in detail another family of invariants, usually finer than the volume growth, and referred to as the $\ell^{p}-$isoperimetric profiles. Such profiles have a more elaborate definition and therefore are harder to compute than the volume growth, but, in return, they also encode more information about the large-scale geometry of the group. See Section~\ref{sec:isoprof} for details. 

Another direction of research is to find groups that are~\textit{quasi-isometrically rigid}. A group $G$ is \textit{quasi-isometrically rigid} if, for any group $G'$ quasi-isometric to $G$, $G'$ is either isomorphic to a finite-index subgroup of $G$ or to a quotient of $G$ by a finite normal subgroup. For instance, free abelian~\cite{Sha04} and free non-abelian groups are quasi-isometrically rigid~\cite{Sta68, Dun85}, as well as mapping class groups~\cite{Beh+12} or solvable Baumslag-Solitar groups~\cite{FM98}. There are also a wide range of classes of groups that are conjectured to be quasi-isometrically rigid, such as polycyclic groups, right-angled Artin groups and outer automorphism groups of free groups. 

Going even further, one may also try to describe all groups that are quasi-isometric to a given metric space. In this direction, Pekka Tukia proved the following: 

\begin{theorem}[{\cite{Tuk86, Tuk94}}]
Let $n\ge 2$. If $G$ is a finitely generated group quasi-isometric to $\mathbb{H}^{n+1}$, then $G$ acts geometrically on $\mathbb{H}^{n+1}$. In particular, $G$ is virtually isomorphic to a uniform lattice in $\text{Isom}(\mathbb{H}^{n+1})$. 
\end{theorem}

Tukia's work also led to impressive quasi-isometric rigidity results for certain solvable Lie groups; see~\cite{DFX23, DFX25}. 

Lastly, the literature also shows that, even in a specific class, tracking quasi-isometries that relate two groups is hard, and requires involved techniques. A good illustrative instance is the work of Farb and Mosher on solvable Baumslag-Solitar groups, who managed to prove the following classification criterion: 

\begin{theorem}[{\cite[Theorem~7.1]{FM98}}]
Let $n,m\ge 2$. Then $\text{BS}(1,n)$ and $\text{BS}(1,m)$ are quasi-isometric if and only if $n$ and $m$ are powers of a common integer. 
\end{theorem}

This result in particular shows that the large-scale geometries of Baumslag-Solitar groups heavily depend on the parameters that define them. 

So far, it appears that quasi-isometric rigidity problems have been investigated in the literature for two main classes of groups: those that exhibit non-positive curvature features, and those that are amenable (and often even more, for instance solvable). One goal of this thesis is to explore such problems in a third direction: the class of~\textit{wreath products} and all their related variants, that we present now. 

\paragraph{Our fellow partners for the thesis.} Given two finitely generated groups $A$ and $B$, their~\textit{wreath product} is the finitely generated group defined by 
\begin{equation*}
    A\wr B \defeq \left(\bigoplus_{B}A\right)\rtimes B
\end{equation*}
where $B$ acts on the direct sum by permutation of the coordinates. We usually call $A$ the~\textit{lamp group}, and $B$ the~\textit{base group}, mostly because of a nice and visual description of elements of a wreath product in terms of~\textit{colourings of graphs} and~\textit{arrows} moving along this graph.

Indeed, consider a group element $(c,p)\in \Z/2\Z\wr\Z$. Think of $c\in\bigoplus_{\Z}\Z/2\Z$ as a colouring of the vertices of the canonical Cayley graph of $\Z$, which looks like an infinite street, with only two possible colours: $0$ or $1$. Think of $p\in\Z$ as being the position of a lamplighter. As explained below, a natural generating set $U$ of $\Z/2\Z\wr\Z$ can be made of two parts, one part corresponding to changing the colourings, and the other part corresponding to changing the lamplighter position. Then, there are two possible moves to go from $(c,p)$ to a neighbouring vertex in $\text{Cay}(\Z/2\Z\wr\Z,U)$: 
\begin{itemize}
    \item either we only move the arrow to a neighbouring vertex in $\Z$, and the colouring stays the same;
    \item or the arrow stays on the vertex where it stands, but changes the colour of this vertex, that goes to $1$ if it was $0$, and to $0$ if it was $1$. 
\end{itemize}

See Example~\ref{ex:descriptionoflamplighters} for a more general explanation. This description led some authors to use the word~\textit{lamplighter} to refer to a wreath product having a~\textit{finite} lamp group. We follow this convention throughout the entire text. 

In the light of Gromov's program, one may ask:

\begin{question}\label{questionGromov}
Let $A_{1}$, $A_{2}$, $B_{1}$, $B_{2}$ be finitely generated groups. When are $A_{1}\wr B_{1}$ and $A_{2}\wr B_{2}$ quasi-isometric? 
\end{question}

A first piece of answer was provided by Anna Erschler in 2000, when she proved that the wreath product construction is compatible with~\textit{bijective} quasi-isometries (i.e. biLipschitz equivalences).

\begin{theorem}[\cite{Dyu00}]
Let $A$, $B$ and $C$ be finitely generated groups. If $A$ and $B$ are biLipschitz equivalent, then $A\wr C$ and $B\wr C$ are biLipschitz equivalent. 
\end{theorem}

Using this stability result, Erschler managed to prove that virtual solvability and virtual torsion-freeness are not quasi-isometry invariants, a question that was left open since Gromov's result on polynomial growth groups; see~\cite[Proposition~1]{Dyu00}. Let us notice here that it is still an open question whether virtual solvability is a quasi-isometry invariant among finitely presented groups. 

\paragraph{Quasi-isometric classifications.} In 2013, Alex Eskin, David Fisher and Kevin Whyte managed to provide a complete answer to Question~\ref{questionGromov} for lamplighters over cyclic groups (that are amenable groups).

\begin{theorem}[{\cite[Theorem~1.2]{EFW13}}]
Let $F_{1}$, $F_{2}$ be non-trivial finite groups. Then $F_{1}\wr\Z$ and $F_{2}\wr\Z$ are quasi-isometric if and only if there exist integers $a,r,s\ge 1$ such that $|F_{1}|=a^{r}$ and $|F_{2}|=a^{s}$. 
\end{theorem}

In particular, this result is the first one to exhibit a dependency between the presence of a quasi-isometry between two lamplighters and the cardinalities of the involved lamp groups.

A second major contribution to Question~\ref{questionGromov} was recently made by Anthony Genevois and Romain Tessera: 

\begin{theorem}[{\cite[Corollary~1.6]{GT24b}}]\label{thm:INTROclassificationGT21}
Let $F_{1}$ and $F_{2}$ be non-trivial finite groups. Let $H_{1}$ and $H_{2}$ be finitely presented groups. Suppose that $H_{1}$ is one-ended.
\begin{enumerate}[label=(\roman*)]
    \item If $H_{1}$ is not amenable, then $F_{1}\wr H_{1}$ and $F_{2}\wr H_{2}$ are quasi-isometric if and only if $|F_{1}|$, $|F_{2}|$ have the same prime divisors and there exists a quasi-isometry $H_{1}\rightarrow H_{2}$.
    \item If $H_{1}$ is amenable, then $F_{1}\wr H_{1}$ and $F_{2}\wr H_{2}$ are quasi-isometric if and only if there exist integers $a,r,s\ge 1$ such that $|F_{1}|=a^{r}$, $|F_{2}|=a^{s}$ and there exists a quasi-$\frac{s}{r}$-to-one quasi-isometry $H_{1}\rightarrow H_{2}$.
\end{enumerate}
\end{theorem}

We notice two important points here: firstly, this result applies to base groups that are finitely presented and one-ended, which already covers a large number of new cases, such properties being encountered in many classes of infinite finitely generated groups. Secondly, it makes even more explicit the dependency between the quasi-isometry between the base groups and the cardinalities of the lamp groups, appealing, for the amenable case, to the notion of~\textit{scaling quasi-isometries} (which is trivial if groups are not amenable). The authors of~\cite{GT24b} have formally introduced the concept in an earlier work~\cite{GT22}, that we explain in detail in Section~\ref{sec:ScalingQI}. 

There exists in fact a more general version of Theorem~\ref{thm:INTROclassificationGT21}, with weaker assumptions, proved in the more recent paper~\cite{GT24a}. We state and explain this result in details in Chapter~\ref{chap:chapter2}, but for the purpose of this introduction, Theorem~\ref{thm:INTROclassificationGT21} will be sufficient. 

\paragraph{Quasi-isometries vs biLipschitz equivalences.} An important direction of research that emerged is the understanding of the difference between being quasi-isometric and being biLipschitz equivalent. Following Gromov's work, it was indeed an open question to know whether there are amenable quasi-isometric groups that are not biLipschitz equivalent. This problem has been solved by Tullia Dymarz in 2010, using lamplighters over $\Z$: 

\begin{theorem}[{\cite[Theorem~1]{Dym10}}]\label{thm:INTROBiLipvsQIDymarz}
Let $F$ be a non-trivial finite group, and let $k\ge 2$ be an integer. Then $F\wr\Z$ and $F^{k}\wr\Z$ are not biLipschitz equivalent if $k$ is not a product of the prime numbers appearing in the decomposition of $|F|$. 
\end{theorem}

Between non-amenable groups, it is known since 1999 and a famous result of Kevin Whyte that any quasi-isometry can be turned into a biLipschitz equivalence, up to a finite distance change; see Theorem~\ref{thm:Whytethm} and the references therein. In particular, a non-amenable group has all its self-quasi-isometries lying at bounded distance from bijections.

On the other hand, as a nice consequence, one can deduce from Theorem~\ref{thm:INTROclassificationGT21} the first examples of \textit{amenable} groups having all their self-quasi-isometries at bounded distance from bijections. 

\begin{corollary}[{\cite[Corollary~1.16]{GT24b}}]
Let $F$ be a non-trivial finite group, and let $K$ be an amenable one-ended finitely presented group. Then any quasi-isometry $F\wr K\longrightarrow F\wr K$ is at bounded distance from a bijection.
\end{corollary}

This result is in fact a straightforward consequence of a result appearing implicitly in Theorem~\ref{thm:INTROclassificationGT21}: any quasi-isometry between two lamplighters (with our running assumptions) must actually be scaling. In particular, the scaling group of $F\wr K$ must be reduced to $\lbrace 1\rbrace$.

\paragraph{Halo products.} In a subsequent article~\cite{GT24a}, Genevois and Tessera have pushed their techniques even further, generalizing them to a more general framework, that of~\textit{halo products}:

\begin{definition}[{see Definition~\ref{def:haloproducts}}]
Let $H$ be a group. A~\textit{halo of groups $\halo$ over $H$} is the data, for any subset $S\subset H$, of a group $L(S)$ such that:
\begin{itemize}
    \item for all $R,S\subset H$, if $R\subset S$ then $L(R)\leqslant L(S)$;
    \item $L(\emptyset)=\lbrace 1\rbrace$ and $L(H)=\langle L(S) : S\subset H \;\text{finite}\rangle$;
    \item for all $R,S\subset H$, $L(R\cap S)=L(R)\cap L(S)$.
\end{itemize}
\end{definition}
Given an action $\alpha\colon H\curvearrowright L(H)$, we then define the~\textit{halo product} $\halo H$ over $H$ as the semi-direct product
\begin{equation*}
    \halo H \defeq L(H)\rtimes_{\alpha} H.
\end{equation*}

For instance, if $F$ is a non-trivial group, choosing $L(S)=\bigoplus_{S}F$ for all $S\subset H$ yields $\halo H = F\wr H$, while choosing $L(S)=\fsym{S}$ the group of finitely supported permutations $S\rightarrow S$ allows to define the~\textit{lampshuffler over $H$}:
\begin{equation*}
    \shuf{H} \defeq \fsym{H}\rtimes H
\end{equation*}
where $\fsym{H}$ stands for the group of finitely supported permutations $H\rightarrow H$. Many other constructions are possible; we refer to Section~\ref{sec:halo} for a richer list. 

The initial motivation to introduce such a formalism is to notice that the semi-direct product structure is in fact the key ingredient for proving Theorem~\ref{thm:INTROclassificationGT21}, as left cosets of the base groups of lamplighters must be quasi-preserved by quasi-isometries. A similar approach therefore yields many other classification results for lampshufflers, lampcloners or lampdesigners. As an illustration: 

\begin{theorem}[{\cite[Corollary~8.9]{GT24a}}]
Let $H$ and $K$ be finitely presented one-ended groups. Then $\shuf{H}$ and $\shuf{K}$ are quasi-isometric if and only if $H$ and $K$ are biLipschitz equivalent. 
\end{theorem}

We refer to Chapter~\ref{chap:chapter2} for a detailed explanation of the techniques and the results exposed in~\cite{GT24a}.

\paragraph{F\o lner functions and isoperimetric profiles.} Let us now turn to a more analytic study of the large-scale geometry of wreath products. Such groups have exponential growth as soon as the base group is infinite. We therefore need more subtle geometric invariants to distinguish wreath products from one another. One such invariant, that carries more information than the volume growth, is the~\textit{F\o lner function}:

\begin{definition}
Let $G$ be a finitely generated group with $S_{G}$ a finite generating set. The~\textit{F\o lner function} of $G$ with respect to $S_{G}$ is the map
\begin{align*}
    \text{F\o l}_{G}\colon \N &\longrightarrow \R_{+} \\
    n&\longmapsto \inf\left\lbrace |A| : \frac{|\partial_{G}A|}{|A|}\le \frac{1}{n}\right\rbrace
\end{align*}
where the~\textit{boundary} of a finite subset $A\subset G$ is given by 
\begin{equation*}
    \partial_{G}A\defeq \lbrace y\in G\setminus A : \exists s\in S_{G},\exists a\in A, y=as\rbrace. 
\end{equation*}
\end{definition}

Intuitively, the F\o lner function of a group indicates how amenable it is; the faster it tends to infinity, the less the group is amenable. 

Rather than the F\o lner function, it is also possible to work with its generalized inverse, called the~\textit{isoperimetric profile}, which is the function $\prof{G}\colon\N\rightarrow \R_{+}$ defined by
\begin{equation*}
    \prof{G}(n) \defeq \sup\left\lbrace \frac{|A|}{|\partial_{G}A|} : |A|\le n\right\rbrace. 
\end{equation*}
We refer to Section~\ref{sec:isoprof} for further details on the isoperimetric profile and its more general $\ell^p-$version. At this point, it will be enough for us to stick to the following result, that gives a precise formula to compute F\o lner functions of wreath products in terms of F\o lner functions of the factors. 

\begin{theorem}[{\cite[Theorem~1]{Ers03}}]\label{thm:INTROformulaforprofilewreathproducts}
Let $G$ and $H$ be finitely generated amenable groups. Suppose that the F\o lner function of $H$ satisfies:
\begin{equation*}
    \forall C>0,\;\exists K>0,\; C\cdot\text{F\o l}_{H}(n) \le \text{F\o l}_{H}(Kn) \;\text{for all $n$ large enough}.
\end{equation*}
Then one has 
\begin{equation*}
    \text{F\o l}_{G\wr H}(n) \simeq \text{F\o l}_{G}(n)^{\text{F\o l}_{H}(n)}.
\end{equation*}
\end{theorem}

We discuss further the assumption of Theorem~\ref{thm:INTROformulaforprofilewreathproducts} in Section~\ref{sec:assumptionstar}. Let us simply mention that, in practice, it holds for any finitely generated amenable group for which we know the F\o lner function. 

\paragraph{Quantitative orbit equivalence...} In addition to being a rather thin quasi-isometry invariant, isoperimetric profiles are in fact tightly connected to a~\textit{measured} way of comparing the large-scale geometries of finitely generated groups. 

Given finitely generated groups $\Gamma$ and $\Lambda$, an~\textit{orbit equivalence coupling} between them is the data of a standard probability space $(X,\mu)$ and two free p.m.p. actions $\Gamma,\Lambda\act (X,\mu)$ sharing the same orbits: 
\begin{equation*}
    \Lambda\cdot x = \Gamma\cdot x
\end{equation*}
for $\mu-$almost every $x\in X$. 

A famous result of Ornstein-Weiss says that, in fact, such couplings always exist for amenable groups. 

\begin{theorem}[{\cite{OW80}}]
Any two infinite amenable groups are orbit equivalent. 
\end{theorem}

Thus, this form of orbit equivalence is too weak to distinguish amenable groups from a geometric point of view. It is therefore necessary to strengthen the definition of orbit equivalence, by imposing constraints on maps that derive from the coupling. 

Indeed, an orbit equivalence coupling naturally gives rise to two maps $c_{\Gamma,\Lambda}\colon \Gamma\times X\rightarrow \Lambda$, $c_{\Lambda,\Gamma}\colon \Lambda\times X\rightarrow \Gamma$, uniquely defined by the equations
\begin{equation*}
    \gamma\cdot x = c_{\Gamma,\Lambda}(\gamma,x)\cdot x, \; \lambda\cdot x = c_{\Lambda,\Gamma}(\lambda,x)\cdot x
\end{equation*}
for all $\gamma\in\Gamma$, $\lambda\in\Lambda$ and $\mu-$almost every $x\in X$. These maps are called the~\textit{cocycles} associated to the coupling, because they satisfy the identity
\begin{equation*}
    c_{\Gamma,\Lambda}(\gamma\gamma',x)=c_{\Gamma,\Lambda}(\gamma, \gamma'\cdot x)c_{\Gamma,\Lambda}(\gamma',x), \; \gamma\in\Gamma, \;\lambda\in\Lambda, \;x\in X
\end{equation*}
and similarly for $c_{\Lambda,\Gamma}$. 

As it turns out, these maps contain valuable information about geometries of $\Gamma$ and $\Lambda$. Indeed, orbits of a free action of $\Gamma$ (resp. $\Lambda$) on a probability space look like a Cayley graph of $\Gamma$ (resp. $\Lambda)$ with respect to a fixed generating set $S_{\Gamma}$ (resp. $S_{\Lambda}$). In every orbit, we put an edge between $x$ and $\gamma\cdot x$ (resp. $x$ and $\lambda\cdot x$) if $\gamma\in S_{\Gamma}\cup S_{\Gamma}^{-1}$ (resp. $\lambda\in S_{\Lambda}\cup S_{\Lambda}^{-1}$). Set-wisely, $\Gamma\cdot x$ equals $\Lambda \cdot x$, and the cocycles encode the distortion needed to transform the edges of one orbit into those of the other. This distortion is measured by computing $|c_{\Gamma,\Lambda}(\gamma,x)|_{S_{\Lambda}}$ and $|c_{\Lambda,\Gamma}(\lambda,x)|_{S_{\Gamma}}$, for generators $\gamma\in S_{\Gamma}$, $\lambda\in S_{\Lambda}$. The bigger these values are, the bigger is the distortion from one orbit to the other. This leads us to define~\textit{$\ld^p-$integrability} of cocycles: for $p\ge 0$, we say that $c_{\Gamma,\Lambda}$ is $\ld^p$ if 
\begin{equation*}
    \int_{X} |c_{\Gamma,\Lambda}(\gamma,x)|_{S_{\Lambda}}^{p}\;\mathrm{d}\mu(x) < +\infty
\end{equation*}
for every $\gamma\in S_{\Gamma}$, and we say it is $\ld^{\infty}$ if for any $\gamma\in\Gamma$, the map $|c_{\Gamma,\Lambda}(\gamma,\cdot)|_{S_{\Lambda}}\colon X\rightarrow\N$ is essentially bounded. We make similar definitions for $c_{\Lambda,\Gamma}$. Then, for $p,q\ge 0$, an $(\ld^p,\ld^q)-$coupling from $\Gamma$ to $\Lambda$ is therefore a coupling $(X,\mu)$ for which $c_{\Gamma,\Lambda}$ is $\ld^p$ and $c_{\Lambda,\Gamma}$ is $\ld^q$. It is possible to consider similar integrability conditions for more general functions, which lead to the notion of $(\varphi,\psi)-$integrable orbit equivalence couplings; see Definition~\ref{def:integrabilityofcocycles}.

In addition to existence results, we need~\textit{rigidity} statements, i.e. statements saying that the existence of a quantified orbit equivalence coupling between two finitely generated groups ensures an upper bound on the possible quantifications. This is actually the main observation of~\cite{DKLMT22}, namely the quantified version of orbit equivalence is suitable for detecting geometric differences between amenable groups. The most important obstruction implies $\ell^p-$isoperimetric profiles. 

\begin{theorem}[{\cite[Theorem~1.1]{DKLMT22}}]\label{thm:INTROObstructionDKLMT}
Let $G$ and $H$ be finitely generated amenable groups. Let $\varphi\colon \R_{+}\rightarrow\R_{+}$. Assume that there is a $(\varphi, \ld^0)$-integrable orbit equivalence coupling from $G$ to $H$. 
\begin{enumerate}[label=(\roman*)]
    \item If $\varphi$ and $t\longmapsto \frac{t}{\varphi(t)}$ are increasing, then $\varphi\circ\prof{H}(n) \preccurlyeq \prof{G}(n)$.
    \item If $\varphi(x)=x^p$ for some $p\ge 1$, then $\profp{H}(n)\preccurlyeq \profp{G}(n)$.
    \end{enumerate}
\end{theorem}

We refer to~\cite[Theorem~3.1]{DKLMT22} for another obstruction implying rather the volume growth. 

\begin{remark}
Orbit equivalence is in fact a particular instance of~\textit{measured equivalence}, introduced by Gromov as a measured analog of the notion of quasi-isometry. There are important connections between the two points of view: for instance, two amenable groups are biLipschitz equivalent if and only if they are $\ld^{\infty}$ orbit equivalent~\cite{Sha04}, and they are quasi-isometric if and only if there exists a~\textit{mutually cobounded} $\ld^{\infty}$ measure equivalence between them. This thesis however only presents and deals with orbit equivalence.
\end{remark}

\paragraph{...between wreath products.} It is proved in~\cite{DKLMT22} that the wreath product construction is also compatible with quantitative orbit equivalence, in the following sense: 

\begin{theorem}[{\cite[Corollary~7.3]{DKLMT22}}]\label{thm:INTROstabilityforOEbetweenlamplighters}
Let $\Lambda$ be a finite group. Let $H$ and $K$ be finitely generated groups and $\varphi,\psi\colon\R_{+}\rightarrow\R_{+}$ be increasing maps. If there exists a $(\varphi,\psi)-$integrable orbit equivalence coupling from $H$ to $K$, then the same holds from $\Lambda\wr H$ to $\Lambda\wr K$.
\end{theorem}

Hence, the knowledge of the isoperimetric profiles of lamplighters (Theorem~\ref{thm:INTROformulaforprofilewreathproducts}) allows us to deduce the next characterisation:

\begin{corollary}\label{cor:INTROquantitativecomparisonbetweenlamplighters}
Let $k,\ell\ge 1$ be integers. Let $p\ge 0$. Let $F$ be a finite group. Then there exists an $(\ld^p, \ld^0)-$orbit equivalence from $F\wr\Z^{k+\ell}$ to $F\wr\Z^{k}$ if and only if $p<\frac{k}{k+\ell}$. 
\end{corollary}

Indeed, if $p<\frac{k}{k+\ell}$, it is proved in~\cite{DKLMT22} that there exists an $(\ld^p,\ld^0)-$orbit equivalence from $\Z^{k+\ell}$ to $\Z^k$, whence the same claim from $F\wr\Z^{k+\ell}$ to $F\wr\Z^{k}$ by Theorem~\ref{thm:INTROstabilityforOEbetweenlamplighters}. Conversely, if such a coupling exists, then Theorem~\ref{thm:INTROObstructionDKLMT}\textit{(i)} implies 
\begin{equation*}
    \ln(n)^{\frac{p}{k}}\simeq \prof{F\wr\Z^{k}}(n)^p \preccurlyeq \prof{F\wr\Z^{k+\ell}}(n) \simeq \ln(n)^{\frac{1}{k+\ell}}
\end{equation*}
whence $p\le \frac{k}{k+\ell}$, and the equality case is excluded by~\cite[Theorem~B]{Cor25}. 

Corollary~\ref{cor:INTROquantitativecomparisonbetweenlamplighters} gives an explicit upper bound on the degree of $\ld^p-$integrability of a cocycle from the bigger group $\Z^{k+\ell}$ to the smaller $\Z^{k}$. Since cocycles measure the required distortion to pass from a Cayley graph of one group to the other, the last statement can be seen as a precise comparison of the geometries of the two groups. 

\bigskip
\noindent \textbf{New results presented in the thesis} 
\bigskip

After this brief overview of the state of the art, let us present the contributions of this thesis. 

\paragraph{Quasi-isometric rigidity for permutational lamplighters.} Chapter~\ref{chap:chapter3} aims at proving a quasi-isometric classification of some permutational wreath products whose lamp groups are finite. The precise statement is the following. 

\begin{theorem}[{see Theorem~\ref{thm:classificationPWPuptoQI}}]\label{thm:INTROclassificationPWPuptoQI}
Let $E$, $F$ be two non-trivial finite groups. Let $G$, $H$ be finitely presented groups, with finitely generated normal infinite subgroups $M\lhd G$, $N\lhd H$. Suppose that $M$ has infinite index in $G$, and that $G$ (resp. $H$) is not coarsely separable by any collection of subspaces that uniformly quasi-isometrically embed into $M$ (resp. $N$). The following claims hold.
\begin{itemize}
    \item If $M$ is not co-amenable in $G$, then $E\wr_{G/M}G$ and $F\wr_{H/N}H$ are quasi-isometric if and only if \;$|E|$ and $|F|$ have the same prime divisors and there exists a quasi-isometry of pairs $(G,M)\longrightarrow (H,N)$. 
    \item If $M$ is co-amenable in $G$, then $E\wr_{G/M}G$ and $F\wr_{H/N}H$ are quasi-isometric if and only if \;$|E|=n^{r}$, $|F|=n^{s}$ for some $n,r,s\ge 1$ and there exists a quasi-isometry of pairs $(G,M)\longrightarrow (H,N)$ inducing a quasi-$\frac{s}{r}$-to-one quasi-isometry $G/M\rightarrow H/N$.
    %\vspace{0.05cm}
\end{itemize}
\end{theorem}

As for Theorem~\ref{thm:INTROclassificationGT21}, the non-amenable situation must be distinguished from the amenable one, the latter bringing more rigidity. Moreover, the undertaken strategy for proving Theorem~\ref{thm:INTROclassificationPWPuptoQI} is similar to that of~\cite{GT24b}, explained in the preliminary Chapter~\ref{chap:chapter2}. 

The main difference with the standard case is that, additionally to preserving cosets of the base groups, a quasi-isometry between two such permutational lamplighters must also preserve the cosets of the involved normal subgroups. We refer to Section~\ref{sec:QIofpairs} for the precise terminology. At this point, let us just emphasize the fact that the presence of a quasi-isometry $E\wr_{G/M} G\longrightarrow F\wr_{H/N} H$ is sufficient to detect the large-scale geometries of the quotients $G/M$ and $H/N$, since a quasi-isometry of pairs $(G,M)\longrightarrow (H,N)$ induces a quasi-isometry $G/M\rightarrow H/N$. 

As an application, we obtain the next corollary, which was in fact the initial motivation for proving the more general Theorem~\ref{thm:INTROclassificationPWPuptoQI}:

\begin{corollary}[{see Corollary~\ref{cor:classificationofPWPoverZ^duptoQI}}]
Let $E$, $F$ be two non-trivial finite groups. Let $m,m',n,n'$ be four integers such that $m\ge n\ge 2$, $m'\ge n'\ge 2$. Then $E\wr_{\Z^n}\Z^m$ and $F\wr_{\Z^{n'}}\Z^{m'}$ are quasi-isometric if and only if $m=m'$, $n=n'$ and $|E|$, $|F|$ are powers of a common number. 
\end{corollary}

In specific cases, where explicit constructions are accessible, it is possible to have more information on the induced quasi-isometry between the quotients. As an application, we show that, for permutational lamplighters over free abelian groups, being quasi-isometric is the same as being biLipschitz equivalent. 

\begin{corollary}[{see Corollary~\ref{cor:classificationPWPoverZ^duptoBILIP}}]
Let $E$ and $F$ be two non-trivial finite groups, and let $m,n$ be two integers such that $m>n\ge 2$. Then $E\wr_{\Z^{n}}\Z^{m}$ and $F\wr_{\Z^{n}}\Z^{m}$ are biLipschitz equivalent if and only if they are quasi-isometric, which happens if and only if\;$|E|$ and $|F|$ are powers of a common number. 
\end{corollary}

In contrast, two standard wreath products can be quasi-isometric without being biLipschitz equivalent (\textit{e.g.} $\Z_{2}\wr\Z^2$ and $\Z_{4}\wr\Z^2$; see~\cite[Corollary~1.15]{GT24b}).

\paragraph{More scaling groups.} A nice consequence of these classification results is that they allow us to identify the scaling group of our spaces, information available only for very specific families of amenable groups. Theorem~\ref{thm:INTROclassificationGT21} extends this list, as it implies that the scaling group of $F\wr K$ is $\lbrace 1\rbrace$ when $F$ is finite and $K$ is amenable, finitely presented and one-ended. As explained in Section~\ref{sec:Scalinggroupsofhaloproducts}, this computation is a consequence of the general form of quasi-isometries between lamplighters, that must be~\textit{aptolic} (see Definition~\ref{def:aptolicity} for the precise definition).  

As proved in~\cite{BGT24}, the aptolicity phenomenon also holds for quasi-isometries between some wreath products having an infinite lamp group. We exploit this in Chapter~\ref{chap:chapter4} to prove the following rigidity property. 

\begin{theorem}[{see Theorem~\ref{thm:maintheoremfortheclassMexp}}]\label{thm:INTROmaintheoremfortheclassMexp}
Let $N$ and $M$ be finitely generated groups of polynomial growth, with growth degrees $n$ and $m$ respectively. Let $G$ and $H$ be finitely presented amenable groups from $\mathcal{M}_{\text{exp}}$. Then any quasi-isometry $N\wr G \longrightarrow M\wr H$ is quasi-$\frac{m}{n}$-to-one. 
\end{theorem}

We refer to Section~\ref{sec:introchapter4} for explanations on the assumptions and the class $\mathcal{M}_{\text{exp}}$. Notice that, in the case of polynomial growth lamp groups, the scaling factor of the quasi-isometry is parametrized by the  involved growth degrees. In particular, it follows that:

\begin{corollary}[{see Corollary~\ref{cor:scalinggroupsdesired}}]\label{cor:INTROscalinggroupsofwreathproducts}
Let $N$ be a finitely generated group with polynomial growth. Let $G$ be a finitely presented amenable group from $\mathcal{M}_{\text{exp}}$. Then any quasi-isometry $N\wr G\longrightarrow N\wr G$ lies at bounded distance from a bijection. In particular, $\text{Sc}(N\wr G)=\lbrace 1\rbrace$.
\end{corollary}

These results have two interesting consequences. The first one is that wreath products of the form $N\wr G$ are~\textit{lamplighter-rigid}, in the sense that when forming lamplighters over such groups, no flexibility with respect to quasi-isometries is possible. 

\begin{proposition}[{see Proposition~\ref{prop:lamplighterrigidity}}]\label{prop:INTROlamplighterrigidity}
Let $n,m\ge 2$ be two integers. Let $N$ be a finitely generated group with polynomial growth, and let $G$ be a finitely presented amenable group from $\mathcal{M}_{\text{exp}}$. Then the lamplighters $\Z_{n}\wr(N\wr G)$ and $\Z_{m}\wr(N\wr G)$ are quasi-isometric if and only if $n=m$. 
\end{proposition}

Pushing the idea further, we derive an arithmetic condition which is able to distinguish some iterated wreath products up to quasi-isometry: 

\begin{proposition}[{see Proposition~\ref{prop:mixingofscalingconditions}}]\label{prop:INTROiteratedlamplighters}
Let $n,m\ge 2$ be two integers. Let $N_{1}$ and $N_{2}$ be finitely generated groups of polynomial growth, with growth degrees $n_{1}$ and $n_{2}$ respectively. Let $G$ and $H$ be finitely presented amenable groups from $\mathcal{M}_{\text{exp}}$. If \;$\Z_{n}\wr(N_{1}\wr G)$ and $\Z_{m}\wr(N_{2}\wr H)$ are quasi-isometric, then there exist $a,r,s\ge 1$ such that $n=a^{r}$, $m=a^{s}$, and $\frac{s}{r}=\frac{n_{2}}{n_{1}}$.
\end{proposition}

For instance, there is no quasi-isometry 
\begin{equation*}
\Z_{2}\wr(\Z^{2}\wr \text{SOL}(\Z)) \longrightarrow \Z_{4}\wr(\Z^{3}\wr \text{SOL}(\Z)).
\end{equation*}
The key observation for getting the last equality in Proposition~\ref{prop:INTROiteratedlamplighters} is that the aptolicity phenomenon implies a mixing of two scaling conditions, so the proposition follows from Theorems~\ref{thm:INTROmaintheoremfortheclassMexp} and~\ref{thm:INTROclassificationGT21}. In fact, we can get a complete classification when $N_{1}$ and $N_{2}$ are virtually abelian; see Corollary~\ref{cor:classificationofiteratedwreathproducts}.

In the remaining part of this introduction, we present the two articles~\cite{cordum25, cordum26}, written jointly with Corentin Correia, in which we show that the halo structure is also particularly well-suited for tracking isoperimetric profiles and orbit equivalence couplings for these spaces. 

\paragraph{Isoperimetric profiles of lampshufflers.}
Chapter~\ref{chap:chapter5} contains results for many classes of halo products, but for conciseness let us focus in this introduction on the case of lampshufflers.

For such groups, a general lower bound on their F\o lner function has been proved by Erschler and Zheng. 

\begin{theorem}[{\cite[Corollary~1.4]{EZ21}}]
Let $H$ be a finitely generated amenable group. Then we have
\begin{equation*}
    \text{F\o l}_{\shuf{H}}(x)\succcurlyeq \beta_{H}(x)^{\beta_{H}(x)}
\end{equation*}
where $\beta_{H}$ is the growth function of $H$. Moreover, this lower bound is optimal if $H$ has polynomial growth. 
\end{theorem}

In particular, when $H$ has exponential growth, inverting this lower bound shows that 
\begin{equation*}
    \prof{\shuf{H}}(x) \preccurlyeq \ln(\ln(x)). 
\end{equation*}
However, this bound cannot be optimal in many cases, since for instance two iterated lamplighters over $\Z$ with different numbers of iterations (e.g. $F\wr\Z$ and $F\wr(F\wr\Z)$) both have exponential growth, whereas we expect the isoperimetric profile of $\shuf{F\wr(F\wr\Z)}$ to be slower than the one of $\shuf{F\wr\Z}$. 

One of the main results of~\cite{cordum25} provides finer bounds, connecting the isoperimetric profile of $\shuf{H}$ directly to the profile of $H$:

\begin{theorem}[{see Theorem~\ref{thm:boundsForProfile intro}}]\label{thm:INTROboundsForProfile intro}
Let $p\ge 1$. Let $H$ be a finitely generated amenable group whose $\ell^p-$isoperimetric profile $\profp{H}$ satisfies Assumption~$(\star)$. Then the $\ell^p-$isoperimetric profile $\profp{\shuf{H}}$ of $\shuf{H}$ satisfies
\begin{equation*}
    \profp{H}\left(\frac{\ln(x)}{\ln(\ln(x))}\right) \preccurlyeq \profp{\shuf{H}}(x) \preccurlyeq \prof{H}(\ln(x)). 
\end{equation*}
\end{theorem}

In this statement, Assumption~$(\star)$ is the one appearing in Theorem~\ref{thm:INTROformulaforprofilewreathproducts}, and indeed the technique for achieving the upper bound is to find, inside $\shuf{H}$, copies of lamplighter graphs and to use the monotonicity of the $\ell^p-$isoperimetric profile combined with (a slight extension of) Erschler's formula (Theorem~\ref{thm:INTROformulaforprofilewreathproducts}). 

From Theorem~\ref{thm:INTROboundsForProfile intro}, it is also possible to derive exact estimates for the $\ell^p-$profile of iterated lampshufflers, defined inductively by
\begin{equation*}
    \shufn{n}{H}\defeq \shuf{\shufn{n-1}{H}}
\end{equation*}
if $n\ge 1$ and $\shufn{0}{H}\defeq H$. For instance, in the polynomial growth range, this gives:

\begin{proposition}[{see Proposition~\ref{prop:profileofshufnofpolynomialgrowthgroupsINTRO}}]\label{prop:INTROprofileofshufnofpolynomialgrowthgroupsINTRO}
Let $H$ be a finitely generated group of polynomial growth of degree $d\ge 1$. Then one has 
\begin{equation*}
    \profp{\shufn{n}{H}}(x) \simeq \left(\frac{\ln^{\circ n}(x)}{\ln^{\circ (n+1)}(x)}\right)^{\frac{1}{d}}
\end{equation*}
for any integer $n\ge 1$ and any real number $p\ge 1$.
\end{proposition}

In turn, these computations also contribute to the quasi-isometric classification program initiated in~\cite{GT24a}, and we are also able to prove that, in some cases, being quasi-isometric is the same as being biLipschitz equivalent. As an illustration:

\begin{corollary}[{see Theorem~\ref{thm:IteratedShufflersPolynomialQIBiLip intro}}]\label{cor:INTROiteratedlampshufflersoverZ^dQI}
Let $d,k\ge 1$ and $n,m\ge 0$ be integers. Then the iterated lampshufflers $\shufn{n}{\Z^d}$ and $\shufn{m}{\Z^k}$ are quasi-isometric if and only if they are biLipschitz equivalent, which happens if and only if $n=m$ and $d=k$. 
\end{corollary}

\paragraph{Quantitative orbit equivalence for lampshufflers.} As we explained above, when two groups are not quasi-isometric, quantitative orbit equivalence is suitable for comparing a bit more precisely their geometries. The existence of optimal orbit equivalence couplings for lamplighters has been established in~\cite{DKLMT22} and~\cite{Cor25}. In~\cite{cordum26}, we undertake a similar study for other halo products. 

For instance, for lampshufflers, we establish the following stability result: 
\begin{theorem}[{see Theorem~\ref{thm:stabilityofcouplings+quantificationLampjugglersINTRO}}]
If two groups $H$ and $K$ are orbit equivalent, then $\shuf{H}$ and $\shuf{K}$ are orbit equivalent. Moreover, if $H$ and $K$ are finitely generated, if $\varphi,\psi\colon\R_{+}\rightarrow\R_{+}$ are non-decreasing maps, and if there exists a $(\varphi,\psi)-$integrable orbit equivalence coupling from $H$ to $K$, then the same holds from $\shuf{H}$ to $\shuf{K}$.
\end{theorem}

This statement should be compared to Theorem~\ref{thm:INTROstabilityforOEbetweenlamplighters}, and in fact we extract from the proof of the latter a general method that can, in practice, be applied for any halo product satisfying some mild assumptions. See Section~\ref{sec:generalmethodstability} for more details. 

Combined with our computations of isoperimetric profiles of iterated lampshufflers (cf. Proposition~\ref{prop:INTROprofileofshufnofpolynomialgrowthgroupsINTRO}), we can deduce optimal orbit equivalence couplings between iterated lampshufflers and lampjugglers over free abelian groups using to Theorem~\ref{thm:INTROObstructionDKLMT}.

\begin{corollary}[{see Theorem~\ref{thm:optimalitypolynomialgrowthJugglersINTRO}}]
Let $k,d\ge 1$ be positive integers such that $k> d$. Let $n\ge 0$ be an integer.
Then $\shufn{n}{\Z^k}$ and $\shufn{n}{\Z^d}$ are $\ld^{p}$ orbit equivalent if and only if $p<\frac{d}{k}$.
\end{corollary}

\bigskip

\noindent \textbf{Structure of the thesis} 

\bigskip

\thispagestyle{empty}
Chapter~\ref{chap:chapter1} contains all the required material on geometric group theory to dive into the subsequent chapters. Chapter~\ref{chap:chapter2} aims at presenting, in a concise manner, the content of the articles~\cite{GT24a, GT24b}, since the first part of the thesis relies heavily on results from these papers. Chapter~\ref{chap:chapter3} presents the content of~\cite{Dum24} and establishes in particular Theorem~\ref{thm:INTROclassificationPWPuptoQI}. Chapter~\ref{chap:chapter4} presents the contributions of~\cite{Dum26} about wreath products with infinite lamp groups (Theorem~\ref{thm:INTROmaintheoremfortheclassMexp}, Corollary~\ref{cor:INTROscalinggroupsofwreathproducts}, Propositions~\ref{prop:INTROlamplighterrigidity} and~\ref{prop:INTROiteratedlamplighters}). We study isoperimetric profiles of halo products in Chapter~\ref{chap:chapter5}, which also contains our applications to questions concerning the existence of regular maps and quasi-isometries between such halo products (among which Corollary~\ref{cor:INTROiteratedlampshufflersoverZ^dQI}). Our results on quantitative orbit equivalence are contained in Chapter~\ref{chap:chapter6} and rely on the paper~\cite{cordum26}. 

\endgroup

\afterpage{\blankpage} %partie introduction ou avant-propos

\fancyhead{} % clear all header fields

\fancyhead[RO]{\leftmark}
\fancyhead[LE]{\textsc{\chaptername~\thechapter}}

\chapter{Generalities on geometric group theory}\label{chap:chapter1}
The goal of this chapter is to develop the necessary tools to study finitely generated groups as metric spaces.

\vspace{0.3cm}

\minitoc

\section{Finitely generated groups}\label{subsection1.1}

The goal of this subsection is to introduce several classes of finitely generated groups of interest. They will constitute our running examples that will follow us for the rest of this text. 

\begin{example}\label{ex:finitelygeneratedgroups}
\textit{(i)} Finite groups are finitely generated.

\noindent \textit{(ii)} The group $(\Z,+)$ is finitely generated, and $S=\lbrace -1,1\rbrace$ is a finite generating set. More generally, if $p,q\in\Z$ are coprime, then $\lbrace \pm p, \pm q\rbrace$ is a symmetric generating set for $\Z$.

\noindent \textit{(iii)} In fact, for any $d\ge 1$, the group $\Z^d$ is finitely generated, and a symmetric generating set is given by the canonical basis 
\begin{equation*}
    \lbrace \pm(1,0,\dots,0), \pm (0,1,0,\dots,0), \dots, \pm (0,0,\dots,0,1)\rbrace.
\end{equation*}

\noindent \textit{(iv)} If $d\ge 1$, the non-abelian free group $F_{d}$ of rank $d$ is finitely generated, a generating set being given by the equivalence classes of words of length one over a set $S$ of cardinality $d$. 

\noindent \textit{(v)} It is not hard to check that the \textit{Heisenberg group} 
\begin{equation*}
    H(\Z) \defeq \Bigg\lbrace \begin{pmatrix}      1 & a & c \\
   0 & 1 & b \\
   0 & 0 & 1
\end{pmatrix} : a,b,c\in\Z\Bigg\rbrace 
\end{equation*}
is generated by the three matrices 
\begin{equation*}
    x=\begin{pmatrix}       
   1 & 1 & 0 \\
   0 & 1 & 0 \\
   0 & 0 & 1
   \end{pmatrix}, \; y=\begin{pmatrix}       
   1 & 0 & 0 \\
   0 & 1 & 1 \\
   0 & 0 & 1
   \end{pmatrix}, \; z=\begin{pmatrix}       
   1 & 0 & 1 \\
   0 & 1 & 0 \\
   0 & 0 & 1
   \end{pmatrix}.
\end{equation*} 

\noindent \textit{(vi)} The group $(\Q,+)$ is not finitely generated. Indeed, suppose for a contradiction that $\Q$ is generated by finitely many rationals $\frac{p_{1}}{q_{1}},\dots,\frac{p_{n}}{q_{n}}$. Any finite sum of these fractions or their inverses is a rational number with denominator at most $q_{1}\dots q_{n}$. Letting $N\defeq q_{1}\dots q_{n}+1$, it follows that $\frac{1}{N}$ cannot be written using $\frac{p_{1}}{q_{1}},\dots,\frac{p_{n}}{q_{n}}$ or their inverses, a contradiction. 

\noindent \textit{(vii)} The group $D_{\infty}\defeq \langle a,t : a^2=1, ata^{-1}=t^{-1}\rangle$ is called the \textit{infinite dihedral group}, and generalises finite dihedral groups allowing a rotation of infinite order. In fact, this group is isomorphic to 
\begin{equation*}
    \langle a,b : a^2=b^2=1\rangle=\Z_{2} * \Z_{2}
\end{equation*}
the free product of two cyclic groups of order $2$.  
\end{example}

Note that a finitely generated group is always countable, but the converse does not hold, as shown by Example~\ref{ex:finitelygeneratedgroups}\textit{(vi)}.

The following proposition shows that being finitely generated is stable under group extensions. 

\begin{proposition}\label{prop:finitegenerationextension}
Let $G$ be a group and $N\lhd G$. If $G$ is finitely generated, then $G/N$ is finitely generated. Conversely, if $N, G/N$ are finitely generated, then $G$ is finitely generated.
\end{proposition}

\begin{proof}
If $G$ is finitely generated and $\pi\colon G\rightarrow G/N$ is the natural surjection, then the image under $\pi$ of a generating set for $G$ is a generating set for $G/N$. 

\noindent Conversely, let $\lbrace g_{1},\dots,g_{n}\rbrace$ be a generating set for $N$ and $\lbrace h_{1}N,\dots,h_{m}N\rbrace$ a generating set for $G/N$. Fix $g\in G$. Then there exists $\varepsilon_{1},\dots,\varepsilon_{m}\in\lbrace -1,1\rbrace$ such that
\begin{equation*}
    gN=(h_{1}N)^{\varepsilon_{1}}\dots(h_{m}N)^{\varepsilon_{m}}=(h_{1}^{\varepsilon_{1}}\dots h_{m}^{\varepsilon_{m}})N
\end{equation*}
and it follows that $g(h_{1}^{\varepsilon_{1}}\dots h_{m}^{\varepsilon_{m}})^{-1}\in N$. We can then write 
\begin{equation*}
    g(h_{1}^{\varepsilon_{1}}\dots h_{m}^{\varepsilon_{m}})^{-1}=g_{1}^{\delta_{1}}\dots g_{n}^{\delta_{n}}
\end{equation*}
for some $\delta_{1},\dots,\delta_{n}\in\lbrace -1,1\rbrace$, and thus $g=g_{1}^{\delta_{1}}\dots g_{n}^{\delta_{n}}h_{1}^{\varepsilon_{1}}\dots h_{m}^{\varepsilon_{m}}$. This proves that 
\begin{equation*}
\lbrace g_{1},\dots,g_{n},h_{1},\dots,h_{m}\rbrace
\end{equation*}
is a finite generating set for $G$, and the proof is complete.
\end{proof}

On the other hand, it is in general not true that subgroups of finitely generated groups are themselves finitely generated. To produce such an example, we introduce an additional group construction, which is the central definition for this thesis. 

\begin{definition}\label{def:wreathproduct}
Let $A$ and $B$ be two groups. Their wreath product $A\wr B$ is the group defined by 
\begin{equation*}
    \bigg(\bigoplus_{B}A\bigg)\rtimes B
\end{equation*}
where $B$ acts on the direct sum via
\begin{equation*}
    (b\cdot f)(b') \defeq f(b^{-1}b')
\end{equation*}
for any $b,b'\in B$ and any $f\in \bigoplus_{B}A$.
\end{definition}

Hence, elements of $A\wr B$ are pairs $(f,b)$ where $f$ is a~\textit{finitely supported} function on $B$ (i.e. $f(b)=e_{A}$ for all but finitely many $b\in B$) and $b\in B$. The multiplication law is given by 
\begin{equation*}
    (f,b)(f',b') = (f+b\cdot f', bb')
\end{equation*}
for all $f,f'\in \bigoplus_{B}A$, $b,b'\in B$, where “$+$” stands for the composition law in the direct sum. 

We then prove that this construction preserves finite generation.

\begin{proposition}\label{prop:wreathproductspreservefinitegeneration}
If $A$ and $B$ are finitely generated, then $A\wr B$ is finitely generated.
\end{proposition}

\begin{proof}
Let $S=\lbrace a_{1},\dots,a_{n}\rbrace$ be a generating set for $A$, and let $T=\lbrace b_{1},\dots,b_{m}\rbrace$ be a generating set for $B$. For any $a\in A$, let $\delta_{a} \in \bigoplus_{B}A$ be defined by $\delta_{a}(e_{B})=a$ and $\delta_{a}(b)=e_{A}$ for any $b\neq e_{B}$. Let also $\textbf{1}$ denote the neutral element of the direct sum, defined as $\textbf{1}(b)=e_{A}$ for any $b\in B$. We claim that the finite set 
\begin{equation*}
    U\defeq \lbrace (\delta_{a_{i}},e_{B}), (\textbf{1}, b_{j}) : 1\le i\le n, 1\le j\le m\rbrace
\end{equation*}
is a generating set for $A\wr B$. 

\noindent First, as the multiplication in $A\wr B$ is the multiplication of $B$ in the second component, and as $T$ generates $B$, it is enough to prove that any pair of the form $(f,e_{B})$ is a product of
\begin{equation*}
    (\delta_{a_{1}},e_{B}),\dots, (\delta_{a_{n}},e_{B}).
\end{equation*}
Since $f$ is finitely supported, it is enough to prove that any pair of the form $(\delta_{a},e_{B})$, $a\in A$, is a product of $(\delta_{a_{1}},e_{B}),\dots, (\delta_{a_{n}},e_{B})$. For $a\in A$, write 
\begin{equation*}
    a=a_{i_{1}}\dots a_{i_{k}}
\end{equation*}
for some $i_{1},\dots,i_{k}\in \lbrace 1,\dots,n\rbrace$, and then it follows that 
\begin{equation*}
    (\delta_{a},e_{B})=(\delta_{a_{i_{1}}},e_{B})\dots (\delta_{a_{i_{k}}},e_{B}).
\end{equation*}
Thus $A\wr B$ is finitely generated. 
\end{proof}

From Proposition~\ref{prop:wreathproductspreservefinitegeneration}, it follows that $\Z/2\Z \wr \Z$ is finitely generated, but it contains $\bigoplus_{\Z}\Z/2\Z$ as a subgroup, and the latter is not finitely generated. 

In Sections~\ref{sec:QIrigidityforhaloproducts} and~\ref{sec:halo}, we introduce similar group constructions, called~\textit{halo products}, that also provide examples of finitely generated groups with non-finitely generated subgroups. 

\section{The metric coarse category}

Now the objects are defined, let us introduce morphisms that relate them. In this section, we give definitions in a slightly more general framework, namely the one of pseudo-metric spaces, and we explain in Section~\ref{sec:groupsaspseudometricspaces} how finitely generated groups fit in this setup. We follow mainly the exposition of~\cite{CH16}. 

\subsection{Coarsely Lipschitz maps and large-scale Lipschitz maps}

An~\textit{upper control} is a non-decreasing function $\Phi_{+}\colon \R_{+}\rightarrow \R_{+}$, and a~\textit{lower control} is a non-decreasing function $\Phi_{-}\colon \R_{+}\rightarrow \R_{+}\cup\lbrace\infty\rbrace$ such that $\lim\limits_{t\rightarrow \infty}\Phi_{-}(t)=\infty$.

If $X$ and $Y$ are pseudo-metric spaces and $f\colon X\rightarrow Y$ is a map, an~\textit{upper control for $f$} is an upper control $\Phi_{+}$ such that 
\begin{equation*}
        d_{Y}(f(x),f(x')) \le \Phi_{+}(d_{X}(x,x'))
\end{equation*}
for any $x,x'\in X$, and dually a~\textit{lower control for $f$} is a lower control $\Phi_{-}$ such that 
\begin{equation*}
        \Phi_{-}(d_{X}(x,x')) \le d_{Y}(f(x),f(x')) 
\end{equation*}
for any $x,x'\in X$.

\begin{definition}\label{def:familiesofmaps}
Let $X,Y$ be pseudo-metric spaces and let $f\colon X\rightarrow Y$ be a map. We say that $f$ is 
\begin{enumerate}[label=(\roman*)]
    \item \textit{coarsely Lipschitz} if there exists an upper control for $f$.
    \item \textit{coarsely expansive} if there exists a lower control for $f$.
    \item a \textit{coarse embedding} if it is coarsely Lipschitz and coarsely expansive.
    \item \textit{essentially surjective} if $f(X)$ is co-bounded in $Y$.
    \item a \textit{metric coarse equivalence} if it is an essentially surjective coarse embedding. 
\end{enumerate}
\end{definition}

If there exists a metric coarse equivalence $f\colon X\rightarrow Y$, we say that $X$ and $Y$ are~\textit{coarsely equivalent}. 

Let us start with equivalent reformulations of the above conditions.

\begin{proposition}\label{prop:formulationofcoarselyLipschitz}
Let $X,Y$ be pseudo-metric spaces and $f\colon X\rightarrow Y$ a map. The following are equivalent.
\begin{enumerate}[label=(\roman*)]
    \item The map $f$ is coarsely Lipschitz.
    \item For all $R\ge 0$, there exists $S\ge 0$ such that if $x,x'\in X$ satisfy $d_{X}(x,x')\le R$, then $d_{Y}(f(x),f(x'))\le S$.
    \item For any sequence of points $(x_{n})_{n\in\N}, (x_{n}')_{n\in\N}$ in $X$ with $\dis\sup_{n\in\N}d_{X}(x_{n},x_{n}')<\infty$, we have
    \begin{equation*}
        \sup_{n\in\N}d_{Y}(f(x_{n}), f(x_{n}'))<\infty.
    \end{equation*}
\end{enumerate}
\end{proposition}

\begin{proof}
\textit{(i)} $\Longrightarrow$ \textit{(ii)}: Suppose that $f$ is coarsely Lipschitz, and denote $\Phi_{+}$ an upper control for $f$. Let $R\ge 0$, and set $S\defeq \Phi_{+}(R) \ge 0$. Then, if $x,x'\in X$ are such that $d_{X}(x,x') \le R$, it follows that 
\begin{equation*}
    d_{Y}(f(x),f(x')) \le \Phi_{+}(d_{X}(x,x')) \le \Phi_{+}(R) = S
\end{equation*}
since $f$ is coarsely Lipschitz and $\Phi_{+}$ is non-decreasing. Thus~\textit{(ii)} holds. 

\noindent \textit{(ii)} $\Longrightarrow$ \textit{(iii)}: Let $(x_{n})_{n\in\N}, (x_{n}')_{n\in\N}$ be two sequences of points in $X$ with 
\begin{equation*}
    C\defeq \dis\sup_{n\in\N}d_{X}(x_{n},x_{n}') < \infty.
\end{equation*}
Using~\textit{(ii)}, there is $S\ge 0$ such that 
\begin{equation*}
    d_{Y}(f(x),f(x')) \le S
\end{equation*}
if $d_{X}(x,x')\le C$. As $d_{X}(x_{n},x_{n}')\le C$ for any $n\in\N$, we thus have
\begin{equation*}
    d_{Y}(f(x_{n}), f(x_{n}')) \le S
\end{equation*}
for any $n\in\N$, hence $\sup_{n\in\N}d_{Y}(f(x_{n}), f(x_{n}')) \le S <\infty$, which shows~\textit{(iii)}. 

\noindent \textit{(iii)} $\Longrightarrow$ \textit{(i)}: For $c\in \R_{+}$, define 
\begin{equation*}
    \Phi_{+}(c) \defeq \sup\lbrace d_{Y}(f(x),f(x')) : x,x'\in X,\; d_{X}(x,x')\le c\rbrace.
\end{equation*}
Then $\Phi_{+}$ is positive and non-decreasing. Towards a contradiction, suppose that $\Phi_{+}(c)=\infty$ for some $c\in\R_{+}$. This implies that there exist two sequences $(x_{n})_{n\in\N}, (x_{n}')_{n\in\N}\subset X$ with $d_{X}(x_{n},x_{n}')\le c$ for any $n\in \N$ and 
\begin{equation*}
    \lim\limits_{n\rightarrow\infty} d_{Y}(f(x_{n}), f(x_{n}')) = \Phi_{+}(c)=\infty
\end{equation*}
which is excluded by \textit{(iii)}. Hence $\Phi_{+}$ takes only finite values, and thus is indeed an upper control for $f$. 
\end{proof}

Dualising the above proof, one gets the same statement for coarsely expansive maps. 

\begin{proposition}\label{prop:formulationofcoarselyexpensive}
Let $X,Y$ be pseudo-metric spaces and $f\colon X\rightarrow Y$ a map. The following are equivalent.
\begin{enumerate}[label=(\roman*)]
    \item The map $f$ is coarsely expansive.
    \item For all $r\ge 0$, there exists $s\ge 0$ such that if $x,x'\in X$ have $d_{X}(x,x')\ge r$, then one has $d_{Y}(f(x),f(x'))\ge s$. 
    \item For any sequence of points $(x_{n})_{n\in\N}, (x_{n}')_{n\in\N}$ in $X$ with $\lim\limits_{n\rightarrow\infty}d_{X}(x_{n},x_{n}')=\infty$, we have
    \begin{equation*}
        \lim\limits_{n\rightarrow\infty}d_{Y}(f(x_{n}), f(x_{n}'))=\infty.
    \end{equation*}
\end{enumerate}
\end{proposition}

Given two maps $f,f'\colon X\rightarrow Y$ between pseudo-metric spaces, we say that $f'$ is~\textit{close} to $f$, or also that $f'$ is~\textit{at bounded distance} from $f$, if there exists $C>0$ such that 
\begin{equation*}
    d_{Y}(f(x),f'(x)) \le C
\end{equation*}
for any $x\in X$. This is the same as requiring that 
\begin{equation*}
    \sup_{x\in X}d_{Y}(f(x),f'(x)) < \infty
\end{equation*}
and, in this case, we write $f\sim_{c} f'$.

\begin{lemma}\label{lem:closenessER}
Closeness is an equivalence relation.
\end{lemma}

\begin{proof}
Clearly $f\sim_{c} f$ as $d_{Y}(f(x),f(x))=0$ for any $x\in X$. Symmetry of $\sim_{c}$ follows from symmetry of $d_{Y}$, and transitivity follows from the triangle inequality for $d_{Y}$. 
\end{proof}

The next result shows that properties from Definition~\ref{def:familiesofmaps} are invariant by finite distance changes. 

\begin{proposition}\label{prop:propertiesstableundercomposition/closeness}
Let $X,Y,Z$ be pseudo-metric spaces, $f,f'\colon X\rightarrow Y$ two close maps, and $g,g'\colon Y\rightarrow Z$ two close maps. 

\noindent \textit{(i)} The map $f$ is coarsely Lipschitz (resp. coarsely expansive, a coarse embedding, essentially surjective, a metric coarse equivalence) if and only if $f'$ is coarsely Lipschitz (resp. coarsely expansive, a coarse embedding, essentially surjective, a metric coarse equivalence).

\noindent \textit{(ii)} If $f,g$ are both coarsely Lipschitz (resp. coarsely expansive, coarse embeddings, essentially surjective, metric coarse equivalences), then $g\circ f$ is coarsely Lipschitz (resp. coarsely expansive, a coarse embedding, essentially surjective, a metric coarse equivalence).

\noindent \textit{(iii)} If $g$ is coarsely Lipschitz, then $g\circ f$ and $g'\circ f'$ are close.
\end{proposition}

\begin{proof}
\textit{(i)} Suppose $f$ is coarsely Lipschitz, and let $C>0$ be such that $d_{Y}(f(x), f'(x))\le C$ for any $x\in X$. Let $R\ge 0$. As $f$ is coarsely Lipschitz, we find $K\ge 0$ such that
\begin{equation*}
    d_{X}(x,x') \le R \Longrightarrow d_{Y}(f(x),f(x')) \le K.
\end{equation*}
Set $S\defeq K+2C$, and let $x,x'\in X$ with $d_{X}(x,x') \le R$. Then it follows that 
\begin{align*}
    d_{Y}(f'(x), f'(x')) &\le d_{Y}(f'(x), f(x))+d_{Y}(f(x),f(x'))+d_{Y}(f(x'),f'(x')) \\
    & \le K+2C \\
    &=S.
\end{align*}
As $R\ge 0$ was arbitrary, Proposition~\ref{prop:formulationofcoarselyLipschitz} guarantees that $f'$ is coarsely Lipschitz as well. The converse follows swapping the roles of $f$ and $f'$.

\noindent Now, suppose that $f$ is essentially surjective. As above, let $C>0$ be such that 
\begin{equation*}
    d_{Y}(f(x), f'(x))\le C
\end{equation*}
for any $xl\in X$, and let $C'>0$ be such that any point in $Y$ is at distance at most $C'$ from the image of $f$. Let $y\in Y$, and choose $x\in X$ with $d_{Y}(f(x),y) \le C'$. Then we get that 
\begin{align*}
    d_{Y}(f'(x),y) \le d_{Y}(f'(x),f(x))+d_{Y}(f(x),y) \le C+C'.
\end{align*}
Hence any point of $Y$ is at distance at most $C+C'$ from the image of $f'$, i.e. $f'$ is essentially surjective. Once again, the converse follows by symmetry, and the proofs for the other properties are completely similar.

\noindent \textit{(ii)} Here also we only do the proof for one of the properties, and similar arguments apply for the others. Suppose for instance that $f$ and $g$ are both coarsely expansive. Let $r\ge 0$. Applying Proposition~\ref{prop:formulationofcoarselyexpensive}\textit{(ii)} to $f$, we find $s\ge 0$ such that 
\begin{equation*}
    d_{X}(x,x')\ge r \Longrightarrow d_{Y}(f(x),f(x'))\ge s
\end{equation*}
and applying now Proposition~\ref{prop:formulationofcoarselyexpensive}\textit{(ii)} to $g$, there is $t\ge 0$ such that 
\begin{equation*}
    d_{Y}(y,y') \ge s \Longrightarrow d_{Z}(g(y),g(y')) \ge t. 
\end{equation*}
Combining these two implications, we conclude that if $x,x'\in X$ are such that $d_{X}(x,x') \ge r$, then 
\begin{equation*}
    d_{Z}\big((g\circ f)(x),(g\circ f)(x')\big)\ge t
\end{equation*}
proving that $g\circ f$ is coarsely expansive.

\noindent \textit{(iii)} Let $C>0$ be such that $d_{Y}(f(x),f'(x))\le C$ for any $x\in X$, and let $C'>0$ playing the same role for $g$ and $g'$. As $g$ is coarsely Lipschitz, there is $K\ge 0$ such that $d_{Z}(g(y),g(y'))\le K$ if $d_{Y}(y,y')\le C$. Then for any $x\in X$ one has
\begin{equation*}
    d_{Z}\big(g(f(x)),g'(f'(x))\big) \le d_{Z}\big(g(f(x)), g(f'(x))\big)+d_{Z}\big(g(f'(x)), g'(f'(x))\big) \le K+C'
\end{equation*}
whence $g\circ f$ and $g'\circ f'$ are close.
\end{proof}

This proposition motivates then the next definition.

\begin{definition}\label{def:metriccoarsecategory}
Let $X,Y$ be pseudo-metric spaces. A~\textit{coarse morphism} from $X$ to $Y$ is a closeness class of coarsely Lipschitz maps from $X$ to $Y$. The~\textit{metric coarse category} is the category whose objects are pseudo-metric spaces and whose morphisms are coarse morphisms. 
\end{definition}

\begin{definition}
Let $X,Y$ be pseudo-metric spaces and $f\colon X\rightarrow Y$ be a map. We say that $f$ is 
\begin{enumerate}[label=(\roman*)]
    \item~\textit{large-scale Lipschitz} if it has an affine upper control, i.e. there exist $c_{+}>0$, $c_{+}'\ge 0$ such that 
    \begin{equation*}
        d_{Y}(f(x),f(x'))\le c_{+}\cdot d_{X}(x,x')+c_{+}'
    \end{equation*}
    for any $x,x'\in X$.
    \item~\textit{large-scale expansive} if it has an affine lower control, i.e. there exist $c_{-}>0$, $c_{-}'\ge 0$ such that
    \begin{equation*}
        d_{Y}(f(x), f(x'))\ge c_{-}\cdot d_{X}(x,x')-c_{-}'
    \end{equation*}
    for any $x,x'\in X$.
    \item a~\textit{quasi-isometric embedding} if it is large-scale Lipschitz and large-scale expansive.
    \item a~\textit{quasi-isometry} if it is an essentially surjective quasi-isometric embedding.
    \end{enumerate}
\end{definition}

If there is a quasi-isometry $f\colon X\rightarrow Y$, we say that $X$ and $Y$ are~\textit{quasi-isometric}, and we denote $X\sim_{Q.I.}Y$. As we will see below, $\sim_{Q.I.}$ is an equivalence relation among pseudo-metric spaces.

In particular, notice that any large-scale Lipschitz map is coarsely Lipschitz, and any large-scale expansive map is coarsely expansive.

\begin{example}\label{ex:QI}
\textit{(i)} For any $n\ge 1$, the natural inclusion $\Z^n\hookrightarrow \R^n$ is a quasi-isometry, since it is an isometry and since any $n-$tuple of real numbers $(x_{1},\dots,x_{n})$ is at distance at most $\sqrt{n}$ from an $n-$tuple of integers, namely $(\lfloor x_{1} \rfloor,\dots, \lfloor x_{n} \rfloor)$.

\noindent \textit{(ii)} Likewise, the natural inclusion $2\Z \hookrightarrow \Z$ is also a quasi-isometry, since it is an isometric map and since any integer is at distance at most $1$ from an even integer. 
\end{example}

Here goes the natural analog of Proposition~\ref{prop:propertiesstableundercomposition/closeness} for large-scale Lipschitz and expansive maps. 

\begin{proposition}
Let $X,Y,Z$ be pseudo-metric spaces, $f,f'\colon X\rightarrow Y$ two close maps, and $g,g'\colon Y\rightarrow Z$ two close maps.

\noindent \textit{(i)} The map $f$ is large-scale Lipschitz (resp. large-scale expansive, a quasi-isometric embedding, a quasi-isometry) if and only if $f'$ is large-scale Lipschitz (resp. large-scale expansive, a quasi-isometric embedding, a quasi-isometry).

\noindent \textit{(ii)} If $f,g$ are both large-scale Lipschitz (resp. large-scale expansive, quasi-isometric embeddings, quasi-isometries), then $g\circ f$ is large-scale Lipschitz (resp. large-scale expansive, a quasi-isometric embedding, a quasi-isometry).
\end{proposition}

This in turn leads to a natural analog of the metric coarse category for large-scale Lipschitz maps.

\begin{definition}
Let $X$ and $Y$ be pseudo-metric spaces. A~\textit{large-scale morphism} from $X$ to $Y$ is a closeness class of large-scale Lipschitz maps from $X$ to $Y$. The~\textit{large-scale category} is the subcategory of the metric coarse category whose objects are pseudo-metric spaces and whose morphisms are large-scale Lipschitz morphisms. 
\end{definition}

\begin{definition}\label{def:largescalecategory}
Let $X,Y$ be pseudo-metric spaces and $f\colon X\rightarrow Y$ be a map. We say that $f$ is 
\begin{enumerate}[label=(\roman*)]
    \item~\textit{Lipschitz} if there is $c_{+}>0$ such that 
    \begin{equation*}
        d_{Y}(f(x),f(x'))\le c_{+}\cdot d_{X}(x,x')
    \end{equation*}
    for any $x,x'\in X$.
    \item~\textit{biLipschitz} if there exist $c_{+}>0$, $c_{-}> 0$ such that
    \begin{equation*}
        c_{-}\cdot d_{X}(x,x')\le d_{Y}(f(x), f(x'))\le c_{+}\cdot d_{X}(x,x')
    \end{equation*}
    for any $x,x'\in X$.
    \item a~\textit{biLipschitz equivalence} if it is biLipschitz and surjective. 
\end{enumerate}
\end{definition}

Let us now give additional examples of such maps. 

\begin{example}\label{ex:examplesofmaps}
\textit{(i)} Let $f\colon X\rightarrow Y$ be a map between two pseudo-metric spaces. If $X$ has finite diameter, then $f$ is large-scale expansive, since 
\begin{equation*}
    d_{Y}(f(x),f(x')) \ge 0 \ge d_{X}(x,x')-\text{diam}(X)
\end{equation*}
for all $x,x'\in X$. If rather $f(X)$ has finite diameter, then $f$ is large-scale Lipschitz, since 
\begin{equation*}
    d_{Y}(f(x),f(x')) \le \text{diam}(f(X))
\end{equation*}
for all $x,x'\in X$. Lastly, if $Y$ has finite diameter and $X\neq\emptyset$, then $f$ is essentially surjective. Combining these three facts, it follows that any non-empty pseudo-metric space of finite diameter is quasi-isometric to the one point space. 

\noindent \textit{(ii)} For any $p\ge 1$, the map $\text{Id}_{\R^n}\colon (\R^{n},d_{\infty})\rightarrow (\R^n,d_{p})$ is a biLipschitz equivalence, since 
\begin{equation*}
    d_{\infty}(x,y) \le d_{p}(x,y) \le n^{\frac{1}{p}}\cdot d_{\infty}(x,y)
\end{equation*}
for any $x,y\in \R^n$, where the metrics $d_{p}$, $d_{\infty}$ are defined as 
\begin{equation*}
    d_{p}(x,y)\defeq \bigg(\sum_{i=1}^{n}|x_{i}-y_{i}|^{p}\bigg)^{\frac{1}{p}}, \; d_{\infty}(x,y) \defeq \max_{1\le i\le n}|x_{i}-y_{i}|
\end{equation*}
for any $x=(x_{1},\dots,x_{n}), y=(y_{1},\dots, y_{n})\in\R^n$. 

\noindent \textit{(iii)} Let $(X,d)$ be a pseudo-metric space. Define a metric $d_{1}$ by setting 
\begin{equation*}
    d_{1}(x,x')\defeq \max\big(1, d(x,x')\big)
\end{equation*}
for all $x\neq x'\in X$ and $d_{1}(x,x)=0$ for all $x\in X$. Then the map $\text{Id}_{X}\colon (X,d) \rightarrow (X,d_{1})$ is a quasi-isometry. Define now another metric $d_{\ln}$ on $X$ by the formula
\begin{equation*}
    d_{\ln}(x,x')\defeq \ln\big(1+d(x,x')\big), \;x,x'\in X.
\end{equation*}
The map $\text{Id}_{X}\colon (X,d)\rightarrow (X,d_{\ln})$ is a metric coarse equivalence, since it is surjective and the functions $\Phi_{-}(t)=\Phi_{+}(t)=\ln(1+t)$ are lower and upper controls for $\text{Id}_{X}$. We claim that $\text{Id}_{X}\colon (X,d)\rightarrow (X,d_{\ln})$ is large-scale expansive if and only if $(X,d)$ has finite diameter. 

\begin{proof}
If $(X,d)$ has finite diameter, $\text{Id}_{X}$ is large-scale expansive by \textit{(i)} above. Conversely, assume there exist $c_{-}>0$, $c_{-}'\ge 0$ with 
\begin{equation*}
    d_{\ln}(x,x')=\ln\big(1+d(x,x')\big)\ge c_{-}\cdot d(x,x')-c_{-}'
\end{equation*}
for all $x,x'\in X$. Towards a contradiction, assume that $\text{diam}(X,d)=\infty$, and pick two sequences $(x_{n})_{n\in\N}$, $(x_{n}')_{n\in\N}$ in $X$ such that $d(x_{n},x_{n}') \rightarrow \infty$ as $n\rightarrow \infty$. It follows from the assumption that 
\begin{equation*}
   \frac{\ln\big(1+d(x_{n},x_{n}')\big)}{d(x_{n},x_{n}')} \ge c_{-}-\frac{c_{-}'}{d(x_{n},x_{n}')}
\end{equation*}
for all $n\in\N$ large enough. Letting $n\rightarrow\infty$ in this inequality provides $c_{-}\le 0$, a contradiction. Thus $(X,d)$ must have finite diameter.  
\end{proof}

\noindent \textit{(iv)} A pseudo-metric space $X$ is \textit{hyperdiscrete} if the set 
\begin{equation*}
    \big\lbrace (x,x')\in X^2 : d_{X}(x,x') \le c, x\neq x'\big\rbrace
    \end{equation*}
is finite for any $c>0$. Then any map from a hyperdiscrete pseudo-metric space to any pseudo-metric space is coarsely Lipschitz. 

\noindent \textit{(v)} Let $X$ be a metric space, $Y$ a pseudo-metric space, and suppose there is $c>0$ so that $d_{X}(x,x')\ge c$ for all $x\neq x'\in X$. Then a map $f\colon X\rightarrow Y$ is large-scale Lipschitz if and only if it is Lipschitz. Indeed, suppose that $f$ is large-scale Lipschitz, which means there exist $c_{+}>0$, $c_{+}'\ge 0$ such that 
\begin{equation*}
    d_{Y}(f(x),f(x'))\le c_{+}\cdot d_{X}(x,x')+c_{+}'
\end{equation*}
for any $x,x'\in X$. Now $1\le \frac{d_{X}(x,x')}{c}$ for all $x,x'\in X$, and it follows that 
\begin{equation*}
    d_{Y}(f(x),f(x')) \le \left(c_{+}+\frac{c_{+}'}{c}\right)\cdot d_{X}(x,x')
\end{equation*}
for any $x,x'\in X$. 
\end{example}

The next lemma ensures that control functions are almost invertible. 

\begin{lemma}\label{lem:controlfunctionsinvertibility}
(i) Let $\Phi_{+}\colon \R_{+}\rightarrow \R_{+}$ be an upper control. The function $\Psi_{-}\colon \R_{+}\rightarrow \R_{+}\cup\lbrace\infty\rbrace$ defined by 
\begin{equation*}
    \Psi_{-}(s)\defeq \inf\lbrace r\ge 0 : \Phi_{+}(r)\ge s\rbrace, \; s\ge 0 
\end{equation*}
is a lower control such that $\Psi_{-}(\Phi_{+}(t)) \le t$ for all $t\ge 0$.

\noindent (ii) Let $\Phi_{-}\colon \R_{+}\rightarrow \R_{+}\cup\lbrace\infty\rbrace$ be a lower control. The function $\Psi_{+}\colon \R_{+}\rightarrow \R_{+}$ defined by 
\begin{equation*}
    \Psi_{+}(s)\defeq \sup\lbrace r\ge 0 : \Phi_{-}(r)\le s\rbrace, \; s\ge 0 
\end{equation*}
is an upper control such that $\Psi_{+}(\Phi_{-}(t)) \ge t$ for all $t\ge 0$.
\end{lemma}

\begin{proof}
We only prove~\textit{(i)} since the proof of~\textit{(ii)} is identical. Clearly $\Psi_{-}(s)\ge 0$ for any $s\ge 0$. Next, if $s_{1}\le s_{2}$, and if $r\ge 0$ is so that $\Phi_{+}(r)\ge s_{2}$, then also $\Phi_{+}(r)\ge s_{1}$, whence $\Psi_{-}(s_{1})\le r$. It follows that $\Psi_{-}(s_{1})\le \inf\lbrace r\ge 0 : \Phi_{+}(r)\ge s_{2}\rbrace = \Psi_{-}(s_{2})$, and $\Psi_{-}$ is non-decreasing. Lastly, if $t\ge 0$, we have 
\begin{equation*}
    \Psi_{-}(\Phi_{+}(t))=\inf\lbrace r\ge 0 : \Phi_{+}(r)\ge \Phi_{+}(t)\rbrace \le t
\end{equation*}
since $t\in \lbrace r\ge 0 : \Phi_{+}(r)\ge \Phi_{+}(t)\rbrace$. 
\end{proof}

\begin{proposition}\label{prop:morphismsincoarsecategory}
Let $X$ and $Y$ be pseudo-metric spaces, $f\colon X\rightarrow Y$ a coarsely Lipschitz map, and $\overline{f}$ the corresponding morphism in the metric coarse category. The following claims hold. 
\begin{enumerate}[label=(\roman*)]
    \item If $X\neq\emptyset$, $\overline{f}$ is an epimorphism if and only if $f$ is essentially surjective. 
    \item The morphism $\overline{f}$ is a monomorphism if and only if $f$ is coarsely expansive. 
    \item The morphism $\overline{f}$ is an isomorphism if and only if $f$ is a metric coarse equivalence. Moreover, if $X\neq\emptyset$, this holds if and only if $\overline{f}$ is an epimorphism and a monomorphism.
\end{enumerate}
\end{proposition}

\begin{proof}
\textit{(i)} Assume first that $f$ is essentially surjective, and let $c\defeq \dis\sup_{y\in Y}d_{Y}(y, f(X))$. Consider a pseudo-metric space $Z$ and two coarsely Lipschitz maps $h_{1},h_{2}\colon Y\rightarrow Z$ such that $h_{1}\circ f \sim_{c} h_{2}\circ f$. Hence there is $c'>0$ such that 
\begin{equation*}
    d_{Z}\big(h_{1}(f(x)), h_{2}(f(x))\big) \le c'
\end{equation*}
for any $x\in X$. Moreover, as $h_{1},h_{2}$ are coarsely Lipschitz, we can find $c_{1}, c_{2}>0$ so that 
\begin{equation*}
    d_{Y}(y,y') \le c \Longrightarrow d_{Z}\big(h_{1}(y),h_{1}(y')\big) \le c_{1}, \; d_{Y}(y,y') \le c \Longrightarrow d_{Z}\big(h_{2}(y),h_{2}(y')\big) \le c_{2}
\end{equation*}
for all $y,y'\in Y$. Let now $y\in Y$, and choose $x\in X$ with $d_{Y}(y,f(x)) \le c$. It then follows from the above implications that 
\begin{align*}
    d_{Z}\big(h_{1}(y), h_{2}(y)\big) &\le d_{Z}\big(h_{1}(y), h_{1}(f(x))\big)+d_{Z}\big(h_{1}(f(x)), h_{2}(f(x))\big)+d_{Z}\big(h_{2}(f(x)), h_{2}(y)\big) \\
    &\le c_{1}+c'+c_{2}.
\end{align*}
As $y\in Y$ was arbitrary, this proves that $h_{1}\sim_{c} h_{2}$, so $\overline{f}$ is an epimorphism. 

\noindent Conversely, suppose $f$ is not essentially surjective. Define $h_{1},h_{2}\colon Y\rightarrow \R_{+}$ by $h_{1}(y)=0$, and $h_{2}(y)=d_{Y}(y,f(X))$. Then $h_{1}\circ f=h_{2}\circ f=0$, so $h_{1}\circ f \sim_{c} h_{2}\circ f$, but $h_{1}\nsim_{c} h_{2}$, as $h_{1}(Y)=\lbrace 0\rbrace$ is bounded in $\R_{+}$ while $h_{2}(Y)$ is not. As $h_{1},h_{2}$ are coarsely Lipschitz (and thus morphisms in the metric coarse category), we deduce that $\overline{f}$ is not an epimorphism. 

\noindent \textit{(ii)} Suppose now that $f$ is coarsely expansive. Let $\Phi_{-}$ be a lower control for $f$ and $\Psi_{+}$ an upper control as in Lemma~\ref{lem:controlfunctionsinvertibility}. Let $W$ be a pseudo-metric space and let $h_{1},h_{2}\colon W \rightarrow X$ be two coarsely Lipschitz maps such that $f\circ h_{1}\sim_{c} f\circ h_{2}$. Hence there is $c>0$ such that 
\begin{equation*}
    d_{Y}\big(f(h_{1}(w)), f(h_{2}(w))\big) \le c
\end{equation*}
for any $w\in W$. It follows that 
\begin{align*}
    d_{X}(h_{1}(w), h_{2}(w)) &\le \Psi_{+}\big(\Phi_{-}(d_{X}(h_{1}(w), h_{2}(w)))\big) \\
    &\le \Psi_{+}\big(d_{Y}(f(h_{1}(w)), f(h_{2}(w)))\big) \\
    &\le \Psi_{+}(c)
\end{align*}
for any $w\in W$, which shows that $h_{1}\sim_{c} h_{2}$. Thus $\overline{f}$ is a monomorphism. 

\noindent Conversely, suppose that $f$ is not coarsely expansive. This implies there exist $c>0$ and $(x_{n})_{n\in \N}, (x_{n}')_{n\in\N} \subset X$ with $\lim\limits_{n\rightarrow\infty}d_{X}(x_{n},x_{n}')=\infty$ and 
\begin{equation*}
    d_{Y}(f(x_{n}),f(x_{n}')) \le c
\end{equation*}
for all $n\in\N$. Consider now $W\defeq \lbrace n^2 : n\in\N\rbrace$ endowed with the usual metric $(d_{W}(n^2,m^2)=|n^2-m^2|, n,m\in\N)$, and the maps $h_{1},h_{2}\colon W\rightarrow X$, $h_{1}(n^2)=x_{n}, h_{2}(n^2)=x_{n}'$. As $W$ is a hyperdiscrete pseudo-metric space, Example~\ref{ex:examplesofmaps}\textit{(iv)} guarantees that $h_{1}, h_{2}$ are coarsely Lipschitz. Now $f\circ h_{1} \sim_{c} f\circ h_{2}$, but $h_{1}\nsim_{c} h_{2}$ as $d_{X}(x_{n},x_{n}')\rightarrow\infty$ when $n\rightarrow\infty$. Hence $\overline{f}$ is not a monomorphism. 

\noindent \textit{(iii)} Suppose that $f$ is a metric coarse equivalence. If $X=Y=\emptyset$, there is nothing to prove. If $Y\neq \emptyset$, then so is $X$ (since the only map $\emptyset \rightarrow Y$ is not essentially surjective). Let $\Phi_{-}$, $\Phi_{+}$ be lower and upper controls for $f$ and let $c>0$ be such that $d_{Y}(y,f(X)) \le c$ for any $y\in Y$. Let $\Psi_{+}$ be as in Lemma~\ref{lem:controlfunctionsinvertibility} and let $\Psi_{-}$ be the lower control given by 
\begin{equation*}
    \Psi_{-}(s)=\inf\big\lbrace r\ge 0 : \Phi_{+}(r)+2c\ge s\big\rbrace, \; s\ge 0.
\end{equation*}
For each $y\in Y$, pick $x_{y}\in X$ such that $d_{Y}(y,f(x_{y}))\le c$, and define $g\colon Y\rightarrow X$ by $g(y)\defeq x_{y}$. Let $y,y'\in Y$. Then we have 
\begin{align*}
    d_{X}(g(y), g(y')) &\le \Psi_{+}\big(\Phi_{-}\big(d_{X}(g(y), g(y'))\big)\big) \\
    &\le \Psi_{+}\big(d_{Y}\big(f(g(y)), f(g(y'))\big)\big) \\
    &= \Psi_{+}\big(d_{Y}\big(f(x_{y}), f(x_{y'})\big)\big) \\
    &\le \Psi_{+}\big(d_{Y}(y,y')+2c\big)
\end{align*}
which proves that $s\longmapsto \Psi_{+}(s+2c)$ is an upper control for $g$, which is then coarsely Lipschitz. On the other hand, we have 
\begin{align*}
    d_{Y}(y,y') &\le d_{Y}\big(y,f(g(y))\big)+d_{Y}\big(f(g(y)), f(g(y'))\big)+d_{Y}\big(f(g(y')), y'\big) \\
    &\le \Phi_{+}\big(d_{X}(g(y),g(y')\big)+2c \\
    &=\widetilde{\Phi_{+}}\big(d_{X}(g(y),g(y'))\big)
\end{align*}
for any $y,y'\in Y$, where $\widetilde{\Phi_{+}}(s) \defeq \Phi_{+}(s)+2c$, $s\ge 0$. It follows that 
\begin{equation*}
    \Psi_{-}(d_{Y}(y,y')) \le \Psi_{-}\big(\widetilde{\Phi_{+}}\big(d_{X}(g(y),g(y'))\big)\big) \le d_{X}(g(y),g(y'))
\end{equation*}
for all $y,y'\in Y$. Therefore $g$ is coarsely expansive as well. By construction, we have $f\circ g \sim_{c} \text{Id}_{Y}$, so $f\circ g\circ f \sim_{c} f$, and as $\overline{f}$ is a monomorphism by~\textit{(ii)}, we conclude that $g\circ f\sim_{c} \text{Id}_{X}$, and finally that $\overline{f}$ is an isomorphism with inverse $\overline{g}$. 

\noindent Conversely, if $\overline{f}$ is an isomorphism, then $f$ is essentially surjective by~\textit{(i)} and coarsely expansive by~\textit{(ii)}, thus it is a metric coarse equivalence. 
\end{proof}

It follows from this result, and the fact that compositions of metric coarse equivalences, are metric coarse equivalences that being coarsely equivalent is an equivalence relation among pseudo-metric spaces. 

\begin{definition}\label{def:coarseinvertibility}
Let $X,Y$ be pseudo-metric spaces and let $f\colon X\rightarrow Y$ be a coarsely Lipschitz map. We say that $f$ is 
\begin{enumerate}[label=(\roman*)]
    \item~\textit{coarsely right-invertible} if there exists a coarsely Lipschitz map $g\colon Y\rightarrow X$ such that $f\circ g\sim_{c}\text{Id}_{Y}$.
    \item~\textit{coarsely left-invertible} if there exists a coarsely Lipschitz map $g\colon Y\rightarrow X$ such that $g\circ f\sim_{c}\text{Id}_{X}$.
    \item~\textit{coarsely invertible} if it is both left-invertible and right-invertible. 
\end{enumerate}
\end{definition}

\begin{remark}\label{rem:remark1.20}
If $f\colon X\rightarrow Y$ is coarsely Lipschitz and coarsely right-invertible, then $\overline{f}$ is an epimorphism. Dually, if $f$ is coarsely left-invertible, then $f$ is a monomorphism. In particular, $f$ is coarsely invertible if and only if $\overline{f}$ is an isomorphism, i.e. if and only if $f$ is a metric coarse equivalence. 
\end{remark}

\begin{definition}
Let $Y$ be a pseudo-metric space. A subspace $Z$ of $Y$ is a~\textit{coarse retract} of $Y$ if the inclusion map $i\colon Z\hookrightarrow Y$ is left-invertible, i.e. there exists a coarsely Lipschitz map $r\colon Y\rightarrow Z$ such that $r\circ i\sim_{c}\text{Id}_{Z}$. 
\end{definition}

Retractions provide an alternative characterization of coarse expansiveness and coarse left-invertibility.

\begin{proposition}
Let $X,Y$ be pseudo-metric spaces and $f\colon X\rightarrow Y$ be a coarsely Lipschitz map. Denote $f_{\text{im}}\colon X\rightarrow f(X)$ the map induced by $f$. Then the following holds.
\begin{enumerate}[label=(\roman*)]
    \item The map $f$ is coarsely expansive if and only if $f_{\text{im}}$ is a metric coarse equivalence. 
    \item The map $f$ is coarsely left-invertible if and only if it is coarsely expansive and $f(X)$ is a coarse retract.
\end{enumerate}
\end{proposition}

\begin{proof}
\textit{(i)} directly follows from the definitions.

\noindent~\textit{(ii)} Suppose that $f$ is coarsely left-invertible and let $g\colon Y\rightarrow X$ be a coarsely Lipschitz map such that $g\circ f\sim_{c} \text{Id}_{X}$. Denote $i$ the natural inclusion of $f(X)$ into $Y$. It follows from Remark~\ref{rem:remark1.20} that $\overline{f}$ is a monomorphism, so $f$ is coarsely expansive by Proposition~\ref{prop:morphismsincoarsecategory}. Now we have 
\begin{equation*}
    \text{Id}_{f(X)}\circ f_{\text{im}}=f\circ \text{Id}_{X} \sim f_{\text{im}}\circ (g\circ i\circ f_{\text{im}})
\end{equation*}
and since $f_{\text{im}}$ is a metric coarse equivalence by (i) we conclude that $\text{Id}_{X}\sim_{c} (f_{\text{im}}\circ g) \circ i$, and $f_{\text{im}}\circ g$ is a coarse retraction from $Y$ to $f(X)$. 

\noindent Conversely, suppose $f$ is coarsely expansive and that $f(X)$ is a coarse retract. Let $r\colon Y\rightarrow f(X)$ be a coarse retraction. By~\textit{(i)}, $f_{\text{im}}$ is a metric coarse equivalence, so let $j\colon f(X)\rightarrow X$ be a coarsely Lipschitz map such that $\overline{j}$ and $\overline{f}$ are inverses of each other. Then one has
\begin{align*}
    (j\circ r)\circ f &= j\circ (r\circ i)\circ f_{\text{im}} \\
    &\sim_{c} j\circ \text{Id}_{f(X)}\circ f_{\text{im}} \\
    &=j\circ f_{\text{im}}\\
    &\sim_{c} \text{Id}_{X}
\end{align*}
whence $j\circ r$ is a coarse left inverse for $f$. 
\end{proof}

\begin{definition}
Let $X,Y$ be pseudo-metric spaces. We say that $Y$ is~\textit{coarsely retractable} on $X$ if there is a coarsely right-invertible coarsely Lipschitz map from $Y$ to $X$, or equivalently if there is a coarsely left-invertible coarsely Lipschitz map from $X$ to $Y$.
\end{definition}

We conclude this part by the analog of Proposition~\ref{prop:morphismsincoarsecategory} in the large-scale category. 

\begin{proposition}\label{prop:morphismsinlargescalecategory}
Let $X,Y$ be pseudo-metric spaces, $f\colon X\rightarrow Y$ a large-scale Lipschitz map, and $\overline{f}$ the corresponding morphism in the large-scale category. The following equivalences hold.
\begin{enumerate}[label=(\roman*)]
    \item If $X\neq\emptyset$, $\overline{f}$ is an epimorphism if and only if $f$ is essentially surjective. 
    \item The morphism $\overline{f}$ is a monomorphism if and only if $f$ is large-scale expansive. 
    \item The morphism $\overline{f}$ is an isomorphism if and only if $f$ is a quasi-isometry. Moreover, if $X\neq\emptyset$, this holds if and only if $\overline{f}$ is an epimorphism and a monomorphism.
    \end{enumerate}
\end{proposition}

\subsection{Coarse and large-scale properties}

Let $(X,d_{X})$ be a pseudo-metric space and $c>0$. If $x,x'\in X$ and $n\ge 0$, a~\textit{$c-$path of $n$ steps from $x$ to $x'$ in $X$} is a sequence 
\begin{equation*}
    x=x_{0},x_{1},\dots,x_{n-1},x_{n}=x'
\end{equation*}
of points in $X$ such that $d_{X}(x_{i-1},x_{i}) \le c$ for all $i=1,\dots,n$.

\begin{definition}\label{def:coarseconnectednessinspaces}
Let $(X,d_{X})$ be a pseudo-metric space  and $c>0$. We say that $X$ is 
\begin{enumerate}[label=(\roman*)]
    \item~\textit{$c-$coarsely connected} if for any pair of points $x,x'\in X$, there is a $c-$path from $x$ to $x'$.
    \item~\textit{$c-$coarsely geodesic} if there exists an upper control $\Phi$ such that, for any pair of points $x,x'\in X$, there is a $c-$path of at most $\Phi(d_{X}(x,x'))$ steps from $x$ to $x'$. 
    \item~\textit{$c-$large-scale geodesic} if there exist $a>0$, $b\ge 0$ such that for any pair of points $x,x'\in X$, there is a $c-$path of at most $ad_{X}(x,x')+b$ steps from $x$ to $x'$.
    \item~\textit{$c-$geodesic} if for any pair of points $x,x'\in X$, there is a $c-$path $x=x_{0},x_{1},\dots,x_{n}=x'$ such that 
    \begin{equation*}
        d_{X}(x,x') = \sum_{i=1}^{n}d_{X}(x_{i-1},x_{i}).
    \end{equation*}
    \item~\textit{geodesic} if for any pair of points $x,x'\in X$ with $d_{X}(x,x')>0$, there exists an isometric map $\sigma\colon [0,d_{X}(x,x')]\rightarrow X$ such that $\sigma(0)=x$ and $\sigma(d_{X}(x,x'))=x'$.
\end{enumerate}
We say that $X$ is~\textit{coarsely connected} (resp.~\textit{coarsely geodesic},~\textit{large-scale geodesic}) if it is $c-$coarsely connected (resp. $c-$coarsely geodesic, $c-$large-scale geodesic) for some $c>0$. 
\end{definition}

\begin{remark}\label{rem:remark1.26}
\textit{(i)} Clearly, we have 
\begin{align*}
    X \; \text{geodesic} &\Longrightarrow X \;\text{$c-$geodesic} \\
    &\Longrightarrow X \;\text{large-scale geodesic} \\
    &\Longrightarrow X \;\text{coarsely geodesic} \\
    &\Longrightarrow X \;\text{coarsely connected}.
\end{align*}

\noindent \textit{(ii)} If $X$ is $c-$coarsely connected (resp. $c-$coarsely geodesic, $c-$large-scale geodesic or $c-$geodesic) for some $c>0$, then $X$ is $C-$coarsely connected (resp. $C-$coarsely geodesic, $C-$large-scale geodesic or $C-$geodesic) for any $C\ge c$. 
\end{remark}

\begin{proposition}\label{prop:invarianceofcoarseconnectedness}
Coarse connectedness, coarse geodesicity (resp. large-scale geodesicity) are invariant under metric coarse equivalences (resp. quasi-isometries).
\end{proposition}

\begin{proof}
We show the proof for large-scale geodesicity, and the others are completely similar. Suppose that $f\colon X\rightarrow Y$ is a $(C,K)-$quasi-isometry, with $C\ge 1$ and $K\ge 0$, and let $c>0$ be such that $X$ is $c-$large scale geodesic and any point of $Y$ is at distance at most $c$ from $f(X)$. Let $y,y'\in Y$ and let $x,x'\in X$ be such that
\begin{equation*}
    d_{Y}(y,f(x)), \; d_{Y}(y',f(x')) \le c. 
\end{equation*}
As $X$ is $c-$large scale geodesic, there exist $a>0$, $b\ge 0$ and a $c-$path 
\begin{equation*}
    x=x_{0},x_{1},\dots,x_{n}=x'
\end{equation*}
such that $n\le a\cdot d_{X}(x,x')+b$. Set 
\begin{equation*}
    y_{0}\defeq y, y_{1}\defeq f(x_{1}),y_{2}\defeq f(x_{2}),\dots, y_{n-1}\defeq f(x_{n-1}), y_{n}\defeq y'.
\end{equation*}
Then one has
\begin{align*}
    d_{Y}(y_{i-1},y_{i})&=d_{Y}(f(x_{i-1}), f(x_{i})) \\
    &\le C\cdot d_{X}(x_{i-1},x_{i})+K \\
    &\le C\cdot c+K
\end{align*}
for all $i=2,\dots,n-1$, and also
\begin{align*}
    d_{Y}(y_{0},y_{1})&=d_{Y}(y,f(x_{1})) \\
    &\le d_{Y}(y,f(x))+d_{Y}(f(x_{0}),f(x_{1})) \\ 
    &\le c+C\cdot c+K \\
    &=(C+1)c+K
\end{align*}
and 
\begin{align*}
    d_{Y}(y_{n-1},y_{n})&=d_{Y}(f(x_{n-1}),y') \\
    &\le d_{Y}(f(x_{n-1}), f(x_{n}))+d_{Y}(f(x'),y') \\ 
    &\le c+C\cdot c+K \\
    &=(C+1)c+K.
\end{align*}

Thus $y=y_{0},y_{1},\dots,y_{n}=y'$ is a $\big((C+1)c+K\big)-$path between $y$ and $y'$ in $Y$: 

\begin{figure}[H]
  \centering
  \includegraphics[
    width=0.75\linewidth
  ]{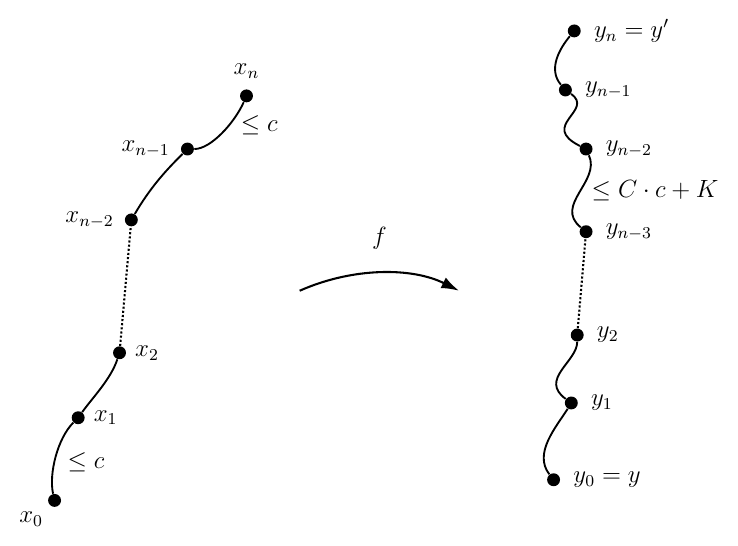}
  \caption{Configuration from the proof of Proposition~\ref{prop:invarianceofcoarseconnectedness}.}
  \label{fig:cayley}
\end{figure}

It is made of at most 
\begin{align*}
    n &\le a\cdot d_{X}(x,x')+b \le a\big(Cd_{Y}(f(x),f(x'))+CK\big)+b\\
    &\le aC\big(2c+d_{Y}(y,y')\big)+aCK+b
\end{align*}
steps, and the latter is indeed an affine upper bound on the length $n$ of the path in term of the distance $d_{Y}(y,y')$ between $y$ and $y'$. As $y,y'\in Y$ were arbitrary, it follows that $Y$ is large-scale geodesic as well. 
\end{proof}

Large-scale geodesicity can also be used to boost coarse properties for maps to large-scale properties.

\begin{proposition}\label{prop:boostingcoarsemapstolargescalemaps}
Let $X,Y$ be pseudo-metric spaces and $f\colon X\rightarrow Y$ a map. 
\begin{enumerate}[label=(\roman*)]
    \item If $X$ is large-scale geodesic and $f$ is coarsely Lipschitz, then $f$ is large-scale Lipschitz.
    \item If $X,Y$ are large-scale geodesic and $f$ is a metric coarse equivalence, then $f$ is a quasi-isometry.
\end{enumerate}
\end{proposition}

\begin{proof}
\textit{(i)} Assume that $X$ is $c-$large-scale geodesic, and let $a>0$, $b\ge 0$ be such that any pair of points $x,x'\in X$ can be joined by a $c-$path of at most $a\cdot d_{X}(x,x')+b$ steps. As $f$ is coarsely Lipschitz, we find $C\ge 0$ such that 
\begin{equation}\label{eq3}
    d_{X}(x,x') \le c \Longrightarrow d_{Y}(f(x),f(x'))\le C.
\end{equation}
Let $x,x'\in X$ and choose a $c-$path $x=x_{0},x_{1},\dots,x_{n}=x'$ from $x$ to $x'$ of at most $n\le a\cdot d_{X}(x,x')+b$ steps. Then 
\begin{align*}
    d_{Y}(f(x),f(x'))&\le \sum_{i=1}^{n}d_{Y}(f(x_{i-1}),f(x_{i})) \\
    &\le C\cdot n \\
    &\le C\big(a\cdot d_{X}(x,x')+b\big)\\
    &=(aC)\cdot d_{X}(x,x')+bC
\end{align*}
where the first inequality follows from the triangle inequality and the second one follows from (\ref{eq3}). Hence $f$ is large-scale Lipschitz. 

\noindent \textit{(ii)} follows directly from~\textit{(i)} applied to $f$ and to $g\colon Y\rightarrow X$ a metric coarse equivalence so that $g\circ f\sim_{c} \text{Id}_{X}$ and $f\circ g\sim_{c}\text{Id}_{Y}$.
\end{proof}

We also notice that checking a given map is Lipschitz is simpler in the particular case of connected graphs. 

\begin{lemma}\label{lem:Lipschitzbetweengraphs}
Let $f\colon X\rightarrow Y$ be a map between two connected graphs. If there exists $C>0$ such that $d(f(x),f(y)) \le C$ for any pair $(x,y)$ of adjacent vertices, then $f$ is $C-$Lipschitz. 
\end{lemma}

\begin{proof}
Let $x,y\in X$ be two vertices, and pick $x=x_{0},x_{1},\dots,x_{n-1},x_{n}=y$ a geodesic connecting $x$ and $y$ in $X$, so that $n=d(x,y)$. From the triangle inequality and the assumption, we get 
\begin{equation*}
    d(f(x),f(y))=d(f(x_{0}),f(x_{n})) \le \sum_{i=0}^{n-1}d(f(x_{i}), f(x_{i+1})) \le C\cdot n=C\cdot d(x,y).
\end{equation*}
Thus $f$ is $C-$Lipschitz, as claimed. 
\end{proof}

Let us conclude this part defining neighbourhoods of subsets in metric spaces and how they behave when applying coarse maps. 

Given a metric space $(X,d_{X})$ and a real number $R\ge 0$, the $R-$\textit{neighbourhood} of a subset $A\subset X$ is the subset $A^{+R}$ given by
\begin{equation*}
    A^{+R} \defeq \bigcup_{a\in A} \lbrace x\in X : d_{X}(x,a) \le R\rbrace=\bigcup_{a\in A}B_{X}(a,R).
\end{equation*}
Given two subsets $A,B\subset X$, we define then the \textit{Hausdorff distance} between $A$ and $B$, denoted $d_{\text{Haus}}(A,B)$, as 
\begin{equation*}
    d_{\text{Haus}}(A,B) \defeq \inf\lbrace R\ge 0 : A\subset B^{+R}, B\subset A^{+R}\rbrace
\end{equation*}
where, by convention, the infimum is $+\infty$ if there is no such $R$. In this way, the set $\mathcal{P}(X)$ with the Hausdorff distance $d_{\text{Haus}}$ is an extended pseudo-metric space.

\begin{lemma}\label{lem:neighborhoodsandQI}
Let $X, Y$ be metric spaces, and let $f\colon X\rightarrow Y$ be a $(C,K)-$quasi-isometry. Let $R,S\ge 0$. The following properties hold:
\begin{enumerate}[label=(\roman*)]
    \item If $R\le S$, then $A^{+R}\subset A^{+S}$ for any $A\subset X$;
    \item $(A^{+R})^{+S}\subset A^{+(R+S)}$ for any $A\subset X$; 
    \item For any $A\subset X$, $f(A^{+R}) \subset f(A)^{+(CR+K)}$;
    \item For any $A,B\subset X$, $d_{\text{Haus}}(f(A),f(B))\le C\cdot d_{\text{Haus}}(A,B)+K$.
\end{enumerate}
\end{lemma}

\begin{proof}
Claim~\textit{(i)} is immediate. For~\textit{(ii)}, if we take $x\in (A^{+R})^{+S}$, then there is $y\in A^{+R}$ such that $d_{X}(x,y)\le S$, and there is $a\in A$ such that $d_{X}(y,a)\le R$, whence
\begin{equation*}
    d_{X}(x,a) \le d_{X}(x,y)+d_{X}(y,a)\le R+S
\end{equation*}
and thus $x\in A^{+(R+S)}$.

\noindent For~\textit{(iii)}, fix $A\subset X$, and let $y\in f(A^{+R})$. Write $y=f(x)$ with $x\in A^{+R}$. Let thus $a\in A$ be such that $d_{X}(x,a)\le R$. Then one has 
\begin{equation*}
    d_{Y}(y,f(a)) = d_{Y}(f(x), f(a)) \le C\cdot d_{X}(x,a)+K \le C\cdot R+K
\end{equation*}
and $f(a)\in f(A)$, whence $y\in f(A)^{+(CR+K)}$. 

\noindent Lastly, for~\textit{(iv)}, write $Q\defeq d_{\text{Haus}}(A,B)$, so that $B\subset A^{+Q}, A\subset B^{+Q}$, and thus by~\textit{(iii)}, one gets 
\begin{equation*}
    f(B) \subset f(A^{+Q}) \subset f(A)^{+(CQ+K)}, \; f(A) \subset f(B^{+Q}) \subset f(B)^{+(CQ+K)} 
\end{equation*}
and the claim follows. 
\end{proof}

\subsection{Groups as pseudo-metric spaces}\label{sec:groupsaspseudometricspaces}

In this part, we explain how finitely generated groups can be seen as objects in the metric coarse category. Additionally, this provides us with numerous examples of metric coarse equivalences and quasi-isometries. 

\begin{definition}\label{def:lengthandmetricongroups}
Let $G$ be a finitely generated group with a finite generating set $S_{G}$. The~\textit{length} of an element $g\in G$ with respect to $S_{G}$ is 
\begin{equation*}
    |g|_{S_{G}} \defeq \min\big\lbrace n\in \N : \exists s_{1},\dots, s_{n}\in S_{G}\cup S_{G}^{-1}, g=s_{1}\dots s_{n}\big\rbrace
\end{equation*}
and the \textit{word metric} associated to $S_{G}$ is the map $d_{S_{G}}(\cdot,\cdot)\colon G\times G\longrightarrow [0,\infty)$ defined as
\begin{equation*}
    d_{S_{G}}(g,h) \defeq |g^{-1}h|_{S_{G}}.
\end{equation*}
\end{definition}

The next proposition justifies the terminology. 

\begin{proposition}\label{prop:fggroupismetricspace}
Let $G$ be a finitely generated group with a finite generating set $S_{G}$. The pair $(G,d_{S_{G}})$ is a metric space. Moreover, $d_{S_{G}}$ is left-invariant.
\end{proposition}

Interestingly, such an abstract group has a visual representation, called its~\textit{Cayley graph}. 

\begin{definition}\label{def:Cayleygraph}
Let $G$ be a finitely generated group with a finite generating set $S_{G}$. Its \textit{Cayley graph with respect to $S_{G}$} is the graph $\text{Cay}(G,S_{G})$ 
\begin{enumerate}[label=(\roman*)]
    \item whose vertices are elements of $G$;
    \item whose edges connect two distinct elements $g,h\in G$ if there exists $s\in S_{G}\cup S_{G}^{-1}$ such that $h=gs$. 
\end{enumerate}
\end{definition}

The fact that $S_{G}$ generates $G$ ensures then that $\text{Cay}(G,S_{G})$ is connected. 

\begin{figure}[H]
  \centering
  \includegraphics[width=0.7\linewidth]{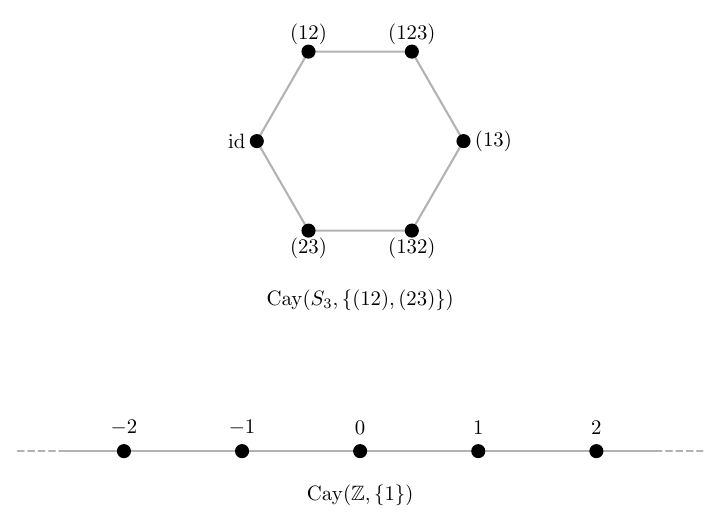}
  \caption{Two examples of Cayley graphs.}
  \label{fig:cayley2}
\end{figure}

\begin{figure}[H]
  \centering
  \includegraphics[width=0.7\linewidth]{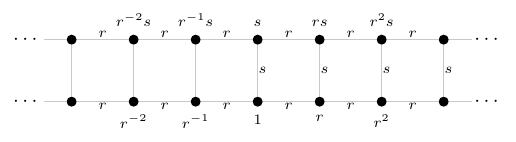}
  \caption{$\text{Cay}(D_{\infty}, \lbrace s,r\rbrace)$.}
  \label{fig:cayley2}
\end{figure}

\begin{example}\label{ex:descriptionoflamplighters}
Recall from Proposition~\ref{prop:wreathproductspreservefinitegeneration} that if $A$ is generated by $S=\lbrace a_{1},\dots, a_{n}\rbrace$ and $B$ is generated by $T=\lbrace b_{1},\dots, b_{m}\rbrace$, then 
\begin{equation*}
    U \defeq \lbrace (\delta_{a_{i}},e_{B}), (\textbf{1}, b_{j}) : 1\le i\le n, 1\le j\le m\rbrace
\end{equation*}
generates the wreath product $A\wr B$. The classical interpretation of such generating sets is as follows. First, think of an element $c\in\bigoplus_{B}A$ as a colouring of the vertices of $\text{Cay}(B,T)$ with colors coming from $A$, with only finitely many vertices having a non-trivial color (and these vertices form the~\textit{support} of $c$, denoted $\text{supp}(c)$). Second, think of an element $(c,p)\in A\wr B$ as a pair made of a colouring $c\in\bigoplus_{B}A$ together with an arrow pointing at some vertex $p\in B$.

\begin{figure}[H]
  \centering
  \includegraphics[width=0.75\linewidth, trim=0 0 0 4mm, clip]{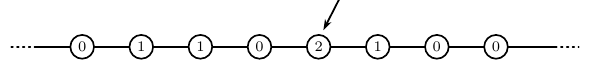}
  \caption{An element of $\Z \wr \Z$.}
  \label{fig:cayley}
\end{figure}

There are two possible moves to go from $(c,p)$ to a neighbouring vertex in $\text{Cay}(A\wr B,U)$: 
\begin{itemize}
    \item either we only move the arrow to a neighbouring vertex in $B$, and the colouring stays the same;
    \item or the arrow stays on the vertex where it stands, but changes the color of this vertex, and replaces it with an adjacent color in $\text{Cay}(A,S)$.
\end{itemize}

\begin{figure}[H]
  \centering
  \includegraphics[width=0.8\linewidth, trim=0 0 0 4mm, clip]{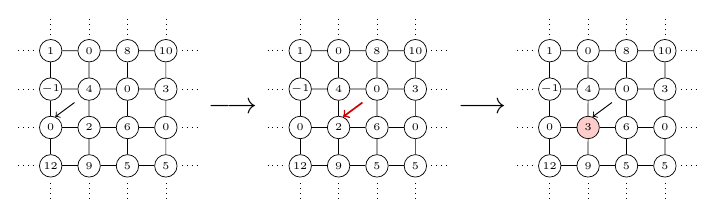}
  \caption{A path of length two in $\text{Cay}(\Z\wr\Z^2,U)$.}
  \label{fig:cayley}
\end{figure}

Hence to go from $(c,p)$ to another vertex $(d,q)$ in $A\wr B$, the arrow has to move in $B$ from $p$ to $q$ by visiting all vertices of $B$ where $c$ and $d$ differ. At each of these vertices $t\in B$, the color is changed, and goes from $c(t)\in A$ to $d(t)\in A$.
\end{example}

Next, we want to ensure that the choice of a finite generating set for a finitely generated group does not affect its large-scale geometry.

\begin{lemma}\label{lem:changeofgensetsBILIP}
Let $G$ be a finitely generated group, and let $S,T\subset G$ be two finite generating sets of $G$. Then the identity map $(G,d_{S})\rightarrow (G,d_{T})$ is a biLipschitz equivalence.
\end{lemma}

\begin{proof}
Let $g,h\in G$, and let $n\defeq d_{S}(g,h)$. We write then $g^{-1}h=s_{1}\dots s_{n}$ with $s_{1},\dots,s_{n}\in S\cup S^{-1}$. Then, by left-invariance of $d_{T}$, one has 
\begin{align*}
    d_{T}(g,h)&=d_{T}(g,gs_{1}\dots s_{n}) \\
    &=d_{T}(e_{G},s_{1}\dots s_{n}) \\
    &\le d_{T}(e_{G},s_{1})+d_{T}(s_{1}, s_{1}\dots s_{n}) \\
    &=d_{T}(e_{G},s_{1})+d_{T}(e_{G}, s_{2}\dots s_{n}) \\
    &\le \dots \\
    &\le \sum_{i=1}^{n}d_{T}(e_{G},s_{i}) \\
    &\le C\cdot n \\
    &= C\cdot d_{S}(g,h)
\end{align*}
where $C\defeq \displaystyle\max_{s\in S}d_{T}(e_{G},s)>0$. By symmetry, one proves that $d_{S}(g,h) \le C'\cdot d_{T}(g,h)$ where $C'\defeq \displaystyle\max_{t\in T}d_{S}(e_{G},t)>0$. Thus we conclude that 
\begin{equation*}
    \frac{1}{C'}\cdot d_{S}(g,h) \le d_{T}(g,h) \le C\cdot d_{S}(g,h)
\end{equation*}
and the claim follows. 
\end{proof}

\subsection{The Milnor-Schwarz lemma}

The next result is an efficient tool to prove the existence of quasi-isometries between groups and geodesic spaces on which they act. To state it, we recall additional terminologies related to group actions. 

\begin{definition}\label{def:propertiesofgroupactions}
Let $G$ be a group acting on a metric space $(X,d_{X})$. The action is called
\begin{enumerate}[label=(\roman*)]
    \item~\textit{isometric} if $d_{X}(g\cdot x, g\cdot y)=d_{X}(x,y)$ for any $g\in G$ and $x,y\in X$;
    \item~\textit{properly discontinuous} if for any compact subset $K\subset X$, the set 
    \begin{equation*}
        \lbrace g\in G : gK\cap K\neq\emptyset\rbrace
    \end{equation*}
    is finite;
    \item~\textit{cocompact} if there exists $x_{0}\in X$, there is $R\ge 0$ such that any $x\in X$ lies at distance $\le R$ from a point of $G\cdot x_{0}$.
    \item~\textit{geometric} if it is isometric, properly discontinuous and cocompact.
\end{enumerate}
\end{definition}

Recall also that a metric space $(X,d_{X})$ is called~\textit{proper} if closed balls are compact.

\begin{theorem}[{Milnor-Schwarz lemma}]\label{lem:Milnor-Schwarz}
Let $G$ be a group and $X$ be a proper geodesic metric space. If $G$ admits a geometric action on $X$, then $G$ is finitely generated and the orbit map \begin{align*}
    i_{x_{0}}\colon G &\longrightarrow X \\
    g&\longmapsto g\cdot x_{0}
\end{align*}
is a quasi-isometry for every $x_{0}\in X$. 
\end{theorem}

\begin{proof}
Let $x_{0}\in X$ be a base point. Let $R\ge 0$ be such that $G-$translates of $B_{X}(x_{0},R)$ cover $X$. Let 
\begin{equation*}
    S \defeq \big\lbrace g\in G : B_{X}(gx_{0},R)\cap B_{X}(x_{0},R)\neq\emptyset\big\rbrace.
\end{equation*}
Since the action is properly discontinuous, and that $X$ is proper, $S$ is finite. We claim that $S$ generates $G$. 

\noindent Write $B\defeq B(x_{0},R)$ and $c\defeq \inf\lbrace d(B,gB) : g\neq e_{G}, g\notin S\rbrace >0$. Let $g\in G$ with $g\notin S\cup\lbrace e_{G}\rbrace$. Then $d(x_{0}, gx_{0}) \ge 2R+c \ge R+c$, so there is $k\ge 2$ such that 
\begin{equation*}
    R+(k-1)c \le d(x_{0}, gx_{0}) < R+kc.
\end{equation*}
Thus we can choose points $x_{1}, x_{2},\dots, x_{k+1}=gx_{0}$ on a geodesic segment from $x_{0}$ to $gx_{0}$ such that $d(x_{0},x_{1}) \le R$ and $d(x_{i}, x_{i+1})<c$ for $1\le i\le k$. Then there are elements $e_{G}=g_{0}, g_{1}, \dots, g_{k}=g$ in $G$ such that $x_{i+1}\in g_{i}B$ for $0\le i\le k$:

\begin{figure}[H]
  \centering
  \includegraphics[width=0.75\linewidth]{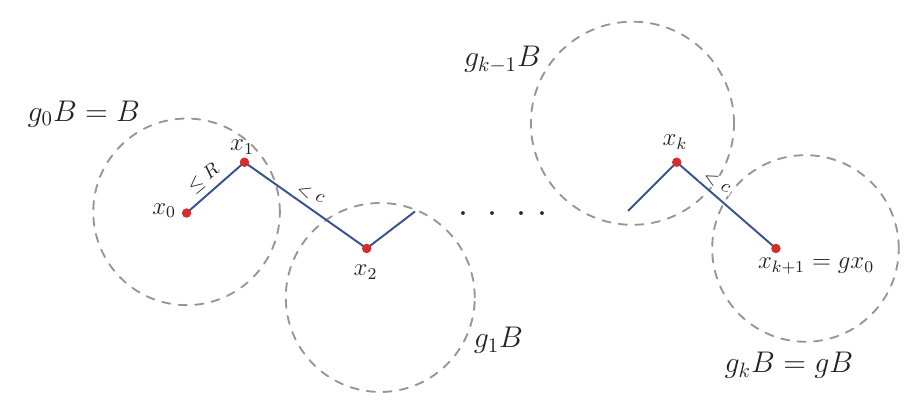}
  \caption{Configuration of the sequence $x_{0},x_{1}, x_{2},\dots, x_{k+1}$ in the space $X$}
\end{figure}
Let then $s_{i}\defeq g_{i-1}^{-1}g_{i}$ for $1\le i\le k$. This way, we get
\begin{equation*}
    d(B,s_{i}B)=d(B, (g_{i-1}^{-1}g_{i})B)=d(g_{i-1}B, g_{i}B) \le d(x_{i}, x_{i+1})<c
\end{equation*}
whence $s_{i}\in S$. Since $s_{1}s_{2}\dots s_{k}=g_{k}=g$, it follows that $S$ generates $G$, so $G$ is finitely generated. 

\noindent Now we show that $i_{x_{0}}$ is a quasi-isometry, when $G$ is equipped with the word metric $d_{S}$. Since any point of $X$ is at distance $\le R$ from a point of $G\cdot x_{0}$ (by cocompactness), $i_{x_{0}}$ is essentially surjective. We only need to show that there are constants $C\ge 1$ and $K\ge 0$ such that 
\begin{equation*}
    \frac{1}{C}\cdot d_{S}(g,h)-K \le d(gx_{0}, hx_{0}) \le C\cdot d_{S}(g,h)+K
\end{equation*}
for any $g,h\in G$. Since one has $d_{S}(g,h)=d_{S}(e_{G}, g^{-1}h)$ (by left-invariance of $d_{S}$) and $d(gx_{0}, hx_{0})=d(x_{0}, (g^{-1}h)x_{0})$ (the action $G\curvearrowright X$ is isometric), it is enough to prove that 
\begin{equation}\label{eq:inequalityQI}
    \frac{1}{C}\cdot d_{S}(e_{G},g)-K \le d(x_{0}, gx_{0}) \le C\cdot d_{S}(e_{G},g)+K
\end{equation}
for all $g\in G$. Let $L\defeq \max\lbrace d(x_{0}, sx_{0}) : s\in S\rbrace$ and set 
\begin{equation*}
    C \defeq \max\left\lbrace \frac{1}{c}, L, 2R\right\rbrace, \; K \defeq \max\left\lbrace \frac{1}{C},c\right\rbrace.
\end{equation*}
If $g=e_{G}$, then (\ref{eq:inequalityQI}) clearly holds. If rather $g=s\in S$, then we have $0\le d(x_{0}, sx_{0})\le 2R$ and $d_{S}(e_{G},s)=1$, and thus 
\begin{equation*}
    \frac{1}{C}\cdot d_{S}(e_{G},s)-K = \frac{1}{C}-K \le 0 \le d(x_{0},sx_{0})\le 2R \le C \le C\cdot d_{S}(e_{G}, s)+K
\end{equation*}
which is exactly (\ref{eq:inequalityQI}). Suppose lastly that $g\notin S\cup\lbrace e_{G}\rbrace$. From the proof that $S$ generates $G$, we get that, if $k$ is the largest integer such that $R+(k-1)c \le d(x_{0},gx_{0})$ then $d_{S}(e_{G},g) \le k$. In particular 
\begin{equation*}
    R+(d_{S}(e_{G},g)-1)c \le d(x_{0}, gx_{0})
\end{equation*}
and it follows that 
\begin{equation*}
    c\cdot d_{S}(e_{G},g)-c \le d(x_{0}, gx_{0}) - R \le d(x_{0}, gx_{0}).
\end{equation*}
Recalling that $L$ is the maximum of $d(x_{0},sx_{0})$ as $s$ runs over $S$, the same argument as in the proof of Lemma~\ref{lem:changeofgensetsBILIP} shows that $d(x_{0}, gx_{0}) \le L\cdot d_{S}(e_{G},g)$, and thus we obtain 
\begin{equation*}
    c\cdot d_{S}(e_{G},g)-c \le d(x_{0}, gx_{0}) \le L\cdot d_{S}(e_{G},g).
\end{equation*}
Since $C\ge L$ and $C\ge \frac{1}{c}$, $K\ge c$, inequality (\ref{eq:inequalityQI}) holds. This finishes the proof. 
\end{proof}

\begin{example}\label{ex:applicationsofMilnor-Schwarz}
\textit{(i)} For any $n\ge 1$, the natural action of $\Z^n$ on $\R^n$ is isometric, proper and cocompact. As $\R^n$ is geodesic and proper, it follows that $\Z^n\sim_{Q.I.} \R^{n}$, as in Example~\ref{ex:QI}\textit{(i)}.

\noindent \textit{(ii)} If $G$ is finitely generated and $S$ is a finite symmetric generating set for $G$, then the natural action of $G$ on its Cayley graph $\text{Cay}(G,S)$ is isometric, proper, and cocompact since it is transitive on the vertices and there are $|S|$ equivalence classes of edges. Thus $G$ is quasi-isometric to $\text{Cay}(G,S)$. 

\noindent \textit{(iii)} The Cayley graph of $\Z_{2}*\Z_{2}*\Z_{2}*\Z_{2}=\langle a,b,c,d : a^2=b^2=c^2=d^2=1\rangle$ with respect to $S=\lbrace a,b,c,d\rbrace$ is a $4-$regular tree, and is therefore quasi-isometric to any Cayley graph of $F_{2}$. Thus $\Z_{2}*\Z_{2}*\Z_{2}*\Z_{2}$ is quasi-isometric to $F_{2}$. 
\end{example}

We now deduce corollaries of interest to get more examples of pairs of quasi-isometric groups. 

\begin{corollary}\label{cor:QIsubgroupsoffiniteindex}
Let $G$ be a finitely generated group and $H\leqslant G$ a finite index subgroup. Then $H$ is finitely generated and quasi-isometric to $G$.   
\end{corollary}

\begin{proof}
Consider $S$ a finite generating set for $G$ and the metric space $(G,d_{S})$. Let $H$ act on $(G,d_{S})$ by left-multiplication. This action is isometric, proper, and cocompact since a finite set of representatives of left $H-$cosets is a compact subset of $(G,d_{S})$ whose translates by $H$ cover $G$. Moreover, $(G,d_{S})$ is geodesic and proper (balls of finite radius centered at $e_{G}\in G$ are finite), whence $H$ is finitely generated and quasi-isometric to $G$ by Theorem~\ref{lem:Milnor-Schwarz}. Moreover, a quasi-isometry is given by an arbitrary orbit map, for an arbitrary choice of base point in $G$. The choice $e_{G}\in G$ shows that the natural inclusion $H\hookrightarrow G$ is a quasi-isometry. 
\end{proof}

\begin{example}\label{ex:finiteindexsubgroupsQItogroup}
\textit{(i)} The dihedral group $D_{\infty}=\langle a, t : a^2=1, ata^{-1}=t^{-1}\rangle \cong \Z\rtimes \Z/2\Z$ contains $\Z$ as a finite index subgroup, and thus is quasi-isometric to $\Z$. 

\noindent \textit{(ii)} The group $\text{SL}_{2}(\Z)$ contains a finite index subgroup isomorphic to $F_{2}$ (see e.g.~\cite[Proposition~4.4.2]{Loh17}), so $\text{SL}_{2}(\Z)$ is quasi-isometric to $F_{2}$. 
\end{example}

\begin{corollary}\label{cor:groupsQItoquotients}
Let $G$ be a finitely generated group and $N\lhd G$ be a finite normal subgroup. Then $G$ is quasi-isometric to $G/N$. 
\end{corollary}

\begin{proof}
The natural action of $G$ on $G/N$ satisfies all assumptions of the Milnor-Schwarz lemma, whence the claim.
\end{proof}

This implies for instance that $\text{SL}_{2}(\Z)$ is quasi-isometric to $\text{PSL}_{2}(\Z)$ since the latter is the quotient of $\text{SL}_{2}(\Z)$ by its center $\lbrace \pm I_{2}\rbrace$.

For the last application, we need a terminology.

\begin{definition}\label{def:commensurability}
Let $G$ and $H$ be two groups. We say that they are
\begin{enumerate}[label=(\roman*)]
    \item~\textit{commensurable} if they contain finite index subgroups $G'\leqslant G$, $H'\leqslant H$ such that $G'\cong H'$.
    \item~\textit{weakly commensurable} if they contain finite index subgroups $G'\leqslant G$, $H'\leqslant H$ with finite normal subgroups $N\lhd G'$, $M\lhd H'$ such that $G'/N\cong H'/M$. 
\end{enumerate}
\end{definition}

The next statement is then a direct consequence of our previous results.

\begin{corollary}\label{cor:weakcommensurabilityQI}
Let $G$ be a finitely generated group. If $H$ is weakly commensurable to $G$, then $H$ is finitely generated and quasi-isometric to $G$.
\end{corollary}

\begin{proof}
Assume $H$ is weakly commensurable to $G$, and let $G',H',N,M$ be as in Definition~\ref{def:commensurability}. As $G$ is finitely generated, we deduce from Corollary~\ref{cor:QIsubgroupsoffiniteindex} that $G'$ is finitely generated, and thus $G'/N$ is finitely generated (Proposition~\ref{prop:finitegenerationextension}). Hence $H'/M$ is finitely generated, so that $H'$ is finitely generated. Thus $H$ is finitely generated, and we have 
\begin{equation*}
    H \sim_{Q.I} H' \sim_{Q.I.} H'/M \cong G'/N \sim_{Q.I.} G' \sim_{Q.I.} G
\end{equation*}
where the first and last quasi-isometries are given by Corollary~\ref{cor:QIsubgroupsoffiniteindex}, and the second and the third quasi-isometries are given by Corollary~\ref{cor:groupsQItoquotients}. The proof is complete.
\end{proof}

\section{Coarse geometric invariants}\label{subsection1.3}

The goal of this section is to develop powerful tools to be able to distinguish pseudo-metric spaces up to quasi-isometry. 

\subsection{Growth of finitely generated groups}\label{subsubsectionGROWTH}

\begin{definition}\label{def:growth}
Let $G$ be a finitely generated group with a finite generating set $S$. The~\textit{growth function} of $G$ with respect to $S$ is the map $\beta_{(G,S)}\colon \N\rightarrow \N$ defined as 
\begin{equation*}
    \beta_{(G,S)}(n)\defeq |B_{d_{S}}(e_{G},n)|. 
\end{equation*}
\end{definition}

\begin{example}\label{ex:examplesofgrowthfunctions}
\textit{(i)} Let $G=\Z$ equipped with the metric $d_{S}$ where $S=\lbrace -1,1\rbrace$. Then, if $n\ge 1$, $B_{d_{S}}(0,n)=\lbrace -n,-n+1,\dots, n-1,n\rbrace$, whence $\beta_{(\Z,S)}(n)=2n+1$ for any $n\ge 1$.

\noindent \textit{(ii)} Let $G=\Z^2$ equipped with the generating set $S=\lbrace (1,0),(0,1),(-1,0),(0,-1)\rbrace$. For $n\ge 0$, the ball of radius $n$ centered at $(0,0)$ is 
\begin{equation*}
    \lbrace (i,j)\in\Z^2 : |i|+|j| \le n\rbrace
\end{equation*}
the diagonal square containing the vertices $(0,r), (0,-r),(r,0),(-r,0)$, $0\le r\le n$. Thus 
\begin{equation*}
    \beta_{(\Z^2,S)}(n)=1+\sum_{r=1}^{n}4r=2n^2+2n+1
\end{equation*}
for all $n\ge 0$. 

\noindent \textit{(iii)} If $G=\Z^2$ is rather endowed with $S'=S\cup \lbrace(1,1),(-1,-1),(1,-1),(-1,1)\rbrace$, then for any $n\ge 0$, the ball of radius $n$ centered at $(0,0)$ is now 
\begin{equation*}
    \lbrace (i,j)\in\Z^2 : |i|\le n, |j|\le n\rbrace=\lbrace -n,\dots,n\rbrace^{2}
\end{equation*}
so that $\beta_{(\Z^2, S')}(n)=(2n+1)^{2}=4n^2+4n+1$ for any $n\ge 0$.

\noindent \textit{(iv)} Let $G=F_{2}$ equipped with $S=\lbrace a,b,a^{-1},b^{-1}\rbrace$.

\begin{figure}[H]
  \centering
  \includegraphics[width=0.6\linewidth]{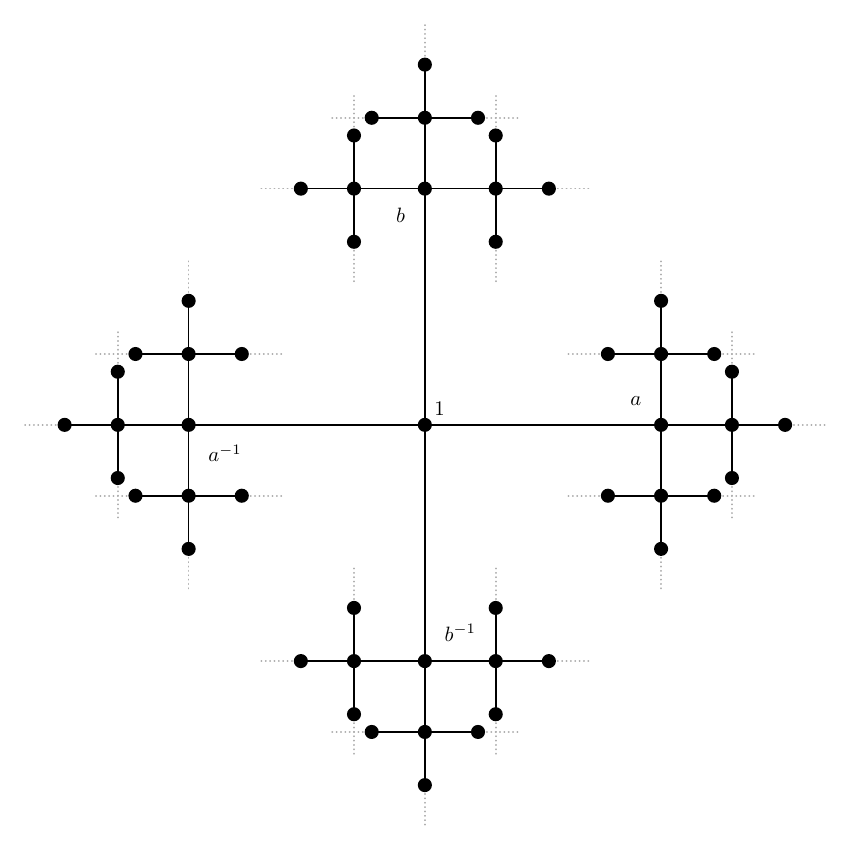}
  \caption{Cay$(F_2, S)$.}
  \label{fig:f2_cayley}
\end{figure}

For all $n\ge 1$, the ball of radius $n$ centered at the identity element has cardinality 
\begin{equation*}
    \beta_{(F_{2},S)}(n)=1+4\cdot \sum_{j=0}^{n-1}3^{j}=2\cdot 3^{n}-1.
\end{equation*}
\end{example}

The previous examples show that the exact values of the growth function depend on the choice of $S$. However, we know that, for different finite generating sets of the same group, the two corresponding metric spaces are related by a biLipschitz equivalence (cf. Lemma~\ref{lem:changeofgensetsBILIP}), and therefore have the same large-scale geometry. In order to make the growth function a geometric invariant, we need to introduce an equivalence relation on the set of functions defined on $\N$.

\begin{definition}\label{def:asymptoticequivalence}
Let $f,g\colon \N\rightarrow \R_{+}$. We say that $g$~\textit{dominates} $f$, and we write $f\preccurlyeq g$, if there is a constant $C>0$ such that $f(n) \le C\cdot g(Cn)$ for any $n$ large enough. We say that $f$ and $g$ are \textit{asymptotically equivalent}, written $f(n)\simeq g(n)$, if $f(n)\preccurlyeq g(n)$ and $g(n)\preccurlyeq f(n)$.
\end{definition}

It is easy to check that $\simeq$ is indeed an equivalence relation, and given a function $f$, its~\textit{class} refers to its equivalence class modulo $\simeq$. If no confusion is possible, we often write $f$ for the class of a function $f$. 

\begin{example}\label{ex:examplesofasymptoticequivalence}
\textit{(i)} Let $a,b>0$ and $c,d>1$. One has $n^{a} \preccurlyeq n^{b}$ if and only if $a\le b$ (and thus $n^{a}\simeq n^{b}$ if and only if $a=b$), and $c^{n}\simeq d^{n}$. 

\begin{proof}
We start with the first equivalence. If $a\le b$, then $n^{a}\le n^{b}$ for any $n\ge 0$, so $n^{a}\preccurlyeq n^{b}$.  Conversely, assume that $n^{a}\preccurlyeq n^{b}$, meaning there is $C>0$ so that
\begin{equation*}
    n^{a} \le C\cdot (Cn)^{b}=C^{b+1}n^{b}
\end{equation*}
for any $n$ large enough. Hence $n^{a-b}$ is bounded for $n$ large enough, which forces $a-b\le 0$, i.e. $a\le b$. 

\noindent For the second claim, note that it is enough to prove that if $1<d\le c$, then $c^{n}\preccurlyeq d^{n}$. 
\end{proof}

\noindent \textit{(ii)} For any $a>0,b>1$, $n^{a} \preccurlyeq b^{n}$ and $n^{a} \nsimeq b^{n}$. 
The first part is proved exactly in the same spirit as the previous example. On the other hand, if we suppose that $n^{a}\simeq b^{n}$, then we find $C>0$ such that 
\begin{equation*}
    b^{n} \le C\cdot (Cn)^{a}=C^{a+1}n^{a}
\end{equation*}
for all $n$ large enough, so $\frac{b^{n}}{n^{a}}$ is bounded as $n\rightarrow\infty$, a contradiction.
\end{example}

We are now ready to prove that the asymptotic behaviour of the growth function is a geometric invariant.

\begin{proposition}\label{prop:invarianceofgrowth}
Let $G$ (resp. $H$) be a finitely generated group with a finite generating set $S$ (resp. $T$). If there exists a quasi-isometry $G\rightarrow H$, then $\beta_{(G,S)}(n) \simeq \beta_{(H,T)}(n)$. 
\end{proposition}

\begin{proof}
Let $C\ge 1$, $K\ge 0$ be the parameters of a quasi-isometry $f\colon G\rightarrow H$. Let $n\in\N$ and fix $g\in G$ with length $\le n$. Then we have 
\begin{align*}
    d_{T}(e_{H},f(g)) &\le d_{T}(e_{H}, f(e_{G}))+d_{T}(f(e_{G}), f(g)) \\
    &\le d_{T}(e_{H}, f(e_{G}))+C\cdot d_{S}(e_{G},g)+K \\
    &\le C\cdot n+K+E \\
    &\le (C+K+E)\cdot n
\end{align*}
where $E\defeq d_{T}(e_{H}, f(e_{G}))$. Thus $f(B_{d_{S}}(e_{G},n)) \subset B_{d_{T}}(e_{H}, D\cdot n)$, where $D\defeq C+K+E$. Additionally, if we let $Q\ge 0$ be such that the pre-image of any point of $H$ under $f$ has cardinality at most $Q$, then we get 
\begin{align*}
    \beta_{(G,S)}(n)=|B_{d_{S}}(e_{G},n)|&\le \big|f^{-1}\big(f(B_{d_{S}}(e_{G},n))\big)\big| \\
    &=\left|\bigcup_{y\in f(B_{d_{S}}(e_{G},n))} f^{-1}(\lbrace y\rbrace)\right| \\
    &\le Q\cdot \big|f(B_{d_{S}}(e_{G},n))\big| \\
    &\le Q\cdot |B_{d_{T}}(e_{H}, D\cdot n)| \\
    &=Q\cdot \beta_{(H,T)}(Dn) \\
    &\le \max(Q,D)\cdot \beta_{(H,T)}\big(\max(Q,D)n\big)
\end{align*}
which proves that $\beta_{(G,S)}(n) \preccurlyeq \beta_{(H,T)}(n)$. Considering a quasi-inverse $\overline{f}\colon H\rightarrow G$ of $f$, one obtains $\beta_{(H,T)}(n) \preccurlyeq \beta_{(G,S)}(n)$ by symmetry. Thus $\beta_{(G,S)}(n) \simeq \beta_{(H,T)}(n)$ as claimed. 
\end{proof}

In particular, if $S$ and $S'$ are different generating sets of the same group $G$, $(G,d_{S})$ and $(G,d_{S'})$ are biLipschitz equivalent, and thus $\beta_{(G,S)}(n)\simeq \beta_{(G,S')}(n)$. We therefore denote by $\beta_{G}$ the (equivalence class of the) growth function of $G$. 

If $\beta_{G}(n)\simeq n^{d}$ for some $d\in\N$, we say that $G$ has~\textit{polynomial growth of degree $d$}, and we say it has~\textit{exponential growth} if $\beta_{G}(n)\simeq e^{n}$. It has~\textit{subexponential growth} if $\beta_{G}(n) \preccurlyeq e^{n}$ and $\beta_{G}(n)\nsimeq e^{n}$, and it has~\textit{superpolynomial growth} if, for any $d\in\N$, $n^{d}\preccurlyeq \beta_{G}(n)$ and $n^{d}\nsimeq \beta_{G}(n)$. Lastly, $G$ has~\textit{intermediate growth} if it has superpolynomial growth and subexponential growth. 

Equivalently, $G$ has subexponential growth if 
\begin{equation*}
    \lim\limits_{n\rightarrow\infty}\frac{\ln(\beta_{G}(n))}{n}=0
\end{equation*}
and exponential growth if this limit does not vanish.

For instance, we know from the above examples that $\Z$ has linear growth, and $\Z^2$ has polynomial growth of degree $2$. More generally, $\Z^d$ has polynomial growth of degree $d\ge 1$, and $F_{d}$ has exponential growth for any $d\ge 2$.

The growth function behaves well when taking subgroups and quotients. 

\begin{proposition}\label{prop:growthofsubgroupsandquotients}
Let $G$ be a finitely generated group. 
\begin{enumerate}[label=(\roman*)]
    \item If $H\leqslant G$ is a finitely generated subgroup of $G$, then $\beta_{H}(n)\preccurlyeq \beta_{G}(n)$. 
    \item If $H\leqslant G$ has finite index, then $\beta_{H}(n)\simeq \beta_{G}(n)$. 
    \item If $N\lhd G$, then $\beta_{G/N}(n)\preccurlyeq \beta_{G}(n)$. In addition, if $N$ is finite, then $\beta_{G/N}(n)\simeq \beta_{G}(n)$.
    \item If $H$ is weakly commensurable to $G$, then $\beta_{G}(n)\simeq\beta_{H}(n)$.
\end{enumerate}
\end{proposition}

In particular, it follows from this proposition that any group containing a non-abelian free group has exponential growth. 

\begin{proof}
\textit{(i)} Let $S$ be a finite symmetric generating set for $G$ and $T$ be a finite symmetric generating set for $H$. Then $U\defeq S\cup T$ is a finite generating set for $G$, and as $T\subset U$, we have $B_{\mathrm{d}_{T}}(e_{H},n)\subset B_{\mathrm{d}_{U}}(e_{G},n)$ for any $n\ge 1$, whence $\beta_{(H,T)}(n) \le \beta_{(G,U)}(n)$ for any $n\ge 1$. This proves that $\beta_{(H,T)}(n)\preccurlyeq\beta_{(G,U)}(n)$.

\noindent \textit{(ii)} If $H\lhd G$ has moreover finite index, it is quasi-isometric to $G$ (Corollary~\ref{cor:QIsubgroupsoffiniteindex}) so the claim follows from Proposition~\ref{prop:invarianceofgrowth}.

\noindent \textit{(iii)} If $S$ generates $G$, then $\pi(S)$ generates $G/N$ and $B_{d_{\pi(S)}}(e_{G/N},n)=\pi(B_{d_{S}}(e_{G},n))$ for any $n\ge 1$, where $\pi\colon G\rightarrow G/N$ is the natural projection. Thus
\begin{equation*}
    \beta_{G/N}(n)=|B_{d_{\pi(S)}}(e_{G/N},n)|=\big|\pi(B_{d_{S}}(e_{G},n))\big| \le |B_{d_{S}}(e_{G},n)| = \beta_{G}(n)
\end{equation*}
for any $n\ge 1$, whence the first claim. Additionally, if $N$ is finite, then $\pi$ is $|N|$-to-one, and one has then 
\begin{align*}
    \beta_{G}(n)=|B_{d_{S}}(e_{G},n)| &\le \big|\pi^{-1}\big(\pi(B_{d_{S}}(e_{G},n))\big)\big| \\
    &\le |N|\cdot \big|\pi(B_{d_{S}}(e_{G},n))\big| \\
    &\le |N|\cdot \beta_{G/N}(|N|\cdot n)
\end{align*}
for any $n\ge 1$, so $\beta_{G}(n)\preccurlyeq \beta_{G/N}(n)$. 

\noindent \textit{(iv)} If $H$ is weakly commensurable to $G$, it is finitely generated and quasi-isometric to $G$ by Corollary~\ref{cor:weakcommensurabilityQI}, and we get the result applying Proposition~\ref{prop:invarianceofgrowth}.
\end{proof}

We also deduce from this proposition the next upper bound for the growth of a finitely generated group.

\begin{corollary}
Let $G$ be a finitely generated group. Then $G$ grows at most exponentially fast.
\end{corollary}

\begin{proof}
Since $G$ is finitely generated, it is a quotient of a free group of finite rank. Such a group grows exponentially fast, and the conclusion follows from Proposition~\ref{prop:growthofsubgroupsandquotients}\textit{(iii)}. 
\end{proof}

For instance, $\text{SL}_{2}(\Z)$, $\text{PSL}_{2}(\Z)$, $\Z_{2}*\Z_{2}*\Z_{2}*\Z_{2}$ all have exponential growth, since all are quasi-isometric to $F_{2}$. On the other hand, $D_{\infty}$ has polynomial growth of degree $1$ since it is quasi-isometric to $\Z$. 

Let us also mention the following. 

\begin{proposition}
Let $G$ (resp. $H$) be a finitely generated group with a finite generating set $S$ (resp. $T$). Then $G\times H$ is finitely generated and $\beta_{G\times H}(n)\simeq \beta_{(G,S)}(n)\beta_{(H,T)}(n)$.
\end{proposition}

\begin{proof}
The fact that $G\times H$ is finitely generated is a consequence of Proposition~\ref{prop:finitegenerationextension}, and the proof of the latter shows that 
\begin{equation*}
    U\defeq \lbrace (s,e_{H}) : s\in S\rbrace\cup\lbrace (e_{G},t) : t\in T\rbrace
\end{equation*}
is a finite generating set for $G\times H$. As the growth type of $G\times H$ is independent of the choice of the generating set (see the remark right after Proposition~\ref{prop:invarianceofgrowth}), we now show that
\begin{equation*}
    \beta_{(G\times H, U)}(n)\simeq \beta_{(G,S)}(n)\beta_{(H,T)}(n).
\end{equation*}

\noindent First of all, if $g\in B_{\mathrm{d}_{S}}(e_{G},n)$ and if $h\in B_{\mathrm{d}_{T}}(e_{H},n)$, then $(g,h)\in B_{\mathrm{d}_{U}}(e_{G\times H},2n)$, and thus 
\begin{align*}
    \beta_{(G,S)}(n)\beta_{(H,T)}(n)&=|B_{\mathrm{d}_{S}}(e_{G},n)||B_{\mathrm{d}_{T}}(e_{H},n)| \\
    &\le |B_{\mathrm{d}_{U}}(e_{G\times H},2n)|\\
    &=\beta_{(G\times H, U)}(2n) \\
    &\le 2\beta_{(G\times H, U)}(2n) 
\end{align*}
for any $n\ge 1$. Hence $\beta_{(G,S)}(n)\beta_{(H,T)}(n) \preccurlyeq \beta_{(G\times H, U)}(n)$. Conversely, if $n\ge 1$ and $(g,h)\in B_{\mathrm{d}_{U}}(e_{G\times H}, n)$, there are $p,r\in\mathbb{N}$ with $p+r\le n$ and group elements $x_{1},\dots,x_{p}\in S\cup S^{-1}$, $y_{1},\dots,y_{r}\in T\cup T^{-1}$ such that 
\begin{equation*}
    (g,h)=(x_{1},e_{H})\dots(x_{p},e_{H})(e_{G}, y_{1})\dots(e_{G},y_{r})=(x_{1}\dots x_{p}, y_{1}\dots y_{r}).
\end{equation*}
Hence $B_{\mathrm{d}_{U}}(e_{G\times H}, n) \subset B_{\mathrm{d}_{S}}(e_{G}, n)\times B_{\mathrm{d}_{T}}(e_{H},n)$, and it follows that 
\begin{equation*}
    \beta_{(G\times H, U)}(n) \le \beta_{(G,S)}(n)\beta_{(H,T)}(n)
\end{equation*}
for any $n\ge 1$. We conclude that $\beta_{(G\times H, U)}(n)\preccurlyeq \beta_{(G,S)}(n)\beta_{(H,T)}(n)$, and thus that 
\begin{equation*}
    \beta_{(G\times H, U)}(n)\simeq \beta_{(G,S)}(n)\beta_{(H,T)}(n).
\end{equation*}
\end{proof}

This implies the following.

\begin{corollary}\label{cor3.29}
Any finitely generated abelian group has polynomial growth. 
\end{corollary}

\begin{proof}
Such a group $G$ splits as a product $\Z^{d}\times F$, where $d\in\N$ and $F$ is a finite group (see e.g.~\cite[Corollary~1.30]{DK18}). As the growth function of $F$ is constant for $n$ large enough, we conclude that $\beta_{G}(n)\simeq \beta_{\Z^d}(n)\simeq n^{d}$, as claimed. 
\end{proof}

These results already allow us to distinguish euclidean spaces up to quasi-isometry.

\begin{corollary}
For any $d\neq d'\in\N$, $\Z^{d}$ is not quasi-isometric to $\Z^{d'}$. As a consequence, $\R^{d}$ is not quasi-isometric to $\R^{d'}$.
\end{corollary}

\begin{proof}
If $d\neq d'$ and $\Z^{d}\sim_{Q.I.}\Z^{d'}$, then $n^{d}\simeq n^{d'}$, so $d=d'$ by Example~\ref{ex:examplesofasymptoticequivalence}\textit{(i)}. This contradiction shows that $\Z^d \nsim_{Q.I.} \Z^{d'}$. In particular, as $\R^d$ (resp. $\R^{d'}$) is quasi-isometric to $\Z^{d}$ (resp. $\Z^{d'}$), we deduce also $\R^d \nsim_{Q.I.} \R^{d'}$.
\end{proof}

Amazingly, it turns out that the behaviour of the growth function of a finitely generated group encodes a lot of information about its algebraic structure. A first major result highlighting this connection is Milnor's theorem. 

\begin{theorem}[{\cite{Mil68}}]\label{thm:Milnor}
A finitely generated solvable group of subexponential growth is polycyclic. 
\end{theorem}

Recall that a group $G$ is~\textit{solvable} (resp.~\textit{polycyclic}) if one can find a sequence of subgroups
\begin{equation*}
    \lbrace 1\rbrace=H_{s} \leqslant H_{s-1}\leqslant \dots\leqslant H_{1}\leqslant H_{0}=G
\end{equation*}
such that, for any $0\le i\le s-1$, $H_{i+1}\lhd H_{i}$ and the quotient $H_{i}/H_{i+1}$ is abelian (resp. cyclic). 

As an application, a wreath product $A\wr B$ has exponential growth as soon as $A$ is a non-trivial finitely generated solvable group, and $B$ is an infinite finitely generated solvable group.

\begin{proof}
Indeed, if $A$ and $B$ are both solvable, then $A\wr B$ is solvable, and it is finitely generated by Proposition~\ref{prop:wreathproductspreservefinitegeneration}. Thus, if $A\wr B$ has subexponential growth, it follows from Milnor's theorem that it is polycyclic, and in particular, all its subgroups are finitely generated~\cite[Proposition~5.5]{CecAdd21}. Thus $\bigoplus_{B}A$ is finitely generated, a contradiction since $B$ is infinite. Hence $A\wr B$ cannot have subexponential growth. 
\end{proof}

On the other hand, Wolf's theorem improves polycyclicity to virtual nilpotency when the group has subexponential growth. 

\begin{theorem}[{\cite{Wol68}}]\label{thm:Wolf}
A finitely generated polycyclic group of subexponential growth is virtually nilpotent.
\end{theorem}

A proof is presented for instance in~\cite[Theorem~7.37]{CecAdd21}.

These two theorems combined therefore show intermediate growth is impossible among solvable groups.

\begin{corollary}
A finitely generated solvable group either has exponential or polynomial growth. In the latter case, it is virtually nilpotent. 
\end{corollary}

In this statement, the first sentence is an immediate consequence of Theorems~\ref{thm:Milnor} and~\ref{thm:Wolf}. The second sentence is not at all obvious, and is in fact a celebrated result of Gromov.

\begin{theorem}[{\cite{Gro81}}]\label{thm:Gromov}
A finitely generated group has polynomial growth if and only if it is virtually nilpotent. 
\end{theorem}

On the other side, the existence of intermediate growth groups has been proved by Grigorchuk in the 80's. 

\begin{theorem}[{\cite{Gri84}}]\label{thm:Grigorchuk}
There exist finitely generated groups of intermediate growth. 
\end{theorem}

\subsection{Amenability}\label{subsection:amenability}

Let us now present another geometric invariant of groups: amenability. 

\begin{definition}
Let $G$ be a finitely generated group with a finite generating set $S$. We say that $G$ is~\textit{amenable} if there exists a sequence $(F_{n})_{n\in\N}$ of finite subsets of $G$ such that 
\begin{equation*}
    \lim\limits_{n\rightarrow\infty}\frac{|\partial_{G}F_{n}|}{|F_{n}|}=0
\end{equation*}
where $\partial_{G}F_{n} \defeq \lbrace g\in G\setminus F_{n} : \exists f\in F_{n}, \exists s\in S, g=fs\rbrace$ is the~\textit{boundary} of $F_{n}$ in $G$.
\end{definition}

When it exists, such a sequence is usually called a~\textit{F\o lner sequence}.

Amenability is known to admit plenty of other equivalent definitions, that relates it to many other fields of study, such as for instance random walks and probability theory, fixed point theorems and functional analysis, or paradoxical decompositions and the Banach-Tarski paradox to name a few. 

\begin{theorem}\label{thm:equivalentdefofamenability}
Let $G$ be a finitely generated group. The following are equivalent. 
\begin{enumerate}[label=(\roman*)]
    \item $G$ is amenable. 
    \item For all $\varepsilon>0$, there exists a finite set $A\subset G$ such that $|\partial_{G}A|\le \varepsilon\cdot |A|$. 
    \item $G$ has a left-invariant mean, i.e. a function $\mu\colon \mathcal{P}(G)\longrightarrow \left[0,1\right]$ such that $\mu(G)=1$, $\mu(A\cup B)=\mu(A)+\mu(B)$ where $A,B\subset G$ are disjoint, and $\mu(gA)=\mu(A)$ for any $A\subset G$ and $g\in G$.
    \item Any action of $G$ on a compact set $K\subset V$ in a locally convex topological vector space $V$ has a fixed point. 
\end{enumerate}
\end{theorem}

\begin{example}\label{ex:examplesofAgroups}
\textit{(i)} Finite groups are amenable. For instance, if $F$ is finite, it suffices to consider the normalized counting measure $\mu(A)=\frac{|A|}{|F|}$, $A\subset F$, to get an invariant mean. 

\noindent \textit{(ii)} $\Z$ is amenable: if $F_{n}=\lbrace -n,\dots,n\rbrace$ and $S=\lbrace -1,1\rbrace$, then $\frac{|\partial_{\Z}F_{n}|}{|F_{n}|}=\frac{2}{2n+1}$ goes to $0$ as $n\rightarrow\infty$. More generally, $\Z^d$ is amenable for any $d\ge 1$, consider the sequence $F_{n}\defeq \lbrace -n,\dots, n\rbrace^{d}$ and the canonical generating set of $\Z^d$. 

\noindent \textit{(iii)} The free group $F_{2}$ is not amenable. 

\begin{proof}
Suppose for a contradiction that there exists an invariant mean $\mu$ on $F_{2}$. Write 
\begin{equation*}
    F_{2}=\lbrace e\rbrace\sqcup A_{+}\sqcup A_{-}\sqcup B_{+}\sqcup B_{-}
\end{equation*}
where $A_{+}$ (resp. $A_{-}$) consists of reduced words starting with an $a$ (resp. $a^{-1}$) and $B_{+}$ (resp. $B_{-}$) consists of reduced words starting with a $b$ (resp. $b^{-1}$). Since the second letter of an element of $A_{+}$ can be an $a$, a $b$ or a $b^{-1}$, multiplying this element by $a^{-1}$ produces an element either of $A_{+}$, $B_{+}$ or $B_{-}$. It follows that 
\begin{equation*}
    a^{-1}A_{+}=\lbrace e \rbrace\sqcup A_{+}\sqcup B_{+} \sqcup B_{-}.
\end{equation*}
and properties of $\mu$ then imply  
\begin{equation*}
    \mu(A_{+})=\mu(a^{-1}A_{+})=\mu\big(\lbrace e \rbrace\sqcup A_{+}\sqcup B_{+} \sqcup B_{-}\big)=\mu(\lbrace e \rbrace)+\mu(A_{+})+\mu(B_{+})+\mu(B_{-}).
\end{equation*}
Thus $\mu(\lbrace e \rbrace)+\mu(B_{+})+\mu(B_{-})=0$. Since $\mu$ is positive, this forces $\mu(\lbrace e\rbrace)=\mu(B_{+})=\mu(B_{-})=0$. Likewise, we get $\mu(A_{+})=\mu(A_{-})=0$. We conclude that
\begin{align*}
    1=\mu(F_{2})&=\mu\big(\lbrace e\rbrace\sqcup A_{+}\sqcup A_{-}\sqcup B_{+}\sqcup B_{-}\big)\\
    &=\mu(\lbrace e\rbrace)+\mu( A_{+})+\mu(A_{-})+\mu(B_{+})+\mu(B_{-}) \\
    &=0
\end{align*}
a contradiction. Therefore such a $\mu$ cannot exist.
\end{proof}
\end{example}

One of the questions that arise from these examples is whether amenability of a group depends on the choice of a finite generating set for that group. The following theorem, asserting as promised that amenability is a geometric invariant, answers the question.

\begin{theorem}[{\cite[Theorem~3.1.5]{NG23}}]
Let $G$ and $H$ be finitely generated groups. If there exists a quasi-isometry $f\colon G\rightarrow H$ and if $H$ is amenable, then $G$ is amenable. 
\end{theorem}

Amenability also enjoys several stability properties under classical group operations, that makes it very robust. 

\begin{theorem}\label{thm:stabilitypropertiesamenability}
Let $G$ be a finitely generated group. 
\begin{enumerate}[label=(\roman*)]
    \item \textit{(Subgroups)} If $G$ is amenable, and $H\leqslant G$, then $H$ is amenable. 
    \item \textit{(Quotients)} If $G$ is amenable, and $N\lhd G$, then $G/N$ is amenable. 
    \item \textit{(Extensions)} If $N\lhd G$, and $N, G/N$ are amenable, then $G$ is amenable. 
\end{enumerate}
\end{theorem}

From the first statement, we can in particular deduce that: 

\begin{corollary}
Any group containing a subgroup isomorphic to $F_{2}$ is not amenable. In particular, $F_{d}$ is not amenable for all $d\ge 2$. 
\end{corollary}

Thus, for instance, $\text{SL}_{n}(\Z)$ is not amenable, as well as $\text{PSL}_{n}(\Z)$, for all $n\ge 2$. 

On the other hand, stability under extensions has the following corollaries, that provide additional examples of amenable groups.

\begin{corollary}\label{cor:fgabeliangroupsareamenable}
Any finitely generated abelian group is amenable.
\end{corollary}

\begin{proof}
If $G$ is finitely generated and abelian, it splits as a direct product $G\cong \Z^d\times F$ where $d\ge 1$ and $F$ is finite. Since both $F$ and $\Z^d$ are amenable, $G$ is amenable by Theorem~\ref{thm:stabilitypropertiesamenability}\textit{(iii)}. 
\end{proof}

\begin{corollary}
Let $A$ and $B$ be finitely generated groups. Then $A\wr B$ is amenable if and only if $A$ and $B$ are amenable. 
\end{corollary}

\begin{proof}
If $A\wr B$ is amenable, amenability of $A$ and $B$ follows from Theorem~\ref{thm:stabilitypropertiesamenability}\textit{(i)} and Theorem~\ref{thm:stabilitypropertiesamenability}\textit{(ii)} respectively. Conversely, if $A$ and $B$ are amenable, then so is $\bigoplus_{B}A$, and thus $A\wr B$ is amenable by Theorem~\ref{thm:stabilitypropertiesamenability}\textit{(iii)}. 
\end{proof}

Let us also emphasize on a fourth stability result: recall that given a group $G$ and a collection $\mathcal{F}$ of subgroups of $G$, $\mathcal{F}$ is called \textit{directed} if for any $H,H'\in \mathcal{F}$, there is $H''\in \mathcal{F}$ such that $H,H' \leqslant H''$. We then say that $G$ is \textit{directed} if $\mathcal{F}$ is directed and if 
\begin{equation*}
    G=\bigcup_{H\in \mathcal{F}}H.
\end{equation*}

\begin{proposition}\label{prop:stabilityunderdirectedunions}
If $G$ is the directed union of $\mathcal{F}$ and if any $H\in\mathcal{F}$ is amenable, then $G$ is amenable. 
\end{proposition}

\begin{corollary}\label{cor:amenabilitybyfgsubgroups}
A group is amenable if and only if all its finitely generated subgroups are amenable. 
\end{corollary}

\begin{proof}
If $G$ is amenable, any of its subgroups is amenable by Theorem~\ref{thm:stabilitypropertiesamenability}\textit{(i)}. Conversely, suppose any finitely generated subgroup of $G$ is amenable, and consider 
\begin{equation*}
    \mathcal{F} \defeq \lbrace H\leqslant G : H \;\text{is finitely generated}\rbrace.
\end{equation*}
The collection $\mathcal{F}$ is directed, as if $H=\la S\ra$, $H'=\la S'\ra$ are both finitely generated, the subgroup $H''=\la S\cup S'\ra$ is in $\mathcal{F}$ and contains both $H$ and $H'$ as subgroups. By assumption, any $H\in\mathcal{F}$ is amenable, and $G$ is the directed union of $\mathcal{F}$, whence $G$ is amenable by Proposition~\ref{prop:stabilityunderdirectedunions}.
\end{proof}

We can then strengthen Corollary~\ref{cor:fgabeliangroupsareamenable}. 

\begin{corollary}\label{cor:abeliangroupsareamenable}
Any abelian group is amenable. 
\end{corollary}

\begin{proof}
Let $G$ be an abelian group. By Corollary~\ref{cor:amenabilitybyfgsubgroups}, it is enough to prove that any finitely generated subgroup of $G$ is amenable. Such a subgroup is then finitely generated and abelian, whence amenable by Corollary~\ref{cor:fgabeliangroupsareamenable}. 
\end{proof}

\begin{corollary}
Any solvable group is amenable. 
\end{corollary}

\begin{proof}
Such a group is obtained from the trivial group by doing finitely many extensions by abelian groups. As these are amenable by Corollary~\ref{cor:abeliangroupsareamenable}, and as amenability is preserved by extensions (cf. Theorem~\ref{thm:stabilitypropertiesamenability}\textit{(iii)}), the claim follows. 
\end{proof}

We conclude our list of examples by a straightforward implication between volume growth and amenability. 

\begin{proposition}
Any finitely generated subexponential growth group is amenable. 
\end{proposition}

\begin{proof}
Assume that a finitely generated group $G=\langle S\rangle$ is not amenable. Then there exists $\varepsilon>0$ such that $|\partial_{G}A|>\varepsilon\cdot |A|$ for any finite subset $A\subset G$. For any $n\ge 1$, the ball $B_{d_{S}}(e_{G},n)$ can be decomposed as the disjoint union 
\begin{equation*}
    B_{d_{S}}(e_{G},n) = B_{d_{S}}(e_{G},n-1)\sqcup S(e_{G}, n)
\end{equation*}
where $S(e_{G}, n)$ is the~\textit{sphere centered at $e_{G}$ of radius $n$}, that is the set of elements of the group having word length exactly $n$. Thus 
\begin{align*}
    \beta_{G}(n) &= |B_{d_{S}}(e_{G},n)| = |B_{d_{S}}(e_{G},n-1)|+|S(e_{G}, n)|\\
    &=|B_{d_{S}}(e_{G},n-1)|\cdot \left(1+\frac{|S(e_{G}, n)|}{|B_{d_{S}}(e_{G},n-1)|}\right)
\end{align*}
and iterating, one obtains 
\begin{align*}
    \beta_{G}(n)&=|B_{d_{S}}(e_{G},n-1)|\left(1+\frac{|S(e_{G}, n)|}{|B_{d_{S}}(e_{G},n-1)|}\right) \\
    &=|B_{d_{S}}(e_{G},n-2)|\left(1+\frac{|S(e_{G}, n-1)|}{|B_{d_{S}}(e_{G},n-2)|}\right)\left(1+\frac{|S(e_{G}, n)|}{|B_{d_{S}}(e_{G},n-1)|}\right)\\
    &= \dots \\
    &=\prod_{k=1}^{n}\left(1+\frac{|S(e_{G}, k)|}{|B_{d_{S}}(e_{G}, k-1)|}\right).
\end{align*}
As $|S(e_{G}, k)|=|\partial_{G}B_{d_{S}}(e_{G}, k-1)| > \varepsilon\cdot |B_{d_{S}}(e_{G}, k-1)|$, one gets $\beta_{G}(n) > (1+\varepsilon)^{n}$, and thus 
\begin{equation*}
    \beta_{G}(n) \succcurlyeq (1+\varepsilon)^{n}. 
\end{equation*}
Hence $G$ has exponential growth as claimed.
\end{proof}

For instance, the Grigorchuk group from Theorem~\ref{thm:Grigorchuk} has subexponential growth, and thus is amenable. 

\subsection{Isoperimetric profiles}\label{sec:isoprof}

We now focus on another geometric invariant of finitely generated groups, which is in general finer than the volume growth, since it also takes into account how~\textit{boundaries} of finite subsets grow. 

For a finitely generated group $G$ and a finite generating set $S_{G}$, its $\ell^{p}-$\textit{isoperimetric profile}, for $p\ge 1$, is the function $\profp{G}\colon \N\rightarrow \R_{+}$ given by
\begin{equation*}
    \profp{G}(n)\defeq\sup_{\substack{f\colon G\to\R_+\\|\supp{f}|\leq n}}{\frac{\|f\|_{p}}{\|\nabla f\|_{p}}}
\end{equation*}
where the support of $f\colon G\rightarrow\R_+$ is $\supp{f}\defeq\lbrace g\in G : f(g)\not=0\rbrace$ and the $\ell^p-$norm of its gradient is defined by
\begin{equation*}
    \|\nabla f\|_{p}^{p}\defeq\sum_{g\in G,\; s\in S_{G}}{|f(g)-f(gs)|^{p}}. 
\end{equation*}
For $p=1$, the $\ell^{1}-$isoperimetric profile is merely called \textit{isoperimetric profile} and one has
\begin{equation*}
    \prof{G}(n) \simeq \sup_{|A|\le n}\frac{|A|}{|\partial_{G} A|}
\end{equation*}
with $\partial_{G}A = \lbrace g\in G\setminus A : \exists s\in S_{G}, \exists h\in A, g=hs\rbrace$ the boundary of $A$ in $G$.

Recall also that the isoperimetric profile of a group $G$ is the generalized inverse of its~\textit{F\o lner function} $\text{F\o l}_{G}\colon \N\rightarrow\R_{+}$, defined as
 \begin{equation*}
     \text{F\o l}_{G}(n) \defeq \inf\left\lbrace |A| : \frac{|\partial_{G}A|}{|A|} \le \frac{1}{n} \right\rbrace.
 \end{equation*}
 
In fact, the $\ell^{p}-$isoperimetric profile of $G$ is the generalized inverse of the~\textit{$\ell^{p}-$F\o lner function} of $G$, that we define in Chapter~\ref{chap:chapter5}. 
 
The next theorem states that, as for the volume growth, the asymptotic behaviour of the $\ell^{p}-$isoperimetric profile is invariant under quasi-isometries. In particular, when studying asymptotics of this function, the choice of the finite generating set is not important.
 
\begin{theorem}[{\cite[Theorem~1]{Tes08}}]\label{thm:profilQIinvariant}
Let $G$ and $H$ be finitely generated groups. Let $p\ge 1$. If there exists a quasi-isometry $f\colon G\rightarrow H$, then $\profp{G}(n)\simeq \profp{H}(n)$.
\end{theorem}

Erschler also noticed that the $\ell^{p}-$profile is monotonically decreasing when taking subgroups:

\begin{lemma}[{\cite[Lemma~4]{Ers03}}]\label{lem:monotonieprofilErschler}
Let $G$ be a finitely generated amenable group, and let $H\leqslant G$ be a finitely generated subgroup. Then, for any $p\ge 1$, $\profp{G}(n) \preccurlyeq \profp{H}(n)$. 
\end{lemma}

We shall also notice that, if $G$ is not amenable, then there is $\varepsilon>0$ such that $|\partial_{G}A| > \varepsilon \cdot |A|$ (Theorem~\ref{thm:equivalentdefofamenability}) for any finite subset $A\subset G$, so the $\ell^{1}-$isoperimetric profile of $G$ is bounded from above. The same holds for the $\ell^{p}-$profile of $G$ for $p\ge 1$, since the profile is monotonous decreasing into the variable $p$ (see Lemma~\ref{lem:MonotonuousProfile} below). Thus it is irrelevant to study the behaviour of these functions for non-amenable groups. 

On the other hand, the asymptotics of the profile are now known for several classical classes of amenable groups, among which:
\begin{itemize}
    \item $\profp{G}(n) \simeq n^{\frac{1}{d}}$ if $G$ has polynomial growth of degree $d\ge 1$;
    \item $\profp{G}(n)\simeq \ln(n)$ for $G=\text{BS}(1,k)$ for $k\ge 2$, or $G=F\wr\Z$, where $F$ is a non-trivial finite group;
    \item $\profp{G}(n) \simeq \ln(n)$ for any polycyclic group $G$ with exponential growth~\cite{Pit95, Pit00}, or more generally any exponential growth group within the class GES of Tessera~\cite[Corollary 5]{Tes13}.
\end{itemize}

In fact, the main result of~\cite{Ers03} provides an explicit formula to compute the F\o lner function of a wreath product in terms of F\o lner functions of both factors. 

\begin{theorem}[{\cite[Theorem~1]{Ers03}}]\label{thm:formulaforprofilewreathproducts}
Let $G$ and $H$ be finitely generated amenable groups. Suppose that the F\o lner function of $H$ satisfies:
\begin{equation*}
    \forall C>0,\;\exists K>0,\; C\cdot\text{F\o l}_{H}(n) \le \text{F\o l}_{H}(Kn) \quad\text{for all $n$ large enough}.
\end{equation*}
Then one has 
\begin{equation*}
    \text{F\o l}_{G\wr H}(n) \simeq \text{F\o l}_{G}(n)^{\text{F\o l}_{H}(n)}.
\end{equation*}
\end{theorem}

The assumption of Theorem~\ref{thm:formulaforprofilewreathproducts} can also be formulated in terms of the $\ell^{1}-$isoperimetric profile: $\prof{H}(Cn)=O(\prof{H}(n))$ for any $C>0$. We make further comments on it in Chapter~\ref{chap:chapter5}, since it also appears over there.

\section{Orbit equivalence}\label{sec:OE}

This part is devoted to introduce a measured framework for finitely generated groups, which is crucial for Chapter~\ref{chap:chapter6}. As we will explain below, the tools we develop now can be used to understand more precisely the large-scale geometry of finitely generated groups. 

The starting point is the notion of~\textit{orbit equivalence}: two groups $G$ and $H$ are~\textit{orbit equivalent} if there exist a standard probability space $(X,\mu)$ and two free probability measure-preserving (p.m.p.) actions $G,H\act (X,\mu)$ that share the same orbits, i.e. $G\cdot x=H\cdot x$ for almost every $x\in X$. In this case, the data of the space $(X,\mu)$ and these two actions is called an~\textit{orbit equivalence coupling} between $G$ and $H$.

The following observation is useful in practice for simplifying computations.

\begin{lemma}\label{lem:EqOrbitsGenerating}
Let $G,H\act (X,\mu)$ be two free p.m.p. actions of countable groups and denote by $S_{G}$ (resp. $S_{H}$) a generating set of $G$ (resp. $H$). Assume that 
\begin{equation*}
    S_{G}\cdot x\subset H\cdot x\;\text{and}\;S_{H}\cdot x\subset G\cdot x
\end{equation*}
for almost every $x\in X$. Then $G$ and $H$ have the same orbits.
\end{lemma}

\begin{proof}
For every $x\in X$ and every $g\in S_{G}$, let us denote by $c_{G,H}(g,x)$ the element of $H$ satisfying
\begin{equation*}
    g\cdot x=c_{G,H}(g,x)\cdot x.
\end{equation*}
We then have
\begin{equation*}
    x=g\cdot(g^{-1}\cdot x)=c_{G,H}(g,g^{-1}\cdot x)\cdot (g^{-1}\cdot x)
\end{equation*}
which implies that
\begin{equation*}
    g^{-1}\cdot x=c_{G,H}(g,g^{-1}\cdot x)^{-1}\cdot x
\end{equation*}
and we define $c_{G,H}(g^{-1},x)\defeq c_{G,H}(g,g^{-1}\cdot x)^{-1}\in H$. Finally, given $g\in G$, we write it as $g=g_{1}\ldots g_{n}$ with $g_{1},\dots,g_{n}\in S_{G}\cup S_G^{-1}$, and it follows that 
\begin{align*}
    g\cdot x&=g_{1}\cdot (g_{2}\dots g_{n} \cdot x)\\
    &=c_{G,H}(g_{1},g_{2}\dots g_{n} \cdot x)g_{2}\dots g_{n} \cdot x\\
    &=\dots \\
    &=c_{G,H}(g_{1},g_{2}\dots g_{n}\cdot x)c_{G,H}(g_{2},g_{3}\dots g_{n} \cdot x)\dots c_{G,H}(g_{n},x)\cdot x
\end{align*}
whence $g\cdot x\in H\cdot x$. This holds for every $g\in G$, so we get $G\cdot x\subset H\cdot x$. We similarly prove the reverse inclusion to conclude the proof. 
\end{proof}

As proved by Ornstein and Weiss in 1980, the notion of orbit equivalence is too flexible to distinguish amenable groups.

\begin{theorem}[{\cite{OW80}}]\label{thm:OW}
Any two infinite amenable groups are orbit equivalent. 
\end{theorem}

Note that an orbit equivalence coupling comes with two maps $c_{G,H}\colon G\times X\rightarrow H$, $c_{H,G}\colon H\times X\rightarrow G$, called the~\textit{cocycles} of the coupling and defined by the equations
\begin{equation*}
    g\cdot x = c_{G,H}(g,x)\cdot x, \; h\cdot x = c_{H,G}(h,x)\cdot x
\end{equation*}
for all $g\in G$, $h\in H$ and $\mu-$almost every $x\in X$. These maps are well-defined by freeness of the initial actions, and satisfy a~\textit{cocycle identity}:
\begin{equation*}
    c_{G,H}(gg',x)=c_{G,H}(g, g'\cdot x)c_{G,H}(g',x)
\end{equation*}
and similarly for $c_{H,G}$. Moreover, for $\mu-$almost every $x\in X$, the maps $c_{G,H}(\cdot,x)\colon G\rightarrow H$ and $c_{H,G}(\cdot,x)\colon H\rightarrow G$ are bijections and inverses of each other. 

Intuitively, cocycles can be used to detect geometric differences between $G$ and $H$, because they measure the distortion we make when passing from a $G-$orbit to an $H-$orbit, and such orbits are copies of Cayley graphs of $G$ and $H$. 

\begin{figure}[H]
  \centering
  \includegraphics[
    width=0.75\linewidth
  ]{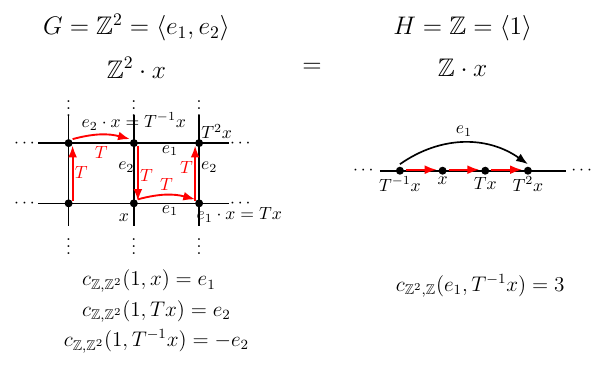}
  \caption{An orbit equivalence between $\Z$ and $\Z^2$.}
  \label{fig:cayley}
\end{figure}

In order to make this formal, we~\textit{quantify} the cocycles:

\begin{definition}\label{def:integrabilityofcocycles}
Let $G$ and $H$ be finitely generated groups with finite generating sets $S_{G}$, $S_{H}$. Let $(X,\mu)$ be an orbit equivalence coupling between $G$ and $H$. Given a map $\varphi\colon \R_{+}\rightarrow \R_{+}$, we say that $c_{G,H}$ is~\textit{$\varphi-$integrable} if, for every $g\in G$, there is a constant $C_{g}>0$ such that 
\begin{equation*}
    \int_{X}\varphi\big(C_{g}|c_{G,H}(g,x)|_{S_{H}}\big)\;\mathrm{d}\mu(x) < +\infty.
\end{equation*}
Then, given two maps $\varphi,\psi\colon\R_{+}\rightarrow \R_{+}$, we say that $(X,\mu)$ is a~\textit{$(\varphi,\psi)-$integrable orbit equivalence coupling} if $c_{G,H}$ is $\varphi-$integrable and $c_{H,G}$ is $\psi-$integrable. 
\end{definition}

In this case, we also say that $G$ and $H$ are~\textit{$(\varphi,\psi)-$integrably orbit equivalent}. If $\psi=\varphi$, we simply say that $G$ and $H$ are~\textit{$\varphi-$integrably orbit equivalent}. Moreover, we say that $c_{G,H}$ is~\textit{$\ld^p-$integrable} if $\varphi(x)=x^p$, for $p>0$, and $\ld^0$ means that there is no requirement on the cocycle. We say that $c_{G,H}$ is $\ld^{\infty}$ if, for any $g\in G$, the map $|c_{G,H}(g,\cdot)|_{S_{H}}\colon X\rightarrow \N$ is essentially bounded. Lastly, $c_{G,H}$ is $\ld^{<\infty}$ if it is $\ld^p$ for any $p<\infty$. 

\begin{remark}\label{rem:checkongenerators}
The use of the constants $C_{g}$ in the definition of $\varphi-$integrability of $c_{G,H}$ is necessary because we want the following properties:
\begin{itemize}
    \item The notion of $\varphi-$integrability does not depend on the particular choices of finite generating sets $S_{G}$, $S_{H}$. Indeed, by Lemma~\ref{lem:changeofgensetsBILIP}, two different generating sets and the associated word metrics on the same group yield biLipschitz equivalent metric spaces;
    \item If $\varphi\simeq \psi$, then $c_{G,H}$ is $\varphi-$integrable if and only if $c_{G,H}$ is $\psi-$integrable;
    \item To prove that $c_{G,H}$ is $\varphi-$integrable, it suffices to check the finiteness of 
    \begin{equation*}
        \int_{X}\varphi\big(C_{g}|c_{G,H}(g,x)|_{S_{H}}\big)\;\mathrm{d}\mu(x)
    \end{equation*}
    for $g$ in a finite generating set of $G$ (see~\textit{e.g.}~\cite[Proposition~2.22]{DKLMT22}). The same remark holds for $\ld^{\infty}$. 
\end{itemize}
\end{remark}

This notion of orbit equivalence is indeed more rigid, and starts to detect the large-scale geometry of the groups we are coupling. One of the first results in this direction is due to Bowen, proved in the appendix of Austin's paper~\cite{Aus16}, and states that volume growth is an invariant of $\ld^{1}-$orbit equivalence:

\begin{theorem}[{\cite[Theorem~B.2]{Aus16}}]
Let $G$ and $H$ be finitely generated groups. If $G$ and $H$ are $\ld^{1}-$orbit equivalent, then $\beta_{G}(n)\simeq \beta_{H}(n)$.
\end{theorem}

This statement has been generalised in~\cite[Theorem~3.1]{DKLMT22}, where the authors proved more generally that a $(\varphi, \ld^0)-$integrable orbit equivalence coupling from $G$ to $H$, where $\varphi$ is increasing and subadditive such that $\varphi(0)=0$, implies that 
\begin{equation*}
    \varphi\circ \beta_{H}(n) \preccurlyeq \beta_{G}(n).
\end{equation*}

In the same paper, they in fact proved the following finer asymptotic inequality, with the isoperimetric profile:

\begin{theorem}[{\cite[Theorem~1.1]{DKLMT22}}]\label{thm:ObstructionDKLMT}
Let $G$ and $H$ be finitely generated amenable groups. Let $\varphi\colon \R_{+}\rightarrow\R_{+}$. Assume that there is a $(\varphi, \ld^0)-$integrable orbit equivalence coupling from $G$ to $H$. 
\begin{enumerate}[label=(\roman*)]
    \item If $\varphi$ and $t\longmapsto \frac{t}{\varphi(t)}$ are increasing, then $\varphi\circ\prof{H}(n) \preccurlyeq \prof{G}(n)$.
    \item If $\varphi(x)=x^p$ for some $p\ge 1$, then $\profp{H}(n)\preccurlyeq \profp{G}(n)$.
    \end{enumerate}
\end{theorem}

Before starting the proof, let us notice here that the assumptions on $\varphi$ in point~\textit{(i)} imply 
\begin{equation}\label{eq:eq1.3}
    \varphi(Ct)=O(\varphi(t))
\end{equation}
for every $C>0$. Indeed:
\begin{itemize}
    \item If $0<C\le 1$, then $\varphi(Ct)\le \varphi(t)$ by increasingness of $\varphi$; 
    \item If $C>1$, then $\frac{t}{\varphi(t)}\le \frac{Ct}{\varphi(Ct)}$ by increasingness of $t\longmapsto\frac{t}{\varphi(t)}$, so $\varphi(Ct) \le C\varphi(t)$.
\end{itemize}
Hence one also has $\varphi(t)=O(\varphi(Ct))$ for every $C>0$, since 
\begin{equation*}
    \varphi(t)=\varphi\left(\frac{1}{C}\cdot Ct\right)=O(\varphi(Ct))
\end{equation*}
using (\ref{eq:eq1.3}) with the constant $\frac{1}{C}$. In conclusion, under assumptions of point~\textit{(i)}, we can forget about the constants $C_{g}$ in the definition of $\varphi-$integrability.

\begin{proof} Let $S_{G}$ and $S_{H}$ be finite generating sets of $G$ and $H$ respectively. 

\noindent \textit{(i)} Assume that $(X,\mu)$ is an orbit equivalence coupling from $G$ to $H$, where $\varphi$ and $t\longmapsto\frac{t}{\varphi(t)}$ are increasing. Consider the cocycles $c_{G,H}\colon G\times X\rightarrow H$, $c_{H,G}\colon H\times X\rightarrow G$. 

\noindent Let $f\colon H\rightarrow \R$ with $|\text{supp}(f)|\le n$ and $\|\nabla_{S_{H}}f\|_{1}=1$. From this function, we define a random map as follows: let $x\in X$ and let 
\begin{align*}
    f_{x}\colon G &\longrightarrow \R \\
    g &\longmapsto f(c_{G,H}(g,x)).
\end{align*}
From the fact that $c_{G,H}(\cdot,x)\colon G\rightarrow H$ is a bijection for almost every $x\in X$, we already have 
\begin{equation*}
    |\text{supp}(f_{x})|=|\text{supp}(f)|, \; \|f_{x}\|_{1}=\|f\|_{1}. 
\end{equation*}
We now look at the norm of the $\ell^1-$gradient of $f_{x}$. We first have 
\begin{align*}
    &\int_{X} \|\nabla_{S_{G}}f_{x}\|_{1}\mathrm{d}\mu(x)=\sum_{s\in S_{G}}\sum_{g\in G}\int_{X} \big|f_{x}(sg)-f_{x}(g)\big|\;\mathrm{d}\mu(x) \\
    &=\sum_{s\in S_{G}}\sum_{g\in G}\int_{X} \big|f(c_{G,H}(sg,x))-f(c_{G,H}(g,x))\big|\;\mathrm{d}\mu(x)\\
    &=\sum_{s\in S_{G}}\sum_{g\in G}\int_{X} \big|f\big(c_{G,H}(s,g\cdot x)c_{G,H}(g,x)\big)-f(c_{G,H}(g,x))\big|\;\mathrm{d}\mu(x)\\
    &=\sum_{s\in S_{G}}\sum_{g\in G}\int_{X} \big|f\big(c_{G,H}(s, x)c_{G,H}(g,g^{-1}\cdot x)\big)-f\big(c_{G,H}(g,g^{-1}\cdot x)\big)\big|\;\mathrm{d}\mu(x)\\
    &=\sum_{s\in S_{G}}\sum_{g\in G}\int_{X} \big|f\big(c_{G,H}(s, x)c_{G,H}(g^{-1},x)^{-1}\big)-f\big(c_{G,H}(g^{-1},x)^{-1}\big)\big|\;\mathrm{d}\mu(x)\\
    &=\sum_{s\in S_{G}}\sum_{h\in H}\int_{X} \big|f\big(c_{G,H}(s, x)h\big)-f(h)\big|\;\mathrm{d}\mu(x)\\
    &=\sum_{s\in S_{G}}\int_{X}\frac{\sum_{h\in H}\big|f(c_{G,H}(s,x)h)-f(h)\big|}{\varphi\big(\sum_{h\in H}\big|f(c_{G,H}(s,x)h)-f(h)\big|\big)}\cdot \varphi\left(\sum_{h\in H}\big|f(c_{G,H}(s,x)h)-f(h)\big|\right)\; \mathrm{d}\mu(x)
\end{align*}
using successively Fubini's theorem, the definition of $f_{x}$, the cocycle identity, the change of variable $x\in X\longmapsto g\cdot x\in X$ and the fact that the $G-$action on $X$ preserves $\mu$, the cocycle identity and the change of variable $c_{G,H}(\lbrace\cdot\rbrace^{-1},x)^{-1}\colon G\rightarrow H$. Since the sum $\sum_{h\in H}\left|f(c_{G,H}(s,x)h)-f(h)\right|$ is bounded from above by $2\|f\|_{1}$, the monotonicity of $t\longmapsto \frac{t}{\varphi(t)}$ implies that  
\begin{equation}\label{eq:eq1.4}
    \frac{\sum_{h\in H}\left|f(c_{G,H}(s,x)h)-f(h)\right|}{\varphi\left(\sum_{h\in H}\left|f(c_{G,H}(s,x)h)-f(h)\right|\right)} \le \frac{2\|f\|_{1}}{\varphi(2\|f\|_{1})}.
\end{equation}
For bounding $\varphi\left(\displaystyle\sum_{h\in H}\big|f(c_{G,H}(s,x)h)-f(h)\big|\right)$, write
\begin{equation*}
    c_{G,H}(s,x)=t_{1}\dots t_{k}
\end{equation*}
where $k=|c_{G,H}(s,x)|_{S_{H}}$ and $t_{i}\in S_{H}$ for every $1\le i\le k$. Then the triangle inequality provides
\begin{align*}
    \big|f(c_{G,H}(s,x)h)-f(h)\big|&=\big|f(t_{1}\dots t_{k}h)-f(h)\big| \\
    &\le \sum_{i=1}^{k}\big|f(t_{i}\dots t_{k}h)-f(t_{i+1}\dots t_{k}h)\big|. 
\end{align*}
Then the change of variables $h\in H \longmapsto t_{i+1}\dots t_{k}h\in H$ provides 
\begin{align*}
    \sum_{h\in H}\big|f(c_{G,H}(s,x)h)-f(h)\big|&\le \sum_{i=1}^{k}\sum_{h\in H}\big|f(t_{i}\dots t_{k}h)-f(t_{i+1}\dots t_{k}h)\big|\\
    &=\sum_{i=1}^{k}\sum_{h\in H}\left|f(t_{i}h)-f(h)\right|\\
    &\le \sum_{i=1}^{k}\|\nabla_{S_{H}}f\|_{1} \\
    &\le k\cdot \|\nabla_{S_{H}}f\|_{1} \\
    &=k\\
    &=|c_{G,H}(s,x)|_{S_{H}}.
\end{align*}
This inequality combined with Equation (\ref{eq:eq1.4}) then yields
\begin{equation*}
    \int_{X} \|\nabla_{S_{G}}f_{x}\|_{1}\mathrm{d}\mu(x) \le \frac{2\|f\|_{1}}{\varphi(2\|f\|_{1})}\cdot \sum_{s\in S_{G}}\int_{X}\varphi\big(|c_{G,H}(s,x)|_{S_{H}}\big)\mathrm{d}\mu(x)=\frac{2\|f\|_{1}}{\varphi(2\|f\|_{1})}\cdot K
\end{equation*}
where $K\defeq \displaystyle\sum_{s\in S_{G}}\int_{X}\varphi\big(|c_{G,H}(s,x)|_{S_{H}}\big)\mathrm{d}\mu(x)$ is finite by $\varphi-$integrability of the cocycle $c_{G,H}$ (recall from the remark before the proof that, under our assumptions, constants in the definition of $\varphi-$integrability can be removed). Thus we deduce that there is $x_{0}\in X$ such that 
\begin{equation*}
    \|\nabla_{S_{G}}f_{x_{0}}\|_{1} \le \frac{2\|f\|_{1}}{\varphi(2\|f\|_{1})}\cdot K.
\end{equation*}
Since $|\text{supp}(f_{x_{0}})| \le n$ and $\|f_{x_{0}}\|_{1}=\|f\|_{1}$, it follows that 
\begin{align*}
    \prof{G}(n) \ge \frac{\|f_{x_{0}}\|_{1}}{\|\nabla f_{x_{0}}\|_{1}} \ge \frac{1}{2K}\varphi(2\|f\|_{1}) \ge \frac{1}{2K}\varphi(\|f\|_{1}).
\end{align*}
Lastly, taking the supremum over all maps $f\colon H\rightarrow \R$ such that $|\text{supp}(f)|\le n$ and $\|\nabla_{S_{H}} f\|_{1}=1$, we get indeed 
\begin{equation*}
    \prof{G}(n) \ge \frac{1}{2K}\varphi(\prof{H}(n))
\end{equation*}
whence $\varphi\circ \prof{H}(n) \preccurlyeq \prof{G}(n)$. 

\noindent \textit{(ii)} Fix $\varphi(t)=t^p$. The strategy is similar to the previous one. Fix $f\colon H\longrightarrow\R$ with $|\text{supp}(f)|\le n$ and $\|\nabla_{S_{H}}f\|_{p}=1$. For $x\in X$ we define $f_{x}$ in the same way as above, so we still have 
\begin{equation*}
    |\text{supp}(f_{x})|=|\text{supp}(f)|, \; \|f_{x}\|_{p}=\|f\|_{p}
\end{equation*}
but we now average rather $\|\nabla_{S_{G}}f_{x}\|_{p}^{p}$ to deduce information on the $\ell^{p}-$norm of the gradient of $f_{x}$. With the same steps, one proves that 
\begin{equation*}
    \int_{X}\|\nabla_{S_{G}}f_{x}\|_{p}^{p}\;\mathrm{d}\mu(x)=\sum_{s\in S_{G}}\int_{X}\sum_{h\in H}\big|f(c_{G,H}(s,x)h) -f(h)\big|^{p}\;\mathrm{d}\mu(x).
\end{equation*}
Fixing $s\in S_{G}$ and writing $c_{G,H}(s,x)=t_{1}\dots t_{k}$ with $k=|c_{G,H}(s,x)|_{S_{H}}$, $t_{i}\in S_{H}$, the triangle inequality gives 
\begin{align*}
    \sum_{h\in H}\big|f\big(c_{G,H}(s,x)h\big) -f(h)\big|^{p}&=\|f(c_{G,H}(s,x)\cdot)-f(\cdot)\|_{p}^{p} \\
    &\le \left(\sum_{i=1}^{k}\|f(t_{i}\dots t_{k}\cdot)-f(t_{i+1}\dots t_{k}\cdot)\|_{p}\right)^{p} \\
    &=\left(\sum_{i=1}^{k}\|f(t_{i}\cdot)-f(\cdot)\|_{p}\right)^{p} \\
    &\le \big(k\cdot \|\nabla_{S_{H}}f\|_{1}\big)^{p} \\
    &=k^p \\
    &=|c_{G,H}(s,x)|_{S_{H}}^{p}.
\end{align*}
Thus we deduce that 
\begin{align*}
    \int_{X}\|\nabla_{S_{G}}f_{x}\|_{p}^{p}\;\mathrm{d}\mu(x)&=\sum_{s\in S_{G}}\int_{X}\sum_{h\in H}\big|f\big(c_{G,H}(s,x)h\big) -f(h)\big|^{p}\;\mathrm{d}\mu(x) \\
    &\le \sum_{s\in S_{G}}\int_{X}|c_{G,H}(s,x)|_{S_{H}}^{p}\;\mathrm{d}\mu(x)
\end{align*}
and the latter quantity, call it $K$, is finite since $c_{G,H}$ is $\ld^{p}-$integrable. Hence there is $x_{0}\in X$ such that $\|\nabla_{S_{G}}f_{x_{0}}\|_{p}^{p} \le K$, and thus 
\begin{equation*}
    \profp{G}(n) \ge \frac{\|f_{x_{0}}\|_{p}}{\|\nabla_{S_{G}}f_{x_{0}}\|_{p}} \ge \frac{1}{K^{\frac{1}{p}}}\|f\|_{p}.
\end{equation*}
It remains to take the supremum over all $f\colon H\rightarrow \R$ having $|\text{supp}(f)|\le n$ and $\|\nabla_{S_{H}}f\|_{p}=1$ to conclude that 
\begin{equation*}
    \profp{G}(n) \ge \frac{1}{K^{\frac{1}{p}}}\profp{H}(n). 
\end{equation*}
\end{proof}

As an application of this asymptotic inequality, the authors of~\cite{DKLMT22} deduce that:

\begin{corollary}
Let $k,\ell \ge 1$ be two integers. Let $p\ge 0$. If there is a $(\ld^p, \ld^0)-$orbit equivalence coupling from $\Z^{k+\ell}$ to $\Z^{k}$, then $p\le \frac{k}{k+\ell}$. 
\end{corollary}

\begin{proof}
We consider two cases. Suppose first that $0\le p\le 1$ and that there is a $(\ld^p, \ld^0)-$orbit equivalence coupling from $\Z^{k+\ell}$ to $\Z^{k}$. Then Theorem~\ref{thm:ObstructionDKLMT}\textit{(i)} implies that 
\begin{equation*}
    n^{\frac{p}{k}} \simeq \prof{\Z^{k}}(n)^{p} \preccurlyeq \prof{\Z^{k+\ell}}(n) \simeq n^{\frac{1}{k+\ell}}
\end{equation*}
whence $\frac{p}{k} \le \frac{1}{k+\ell}$, i.e. $p\le \frac{k}{k+\ell}$.

\vspace{0.15cm}

\noindent Assume now that $p>1$. Then a $(\ld^p, \ld^0)$ coupling from $\Z^{k+\ell}$ to $\Z^{k}$ would be in particular $(\ld^q,\ld^0)$ for any $\frac{k}{k+\ell}<q\le 1$, which contradicts the first case proved above. 
\end{proof}

In addition, it is also proved in~\cite[Theorem~6.12]{DKLMT22} that such a coupling does exist if $p<\frac{k}{k+\ell}$, using a particular type of F\o lner sequences of amenable groups, called~\textit{F\o lner tilings}. These sequences are particularly useful for constructing free p.m.p. actions on product spaces, and thus provide orbit equivalence couplings. If additionally one can control precisely the geometry of the tiles we use, one can then control the integrability of cocycles of the coupling with respect to a given function. We give more details on this particular technique in Chapter~\ref{chap:chapter6}, where we use it to build couplings between lamplighter-like groups.

Thus, for free abelian groups, the only thing 
unclear in~\cite{DKLMT22} was whether an orbit coupling from $\Z^{k+\ell}$ to $\Z^{k}$ could be $(\ld^{\frac{k}{k+\ell}}, \ld^0)$. Correia recently answered the question, showing that the threshold $p=\frac{k}{k+\ell}$ cannot be reached~\cite[Theorem~4.3]{Cor25}. The general statement, that will also be useful below for couplings between halo products, is the following:
\begin{theorem}[{\cite[Theorem~B]{Cor25}}]\label{thm:threshold}
Let $G$ and $H$ be finitely generated groups. Assume there exists a non-decreasing function $h_{G}$ and an increasing function $h_{H}$ satisfying $h_{G}(x) \simeq\prof{G}(x)$, $h_{H}(x)\simeq\prof{H}(x)$ and the following as $x\rightarrow+\infty$:
\begin{equation*}
		 h_{G}(x)=o\left(h_{H}(x)\right),
\end{equation*}
\begin{equation*}
\forall C>0,\; h_{G}(Cx)=O\left(h_{G}(x)\right),
\end{equation*}
\begin{equation*}
\forall C>0,\; h_{G} \circ  h_{H}^{-1}(Cx)=O\left(h_{G}\circ h_{H}^{-1}(x)\right).
\end{equation*}
Then there is no $(h_{G}\circ h_{H}^{-1},\ld^0)-$integrable orbit equivalence coupling from $G$ to $H$.
\end{theorem}

In conclusion, there exists a $(\ld^p,\ld^0)-$orbit equivalence from $\Z^{k+\ell}$ to $\Z^{k}$ if and only if $p<\frac{k}{k+\ell}$.

Turning to another class of groups, there also exists a method to construct couplings between wreath products provided couplings between factors of the wreath products. 

\begin{theorem}[{\cite[Corollary~7.3]{DKLMT22}}]\label{thm:stabilityofOEforwreathproducts}
Let $G_{1},H_{1},G_{2},H_{2}$ be finitely generated groups. Let $\varphi,\psi\colon\R_{+}\rightarrow\R_{+}$.
If there is a $(\varphi,\psi)-$integrable orbit equivalence coupling from $G_{1}$ to $G_{2}$, and if there is a $(\varphi,\psi)-$integrable orbit equivalence coupling from $H_{1}$ to $H_{2}$, then there is a $(\varphi,\psi)-$integrable orbit equivalence coupling from $G_{1}\wr H_{1}$ to $G_{2}\wr H_{2}$.
\end{theorem}

Combined with the computation of the isoperimetric profile of wreath products from Theorem~\ref{thm:formulaforprofilewreathproducts}, the isoperimetric obstruction from Theorem~\ref{thm:ObstructionDKLMT} and the description at the critical threshold given in~\cite{Cor25}, it follows for instance that:

\begin{corollary}
Let $k,\ell\ge 1$ be integers. Let $F$ be a non-trivial finite group. Then there is an $(\ld^p, \ld^0)-$orbit equivalence from $F\wr\Z^{k+\ell}$ to $F\wr\Z^{k}$ if and only if $p<\frac{k}{k+\ell}$.
\end{corollary}

We can also formulate an iterated version of this equivalence, if for a non-trivial finite group $F$, we define iteratively the $n-$th iterated wreath product $H_{n}(k)$ as: $H_{0}(k)\defeq \Z^k$, $H_{n+1}(k)\defeq F\wr H_{n}(k)$.

\begin{corollary}
Let $k,\ell\ge 1$ be integers. Let $n\ge 0$. Then there is an $(\ld^p, \ld^0)-$orbit equivalence from $H_{n}(k+\ell)$ to $H_{n}(k)$ if and only if $p<\frac{k}{k+\ell}$. 
\end{corollary}

These descriptions for quantitative orbit equivalence of (iterated) lamplighters serve as a motivation to try to achieve similar description for variants of lamplighters, such as lampshufflers or lampcloners. See Chapter~\ref{chap:chapter6}.

On the other hand, it is also proved in~\cite{DKLMT22} that other classical quasi-isometry invariants cannot provide any valuable information on the integrability of the cocycles of a coupling. One such geometric invariant is finite presentability:

\begin{theorem}[{\cite[Corollary~8.2]{DKLMT22}}]
Finite presentability is unstable under the equivalence relation of strongly exponential orbit equivalence. In particular, it is unstable under $\ld^{<\infty}$ orbit equivalence. 
\end{theorem}

Among many other consequences, methods from~\cite{DKLMT22} allow one to prove a far-reaching extension of Theorem~\ref{thm:profilQIinvariant} and Lemma~\ref{lem:monotonieprofilErschler}:

\begin{theorem}[{\cite[Corollary~5.6]{DKLMT22}}]\label{thm:profilesmonotonuousregularmaps}
Let $G$ and $H$ be finitely generated amenable groups. Let $p\ge 1$. If there exists a regular map from $G$ to $H$, then $\profp{H}(n) \preccurlyeq \profp{G}(n)$. 
\end{theorem}

The statement is false if we drop the amenability assumption. For instance, it is well-known that $\Z/2\Z\wr\Z$ contains biLipschitz embedded copies of $F_{2}$~\cite{Woe05} (which is not amenable by Example~\ref{ex:examplesofAgroups}\textit{(iii)}), and indeed the isoperimetric profile of $\Z/2\Z\wr\Z$ cannot be bounded by the one of $F_{2}$.

Finally, let us mention the following result of composition of actions, which enables us to have some kind of transitivity of quantitative orbit equivalence among finitely generated amenable groups. This is a reformulation of~\cite[Propositions~2.26,~2.29 and~2.30]{DKLMT22} in terms of orbit equivalence.

\begin{theorem}\label{thm:CompositionCouplings}
Let $G_{1}, G_{2}, G_{3}$ be three finitely generated groups and $\varphi,\varphi',\psi,\psi'\colon\R_{+}\rightarrow\R_{+}$ be non-decreasing maps. Let us assume that there exists a $(\varphi,\psi)-$integrable orbit equivalence coupling from $G_{1}$ to $G_{2}$, and a $(\varphi',\psi')-$integrable orbit equivalence coupling from $G_{2}$ to $G_{3}$. Let us also assume that the pairs $(\varphi,\varphi')$ and $(\psi',\psi)$ lie in the set of pairs of functions $(f,g)$ satisfying one of the following two conditions:
\begin{enumerate}[label=(\roman*)]
    \item\label{item:CompositionHyp1} either $f$ and $g$ are subadditive and $g$ is concave;
    \item\label{item:CompositionHyp2} or $f(x)\simeq g(x)\simeq x^p$ for some $p\ge 1$.
\end{enumerate}
Then there exists a $(\varphi'',\psi'')-$integrable orbit equivalence coupling from $G_{1}$ to $G_{3}$, where
\begin{equation*}
\varphi''(x)=\left\{\begin{array}{ll}
        \varphi\circ\varphi'(x) & \text{if }(\varphi,\varphi')\text{ satisfies Assumption~\textit{\ref{item:CompositionHyp1}}}\\
        x^p & \text{if }(\varphi,\varphi')\text{ satisfies Assumption~\textit{\ref{item:CompositionHyp2}} for some }p\ge 1
\end{array}\right.
\end{equation*}
    and
\begin{equation*}
\psi''(x)=\left\{\begin{array}{ll}
        \psi'\circ\psi(x) & \text{if }(\psi',\psi)\text{ satisfies Assumption~\ref{item:CompositionHyp1}}\\
        x^p & \text{if }(\psi',\psi)\text{ satisfies Assumption~\ref{item:CompositionHyp2} for some }p\ge 1
\end{array}\right..  
\end{equation*}
\end{theorem}

\section{Measure-scaling quasi-isometries}\label{sec:ScalingQI}

This part introduces measure-scaling quasi-isometries in the framework of bounded degree graphs. It follows mainly~\cite{GT22}.

\subsection{Definition and first properties}

The main definition of this section is the following one. 

\begin{definition}\label{def:measurescalingQI}
Let $X$ and $Y$ be bounded degree graphs. Let $k>0$. A quasi-isometry $f\colon X\rightarrow Y$ is~\textit{quasi-$k$-to-one} if there exists $C>0$ such that 
\begin{equation*}
    \left|k|A|-|f^{-1}(A)|\right| \le C\cdot |\partial_{Y}A|
\end{equation*}
for all finite subsets $A\subset Y$, where $\partial_{Y}A \defeq \lbrace y\in Y\setminus A : \exists a\in A, y\sim_{Y}a\rbrace$ is the~\textit{boundary} of $A$ in $Y$.
\end{definition}

In this case, we say that $f$ is~\textit{measure-scaling}, and that the real number $k>0$ is the~\textit{scaling factor}. 

Informally, being quasi-$k$-to-one means that the pre-image of a finite subset of the target space is close to have the cardinality of the finite subset multiplied by $k$. The scaling factor should therefore be thought as a measurement of how far our map is from a bijection. The next theorem, due to Whyte, also serves as a motivation for this viewpoint. 

\begin{theorem}[{\cite{Why99};~\cite[Proposition~4.1]{GT22}}]\label{thm:Whytethm}
Let $X$ and $Y$ be bounded degree graphs. A quasi-isometry $f\colon X\rightarrow Y$ is quasi-one-to-one if and only if it lies at bounded distance from a bijection. 
\end{theorem}

\begin{remark}\label{rem:ScalingQINAspaces}
Observe also that between non-amenable spaces, Definition~\ref{def:measurescalingQI} is satisfied for any quasi-isometry $f\colon X\rightarrow Y$ and any $k>0$. Indeed, if $Y$ is not amenable, there exists $\varepsilon>0$ such that $|\partial_{Y}A|>\varepsilon\cdot |A|$ for any finite subset $A\subset Y$, and thus 
\begin{align*}
    \left|k|A|-|f^{-1}(A)|\right| &\le k|A|+|f^{-1}(A)| \\
    &\le (k+P)\cdot |A| \\
    &\le \frac{k+P}{\varepsilon}\cdot |\partial_{Y}A|
\end{align*}
where $P\ge 1$ is a uniform bound on $|f^{-1}(\lbrace y\rbrace)|$ for $y\in Y$. Thus $f$ is quasi-$k$-to-one. 
\end{remark}

In particular, from this remark and Theorem~\ref{thm:Whytethm}, we deduce that any quasi-isometry $X\rightarrow Y$ between non-amenable graphs lies at bounded distance from a bijection. 

The situation is drastically different for amenable spaces. 

\begin{lemma}\label{lem:uniquenessscalingfactor}
Let $X$ and $Y$ be amenable bounded degree graphs. If $f\colon X\rightarrow Y$ is quasi-$k$-to-one and quasi-$k'$-to-one, then $k=k'$. 
\end{lemma}

\begin{proof}
As $f$ is quasi-$k$-to-one, there is $C>0$ such that 
\begin{equation*}
    \left|k|A|-|f^{-1}(A)|\right| \le C\cdot |\partial_{Y}A|
\end{equation*}
for all finite subsets $A\subset Y$. In particular, applying this inequality to a F\o lner sequence $(F_{n})_{n\in\N}$ of $Y$ leads to 
\begin{equation*}
     \left|k|F_{n}|-|f^{-1}(F_{n})|\right| \le C\cdot |\partial_{Y}F_{n}|
\end{equation*}
for any $n\in\N$, or equivalently 
\begin{equation*}
     \left|k-\frac{|f^{-1}(F_{n})|}{|F_{n}|}\right| \le C\cdot \frac{|\partial_{Y}F_{n}|}{|F_{n}|}
\end{equation*}
for any $n\in\N$. Thus we get that $k=\lim\limits_{n\rightarrow\infty}\frac{|f^{-1}(F_{n})|}{|F_{n}|}$, and by the same reasoning, this last limit also equals $k'$, whence $k=k'$.
\end{proof}

Notice also that in the amenable case, a quasi-isometry may not be measure-scaling: for instance, the map $f\colon \Z\rightarrow\Z$ defined by $f(n)=n$ if $n\ge 0$ and $f(n)=3n$ if $n\le 0$ is not measure-scaling. 

On the other hand, many classical quasi-isometries turn out to be scaling. For instance, if $H$ is a proper finite-index subgroup of an amenable finitely generated group $G$, then the natural inclusion $H\hookrightarrow G$ is quasi-$\frac{1}{[G:H]}$-to-one. 

Let us now analyse how scaling maps behave under classical operations. 

\begin{proposition}\label{prop:stabilitypropertiesforscalingQI}
Let $X$, $Y$ and $Z$ be bounded degree graphs. Let $k_{1},k_{2}>0$ and let $f,h\colon X\rightarrow Y$, $g\colon Y\rightarrow Z$ be three quasi-isometries.
\begin{enumerate}[label=(\roman*)]
    \item If $f$ is quasi-$k_{1}$-to-one and $h$ is at bounded distance from $f$, then $h$ is quasi-$k_{1}$-to-one. 
    \item If $f$ is quasi-$k_{1}$-to-one and $g$ is quasi-$k_{2}$-to-one, then $g\circ f\colon X\rightarrow Z$ is quasi-$k_{1}k_{2}$-to-one. 
    \item If $f$ is quasi-$k_{1}$-to-one, then any of its quasi-inverses is quasi-$\frac{1}{k_{1}}$-to-one. 
\end{enumerate}
\end{proposition}

The proof of Proposition~\ref{prop:stabilitypropertiesforscalingQI} requires several intermediate observations. 

\begin{lemma}\label{lem:Boundingneighborhoodsingraphs}
Let $X$ be a bounded degree graph. 
\begin{enumerate}[label=(\roman*)]
    \item For any finite subset $A\subset X$ and any $S\ge 0$, one has $|A^{+S}|\le N^{S}\cdot|A|$, where $N\ge 3$ is an integer larger than the maximal degree of a vertex of $X$.
    \item For any finite subset $A\subset X$ and any $S\ge 0$, there is a constant $R>0$ (depending only on $S$ and on $X$) such that $|A^{+S}\setminus A| \le R\cdot |\partial_{X} A|$.
\end{enumerate}
\end{lemma}

\begin{proof}
\textit{(i)} Fix any $0\le i\le S-1$. Given any element $a$ of $A$, the fact that any vertex of $X$ has degree $\le N$ implies that there are at most $N^{i}$ paths of length $i$ that start at $a\in A$ and that end at an element of $A^{+i}$. Thus 
\begin{equation*}
    \left|A^{+S}\right|\le \left|A\right|\cdot \sum_{i=0}^{S-1}N^{i} \le N^{S}\cdot \left|A\right|
\end{equation*}
as claimed. 

\noindent \textit{(ii)} It suffices to note that $A^{+S}\setminus A \subset (\partial_{X}A)^{+(S-1)}$ and to apply point \textit{(i)}.
\end{proof}

\begin{lemma}\label{lem:Boundingboundariesofpreimages}
Let $X$ and $Y$ be two graphs of bounded degree. Let $f\colon X\rightarrow Y$ be a quasi-isometry. Then there exists a constant $L>0$ such that 
\begin{equation*}
    |\partial_{X}f^{-1}(A)| \le L \cdot |\partial_{Y}A|
\end{equation*}
for any finite subset $A\subset Y$.
\end{lemma}

\begin{proof}
Let $C\ge 1$, $K\ge 0$ be such that $f$ is a $(C,K)-$quasi-isometry. Notice that $f$ sends two adjacent vertices to two vertices at distance $\le C+K$. Let $P\ge 1$ be a uniform bound on pre-images of points under $f$. 

\noindent Let $A\subset Y$ be finite. Let $x\in \partial_{X}f^{-1}(A)$. Thus there is $y\in f^{-1}(A)$ such that $d_{X}(x,y)=1$. It follows that $d_{Y}(f(x),f(y)) \le C+K$ and $f(y)\in A$, so $f(x)\in A^{+(C+K)}$, and thus $x\in f^{-1}(A^{+(C+K)})$. Additionally, $x\notin f^{-1}(A)$, so we have proved 
\begin{equation*}
    \partial_{X}f^{-1}(A) \subset f^{-1}\big(A^{+(C+K)}\big)\setminus f^{-1}(A).
\end{equation*}
Taking cardinalities, it follows that 
\begin{align*}
    |\partial_{X}f^{-1}(A)| &\le \big|f^{-1}\big(A^{+(C+K)}\big)\setminus f^{-1}(A)\big| \\
    &=\big|f^{-1}\big(A^{+(C+K)}\setminus A\big)\big| \\
    &\le P\cdot \big|A^{+(C+K)}\setminus A\big| \\
    &\le P\cdot R\cdot |\partial_{Y}A|
\end{align*}
where $R>0$ is the constant provided by Lemma~\ref{lem:Boundingneighborhoodsingraphs}\textit{(ii)}, that depends only on $Y$ and the parameters $C$ and $K$ of $f$. This finishes the proof. 
\end{proof}

\begin{lemma}\label{lm:preimagesandquasiinverses}
Let $X$ and $Y$ be bounded degree graphs, and let $f\colon X\rightarrow Y$ be a $(C,K)-$quasi-isometry, with a $(C,K)-$quasi-inverse $\overline{f}\colon Y\rightarrow X$. Then we have 
\begin{equation*}
    d_{\text{Haus}}\left(\overline{f}(A),f^{-1}(A^{+K})\right)\le (C+2)K
\end{equation*}
for any subset $A\subset Y$.
\end{lemma}

\begin{proof}
Let $A\subset Y$, and let $x\in f^{-1}(A^{+K})$. Then $f(x)\in A^{+K}$, and we may pick $y\in A$ such that $d_{Y}(y, f(x))\le K$. It follows that 
\begin{align*}
    d_{X}\big(\overline{f}(y), x\big) &\le d_{X}\big(\overline{f}(y), \overline{f}(f(x))\big)+d_{X}\big(\overline{f}(f(x)), x\big) \\
    &\le C\cdot d_{Y}(y, f(x))+K+K \\
    &\le (C+2)K
\end{align*}
and since $\overline{f}(y)\in \overline{f}(A)$, we get $x\in \overline{f}(A)^{+(C+2)K}$. 

\noindent Conversely, if $x=\overline{f}(y)$ with $y\in A$, then $f(x)=f(\overline{f}(y))$ is at distance at most $K$ from $y\in A$, i.e. $f(x)\in A^{+K}$, so $x\in f^{-1}(A^{+K})$. We conclude that
\begin{equation*}
    d_{\text{Haus}}\left(\overline{f}(A),f^{-1}(A^{+K})\right) \le (C+2)K
\end{equation*}
as claimed. 
\end{proof}

\begin{proof}[Proof of Proposition~\ref{prop:stabilitypropertiesforscalingQI}]
\textit{(i)} Assume that the distance between $f$ and $h$ is at most $Q$, and observe that this implies that $f^{-1}(B)\subset h^{-1}(B^{+Q})$ and $h^{-1}(B)\subset f^{-1}(B^{+Q})$ for any finite subset $B\subset Y$. 

\noindent Now, fix any finite subset $A\subset Y$. We first estimate 
\begin{equation}\label{eq:finitedistancechange}
    \left|k_{1}|A|-|h^{-1}(A)|\right| \le \left|k_{1}|A|-|f^{-1}(A)|\right| + \left||f^{-1}(A)|-|h^{-1}(A)|\right|.
\end{equation}
We then have 
\begin{align*}
    |f^{-1}(A)|-|h^{-1}(A)| &\le |h^{-1}(A^{+Q})|-|h^{-1}(A)| \\
    &\le \big|h^{-1}\big(A^{+Q}\setminus A\big)\big| \\
    &\le P\cdot \big|A^{+Q}\setminus A\big| \\
    &\le P\cdot R\cdot |\partial_{Y}A|
\end{align*}
where $P\ge 1$ is a uniform bound on sizes of pre-images of points under $h$ and $R>0$ is the constant from Lemma~\ref{lem:Boundingneighborhoodsingraphs}\textit{(ii)}. On the other hand, we have similarly 
\begin{align*}
    |h^{-1}(A)|-|f^{-1}(A)| &\le |f^{-1}(A^{+Q})|-|f^{-1}(A)| \\
    &\le |f^{-1}\big(A^{+Q}\setminus A\big)| \\
    &\le P'\cdot \big|A^{+Q}\setminus A\big| \\
    &\le P'\cdot R\cdot |\partial_{Y}A|
\end{align*}
where $P'$ is a uniform bound on sizes of pre-images of points under $f$. Thus, combining (\ref{eq:finitedistancechange}) and the fact that $f$ is quasi-$k_{1}$-to-one, we get the existence of $C>0$ such that 
\begin{align*}
    \left|k_{1}|A|-|h^{-1}(A)|\right| &\le \left|k_{1}|A|-|f^{-1}(A)|\right| + \left||f^{-1}(A)|-|h^{-1}(A)|\right| \\
    &\le C\cdot |\partial_{Y}A| + \max(P,P')\cdot R\cdot |\partial_{Y}A| \\
    &=\big(C+\max(P,P')\cdot R\big)\cdot |\partial_{Y}A|
\end{align*}
and this proves that $h$ is quasi-$k_{1}$-to-one.

\noindent \textit{(ii)} Fix a finite subset $A\subset Z$. By our assumptions, we know there exist $C>0$ and $D>0$ such that
\begin{align*}
    \left|k_{1}k_{2}|A|-|(g\circ f)^{-1}(A)|\right| &\le \left|k_{1}k_{2}|A|-k_{1}|g^{-1}(A)|\right|+\left| k_{1}|g^{-1}(A)|-|f^{-1}(g^{-1}(A))|\right| \\
    &\le k_{1}\cdot D\cdot |\partial_{Z}A|+|\partial_{Y}g^{-1}(A)| \\
    &\le k_{1}\cdot D\cdot |\partial_{Z}A|+L\cdot |\partial_{Z}A| \\
    &\le (k_{1}\cdot D+L)\cdot |\partial_{Z}A|
\end{align*}
where $L>0$ is the constant provided by Lemma~\ref{lem:Boundingboundariesofpreimages}. This proves that $g\circ f$ is quasi-$k_{1}k_{2}$-to-one. 

\noindent \textit{(iii)} Let $C\ge 1$ and $K\ge 0$ be such that $f$ and one of its quasi-inverses $\overline{f}\colon Y\rightarrow X$ are $(C,K)-$quasi-isometries, and such that $f\circ \overline{f}$, $\overline{f}\circ f$ lie at distance $\le K$ from $\text{Id}_{Y}$, $\text{Id}_{X}$ respectively. Fix a finite subset $A\subset X$. 

\noindent We first rewrite 
\begin{align*}
    \left|\frac{1}{k_{1}}|A|-\big|\overline{f}^{-1}(A)\big|\right|&=\frac{1}{k_{1}}\left||A|-k_{1}\big|\overline{f}^{-1}(A)\big|\right| \\
    &\le \frac{1}{k_{1}}\left|k_{1}\big|\overline{f}^{-1}(A)\big|-\big|f^{-1}\big(\overline{f}^{-1}(A)\big)\big|\right|+\frac{1}{k_{1}}\left|\big|f^{-1}\big(\overline{f}^{-1}(A)\big)\big|-|A|\right|.
\end{align*}
In this sum, we control the first term using that $f$ is quasi-$k_{1}$-to-one: there exists a constant $T>0$ such that 
\begin{align*}
    \left|k_{1}\big|\overline{f}^{-1}(A)\big|-\big|f^{-1}\big(\overline{f}^{-1}(A)\big)\big|\right| &\le T\cdot \big|\partial_{Y}\overline{f}^{-1}(A)\big| \\
    &\le T\cdot L\cdot |\partial_{X}A|
\end{align*}
where $L>0$ is provided by Lemma~\ref{lem:Boundingboundariesofpreimages} (and depends only on $X$ and the parameters $C$ and $K$ of $\overline{f})$. For the second term, we write
\begin{equation*}
    \left|\big|f^{-1}\big(\overline{f}^{-1}(A)\big)\big|-|A|\right| \le \big|A^{+K}\big|-|A| = \big|A^{+K}\setminus A\big| \le R\cdot |\partial_{X}A|
\end{equation*}
since $\overline{f}\circ f$ is at distance $\le K$ from $\text{Id}_{X}$ and where $R>0$ comes from Lemma~\ref{lem:Boundingneighborhoodsingraphs}\textit{(ii)} (and depends only on $X$ and the parameters $C$ and $K$ of $f$). By symmetry, one also gets that 
\begin{equation*}
    |A|-\big|f^{-1}\big(\overline{f}^{-1}(A)\big)\big| \le \big|\big(\overline{f}\circ f\big)^{-1}(A^{+K}\setminus A)\big| \le P\cdot \big|A^{+K}\setminus A\big| \le P\cdot R\cdot |\partial_{X}A| 
\end{equation*}
where $P\ge 1$ is a uniform bound on the sizes of pre-images of points under $\overline{f}\circ f$. Combining all these inequalities, it follows that 
\begin{align*}
    \left|\frac{1}{k_{1}}|A|-\big|\overline{f}^{-1}(A)\big|\right|&\le  \frac{1}{k_{1}}\left|k_{1}\big|\overline{f}^{-1}(A)\big|-\big|f^{-1}\big(\overline{f}^{-1}(A)\big)\big|\right|+\frac{1}{k_{1}}\left|\big|f^{-1}\big(\overline{f}^{-1}(A)\big)\big|-|A|\right| \\
    &\le \frac{T\cdot L}{k_{1}}|\partial_{X}A|+\frac{P}{k_{1}}\cdot |\partial_{X}A| \\
    &=\frac{T\cdot L + P}{k_{1}}\cdot |\partial_{X}A|
\end{align*}
and we conclude that $\overline{f}\colon Y\rightarrow X$ is quasi-$\frac{1}{k_{1}}$-to-one. This finishes the proof. 
\end{proof}

As an application of these properties,~\cite[Theorem~5.1]{GT22} provides a version of Milnor-Schwarz lemma (cf. Lemma~\ref{lem:Milnor-Schwarz}) in the context of scaling quasi-isometries. As an interesting special case, let us quote:

\begin{corollary}[{\cite[Corollary~5.8]{GT22}}]
Let $G$ be a finitely generated group, and let $H_{1}$ and $H_{2}$ be two finite-index subgroups of $G$. Then there exists a measure-scaling quasi-isometry $H_{1}\rightarrow H_{2}$, of scaling factor $\frac{[G:H_{2}]}{[G:H_{1}]}$.
\end{corollary}

In particular, from this corollary and Whyte's theorem (Theorem~\ref{thm:Whytethm}), one deduces the following nice fact: in a given finitely generated amenable group, if two finite-index subgroups have the same index, then they must be biLipschitz equivalent. See~\cite[Proposition~1.5]{GT22} for a more general statement.

We refer to~\cite[Section~5]{GT22} for other examples of scaling quasi-isometries, in the framework of locally compact groups and metric measure spaces.

\subsection{Rational scaling factors}

The following result gives an alternative description of measure-scaling quasi-isometries when the scaling factor is rational. Since we will use it several times in the sequel to construct quasi-isometries between wreath products, we provide a proof below. 

\begin{theorem}[{\cite[Proposition~4.2]{GT22}}]\label{thm:rationalscalingfactor}
Let $n,m\ge 1$ be natural integers and let $X$ and $Y$ be two bounded degree graphs. Let $f\colon X\rightarrow Y$ be a quasi-isometry. The following statements are equivalent:
\begin{enumerate}[label=(\roman*)]
    \item The map $f$ is quasi-$\frac{m}{n}$-to-one;
    \item There exist a partition $\mathcal{P}_{X}$ (resp. $\mathcal{P}_{Y}$) with uniformly bounded pieces of size $m$ (resp. $n$) and a bijection $\psi\colon \mathcal{P}_{X}\rightarrow \mathcal{P}_{Y}$ such that $f$ is at bounded distance from a map $g\colon X\rightarrow Y$ satisfying $g(P)\subset \psi(P)$ for any $P\in\mathcal{P}_{X}$. 
\end{enumerate}
\end{theorem}

\begin{proof}
We start by showing that~\textit{(ii)} $\Longrightarrow$~\textit{(i)}. Since $f$ and $g$ are at bounded distance, it suffices to prove that $g$ is quasi-$\frac{m}{n}$-to-one and apply Proposition~\ref{prop:stabilitypropertiesforscalingQI} to conclude that $f$ is also quasi-$\frac{m}{n}$-to-one.

\smallskip 

\noindent Let $A\subset Y$ be finite. Let $A^{+}$ denote the union of all pieces of $\mathcal{P}_{Y}$ that contain at least one point of $A$. If $A^{+}$ consists of $k$ pieces, then $g^{-1}(A^{+})$ is a union of $k$ pieces of $\mathcal{P}_{X}$, hence $\frac{|g^{-1}(A^{+})|}{m}=k=\frac{|A^{+}|}{n}$. Thus we get 
\begin{align*}
    \left|\frac{m}{n}|A|-|g^{-1}(A)|\right| &\le \frac{m}{n}\left||A|-|A^{+}|\right|+\left|\frac{m}{n}|A^{+}|-|g^{-1}(A^{+})|\right|+\left||g^{-1}(A^{+})|-|g^{-1}(A)|\right| \\
    &=\frac{m}{n}|A^{+}\setminus A|+|g^{-1}(A^{+}\setminus A)| \\
    &\le \left(\frac{m}{n}+m\right)\cdot |A^{+}\setminus A|
\end{align*}
using for the last inequality that $|g^{-1}(B)| \le m\cdot |B|$ for any finite subset $B\subset Y$. Since pieces of $\mathcal{P}_{Y}$ are uniformly bounded, there is $K\ge 0$ such that $A^{+}\subset A^{+K}$, whence 
\begin{align*}
     \left|\frac{m}{n}|A|-|g^{-1}(A)|\right| &\le \left(\frac{m}{n}+m\right)\cdot |A^{+}\setminus A| \\
     &\le \frac{m}{n}(n+1)\cdot|A^{+K}\setminus A| \\
     &\le \frac{m}{n}(n+1)\cdot R \cdot |\partial_{Y}A|
\end{align*}
using Lemma~\ref{lem:Boundingneighborhoodsingraphs}\textit{(ii)} in the last line. Thus $g$ is quasi-$\frac{m}{n}$-to-one, hence so is $f$, and~\textit{(i)} is proved. 

\smallskip 

\noindent~\textit{(i)} $\Longrightarrow$~\textit{(ii)}: Let $C\ge 1$, $K\ge 0$ be such that $f$ is a $(C,K)-$quasi-isometry. 

\smallskip

\noindent According to~\cite[Lemma~4.3]{GT22}, there exists a partition $\mathcal{P}_{X}$ (resp. $\mathcal{P}_{Y}$) of $X$ (resp. of $Y$) with uniformly bounded pieces of size $m$ (resp. $n$). Let $L$ denote a uniform upper bound on diameters of pieces of $\mathcal{P}_{X}$ and $\mathcal{P}_{Y}$. We define a graph structure on $\mathcal{P}_{X}$ (resp. $\mathcal{P}_{Y}$) by connecting two pieces of $\mathcal{P}_{X}$ (resp. $\mathcal{P}_{Y}$) by an edge whenever they contain adjacent vertices of $X$ (resp. of $Y$). For each $P\in\mathcal{P}_{X}$ (resp. $Q\in\mathcal{P}_{Y}$), we fix a basepoint $x_{P}\in P$ (resp. $y_{Q}\in Q$). Then we define the maps $i\colon \mathcal{P}_{X}\longrightarrow X$, $i(P)=x_{P}$ and 
\begin{align*}
    j\colon Y&\longrightarrow \mathcal{P}_{Y} \\
    y&\longmapsto \text{piece containing} \;y.
\end{align*} 
These maps are quasi-isometries and, applying the implication \textit{(ii)} $\Longrightarrow$ \textit{(i)} already proved, $i$ is quasi-$\frac{1}{m}$-to-one and $j$ is quasi-$n$-to-one. Therefore, Proposition~\ref{prop:stabilitypropertiesforscalingQI}\textit{(ii)} implies that $j\circ f\circ i\colon \mathcal{P}_{X}\longrightarrow \mathcal{P}_{Y}$ is quasi-one-to-one and thus, by Theorem~\ref{thm:Whytethm}, it lies at bounded distance, say $L\ge 0$, from a bijection $\psi\colon \mathcal{P}_{X}\longrightarrow\mathcal{P}_{Y}$. Finally, we define the map 
\begin{align*}
    g\colon X&\longrightarrow Y \\
    x&\longmapsto y_{\psi(P)}
\end{align*}
where $P\in\mathcal{P}_{X}$ contains $x$. By construction, $g(P)\subset \psi(P)$ for any $P\in\mathcal{P}_{X}$. Let $x\in X$ and let $P\in\mathcal{P}_{X}$ be the piece containing $x$. Then one has 
\begin{align*}
    d_{Y}(f(x),g(x)) &\le d_{Y}(f(x),f(x_{P}))+d_{Y}(f(x_{P}),g(x)) \\
    &\le C\cdot d_{X}(x,x_{P})+K+d_{Y}(f(x_{P}),g(x))\\
    &\le C\cdot D+K+d_{Y}(f(x_{P}),g(x)).
\end{align*}
It remains to estimate $d_{Y}(f(x_{P}),g(x))$. If we let $A\ge 1$ and $B\ge 0$ denote the parameters of $j$, then it follows that
\begin{align*}
    d_{Y}(f(x_{P}),g(x))&\le A\cdot \big(d_{\mathcal{P}_{Y}}(j(f(x_{P}), j(g(x)))+B\big) \\
    &\le A\cdot \big(d_{\mathcal{P}_{Y}}(j(f(i(P)), \psi(P))+B\big) \\
    &\le A\cdot L + B.
\end{align*}
Thus $d_{Y}(f(x),g(x)) \le C\cdot D +K+A\cdot L+B$, and $g$ is at bounded distance from $f$. This concludes the proof. 
\end{proof}

\subsection{Scaling groups}\label{sec:Scalinggroups}

Let $X$ be an amenable bounded degree graph, and let us define its~\textit{scaling quasi-isometry group} as 
\begin{equation*}
    \text{QI}_{\text{sc}}(X) \defeq \lbrace \text{measure-scaling quasi-isometries $X\rightarrow X$}\rbrace/\text{bounded distance}.
\end{equation*}

Then, it follows from Lemma~\ref{lem:uniquenessscalingfactor} and Proposition~\ref{prop:stabilitypropertiesforscalingQI} that there is a well-defined group morphism 
\begin{align*}
    \text{Sc}\colon  \text{QI}_{\text{sc}}(X)&\longrightarrow \R_{>0} \\
    f\;\text{quasi-$k$-to-one}&\longmapsto k
\end{align*}
taking any scaling quasi-isometry of $X$ to its scaling factor. We call this map the \textit{scale morphism}.

\begin{definition}\label{def:scalinggroup}
Let $X$ be an amenable bounded degree graph. Its~\textit{scaling group}, denoted $\text{Sc}(X)$, is the image of the scale morphism $\text{Sc}\colon  \text{QI}_{\text{sc}}(X)\longrightarrow \R_{>0}$. 
\end{definition}

The scaling group $\text{Sc}(X)$ is a measure-scaling quasi-isometry invariant of the space $X$. In particular, given a finitely generated group $G$, its scaling group $\text{Sc}(G)$ is defined as the scaling group of its Cayley graph with respect to an arbitrary finite generating set, and $\text{Sc}(G)$ does not depend on the choice of such a generating set. 

The object being introduced, we should keep in mind two central questions about the behaviour of quasi-isometries of a given finitely generated group:
\begin{itemize}
    \item Is it true that any self-quasi-isometry $G\rightarrow G$ is measure-scaling?
    \item In general, what is $\text{Sc}(G)$?
\end{itemize}

Let us now give few examples of scaling groups that are known. 

\begin{proposition}[{\cite[Corollary~6.6]{GT22}}]\label{prop:examplesofscalinggroups}
We have $\text{Sc}(G)=\R_{>0}$ if:
\begin{enumerate}[label=(\roman*)]
    \item $G$ is a Carnot group or a lattice in a Carnot group;
    \item $G=\text{SOL}(\R)$ or a lattice in $\text{SOL}(\R)$;
    \item $G=\text{BS}(1,n)$ for all $n\ge 2$.
\end{enumerate}
\end{proposition}

On the other hand, for point~\textit{(i)}, one should keep in mind the fact that $\text{Sc}(\Z^d)=\text{Sc}(\R^d)=\R_{>0}$ for any $d\ge 1$, but a self-quasi-isometry of $\Z^d$ or $\R^d$ need not be scaling. 

The scaling group is also known to be smaller in a number of other situations. However, the proof of these facts are in general much harder, because showing a restriction on the scaling group of a space requires to understand precisely the behaviour of all quasi-isometries of the space. Such information is known only for very specific classes of groups. 

\begin{proposition}[{\cite{Dym10}; \cite{EFW13}}]\label{prop:intermediatescalinggroups}
Let $F$ be a non-trivial finite group. Then one has 
\begin{equation*}
    \text{Sc}(F\wr\Z)=\langle p_{1},\dots,p_{k}\rangle
\end{equation*}
where $p_{1},\dots,p_{k}$ are the prime numbers appearing in the decomposition of $|F|$.
\end{proposition}

Additionally, it follows from results of~\cite{DPT15} that higher rank lamplighters $\Gamma_{d}(q)$ introduced in this article have the same scaling groups, namely 
\begin{equation*}
    \text{Sc}(\Gamma_{d}(q))=\langle p_{1},\dots,p_{k}\rangle
\end{equation*}
where $p_{1},\dots,p_{k}$ are the prime numbers appearing in the decomposition of $q\ge 2$.

As we will explain in Chapter~\ref{chap:chapter2}, it follows from the work of Genevois and Tessera on the quasi-isometric rigidity of lamplighters that $\text{Sc}(F\wr K)=\lbrace 1\rbrace$ if $F$ is a non-trivial finite group and $K$ is amenable, finitely presented and one-ended. We prove in Chapter~\ref{chap:chapter4} a similar result, when $F$ is allowed to be infinite and of polynomial growth. 

\paragraph{Algebraic and geometric consequences.} Let us now outline some nice consequences that can be deduced from the computation of the scaling group of a finitely generated group. 

The first one is algebraic and concerns finite-index subgroups of a given finitely generated group.

\begin{corollary}\label{cor:nopropersubgroupisotothegroup}
Let $G$ be a finitely generated group. If $\text{Sc}(G)=\lbrace 1\rbrace$, then $G$ has no proper finite-index subgroups isomorphic to itself. 
\end{corollary}

\begin{proof}
If $H$ is a finite-index subgroup isomorphic to $G$, composing this isomorphism with the natural inclusion $H\hookrightarrow G$ provides a quasi-$\frac{1}{[G:H]}$-to-one self-quasi-isometry $G\rightarrow G$ (using Proposition~\ref{prop:stabilitypropertiesforscalingQI}\textit{(ii)}), and thus $[G:H]=1$ since $\text{Sc}(G)=\lbrace 1\rbrace$, i.e. $H=G$ cannot be proper. 
\end{proof}

More generally, the same proof yields:

\begin{corollary}\label{cor:biLipsubgroupshavesameindex}
Let $G$ be a finitely generated group. If $\text{Sc}(G)=\lbrace 1\rbrace$, and if $H_{1}, H_{2}$ are two biLipschitz equivalent finite-index subgroups of $G$, then $[G:H_{1}]=[G:H_{2}]$.
\end{corollary}

Even though these two facts are elementary, they illustrate well the interaction between the algebraic structure of a group with its large-scale geometry. From a complete description of its quasi-isometries, one can deduce non-trivial algebraic facts about subgroups of a given group. 

The second consequence is geometric, and deals with the difference between quasi-isometries and biLipschitz equivalences. While finite-index extensions of a finitely generated group yield groups that are quasi-isometric, one can wonder whether these extensions are in fact biLipschitz equivalent to the initial group. In view of what has been said above in the non-amenable setting, the question is only relevant for amenable groups.

\begin{corollary}\label{cor:BiLip.eqandSc}
Let $G$ be a finitely generated amenable group, and let $H$ be a group having $G$ as a finite-index subgroup. Then $G$ and $H$ are biLipschitz equivalent if and only if $[H:G]\in \text{Sc}(G)$.
\end{corollary}

\begin{proof}
Suppose that there is a biLipschitz equivalence $H\rightarrow G$. Such a map is quasi-one-to-one, and pre-composing it with the natural inclusion $G\hookrightarrow H$ provides a quasi-$\frac{1}{[H:G]}$-to-one quasi-isometry $G\rightarrow G$, by Proposition~\ref{prop:stabilitypropertiesforscalingQI}. Thus $\frac{1}{[H:G]}\in \text{Sc}(G)$, whence $[H:G]\in \text{Sc}(G)$.

\noindent Conversely, assume that $[H:G]\in \text{Sc}(G)$, and fix a quasi-$[H:G]$-to-one map $f\colon G\rightarrow G$. Post-composing it with the natural inclusion $G\hookrightarrow H$ gives a quasi-one-to-one map $G\rightarrow H$, and such a map lies at bounded distance from a bijection by Theorem~\ref{thm:Whytethm}. The latter is the desired biLipschitz equivalence. 
\end{proof}

\subsection{Open questions} As already pointed out, the systematic study of scaling quasi-isometries of amenable groups has been initiated recently, in~\cite{GT22}, and many challenging questions arise. We propose here a selection of such questions. 

\smallskip 

The scaling group is invariant under measure-scaling quasi-isometries. Hence:

\begin{question}
Are there pairs of quasi-isometric amenable groups with different scaling groups? 
\end{question}

In the light of Proposition~\ref{prop:examplesofscalinggroups}\textit{(i)}, it is tempting to conjecture that the scaling group of a group $G$ tends to be all of $\R_{>0}$ if $G$ is “small” in some sense. We then ask: 

\begin{question}
Is it true that $\text{Sc}(G)=\R_{>0}$ for any finitely generated group $G$ with polynomial growth? 
\end{question}

Lastly, there is a lack of known subgroups of $\R_{>0}$ that can be obtained as the scaling group of an amenable finitely generated group. Apart from lamplighters over $\Z$ (Proposition~\ref{prop:intermediatescalinggroups}) and higher rank lamplighters from~\cite{DPT15}, there is no other instance of a proper subgroup of $\R_{>0}$ being the scaling group of an amenable group $G$.

\begin{question}
Is it true that any finitely generated subgroup of $\R_{>0}$ can be realised as the scaling group of an amenable group?
\end{question}

In~\cite{Lev24}, it is proved that any finitely generated subgroup of $\R_{>0}$ can be realised as the scaling group of a bounded degree graph. However, the graphs considered are far from being Cayley graphs of amenable finitely generated groups. 

\section{Quasi-isometries of pairs}\label{sec:QIofpairs}

In our study of quasi-isometries of permutational lamplighters below, we will be interested in quasi-isometries that moreover quasi-preserve some subspaces of their source space. This is the next refinement we introduce. 

\begin{definition}\label{def:QIofpairs}
Let $X,Y$ be metric spaces and let $\mathcal{A}$ (resp. $\mathcal{B}$) be a collection of subspaces of $X$ (resp. $Y$). A quasi-isometry $f\colon X\rightarrow Y$ is a~\textit{quasi-isometry of pairs} $f\colon (X,\mathcal{A})\longrightarrow (Y,\mathcal{B})$ if there exists $Q>0$ such that:
\begin{itemize}
    \item For any $A\in\mathcal{A}$, there exists $B\in\mathcal{B}$ such that $d_{\text{Haus}}(f(A),B)\le Q$.
    \item For any $B\in\mathcal{B}$, there exists $A\in\mathcal{A}$ such that $d_{\text{Haus}}(f(A),B)\le Q$.
\end{itemize}
In this case, if $f$ is a $(C,K)-$quasi-isometry, then we say that $f$ is a $(C,K,Q)-$\textit{quasi-isometry of pairs}, and that $(X,\mathcal{A}), (Y,\mathcal{B})$ are~\textit{quasi-isometric pairs}. 
\end{definition}

In the particular case where $X=G$ and $Y=H$ are finitely generated groups and $\mathcal{A}=\mathcal{C}_{M}$, $\mathcal{B}=\mathcal{C}_{N}$ are collections of left cosets of subgroups $M\leqslant G$, $N\leqslant H$, we simply write $(G,M)\longrightarrow (H,N)$ to denote a quasi-isometry of pairs $(G,\mathcal{C}_{M})\longrightarrow (H,\mathcal{C}_{N})$.

\smallskip

In this part, we record basic properties of quasi-isometries of pairs that will be useful below, and we specify then to quasi-isometries of pairs between finitely generated groups quasi-preserving cosets of a fixed normal subgroup. For a more general approach on the notion and many geometric properties with respect to such quasi-isometries, we refer to~\cite{HM23, AM24} and the references therein.  

\smallskip

We start by gathering in a lemma some elementary observations for future use.

\begin{lemma}\label{lem:stabilitypropertiesQIofpairs}
Let $X,Y,Z$ be three metric spaces with collections of subspaces $\mathcal{A}, \mathcal{B}, \mathcal{D}$ respectively. 
\begin{enumerate}[label=(\roman*)]
    \item If \;$f\colon (X,\mathcal{A})\longrightarrow (Y,\mathcal{B})$ is a quasi-isometry of pairs and if $h\colon X\rightarrow Y$ lies at finite distance from $f$, then $h\colon (X,\mathcal{A})\longrightarrow (Y,\mathcal{B})$ is a quasi-isometry of pairs.
    \item If \;$f\colon (X, \mathcal{A})\longrightarrow (Y,\mathcal{B})$, $g\colon (Y,\mathcal{B})\longrightarrow (Z,\mathcal{D})$ are quasi-isometries of pairs, then $g\circ f\colon (X,\mathcal{A})\longrightarrow (Z,\mathcal{D})$ is a quasi-isometry of pairs. 
    \item If \;$f\colon (X,\mathcal{A})\longrightarrow (Y,\mathcal{B})$ is a quasi-isometry of pairs, then any of its quasi-inverses $f'\colon (Y,\mathcal{B})\longrightarrow (X,\mathcal{A})$ is a quasi-isometry of pairs.  
\end{enumerate}
\end{lemma}

\begin{proof}
\textit{(i)} Let $f\colon (X,\mathcal{A})\longrightarrow (Y,\mathcal{B})$ be a $(C,K,Q)-$quasi-isometry of pairs, and suppose that $h\colon X\rightarrow Y$ lies at distance $\le L$ from $f$. Then $h$ is a $(C,K+2L)-$quasi-isometry. Now, let $A\in\mathcal{A}$. By assumption, we find $B\in\mathcal{B}$ such that $d_{\text{Haus}}(f(A),B)\le Q$, so it follows that
\begin{equation*}
    d_{\text{Haus}}(h(A), B) \le d_{\text{Haus}}(h(A), f(A))+d_{\text{Haus}}(f(A),B)\le Q+L.
\end{equation*}
Similarly, given any $B\in\mathcal{B}$ there exists $A\in \mathcal{A}$ with $d_{\text{Haus}}(h(A),B)\le Q+L$, thus $h\colon (X,\mathcal{A})\longrightarrow (Y,\mathcal{B})$ is a $(C,K+2L, Q+L)-$quasi-isometry of pairs. 

\noindent\textit{(ii)} Let $f$ be a $(C,K,Q)-$quasi-isometry of pairs, and let $g$ be a $(C',K',Q')-$quasi-isometry of pairs. Then a computation shows that $g\circ f$ is a $(CC',C'K+2K')-$quasi-isometry. Let $A\in \mathcal{A}$. By assumption, we may find some $B\in\mathcal{B}$ such that $d_{\text{Haus}}(f(A),B)\le Q$, and we may find some $D\in\mathcal{D}$ such that $d_{\text{Haus}}(g(B),D)\le Q'$. Thus it follows 
\begin{align*}
    d_{\text{Haus}}(g(f(A)), D) &\le d_{\text{Haus}}\big(g(f(A)), g(B)\big)+d_{\text{Haus}}(g(B),D) \\ 
    &\le C'\cdot d_{\text{Haus}}(f(A),B)+K'+Q' \\
    &\le C'\cdot Q+K'+Q'.
\end{align*}
Conversely, given any $D\in\mathcal{D}$, we pick $B\in\mathcal{B}$ such that $d_{\text{Haus}}(g(B),D)\le Q'$, and we pick $A\in \mathcal{A}$ such that $d_{\text{Haus}}(f(A),B)\le Q$. We then obtain
\begin{align*}
    d_{\text{Haus}}(g(f(A)), D) &\le d_{\text{Haus}}\big(g(f(A)), g(B)\big)+d_{\text{Haus}}(g(B),D) \\ 
    &\le C'\cdot d_{\text{Haus}}(f(A),B)+K'+Q' \\
    &\le C'\cdot Q+K'+Q'.
\end{align*}
We conclude that $g\circ f$ is an $(CC', C'K+2K',C'Q+K'+Q')-$quasi-isometry of pairs from $(X,\mathcal{A})$ to $(Z,\mathcal{D})$. 

\noindent\textit{(iii)} Assume that  $f\colon (X,\mathcal{A})\longrightarrow (Y,\mathcal{B})$ is a $(C,K,Q)-$quasi-isometry of pairs and let $f'\colon Y\rightarrow X$ be a quasi-inverse. Up to increasing $C$ and $K$, we assume that $f'$ is also a $(C,K)-$quasi-isometry. Fix now $B\in\mathcal{B}$. By assumption, we find $A\in \mathcal{A}$ such that
\begin{equation*}
    d_{\text{Haus}}(f(A),B)\le Q.
\end{equation*}
This implies that $d_{\text{Haus}}\big(f'(f(A)),f'(B)\big)\le C\cdot Q+K$, so that
\begin{equation*}
    d_{\text{Haus}}(f'(B),A)\le d_{\text{Haus}}\big(f'(B), f'(f(A))\big)+d_{\text{Haus}}(f'(f(A)), A) \le C\cdot Q+2K
\end{equation*}
and similarly one shows that given any $A\in\mathcal{A}$ there exists $B\in\mathcal{B}$ such that the Hausdorff distance between $f'(B)$ and $A$ is at most $C\cdot Q+2K$. Thus $f'\colon (Y,\mathcal{B})\longrightarrow (X,\mathcal{A})$ is a $(C,K,C\cdot Q+2K)-$quasi-isometry of pairs. The proof is complete.
\end{proof}

We now focus on finitely generated groups with collections given by cosets of normal subgroups.

\begin{proposition}\label{prop:inducedmapbetweenquotients}
Let $G$ and $H$ be finitely generated groups with normal subgroups $M\lhd G$, $N\lhd H$. A quasi-isometry of pairs $f\colon (G,M)\longrightarrow (H,N)$ induces a coarsely defined quasi-isometry $\overline{f}\colon G/M\rightarrow H/N$. Moreover, for any $p\in G$, $f(pM)$ is at a uniform finite Hausdorff distance from $\overline{f}(pM)$.
\end{proposition}

\begin{proof}
Fix a finite generating set $T$ (resp. $S$) of $G$ (resp. of $H$), and let $\pi_{G}$ (resp. $\pi_{H}$) denote the canonical projection of $G$ (resp. of $H$) onto $G/M$ (resp. $H/N$). We equip $G/M$ (resp. $H/N$) with the word metric provided by the generating set $\pi_{G}(T)$ (resp. $\pi_{H}(S)$).
Let $C\ge 1$, $K,Q\ge 0$ be constants such that $f$ and a quasi-inverse $g\colon (H,N)\longrightarrow (G,M)$ are $(C,K,Q)-$quasi-isometry of pairs. 

\noindent We define $\overline{f}$ as follows: if $pM\in G/M$, then this coset is mapped by $f$ at distance $\le Q$ from some $N-$coset $y_{p}N$. We then let $\overline{f}(pM) \defeq y_{p}N$. By construction, the value of $\overline{f}$ on a given coset does not depend on the choice of a representative of this specific coset (note however that the choice of $y_{p}$ is not unique). We now check that $\overline{f}$ is a quasi-isometry.

\noindent Let $a=pM$, $b=p'M$ be two adjacent vertices in $\text{Cay}(G/M, \pi_{G}(T))$. This means that $b=a\pi_{G}(t)=(pM)(tM)=ptM$ for some $t\in T$, so that 
\begin{align*}
    d_{H/N}\big(\overline{f}(a), \overline{f}(b)\big)&=d_{H/N}(y_{p}N, y_{pt}N) \\
    &\le d_{H/N}(y_{p}N, f(p)N)+d_{H/N}(f(p)N,f(pt)N)+d_{H/N}(f(pt)N, y_{pt}N).
\end{align*}
In this expression, the second term is bounded from above by $C+K$, since
\begin{equation*}
    d_{H/N}\big(f(p)N,f(pt)N\big) \le d_{H}(f(p),f(pt))\le C\cdot d_{G}(p,pt)+K=C+K.
\end{equation*}
For the first term, we know by assumption that $d_{\text{Haus}}(f(pM),y_{p}N)\le Q$, which implies in particular that $f(p)\in (y_{p}N)^{+Q}$, so $f(p)$ is connected to a point in $y_{p}N$ by a path in $\text{Cay}(H,S)$ of length $\le Q$. The projection of such a path to $H/N$ is a path in the Cayley graph $\text{Cay}(H/N,\pi_{H}(S))$ of length $\le Q$ connecting $f(p)N$ to $y_{p}N$, so that 
\begin{equation*}
    d_{H/N}(y_{p}N, f(p)N) \le Q. 
\end{equation*}
We show similarly that $d_{H/N}(f(pt)N, y_{pt}N) \le Q$, whence 
\begin{equation*}
    d_{H/N}\big(\overline{f}(a), \overline{f}(b)\big) \le C+K+2Q.
\end{equation*}
We conclude from Lemma~\ref{lem:Lipschitzbetweengraphs} that $\overline{f}\colon G/M\rightarrow H/N$ is $(C+K+2Q)-$Lipschitz. 

\noindent Now, consider the map $\overline{g}\colon H/N\rightarrow G/M$ induced by $g$, defined by $\overline{g}(qN) \defeq z_{q}M$, where $z_{q}M$ is an $M-$coset at Hausdorff distance at most $Q$ from $g(qN)$. We can reproduce the above argument to prove that $\overline{g}$ is also $(C+K+2Q)-$Lipschitz. Now, if $p\in G$, note that
\begin{align*}
    d_{\text{Haus}}(z_{y_{p}}M, pM) &\le d_{\text{Haus}}(z_{y_{p}}M, g(y_{p}N))+d_{\text{Haus}}\big(g(y_{p}N), g(f(pM))\big)\\
    &+d_{\text{Haus}}\big(g(f(pM)), pM\big) \\
    &\le Q+C\cdot d_{\text{Haus}}(y_{p}N, f(pM))+K+K \\
    &\le (C+1)\cdot Q+2K
\end{align*}
since $g\circ f$ lies at distance $\le K$ from $\text{Id}_{G}$ and since $g$ is a $(C,K)-$quasi-isometry. In particular, we get that $p\in (z_{y_{p}}M)^{+((C+1)Q+2K)}$, so there is a path in $\text{Cay}(G,T)$ of length $\le (C+1)Q+2K$ connecting $p$ to a point of $z_{y_{p}}M$. The projection of such a path to $G/M$ provides a path in $\text{Cay}(G/M,\pi_{G}(T))$ of length $\le (C+1)Q+2K$ connecting $pM$ to $z_{y_{p}}M=\overline{g}(\overline{f}(pM))$, so that 
\begin{equation*}
    d_{G/M}\left(\overline{g}(\overline{f}(pM)), pM\right)\le (C+1)\cdot Q+2K. 
\end{equation*}
We conclude that $\overline{g}\circ\overline{f}$ lies at distance at most $(C+1)\cdot Q+2K$ from $\text{Id}_{G/M}$, and similarly $\overline{f}\circ\overline{g}$ lies at distance at most $(C+1)\cdot Q+2K$ from $\text{Id}_{H/N}$. We can now conclude that 
\begin{align*}
    d_{H/N}\big(\overline{f}(a),\overline{f}(b)\big) &\ge \frac{1}{C+K+2Q}d_{G/M}\big(\overline{g}(\overline{f}(a)), \overline{g}(\overline{f}(b))\big) \\
    &\ge \frac{1}{C+K+2Q}\cdot d_{G/M}(a,b)-\frac{2((C+1)Q+2K)}{C+K+2Q}
\end{align*}
for any $a,b\in G/M$, so $\overline{f}$ is a quasi-isometry with $\overline{g}$ as a quasi-inverse. 
\end{proof}

\begin{remark}
We should emphasize here that, given a $(C,K,Q)-$quasi-isometry of pairs $f\colon (G,M)\longrightarrow (H,N)$, the map $\overline{f}$ given by this statement is not unique, as given an $M-$coset $pM$ in $G$, there could be many $N-$cosets in $H$ lying at Hausdorff distance at most $Q$ from $f(pM)$. However, any two such cosets would be at Hausdorff distance $\le 2Q$, so that any two such maps induced by $f$ always lie at distance at most $2Q$ from each other. 
\end{remark}

We record now some elementary properties of the correspondence $f\longmapsto \overline{f}$. 

\begin{proposition}\label{prop:compatibilityoperationsquotients}
Let $G,H,I$ be finitely generated groups with normal subgroups $M,N,J$. Let $f\colon (G,M)\longrightarrow (H,N)$ be a quasi-isometry of pairs, and fix an induced map $\overline{f}\colon G/M\rightarrow H/N$. The following properties hold. 
\begin{enumerate}[label=(\roman*)]
    \item If \;$h\colon G\rightarrow H$ is a map at bounded distance from $f$, then any quasi-isometry $\overline{h}$ induced by $h$ lies at bounded distance from $\overline{f}$, and there is at least one choice of $\overline{h}$ that coincides with $\overline{f}$.
    \item If $f'\colon H\rightarrow G$ is a quasi-inverse of $f$, then $\overline{f'}$ is a quasi-inverse of $\overline{f}$. 
    \item If $g\colon (H,N)\longrightarrow (I,J)$ is another quasi-isometry of pairs with an induced quasi-isometry $\overline{g}\colon H/N\longrightarrow I/J$, then $\overline{g\circ f}$ lies at bounded distance from $\overline{g}\circ\overline{f}$, and there is at least one choice of $\overline{g\circ f}$ that coincides with $\overline{g}\circ \overline{f}$.  
\end{enumerate}
\end{proposition}

\begin{proof}
\textit{(i)} Let $f$ be a $(C,K,Q)-$quasi-isometry of pairs. By Lemma~\ref{lem:stabilitypropertiesQIofpairs}\textit{(i)} and its proof, $h$ is then a $(C,K+2L, Q+L)-$quasi-isometry of pairs, where $L\ge 0$ is a constant that controls the distance between $h$ and $f$. Moreover, the proof of Lemma~\ref{lem:stabilitypropertiesQIofpairs}\textit{(i)} shows that, if $pM$ is an $M-$coset, then $h(pM)$ lies at Hausdorff distance at most $Q+L$ from the $N-$coset $y_{p}N$ lying at Hausdorff distance at most $Q$ from $f(pM)$, given by the assumption that $f$ is a quasi-isometry of pairs. Since this is true for any $M-$coset $pM$ of $G$, we indeed have that $h$ induces a quasi-isometry $\overline{h}$ equal to $\overline{f}$. 

\noindent\textit{(ii)} has been proved in Proposition~\ref{prop:inducedmapbetweenquotients}, and \textit{(iii)} is proved similarly to \textit{(i)}. 
\end{proof}

\begin{remark}\label{rm2.12}
Let $f\colon (G,M)\longrightarrow (H,N)$ be a $(C,K,Q)-$quasi-isometry of pairs. By assumption, for any $p\in G$, $f$ sends a coset $pM\subset G$ at Hausdorff distance $\le Q$ from a coset $y_{p}N$, used to define $\overline{f}\colon G/M\rightarrow H/N$. We can thus modify $f$ to get a map $f'\colon G\rightarrow H$ at distance $\le Q$ from $f$, sending the coset $pM$ into the coset $y_{p}N$. This in turn implies that $y_{p}N=f'(p)N$ for any $p\in G$, which means that $f'$ has an induced quasi-isometry $G/M\longrightarrow H/N$ with a simpler formula: $\overline{f'}(pM) \defeq f'(p)N$. Note that this induced quasi-isometry agrees with $\overline{f}$. Hence, in the sequel, when we will be given a quasi-isometry of pairs $f\colon (G,M)\longrightarrow (H,N)$ with an induced quasi-isometry $\overline{f}\colon G/M\rightarrow H/N$ satisfying a property $\mathcal{P}$ (for instance being quasi-$k$-to-one for some $k>0$), we will always be able to assume (up to finite distance) that $f$ sends cosets into cosets and that the new induced quasi-isometry we get after this change still has property $\mathcal{P}$.
\end{remark}

\chapter{Quasi-isometric rigidity of lamplighters and halo products}\label{chap:chapter2}

The main objects of this thesis being defined, we now present the main questions we are concerned with and the important results on the topic already recorded in the literature. All these results constitute the basis on which we present our contributions in the second part. 

\vspace{0.3cm}

\minitoc

\section{Generalities}

The wreath product of two finitely generated groups being finitely generated itself (cf. Proposition~\ref{prop:wreathproductspreservefinitegeneration}), it naturally inherits a metric structure. From the point of view of large-scale geometry, it is therefore natural to ask:

\begin{question}\label{question:wreathproductsQI}
Let $G_{1},G_{2},H_{1},H_{2}$ be finitely generated groups. When are $G_{1}\wr H_{1}$ and $G_{2}\wr H_{2}$ quasi-isometric?
\end{question}

This question can be seen as part of the ambitious program initiated by Gromov in the 1980s~\cite{Gro81}, aiming at classifying, or at least “understanding”, finitely generated groups up to quasi-isometry. 

Question~\ref{question:wreathproductsQI} is still open in full generality, but much progress has been made on the understanding of quasi-isometry classes of wreath products. An important first piece of work was carried out by Dyubina in 2000, showing that the wreath product construction is compatible with biLipschitz equivalences (i.e. bijective quasi-isometries).

\begin{proposition}[{\cite[Lemma~1]{Dyu00}}]\label{prop:BiLipwreathproductscolor}
Let $A$, $B$ and $C$ be finitely generated groups. If $A$ and $B$ are biLipschitz equivalent, then $A\wr C$ and $B\wr C$ are biLipschitz equivalent. 
\end{proposition}

A similar statement holds if the base groups are biLipschitz equivalent.

\begin{proposition}\label{prop:BiLipwreathproductsbase}
Let $A$, $B$ and $C$ be finitely generated groups. If $B$ and $C$ are biLipschitz equivalent, then $A\wr B$ and $A\wr C$ are biLipschitz equivalent.
\end{proposition}

\begin{proof}
Fix finite generating sets $S_{A}$, $S_{B}$, $S_{C}$ of $A$, $B$ and $C$ respectively. Let $f\colon B\rightarrow C$ be a biLipschitz equivalence. This means that $f$ is bijective, and that $f$ and its inverse $f^{-1}\colon C\rightarrow B$ are $L-$Lipschitz for some $L\ge 1$. Then define the map
\begin{align*}
    \varphi\colon A\wr B &\longrightarrow A\wr C \\
    (c,p)&\longmapsto (c\circ f^{-1}, f(p))
\end{align*}
whose inverse is 
\begin{align*}
    \psi \colon A\wr C&\longrightarrow A\wr B \\
    (d,q)&\longmapsto (d\circ f, f^{-1}(q)).
\end{align*}
We claim that $\varphi$ is Lipschitz. Indeed, fix two adjacent vertices $a=(c,p)$, $b=(c',p')$ in $\text{Cay}(A\wr B, S_{A\wr B})$ (here $S_{A\wr B}$ denotes the generating set of $A\wr B$ exhibited in Example~\ref{ex:descriptionoflamplighters}). There are two cases to consider. 
\begin{itemize}
    \item If $c'=c$ and $p'$ is adjacent to $p$ in $\text{Cay}(B,S_{B})$, then $p'=ps$ for some $s\in S_{B}$, and thus one has 
    \begin{align*}
        d_{A\wr C}(\varphi(a), \varphi(b))&=d_{A\wr C}\big((c\circ f^{-1}, f(p)), (c\circ f^{-1}, f(ps)\big) \\
        &=d_{C}(f(p), f(ps)) \\
        &\le L\cdot d_{B}(p,ps) \\
        &=L.
    \end{align*}
    \item If rather $p=p'$ and $c$, $c'$ differ only on this vertex, with $c(p)$ and $c'(p)$ being adjacent in $\text{Cay}(A,S_{A})$, then $c\circ f^{-1}$, $c'\circ f^{-1}$ only differ on $f(p)$ and take adjacent colors on this vertex, so that $\varphi(a)=(c\circ f^{-1}, f(p))$, $\varphi(b)=(c'\circ f^{-1}, f(p))$ are neighbours in $\text{Cay}(A\wr C, S_{A\wr C})$, i.e.
    \begin{equation*}
        d_{A\wr C}(\varphi(a), \varphi(b))=1\le L.
    \end{equation*}
\end{itemize}
In any case, $\varphi$ sends adjacent vertices to vertices at distance $\le L$ in $\text{Cay}(A\wr C, S_{A\wr C})$. Lemma~\ref{lem:Lipschitzbetweengraphs} implies then that $\varphi$ is $L-$Lipschitz. A similar computation, using that $f^{-1}$ is $L-$Lipschitz, implies that $\psi$ is also $L-$Lipschitz. This proves that $\varphi$ is a biLipschitz equivalence between $A\wr B$ and $A\wr C$. 
\end{proof}

Proposition~\ref{prop:BiLipwreathproductscolor} has several consequences. The first one is that the quasi-isometry class of a lamplighter $A\wr B$ (i.e. $A$ is finite) only depends on $|A|$. For instance, $\Z/6\Z\wr\Z$ and $S_{3}\wr \Z$ are quasi-isometric. But, most notably, Dyubina used this stability property in~\cite{Dyu00} to deduce that some algebraic properties of groups, such as being virtually solvable or virtually torsion-free, are not stable under quasi-isometry. Indeed, for virtual solvability, $\Z/60\Z\wr\Z$ and $A_{5}\wr\Z$, where $A_{5}$ is the alternating group over five elements, are quasi-isometric, but the first one is solvable, while the second one is not virtually solvable. 

For virtual torsion-freeness, note that $\Z\wr\Z$ and $(\Z\times A_{5})\wr\Z$ are biLipschitz equivalent, since $\Z$ and $\Z\times A_{5}$ are biLipschitz equivalent (use for instance Corollary~\ref{cor:BiLip.eqandSc} and $\text{Sc}(\Z)=\R_{>0}$). But $\Z\wr\Z$ is torsion-free, while any finite index subgroup of $(\Z\times A_{5})\wr\Z$ has elements of order $|A_{5}|$. 

However, all these examples are infinitely presented groups, and it is then also natural to wonder whether being virtually solvable or being virtually torsion-free is invariant under quasi-isometry among finitely presented groups. For virtual torsion-freeness, the answer remains negative, and examples have been provided by Cornulier in~\cite[Proposition~2.12]{Cor06}, using permutational wreath products, and by Le Boudec~\cite[Corollary~6]{LB25}. The question for virtual solvability is still open.

In another direction, Dymarz also used lamplighters over $\Z$ in~\cite{Dym10} to prove that, among amenable groups, two groups can be quasi-isometric without being biLipschitz equivalent. 

\begin{proposition}[{\cite[Theorem~1]{Dym10}}]\label{prop:BiLipvsQIDymarz}
Let $F$ be a non-trivial finite group, and let $k\ge 2$ be an integer. Then $F\wr\Z$ and $F^{k}\wr\Z$ are not biLipschitz equivalent if $k$ is not a product of the prime numbers appearing in the decomposition of $|F|$. 
\end{proposition}

On the other hand, it is not hard to construct explicitly a quasi-isometry between $F\wr\Z$ and $F^{k}\wr \Z$ for any choice of $k\ge 2$. See for instance~\cite[Proposition~3.12]{GT24b}. 

Here also, the groups involved are infinitely presented, and a few years later, building on these examples, Dymarz, Peng and Taback introduced in~\cite{DPT15} a family of finitely presented amenable groups, called~\textit{higher rank lamplighters}, to prove that there are also finitely presented groups that can be quasi-isometric without being biLipschitz equivalent. 

\paragraph{Lamplighters over two-ended groups.} The first classification result that provides a partial answer to Question~\ref{question:wreathproductsQI} goes back to 2013 and is due to Eskin, Fisher and Whyte.

\begin{theorem}[{\cite[Theorem~1.2]{EFW13}}]\label{thm:Eskin-Fisher-Whyte}
Let $F_{1}$, $F_{2}$ be non-trivial finite groups. Then $F_{1}\wr\Z$ and $F_{2}\wr\Z$ are quasi-isometric if and only if there exist integers $a,r,s\ge 1$ such that $|F_{1}|=a^{r}$ and $|F_{2}|=a^{s}$. 
\end{theorem}

From this result, one deduces that, if $H_{1}$, $H_{2}$ are both two-ended, then $F_{1}\wr H_{1}$ and $F_{2}\wr H_{2}$ are quasi-isometric if and only if there exist $a,r,s\ge 1$ such that $|F_{1}|=a^{r}$ and $|F_{2}|=a^{s}$. Indeed, if $H_{1}$ is two-ended, it is virtually cyclic, and thus biLipschitz equivalent to $\Z$. The wreath product $F_{1}\wr H_{1}$ is therefore biLipschitz equivalent to $F_{1}\wr\Z$ using Proposition~\ref{prop:BiLipwreathproductsbase}, and similarly $F_{2}\wr H_{2}$ is biLipschitz equivalent to $F_{2}\wr\Z$. The existence of a quasi-isometry $F_{1}\wr H_{1}\longrightarrow F_{2}\wr H_{2}$ is thus equivalent to the existence of a quasi-isometry $F_{1}\wr\Z\longrightarrow F_{2}\wr\Z$, hence the claim by Theorem~\ref{thm:Eskin-Fisher-Whyte}.

\section{Lamplighters over one-ended groups} After Eskin-Fisher-Whyte's breakthrough, a second classification result has been obtained by Genevois and Tessera in~\cite{GT24b}.

\begin{theorem}[{\cite[Corollary~1.6]{GT24b}}]\label{thm:classificationGT21}
Let $F_{1}$ and $F_{2}$ be non-trivial finite groups. Let $H_{1}$ and $H_{2}$ be finitely presented groups. Suppose that $H_{1}$ is one-ended.
\begin{enumerate}[label=(\roman*)]
    \item If $H_{1}$ is not amenable, then $F_{1}\wr H_{1}$ and $F_{2}\wr H_{2}$ are quasi-isometric if and only if $|F_{1}|$, $|F_{2}|$ have the same prime divisors and there exists a quasi-isometry $H_{1}\rightarrow H_{2}$.
    \item If $H_{1}$ is amenable, then $F_{1}\wr H_{1}$ and $F_{2}\wr H_{2}$ are quasi-isometric if and only if there exist integers $a,r,s\ge 1$ such that $|F_{1}|=a^{r}$, $|F_{2}|=a^{s}$ and there exists a quasi-$\frac{s}{r}$-to-one quasi-isometry $H_{1}\rightarrow H_{2}$.
\end{enumerate}
\end{theorem}

Several comments are in order. First, the classification of wreath products is actually a corollary of the main result of~\cite{GT24b}, namely~\cite[Theorem~1.5]{GT24b}, which is a classification result for more general lamplighter graphs. The corollary follows directly from this classification since, as mentioned above, Cayley graphs of lamplighter groups coincide with lamplighter graphs over Cayley graphs of groups. We refer the reader to~\cite[Theorem~1.5]{GT24b} for the more general statement and here we focus on the case of wreath products of groups. 

Second, observe that the amenable case is more rigid, as it involves a compatibility condition between the quasi-isometry between the base groups and the cardinality of the lamp groups. Such a condition does not appear in the non-amenable case, but it is in fact hidden: as said earlier, between non-amenable spaces, a quasi-isometry is quasi-$k$-to-one for any real number $k>0$ (cf. Remark~\ref{rem:ScalingQINAspaces}). 

Before turning to the strategy of proof, let us highlight several nice consequences. 

\paragraph{Consequences of Theorem~\ref{thm:classificationGT21}.} Theorem~\ref{thm:classificationGT21} can be reformulated using the scaling group in the case where base groups are the same.

\begin{corollary}[{\cite[Corollary~1.9]{GT24b}}]\label{cor:samebasegroups}
Let $F_{1}$ and $F_{2}$ be non-trivial finite groups. Let $H$ be a finitely presented amenable one-ended group. Then $F_{1}\wr H$ and $F_{2}\wr H$ are quasi-isometric if and only if there exist integers $a,r,s\ge 1$ such that $|F_{1}|=a^{r}$, $|F_{2}|=a^{s}$ and $\frac{s}{r}\in\text{Sc}(H)$.
\end{corollary}

In this case, the classification of lamplighters over $H$ up to quasi-isometry reduces to the understanding of the scaling group of $H$. Such information is available in some cases, for instance:

\begin{corollary}[{\cite[Corollary~1.10]{GT24b}}]\label{cor:lamplightersoverZ^d}
Let $F_{1}$ and $F_{2}$ be non-trivial finite groups. Let $d\ge 2$. Then $F_{1}\wr \Z^d$ and $F_{2}\wr \Z^d$ are quasi-isometric if and only if there exist integers $a,r,s\ge 1$ such that $|F_{1}|=a^{r}$ and $|F_{2}|=a^{s}$.  \end{corollary}

Genevois and Tessera's classification provides additional examples of pairs of quasi-isometric groups that are not biLipschitz equivalent:

\begin{corollary}[{\cite[Corollary~1.15]{GT24b}}]\label{cor:BiLiplamplighters}
Let $F_{1}$ and $F_{2}$ be non-trivial finite groups. Let $H$ be an amenable finitely presented one-ended group. Then $F_{1}\wr H$ and $F_{2}\wr H$ are biLipschitz equivalent if and only if $|F_{1}|=|F_{2}|$.    
\end{corollary}

For instance, $\Z/2\Z\wr\Z^2$ and $\Z/4\Z\wr\Z^2$ are quasi-isometric but not bijectively quasi-isometric. 

\paragraph{Strategy for the proof of Theorem~\ref{thm:classificationGT21}.} The proof of Theorem~\ref{thm:classificationGT21} is decomposed into four main steps. We explain them here since they constitute the basis for the strategy developed in Chapter 3 in order to achieve similar classification results for permutational wreath products. 

\paragraph{The embedding theorem.} The first step is perhaps the hardest one, and reflects a strong rigidity phenomenon: it consists in proving that any quasi-isometry $F_{1}\wr H_{1}\longrightarrow F_{2}\wr H_{2}$ must actually be a quasi-isometry of pairs $(F_{1}\wr H_{1}, H_{1})\longrightarrow (F_{2}\wr H_{2}, H_{2})$, namely it must send $H_{1}-$cosets at bounded Hausdorff distance from $H_{2}-$cosets, and the same property must hold for a quasi-inverse. This fact is deduced from the following~\textit{embedding theorem}:

\begin{theorem}[{\cite[Theorem~1.4]{GT24b}}]\label{thm:embedding1}
Let $G$ be a finitely presented one-ended group. Let $F$ be a non-trivial finite group and let $H$ be an infinite finitely generated group. Then every coarse embedding $\rho\colon G\longrightarrow F\wr H$ has its image contained in the neighbourhood of an $H-$coset. 
\end{theorem}

The idea for proving such a statement is the following. For simplicity, let us focus on $F=\Z_{2}$ and let $\rho\colon G\longrightarrow  \Z_{2}\wr H$ be a coarse embedding. First of all, as 
\begin{equation*}
    \Z_{2}\wr H = \langle a, S : a^2=1, [a, hah^{-1}]=1, \; h\in H\rangle
\end{equation*}
is infinitely presented, one can consider, for finite subsets $S\subset H$, the group given by the~\textit{truncated} presentation 
\begin{equation*}
    \Z_{2}\square_{S} H \defeq \langle a, S : a^2=1, [a, hah^{-1}]=1, \;h\in S\rangle
\end{equation*}
and it is a well-known fact (see e.g.~\cite[Lemma~6.21]{BGT24}) from coarse topology that $\rho$ can be factored through $\Tilde{\rho}\colon G\longrightarrow \Z_{2}\square_{S} H$, such that $\pi_{S}\circ\Tilde{\rho}=\rho$, where $\pi_{S}\colon \Z_{2}\square_{S} H \longrightarrow  \Z_{2}\wr H$ is the natural surjection. The interest of doing this is that, in fact, $\Z_{2}\square_{S} H$ splits algebraically as a semi-direct product
\begin{equation*}
    C(\Gamma)\rtimes H
\end{equation*}
where $\Gamma\defeq \text{Cay}(H,S)$ and where $C(\Gamma)$ is the right-angled Coxeter group over the graph $\Gamma$. Moreover, a computation shows that $\pi_{S}$ sends $H-$cosets to $H-$cosets, so in order to prove that $\rho(G)$ lies in a neighbourhood of an $H-$coset in $F\wr H$, it is enough to prove that $\Tilde{\rho}(G)$ lies in a neighbourhood of an $H-$coset in $\Z_{2}\square_{S} H$.

To prove the latter, we exploit the~\textit{wallspace} structure of $C(\Gamma)$, because the Cayley graph of a right-angled Coxeter group is~\textit{median} (i.e. it is the $1-$skeleton of a $\text{CAT}(0)$ cube complex), and naturally has~\textit{hyperplanes}. Each hyperplane delimits two~\textit{halfspaces}, and the key observation is now that $\Tilde{\rho}(G)$ must necessarily lie in one of these two halfspaces. Indeed, if the hyperplane “cuts” $\Tilde{\rho}(G)$ too heavily, then, projecting back to the lamplighter, $\rho(G)=\pi_{S}(\Tilde{\rho}(G))$ would be separated by $\pi_{S}(\text{hyperplane})$, which is bounded. This would contradict one-endedness of $\rho(G)$. Since this holds for~\textit{any} hyperplane, it follows that $\Tilde{\rho}(G)$ lies in the intersection of all these halfspaces, which consists of a single vertex in the $\text{CAT}(0)$ cube complex. This implies that, in $C(\Gamma)\rtimes H$, $\Tilde{\rho}(G)$ lies in the neighbourhood of an $H-$coset. 

For larger finite groups, Coxeter groups are replaced with graph products of groups, and median geometry is replaced by~\textit{quasi-median geometry}, but the idea is the same. We undertake this strategy below in Chapter~\ref{chap:chapter3} in the case of permutational wreath products with infinite stabilizers. 

Once Theorem~\ref{thm:embedding1} is proved, it remains to observe that, given our quasi-isometry $F_{1}\wr H_{1}\longrightarrow F_{2}\wr H_{2}$, two $H_{1}-$cosets cannot be sent close to the same $H_{2}-$coset, because the coarse intersection of two $H_{1}-$cosets in $F_{1}\wr H_{1}$ is bounded:

\begin{lemma}\label{lem:coarseintersectionsofleavesarebounded}
Let $F$ be a non-trivial finite group, and $H$ be an infinite finitely generated group. Let $R\ge 0$. Then, for any distinct colourings $c,d\in F^{(H)}$, $(cH)\cap (dH)^{+R}$ is bounded.
\end{lemma}

\begin{proof}
Assume, towards a contradiction, that $(cH)\cap (dH)^{+R}$ is infinite. Fix a point $(c,q)\in (cH)\cap (dH)^{+R}$. Then there is a point $(d,p)$ such that $d_{F\wr H}\big((c,q),(d,p)\big) \le R$, which in particular implies that $c$ and $d$ only differ on points at distance $\le R$ from $q$. It follows that $\text{supp}(c^{-1}d)$ is contained in $B_{H}(q,R)$ for infinitely many $q$, which is impossible as it is bounded. 
\end{proof}

\paragraph{Leaf-preserving quasi-isometries.} From the embedding theorem and Lemma~\ref{lem:coarseintersectionsofleavesarebounded}, it follows that any quasi-isometry $\varphi\colon F_{1}\wr H_{1}\longrightarrow F_{2}\wr H_{2}$ is~\textit{leaf-preserving}, in the following sense:

\begin{definition}
Let $F_{1}, F_{2}$ be non-trivial finite groups and let $H_{1}, H_{2}$ be finitely generated groups. Let $\varphi\colon F_{1}\wr H_{1}\longrightarrow F_{2}\wr H_{2}$ be a quasi-isometry and let $\overline{\varphi}\colon F_{2}\wr H_{2}\longrightarrow F_{1}\wr H_{1}$ be a quasi-inverse. We say that $\varphi$ is~\textit{leaf-preserving} if there is a constant $C\ge 0$ such that, for any colouring $c\in F_{1}^{(H_{1})}$ (resp. $d\in F_{2}^{(H_{2})}$), $\varphi(cH_{1})$ (resp. $\overline{\varphi}(dH_{2})$) is at Hausdorff distance at most $C$ from an $H_{2}-$coset (resp. an $H_{1}-$coset).
\end{definition}

This property can be expressed as follows: 

\begin{lemma}[{\cite[Lemma~4.2]{GT24b}}]\label{lem:characterisationofleafpreservingness}
Let $F_{1}, F_{2}$ be non-trivial finite groups and let $H_{1}. H_{2}$ be finitely generated groups. A quasi-isometry $\varphi\colon F_{1}\wr H_{1}\longrightarrow F_{2}\wr H_{2}$ is leaf-preserving if and only if there is a constant $C\ge 0$ and a bijection $\alpha\colon F_{1}^{(H_{1})}\longrightarrow F_{2}^{(H_{2})}$ such that, for any $c\in F_{1}^{(H_{1})}$, the Hausdorff distance between $\varphi(cH_{1})$ and $\alpha(c)H_{2}$ is $\le C$.
\end{lemma}

\begin{proof}
Fix a quasi-inverse $\overline{\varphi}\colon F_{2}\wr H_{2}\longrightarrow F_{1}\wr H_{1}$ of $\varphi$ and constants $A\ge 1$, $B\ge 0$ such that $\varphi$, $\overline{\varphi}$ are $(A,B)-$quasi-isometries. 

\noindent Suppose first that $\varphi$ is leaf-preserving. Then there is a constant $C\ge 0$ such that, for any $c\in F_{1}^{(H_{1})}$, there is $\alpha(c)\in F_{2}^{(H_{2})}$ such that the Hausdorff distance between $\varphi(cH_{1})$ and $\alpha(c)H_{2}$ is $\le C$. Similarly, for any $d\in F_{2}^{(H_{2})}$, there is $\alpha'(d)\in F_{1}^{(H_{1})}$ such that the Hausdorff distance between $\overline{\varphi}(dH_{2})$ and $\alpha'(d)H_{1}$ is $\le C$. Thus one obtains
\begin{align*}
    \overline{\varphi}(\varphi(cH_{1})) &\subset \overline{\varphi}\big((\alpha(c)H_{2})^{+C}\big) \\
    & \subset \overline{\varphi}(\alpha(c)H_{2})^{+(AC+B)} \\
    & \subset \big(\alpha'(\alpha(c))H_{1}\big)^{+(AC+C+B)}
\end{align*}
and since $\overline{\varphi}\circ\varphi$ is at distance $\le B$ from the identity, it follows that
\begin{equation*}
    cH_{1}\subset \overline{\varphi}(\varphi(cH_{1}))^{+B} \subset (\alpha'(\alpha(c))H_{1})^{+((A+1)C+2B)}.
\end{equation*}
Hence, if by contradiction $\alpha'(\alpha(c)) \neq c$, one gets that the intersection 
\begin{equation*}
    cH_{1}\cap (\alpha'(\alpha(c))H_{1})^{+((A+1)C+2B)}
\end{equation*}
is infinite as it contains $cH_{1}$. This contradicts Lemma~\ref{lem:coarseintersectionsofleavesarebounded}, so in fact $\alpha'(\alpha(c))=c$, and $\alpha'\circ \alpha=\text{Id}_{F_{1}^{(H_{1})}}$. The same reasoning shows that $\alpha\circ \alpha'=\text{Id}_{F_{2}^{(H_{2})}}$, hence $\alpha$ and $\alpha'$ are bijections, inverses of each other. 

\noindent Conversely, if $\varphi$ has the desired property, then there is $C\ge 0$ such that, for any $c\in F_{1}^{(H_{1})}$, the Hausdorff distance between $\varphi(cH_{1})$ and $\alpha(c)H_{2}$ is $\le C$, and for any $d\in F_{2}^{(H_{2})}$, the Hausdorff distance between $\overline{\varphi}(dH_{2})$ and $\alpha^{-1}(d)H_{1}$ is $\le C$. In other words, $\varphi$ is leaf-preserving. 
\end{proof}

At this stage, we can therefore say that, combining Lemma~\ref{lem:characterisationofleafpreservingness} and the embedding theorem, any quasi-isometry $\varphi\colon F_{1}\wr H_{1}\longrightarrow F_{2}\wr H_{2}$ can be written, up to a bounded perturbation, as 
\begin{equation*}
    \varphi(c,p)=(\alpha(c),\beta_{c}(p))
\end{equation*}
for some bijection $\alpha\colon F_{1}^{(H_{1})}\longrightarrow F_{2}^{(H_{2})}$ and some collection of maps $\beta_{c}\colon H_{1}\rightarrow H_{2}$, indexed by $c\in F_{1}^{(H_{1})}$. 

\paragraph{Aptolic quasi-isometries.} The key observation, recorded first in~\cite[Theorem~4.3]{GT24b}, is now that, in fact, all $\beta_{c}$ lie at bounded distance from each other, thus showing that $\varphi$ is in fact at bounded distance from an~\textit{aptolic} quasi-isometry. 

\begin{definition}\label{def:aptolicity}
A quasi-isometry $\varphi\colon F_{1}\wr H_{1}\longrightarrow F_{2}\wr H_{2}$ is~\textit{aptolic} if there exist four maps $\alpha\colon F_{1}^{(H_{1})}\longrightarrow F_{2}^{(H_{2})}$, $\alpha'\colon F_{2}^{(H_{2})}\longrightarrow F_{1}^{(H_{1})}$, $\beta\colon H_{1}\rightarrow H_{2}$, $\beta'\colon H_{2}\rightarrow H_{1}$ such that 
\begin{equation*}
    \varphi(c,p)=(\alpha(c),\beta(p))
\end{equation*}
for any $(c,p)\in F_{1}\wr H_{1}$, and 
\begin{equation*}
    \varphi'(c,p)=(\alpha'(c),\beta'(p)), \; (c,p)\in F_{2}\wr H_{2}
\end{equation*}
is a quasi-inverse of $\varphi$. 
\end{definition}

At this point, it is unnecessary to sketch the strategy for proving~\cite[Theorem~4.3]{GT24b}, since a similar result was later proved by Genevois and Tessera in greater generality, that we explain below (see Theorem~\ref{thm:boundeddistancefromaptolicQI}). The proof of the latter is also simpler than that of~\cite[Theorem~4.3]{GT24b}, both technically and conceptually. 

The next step of the strategy is therefore to prove the following three properties that characterise aptolic quasi-isometries:

\begin{proposition}[{\cite[Proposition~3.1]{GT24b}}]\label{prop:characterisationofaptolicQI:finitecase}
Let $F_{1}, F_{2}$ be two non-trivial finite groups, and let $H_{1}, H_{2}$ be infinite finitely generated groups. Let $\alpha\colon F_{1}^{(H_{1})}\longrightarrow F_{2}^{(H_{2})}$, $\beta\colon H_{1}\rightarrow H_{2}$ be two maps. Then the map 
\begin{align*}
    \varphi\colon F_{1}\wr H_{1}&\longrightarrow F_{2}\wr H_{2} \\
    (c,p)&\longmapsto (\alpha(c),\beta(p))
\end{align*}
is a quasi-isometry if and only if:
\begin{enumerate}[label=(\roman*)]
    \item $\alpha$ is a bijection; 
    \item $\beta$ is a quasi-isometry;
    \item There exists a constant $Q\ge 0$ such that, for any colourings $c_{1},c_{2}\in F_{1}^{(H_{1})}$, the Hausdorff distance between $\beta(\text{supp}(c_{1}^{-1}c_{2}))$ and $\text{supp}(\alpha(c_{1})^{-1}\alpha(c_{2}))$ is $\le Q$. 
\end{enumerate}
If so, every quasi-inverse of $\varphi$ is of the form 
\begin{align*}
    \overline{\varphi}\colon F_{2}\wr H_{2}&\longrightarrow F_{1}\wr H_{1} \\
    (c,p)&\longmapsto (\alpha^{-1}(c),\overline{\beta}(p))
\end{align*}
where $\overline{\beta}\colon H_{2}\rightarrow H_{1}$ is a quasi-inverse of $\beta$. 
\end{proposition}

\begin{proof} Fix finite generating sets $S_{1}$, $S_{2}$ of $H_{1}$, $H_{2}$ respectively.

\noindent Fix constants $A\ge 1$, $B\ge 0$ such that $\varphi$ and a quasi-inverse $\overline{\varphi}$, of the form 
\begin{equation*}
    \overline{\varphi}(c,p)=(\overline{\alpha}(c), \overline{\beta}(p)), \; (c,p)\in F_{2}\wr H_{2}
\end{equation*}
are $(A,B)-$quasi-isometries, with $\overline{\varphi}\circ\varphi$, $\varphi\circ\overline{\varphi}$ being at distance $\le B$ from the identities. 

\noindent We prove~\textit{(ii)} first. Let $p,q\in H_{1}$. Then one has 
\begin{align*}
    d(\beta(p), \beta(q))&=d\big((\alpha(\mathbf{1}), \beta(p)), (\alpha(\mathbf{1}), \beta(q)\big) \\
    &=d(\varphi(\mathbf{1}, p), \varphi(\mathbf{1},q)) \\
    &\le A\cdot d((\mathbf{1}, p), (\mathbf{1},q))+B \\
    &=A\cdot d(p,q)+B
\end{align*}
and similarly one proves that $d(\beta(p), \beta(q)) \ge \frac{1}{A}\cdot d(p,q)-B$. Thus $\beta$ is a quasi-isometric embedding. Now, given an arbitrary $h\in H_{2}$, there exists a point $(c,p)\in F_{1}\wr H_{1}$ such that 
\begin{equation*}
    d\big((\alpha(c),\beta(p)), (\mathbf{1}, h)\big)=d(\varphi(c,p), (\mathbf{1}, h))\le B
\end{equation*}
and in particular $d(\beta(p), h)\le B$. Hence $\beta$ is coarsely surjective, and thus it is a quasi-isometry, with the same parameters as $\varphi$.

\noindent We focus now on~\textit{(i)}. Let $(c,p)\in F_{2}\wr H_{2}$. We have already 
\begin{equation*}
    d\big((\alpha\circ \overline{\alpha}(c), \beta\circ\overline{\beta}(p)), (c,p)\big)=d(\varphi\circ\overline{\varphi}(c,p), (c,p)) \le B
\end{equation*}
so $d(\beta\circ\overline{\beta}(p), p)\le B$ and the colourings $c, \alpha\circ \overline{\alpha}(c)$ can only differ on the ball $B_{H_{2}}(p, B)$. It follows from the first observation that $\beta\circ\overline{\beta}$ is at distance $\le B$ from $\text{Id}_{H_{2}}$, and from the second one that $\alpha\circ \overline{\alpha}(c)=c$ by letting $p$ tend to infinity in $H_{2}$. By symmetry, using that $\overline{\varphi}\circ\varphi$ is at distance $\le B$ from the $\text{Id}_{F_{1}\wr H_{1}}$, one gets that $\overline{\beta}\circ\beta$ is at distance $\le B$ from $\text{Id}_{H_{1}}$ and that $\overline{\alpha}\circ\alpha(c)=c$ for any $c\in F_{1}^{(H_{1})}$. In conclusion, $\alpha$ is a bijection with $\overline{\alpha}$ as inverse, $\overline{\beta}$ is a quasi-inverse of $\beta$. This also shows the last claim of the statement. 

\noindent Let us now prove~\textit{(iii)}. Fix two colourings $c_{1},c_{2}\in F_{1}^{(H_{1})}$, and consider a sequence of colourings 
\begin{equation*}
    a_{0}=c_{1}, a_{1},\dots, a_{n-1},a_{n}=c_{2}
\end{equation*}
such that, for any $0\le i\le n-1$, $a_{i}$ and $a_{i+1}$ only differ on a point $p_{i}\in H_{1}$. Then, for $0\le i\le n-1$, one gets 
\begin{align*}
    d\big((\alpha(a_{i}), \beta(p_{i})), (\alpha(a_{i+1}),\beta(p_{i}))\big)&=d(\varphi(a_{i},p_{i}), \varphi(a_{i+1}, p_{i})) \\
    &\le A\cdot d((a_{i}, p_{i}), (a_{i+1},p_{i}))+B \\
    &=A+B 
\end{align*}
so $\alpha(a_{i})$ and $\alpha(a_{i+1})$ can only differ on $B_{H_{2}}(\beta(p_{i}), A+B)$. Thus $\alpha(a_{0})=\alpha(c_{1})$ and $\alpha(a_{n})=\alpha(c_{2})$ can only differ on 
\begin{equation*}
    \bigcup_{i=0}^{n-1}B_{H_{2}}(\beta(p_{i}), A+B)=\big\lbrace \beta(p_{0}),\beta(p_{1}),\dots,\beta(p_{n-1}) \big\rbrace^{+(A+B)} = \beta\big(\text{supp}\big(c_{1}^{-1}c_{2})\big)^{+(A+B)},
\end{equation*}
in other words $\text{supp}\big(\alpha(c_{1})^{-1}\alpha(c_{2})\big) \subset \beta\big(\text{supp}(c_{1}^{-1}c_{2})\big)^{+(A+B)}$. The same reasoning with $\overline{\varphi}$ shows that, for any $d_{1},d_{2}\in F_{2}^{(H_{2})}$, the support of $\alpha^{-1}(d_{1})^{-1}\alpha^{-1}(d_{2})$ lies in the $(A+B)-$neighbourhood of $\overline{\beta}\big(\text{supp}\big(d_{1}^{-1}d_{2})\big)$. Applying this claim with $d_{1}=\alpha(c_{1})$, $d_{2}=\alpha(c_{2})$, we get 
\begin{equation*}
    \text{supp}(c_{1}^{-1}c_{2}) \subset \overline{\beta}\big(\text{supp}(\alpha(c_{1})^{-1}\alpha(c_{2}))\big)^{+(A+B)}
\end{equation*}
and applying $\beta$ to both sides and using Lemma~\ref{lem:neighborhoodsandQI}\textit{(iv)}, and that $\beta\circ\overline{\beta}$ lies at distance $\le B$ from $\text{Id}_{H_{2}}$, it follows that 
\begin{align*}
    \beta(\text{supp}(c_{1}^{-1}c_{2})) &\subset \beta\left(\overline{\beta}(\text{supp}(\alpha(c_{1})^{-1}\alpha(c_{2})))^{+(A+B)}\right) \\
    &\subset \beta\left(\overline{\beta}(\text{supp}(\alpha(c_{1})^{-1}\alpha(c_{2})))\right)^{+(A(A+B)+B)} \\
    &\subset \text{supp}\big(\alpha(c_{1})^{-1}\alpha(c_{2})\big)^{+(A(A+B)+2B)}.
\end{align*}
We conclude that the Hausdorff distance between $\text{supp}(\alpha(c_{1})^{-1}\alpha(c_{2}))$ and $\beta\big(\text{supp}(c_{1}^{-1}c_{2})\big)$ is $\le Q\defeq A(A+B)+2B$, which proves~\textit{(iii)}. 

\begin{figure}[H]
  \centering
  \includegraphics[width=0.9\linewidth]{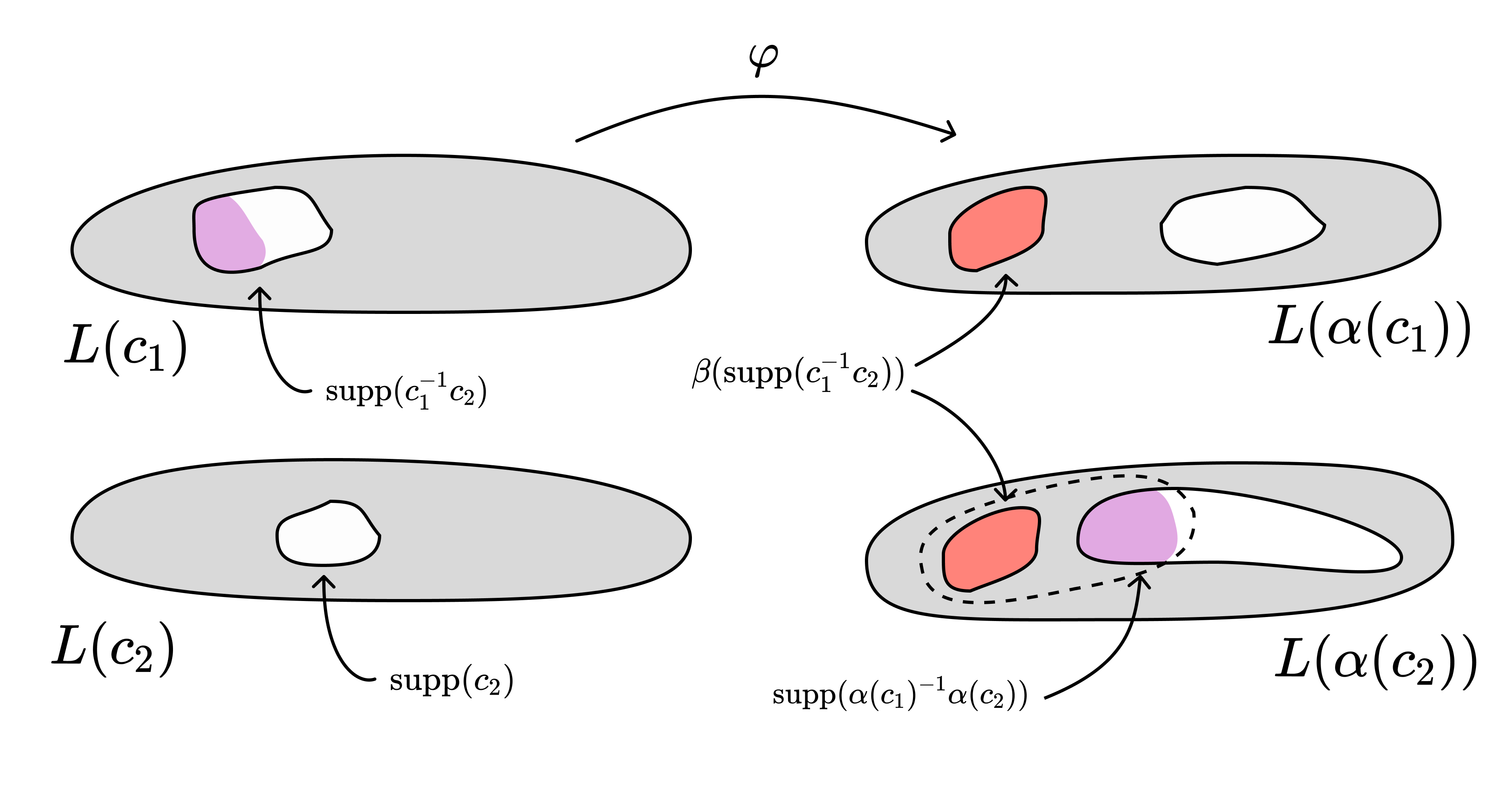}
  \caption{Configuration in the proof of~\textit{(iii)}}
\end{figure}

\noindent Conversely, let us assume~\textit{(i)}-\textit{(iii)} and let us show that $\varphi$ is a quasi-isometry. Let $A\ge 1$, $B\ge 0$ be the parameters of $\beta$ and one of its quasi-inverses $\overline{\beta}$. Let $a$, $b$ be two adjacent vertices in $\text{Cay}(F_{1}\wr H_{1}, F_{1}\cup S_{1})$. There are two cases: 
\begin{itemize}
    \item If $a=(c,p)$ and $b=(c,q)$ for some neighbour $q$ of $p$ in $H_{1}$, then 
    \begin{equation*}
        d(\varphi(a),\varphi(b))=d\big((\alpha(c), \beta(p)), (\alpha(c), \beta(q))\big)=d(\beta(p), \beta(q)) \le A\cdot d(p,q)+B=A+B.
    \end{equation*}
    \item If $a=(c,p)$ and $b=(c',p)$ where $c$ and $c'$ only differ on $p\in H_{1}$, we rather get 
    \begin{align*}
        d(\varphi(a),\varphi(b))&=d\big((\alpha(c), \beta(p)), (\alpha(c'), \beta(p))\big) \\
        &\le \big|\text{supp}(\alpha(c)^{-1}\alpha(c'))\big| \\
        &\le \big|B_{H_{2}}(\beta(p), Q)\big|\\
        &\le D^{Q}
    \end{align*}
    where $D\ge 2$ is a fixed integer at least as large as the maximal degree of a vertex in $\text{Cay}(H_{2},S_{2})$, since, by assumption~\textit{(iii)}, $\text{supp}(\alpha(c)^{-1}\alpha(c'))$ is contained in 
    \begin{equation*}
        \beta\big(\text{supp}(c^{-1}c')\big)^{+Q}=\lbrace \beta(p)\rbrace^{+Q}=B_{H_{2}}(\beta(p), Q). 
    \end{equation*}
\end{itemize}
In any case, we get that $d(\varphi(a),\varphi(b)) \le \max(A+B, D^{Q})$ for any two adjacent vertices $a,b$, so $\varphi$ is $\max(A+B, D^{Q})-$Lipschitz by Lemma~\ref{lem:Lipschitzbetweengraphs}. 

\noindent If we now define the map
\begin{align*}
    \psi\colon F_{2}\wr H_{2}&\longrightarrow F_{1}\wr H_{1} \\
    (d,q)&\longmapsto (\alpha^{-1}(d), \overline{\beta}(q))
\end{align*}
then it also satisfies assumption~\textit{(iii)}. Indeed, for any $d_{1},d_{2}\in F_{2}^{(H_{2})}$, assumption~\textit{(iii)} for $\varphi$ says that the Hausdorff distance between $\beta\big(\text{supp}(\alpha^{-1}(d_{1})^{-1}\alpha^{-1}(d_{2}))\big)$ and $\text{supp}(d_{1}^{-1}d_{2})$ is at most $Q$. Thus the Hausdorff distance between $\overline{\beta}\big(\beta\big(\text{supp}(\alpha^{-1}(d_{1})^{-1}\alpha^{-1}(d_{2}))\big)\big)$ and $\overline{\beta}\big(\text{supp}(d_{1}^{-1}d_{2})\big)$ is at most $A\cdot Q+B$ by Lemma~\ref{lem:neighborhoodsandQI}\textit{(iii)}. Since the first set is at Hausdorff distance at most $B$ from $\text{supp}\big(\alpha^{-1}(d_{1})^{-1}\alpha^{-1}(d_{2})\big)$, it follows that the Hausdorff distance between $\text{supp}\big(\alpha^{-1}(d_{1})^{-1}\alpha^{-1}(d_{2})\big)$ and $\overline{\beta}\big(\text{supp}(d_{1}^{-1}d_{2})\big)$ is at most $A\cdot Q+2B$. Thus~\textit{(iii)} holds for $\psi$ with $Q'\defeq A\cdot Q+2B$. One can therefore reproduce the argument above to conclude that $\psi$ is also $\max(A+B, D'^{Q'})-$Lipschitz, where $D'\ge 2$ is a fixed integer larger than the maximal degree of a vertex in $\text{Cay}(H_{1},S_{1})$. 

\noindent Finally, since 
\begin{equation*}
    d\big(\psi\circ\varphi(c,p), (c,p)\big)=d\big(\overline{\beta}\circ\beta(p),p\big)
\end{equation*}
for any $(c,p)\in F_{1}\wr H_{1}$, and since $\overline{\beta}\circ\beta$ is at distance at most $B$ from $\text{Id}_{H_{1}}$, $\psi\circ\varphi$ is at distance at most $B$ from $\text{Id}_{F_{1}\wr H_{1}}$. Likewise, $\varphi\circ\psi$ is at distance at most $B$ from $\text{Id}_{F_{2}\wr H_{2}}$. 

\noindent Putting all this together, we deduce that
\begin{align*}
    d(\varphi(a), \varphi(b)) &\ge \frac{1}{\max(A+B, D'^{Q'})}\cdot d\big(\psi(\varphi(a)), \psi(\varphi(b))\big) \\
    &\ge \frac{1}{\max(A+B, D'^{Q'})}\big(d(a,b)-d(\psi(\varphi(a)),a)-d(\psi(\varphi(b)),b)\big) \\
    &\ge \frac{1}{\max(A+B, D'^{Q'})}\cdot d(a,b) - \frac{2B}{\max(A+B, D'^{Q'})}
\end{align*}
for any $a,b\in F_{1}\wr H_{1}$. Thus $\varphi$ is a quasi-isometry, with $\psi$ as a quasi-inverse.
\end{proof}

Since an aptolic quasi-isometry $F_{1}\wr H_{1}\longrightarrow F_{2}\wr H_{2}$ provides a constant $Q\ge 0$ with the property of Proposition~\ref{prop:characterisationofaptolicQI:finitecase}\textit{(iii)}, the key observation to deduce an arithmetic constraint on $|F_{1}|$ and $|F_{2}|$ is the following:
\begin{proposition}
Let $F_{1}$, $F_{2}$ be two non-trivial finite groups, and let $H_{1}$, $H_{2}$ be infinite finitely generated groups. Let $\varphi\colon F_{1}\wr H_{1}\longrightarrow F_{2}\wr H_{2}$ be an aptolic quasi-isometry, of the form $\varphi(c,p)=(\alpha(c),\beta(p))$, $(c,p)\in F_{1}\wr H_{1}$. Then, for every quasi-inverse $\overline{\beta}$ of $\beta$, there is a constant $Q\ge 0$ such that: 
\begin{align*}
    &\text{For every $Q'\ge Q$ and any $A_{1}\subset H_{1}$, $\alpha^{-1}\big(\mathcal{L}(\beta(A_{1})^{+Q'})\big)$ is a union of cosets of $\mathcal{L}(A_{1})$};\\
    &\text{conversely, for any $Q'\ge Q$ and any $A_{2}\subset H_{2}$, $\alpha\big(\mathcal{L}(\overline{\beta}(A_{2})^{+Q'})\big)$ is a union of cosets of $\mathcal{L}(A_{2})$}.
\end{align*}
As a consequence, $|F_{1}|$ and $|F_{2}|$ have the same prime divisors. 
\end{proposition}

We do not prove here the first part of this statement, since we establish an analogous result in Chapter~\ref{chap:chapter4} (Proposition~\ref{prop:finiteunionofcosets}), in the context of arbitrary lamp groups. But let us explain how to deduce from it the second part of the statement. 

Given a subset $A_{1}$ of $H_{1}$, $\alpha^{-1}\big(\mathcal{L}(\beta(A_{1})^{+Q})\big)$ is a union of cosets of $\mathcal{L}(A_{1})$, so there are $k\ge 1$ and colourings $c_{1},\dots,c_{k}$ such that 
\begin{equation*}
    \alpha^{-1}\big(\mathcal{L}(\beta(A_{1})^{+Q})\big) = \bigsqcup_{i=1}^{k}c_{i}\mathcal{L}(A_{1})
\end{equation*}
which implies, since $\alpha$ is a bijection, that
\begin{equation*}
    |\mathcal{L}(\beta(A_{1})^{+Q})|=|\alpha^{-1}\big(\mathcal{L}(\beta(A_{1})^{+Q})\big)|=\left|\bigsqcup_{i=1}^{k}c_{i}\mathcal{L}(A_{1})\right|=\sum_{i=1}^{k}|c_{i}\mathcal{L}(A_{1})|=k\cdot|\mathcal{L}(A_{1})|.
\end{equation*}
Since $\mathcal{L}(A_{1})=\bigoplus_{A_{1}}F_{1}$, it follows that
\begin{equation*}
    k\cdot |F_{1}|^{|A_{1}|}=k\cdot |\mathcal{L}(A_{1})|=|\mathcal{L}(\beta(A_{1})^{+Q})|=|F_{2}|^{|\beta(A_{1})^{+Q}|}
\end{equation*}
and this equality implies that $|F_{1}|$ and $|F_{2}|$ share the same prime divisors. 

Finally, in the case where $H_{1}$ and $H_{2}$ are amenable, this conclusion can be strengthened to see that $|F_{1}|$ and $|F_{2}|$ are in fact powers of a common number. The details of the computations leading to this conclusion are provided in~\cite[Theorem~3.9]{GT24b}, or alternatively in Theorem~\ref{thm:AlamplightersoverTBPgroups-rigiditypart} below. 

\paragraph{Constructing quasi-isometries.} Conversely, let us now explain how to deduce the flexibility part of Theorem~\ref{thm:classificationGT21}. We consider first the non-amenable case, in which we have the following existence result:
\begin{theorem}[{\cite[Proposition~3.13]{GT24b}}]\label{thm:constructingQIbetweenNAlamplightersGT21}
Let $n,m\ge 2$. Let $G$ be a non-amenable finitely generated group. If $n$ and $m$ have the same prime divisors, then there exists an aptolic quasi-isometry from $\Z_{n}\wr G$ to $\Z_{m}\wr G$.  
\end{theorem}

We only sketch the main idea of the proof on an example, namely $\Z_{6}\wr F_{2}$ and $\Z_{54}\wr F_{2}$, and we refer to~\cite[Proposition~3.13]{GT24b} for the complete proof. The first thing to do is to replace $\Z_{54}\wr F_{2}$ with $(\Z_{2}\oplus\Z_{3}^{3})\wr F_{2}$ (which changes nothing since this group is biLipschitz equivalent to the first one, thanks to Proposition~\ref{prop:BiLipwreathproductscolor}). Geometrically, this amounts to splitting each lamp on the Cayley graph of $F_{2}$ in two half-lamps. Each finitely supported colouring $c\colon F_{2}\rightarrow \Z_{54}$ becomes the “sum” $c_{1}\oplus c_{2}$ of two colourings $c_{1}\colon F_{2}\rightarrow \Z_{2}$, $c_{2}\colon F_{2}\rightarrow \Z_{3}^{3}$. We do the same with $\Z_{6}\wr F_{2}$ that becomes $(\Z_{2}\oplus\Z_{3})\wr F_{2}$. Now, to any colouring $c\colon F_{2}\rightarrow \Z_{3}$, one can create a new colouring $\overline{c}\colon F_{2}\rightarrow \Z_{3}^3$ as follows: fix a direction at infinity $\xi\in\partial F_{2}$, and for any $p\in F_{2}$, let $\overline{c}(p)\in \Z_{3}^3$ be given by $(c(p_{1}), c(p_{2}), c(p_{3}))$, where $p_{1},p_{2},p_{3}$ are the three vertices separated from $\xi$ by $p$. Define then the map 
\begin{align*}
    \varphi\colon (\Z_{2}\oplus\Z_{3})\wr F_{2} &\longrightarrow (\Z_{2}\oplus\Z_{3}^3)\wr F_{2} \\
(c_{1}\oplus c_{2},p)&\longmapsto (c_{1}\oplus\overline{c_{2}}, p)
\end{align*}
namely we keep unchanged the first half of the lamps and the position of the arrow, and we only change the second half of the lamps: 

\begin{figure}[H]
  \centering
  \includegraphics[
    width=0.9\linewidth
  ]{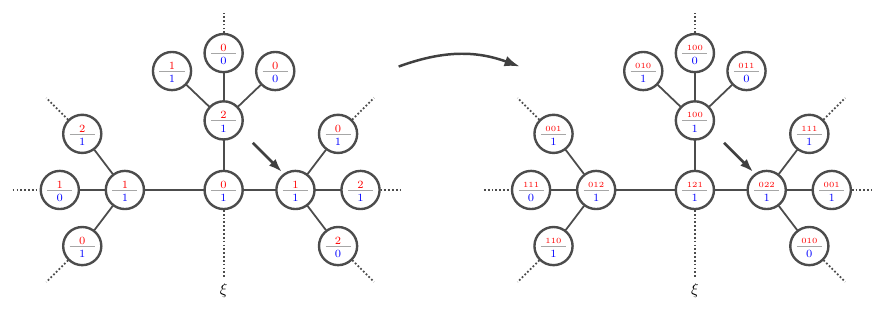}
  \caption{A quasi-isometry $(\Z_{2}\oplus\Z_{3})\wr F_{2} \longrightarrow (\Z_{2}\oplus\Z_{3}^3)\wr F_{2}$}
  \label{fig:cayley}
\end{figure}

In this example, the key point is that there is a $3$-to-one map $F_{2}\rightarrow F_{2}$ that lies at bounded distance from the identity, namely the map that sends any vertex of $F_{2}$ to its neighbour according to $\xi$. Such maps cannot exist between amenable spaces (as a consequence of Lemma~\ref{lem:uniquenessscalingfactor}), hence the appearance of rigidity and the need for another strategy for constructing quasi-isometries. 

With Theorem~\ref{thm:constructingQIbetweenNAlamplightersGT21}, the proof of the right-to-left direction of the first point of Theorem~\ref{thm:classificationGT21} is straightforward: indeed, if $|F_{1}|$ and $|F_{2}|$ have the same prime divisors, then there is a quasi-isometry $f$ from $F_{1}\wr H_{1}$ to $F_{2}\wr H_{1}$. The quasi-isometry between $H_{1}$ and $H_{2}$ that we have by assumption is quasi-one-to-one, so it can be replaced by a bijective quasi-isometry according to Theorem~\ref{thm:Whytethm}, and thus $F_{2}\wr H_{1}$ and $F_{2}\wr H_{2}$ are biLipschitz equivalent thanks to Proposition~\ref{prop:BiLipwreathproductsbase}. It remains to pre-compose this biLipschitz equivalence with $f$ to conclude. 

Between amenable spaces, one also has to make local modifications of the colouring that takes into account the exponents in the cardinalities of the lamp groups. Here the trick is to use instead the criterion of Theorem~\ref{thm:rationalscalingfactor} for scaling quasi-isometries and modify colors on pieces of the partition thus obtained. We refer to~\cite[Proposition~3.12]{GT24b} for the detailed proof, and to the proof of Proposition~\ref{prop:aptolicQIfromscalingQIbetweenquotients-Acase} for a similar construction but in the context of permutational lamplighters. 

\section{Quasi-isometric rigidity for halo products}\label{sec:QIrigidityforhaloproducts} Later, Genevois and Tessera managed to extend their techniques to other~\textit{halo products}, that we define now. 

\begin{definition}
Let $X$ be a set. A~\textit{halo of groups over $X$} is the data, for any subset $S\subset X$, of a group $L(S)$ such that:
\begin{itemize}
    \item for all $R,S\subset X$, if $R\subset S$, then $L(R)\leqslant L(S)$; 
    \item $L(\emptyset)=\lbrace 1\rbrace$ and $L(X)=\langle L(S) : S\subset X\;\text{finite}\rangle$;
    \item for all $R,S\subset X$, $L(R\cap S)=L(R)\cap L(S)$. 
\end{itemize}
\end{definition}

Given an action $H\act X$ and a morphism $\alpha\colon H \longrightarrow \text{Aut}(L(X))$ satisfying $\alpha(h)(L(S))=L(hS)$ for all $S\subset X$ and all $h\in H$, the \textit{permutational halo product} $\halo_{X,\alpha}H$ is the semi-direct product
\begin{equation*}
    \halo_{X,\alpha}H \defeq L(X)\rtimes H.
\end{equation*}

Let us mention our favorite examples. 

\paragraph{Lampjugglers.} Given a group $H$ and an integer $r\ge 1$, the \textit{lampjuggler over $H$} is the semi-direct product
\begin{equation*}
    \juggler{r}{H} \defeq \fsym{H\times\lbrace 1,\dots,r\rbrace} \rtimes H
\end{equation*}
where $\fsym{H\times\lbrace 1,\dots,r\rbrace}$ is the group of finitely supported permutations of $H\times \lbrace 1,\dots,r\rbrace$ and where $H$ acts on it through its initial action on $H\times\lbrace 1,\dots,r\rbrace$ given by $h\cdot (x,i) \defeq (hx, i)$. It can be described as the halo group $\halo H$ where 
\begin{equation*}
    L(S)\defeq \fsym{S\times\lbrace 1,\dots, r\rbrace}, S\subset H.
\end{equation*}
Lampjugglers over finitely generated groups are finitely generated, and one can check that if $S_{H}$ is a finite generating set for $H$, then the finite set
\begin{equation*}
    \left\lbrace (\tau_{(1_{H},i),(s,j)}, 1_{H}) : s\in S_{H}, 1\le i,j\le r\right\rbrace  \cup \lbrace (\text{id},s) : s\in S_{H}\rbrace
\end{equation*}
generates $\juggler{r}{H}$, where $\tau_{(1_{H},i),(s,j)}$ stands for the transposition swapping $(1_{H},i)$ with $(s,j)$. A particular case is that of \textit{lampshufflers} $\shuf{H}$, i.e. lampjugglers with $r=1$. 

\paragraph{Lampdesigners.} Let $F$ and $H$ be two groups. The \textit{lampdesigner over $H$} is the semi-direct product 
\begin{equation*}
    \designer{H} \defeq (F\wr_{H}\fsym{H})\rtimes H
\end{equation*}
where $H$ acts on $\bigoplus_{H}F$ by permuting the coordinates through its action on itself by left-multiplication and acts on $\fsym{H}$ as described above. It is the halo product $\halo H$ for the collection $L(S)\defeq F\wr_{S}\fsym{S}$, $S\subset H$. 

\paragraph{Lampcloners.} Let $H$ be a group and let $\field$ be a field. Denote $V_{H}$ the $\field-$vector space admitting $H$ as a basis, and denote by $\lbrace e_{u} : u\in H\rbrace$ a formal basis. Let $\text{FGL}(H)$ be the group of linear automorphisms $V_{H}\rightarrow V_{H}$ that fix all but finitely many basis elements. This group can also be seen as the group of finitely supported invertible matrices with coefficients in $\field$ whose entries are indexed by $H\times H$. Once again, the action of $H$ on itself naturally yields an action of $H$ on $\text{FGL}(H)$. The \textit{lampcloner over $H$} is the semi-direct product 
\begin{equation*}
    \cloner{H} \defeq \text{FGL}(H)\rtimes H.
\end{equation*}
It is a halo product, for the collection $L(S)\defeq \text{FGL}(S)$, $S\subset H$, where $\text{FGL}(S)$ is thought of as the subgroup of $\text{FGL}(H)$ of linear automorphisms $V_{H}\rightarrow V_{H}$ that fix $H\setminus S$ and that stabilise the subspace $\langle S\rangle\subset V_{H}$. 

We refer the reader to Section~\ref{sec:halo} for a more precise description of these examples and additional ones. 

At this point, we only emphasize the fact that the aptolicity phenomenon also holds for quasi-isometries between lampshufflers, lampcloners or lampdesigners, as shown in~\cite{GT24a}. The proof of this fact relies on the next statement, that vastly generalises Theorem~\ref{thm:embedding1}. 

\begin{theorem}[{\cite[Theorem~4.17]{GT24a}}]\label{thm:embedding2}
Let $Z$ be a geodesic metric space satisfying the thick bigon property and let $\halo H$ be a finitely generated halo product with $\halo$ full and $L(H)$ locally finite. Then every coarse embedding $\rho\colon Z \rightarrow \halo H$ has its image contained in a neighbourhood of an $H-$coset. Moreover, the size of this neighbourhood only depends on $Z$, $\halo H$, and the parameters of $\rho$. 
\end{theorem}

Here, a halo product $\halo H$ over a group $H$ equipped with a word metric $d$ is~\textit{full} if there is a constant $K\ge 0$ such that, for all subsets $R$, $S$, $T$ and $U$ of $H$, if 
\begin{equation*}
    d(R,S), \; d(R,T), \; d(R,U), \; d(S,T) \ge K
\end{equation*}
then 
\begin{equation*}
    L(R\cup S)L(U) \cap L(R\cup T) \subset L(R)L(T).
\end{equation*}
It is a technical assumption needed in the proof, and for now it is enough to keep in mind that our favorite halo products, namely wreath products, lampjugglers, and lampcloners are full~\cite[Lemma~4.39]{GT24a}.

\paragraph{Main idea for proving Theorem~\ref{thm:embedding2}.} Fix a finitely generated halo product $\halo H$ with $L(H)$ locally finite. The first ingredient of the proof is to realize $\halo H$ as a subcomplex $\mathcal{C}_{k}$ of some huge cubical complex $\mathcal{C}$, for any $k\ge 1$. If we think of an element $(c,p)\in\halo H$ as a colouring $c\in L(H)$ and an arrow pointing at $p\in H$, then moving in $\halo H$ amounts to moving the arrow in $H$ and, along the way, to modifying the colouring by an element of $L(H)$ supported in a small ball around the arrow, and $\mathcal{C}_{k}$ is constructed as follows: replace the arrow by a finite connected subgraph of $H$, that we call a~\textit{crowd}. The~\textit{height} of the crowd is the size of the finite connected subgraph, and moving in $\halo H$ amounts to moving the crowd by adding or removing vertices and, along our way, to modifying the colouring by an element of $L(H)$ supported in the crowd. 

Now, at a fixed size $k\ge 2$, two points $x,y\in \mathcal{C}_{k}$ that are not in a neighbourhood of a leaf must be separated in a topological way by specific bounded subcomplexes, called~\textit{blocks} in \cite[Definition~4.10]{GT24a}. For lamplighters, this comes from the following observation. Fix a wreath product $\Z/2\Z\wr H$, and a path $\gamma$ from $1_{\halo H}$ to a point $(c,h)$. Fix a point $p\in H$ where $c$ is non-trivial (i.e. $p\in\text{supp}(c)$). Following $\gamma$, the arrow must go in $H$ from $1_{H}$ to $h$ and, along the way, visits all points in $\text{supp}(c)$ to change the color. In particular, at some stage, the arrow will pass by $p$. The same applies to any other path $\gamma'$ connecting $1_{\halo H}$ to $(c,h)$. Now, on $\gamma$ and $\gamma'$, consider the points $(c_{n}, p_{n})\in\gamma$, $(c_{m}', p_{m}')\in\gamma'$ where the arrow points at $p$, i.e. $p_{n}=p_{m}'=p$. We compare $c_{n}$ and $c_{m}'$ outside of $p$. If these colourings are very different, this means that the orders according to which lamps are switched on are very different whether we are following $\gamma$ or $\gamma'$. For instance, there is $q\in\text{supp}(c)$ very far from $p$ in $H$ such that $\gamma$ turns on $p$ after $q$, while $\gamma'$ turns on $q$ after $p$. But, if $p$ and $q$ are very far, the cycle in $\Z/2\Z\wr H$ that turns on $p$, turns on $q$, turns off $p$ and finally turns off $q$ to go back to $1$, can be homotopically trivial only at a very large scale. Said differently, if we are given $\gamma$ and a fixed scale at which we can modify $\gamma$ up to coarse homotopy, then, along our new path, we can always find a point where the arrow points at $p$ and where the colouring outside $p$ is fixed. This observation leads the authors of~\cite{GT24a} to introduce the notion of~\textit{essential separation} and to prove that, at a given scale $k\ge 2$, any path connecting $x$ and $y$ in $\mathcal{C}_{k}$ must intersect essentially a block. 

Finally, fix for instance a coarsely simply connected and one-ended graph $Z$, and a coarse embedding $Z\rightarrow \halo H$. Since $\halo H$ sits naturally inside $\mathcal{C}$, we get a coarse embedding $Z\rightarrow \mathcal{C}_{2}\subset\mathcal{C}$. The image of $Z$ must then be simply connected in $\mathcal{C}_{k}$ for some large $k\ge 2$. But then, as a consequence of what we explained above, if the image is not in a neighbourhood of a leaf, then it contains a path $\gamma$ that connects two points $x$ and $y$ that are very far and such that any homotopically equivalent path to $\gamma$ has to intersect a given bounded block. This contradicts one-endedness of the image, which guarantees that any path connecting two points in the image can be deformed up to coarse homotopy to avoid a bounded region. It follows that the image of $Z$, as it is both coarsely simply connected and one-ended, must be contained in a neighbourhood of a leaf in $\halo H$. 

In this strategy, one-endedness is not that important; what really matters is to always be able to coarsely modify paths so that the resulting cycles are homotopically trivial. This property has a name, namely the~\textit{thick bigon property}, recently introduced in~\cite{GT24a} and that we explain in detail now. In particular, finitely presented one-ended groups have the thick bigon property, so that Theorem~\ref{thm:embedding2} is indeed a generalisation of Theorem~\ref{thm:embedding1}.

\paragraph{The thick bigon property.} To define the thick bigon property, we need to recall the notion of~\textit{coarsely homotopic paths} in a metric space.

\begin{definition}
Let $X$ be a metric space. Given $C\ge 0$, two paths $\alpha$ and $\beta$ with the same endpoints are~\textit{$C-$coarsely homotopic} if there exists a sequence of paths 
\begin{equation*}
    \alpha=\gamma_{0}, \gamma_{1},\dots,\gamma_{n-1},\gamma_{n}=\beta
\end{equation*}
such that, for any $0\le i\le n-1$, $\gamma_{i+1}$ is obtained from $\gamma_{i}$ by replacing a subpath $\zeta\subset \gamma_{i}$ with a path $\xi$ between the same endpoints satisfying $\text{diam}(\zeta\cup\xi)\le C$. A path is~\textit{$C-$coarsely homotopically trivial} if it is $C-$coarsely homotopic to a constant path. 
\end{definition}

Then, a metric space $X$ is~\textit{$C-$coarsely simply connected} if every path is $C-$coarsely homotopically trivial. A space $X$ is~\textit{coarsely simply connected} if it is $C-$coarsely simply connected for some $C\ge 0$. 

Recall that a finitely generated group is coarsely simply connected if and only if it is finitely presented (see e.g.~\cite{CH16}).

\begin{definition}\label{def:TBP}
A metric space $X$ has the~\textit{thick bigon property} if there exists some $C\ge 0$ such that, for every $R\ge 0$, there exists $L\ge 0$ such that the following holds: any two points $x,y\in X$ can be connected by some path $\gamma_{1}$ such that, for any point $p\in \gamma_{1}$ satisfying $d(p,x), d(p,y)\ge L$, there is a path $\gamma_{2}$ between $x$ and $y$ that is $C-$coarsely homotopic to $\gamma_{1}$ and that avoids $B_{X}(p,R)$.
\end{definition}

\begin{figure}[H]
  \centering
  \includegraphics[width=0.5\linewidth]{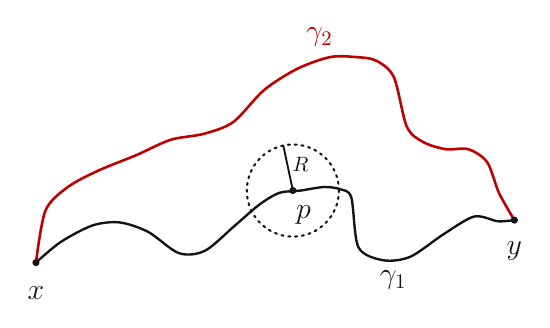}
  \caption{Configuration of points in a space $X$ having the thick bigon property}
\end{figure}

 We need one more terminology to exhibit a first source of spaces with the thick bigon property.

\begin{definition}
Let $X$ be a geodesic metric space. Given a map $\rho$, $X$ is~\textit{$\rho-$uniformly one-ended} if, for every $R\ge 0$, the complement of a ball of radius $R$ has one unbounded connected component and all other components have diameter $\le \rho(R)$. 
\end{definition}

A space $X$ is~\textit{uniformly one-ended} if it is $\rho-$uniformly one-ended for some map $\rho$. 

The first source of examples is then given by coarsely simply connected uniformly one-ended spaces.

\begin{lemma}[{\cite[Lemma~3.4]{GT24a}}]
Let $X$ be a geodesic metric space. If $X$ is uniformly one-ended and coarsely simply connected, then $X$ has the thick bigon property.
\end{lemma}

\begin{proof}
Fix $C\ge 0$ and a map $\rho$ such that $X$ is $C-$coarsely simply connected and $\rho-$uniformly one-ended. Let $R\ge 0$ and set $L\defeq R+\rho(R)+1$.

\noindent Fix two points $x,y\in X$ connected by a path $\gamma_{1}$, and let $p\in \gamma_{1}$ be a point such that $d(p,x),d(p,y)\ge L$. The subpath of $\gamma_{1}$ connecting $x$ to $B(p,R)$ has length at least
\begin{equation*}
    d(x,p)-R \ge L-R =\rho(R)+1>\rho(R).
\end{equation*}
Thus $x$ belongs to the unique unbounded connected component of $X\setminus B(p,R)$, and likewise for $y$. We deduce that there is a path $\gamma_{2}$ connecting $x$ and $y$ that avoids $B(p,R)$. As $X$ is $C-$coarsely simply connected, $\gamma_{1}$ and $\gamma_{2}$ are $C-$coarsely homotopic. This concludes the proof. 
\end{proof}

In particular, Theorem~\ref{thm:embedding2} indeed generalises Theorem~\ref{thm:embedding1}. Note also that finitely presented one-ended groups have the thick bigon property. 

On the other hand, the following notion allows to also consider some infinitely presented groups. 

\begin{definition}\label{def:fleshy}
Let $X$ be a metric space and let $\mathcal{P}$ be a collection of uniformly coarsely embedded subspaces. We say that $X$ is~\textit{fleshy relative to $\mathcal{P}$} if:
\begin{itemize}
    \item There is $C\ge 0$ and a map $\rho$ such that any $P\in\mathcal{P}$ is $C-$coarsely simply connected and $\rho-$uniformly one-ended;
    \item For any $x,y\in X$, there exist $P_{1},\dots, P_{n}\in\mathcal{P}$ such that $x\in P_{1}$, $y\in P_{n}$ and $P_{i}\cap P_{i+1}$ is unbounded for any $1\le i\le n-1$.
\end{itemize}
\end{definition}

A metric space $X$ is~\textit{fleshy} if it is fleshy relative to some collection of subspaces. 

\begin{proposition}[{\cite[Proposition~3.6]{GT24a}}]\label{prop:fleshyimpliesTBP}
Let $X$ be a geodesic metric space. If $X$ is fleshy, then $X$ has the thick bigon property.
\end{proposition}

It is proved in~\cite[Proposition~3.9]{GT24a} that wreath products $E\wr H$ with $E$ infinite and $H$ non-trivial are fleshy. In particular, such wreath products have the thick bigon property. 

On the other hand, wreath products $F\wr H$, with $F$ a non-trivial finite group, do not have the thick bigon property. Indeed, if it were the case, applying Theorem~\ref{thm:embedding2} to $Z=F\wr H$ and $\rho=\text{Id}$ would show that the identity must have its image completely contained in the neighbourhood of an $H-$coset, which is obviously not the case. The same argument shows that any halo product $\halo H$ which is full and whose $L(H)$ is locally finite cannot have the thick bigon property. This encompasses lampshufflers and lampcloners.

Let us also mention another sufficient condition for a group to have the thick bigon property. 

\begin{lemma}[{\cite[Lemma~3.14]{GT24a}}]\label{lem:conjugateofsubgroupimpliesTBP}
Let $G$ be a group with a finite generating set $S$. If $H\leqslant G$ is a finitely presented one-ended subgroup of $G$ such that $|H\cap sHs^{-1}|=\infty$ for every $s\in S$, then $G$ is fleshy. 
\end{lemma}

\begin{proof}
Fix $x,y\in G$ and write $x^{-1}y=s_{1}\dots s_{n}$ as a product of generators and their inverses. Consider the sequence of $H-$cosets
\begin{equation*}
    xH, xs_{1}H,\dots, xs_{1}\dots s_{n-1}H, xs_{1}\dots s_{n}H=yH.
\end{equation*}
For any $0\le i\le n-1$, $xs_{1}\dots s_{i}(H\cap s_{i+1}Hs_{i+1}^{-1})$ is contained in the $1-$neighbourhood of $xs_{1}\dots s_{i}H$ and of $xs_{1}\dots s_{i}s_{i+1}H$, and is unbounded by assumption. It follows that $G$ is fleshy relative to $\lbrace (gH)^{+1} : g\in G\rbrace$. 
\end{proof}

For instance, this result implies that any finitely generated group with a normal one-ended finitely presented subgroup has the thick bigon property. 

Another family of groups that is captured is the one of permutational wreath products, that will be of special interest for us. 

\begin{corollary}\label{cor:PWPhaveTBP}
Let $F$ and $H$ be non-trivial finitely generated groups. If $H$ is one-ended and finitely presented, and if $K\leqslant H$ is infinite, then $F\wr_{H/K}H$ is fleshy. In particular, it has the thick bigon property.
\end{corollary}

\begin{proof}
Fix finite generating sets $R\subset F$ and $S\subset H$. For every $s\in S$, $H\cap sHs^{-1}=H$ is infinite, and for $r\in F$, $H\cap rHr^{-1}$ is infinite since it contains $K$. Thus Lemma~\ref{lem:conjugateofsubgroupimpliesTBP} applies and shows that $F\wr_{H/K}H$ is fleshy. 
\end{proof}

Lemma~\ref{lem:conjugateofsubgroupimpliesTBP} can also be used to prove the following theorem.

\begin{theorem}[{\cite[Theorem~3.12]{GT24a}}]\label{thm:freeabeliansubgroupimpliesTBP}
Let $G$ be a finitely generated group. If $G$ is not virtually cyclic and contains a normal free abelian subgroup of positive rank, then $G$ satisfies the thick bigon property.   
\end{theorem}

In particular, we have:

\begin{corollary}
If a finitely generated torsion-free group is solvable but not cyclic, then it satisfies the thick bigon property. 
\end{corollary}

For instance, the corollary shows that all free solvable groups have the thick bigon property, even though they are infinitely presented~\cite[Corollary~2.14]{Cor06}. 

Lastly, the next criterion is of independent interest, even though it is used in~\cite{GT24a} as a step towards the proof of Theorem~\ref{thm:freeabeliansubgroupimpliesTBP}.

\begin{proposition}[{\cite[Proposition~3.16]{GT24a}}]\label{normalf.g.subgroupimpliesTBP}
Let $G$ be a finitely generated group. If $G$ contains an infinite normal subgroup $N\lhd G$ of infinite index which is finitely generated, then $G$ satisfies the thick bigon property.
\end{proposition}

Here, the assumption of finite generation of the subgroup $N$ cannot be removed: indeed, the lamplighter $G=\Z/2\Z\wr\Z^2$ has $\bigoplus_{\Z^2}\Z/2\Z$ as a normal subgroup of infinite index, and yet $G$ does not have the thick bigon property.

An elementary but still crucial fact, left unproved in~\cite{GT24a}, is that the thick bigon property is a quasi-isometry invariant among geodesic metric spaces. 

\begin{proposition}\label{prop:TBPisaQIinvariant}
Let $X$ and $Y$ be geodesic metric spaces. If $X$ and $Y$ are quasi-isometric, and if $X$ has the thick bigon property, then $Y$ has the thick bigon property.
\end{proposition}

\begin{proof}
Suppose that $f$ is a $(C,K)-$quasi-isometry, with any point of $Y$ being at distance $\le K$ from $f(X)$. It is well-known that if two paths are $C-$coarsely homotopic in $X$ and $f\colon X\rightarrow Y$ is a quasi-isometry, then their images in $Y$ are $C'-$coarsely homotopic for some $C'\ge 0$. See for instance~\cite[Proposition~6.A.7]{CH16}.

\noindent Next, let $R\ge 0$, and set $L\defeq C\cdot L'+K$, where $L'$ is given by the thick bigon property of $X$ for $R'\defeq C\cdot (R+K)$. Fix $y_{1},y_{2}\in Y$, and choose $x_{1},x_{2}\in X$ such that $d(y_{1},f(x_{1})), d(y_{2},f(x_{2}))\le K$. Connect $x_{1},x_{2}$ by some path $\gamma_{1}$ as in the definition, so that $f(x_{1}), f(x_{2})$ are connected by $f(\gamma_{1})$. Fixing a geodesic $[y_{1},f(x_{1})]\subset Y$ and a geodesic $[f(x_{2}), y_{2}]\subset Y$, the concatenation 
\begin{equation*}
    \eta_{1}\defeq [y_{1},f(x_{1})]\cup f(\gamma_{1})\cup [f(x_{2}), y_{2}]
\end{equation*}
is a path from $y_{1}$ to $y_{2}$. Now let $p\in\eta_{1}$ be such that $d(y_{1},p), d(y_{2},p)\ge L$. From the definition of $L$ and the fact that $y_{1}$ (resp. $y_{2}$) is at distance $\le K$ from $f(x_{1})$ (resp. $f(x_{2})$), $p$ necessarily lies on the subpath $f(\gamma_{1})\subset\eta_{1}$. Then $p=f(q)$ for some $q\in\gamma_{1}$, and additionally $d(q,x_{1}), d(q,x_{2})\ge L'$ since $d(y_{1},p), d(y_{2}, p)\ge L$. Hence there exists $\gamma_{2}$ a path between $x_{1}$ and $x_{2}$ that is $C-$coarsely homotopic to $\gamma_{1}$ and that avoids $B(q,R')$. Then $f(\gamma_{2})$ avoids $B(p,R)$
and is $C'-$coarsely homotopic to $f(\gamma_{1})$. We conclude that the concatenation 
\begin{equation*}
    \eta_{2} \defeq  [y_{1},f(x_{1})]\cup f(\gamma_{2})\cup [f(x_{2}), y_{2}]
\end{equation*}
is a path between $y_{1}$ and $y_{2}$ avoiding $B(p,R)$ and $C'-$coarsely homotopic to $\eta_{1}$. Thus $Y$ also has the thick bigon property.
\end{proof}

\begin{figure}[H]
  \centering
  \includegraphics[width=0.9\linewidth]{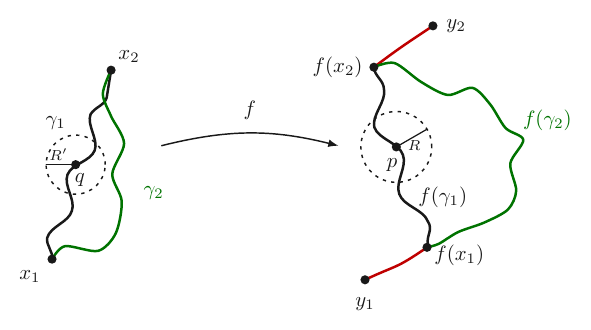}
  \caption{Configuration for the proof of Proposition 2.3.14}
  \label{fig:chemin}
\end{figure}

This result can be useful to distinguish up to quasi-isometry finitely generated groups that have many geometric properties in common. For instance:

\begin{proposition}
Let $d,r\ge 2$. Then the free solvable group $S_{d,r}$ and $\shufn{(d-1)}{\Z^r}$ are not quasi-isometric. 
\end{proposition}

Indeed, $S_{d,r}$ has the thick bigon property, while $\shufn{(d-1)}{\Z^r}$ does not. On the other hand, many other invariants fail to distinguish these two groups: they have the same volume growth, the same number of ends, the same isoperimetric profile (see Proposition~\ref{prop:profileofshufnofpolynomialgrowthgroups} below and~\cite{SCZ15}), and both are infinitely presented.

The embedding theorem and its main assumption being stated, let us now turn to the second step of the strategy: aptolicity.

\paragraph{Aptolicity.} As in the case of lamplighters over one-ended finitely presented groups~\cite{GT24b}, it is a consequence of the (upgraded) embedding theorem (cf. Theorem~\ref{thm:embedding2}) that any quasi-isometry 
\begin{equation*}
    \varphi\colon \mathcal{M}A\longrightarrow \mathcal{N}B
\end{equation*}
with $\mathcal{M}, \mathcal{N}$ full, $M(A)$, $N(B)$ locally finite, and $A$, $B$ having the thick bigon property, is \textit{leaf-preserving}, i.e. there is a constant $R\ge 0$ such that $\varphi$ sends any $A-$coset into the $R-$neighbourhood of a $B-$coset and has a quasi-inverse that does the same with roles of $A$ and $B$ reversed. Thus, up to bounded distance, $\varphi$ can be written as 
\begin{equation*}
    \varphi(c,p) = (\alpha(c), \beta_{c}(p))
\end{equation*}
for all $c\in M(A)$ and $p\in A$, for some bijection $\alpha\colon M(A)\longrightarrow N(B)$ and some collection of maps $\beta_{c}\colon A\rightarrow B$, indexed by $c\in M(A)$. 

At this stage, it is not hard to derive that each $\beta_{c}$ is a quasi-isometry between $A$ and $B$, since $\varphi$ is. Thus, it is already a consequence of the embedding theorem that $\mathcal{M}A$ and $\mathcal{N}B$ being quasi-isometric implies that $A$ and $B$ are quasi-isometric under our running assumptions. 

But more is true: it turns out that all $\beta_{c}$ lie at a uniform bounded distance from each other, and thus $\varphi$ itself lies at a bounded distance from an \textit{aptolic} quasi-isometry, in the sense that:

\begin{definition}
Let $\mathcal{M}A$, $\mathcal{N}B$ be finitely generated halo products. A quasi-isometry $\varphi\colon \mathcal{M}A\longrightarrow \mathcal{N}B$ is \textit{aptolic} if there is a bijection $\alpha\colon M(A)\longrightarrow N(B)$ and a quasi-isometry $\beta\colon A\rightarrow B$ such that 
\begin{equation*}
    \varphi(c,p)=(\alpha(c),\beta(p))
\end{equation*}
for all $(c,p)\in \mathcal{M}A$. 
\end{definition}

This observation on the collection $(\beta_{c})_{c\in M(A)}$ is proved with the help of geometric structures in halo products, called~\textit{squares} and~\textit{ladders} of leaves, that allow to detect geometrically when two points in a halo group have close projections onto the base group. 

To that end, we first construct an~\textit{angular graph} from a given halo product. Recall that such a graph $X$ comes with an angle map $\angle$ that assigns to each vertex $u\in X$ and each pair of neighbours $v,w\in X$ a real number $\angle_{o}(v,w)$. In such a graph, a subgraph $Y$ is called~\textit{$R-$obtuse} for some $R\ge 0$ if all angles in $Y$ are $\ge R$.

The angular graph that we construct from a halo product $\halo H$ encodes how leaves are organised in $\halo H$.

\begin{definition}
Let $\halo H$ be a finitely generated halo product. Let $\varepsilon\ge 0$. The~\textit{graph of leaves} $\mathcal{G}_{\varepsilon}(\halo H)$ is the graph whose vertices are $H-$cosets in $\halo H$ and whose edges connect two cosets whenever they are at distance $\le \varepsilon$. If $P$ is a vertex of $\mathcal{G}_{\varepsilon}(\halo H)$ and $Q,R$ are two neighbours of $P$ in $\mathcal{G}_{\varepsilon}(\halo H)$, the angle $\angle_{P}(Q,R)$ is the smallest distance between two points of $P$ respectively minimising the distance to $Q$ and $R$. 
\end{definition}

An important observation on this construction is that a leaf-preserving quasi-isometry between halo products provides quasi-isometric graphs of leaves.

\begin{proposition}[{\cite[Proposition~4.33]{GT24a}}]\label{prop:QIgraphofleaves}
Let $\mathcal{M}A$, $\mathcal{N}B$ be two finitely generated halo products. For all $C,K,\varepsilon\ge 0$, there exists $\eta\ge 0$ such that any leaf-preserving $(C,K)-$quasi-isometry $\varphi\colon \mathcal{M}A \longrightarrow \mathcal{N}B$ provides a bijective quasi-isometry $\overline{\varphi}\colon \mathcal{G}_{\varepsilon}(\mathcal{M}A) \longrightarrow \mathcal{G}_{\eta}(\mathcal{N}B)$ such that:
\begin{itemize}
    \item $\overline{\varphi}$ sends two adjacent vertices to two adjacent vertices;
    \item there exist $L,S>0$ such that
    \begin{equation*}
        \frac{1}{L}\cdot \angle_{P}(Q,R)-S \le \angle_{\overline{\varphi}(P)}(\overline{\varphi}(Q), \overline{\varphi}(R)) \le L\cdot \angle_{P}(Q,R)+S
    \end{equation*}
    for all vertex $P\in\mathcal{G}_{\varepsilon}(\mathcal{M}A)$ and neighbours $Q,R\in \mathcal{G}_{\varepsilon}(\mathcal{M}A)$. 
\end{itemize}
\end{proposition}

We can now introduce properly squares and ladders of leaves. 

\begin{definition}
Let $\halo H$ be a finitely generated halo product. Let $\varepsilon$, $R>0$. An $(\varepsilon,R)-$\textit{ladder of leaves} is a collection $P_{1},Q_{1},P_{2},Q_{2},\dots, P_{k},Q_{k}$ of leaves such that:
\begin{itemize}
    \item For any $1\le i \le k-1$, $P_{i}$ (resp. $Q_{i}$) lies at distance $\le \varepsilon$ from $P_{i-1}, Q_{i},P_{i+1}$ (resp. $Q_{i-1},P_{i}, Q_{i+1}$);
    \item For any $1\le i \le k-1$, the square $P_{i}, P_{i+1},Q_{i+1},Q_{i}$ is $R-$obtuse. 
\end{itemize}
\end{definition}

\begin{figure}[H]
  \centering
  \includegraphics[
    width=0.75\linewidth
  ]{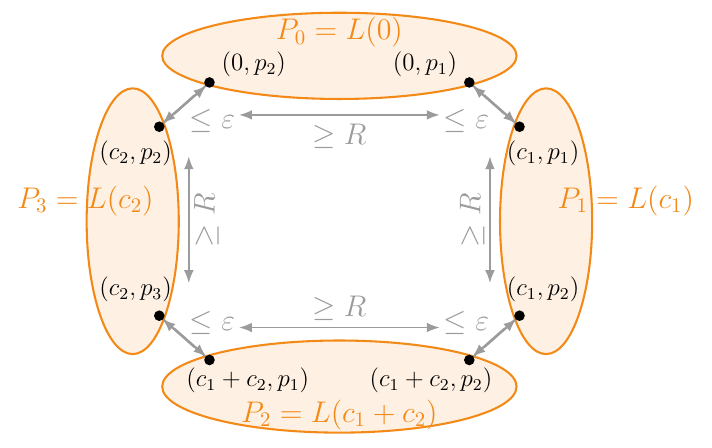}
  \caption{An $(\varepsilon,R)-$square of leaves in $\halo H$.}
  \label{fig:cayley}
\end{figure}

In~\cite{GT24a}, it is then proved that, for such a structure in a large-scale commutative halo product, if a point $(x,p)$ is close to both $P_{1}$ and $Q_{1}$ and a point $(y,q)$ is close to both $P_{k}$ and $Q_{k}$, then $p$ and $q$ must be close in $H$:

\begin{theorem}[{\cite[Theorem~6.3]{GT24a}}]\label{thm:distancebetweenprojections}
Let $\halo H$ be a finitely generated large-scale commutative halo product. Let $D\ge 0$ be the corresponding constant. For all $\varepsilon, \eta>0$ and $R>D+4\varepsilon+2r_{0}$, for every $(\varepsilon, R)-$ladder of leaves $P_{1},Q_{1},\dots, P_{k},Q_{k}$, if $(x,p)\in P_{1}^{+\eta}\cap Q_{1}^{+\eta}$ and $(y,q)\in P_{k}^{+\eta}\cap Q_{k}^{+\eta}$, then $d(p,q)\le 6\eta+\varepsilon$. 
\end{theorem}

A halo product $\halo H$ over a finitely generated group $H=\langle S_{H}\rangle$ is~\textit{large-scale commutative} if there is a constant $D\ge 0$ such that, for any $R,S\subset H$ with $d_{S_{H}}(R,S)\ge D$, the subgroups $L(R)$ and $L(S)$ commute. Moreover, if $\eta>0$ and if two points $(x,p), (y,q)\in \halo H$ are as in the statement, we say that $(x,p)$ and $(y,q)$ are~\textit{$\eta-$connected}. Thus, the previous theorem is saying that, in a large-scale commutative halo group, if two points are $\eta-$connected by a ladder whose squares have big enough angles, then the projections of our two points must be close in the base group. 

\begin{figure}[H]
  \centering
  \includegraphics[width=0.9\linewidth]{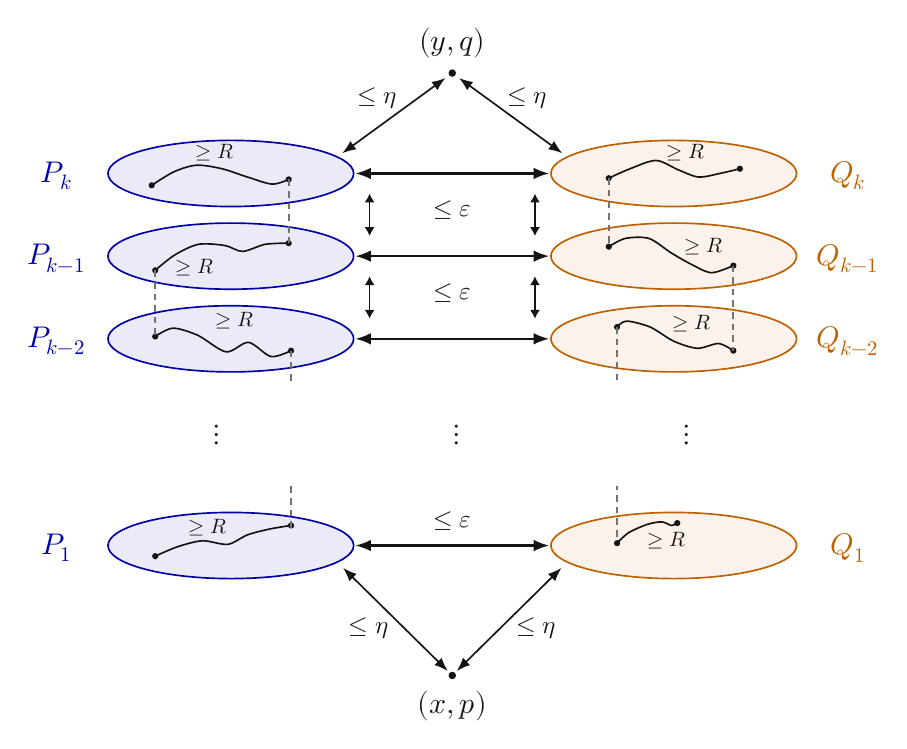}
  \caption{Configuration from Theorem~\ref{thm:distancebetweenprojections}}
\end{figure}

The last terminology we need is that of the~\textit{altitude} of a halo product. Intuitively, it corresponds to the minimal number of ladders we need to connect any point of $\halo H$ to a point of the coset $H\subset \halo H$. 

\begin{definition}
Let $\halo H$ be a finitely generated halo product. Given $\varepsilon, \eta, R>0$, the $(\varepsilon, \eta, R)-$\textit{altitude} of $\halo$ is the minimal number $N\in\N\cup\lbrace\infty\rbrace$ such that, for any $a\in\halo H$, there exist $N+1$ points $x_{0}=a,x_{1},\dots,x_{N}$ successively $\eta-$connected by an $(\varepsilon,R)-$ladder with $x_{N}\in H$.
\end{definition}

The main result of~\cite[Section~6]{GT24a} states then that, between large-scale commutative halo products having finite altitude, leaf-preserving quasi-isometries are at bounded distance from aptolic quasi-isometries. 

\begin{theorem}[{\cite[Corollary~6.11]{GT24a}}]\label{thm:boundeddistancefromaptolicQI}
Let $\mathcal{M}A$, $\mathcal{N}B$ be two finitely generated large-scale commutative halo groups. For all $C,K,\varepsilon, \varepsilon'>0$, there exists $Q>0$ such that, if $\mathcal{M}$ has finite $(\varepsilon, \varepsilon',Q)-$altitude, then every leaf-preserving $(C,K)-$quasi-isometry $\mathcal{M}A\longrightarrow \mathcal{N}B$ lies at bounded distance from an aptolic quasi-isometry.   
\end{theorem}

\begin{proof}
Let $L,S,\eta\ge 0$ be the constants given by Proposition~\ref{prop:QIgraphofleaves}, and let $D$ be the large-scale commutativity constant for $\mathcal{N}$. Fix any $Q>0$ such that $Q>L(D+4\eta+2r_{0}+S)$. Let $\varphi\colon \mathcal{M}A\longrightarrow\mathcal{N}B$ be a leaf-preserving $(C,K)-$quasi-isometry, and write it as 
\begin{equation*}
    \varphi(c,p)=(\alpha(c), \beta_{c}(p))
\end{equation*}
for some bijection $\alpha\colon M(A)\longrightarrow N(B)$ and $(C,K)-$quasi-isometries $\beta_{c}\colon A\rightarrow B$, $c\in M(A)$.

\noindent Let $(c,p)\in\mathcal{M}A$. Since $\mathcal{M}$ has finite $(\varepsilon, \varepsilon',Q)-$altitude $N$, there exist $k\le N+1$ points $z_{0}=(c,p),z_{1},\dots,z_{k}\defeq (1,q)\in A$ such that $z_{i}$ and $z_{i+1}$ are $\varepsilon'-$connected by an $(\varepsilon, Q)-$ladder of leaves, for every $0\le i\le k-1$. Applying Theorem~\ref{thm:distancebetweenprojections} $k$ times, it follows that 
\begin{equation*}
    d(p,q) \le k(6\varepsilon'+\varepsilon).
\end{equation*}
Using Proposition~\ref{prop:QIgraphofleaves}, we also know that, for any $0\le i\le k-1$, $\varphi(z_{i})$ and $\varphi(z_{i+1})$ are $(C\varepsilon'+K)-$connected by a $(\eta, \frac{Q}{L}-S)-$ladder of leaves in $\mathcal{N}B$. Thus, applying once again Theorem~\ref{thm:distancebetweenprojections} $k$ times, we deduce that 
\begin{equation*}
    d\big(\beta_{c}(p), \beta_{\mathbf{1}}(q)\big) \le k(6C\varepsilon'+6K+\eta)
\end{equation*}
and thus it follows that 
\begin{align*}
    d\big(\beta_{c}(p), \beta_{\mathbf{1}}(p)\big) &\le d\big(\beta_{c}(p), \beta_{\mathbf{1}}(q)\big)+d\big(\beta_{\mathbf{1}}(q), \beta_{\mathbf{1}}(p)\big) \\
    &\le k\big(6C\varepsilon'+6K+\eta\big)+C\cdot d(p,q)+K \\
    &\le k\big(6C\varepsilon'+6K+\eta\big)+C\cdot k(6\varepsilon'+\varepsilon)+K.
\end{align*}
Hence any $\beta_{c}$ lies at bounded distance from $\beta_{\mathbf{1}}$, so we conclude that $\varphi$ lies at bounded distance from the aptolic quasi-isometry $(c,p)\longmapsto (\alpha(c),\beta_{\mathbf{1}}(p))$. 
\end{proof}

Thus, to apply Theorem~\ref{thm:boundeddistancefromaptolicQI} in a concrete example, it is enough to check that the halo product we consider has finite altitude. This is done in~\cite[Lemmas~6.12-6.16]{GT24a} for wreath products, lampjugglers, lampdesigners and lampcloners. 

\paragraph{Growth of lamps.} From here, we therefore know that any quasi-isometry between two halo products can be chosen aptolic, under appropriate assumptions on the base groups and the halos. To derive information and notably scaling properties for quasi-isometries between the base groups, the idea is, roughly speaking, to look at “how big” a coset of $L(S)$ in $\halo H$ (for $S\subset H$) is, and then, to quantify the speed at which $L(S)$ grows when $S$ gets bigger and bigger. When our halo product $\halo H$ has its subgroup $L(H)$ locally finite, it suffices to consider cardinalities of cosets of $L(S)$. To introduce the corresponding function, we require additional terminology.

\begin{definition}\label{def:consistencyforhalos}
Let $H$ be a group, and let $\halo H$ be a halo product over $H$. We say that $\halo H$ has~\textit{finite blocks} (resp.~\textit{finitely generated blocks}) if $L(S)$ is finite for any $S\subset H$ (resp. finitely generated). If, furthermore, $|L(S)|$ only depends on $|S|$, we say that $\halo H$ is~\textit{consistent}.
\end{definition}

In practice, consistent halo products encompass all classes we are interested in, namely lamplighters, lampjugglers and lampcloners.

\begin{definition}\label{def:Growthoflamps}
Let $\halo H$ be a consistent halo product over $H$. The~\textit{lamp growth sequence} of $\halo H$ is the function $\Lambda_{\halo H}\colon \N\rightarrow \N$ defined by 
\begin{equation*}
    \Lambda_{\halo H}\colon n\longmapsto |L(S)|, \;\text{where}\; |S|=n.
\end{equation*}
\end{definition}

For instance, if $F$ is a finite group, the lamp growth sequence of $F\wr H$ is $\Lambda_{F\wr H}(n)=|F|^n$, while for lampjugglers it is 
\begin{equation*}
    \Lambda_{\juggler{r}{H}}(n)=(rn)!. 
\end{equation*}
Other examples are contained in~\cite[Facts~7.12-7.16]{GT24a}. 

In our strategy, growth of lamps is relevant because it turns out that, when looking at an aptolic quasi-isometry $\mathcal{M}A\longrightarrow \mathcal{N}B$ of the form $(c,p)\longmapsto (\alpha(c), \beta(p))$, the map $\alpha$ must actually send a coset of a subgroup of the form $M(S)\subset M(A)$ to a coset of a subgroup of the form $N(\Tilde{S}) \subset N(B)$, where $\Tilde{S}$ is constructed from $S$ using $\beta$ and the parameters $C,K$ of the quasi-isometry. This observation is proved in~\cite[Lemma~3.5]{GT24b} for lamplighters and in~\cite[Lemma~7.4]{GT24a} for general halo products. 

This observation then allows to show that the behaviour of the lamp growth sequence is invariant when applying aptolic quasi-isometries. For the statement, recall that a sequence $u\colon \N\rightarrow\N$ \textit{dominates} a sequence $v\colon\N\rightarrow \N$, written $v\prec u$, if there is $C>0$ such that $v(n) \le u(Cn)$, and $u\sim v$ if $u\prec v$ and $v\prec u$. 

\begin{corollary}[{\cite[Corollary~7.11]{GT24a}}]
Let $\mathcal{M}A$, $\mathcal{N}B$ be two finitely generated consistent halo products with $M(A), N(B)$ locally finite. If there exists an aptolic quasi-isometry $\mathcal{M}A\longrightarrow \mathcal{N}B$, then $\Lambda_{\mathcal{M}A} \sim \Lambda_{\mathcal{N}B}$. 
\end{corollary}

The invariance of the asymptotic behaviour of the lamp growth sequence is already powerful enough to distinguish geometrically halos of different nature. 

\begin{corollary}
Let $E$, $F$ be non-trivial finite groups. Let $G,H,I,K$ be finitely generated groups with the thick bigon property. Let $\field$ be a finite field and let $s\ge 1$ be an integer. Then:
\begin{enumerate}[label=(\roman*)]
    \item The lamplighter $E\wr G$ is not quasi-isometric to $\juggler{s}{H}$, $\designer{I}$ or $\cloner{K}$.
    \item The lampjuggler $\juggler{s}{H}$ and the lampdesigner $\designer{I}$ are not quasi-isometric to the lampcloner $\cloner{K}$. 
\end{enumerate}
\end{corollary}

On the other hand, the equivalence class of the lamp growth sequence under the relation $\sim$ does not allow to distinguish halos of the same kind, for instance two lamplighters or two lampjugglers. We must compare more subtly the behaviour of these two sequences. To this end,~\cite{GT24a} introduces the notion of \textit{interlaced} sequences. 

\begin{definition}
Let $\Delta\colon \N\rightarrow \N$. Given two sequences $u,v\colon \N\rightarrow\N$, we say that they are $\Delta-$\textit{interlaced} if there exist $C\ge 0$ and a sequence $(x_{n})_{n\in\N}$ such that $x_{n}=\Theta(n)$ and such that 
\begin{equation*}
    u(n) \;\text{divides}\; v(x_{n}),\;\text{which divides}\; u(n+C\cdot \Delta(n))
\end{equation*}
for any $n\ge 1$. 
\end{definition}

The function $\Delta$ to consider when working with lamp growth sequences of halo groups is the \textit{boundary growth} function $\Delta=\mathbb{B}$.

\begin{definition}
If $X$ is a locally finite graph, its~\textit{boundary growth function} is the function $\mathbb{B}_{X}\colon\N\rightarrow \N$ defined by 
\begin{equation*}
    \mathbb{B}_{X}(n) \defeq \min\lbrace |\partial S| : S\subset X\;\text{connected of size}\; n\rbrace. 
\end{equation*}
\end{definition}

The presence of an aptolic quasi-isometry between two halo products then forces lamp growth sequences to be $\mathbb{B}-$interlaced. 

\begin{proposition}[{\cite[Proposition~7.22]{GT24a}}]\label{prop:interlaced}
Let $\mathcal{M}A$, $\mathcal{N}B$ be finitely generated consistent halo products with $M(A), N(B)$ locally finite. If there exists an aptolic quasi-isometry $\mathcal{M}A\longrightarrow \mathcal{N}B$, then $\Lambda_{\mathcal{M}A}$ and $\Lambda_{\mathcal{N}B}$ are $\mathbb{B}_{A}-$interlaced. 
\end{proposition}

As an application, this property allows to detect an arithmetic condition on the cardinalities of the lamp groups when two lamplighters are quasi-isometric. The next statement generalises the one already proved in~\cite{GT24b}. 

\begin{proposition}[{\cite[Proposition~7.23]{GT24a}}]\label{prop:lamplightersoverTBPgroups-rigiditypart}
Let $E$, $F$ be non-trivial finite groups. Let $G$ and $H$ be finitely generated groups satisfying the thick bigon property. If $E\wr G$ and $F\wr H$ are quasi-isometric, then $|E|$ and $|F|$ have the same prime divisors. 
\end{proposition}

\begin{proof}
By Theorem~\ref{thm:embedding2}, our quasi-isometry $\varphi$ quasi-preserves leaves, hence by Theorem~\ref{thm:boundeddistancefromaptolicQI}, it lies at bounded distance from an aptolic quasi-isometry. It thus follows from Proposition~\ref{prop:interlaced} that $\Lambda_{E\wr G}$ and $\Lambda_{F\wr H}$ are $\mathbb{B}_{G}-$interlaced. In particular, there is a sequence $(x_{n})_{n\in\N}$ such that $\Lambda_{E\wr G}(n)=|E|^n$ divides $\Lambda_{F\wr H}(x_{n})=|F|^{x_{n}}$, and thus any prime divisor of $|E|$ also divides $|F|$. The same reasoning applied to a quasi-inverse of $\varphi$ shows that a prime divisor of $|F|$ also divides $|E|$. This concludes the proof. 
\end{proof}

For lampcloners, one can deduce from Proposition~\ref{prop:interlaced} that the characteristic of the underlying field is a quasi-isometry invariant. 

\begin{proposition}[{\cite[Proposition~7.24]{GT24a}}]
Let $\field$, $\field'$ be two finite fields. Let $G$ and $H$ be finitely generated groups satisfying the thick bigon property. If $\cloner{G}$ and $\mathsf{Cloner}_{\field'}(H)$ are quasi-isometric, then $\field$ and $\field'$ have the same characteristic. 
\end{proposition}

\paragraph{Amenability adds rigidity.} All previous statements hold in a general situation, regardless of (non-)amenability of the groups. However, it turns out that, as already highlighted in~\cite{GT24b}, stronger rigidity statements hold if we additionally know that our groups are amenable. Indeed, thanks to aptolicity, we know that a quasi-isometry between halo products implies the existence of a quasi-isometry between the base groups, and thanks to amenability, this induced quasi-isometry must be measure-scaling. This statement is proved once again by looking at the asymptotic behaviour of lamp growth sequences. 

\begin{definition}
Let $u,v\colon \N\rightarrow \N$ be two non-decreasing sequences. Let $r>0$ be a real number. We write $u \lhd_{r} v$ if there exists $t\in\N$ such that 
\begin{equation*}
    u(k) \le v(\lfloor rk\rfloor+t)
\end{equation*}
for all sufficiently large $k\in\N$. We write $u\bowtie_{r} v$ if $u\lhd_{r} v$ and $v\lhd_{\frac{1}{r}}u$. 
\end{definition}

Then, it is proved in~\cite{GT24a} that such a relation between lamp growth sequences of our halo groups provides scaling properties.

\begin{proposition}[{\cite[Proposition~8.3]{GT24a}}]\label{prop:bowtieimpliesscaling}
Let $\mathcal{M}A$, $\mathcal{N}B$ be finitely generated consistent halo products, with $M(A)$, $N(B)$ locally finite. Let $\alpha\colon M(A)\longrightarrow N(B)$ be a bijection and $\beta\colon A\rightarrow B$ be a quasi-isometry such that
\begin{align*}
    \mathcal{M}A&\longrightarrow \mathcal{N}B \\
    (c,p)&\longmapsto (\alpha(c), \beta(p))
\end{align*}
is a quasi-isometry. If $\Lambda_{\mathcal{N}B} \bowtie_{k}\Lambda_{\mathcal{M}A}$ for some $k>0$, then $\beta$ is quasi-$k$-to-one. 
\end{proposition}

As an illustration, for lamplighters, Theorem~\ref{thm:classificationGT21} can be strengthened as follows. 

\begin{theorem}[{\cite[Theorem~8.6]{GT24a}}]\label{thm:AlamplightersoverTBPgroups-rigiditypart}
Let $E$, $F$ be non-trivial finite groups. Let $G$ and $H$ be finitely generated amenable groups satisfying the thick bigon property. If $E\wr G$ and $F\wr H$ are quasi-isometric, then there exist $a,r,s\ge 1$ such that $|E|=a^{r}$ and $|F|=a^{s}$, and there exists a quasi-$\frac{s}{r}$-to-one quasi-isometry $G\rightarrow H$. 
\end{theorem}

\begin{proof}
Fix a quasi-isometry $\varphi\colon E\wr G\longrightarrow F\wr H$. By Theorem~\ref{thm:embedding2}, $\varphi$ is leaf-preserving, and by Theorem~\ref{thm:boundeddistancefromaptolicQI}, it can be taken aptolic, and we write it as 
\begin{equation*}
    \varphi(c,p) = (\alpha(c), \beta(p)), \; (c,p)\in E\wr G. 
\end{equation*}
Applying Proposition~\ref{prop:interlaced}, it follows that there exist sequences $(x_{n})_{n\in\N}$, $(y_{n})_{n\in\N}$ such that $x_{n}=\Theta(y_{n})$ and $\Lambda_{E\wr G}(x_{n})=|E|^{x_{n}}$ divides $\Lambda_{F\wr H}(y_{n})=|F|^{y_{n}}$, which itself divides $\Lambda_{E\wr G}(y_{n}+o(n))=|E|^{y_{n}+o(n)}$. Thus, if $p\ge 2$ is prime, we get that 
\begin{equation*}
    x_{n}\cdot \text{val}_{p}(|E|) \le y_{n}\cdot\text{val}_{p}(|F|) \le (y_{n}+o(n))\cdot \text{val}_{p}(|E|)
\end{equation*}
whence $\text{val}_{p}(|E|)=0$ if and only if $\text{val}_{p}(|F|)=0$, i.e. $p$ divides $|E|$ if and only if $p$ divides $|F|$. Additionally, the sequence $(\frac{x_{n}}{y_{n}})_{n\in\N}$ converges to $\frac{\text{val}_{p}(|E|)}{\text{val}_{p}(|F|)}$ for every prime $p$ dividing $|E|$ and $|F|$. Thus, the latter cannot depend on $p$, and there exist $r,s\ge 1$ such that 
\begin{equation*}
    \frac{\text{val}_{p}(|E|)}{\text{val}_{p}(|F|)} = \frac{r}{s}
\end{equation*}
for any prime $p$ dividing $|E|$ and $|F|$. Thus there exists $a\ge 1$ such that $|E|=a^{r}$ and $|F|=a^{s}$.

\noindent To conclude, it remains to notice that
\begin{equation*}
    \Lambda_{F\wr H}(n)=|F|^n=a^{sn} \le a^{r(\lfloor \frac{s}{r}n\rfloor+1)} \le \Lambda_{E\wr G}\left(\left\lfloor \frac{s}{r}n\right\rfloor+1\right)
\end{equation*}
for any $n\in\N$, thus $\Lambda_{F\wr H} \lhd_{\frac{s}{r}}\Lambda_{E\wr G}$. Similarly, we prove that $\Lambda_{E\wr G} \lhd_{\frac{r}{s}}\Lambda_{F\wr H}$. Thus $\Lambda_{F\wr H} \bowtie_{\frac{s}{r}}\Lambda_{E\wr G}$, and it follows from Proposition~\ref{prop:bowtieimpliesscaling} that $\beta$ is quasi-$\frac{s}{r}$-to-one. 
\end{proof}

The authors of~\cite{GT24a} derived similar statements for lampjugglers, lampdesigners and lampcloners. For instance: 

\begin{theorem}[{\cite[Theorem~8.7]{GT24a}}]\label{thm:QIrigidityofjugglers}
Let $H$ and $K$ be two finitely generated amenable groups satisfying the thick bigon property. Let $r,s\ge 1$ be two integers. Then $\juggler{r}{H}$ and $\juggler{s}{K}$ are quasi-isometric if and only if there exists a quasi-$\frac{s}{r}$-to-one quasi-isometry $H\rightarrow K$. 
\end{theorem}

As a consequence:

\begin{corollary}[{\cite[Corollary~8.9]{GT24a}}]\label{cor:BILIPbetweenLampshufflers}
Let $H$ and $K$ be finitely generated amenable groups satisfying the thick bigon property. Then $\shuf{H}$ and $\shuf{K}$ are quasi-isometric if and only if $H$ and $K$ are biLipschitz equivalent. 
\end{corollary}

\section{Scaling groups of halo products}\label{sec:Scalinggroupsofhaloproducts} From Genevois and Tessera's results about the quasi-isometric classification of halo products and their method of proof explained above, we can deduce the following:

\begin{corollary}\label{cor:scalinggroupsofhaloproducts}
Let $H$ and $K$ be finitely generated amenable groups with the thick bigon property. If $\halo$ is either:
\begin{itemize}
    \item $\halo(\cdot)=F\wr(\cdot)$ where $F$ is a non-trivial finite group;
    \item $\halo(\cdot)=\juggler{s}{\cdot}$, where $s\ge 1$;
    \item $\halo(\cdot)=\designer{\cdot}$ where $F$ is a non-trivial finite group;
    \item $\halo(\cdot)=\cloner{\cdot}$ where $\field$ is a finite field;
\end{itemize}
then any quasi-isometry $\halo H\longrightarrow \halo K$ is measure-scaling. In particular, $\text{Sc}(\halo H)=\lbrace 1\rbrace$.
\end{corollary}

The proof of the corollary relies on the following intermediate observation, proved in~\cite[Lemma~3.10]{GT24b} for lamplighters over one-ended groups and mentioned in~\cite[Section~9]{GT24a}. 

\begin{lemma}\label{lem:scalingQIbetweenhalos}
Let $H$ and $K$ be finitely generated groups, and let $\halo H$ and $\halo K$ be two halo products over $H$ and $K$. Let $\alpha\colon L(H)\longrightarrow L(K)$ and $\beta\colon H\rightarrow K$ be two maps such that 
\begin{align*}
    \varphi\colon \halo H&\longrightarrow \halo K \\
    (c,p)&\longmapsto (\alpha(c),\beta(p))
\end{align*}
is a quasi-isometry. If $\beta$ is quasi-$k$-to-one for some $k>0$, then $\varphi$ is quasi-$k$-to-one. 
\end{lemma}

\begin{proof}
Denote by $\pi_{H}$ (resp. $\pi_{K}$) the projection of $\halo H$ (resp. $\halo K$) onto the second factor $H$ (resp. $K$). Fix a finite subset $A\subset \halo K$. Let $\mathcal{C}\subset L(K)$ be the collection of all elements of $L(K)$ appearing as first coordinate of an element of $A$, and for any $c\in\mathcal{C}$, let $A_{c}\subset A$ be the subset of all elements of $A$ having $c$ as first coordinate. We then have
\begin{equation*}
    A= \bigsqcup_{c\in\mathcal{C}}A_{c}
\end{equation*}
and then $|A|=\sum_{c\in \mathcal{C}}|A_{c}|$, $|\varphi^{-1}(A)|=\sum_{c\in \mathcal{C}}|\varphi^{-1}(A_{c})|$. Now, if we set $B_{c}\defeq \pi_{K}(A_{c})$ for any $c\in\mathcal{C}$, we have $|A_{c}|=|B_{c}|$ and $|\varphi^{-1}(A_{c})|=|\beta^{-1}(B_{c})|$ since $\pi_{K}\circ \varphi = \beta\circ \pi_{H}$. We then deduce that 
\begin{align*}
    \left|k|A|-|\varphi^{-1}(A)|\right| &\le \sum_{c\in\mathcal{C}}\left|k|A_{c}|-|\varphi^{-1}(A_{c})|\right| \\
    &=\sum_{c\in\mathcal{C}}\left|k|B_{c}|-|\beta^{-1}(B_{c})|\right| \\
    &\le C\cdot \sum_{c\in\mathcal{C}}|\partial_{K}B_{c}|
\end{align*}
where $C>0$ is the constant coming from the quasi-$k$-to-one estimate of $\beta$. It now remains to notice that $\partial_{\halo K}A$ contains $\bigsqcup_{c\in\mathcal{C}}\lbrace c\rbrace\times \partial_{K}B_{c}$ to conclude that 
\begin{equation*}
    \left|k|A|-|\varphi^{-1}(A)|\right| \le C\cdot |\partial_{\halo K}A|
\end{equation*}
as was to be proved. Thus $\varphi$ is quasi-$k$-to-one. 
\end{proof}

\begin{proof}[Proof of Corollary~\ref{cor:scalinggroupsofhaloproducts}]
As explained above, it is proved in~\cite[Corollary~6.11 and Lemmas~6.12-6.16]{GT24a} that any quasi-isometry between $\halo H$ and $\halo K$ is, under our assumptions, at bounded distance from an aptolic quasi-isometry $\varphi\colon \halo H\longrightarrow \halo K$, $\varphi(c,p)=(\alpha(c),\beta(p))$. Since quasi-isometries that are at bounded distance from scaling quasi-isometries are themselves scaling by Proposition~\ref{prop:stabilitypropertiesforscalingQI}, and since aptolic quasi-isometries are scaling if their second components are scaling by Lemma~\ref{lem:scalingQIbetweenhalos}, it only remains to see that $\beta\colon H\rightarrow K$ is measure-scaling. This fact is proved in~\cite[Section~8.1]{GT24a} for lamplighters, in~\cite[Section~8.2]{GT24a} for lampshufflers and lampjugglers, in~\cite[Section~8.3]{GT24a} for lampdesigners and in~\cite[Section~8.4]{GT24a} for lampcloners.
\end{proof}

\afterpage{\blankpage}

\chapter{On the quasi-isometric classification of permutational wreath products}\label{chap:chapter3}

This chapter presents the content of~\cite{Dum24}.

\vspace{0.3cm}

\minitoc

\section{Introduction}\label{subsection3.1}

\setcounter{theoremletter}{0}

\sloppy In his groundbreaking work in the 80's and 90's~\cite{Gro81, Gro93}, Gromov initiated the program of classifying finitely generated groups up to quasi-isometries. The motivation for this program is that the large-scale geometry of such a group may be in fact deeply related to its algebraic structure, as shown for instance by Gromov's celebrated theorem on groups of polynomial growth~\cite{Gro81} or Stallings' theorem about multi-ended groups~\cite{Sta68}. 

\sloppy Much progress has been made since then towards understanding the quasi-isometries of various classes of groups, and some groups are even known to be~\textit{quasi-isometrically rigid}. Here, a group $G$ is quasi-isometrically rigid if any group which is quasi-isometric to $G$ is in fact isomorphic to a finite-index subgroup of $G$, to a quotient of $G$ by a finite normal subgroup, or contains $G$ as a finite-index subgroup. In this case as well, the connection between the geometry of the group and its algebraic structure is particularly strong and explicit. For instance, abelian and non-abelian free groups are quasi-isometrically rigid~\cite{Dun85}, as well as mapping class groups~\cite{Beh+12}, or solvable Baumslag-Solitar groups $\text{BS}(1,n), n\ge 2$~\cite{FM98}.

On the other hand, there are also many known invariants that allow to distinguish groups up to quasi-isometries, including the volume growth, amenability, hyperbolicity, the number of ends, or being finitely presentable. However, for many groups of interest, these invariants are insufficient, and other methods must be developed. 

One such class, whose quasi-isometries are hard to tackle, consists of~\textit{permutational wreath products}. Recall that given two groups $G,H$ and an action of $H$ on a set $X$, the permutational wreath product $G\wr_{X}H$ is defined as 
\begin{equation*}
    G\wr_{X}H \defeq \left(\bigoplus_{X}G\right)\rtimes H
\end{equation*}
where $H$ acts on the direct sum permuting the coordinates through its initial action on $X$. In the particular case $X=H$ and $H$ acts on itself by left-multiplication, we only write $G\wr H$ and we refer to this group as a~\textit{(standard) wreath product}, or a~\textit{lamplighter group} when $G$ is finite. Such groups are well-known and of interest in group theory for several reasons, the main one being that on the one hand the explicit definition makes computations possible, and on the other hand the definition is sufficiently elaborate to reflect unexpected and interesting behaviours. They have been studied extensively in relation with many topics of interests in geometric group theory, for instance amenability and isoperimetric profiles~\cite{Ers03, MO10}, Haagerup property~\cite{CSV08, CSV12}, random walks~\cite{LPP96, PSC02, BE17}, subgroup distortion~\cite{DO11, BLP15, Ril22}, bounded cohomology~\cite{Mon22}, fixed-point properties~\cite{CK11, LS22, LS24} or coarse embeddability into Hilbert spaces and the related compressions~\cite{Li10, Gen22, AT19}.

For the large-scale geometric analysis of lamplighter groups, a first question that arises naturally is then:
\begin{question}\label{question1.1}
Let $F_{1}, F_{2}$ be finite groups and let $H_{1}, H_{2}$ be finitely generated groups. When are $F_{1}\wr H_{1}$ and $F_{2}\wr H_{2}$ quasi-isometric?
\end{question}

An important first piece of answer have been brought by Eskin-Fisher-Whyte in~\cite{EFW12, EFW13}, where they prove that if $H_{1}$ is two-ended, then $F_{1}\wr H_{1}$ and $F_{2}\wr H_{2}$ are quasi-isometric if and only if $H_{1}, H_{2}$ are quasi-isometric and $|F_{1}|, |F_{2}|$ are powers of a common number. However, the proof of this result heavily relies on the two-endedness assumption, and does not seem to be replicable outside of this field. 
In a recent work~\cite{GT24b}, Genevois and Tessera obtained a complete classification of lamplighters over one-ended finitely presented groups. More precisely:
\begin{theorem}\label{thm1.2}
Let $F_{1}, F_{2}$ be non-trivial finite groups and $H_{1}, H_{2}$ be two finitely presented groups. Suppose that $H_{1}$ is one-ended. 
\begin{itemize}
    \item If $H_{1}$ is amenable, then $F_{1}\wr H_{1}$ and $F_{2}\wr H_{2}$ are quasi-isometric if and only if there exist $k,n_{1},n_{2}\ge 1$ such that $|F_{1}|=k^{n_{1}}$, $|F_{2}|=k^{n_{2}}$ and there exists a quasi-$\frac{n_{2}}{n_{1}}$-to-one quasi-isometry $H_{1}\rightarrow H_{2}$.
    \item If $H_{1}$ is not amenable, then $F_{1}\wr H_{1}$ and $F_{2}\wr H_{2}$ are quasi-isometric if and only if $H_{1}$ and $H_{2}$ are quasi-isometric and $|F_{1}|,|F_{2}|$ have the same prime divisors. 
\end{itemize}
\end{theorem}

We refer the reader to Definition~\ref{def:measurescalingQI} for the definition of quasi-$k$-to-one quasi-isometries and the related scaling groups $\text{Sc}(\cdot)$.

The strategy employed for the proof of Theorem~\ref{thm1.2} in~\cite{GT24b} is completely different from that of~\cite{EFW12, EFW13}, and relies crucially on quasi-median geometry and the assumption that groups over which we consider wreath products are finitely presented and one-ended. 

More recently, Genevois and Tessera showed in~\cite{GT24a} that their techniques can also be applied to a wider class of groups generalizing lamplighter groups, and exhibiting a common~\textit{halo structure}. Such a class encompasses standard wreath products, but also lampshufflers, lampdesigners, nilpotent wreath products and many other variations. This allows to extend, at least partially, the classification provided by Theorem~\ref{thm1.2}. In another work with Bensaid~\cite{BGT24}, they also extended their methods to some standard lamplighters with infinite lamp groups, putting into the picture the notion of “coarse separation” that, roughly speaking, generalizes (non-)one-endedness. 

However, one classical variation of standard wreath products remains out of the field covered until now, that of permutational wreath products with infinite stabilizers (the finite stabilizers case is contained in~\cite{GT24b} as a consequence of Theorem~\ref{thm1.2}). Indeed, methods from~\cite{GT24a} do not apply to these groups, according to~\cite[Corollary~4.26]{GT24a}. Our main result is a complete classification of some of these permutational wreath products up to quasi-isometry:

\begin{theoremletter}\label{thm:classificationPWPuptoQI}
Let $E$, $F$ be two non-trivial finite groups. Let $G$, $H$ be finitely presented groups, with finitely generated normal infinite subgroups $M\lhd G$, $N\lhd H$. Suppose that $M$ has infinite index in $G$, and that $G$ (resp. $H$) is not coarsely separable by any collection of subspaces that uniformly quasi-isometrically embed into $M$ (resp. $N$). The following claims hold.
\begin{itemize}
    \item If $M$ is co-amenable in $G$, then $E\wr_{G/M}G$ and $F\wr_{H/N}H$ are quasi-isometric if and only if \;$|E|=n^{r}$, $|F|=n^{s}$ for some $n,r,s\ge 1$ and there exists a quasi-isometry of pairs $(G,M)\longrightarrow (H,N)$ inducing a quasi-$\frac{s}{r}$-to-one quasi-isometry $G/M\rightarrow H/N$.
    %\vspace{0.05cm}
    \item If $M$ is not co-amenable in $G$, then $E\wr_{G/M}G$ and $F\wr_{H/N}H$ are quasi-isometric if and only if \;$|E|$ and $|F|$ have the same prime divisors and there exists a quasi-isometry of pairs $(G,M)\longrightarrow (H,N)$. 
\end{itemize}
\end{theoremletter}

Here, a quasi-isometry $f\colon G\rightarrow H$ is a~\textit{quasi-isometry of pairs} if it sends any $M-$coset at finite Hausdorff distance from an $N-$coset, and we write $f\colon (G,M)\longrightarrow (H,N)$. Moreover, such a map always induces a quasi-isometry between the quotients $G/M$ and $H/N$, unique up to bounded distance. See Definition~\ref{def:QIofpairs} and Proposition~\ref{prop:inducedmapbetweenquotients} for more details.

Recall also that a subgroup $M\lhd G$ is~\textit{co-amenable in} $G$ if the quotient group $G/M$ is amenable. See Section~\ref{subsection:amenability} for further details.

Next, notice that in this statement, the assumption on the index of $M$ in $G$ is not restrictive. Indeed, if $M$ has finite index in $G$, then $E\wr_{G/M}G$ is quasi-isometric to $G$, so $E\wr_{G/M}G$ and $F\wr_{H/N}H$ are quasi-isometric if and only if $G$ and $H$ are. 

The other assumptions in the statement of Theorem~\ref{thm:classificationPWPuptoQI} seem to be, at least at first glance, very restrictive. In particular, the systematic study of coarse separation of spaces, although the notion already appeared at various places in the literature under different formulations (see e.g.~\cite[Section~9.7]{DK18},~\cite{FS96, Mar21}), seems to have been initiated only very recently in~\cite{BGT24}. Moreover, the latter focuses essentially on coarse separation by subspaces of subexponential growth, and the global picture remains to be established. We refer to~\cite[Section~7]{BGT24} for a deeper discussion and several open questions on the subject, and to~\cite{BGT26a, BGT26b} for the relation between coarse separation and splittings in hyperbolic and right-angled Artin groups. 

However, assumptions of Theorem~\ref{thm:classificationPWPuptoQI} are satisfied in a number of classical cases, and combining it with results from~\cite{GT24a, BGT24} already allows us to derive many applications. The following example is the one that motivated the search for a general criterion.

\begin{corollaryletter}\label{cor:classificationofPWPoverZ^duptoQI}
Let $E$ and $F$ be two non-trivial finite groups, and let $m,m',n,n'$ be four integers such that $m\ge n\ge 2$, $m'\ge n'\ge 2$. Then $E\wr_{\Z^n}\Z^m$ and $F\wr_{\Z^{n'}}\Z^{m'}$ are quasi-isometric if and only if $m=m'$, $n=n'$ and $|E|$, $|F|$ are powers of a common number. 
\end{corollaryletter}

For this result (resp. Corollary~\ref{cor:classificationPWPoverZ^duptoBILIP} and Proposition~\ref{prop:classificationAlamplightersoverPWPoverZ^d} below), we refer to Section~\ref{subsection3.7}, more precisely to Question~\ref{question7.6}, for a discussion on the assumption $n,n'\ge 2$ (resp. $n\ge 2$, $k\ge 2$).

In fact, the proof of this corollary shows that we can reformulate Theorem~\ref{thm:classificationPWPuptoQI} in the case of permutational wreath products over base groups that are direct products.

\begin{corollaryletter}\label{cor:classificationofPWPoverdirectproducts}
Let $E$ and $F$ be two non-trivial finite groups, and let $M,K$ be infinite finitely presented groups. Suppose that $K$ is amenable, and that $G\defeq M\times K$ is not coarsely separable by any collection of subspaces that uniformly quasi-isometrically embed into $M$. Then $E\wr_{K}G$ and $F\wr_{K}G$ are quasi-isometric if and only if there exist $n,r,s\ge 1$ such that $|E|=n^{r}$, $|F|=n^{s}$ and $\frac{s}{r}\in \text{Sc}(K)$.
\end{corollaryletter}

Here also, we refer to Definition~\ref{def:scalinggroup} for the definition of scaling groups of finitely generated groups.

Hence one can classify up to quasi-isometry many permutational wreath products over a direct product for which the scaling group is known, such as free abelian groups, solvable Baumslag-Solitar groups or $\text{SOL}(\Z)$ (see Proposition~\ref{prop:examplesofscalinggroups}).

As another consequence of Theorem~\ref{thm:classificationPWPuptoQI}, we may also classify permutational wreath products over a group which is not necessarily a direct product, if the quotient has trivial scaling group:

\begin{corollaryletter}\label{cor:PWPwithquotientstrivialscalinggroup}
Let $E$ and $F$ be two non-trivial finite groups. Let $G$ be a finitely presented group with a finitely generated normal infinite subgroup $M\lhd G$ of infinite index. Suppose that $M$ is co-amenable in $G$, and that $G$ is not coarsely separable by any collection of subspaces that uniformly quasi-isometrically embed into $M$. If $\text{Sc}(G/M)=\lbrace 1\rbrace$, then $E\wr_{G/M}G$ and $F\wr_{G/M}G$ are quasi-isometric if and only if $|E|=|F|$. 
\end{corollaryletter}

On the other hand, Theorem~\ref{thm:classificationPWPuptoQI} also applies to many permutational wreath products built over nilpotent groups. Indeed, such a group $G$ is automatically finitely presented when finitely generated, has all its subgroups finitely generated, and is not coarsely separated by collections of subspaces having growth degree $\le \text{deg}(G)-2$~\cite[Theorem~1.5]{BGT24}, where $\text{deg}(G)$ is the growth degree of $G$. For instance, if $H$ is the Heisenberg group over the integers, and $Z(H)$ denotes its centre, then $\Z_{n}\wr_{H/Z(H)}H$ and $\Z_{m}\wr_{H/Z(H)}H$ are quasi-isometric if and only if $n$ and $m$ are powers of a common number. 

We also emphasize that assumptions of Theorem~\ref{thm:classificationPWPuptoQI} are not necessary for proving both directions of the equivalences. In fact, when constructing quasi-isometries of permutational wreath products from the data of a quasi-isometry of pairs between the base groups, we can get rid of the finite generation and coarse separation assumptions on the subgroups, as well as the finite presentation on the groups. Especially, in the case of non-co-amenable subgroups, we obtain a wide range of quasi-isometric permutational wreath products:

\begin{propositionletter}\label{prop:nonQIPWP}
Let $n,m\ge 2$ be two integers, and let $G$ and $H$ be finitely generated groups with normal subgroups $M\lhd G$, $N\lhd H$. Assume that $M$ is not co-amenable in $G$. If there exists a quasi-isometry of pairs $(G,M)\longrightarrow (H,N)$ and if $n$ and $m$ have the same prime divisors, then there exists a quasi-isometry
\begin{equation*}
    \Z_{n}\wr_{G/M}G\longrightarrow \Z_{m}\wr_{H/N}H.
\end{equation*}
\end{propositionletter}

In addition, in~\cite{GT24b}, the authors use Theorem~\ref{thm1.2} to deduce a classification of wreath products over amenable finitely presented one-ended groups up to biLipschitz equivalences~\cite[Corollary~1.15]{GT24b}. Precisely, they proved, under the above hypotheses, a strong rigidity statement: $\Z_{n}\wr H$ and $\Z_{m}\wr H$ are biLipschitz equivalent if and only if $n=m$. 

For permutational wreath products, we observe more flexibility than in the standard case. For instance:

\begin{corollaryletter}\label{cor:classificationPWPoverZ^duptoBILIP}
Let $E$ and $F$ be two non-trivial finite groups, and let $m,n$ be two integers such that $m>n\ge 2$. Then $E\wr_{\Z^{n}}\Z^{m}$ and $F\wr_{\Z^{n}}\Z^{m}$ are biLipschitz equivalent if and only if\; $|E|$ and $|F|$ are powers of a common number. 
\end{corollaryletter}

This result can be seen as an extension of~\cite[Proposition~A.2]{Cor06}, and provides examples of biLipschitz equivalent permutational wreath products over non-biLipschitz equivalent lamp groups. It is in fact a particular case of Proposition~\ref{prop:BILIPbetweenPLfromscalingQIbetweenquotients} below, proved using the same techniques as for Theorem~\ref{thm:classificationPWPuptoQI}.

Notice also that the assumption of co-amenability of the subgroup in Corollary~\ref{cor:classificationPWPoverZ^duptoBILIP} is not restrictive at all: if $M\lhd G$ is not co-amenable, one directly deduces from Theorem~\ref{thm:classificationPWPuptoQI} and Whyte's theorem (see Theorem~\ref{thm:Whytethm}) that $E\wr_{G/M}G$ and $F\wr_{H/N}H$ are biLipschitz equivalent if and only if $|E|,|F|$ have the same prime divisors and there exists a quasi-isometry of pairs $(G,M)\longrightarrow (H,N)$.

Lastly, our classification can also be used to classify wreath products whose base groups are themselves permutational wreath products. For instance, the non-amenable part of Theorem~\ref{thm:classificationPWPuptoQI} together with results from~\cite{GT24a} provide:

\begin{corollaryletter}\label{cor:classificationNAlamplightersoverPWP}
Let $n,m,p,q\ge 2$ be four integers. Let $H$ be a one-ended finitely presented group, and let $N\lhd H$ be a normal finitely generated infinite subgroup of infinite index. Suppose that $N$ is not co-amenable in $H$, and that $H$ is not coarsely separable by any collection of subspaces that uniformly quasi-isometrically embed into $N$. Then the lamplighters $\Z_{n}\wr(\Z_{p}\wr_{H/N}H)$ and $\Z_{m}\wr(\Z_{q}\wr_{H/N}H)$ are quasi-isometric if and only if $n$ and $m$ have the same prime divisors, and $p$ and $q$ have the same prime divisors.
\end{corollaryletter}

The proof of this statement is decomposed into several intermediate statements, that are themselves of independent interest for the study of rigidity and flexibility properties of such iterated wreath products.

On the amenable side, an additional scaling condition is required. We compute then the scaling groups of some permutational wreath products in Section~\ref{subsection3.6}, thus proving:

\begin{propositionletter}\label{prop:classificationAlamplightersoverPWPoverZ^d}
Let $n,m,p,q,d,k\ge 2$ be integers such that $d>k\ge 2$. Then the lamplighters $\Z_{n}\wr(\Z_{p}\wr_{\Z^{k}}\Z^{d})$ and $\Z_{m}\wr(\Z_{q}\wr_{\Z^{k}}\Z^{d})$ are quasi-isometric if and only if $n$ and $m$ are powers of a common number, and $p$ and $q$ are powers of a common number.
\end{propositionletter}

More generally, one can establish the classification up to quasi-isometry of many finitely generated groups built out of different halos structures. For instance, if $n,m,d,k\ge 2$ are four integers with $d>k$, the combination of Corollary~\ref{cor:classificationPWPoverZ^duptoBILIP} and~\cite[Corollary~8.9]{GT24a} shows that $\shuf{\Z_{n}\wr_{\Z^{k}}\Z^{d}}$ and $\shuf{\Z_{m}\wr_{\Z^{k}}\Z^{d}}$ are quasi-isometric if and only if $n$ and $m$ are powers of a common number, while if $d=k$, $\shuf{\Z_{n}\wr\Z^{d}}$ and $\shuf{\Z_{m}\wr\Z^{d}}$ are quasi-isometric if and only if $n=m$. 

\paragraph{Outlines of the proof of Theorem~\ref{thm:classificationPWPuptoQI}.} Let us now briefly explain the main steps in order to deduce Theorem~\ref{thm:classificationPWPuptoQI}.

The first part of the proof consists in proving that any quasi-isometry between two permutational wreath products always preserves, up to a finite distance, cosets of the base groups. This statement follows from a general~\textit{embedding theorem}:

\begin{theoremletter}\label{thm:embeddingthmPWPintro}
Let $F$ be a finite group. Let $H$ be a finitely generated group and $N\lhd H$ a normal finitely generated subgroup of infinite index. Let $A$ be a coarsely simply connected graph, and let $\rho\colon A\longrightarrow F\wr_{H/N}H$ be a coarse embedding. Assume that $\rho(A)$ is not coarsely separable by any collection of subspaces that uniformly quasi-isometrically embed into $N$.

Then $\rho(A)$ lies into the neighbourhood of an $H-$coset. Moreover, the size of this neighbourhood depends only on $A$, $F$, $H$ and the parameters of $\rho$.
\end{theoremletter}

The proof strategy consists in factorizing our coarse embedding through a group that “approximates” our permutational wreath product $F\wr_{H/N}H$ in some suitable sense. For extracting valuable information from this step, our approximation group needs to be not equal to our original wreath product $F\wr_{H/N}H$. This observation follows from the fact that $N$ is normal and has infinite index in $H$. Then, our approximation group surjects onto $F\wr_{H/N}H$ and sends specific subspaces, called~\textit{leaves}, to $H-$cosets. Thus we only have to prove that the image of our graph under the factorization embedding lies close to a leaf. This is done by noticing that the approximation group has a nice structure, since it splits as a semi-direct product of $H$ and a graph product of groups. This graph product exhibits a quasi-median geometry, that we exploit in order to prove the desired claim.

This embedding theorem is similar in spirit to the ones proved in~\cite{GT24b} and in~\cite{GT24a}. The main difference is in the construction of a geometric model of the approximation group. We must construct such a model taking into account the action of $H$ on the quotient $H/N$.  

As already mentioned, the proof of Theorem~\ref{thm:embeddingthmPWPintro} is the step requiring most of our assumptions, namely the finite presentation of the base groups, their non-coarse separation by subspaces of the considered subgroups, and the finite generation of the subgroups. 

In~\cite{GT24a, GT24b}, the next step consists in using the property of preserving cosets of the base group, referred to as~\textit{leaf-preservingness}, to deduce that all quasi-isometries between two wreath products (with additional assumptions) preserve the lamplighter structure in a strong way, namely they are all~\textit{aptolic}, in the sense:

\begin{definition}\label{def:aptolicQIbetweenPWP}
Let $A,B,C,D$ be four finitely generated groups. A quasi-isometry $q\colon A\wr B \longrightarrow C\wr D$ is of~\textit{aptolic form} if there exist two maps $\alpha\colon \bigoplus_{B}A\longrightarrow \bigoplus_{D}C$ and $\beta\colon B\rightarrow D$ such that 
\begin{equation*}
   \forall(c,p)\in A\wr B,\; q(c,p)=(\alpha(c),\beta(p)).
\end{equation*}
A quasi-isometry $q\colon A\wr B\longrightarrow C\wr D$ is \textit{aptolic} if it is of aptolic form and has a quasi-inverse of aptolic form. 
\end{definition}

In the permutational case, it turns out that such a rigidity statement does not hold, and we construct in Proposition~\ref{prop:non-aptolicQI} many examples of leaf-preserving non-aptolic quasi-isometries between permutational wreath products. Loosely speaking, these examples come from the fact that the lamplighter can control at the same time the color of all lamps in a single coset. 

This lack of aptolicity thus suggests that another strategy than the one followed in~\cite{GT24b} is needed to reach the conclusion. The key observation is that the leaf-preservingness property implies that the base groups must be quasi-isometric through a quasi-isometry quasi-preserving the cosets of the subgroups. The proof of this claim uses in a crucial way the fact that our subgroups are normal and the resulting parallelism between the cosets of the subgroups in the base groups. Once this claim is established, we can cone-off cosets of the subgroups in our permutational wreath products and obtain in this way an induced quasi-isometry between two standard wreath products. This induced quasi-isometry turns out to be, up to finite distance, aptolic. This phenomenon can be thought of as a form of~\textit{relative aptolicity}:

\begin{propositionletter}\label{prop:relativeaptolicity}
Let $n,m\ge2$. Let $G$, $H$ be finitely presented groups with finitely generated normal infinite groups $M\lhd G$, $N\lhd H$. Suppose that $G$ (resp. $H$) is not coarsely separable by any collection of subspaces that uniformly quasi-isometrically embed into $M$ (resp. $N$). Then any quasi-isometry $q\colon \Z_{n}\wr_{G/M}G\longrightarrow \Z_{m}\wr_{H/N}H$ induces a quasi-isometry 
\begin{equation*}
    q^{\text{in}}\colon \Z_{n}\wr G/M \longrightarrow \Z_{m}\wr H/N
\end{equation*}
that lies at finite distance from an aptolic quasi-isometry $\Z_{n}\wr G/M\longrightarrow \Z_{m}\wr H/N$.
\end{propositionletter}

Once we are reduced to an aptolic quasi-isometry between standard lamplighters, we can then deduce our conclusion thanks to results from~\cite{GT24b}. 

\paragraph{Structure of the chapter.} In Section~\ref{subsection3.2}, we fix our notations and we recall many concepts and definitions on coarse and quasi-median geometry that are used throughout the next parts. We also gather the results of other articles that will be useful for our applications.

Section~\ref{subsection3.3} is devoted to the proof of the embedding theorem (cf. Theorem~\ref{thm:embeddingthmPWPintro}) and relies heavily on tools from quasi-median geometry recalled in Section~\ref{subsection3.2}. 

Section~\ref{subsection3.4} is the core of the article, and we prove in Theorem~\ref{thm:rigiditypart}, roughly speaking, the first half of Theorem~\ref{thm:classificationPWPuptoQI}. In subsection~\ref{subsubsection3.4.3}, we show that the leaf-preservingness property guarantees the existence of a quasi-isometry of pairs between the base groups (Proposition~\ref{prop:leafpreservingQIareQIofpairs}). Subsection~\ref{subsubsection3.4.4} introduces cone-offs of graphs and several results regarding quasi-isometries between them. In subsection~\ref{subsubsection3.4.5} we combine all these results to prove Proposition~\ref{prop:relativeaptolicity} and Theorem~\ref{thm:rigiditypart}.

In Section~\ref{subsection3.5}, we prove the other half of Theorem~\ref{thm:classificationPWPuptoQI}. We begin by proving Proposition~\ref{prop:nonQIPWP} in subsection~\ref{subsubsection3.5.1} (cf. Corollary~\ref{cor:aptolicQIfromaQIofpairsbetweenbases}), because the non-amenable case is less rigid than the amenable one and provides substantially more quasi-isometric permutational lamplighters, and we treat the amenable case in subsection~\ref{subsubsection3.5.2}. We emphasize that these results require much less assumptions on the subgroups under consideration. We finish the proof of Theorem~\ref{thm:classificationPWPuptoQI} in subsection~\ref{subsubsection3.5.3}. 

Section~\ref{subsection3.6} is devoted to some consequences of our main theorem. We prove Corollaries~\ref{cor:classificationofPWPoverZ^duptoQI} and~\ref{cor:classificationofPWPoverdirectproducts} in subsection~\ref{subsubsection3.6.1}, Corollary~\ref{cor:classificationPWPoverZ^duptoBILIP} in subsection~\ref{subsubsection3.6.2}, and Corollary~\ref{cor:classificationNAlamplightersoverPWP} and Proposition~\ref{prop:classificationAlamplightersoverPWPoverZ^d} in subsection~\ref{subsubsection3.6.3}. 

Lastly, Section~\ref{subsection3.7} records various questions related to the results of the article. 

\section{Preliminaries on coarse separation and quasi-median geometry}\label{subsection3.2}

\subsection{Notations}\label{subsubsection3.2.1} In this text, all considered graphs are unoriented. If $\Gamma$ is such a graph, $V(\Gamma)$ refers to its vertex set, and $E(\Gamma)$ denotes its edge set.

Given a group $G$, we denote $1_{G}$ its neutral element, and if $G$ is generated by a finite set $S$, $\text{Cay}(G,S)$ refers to the Cayley graph of $G$ with respect to $S$, that is the graph with elements of $G$ as vertices and edges are pairs of the form $(g,gs)$ with $g\in G$ and $s\in S\cup S^{-1}\setminus\lbrace 1_{H}\rbrace$, while $\ell_{S}$ denotes the usual length function on $G$ associated to $S$. 
If $n\ge 1$ is an integer, $\Z_{n}$ stands for the cyclic group of order $n$. 

If $a,b,n$ are three integers, we write $a\equiv b\;[n]$ to say that $a-b$ is a multiple of $n$. 

Lastly, recall that, given $n\ge 2$ and a graph $X$, the~\textit{lamplighter graph over $X$} is the graph $\mathcal{L}_{n}(X)$
\begin{itemize}
    \item whose vertex-set is the set of pairs $(c,p)$ where $c\colon V(X)\rightarrow \Z_{n}$ is finitely supported (i.e. $c(x)=1$ for all but finitely many vertices $x\in V(X)$) and $p\in V(X)$;
    \item whose edges connect $(c_{1},p_{1})$ to $(c_{2},p_{2})$ if either $c_{1}=c_{2}$ and $p_{1}, p_{2}$ are adjacent in $X$, or if $p_{1}=p_{2}$ and $c_{1},c_{2}$ differ only on this vertex.
\end{itemize}

\begin{figure}[H]
  \centering
  \includegraphics[width=0.45\linewidth]{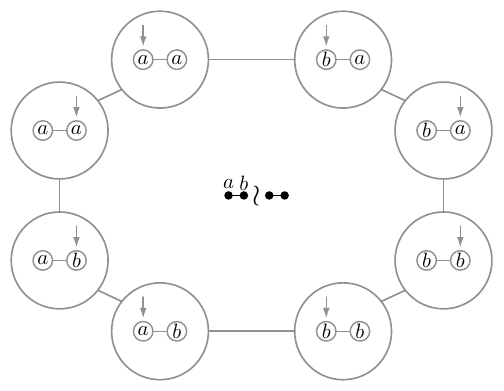}
  \caption{The lamplighter graph $\mathcal{L}_{2}(\text{edge})$.}
  \label{fig:chemin}
\end{figure}

The set of finitely supported colourings $V(X)\rightarrow \Z_{n}$ is denoted $\Z_{n}^{(V(X))}$. In the special case where $X$ is a Cayley graph of a finitely generated group $G$, we denote $\Z_{n}^{(V(X))}$ by $\Z_{n}^{(G)}$. 

\subsection{Coarse separation}\label{subsubsection3.2.2} 

\begin{definition}\label{def:coarseseparation}
Let $(X,d_{X})$ be a metric space, and let $\mathcal{Z}\subset \mathcal{P}(X)$ be a collection of subspaces of $X$. We say that $\mathcal{Z}$~\textit{coarsely separates} $X$ if there exist $k>0$ and $L\ge 0$ such that for any $D\ge 0$, there exists some $Z\in\mathcal{Z}$ such that $X\setminus Z^{+L}$ contains at least two $k-$coarsely connected components with points at distance $\ge D$ from $Z$.
\end{definition}

Naturally, the property of being coarsely separated is preserved by quasi-isometries: 

\begin{proposition}\label{prop:coarseseparationisaQIinvariant}
Let $X$ and $Y$ be two graphs. Let $\varphi\colon X\rightarrow Y$ be a quasi-isometry. If a collection of subgraphs $\mathcal{Z}\subset\mathcal{P}(X)$ coarsely separates $X$, then $\varphi(\mathcal{Z}) \defeq \lbrace \varphi(Z) : Z\in\mathcal{Z}\rbrace\subset \mathcal{P}(Y)$ coarsely separates $Y$.
\end{proposition}

\begin{proof}
Let $\overline{\varphi}$ be a quasi-inverse of $\varphi$ and let $C\ge 1$, $K\ge 0$ be such that $\varphi, \overline{\varphi}$ are $(C,K)-$quasi-isometries and such that $\overline{\varphi}\circ \varphi$, $\varphi\circ\overline{\varphi}$ are at distance at most $K$ from the identities. Since $\mathcal{Z}$ coarsely separates $X$, there is some $L\ge 0$ such that, for any $D\ge 0$, one can find a subgraph $Z\in\mathcal{Z}$ and two vertices $u,v\in V(X)$ at distance $\ge C^2(C+L+K)+C(D+3K)+L$ from $Z^{+L}$ and such that $Z^{+L}$ separates $u$ and $v$. Then one obtains that 
\begin{align*}
    d_{Y}\left(\varphi(u), \varphi(Z)^{+(C(L+C+K)+2K)}\right) &\ge d_{Y}(\varphi(u), \varphi(Z))-C(L+C+K)-2K \\
    &\ge \frac{1}{C}\cdot d_{X}(u, Z)-C(L+C+K)-3K \\
    &\ge \frac{1}{C}\cdot d_{X}(u,Z^{+L})-\frac{L}{C}-C(L+C+K)-3K \\
    &\ge D
\end{align*}
and likewise $d_{Y}\left(\varphi(v), \varphi(Z)^{+\left(C(L+C+K)+2K\right)}\right) \ge D$. Now, we claim that $\varphi(Z)^{+\left(C(L+C+K)+2K\right)}$ separates $\varphi(u)$ and $\varphi(v)$. Indeed, fix a path $\gamma\subset Y$ connecting $\varphi(u)$ to $\varphi(v)$. Then $\overline{\varphi}(\gamma)$ is a path in $X$, connecting $\overline{\varphi}(\varphi(u))$ to $\overline{\varphi}(\varphi(v))$. Fixing paths $[u,\overline{\varphi}(\varphi(u))]$, $[v,\overline{\varphi}(\varphi(v))]$ connecting $u$ to $\overline{\varphi}(\varphi(u))$ and $v$ to $\overline{\varphi}(\varphi(v))$ respectively, of length $\le K$, the concatenation 
\begin{equation*}
    \eta \defeq [u,\overline{\varphi}(\varphi(u))]\cup \overline{\varphi}(\gamma) \cup [v,\overline{\varphi}(\varphi(v))]
\end{equation*}
is a path connecting $u$ and $v$. This path must therefore intersect $Z^{+L}$, and as $d(u, Z^{+L}), d(v, Z^{+L})\ge K$, it follows that the intersection point $q\in \eta\cap Z^{+L}$ must lie on the subpath $\overline{\varphi}(\gamma)\subset \eta$ (see Figure 3.1 below). The point $q$ is at distance $\le C+K$ from a vertex $w=\overline{\varphi}(p)$ for some $p\in\gamma$. Then $w\in \eta\cap Z^{+(L+C+K)}$, hence
\begin{equation*}
    \varphi(w) \in \varphi(\overline{\varphi}(\gamma))\cap \varphi(Z^{+(L+C+K)}) \subset \gamma^{+K} \cap \varphi(Z)^{+(C(L+C+K)+K)}
\end{equation*}
using Lemma~\ref{lem:Boundingneighborhoodsingraphs}\textit{(iv)}, and $p$ being at distance $\le K$ from $\varphi(\overline{\varphi}(p))=\varphi(w)$, we conclude that
\begin{equation*}
    p\in\gamma \cap \varphi(Z)^{+\left(C(L+C+K)+2K\right)}. 
\end{equation*}
Hence $\varphi(Z)^{+\left(C(L+C+K)+2K\right)}$ separates $\varphi(u)$ and $\varphi(v)$.
\end{proof}

\begin{figure}
  \centering
  \includegraphics[width=0.8\linewidth]{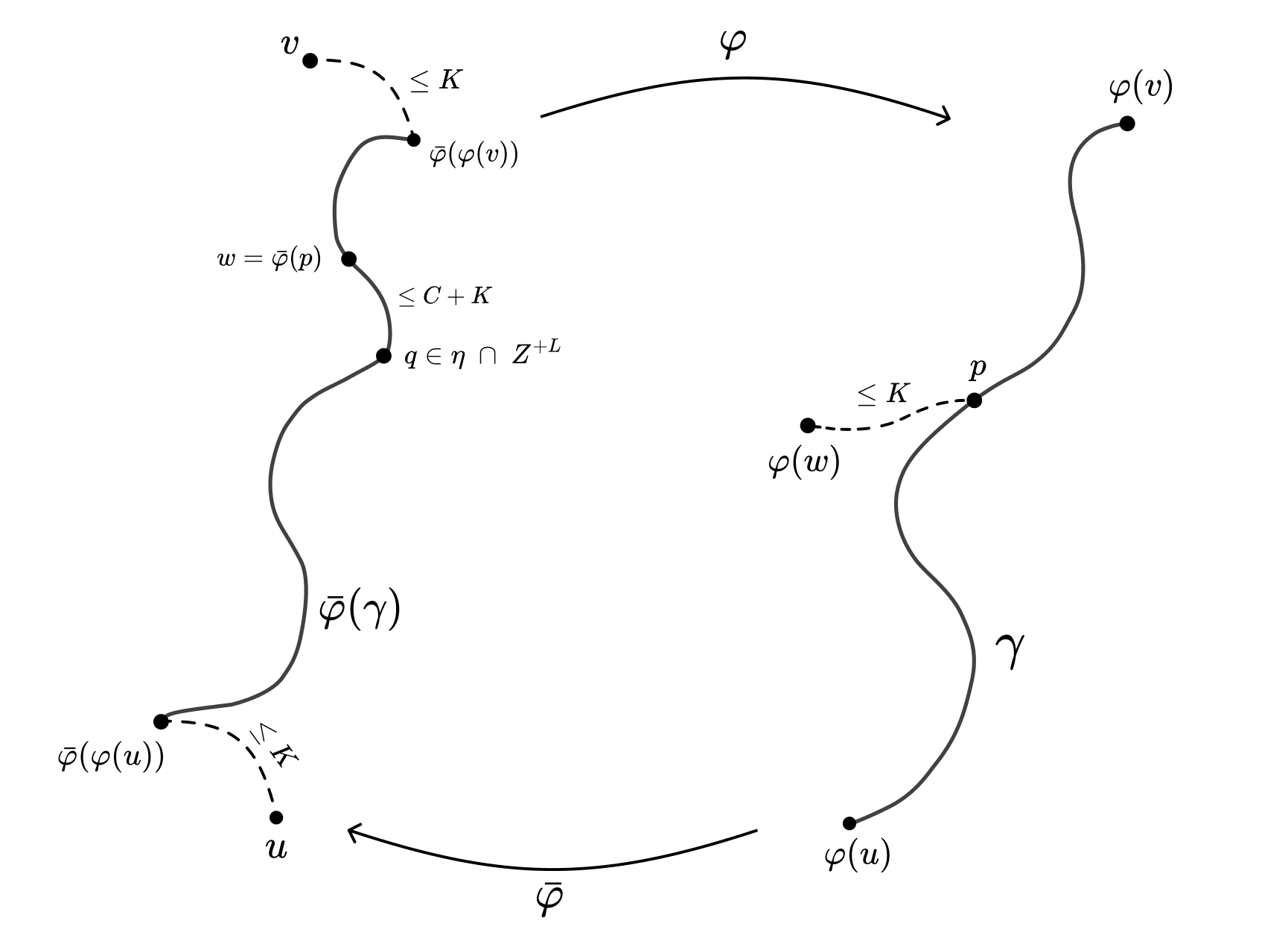}
  \caption{Configuration of the proof of Proposition~\ref{prop:coarseseparationisaQIinvariant}}
\end{figure}

In fact, the property of being coarsely separated is also preserved by coarse equivalences. See~\cite[Lemma~2.3]{BGT24} for a general statement.

We record here one of the main result from~\cite{BGT24} on the coarse separation of groups of polynomial growth, for our future applications.

\begin{theorem}[{\cite[Theorem~1.5]{BGT24}}]\label{thm:noncoarseseparationinpolynomialgrowthgroups}
Let $G$ be a finitely generated group with polynomial growth of degree $d\ge 2$. Then $G$ cannot be coarsely separated by a collection of subspaces having growth degree $\le d-2$.
\end{theorem}

We refer to~\cite[Theorem~1.5]{BGT24} for the proof of this statement.

\subsection{Quasi-median geometry}\label{subsubsection3.2.3} We conclude this section by gathering useful tools for the proof of the embedding theorem. For a more in-depth introduction on quasi-median geometry, we refer to~\cite{Gen17a, GM19}. 

\begin{definition}
A connected graph $X$ is~\textit{quasi-median} if it does not contain $K_{3,2}$ and $K_{4}^{-}$ (the graph obtained from the complete graph $K_{4}$ by removing an edge) as induced subgraphs, and if it satisfies the following two conditions:
\begin{enumerate}[label=(\roman*)]
    \item For all vertices $u,v,w\in X$ such that $d(v,w)=1$, $d(u,v)=d(u,w)=k$, there exists a vertex $x\in X$ which is a common neighbour of $v$ and $w$ and which satisfies $d(u,x)=k-1$. 
    \item For all vertices $u,v,w,z\in X$ such that $d(v,z)=d(w,z)=1$, $d(u,v)=d(u,w)=k$ and $d(u,z)=k+1$, then there exists $x\in X$ a common neighbour of $v$ and $w$ so that $d(u,x)=k-1$.
\end{enumerate}
\end{definition}

These two conditions are usually referred to as the~\textit{triangle condition} and the~\textit{quadrangle condition}. 

\begin{definition}
Let $X$ be a quasi-median graph. A~\textit{hyperplane} $J$ is an equivalence class of edges with respect to the transitive closure of the relation identifying two edges in a common $3-$cycle or two opposite edges in a common $4-$cycle. The~\textit{neighbourhood} of $J$, denoted $N(J)$, is the subgraph generated by the edges of $J$. The connected components of the graph $X\rem J$ obtained from $X$ by removing the interiors of the edges of $J$ are called the~\textit{sectors delimited by $J$}, and the connected components of $N(J)\rem J$ are the~\textit{fibers} of $J$. Two distinct hyperplanes $J_{1}, J_{2}$ are~\textit{transverse} if $J_{2}$ contains an edge in $N(J_{1})\rem J_{1}$. 
\end{definition}

Recall that in a quasi-median graph $X$, a subgraph $Y\subset X$ is~\textit{gated} if for any vertex $x\in X$, there exists a vertex $y\in Y$ so that for any $z\in Y$, there exists a geodesic from $x$ to $z$ passing through $y$. Such a vertex is referred to as the~\textit{projection of $x$ onto $Y$}. Note that a gated subgraph is convex. 

The next statement shows how combinatorics of hyperplanes in a quasi-median graph encode information about the geometry of the graph. See~\cite[Theorem~2.15]{Gen17a} for a self-contained proof. 

\begin{theorem}\label{thm:hyperplanesinQMgraphs}
Let $X$ be a quasi-median graph. The following claims hold.
\begin{enumerate}[label=(\roman*)]
    \item For any hyperplane $J$ of $X$, $X\rem J$ has at least two connected components.
    \item For any hyperplane of $X$, its sectors, neighbourhood and fibers are gated subgraphs of $X$.
    \item A path in $X$ is a geodesic if and only if it crosses each hyperplane at most once.
    \item The distance between two vertices of $X$ coincides with the number of hyperplanes separating them. 
\end{enumerate}
\end{theorem}

The next statement ensures that the Helly property provided above for finite families of gated subgraphs actually extends to infinite families of sectors of hyperplanes, under an additional assumption of finite cubical dimension. Recall that the~\textit{cubical dimension} of a quasi-median graph $X$, denoted $\text{dim}(X)$, is the maximal size of a collection of pairwise transverse hyperplanes.

\begin{lemma}\label{lem:nontrivialintersectioninQMgraphs}
Let $X$ be a quasi-median graph with $\text{dim}(X)<\infty$. For every hyperplane $J$ of $X$, let $J^{+}$ be a sector. Assume that
\begin{enumerate}[label=(\roman*)]
    \item For any hyperplanes $J_{1}, J_{2}$ of $X$, $J_{1}^{+}\cap J_{2}^{+}\neq\emptyset$.
    \item Every non-increasing sequence $J_{1}^{+}\supset J_{2}^{+}\supset \dots$ eventually stabilizes.
\end{enumerate}
Then $\dis\bigcap_{J\;\text{hyperplane}}J^{+}$ is non-empty and reduced to a single vertex. 
\end{lemma}

Another useful result using the cubical dimension is the following lemma~\cite[Lemma~5.8]{GT24b}.

\begin{lemma}\label{lem:boundingdistancesinQMgraphs}
Let $X$ be a quasi-median graph and let $x,y\in X$ be two vertices. If $k$ denotes the maximal number of pairwise non-transverse hyperplanes separating $x$ and $y$, then one has
\begin{equation*}
    d(x,y)\le k\cdot\text{dim}(X).
\end{equation*}
\end{lemma}

\paragraph{Graph products of groups.} Let $\Gamma$ be a simplicial graph and let $\mathcal{C}=\lbrace G_{u} : u\in V(\Gamma)\rbrace$ be a collection of groups indexed by the vertices of $\Gamma$. The~\textit{graph product} is the group denoted $\Gamma\mathcal{C}$ and defined by 
\begin{equation*}
    \Gamma\mathcal{C}\defeq \left. \left(\mystar_{u\in V(\Gamma)}G_{u}\right)\right/\langle\langle [g,h]: g\in G_{u}, h\in G_{v}, \lbrace u,v\rbrace\in E(\Gamma)\rangle\rangle.
\end{equation*}
The groups of the collection $\mathcal{C}$ are called the~\textit{vertex-groups}, and embed naturally into the graph product. If $G_{u}=G$ for all $u\in V(\Gamma)$, we denote the graph product $\Gamma G$ rather than $\Gamma\mathcal{C}$. 

Graph products of groups have recently attracted much attention since they provide a unified way to study several classes of groups with a geometric flavour, such as Coxeter or right-angled Artin groups (see~\cite{Gen17a, GM19, GV20, Gen24}). They also provide good examples of quasi-median graphs, as shown by the following theorem, proved in~\cite[Proposition~8.2]{Gen17a}.

\begin{theorem}\label{thm:CayleygraphsofgraphproductsareQM}
Let $\Gamma$ be a simplicial graph and $\mathcal{C}=\lbrace G_{u} : u\in V(\Gamma)\rbrace$ a collection of groups indexed by the vertices of \;$\Gamma$. Then the graph
\begin{equation*}
    \text{QM}(\Gamma, \mathcal{C}) \defeq \text{Cay}\left(\Gamma\mathcal{C},\bigcup_{u\in V(\Gamma)}G_{u}\setminus\lbrace 1\rbrace\right)
\end{equation*}
is a quasi-median graph of cubical dimension $\text{clique}(\Gamma)$.
\end{theorem}

Recall here that $\text{clique}(\Gamma)$ is the cardinality of a maximal clique of $\Gamma$.

Here also, when the collection $\mathcal{C}$ consists of a single group $G$, we will write $\text{QM}(\Gamma, G)$ instead of $\text{QM}(\Gamma,\mathcal{C})$. 

Additionally, we can relate cliques, prisms (i.e. maximal complete subgraphs) and neighbourhoods of hyperplanes of this Cayley graph to the structure of the graph $\Gamma$. The next two lemmas are the content of~\cite[Lemma~8.6]{Gen17a} and~\cite[Corollary~8.10]{Gen17a}.

\begin{lemma}\label{lem:cliquesinQMgraphs}
Let $\Gamma$ be a simplicial graph and $\mathcal{C}$ be a collection of groups indexed by $V(\Gamma)$. The cliques of $\text{QM}(\Gamma,\mathcal{C})$ coincide with the cosets of the vertex-groups.
\end{lemma}

In particular, for any $u\in V(\Gamma)$, the subgroup $G_{u}$ is a clique in $\text{QM}(\Gamma,\mathcal{C})$, and we denote $J_{u}$ the hyperplane~\textit{dual} to the clique $G_{u}$ (i.e. the edges of the clique $G_{u}$ belong to the class of edges formed by $J_{u}$; cf.~\cite[Definition~2.14]{Gen17a}). As a consequence of~\cite[Lemma~8.5 and Lemma~8.8]{Gen17a}, any hyperplane of $\text{QM}(\Gamma,\mathcal{C})$ is a translate of some $J_{u}$:

\begin{lemma}\label{lem:hyperplanesinQMraphs}
Let $\Gamma$ be a simplicial graph and $\mathcal{C}=\lbrace G_{u} : u\in V(\Gamma)\rbrace$ a collection of groups indexed by the vertices of\; $\Gamma$. Let $J$ be a hyperplane of $\text{QM}(\Gamma,\mathcal{C})$. Then there exist $\alpha\in \Gamma\mathcal{C}$ and $u\in V(\Gamma)$ such that $J=\alpha J_{u}$. 
\end{lemma}

For the next statement, recall that the~\textit{star} of a vertex $u$ of $\Gamma$, denoted $\text{star}(u)$, is the subgraph induced by $u$ and its neighbours, while $\langle \text{star}(u)\rangle$ stands for the subgroup generated by the vertices of $\text{star}(u)$, identified with a subset of $\Gamma$.

\begin{lemma}\label{lem:cliquesinhyperplanesinQMgraphs}
Let $\Gamma$ be a simplicial graph and $\mathcal{C}$ be a collection of groups indexed by $V(\Gamma)$. Let $u\in V(\Gamma)$. The cliques in $J_{u}$ are the cosets $gG_{u}$, where $g\in\langle\text{star}(u)\rangle$. Consequently, $N(J_{u})=\langle\text{star}(u)\rangle$.
\end{lemma}

An immediate consequence is the following useful fact.

\begin{lemma}\label{lem:spanningprismsinQMgraphs}
Let $\Gamma$ be a simplicial graph and $\mathcal{C}$ be a collection of groups indexed by $V(\Gamma)$. Let $J\subset \text{QM}(\Gamma,\mathcal{C})$ be a hyperplane, $x\in N(J)$, and $C_{1},C_{2}\subset N(J)$ two distinct cliques containing $x$. If $C_{1}\subset J$, then $C_{1}$ and $C_{2}$ span a prism.
\end{lemma}

\begin{proof}
Up to right-multiplying $x$ by its inverse, we may assume that $x=1$. Hence there exist $u,v\in V(\Gamma)$ so that $J=J_{u}$, $C_{1}=G_{u}$, $C_{2}=G_{v}$. By Lemma~\ref{lem:cliquesinhyperplanesinQMgraphs}, $v\in\text{star}(u)$. Since $C_{1}\neq C_{2}$, $v$ must be adjacent to $u$, so that $C_{1}$ and $C_{2}$ span the prism $G_{u}\oplus G_{v}$. 
\end{proof}

\begin{definition}\label{def:projectionfromgraphproducts}
Let $\Gamma$ be a simplicial graph and $\mathcal{C}$ be a collection of groups indexed by $V(\Gamma)$. Denote by 
\begin{equation*}
    \xi\colon \Gamma\mathcal{C}\twoheadrightarrow \bigoplus_{u\in V(\Gamma)}G_{u}
\end{equation*}
the canonical projection. For $x\in\Gamma\mathcal{C}$, let $x^{\xi}$ denote its image under $\xi$. 
\end{definition}

\section{The embedding theorem for permutational wreath products}\label{subsection3.3}

This section is the central part of the article, dedicated to the proof of the embedding theorem. We first recall the statement from the introduction, before outlining the strategy of the proof.

\begin{theorem}\label{thm:embeddingthmPWP}
Let $F$ be a finite group, $H$ be a finitely generated group and $N\lhd H$ a normal finitely generated subgroup of infinite index. Let $A$ be a coarsely simply connected graph, and let $\rho\colon A\longrightarrow F\wr_{H/N}H$ be a coarse embedding. Assume that $\rho(A)$ is not coarsely separable by any collection of subspaces that uniformly quasi-isometrically embed into $N$.

Then $\rho(A)$ lies in a neighbourhood of an $H-$coset. Moreover, the size of this neighbourhood depends only on $A$, $F$, $H$ and the parameters of $\rho$.
\end{theorem} 

The strategy towards the proof of Theorem~\ref{thm:embeddingthmPWP} is as follows. First of all, we make use of a well-known approximation theorem to factor our coarse embedding $\rho$ into a group obtained as a truncation of our initial permutational lamplighter $F\wr_{H/N}H$. In order to get a relevant factorization, the approximation needs to be non-trivial, i.e. not equal to the whole $F\wr_{H/N}H$, so we start to check that such permutational wreath products are infinitely presented. 

In a second step, we observe that these truncations are in fact semi-direct products of $H$ with graph products of copies of $F$, that are supported on Cayley graphs of $H/N$. Next, we construct a geometric model of the Cayley graph of the approximation that, moreover, captures the fact that $H\curvearrowright H/N$ has a non-trivial stabilizer. Then, this approximation, that we denote $F\square_{\Gamma}H$ naturally surjects onto $F\wr_{H/N}H$ and this projection sends~\textit{leaves}, i.e. natural copies of $H$ inside $F\square_{\Gamma}H$, to $H-$cosets. Therefore, to conclude the proof of the theorem, it is enough to prove that the image of the factorization of $\rho$ lies into a neighbourhood of a leaf of $F\square_{\Gamma}H$. We show this claim exploiting the median structure of the Cayley graph of this graph product (cf. Theorem~\ref{thm:CayleygraphsofgraphproductsareQM}) and the non-coarse separation assumption on the image $\rho(A)$.

\subsection{The approximation group}\label{subsubsection3.3.1} Here is the first ingredient towards the proof of Theorem~\ref{thm:embeddingthmPWP}. The idea is that a group $G$ having a presentation of the form $\langle S \; | \; r_{1},r_{2},\dots\rangle$ can be coarsely approximated by its finitely presented truncations $G_{p}=\langle S \; | \; r_{1},r_{2},\dots,r_{p}\rangle$, $p\ge 1$ (see~\cite[Fact~A.2]{GT25}, as well as~\cite[Lemma~6.21]{BGT24} and the references therein for more details on this strategy). Denote $\pi_{p}\colon G_{p}\twoheadrightarrow G$ the canonical morphism.

\begin{theorem}\label{thm:factorizationscoarseembeddings}
Let $A$ be a coarsely simply connected graph. For any coarse embedding $\rho\colon A\rightarrow G$, there exists $p\ge 1$ and a coarse embedding $\eta\colon A\rightarrow G_{p}$ such that $\pi_{p}\circ \eta=\rho$.
\end{theorem} 

In our setting, to ensure that this theorem provides a non-trivial factorization of our coarse embedding, we must check that our permutational wreath product $F\wr_{H/N}H$ is infinitely presented. This is the content of the next lemma:
\begin{lemma}\label{lem:permutationallamplightersareinfinitelypresented}
Let $F,H$ be non-trivial finitely presented groups, and let $N\lhd H$ be a normal finitely generated subgroup of $H$ of infinite index. Then $F\wr_{H/N}H$ is not finitely presented.
\end{lemma}

\begin{proof}
To get our conclusion, it is enough to show that at least one of the three conditions of~\cite[Theorem~1.1]{Cor06} does not hold. The first and second conditions are contained in our assumptions, thus we show that the third condition does not hold, i.e. we must check that the product action of $H$ on $H/N\times H/N$ has infinitely many orbits. This follows from the next claim: 

\begin{claim}\label{claim:Productinfinitelymanyorbits}
Let $(X,d_{X})$ be a metric space of infinite diameter, and let $G$ be a group acting isometrically on $X$. Then the product action $G\curvearrowright X^2$ has infinitely many orbits.
\end{claim}

\begingroup
\renewcommand{\qedsymbol}{$\blacksquare$}
\begin{proof}[Proof of Claim~\ref{claim:Productinfinitelymanyorbits}]
As $X$ has infinite diameter, fix two sequences $(a_{k})_{k\in\mathbb{N}}, (b_{k})_{k\in\mathbb{N}}\subset X$ such that $d_{X}(a_{k},b_{k})\ge k$ for any $k\in\mathbb{N}$. Towards a contradiction, suppose that $G\curvearrowright X^2$ has finitely many orbits $\mathcal{O}_{1},\dots,\mathcal{O}_{r}$. Then there exists $i\in\lbrace 1,\dots, r\rbrace$ such that $(a_{k},b_{k})\in \mathcal{O}_{i}$ for infinitely many $k\in\mathbb{N}$ and thus we pick $k,k' \in \mathbb{N}$ such that $d_{X}(a_{k},b_{k})<k'$ and such that $(a_{k},b_{k}), (a_{k'},b_{k'})\in\mathcal{O}_{i}$. Hence, there is $g\in G$ such that $(a_{k'},b_{k'})=g\cdot (a_{k},b_{k})=(g\cdot a_{k}, g\cdot b_{k})$, and it follows that 
\begin{equation*}
    d_{X}(a_{k},b_{k})<k'\le d_{X}(a_{k'},b_{k'})=d_{X}(g\cdot a_{k},g\cdot b_{k})=d_{X}(a_{k},b_{k})
\end{equation*}
using in the last step that $G\curvearrowright X$ is isometric. This contradiction proves that $G\curvearrowright X^2$ has infinitely many orbits.
\end{proof}
\endgroup 

\begingroup
\renewcommand{\qedsymbol}{\openbox}
\noindent In our case, assumptions of Claim~\ref{claim:Productinfinitelymanyorbits} are satisfied, because the action of $H$ on $H/N$ is isometric (when $H/N$ is equipped with a word metric induced by a finite generating set) and $H/N$ has infinite diameter since $N$ has infinite index. Thus $F\wr_{H/N}H$ is not finitely presented.
\end{proof}
\endgroup

We now give an explicit description of the truncations of our wreath product $F\wr_{H/N}H$. The latter admits
\begin{equation*}
    \langle H, F_{g} \;(g\in H/N) \; | \; [F_{1_{H}N},F_{g}] \;(g\in H/N), \;hF_{g}h^{-1}=F_{h\cdot g} \;(h\in H, g\in H/N)\rangle 
\end{equation*}
as a presentation, where $F_{g}$ is a copy of $F$, for any $g\in H/N$. Fixing $S$ a finite subset of $H/N$, we see that the truncation
\begin{equation*}
    \langle H, F_{g} \;(g\in H/N) \; | \; [F_{1_{H}N},F_{g}] \;(g\in S), \;hF_{g}h^{-1}=F_{h\cdot g} \;(h\in H, g\in H/N)\rangle 
\end{equation*}
of $F\wr_{H/N}H$ can be rewritten as 
\begin{equation*}
    \langle H, F_{g}\; (g\in H/N)\;|\; [F_{g},F_{h}] \;(g^{-1}h\in S), \;hF_{g}h^{-1}=F_{h\cdot g} \;(h\in H, g\in H/N)\rangle.
\end{equation*}
This is a presentation of a semi-direct product between $H$ and a graph product of infinitely many copies of $F$, the latter being supported on the graph $\Gamma=\text{Cay}(H/N,S)$. Denoting
\begin{equation*}
    \pi_{S}\colon \Gamma F\rtimes H\longrightarrow F\wr_{H/N}H, \; (x,h)\longmapsto (x^{\xi},h)
\end{equation*}
the natural projection, we deduce from Theorem~\ref{thm:factorizationscoarseembeddings} the following statement:
\begin{corollary}\label{cor:factorizationsofcoarsembeddingsPWP}
Let $A$ be a coarsely simply connected graph. For any coarse embedding $\rho\colon A\longrightarrow F\wr_{H/N}H$, there exists a locally finite Cayley graph $\Gamma$ of $H/N$ and a coarse embedding $\eta\colon A\longrightarrow \Gamma F\rtimes H$ such that $\pi\circ \eta=\rho$.
\end{corollary}

\subsection{A geometric model}\label{sec:geometricmodel} We present now a geometric description of the group $\Gamma F\rtimes H$, where $\Gamma$ is the graph given by Corollary~\ref{cor:factorizationsofcoarsembeddingsPWP}. Fix $T$ a finite generating set of $N$, $S$ a finite generating set of $H/N$. Then $U\defeq \iota(T)\cup s(S)$ is a finite generating set of $H$, where $s\colon H/N\rightarrow H$ is a section of the natural surjection $\pi\colon H\rightarrow H/N$, and $\iota\colon N\hookrightarrow H$ is the natural inclusion.

Let $F\square_{\Gamma}H$ be the graph of~\textit{pointed-marked cliques} of $\text{QM}(\Gamma,F)$, whose vertices are triples of the form $(C,x,u)$ where $C\subset \text{QM}(\Gamma,F)$ is a clique, $x\in C$ is the~\textit{point}, $u\in N$ is the~\textit{mark}, and whose edges have one of the following forms:
\begin{itemize}
    \item (\textit{Slide}) $(C,x,u)-(C,xs,u)$, where $C\subset \text{QM}(\Gamma,F)$ is a clique labelled by $g\in V(\Gamma)=H/N$ (i.e. $C=xF_{g}$, by Lemma~\ref{lem:cliquesinQMgraphs}), $\alpha \in F_{g}$, and $u\in N$;
    \item (\textit{Jump}) $(xF_{g}, x, u)-(xF_{g}, x, us(g)ts(g)^{-1})$, where $x\in \text{QM}(\Gamma,F)$, $g\in H/N$, $t\in T$, $u\in N$;
    \item (\textit{Rotation}) $(xF_{g},x,u)-(xF_{gq}, x, u\delta_{g, q})$ where $x\in \text{QM}(\Gamma,F)$, $g\in H/N$, $q\in S$, $u\in N$
\end{itemize}
and where, for any $a,b\in H/N$, $\delta_{a,b}\defeq s(a)s(b)s(ab)^{-1}$ is the~\textit{defect} of $s$ to be a morphism. 

Note that $\delta_{a,b}\in N$ for any $a,b\in H/N$: indeed one has 
\begin{equation*}
\pi(\delta_{a,b})=\pi(s(a)s(b)s(ab)^{-1})=\pi(s(a))\pi(s(b))\pi(s(ab)^{-1})=ab(ab)^{-1}=1_{H/N}
\end{equation*}
since $\pi$ is a morphism and $\pi\circ s=\text{Id}_{H/N}$. In particular, the third type of edges is well-defined, as well as the second since $N$ is normal in $H$.

As we will see below, these three types of edges are in one-to-one correspondence with the elementary moves in the Cayley graph of $\Gamma F\rtimes H$ endowed with its natural generating set.
The geometric picture to keep in mind is the following. Let $(C,x,u)\subset \text{QM}(\Gamma, F)$ be a pointed-marked clique. By Lemma~\ref{lem:cliquesinQMgraphs}, $C=xF_{g}$ for some $g\in V(\Gamma)=H/N$ and $x\in  \text{QM}(\Gamma, F)$. Given a generator $s\in F_{g}$, we go from $(C,x,u)$ to one of its neighbours in $F\square_{\Gamma}H$ by “sliding” the vertex $x$ through the edge between $x$ and $xs$ in $\text{QM}(\Gamma, F)$. To understand better the third type of edges, observe that cliques of $\text{QM}(\Gamma, F)$ are labelled by $V(\Gamma)=H/N$, and that two such cliques span a prism if and only if their labels in $\Gamma$ are adjacent (see Lemma~\ref{lem:spanningprismsinQMgraphs}). Given a neighbour $g'=gq$ of $g$ in $\Gamma$, we go from $(xF_{g},x,u)$ to one of its neighbours by “rotating” the clique $C=xF_{g}$ around $x$ from $C=xF_{g}$ to $xF_{gq}$. Such a rotation is “elementary”, as it corresponds to an elementary move in $\Gamma$, or equivalently it corresponds to replacing $C$ by a clique $C'$ such that $C\cup C'$ span a prism. Lastly, flexibility of the marks has to be taken into account, because our subgroup $N$ is finitely generated and contributes to the generating set of $\Gamma F\rtimes H$. As we shall see, this contribution results in a conjugation of the marks, that we record as “jumps”, leaving the clique and the point unchanged. 

For $x\in\text{QM}(\Gamma,F)$ a fixed vertex, we denote $\mathcal{L}_{x}$ the subgraph given by 
\begin{equation*}
    \mathcal{L}_{x}\defeq \lbrace (xF_{g},x,u) : g\in H/N, \;u\in N\rbrace.
\end{equation*}
We refer to such a subgraph as a~\textit{leaf} in $F\square_{\Gamma}H$.

Then we define the map 
\begin{align*}
    \psi \colon \text{Cay}(\Gamma F\rtimes H, F_{1_{H}N}\cup U) &\longrightarrow F\square_{\Gamma}H \\
    (x,h)&\longmapsto \big(xF_{\pi(h)}, x, hs(\pi(h))^{-1}\big) 
\end{align*}
as well as 
\begin{align*}
    \varphi\colon F\square_{\Gamma}H &\longrightarrow \text{Cay}(\Gamma F\rtimes H, F_{1_{H}N}\cup U) \\
    (xF_{g}, x, u) &\longmapsto (x,\iota(u)s(g)).
\end{align*}
The fact that $F\square_{\Gamma}H$ is a suitable geometric model for the semi-direct product $\Gamma F\rtimes H$ is expressed as follows:

\begin{lemma}\label{lem:graphisomorphisms}
The maps $\varphi$, $\psi$ are graph isomorphisms, that are inverses of each other. Moreover, $\varphi$ sends leaves to $H-$cosets. 
\end{lemma}

\begin{proof}
First of all, we fix a vertex $(xF_{g},x,u)\in F\square_{\Gamma}H$, and we compute that
\begin{align*}
\psi\circ\varphi(xF_{g},x,u)&=\psi(x,\iota(u)s(g))=\big(xF_{\pi(\iota(u)s(g))}, x,\iota(u)s(g)s(\pi(\iota(u)s(g)))^{-1}\big)\\
&=(xF_{g},x,u)
\end{align*}
since $\pi(\iota(u)s(g))=\pi(s(g))=g$, as $\pi\circ s=\text{Id}_{H/N}$. The other way around, if $(x,h)$ is a vertex of the Cayley graph of $\Gamma F\rtimes H$, then
\begin{equation*}
\varphi\circ\psi(x,h)=\varphi\big(xF_{\pi(h)}, x, hs(\pi(h))^{-1}\big)=\big(x, \iota(hs(\pi(h))^{-1})s(\pi(h))\big)=(x,h)
\end{equation*}
so that $\varphi, \psi$ are bijections on the vertices, inverses of each other. It remains to prove they both preserve adjacency. 

\noindent Let then $a$ and $b$ be adjacent vertices in $\text{Cay}(\Gamma F\rtimes H, F_{1_{H}N}\cup U)$. Assume first that $a=(x,h)$ and that there is $r\in F_{1_{H}N}$ so that $b=(x,h)(r,1_{H})=(xhrh^{-1}, h)$. Then 
\begin{equation*}
    \psi(a)=\big(xF_{\pi(h)}, x, hs(\pi(h))^{-1}\big)
\end{equation*}
and 
\begin{equation*}
    \psi(b)=\big(xhrh^{-1}F_{\pi(h)}, xhrh^{-1}, hs(\pi(h))^{-1}\big).
\end{equation*}
As $hrh^{-1}\in hF_{1_{H}N}h^{-1}=F_{\pi(h)}$, the latter reduces to $\psi(b)=\big(xF_{\pi(h)},xhrh^{-1}, hs(\pi(h))^{-1}\big)$, so that $\psi(a)$ and $\psi(b)$ are adjacent in $F\square_{\Gamma}H$. 

\noindent Next, assume that $a=(x,h)$ and $b=(x,ht)$ for some $t\in T$. Then one has 
\begin{equation*}
    \psi(b)=\big(xF_{\pi(ht)}, x, hts(\pi(ht))^{-1}\big)=\big(xF_{\pi(h)}, x, hts(\pi(h))^{-1}\big)
\end{equation*}
and since 
\begin{align*}
hts(\pi(h))^{-1}=hs(\pi(h))^{-1}s(\pi(h))ts(\pi(h))^{-1}
\end{align*}
we conclude that $\psi(a)$ and $\psi(b)$ are adjacent in $F\square_{\Gamma}H$. 

\noindent Lastly, if $a=(x,h)$ and $b=(x,hs(q))$ for some $q\in S$, then $\psi(a)=(xF_{\pi(h)}, x, hs(\pi(h))^{-1})$ and 
\begin{align*}
    \psi(b)&=\big(xF_{\pi(hs(q))}, x, hs(q)s(\pi(hs(q)))^{-1}\big) \\
    &=\big(xF_{\pi(h)\pi(s(q))}, x, hs(q)s(\pi(h)\pi(s(q)))^{-1}\big) \\
    &=\big(xF_{\pi(h)q}, x, hs(q)s(\pi(h)q)^{-1}\big) 
\end{align*}
and the mark of $\psi(b)$ can in fact be written as 
\begin{equation*}
hs(q)s(\pi(h)q)^{-1}=hs(\pi(h))^{-1}s(\pi(h))s(q)s(\pi(h)q)^{-1}=hs(\pi(h))^{-1}\delta_{\pi(h),q}.
\end{equation*}
It follows that $\psi(a)$ and $\psi(b)$ are adjacent in this case as well, so $\psi$ preserves adjacency. 

\noindent Conversely, let $a$ and $b$ be adjacent vertices in $F\square_{\Gamma}H$. Again, we distinguish three cases. Assume first that $a=(xF_{g},x,u)$ and $b=(xF_{g},x\alpha,u)$. Then 
\begin{equation*}
    \varphi(a)=\varphi(xF_{g},x,u)=(x,\iota(u)s(g)), \; \varphi(b)=\varphi(x\alpha F_{g}, x\alpha,u)=(x\alpha, \iota(u)s(g))
\end{equation*}
are adjacent in $\text{Cay}(\Gamma F\rtimes H, F_{1_{H}N}\cup U)$.

\noindent Assume now that $a=(xF_{g},x,u)$ and $b=(xF_{g}, x, us(g)ts(g)^{-1})$, where $t\in T$. Then $\varphi(a)=(x, \iota(u)s(g))$ and 
\begin{equation*}
    \varphi(b)=\big(x, \iota(us(g)ts(g)^{-1})s(g)\big)=\big(x, \iota(u)s(g)t\big)
\end{equation*}
so $\varphi(a)$ and $\varphi(b)$ are adjacent in $\text{Cay}(\Gamma F\rtimes H, F_{1_{H}N}\cup U)$.

\noindent For the last case, write $a=(xF_{g},x,u)$ and $b=(xF_{gq}, x, u\delta_{g,q})$ for some $q\in S$. Then 
\begin{equation*}
    \varphi(b)=\big(x, \iota(u\delta_{g,q})s(gq)\big)=\big(x, \iota(u)s(g)s(q)\big)
\end{equation*}
with $s(q)\in s(S)$, so $\varphi(a), \varphi(b)$ are adjacent in $\text{Cay}(\Gamma F\rtimes H, F_{1_{H}N}\cup U)$.

\noindent For the last claim of the statement, let $x\in \text{QM}(\Gamma,F)$, and notice that
\begin{equation*}
    \varphi(\mathcal{L}_{x})=\lbrace (x, \iota(u)s(g)) : g\in V(\Gamma)=H/N, u\in N\rbrace = (x,1_{H})H
\end{equation*}
as desired. This concludes the proof. 
\end{proof}

As a consequence, there is a natural projection $p_{\Gamma}\colon F\square_{\Gamma}H \longrightarrow F\wr_{H/N}H$, given by $p_{\Gamma}\defeq \pi\circ \varphi$, i.e.
\begin{equation*}
    p_{\Gamma}(xF_{g},x,u) \defeq (x^{\xi}, \iota(u)s(g)), \; x\in \Gamma F,\; u\in N, \;g\in H/N.
\end{equation*}

Now, since a finitely generated group is quasi-isometric to any of its Cayley graphs, and since composition of coarse embeddings are coarse embeddings, we combine Corollary~\ref{cor:factorizationsofcoarsembeddingsPWP} and Lemma~\ref{lem:graphisomorphisms} to get the following result:
\begin{corollary}\label{cor3.7}
Let $A$ be a coarsely simply connected graph. For any coarse embedding $\rho\colon A\longrightarrow F\wr_{H/N}H$, there exists a locally finite Cayley graph $\Gamma$ of $H/N$ and a coarse embedding $\eta\colon A\longrightarrow F\square_{\Gamma}H$ such that $p_{\Gamma}\circ\eta=\rho$.
\end{corollary}

We conclude this part with an easy observation on the projection $p_{\Gamma}$.

\begin{lemma}\label{lem:projectionsendsleavestocosets}
The projection $p_{\Gamma}\colon F\square_{\Gamma}H \longrightarrow F\wr_{H/N}H$ sends leaves to $H-$cosets. 
\end{lemma}

\begin{proof}
If $x\in\text{QM}(\Gamma, F)$ and $\mathcal{L}_{x}=\lbrace (xF_{g},x,u) : g\in V(\Gamma), u\in N\rbrace$ is such a leaf, then 
\begin{align*}
    p_{\Gamma}(\mathcal{L}_{x}) &= \lbrace (x^{\xi}, \iota(u)s(g)) : g\in V(\Gamma), u\in N\rbrace \\
    & =(x^{\xi}, 1_{H})\lbrace (0,\iota(u)s(g)) : g\in H/N, u\in N\rbrace \\
    &=(x^{\xi},1_{H})H
\end{align*}
and the lemma follows. 
\end{proof}

\subsection{Proof of the embedding theorem}\label{subsubsection3.3.3}

We are now ready to prove our embedding theorem, using tools from quasi-median geometry recalled in the previous section.

\begin{proof}[Proof of Theorem~\ref{thm:embeddingthmPWP}]
Let $A$ be a coarsely simply connected graph, $F,H$ be as in the statement, and let $\rho\colon A\longrightarrow F\wr_{H/N}H$ be a coarse embedding. Without loss of generality, we may assume that $\rho$ is continuous. By Corollary~\ref{cor3.7}, there exists a locally finite Cayley graph $\Gamma$ of $H/N$ and a continuous coarse embedding $\eta\colon A \longrightarrow F\square_{\Gamma}H$ such that $p_{\Gamma}\circ \eta=\rho$:

\begin{figure}[H]
  \centering
  \includegraphics[width=0.5\linewidth]{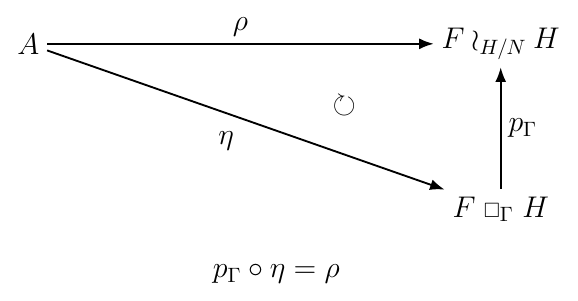}
  \caption{Configuration from the proof of Theorem~\ref{thm:embeddingthmPWP}.}
\end{figure}

\noindent By Lemma~\ref{lem:projectionsendsleavestocosets}, $p_{\Gamma}$ sends leaves to $H-$cosets. It is therefore enough to show that $\eta(A)$ lies in the neighbourhood of a leaf in $F\square_{\Gamma}H$ to conclude that $p_{\Gamma}(\eta(A))=\rho(A)$ lies in the neighbourhood of an $H-$coset in $F\wr_{H/N}H$. 

\begin{claim}\label{claim:hyperplaneseparate}
Let $J$ be a hyperplane of $\text{QM}(\Gamma,F)$, and let 
\begin{equation*}
    \mathcal{H}_{J}\defeq \lbrace (C,x,u) \in \eta(A) : C\subset J, x\in N(J)\rbrace.
\end{equation*}
There exist a sector $W$ and a constant $R\ge 0$, depending only on $A,\Gamma$ and the parameters of $\eta$ and $\rho$ such that $\eta(A)_{W}$ contains points arbitrarily far away from $\mathcal{H}_{J}$ and any vertex of $\eta(A)_{W'}$ is at distance at most $R$ from $\mathcal{H}_{J}$, for any sector $W'\neq W$ delimited by $J$. 
\end{claim}

Here, the notation $\eta(A)_{W}$ (resp. $\eta(A)_{W'}$) is used as a shortcut to denote the subset of vertices of $F\square_{\Gamma}H$ having their first coordinate being cliques lying in $W$ (resp. in $W'$).

\begingroup
\renewcommand{\qedsymbol}{$\blacksquare$}
\begin{proof}[Proof of Claim~\ref{claim:hyperplaneseparate}]
By Lemma~\ref{lem:hyperplanesinQMraphs}, we may write $J=\alpha J_{g}$ for some $g\in V(\Gamma)$ and $\alpha\in \Gamma F$. Given $(C,x,u)\in \mathcal{H}_{J}$, Lemma~\ref{lem:cliquesinhyperplanesinQMgraphs} tells us that $C=\alpha y F_{g}$ for $y\in\langle \text{star}(g)\rangle$, so $x=\alpha y f_{g}$ for some $f_{g}\in F_{g}$, whence $x^{\xi}=\alpha^{\xi}y^{\xi}f_{g}^{\xi}$. Since $c\defeq y^{\xi}f_{g}^{\xi}$ is an element of $\bigoplus_{H/N}F$ whose non-trivial coordinates can only be the ones indexed by vertices of $\text{star}(g)$, it follows that 
\begin{align*}
p_{\Gamma}(\mathcal{H}_{J})&=\lbrace (x^{\xi}, \iota(u)s(g)) \in \rho(A) : x\in N(J),\; u\in N\rbrace \\
&=\lbrace (\alpha^{\xi}c, \iota(u)s(g))\in\rho(A) : \text{supp}(c)\subset \text{star}(g),\;u\in N\rbrace \\
&=\rho(A)\cap \left(\bigsqcup_{c\in\oplus_{H/N} F,\;\text{supp}(c)\subset\text{star}(g)}\big\lbrace (\alpha^{\xi}c,\iota(u)s(g)) : u\in N\big\rbrace\right) \\
&=\rho(A)\cap \left(\bigsqcup_{c\in\oplus_{H/N} F,\;\text{supp}(c)\subset\text{star}(g)}(\alpha^{\xi}c,1_{H})\big\lbrace (\mathbf{1},\iota(u)s(g)) : u\in N\big\rbrace\right)
\end{align*}
where we denote $\text{supp}(c) \defeq \lbrace g\in H/N : c(g) \neq 1_{F}\rbrace$ the~\textit{support} of an element $c\in\bigoplus_{H/N}F$. Since $\Gamma$ is locally finite, $\text{star}(g)$ is finite, so $p_{\Gamma}(\mathcal{H}_{J})$ is a subspace of $\rho(A)$ which is a translate of a subspace which is quasi-isometric to $N$. Indeed, setting $Z_{g}\defeq \lbrace (\mathbf{1},\iota(u)s(g)) : u\in N\rbrace$, the natural map
\begin{align*}
    Z_{g} &\longrightarrow N \\
    (\mathbf{1},\iota(u)s(g))&\longmapsto u
\end{align*}
is a quasi-isometry (actually even a bijective quasi-isometry), when $Z_{g}$ is endowed with the induced metric from $F\wr_{H/N}H$ and $N$ is endowed with the word metric coming from its generating set $T$ (that we fixed once and for all at the beginning of Section~\ref{sec:geometricmodel}). By assumption, $\rho(A)$ cannot be coarsely separated by the collection of its subspaces that uniformly quasi-isometrically embed into $N$, which implies that there exists a constant $D\ge0$ so that $\rho(A)\setminus p_{\Gamma}(\mathcal{H}_{J})$ contains one coarsely connected component with points at arbitrarily large distance from $p_{\Gamma}(\mathcal{H}_{J})$, and any point in another component is at distance $<D$ from $p_{\Gamma}(\mathcal{H}_{J})$. The claim follows. 
\end{proof}
\endgroup

For any hyperplane $J$ of $\text{QM}(\Gamma,F)$, denote $J^{+}$ the sector given by Claim~\ref{claim:hyperplaneseparate}.

\begin{claim}\label{claim:intersectionreducedtoapoint}
The intersection $\dis\bigcap_{J\;\text{hyperplane}}J^{+}$ is reduced to a single vertex.
\end{claim}

\begingroup
\renewcommand{\qedsymbol}{$\blacksquare$}
\begin{proof}[Proof of Claim~\ref{claim:intersectionreducedtoapoint}]
We check the assumptions of Lemma~\ref{lem:nontrivialintersectioninQMgraphs}. First, if $J_{1}$ and $J_{2}$ are transverse, then $J_{1}^{+}$ and $J_{2}^{+}$ clearly intersect. If $J_{1}, J_{2}$ are not transverse and $J_{1}^{+}$ contains $J_{2}$, then $J_{1}^{+}$ and $J_{2}^{+}$ intersect. Lastly, if $J_{1}$ and $J_{2}$ are not transverse and $J_{1}^{+}$ does not contain $J_{2}$, then $J_{2}^{+}$ must be the sector delimited by $J_{2}$ that contains $J_{1}$ as a consequence of Claim~\ref{claim:hyperplaneseparate}. Hence $J_{1}^{+}$ and $J_{2}^{+}$ intersect.

\noindent Now, towards a contradiction, assume there exists a decreasing sequence $J_{1}^{+}\supset J_{2}^{+}\supset \dots$ that does not stabilize. Fix a vertex $(C,x,u)\in F\square_{\Gamma}H$ with $C\subset J_{1}^{+}$. Since any two vertices of $\text{QM}(\Gamma,F)$ are separated by only finitely many hyperplanes, there is $r\ge 1$ so that $C$ is disjoint from $J_{r}^{+}$. Given an $s>r$, we then know that $C$ is disjoint from $J_{s}^{+}$ (since $J_{r}^{+}\supset J_{s}^{+}$), so it follows from Claim~\ref{claim:hyperplaneseparate} that $(C,x,u)$ is at distance at most $R$ from $\mathcal{H}_{J_{s}}$. Hence 
\begin{equation*}
    R \ge d((C,x,u),\mathcal{H}_{J_{s}}) \ge d(x, J_{s}) \ge s-r
\end{equation*}
where the second inequality comes from the definition of the metric on $F\square_{\Gamma}H$, and the last inequality comes from Theorem~\ref{thm:hyperplanesinQMgraphs}\textit{(iv)}, since $J_{r},\dots, J_{s-1}$ separate $x$ from $J_{s}$. We get a contradiction if $s$ is chosen sufficiently large, namely $s>r+R$.
\end{proof}
\endgroup

\noindent Let thus $x\in\text{QM}(\Gamma,F)$ be the vertex given by Claim~\ref{claim:intersectionreducedtoapoint}, and consider the leaf 
\begin{equation*}
    \mathcal{L}\defeq \lbrace (xF_{g},x,u) : g\in V(\Gamma), u\in N\rbrace
\end{equation*}
in $F\square_{\Gamma}H$. We show that the image of $\eta$ lies in a neighbourhood of $\mathcal{L}$. 

\noindent Let $(Q,z,w)$ be a vertex in $\eta(A)$. Fix a geodesic $\gamma$ from $z$ to $x$ in $\text{QM}(\Gamma,F)$, and let $K$ be the clique containing the last edge of $\gamma$. Hence the hyperplane $J$ of $\text{QM}(\Gamma,F)$ containing $K$ separates $z$ and $x$. Since moreover $Q\not\subset J^{+}$, it follows from the definition of $J^{+}$ that $J$ separates $z$ from a pointed-marked clique corresponding to a vertex in $\eta(A)$, in the sense that the point of this vertex is separated from $x$ by $J$. Since $\eta(A)$ is connected, there must exist $(C,y,u)\in\eta(A)$ so that $C\subset J$. Denoting $z'\in C$ the unique vertex of $C$ in the same sector delimited by $J$ as $z$ (uniqueness coming from the fact that, since $J$ separates the vertices of $C$, there cannot be two vertices of $C$ lying in the same sector delimited by $J$ as $z$), and applying the triangle inequality to the triple of vertices $(Q,z,w),(C,z',u),(C,y,u)$, one has
\begin{equation}\label{eq3.1}
    d((Q,z,w),(C,y,u))\le R+2.
\end{equation}

\begin{claim}\label{claim:closetoamark}
There exists a mark $b\in N$ such that $d((C,y,u),(K,x,b)) \le 3d(y,x)+1$.
\end{claim}

\begingroup
\renewcommand{\qedsymbol}{$\blacksquare$}
\begin{proof}[Proof of Claim~\ref{claim:closetoamark}]
Let $x_{0}$ be the unique vertex of $C$ that belongs to $J^{+}$, as $x$. Fix a geodesic $x_{0},\dots,x_{k-1},x_{k}=x$ from $x_{0}$ to $x$ in $\text{QM}(\Gamma,F)$.  By Theorem~\ref{thm:hyperplanesinQMgraphs}(\textit{ii}), this geodesic lies in a fiber of $J$, so for any $i\in\lbrace 0,\dots,k\rbrace$, $x_{i}$ belongs to a clique $C_{i}\subset J$. By construction, $C_{0}=C$ and $C_{k}=K$. For all $i\in\lbrace 0,\dots,k-1\rbrace$, denote $K_{i}$ the clique containing the edge $[x_{i},x_{i+1}]$. By Lemma~\ref{lem:spanningprismsinQMgraphs}, $C_{i}$ spans a prism with both $K_{i-1}$ and $K_{i}$. We can therefore use slides and rotations in $F\square_{\Gamma}H$ to construct a path from $(C,y,u)$ to some $(K,x,b)\in\mathcal{L}$ as follows:
\begin{itemize}
    \item Start from $(C,y,u)=(C_{0},y,u)$, and slide to $(C_{0},x_{0},u)$.
    \item Use a rotation (that will modify the mark) to go from $(C_{0},x_{0},u)$ to some vertex $(K_{0},x_{0},u_{0})$. Then slide to $(K_{0},x_{1},u_{0})$.
    \item Rotate the clique $K_{0}$ to $C_{1}$ to go to some vertex $(C_{1},x_{1},u_{1})$, and then rotate again to go to $(K_{1},x_{1},\Tilde{u_{1}})$. Then slide to $(K_{1}, x_{2},\Tilde{u_{1}})$.
    \item Rotate now to go to $(C_{2},x_{2},u_{2})$. 
    \item More generally, being at $(C_{i},x_{i},u_{i})$ for some $u_{i}\in N$, make a rotation-slide-rotation to go to $(C_{i+1},x_{i+1},u_{i+1})$ for some $u_{i+1}\in N$. 
\end{itemize}
At the end, we reach the vertex $(C_{k},x_{k},u_{k})=(K,x,u_{k})$. Set $b\defeq u_{k}$. Counting each step in the above path, it follows that 
\begin{equation*}
    d((C,y,u),(K,x,b)) \le 3k+1 \le 3d(y,x)+1
\end{equation*}
as claimed.
\end{proof}
\endgroup

\begingroup
\renewcommand{\qedsymbol}{\openbox}
\noindent It thus remains to estimate $d(y,x)$. Let $J_{1},\dots,J_{\ell}$ be a maximal collection of pairwise non-transverse hyperplanes separating $z$ and $x$. Without restrictions, we may assume that $J_{i}$ separates $J_{i-1}$ from $J_{i+1}$ for any $i\in\lbrace 2,\dots,\ell-1\rbrace$, and that $J_{1}$ separates $z$ from $J_{k}$. As $C$ is disjoint from $J_{k}^{+}$ and that $\eta(A)$ is connected, we find a vertex $(M,p,e)\in\eta(A)$ so that $M\subset J_{k}$. Exactly as above, we have
\begin{equation}\label{eq3.2}
    d(z,p)\le d\big((Q,z,w),(M,p,e)\big)\le R+2
\end{equation}
and on the other hand, since $z$ and $p$ are separated by $J_{1},\dots,J_{\ell-1}$, we also have $d(z,p)\ge \ell-1$. Thus 
\begin{align*}
    d(y,x)\le \text{clique}(\Gamma)\cdot\ell \le \text{clique}(\Gamma)\cdot(d(z,p)+1) \le \text{clique}(\Gamma)\cdot(R+3)
\end{align*}
using Lemma~\ref{lem:boundingdistancesinQMgraphs} for the first inequality and (\ref{eq3.2}) for the last inequality. As also $d(y,z)\le d((C,y,u),(Q,z,w))\le R+2$ by (\ref{eq3.1}), we conclude that
\begin{equation}\label{eq3.3}
    d(y,x)\le d(y,z)+d(z,x)\le R+2+\text{clique}(\Gamma)(R+3).
\end{equation}
Combining inequality (\ref{eq3.1}) and Claim \ref{claim:closetoamark}, we deduce 
\begin{align*}
d\big((Q,z,w),(K,x,b)\big) &\le d\big((Q,z,w), (C,y,u)\big)+d\big((C,y,u),(K,x,b)\big) \\
&\le R+2+3d(y,x)+1 \\
&\le R+2+3\big(R+2+\text{clique}(\Gamma)(R+3)\big)+1 \\
&=4R+9+\text{clique}(\Gamma)(3R+9).
\end{align*}
Hence it follows that $d((Q,z,w), \mathcal{L}) \le 4R+9+\text{clique}(\Gamma)(3R+9)$, and we conclude that $\eta(A)$ lies in a neighbourhood of a leaf in $F\square_{\Gamma}H$, as was to be shown. 
\end{proof}
\endgroup 

\section{Large-scale geometry of permutational lamplighters}\label{subsection3.4}

In this section, we use our embedding theorem and several additional tools on cone-offs of graphs in order to prove:

\begin{theorem}\label{thm:rigiditypart}
Let $E$ and $F$ be non-trivial finite groups. Let $G,H$ be finitely presented groups with normal finitely generated infinite subgroups $M\lhd G$, $N\lhd H$. Assume that $M$ has infinite index in $G$ and that $G$ (resp. $H$) is not coarsely separable by any collection of subspaces that uniformly quasi-isometrically embed into $M$ (resp. $N$). If $E\wr_{G/M}G$ and $F\wr_{H/N}H$ are quasi-isometric, then 
\begin{enumerate} [label=(\roman*)]
    \item $|E|$ and $|F|$ have the same prime divisors;
    \item There exists a quasi-isometry of pairs $(G,M)\longrightarrow (H,N)$. In particular, $M$ and $N$ are quasi-isometric, and $G/M$ and $H/N$ are quasi-isometric.
\end{enumerate}
Moreover, if $M$ is co-amenable in $G$, then $N$ is co-amenable in $H$ and there exist $n,r,s\ge 1$ such that $|E|=n^{r}$, $|F|=n^{s}$, and the quasi-isometry of pairs from (ii) induces a quasi-$\frac{s}{r}$-to-one quasi-isometry $G/M\rightarrow H/N$. 
\end{theorem}

\subsection{An alternative description of permutational lamplighters}\label{subsubsection3.4.1} Let us first introduce some notations and a convenient description of permutational wreath products.

\begin{definition}
Let $n\ge 2$, and let $X$ be a connected graph together with a partition $\mathcal{C}$ of $V(X)$. The lamplighter graph over $X$ with respect to $\mathcal{C}$ is the graph $\mathcal{L}_{n}(X,\mathcal{C})$ 
\begin{itemize}
    \item whose vertices are pairs $(c,p)$, where $c\colon X\rightarrow \Z_{n}$ is a colouring of $X$ such that $c$ is constant on the pieces of $\mathcal{C}$ and all but finitely many pieces have a trivial color (the neutral element of $\Z_{n}$), and $p\in V(X)$;
    \item whose edges connect $(c_{1},p_{1}), (c_{2},p_{2})$ either if $c_{1}=c_{2}$ and $p_{1}, p_{2}$ are adjacent in $X$, or if $p_{1}=p_{2}$ and $c_{1},c_{2}$ only differ on the piece containing $p_{1}$. 
\end{itemize}
\end{definition}

The set of colourings as in the first point above is denoted $\Z_{n}^{(X,\mathcal{C})}$.

Hence we think of the vertex $(c,p)\in \mathcal{L}_{n}(X,\mathcal{C})$ as being a colouring of the pieces of $X$ (these pieces also being referred to as~\textit{zones}), with all but finitely many pieces having a trivial color, together with an arrow pointing at $p\in V(X)$. The two types of edges correspond then to:
\begin{itemize}
    \item either the arrow changes its position from $p$ to a neighbour of $p$, and the colouring stays unchanged;
    \item or the arrow stays on the vertex where it stands, but changes the color of the zone containing this vertex. 
\end{itemize}

We define the support of a colouring $c\in\Z_{n}^{(X,\mathcal{C})}$ as being the collection of pieces of $\mathcal{C}$ where $c$ is non-trivial, i.e. 
\begin{equation*}
    \text{supp}(c) \defeq \left\lbrace C\in\mathcal{C} : c(C)\neq 0\right\rbrace. 
\end{equation*}

In the graph $\mathcal{L}_{n}(X,\mathcal{C})$, the distance between $(c_{1},p_{1})$ and $(c_{2},p_{2})$ is thus given by 
\begin{equation*}
    d\big((c_{1},p_{1}), (c_{2},p_{2})\big)=\text{length}(\alpha)+\left|\text{supp}(c_{1}^{-1}c_{2})\right|
\end{equation*}
where $\alpha$ is the shortest path in $X$ starting from $p_{1}$, visiting all zones where $c_{1}$ and $c_{2}$ differ, and ending at $p_{2}$. 

Lastly, we define also~\textit{leaves} of $\mathcal{L}_{n}(X,\mathcal{C})$ as subgraphs of the form 
\begin{equation*}
    L(c) \defeq \lbrace (c,p)\in \mathcal{L}_{n}(X,\mathcal{C}) : p\in V(X)\rbrace
\end{equation*}
where $c\in \Z_{n}^{(X,\mathcal{C})}$ is a fixed colouring and where $(c, p_{1})$ and $(c,p_{2})$ are connected by an edge if $p_{1}$ and $p_{2}$ are connected by an edge in $X$.

Hence, we can think of edges contained in a single leaf as “horizontal” edges, along which the arrow moves, while “vertical” edges are those connecting different leaves in the lamplighter graph, changing only the colouring and keeping the arrow on the same vertex.

An elementary but crucial fact about leaves is that they are convex. 

\begin{lemma}\label{lem:leavesareconvex}
Let $n\ge 2$. Let $X$ be a connected graph together with a partition $\mathcal{C}$ of $V(X)$. In $\mathcal{L}_{n}(X,\mathcal{C})$, leaves are convex. In particular, they are isometrically embedded in $\mathcal{L}_{n}(X,\mathcal{C})$.
\end{lemma}

\begin{proof}
Let $P=L(c)$ be a leaf in $\mathcal{L}_{n}(X,\mathcal{C})$. Fix $a=(c,p)\in P$, $b=(c,q)\in P$ and a geodesic $\gamma$ connecting $a$ and $b$. Towards a contradiction, suppose that some subpath $\sigma\subset\gamma$ stays outside of $P$. Let $v_{0}=(c,p_{0})$ and $v_{r}=(c,p_{r})$ denote the starting point and the ending point of $\sigma$. Consider the subpaths of $\sigma$ consisting of horizontal edges, and denote $q_{1},\dots,q_{\ell}$ the second components of the vertices along those subpaths. Then the concatenation 
\begin{equation*}
    \gamma'\cup [(c,p_{0}), (c,q_{1})]\cup [(c,q_{1}),(c,q_{2})]\cup\dots\cup [(c,q_{\ell-1}), (c,q_{\ell})]\cup\gamma''
\end{equation*}
where $\gamma'$ (resp. $\gamma''$) denotes the subpath of $\gamma$ connecting $a$ (resp. $v_{r}$) and $v_{0}$ (resp. $b$), is a path connecting $a$ and $b$ which is shorter than $\gamma$. This contradicts the fact that $\gamma$ is a geodesic. Thus $P$ is convex, as claimed.  
\end{proof}

\begin{figure}[H]
  \centering
  \includegraphics[width=0.9\linewidth]{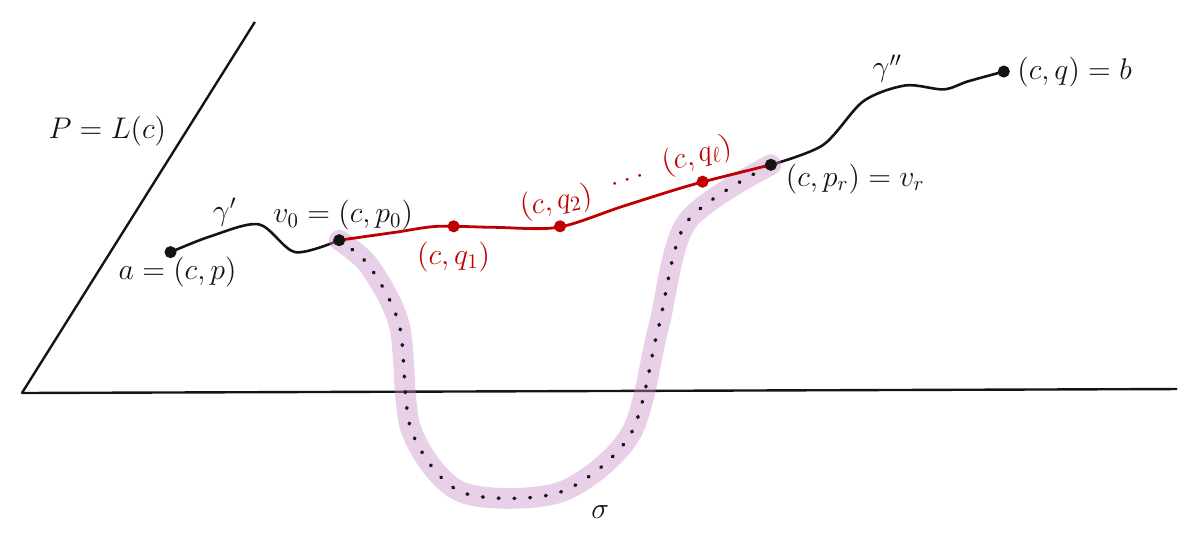}
  \caption{The situation of Lemma~\ref{lem:leavesareconvex}}
\end{figure}

Now, observe that given a non-trivial finite group $F$, a finitely generated group $H$ with a finite generating set $S$ and a subgroup $N\leq H$, one has a graph isomorphism
\begin{equation*}
    \text{Cay}\left(F\wr_{H/N}H, F\cup S\right)\cong \mathcal{L}_{|F|}\left(\text{Cay}(H,S), \mathcal{C}_{N}\right)
\end{equation*}
where $\mathcal{C}_{N}$ denotes the collection of left $N-$cosets. Note that, in $\mathcal{L}_{|F|}\left(\text{Cay}(H,S), \mathcal{C}_{N}\right)$, leaves are simply $H-$cosets, since
\begin{equation*}
    L(c)=\lbrace (c,p) : p\in H\rbrace=(c,1_{H})H.
\end{equation*}

Hence, following our previous description of lamplighters over partitioned graphs, we can see an element $(c,p)\in F\wr_{H/N}H$ as a pair made of a colouring of the left $N-$cosets, supported on finitely many cosets, together with an arrow pointing at some $p\in H$. More precisely, there is a canonical identification between $\bigoplus_{H/N} F$ and the set 
\begin{align*}
    F^{(H,\mathcal{C}_{N})} \defeq \big\lbrace &c\colon H\rightarrow F \;|\; c \;\text{constant on each $N-$coset and}\\
    &c\equiv 1_{F} \;\text{for all but finitely many $N-$cosets}\big\rbrace
\end{align*}
where $\mathcal{C}_{N}$ denotes the collection of left $N-$cosets (i.e. the pieces of the partition). We will denote $\textbf{1}$ the trivial colouring, defined by $\textbf{1}(pN)=1_{F}$ for any $p\in H$. 

The set $F^{(H,\mathcal{C}_{N})}$ has an obvious group structure (written multiplicatively), when identified with $\bigoplus_{H/N}F$. Given now any $S\subset H/N$, the subgroup $\bigoplus_{S}F$ of $\bigoplus_{H/N}F$ is identified to the subgroup $F^{(H,S)}$ of $F^{(H,\mathcal{C}_{N})}$ consisting of colourings supported only on zones from $S$. 

Therefore, the two elementary moves we can do in $\text{Cay}\left(F\wr_{H/N}H, F\cup S\right)$ to go from $(c,p)$ to a neighbour are:
\begin{itemize}[label=\textbullet]
\item either the colouring stays the same, and the arrow moves from $p$ to a neighbour in $\text{Cay}(H,S)$;
\item or the arrow stays where it stands, and changes the color of the coset $pN$.
\end{itemize}

\begin{figure}[H]
  \centering
  \includegraphics[width=0.9\linewidth, trim=0 0 0 4mm, clip]{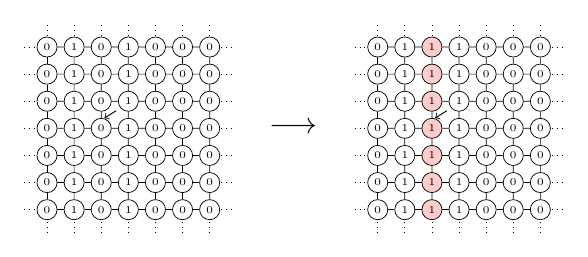}
  \caption{Two neighbours in $\text{Cay}(\Z/2\Z\wr_{\Z}\Z^2, \Z/2\Z\cup\lbrace\pm(1,0),\pm(0,1)\rbrace)$.}
  \label{fig:cayley_path_perm}
\end{figure}

Additionally, if $gN$, $hN$ are two zones in $H$, the distance between them is defined by
\begin{equation*}
    d_{H}(gN, hN) \defeq \min_{a\in gN, \;b\in hN}d_{H}(a,b)
\end{equation*}
and given $x\in H$, $R\ge0$, we let $\lbrace x\rbrace_{+R} \subset \mathcal{C}_{N}$ denote the subset of zones at distance at most $R$ from $x$, i.e.
\begin{equation*}
    \lbrace x\rbrace_{+R}\defeq \lbrace S\in\mathcal{C}_{N} : d_{H}(x,S)\le R\rbrace.
\end{equation*}

In the sequel, to shorten notations, we will often denote the leaf $L(c)$ as $cH$. Also, if $n\ge 1$ and if $H$ is a finitely generated group with a subgroup $N\le H$, we write merely $\mathcal{L}_{n}(H,\mathcal{C}_{N})$ for the graph $\mathcal{L}_{n}(\text{Cay}(H,S), \mathcal{C}_{N})$, where $S$ is a finite generating set of $H$. Note that the large-scale geometry of $\mathcal{L}_{n}(\text{Cay}(H,S), \mathcal{C}_{N})$ is in fact independent of the choice of the generating set $S$, since: 

\begin{fact}
If $X$ and $Y$ are biLipschitz equivalent graphs through $\varphi\colon X\rightarrow Y$ and if $\mathcal{C}$ is a partition of $X$, then $\mathcal{L}_{n}(X, \mathcal{C})$ and $\mathcal{L}_{n}(Y, \varphi(\mathcal{C}))$ are biLipschitz equivalent.
\end{fact}
 
\begin{proof}
Since $\varphi$ is bijective and $\mathcal{C}$ is a partition of $V(X)$, $\varphi(\mathcal{C})$ is a partition of $V(Y)$. Using that $\varphi$ and its inverse $\varphi^{-1}\colon Y\rightarrow X$ are Lipschitz, it is not hard to prove with Lemma~\ref{lem:Lipschitzbetweengraphs} that the map 
\begin{align*}
    f\colon \mathcal{L}_{n}(X, \mathcal{C}) &\longrightarrow \mathcal{L}_{n}(Y, \varphi(\mathcal{C})) \\
    (c,p)&\longmapsto (c\circ \varphi^{-1}, \varphi(p))
\end{align*}
is Lipschitz, as well as its inverse 
\begin{align*}
    f^{-1}\colon \mathcal{L}_{n}(Y, \varphi(\mathcal{C})) &\longrightarrow \mathcal{L}_{n}(X, \mathcal{C}) \\
    (c,p)&\longmapsto (c\circ \varphi, \varphi^{-1}(p)).
\end{align*}
\end{proof}

To conclude this part, we observe that some lamplighters over partitioned graphs may coincide with standard lamplighter graphs. The next criterion will be crucial in our future purposes:

\begin{proposition}\label{prop:PLovergraphswithboundedpiecesreducetoCL}
Let $n\ge 2$ and let $X$ be a graph with a partition $\mathcal{C}$. Assume that there exists $Q\ge 0$ such that $\text{diam}(C)\le Q$ for any $C\in\mathcal{C}$. Then there exists a graph $Y$, which is quasi-isometric to $X$, such that $\mathcal{L}_{n}(X,\mathcal{C})$ is quasi-isometric to $\mathcal{L}_{n}(Y)$.
\end{proposition}

\begin{proof}
For any $C\in\mathcal{C}$, choose a vertex $x_{C}\in C$. Define $Y$ as the graph 
\begin{itemize}
    \item whose vertex set is $\lbrace x_{C} : C\in\mathcal{C}\rbrace$;
    \item whose edges connect $x_{C}$ to $x_{C'}$ whenever $C\neq C'\in \mathcal{C}$ contain adjacent vertices in $X$.
\end{itemize}
Define then the map 
\begin{align*}
    \varphi\colon X &\longrightarrow Y \\
    x&\longmapsto x_{C}, \;\text{where $C\in\mathcal{C}$ contains $x$.} 
\end{align*}
Notice first that $\varphi$ is surjective, in particular coarsely surjective. Next, let $x,y\in X$ be adjacent vertices. There are two possible cases:
\begin{itemize}
    \item Either $x,y\in C$ are in a same piece of the partition, in which case $\varphi(x)=\varphi(y)$, so $d_{Y}(\varphi(x),\varphi(y))=0$;
    \item Or $x\in C$, $y\in C'$, $C\neq C'$. Then $x_{C}$, $x_{C'}$ are adjacent in $Y$, so
    \begin{equation*}
        d_{Y}(\varphi(x),\varphi(y))=d_{Y}(x_{C}, x_{C'})=1.
    \end{equation*}
\end{itemize}
In both cases, we have $d_{Y}(\varphi(x),\varphi(y))\le 1$, so we conclude from Lemma~\ref{lem:Lipschitzbetweengraphs} that $\varphi$ is $1-$Lipschitz. The other way around, define 
\begin{align*}
    \psi\colon Y&\longrightarrow X \\
    x_{C}&\longmapsto x_{C}
\end{align*}
and note that, since any vertex $x\in X$ lies at distance at most $Q$ from a vertex of $Y$, $\psi$ is coarsely surjective. Additionally, if $x_{C}, x_{C'}\in Y$ are adjacent, then $C,C'$ contain adjacent vertices in $X$, say $u\in C$ and $v\in C'$, and thus
\begin{equation*}
    d_{X}(\psi(x_{C}), \psi(x_{C'}))=d_{X}(x_{C},x_{C'})\le d_{X}(x_{C},u)+d_{X}(u,v)+d_{X}(v,x_{C'}) \le 2Q+1.
\end{equation*}
It follows from Lemma~\ref{lem:Lipschitzbetweengraphs} that $\psi$ is $(2Q+1)-$Lipschitz. Lastly, $\varphi\circ \psi=\text{Id}_{Y}$ and for any $x\in X$ one has
\begin{equation*}
    d_{X}(x,\psi(\varphi(x)))=d_{X}(x,x_{C}) \le Q
\end{equation*}
where $C\in\mathcal{C}$ is the piece containing $x$. 
Hence $\psi\circ\varphi$ is at distance $\le Q$ from $\text{Id}_{X}$, and we deduce that 
\begin{align*}
    d_{X}(x,y) &\ge d_{Y}(\varphi(x),\varphi(y)) \\
    &\ge \frac{1}{2Q+1}d_{X}\big(\psi(\varphi(x)), \psi(\varphi(y))\big) \\
    &\ge \frac{1}{2Q+1}\big(d_{X}(x,y)-d_{X}(\psi(\varphi(x)),x)-d_{X}(y,\psi(\varphi(y)))\big)\\
    &\ge \frac{1}{2Q+1}d_{X}(x,y)-\frac{2Q}{2Q+1}
\end{align*}
for any $x,y\in X$, so $\varphi \colon X\rightarrow Y$ is a quasi-isometry, with $\psi$ as a quasi-inverse. This shows the first part of the proposition. 

\noindent To prove the second part, we consider the map 
\begin{align*}
    f\colon \mathcal{L}_{n}(X,\mathcal{C}) &\longrightarrow \mathcal{L}_{n}(Y) \\
    (c,x)&\longrightarrow (c^{\flat}, x_{C}), \; \text{where $C\in\mathcal{C}$ contains $x$}
\end{align*}
and where $c^{\flat}$ is the colouring of $Y$ naturally defined by $c^{\flat}(x_{P})\defeq c(P)$, for any $P\in\mathcal{C}$. The other way around, define 
\begin{align*}
    g\colon \mathcal{L}_{n}(Y) &\longrightarrow  \mathcal{L}_{n}(X,\mathcal{C}) \\
    (c,x_{C})&\longmapsto (\Tilde{c}, x_{C})
\end{align*}
where $\Tilde{c}$ is the colouring of $X$ naturally defined by $\Tilde{c}(x)\defeq c(x_{C})$ where $C\in\mathcal{C}$ is the piece containing $x$. Obviously, one has $\Tilde{c^{\flat}}=c$ (resp. $\Tilde{c}^{\flat}=c$) for any $c\in \Z_{n}^{(X,\mathcal{C})}$ (resp. $c\in\Z_{n}^{(Y)})$, so this already tells us that $f\circ g=\text{Id}_{\mathcal{L}_{n}(Y)}$. Also, for $(c,x)\in\mathcal{L}_{n}(X,\mathcal{C})$ and $C\in\mathcal{C}$ the zone containing $x$, one has 
\begin{equation*}
    d_{\mathcal{L}_{n}(X,\mathcal{C})}\big((c,x), g(f(c,x))\big)=d_{\mathcal{L}_{n}(X,\mathcal{C})}\big((c,x),(c,x_{C})\big)=d_{X}(x,x_{C})\le Q
\end{equation*}
so $g\circ f$ is at distance $\le Q$ from $\text{Id}_{\mathcal{L}_{n}(X,\mathcal{C})}$. We now show that both $f$ and $g$ are Lipschitz.

\noindent Let $a,b\in\mathcal{L}_{n}(X,\mathcal{C})$ be two adjacent vertices. As usual, two cases can occur:
\begin{itemize}
    \item Either $a=(c,x)$ and $b=(c',x)$, where $c,c'$ only differ on the zone $C\in\mathcal{C}$ with $x\in C$. Then
    \begin{equation*}
        d_{\mathcal{L}_{n}(Y)}(f(a),f(b))=d_{\mathcal{L}_{n}(Y)}\big((c^{\flat},x_{C}), ((c')^{\flat}, x_{C})\big)=1
    \end{equation*}
    since $c^{\flat}$, $(c')^{\flat}$ only differ on the vertex $x_{C}\in Y$. 
    \item Either $a=(c,x)$ and $b=(c,y)$ where $x,y$ are adjacent in $X$. If $x$ and $y$ are in the same piece of the partition one has $d_{\mathcal{L}_{n}(Y)}(f(a),f(b))=0$, whereas if $x\in C$ and $y\in C'$ for $C\neq C'$, then 
    \begin{equation*}
        d_{\mathcal{L}_{n}(Y)}(f(a),f(b))=d_{\mathcal{L}_{n}(Y)}\big((c^{\flat},x_{C}),(c^{\flat}, x_{C'})\big)=d_{Y}(x_{C},x_{C'})=1.
    \end{equation*}
\end{itemize}
In both cases we conclude that $d_{\mathcal{L}_{n}(Y)}(f(a),f(b))\le 1$, so $f$ is $1-$Lipschitz according to Lemma~\ref{lem:Lipschitzbetweengraphs}. 

\noindent Let now $a,b\in \mathcal{L}_{n}(Y)$ be adjacent vertices. We treat two cases:
\begin{itemize}
    \item Assume that $a=(c,x_{C})$ and $b=(c',x_{C})$ where $c,c'$ only differ on $x_{C}\in Y$. Then one has 
    \begin{equation*}
        d_{\mathcal{L}_{n}(X,\mathcal{C})}(g(a), g(b))=d_{\mathcal{L}_{n}(X,\mathcal{C})}\big((\Tilde{c},x_{C}),(\Tilde{c'},x_{C})\big)=1
    \end{equation*}
    since $\Tilde{c}, \Tilde{c'}$ only differ on the zone $C$ containing $x_{C}$. 
    \item Assume now $a=(c,x_{C})$ and that $b=(c,x_{C'})$ with $x_{C}, x_{C'}$ adjacent in $Y$. Then $C,C'$ contain adjacent vertices in $X$, say $u\in C$ and $v\in C'$, so that
    \begin{align*}
        d_{\mathcal{L}_{n}(X,\mathcal{C})}(g(a), g(b))&=d_{\mathcal{L}_{n}(X,\mathcal{C})}\big((\Tilde{c},x_{C}), (\Tilde{c},x_{C'})\big) \\
        &=d_{X}(x_{C},x_{C'}) \\
        &\le d_{X}(x_{C},u)+d_{X}(u,v)+d_{X}(v,x_{C'}) \\
        &\le 2Q+1.
    \end{align*}
\end{itemize}
In both cases, we conclude that $g$ sends adjacent vertices in $\mathcal{L}_{n}(Y)$ to vertices at distance at most $2Q+1$ in $\mathcal{L}_{n}(X,\mathcal{C})$. Applying one last time Lemma~\ref{lem:Lipschitzbetweengraphs} it follows that $g$ is $(2Q+1)-$Lipschitz. We have then
\begin{align*}
    d_{\mathcal{L}_{n}(X,\mathcal{C})}(a,b) \ge d_{\mathcal{L}_{n}(Y)}(f(a),f(b)) &\ge \frac{1}{2Q+1}d_{\mathcal{L}_{n}(X,\mathcal{C})}\left(g(f(a)), g(f(b))\right) \\
    &\ge \frac{1}{2Q+1}d_{\mathcal{L}_{n}(X,\mathcal{C})}(a,b)-\frac{2Q}{2Q+1}
\end{align*}
for any $a,b\in \mathcal{L}_{n}(X,\mathcal{C})$. Thus $f$ is a quasi-isometry with $g$ as a quasi-inverse. 
\end{proof}

\subsection{Leaf-preserving, aptolic and non-aptolic quasi-isometries}\label{subsubsection3.4.2} In~\cite{GT24b}, Genevois and Tessera established a strong rigidity result for quasi-isometries between standard wreath products and, later, for more general halo products~\cite[Corollary~6.11]{GT24a}. More precisely, they proved that a quasi-isometry between two halo products (satisfying additional properties) always lies at finite distance from an aptolic quasi-isometry. Perhaps surprisingly, we will show in this section that such a rigidity does not hold in the permutational case (see Proposition~\ref{prop:non-aptolicQI}). This will motivate the strategy developed in subsequent sections.

We begin by defining the leaf-preservingness property for permutational wreath products.

\begin{definition}
A quasi-isometry $q\colon E\wr_{G/M}G\longrightarrow F\wr_{H/N}H$ is~\textit{leaf-preserving} if it sends $G-$cosets to $H-$cosets and admits a quasi-inverse that sends $H-$cosets to $G-$cosets. 
\end{definition}

From our embedding theorem we deduce that:

\begin{corollary}\label{cor:QIareleafpreserving}
Let $E$ and $F$ be non-trivial finite groups. Let $G,H$ be finitely presented groups with normal finitely generated infinite subgroups $M\lhd G$, $N\lhd H$. Suppose that $M$ has infinite index in $G$ and that $G$ (resp. $H$) is not coarsely separable by any collection of subspaces that uniformly quasi-isometrically embed into $M$ (resp. $N$). If $q\colon E\wr_{G/M}G \longrightarrow F\wr_{H/N}H$ is a quasi-isometry, then $q$ is (up to finite distance) leaf-preserving. 
\end{corollary}

The proof of this corollary requires, in addition to the embedding theorem, to understand better the coarse intersection of two leaves in a permutational wreath product. This is the goal of the next two lemmas. 

\begin{lemma}\label{lem:Hausdorffdistancebetweencosets}
Let $G$ be a finitely generated group with an infinite normal subgroup $M\lhd H$. Endow $G$ with a word metric $d_{G}$ coming from a finite generating set $S$. Then we have 
\begin{equation*}
    d_{\text{Haus}}(gM, hM) \le d_{G}(gM, hM)
\end{equation*}
for any $g,h\in G$.
\end{lemma}

\begin{proof}
Without restrictions, we assume that $S$ is symmetric. Let $x\in gM$, $y\in hM$ be such that $d_{G}(x,y)=d_{G}(gM, hM)$, and let $\gamma$ be a geodesic in $\text{Cay}(G,S)$ from $x$ to $y$. Let $u\in gM$, and consider in $\text{Cay}(G,S)$ the path $ux^{-1}\gamma$:  

\begin{figure}[H]
  \centering
  \includegraphics[width=0.7\linewidth]{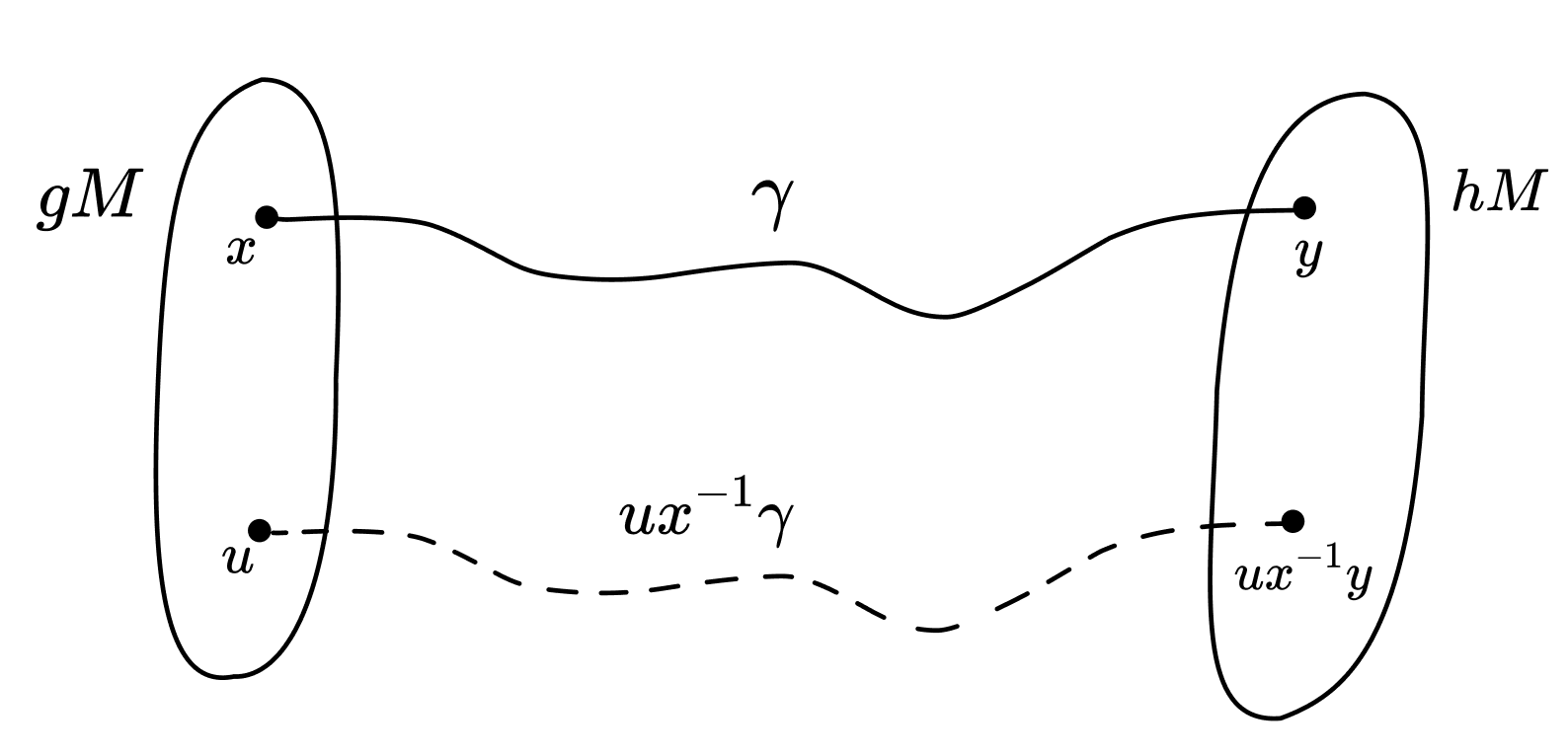}
  \caption{Configuration from the proof of Lemma~\ref{lem:Hausdorffdistancebetweencosets}}
\end{figure}

\noindent It connects $ux^{-1}x=u$ to $ux^{-1}y$, and we claim that the latter belongs to $hM$. Indeed, 
\begin{align*}
    ux^{-1}y\in hM &\Longleftrightarrow ux^{-1}y\in yM \\
    &\Longleftrightarrow ux^{-1}\in yMy^{-1}\\
    &\Longleftrightarrow ux^{-1}\in M
\end{align*}
and this last assertion is true since $uM=xM$. Thus any element of $gM$ is connected to an element of $hM$ by a path of length $\text{length}(\gamma)=d_{G}(gM, hM)$. A symmetric reasoning shows that any element of $hM$ can be connected to an element of $gM$ by a path of length $\text{length}(\gamma)=d_{G}(gM, hM)$. We then conclude that 
\begin{equation*}
    d_{\text{Haus}}(gM, hM) \le d_{G}(gM, hM)
\end{equation*}
and the proof is complete.  
\end{proof}

\begin{lemma}\label{lem:coarseintersectionsofzonesisaneighborhoodof1zone}
Let $E$ be a non-trivial finite group, and let $G$ be a finitely generated group with a normal finitely generated subgroup $M\lhd G$. For every distinct colourings $c,d\in E^{(G,\mathcal{C}_{M})}$, for every $R\ge 0$, there is a constant $Q\ge 0$ such that
\begin{equation*}
    (cG)\cap (dG)^{+R} \subset \left\lbrace (c,p)\in E\wr_{G/M}G : p\in M^{+Q}\right\rbrace.
\end{equation*}
\end{lemma}

\begin{proof}
Fix $S$ a finite generating set of $G$, two distinct colourings $c,d\in E^{(G,\mathcal{C}_{M})}$ and $R\ge 0$. Let $(c,p)\in (cG)\cap (dG)^{+R}$. Thus there exists $(d,q)\in dG$ such that $d_{E\wr_{G/M}G}\big((c,p),(d,q)\big) \le R$, which implies that $c^{-1}d$ is supported on $M-$cosets at distance $\le R$ from $p$. In other words, if we write explicitly $\text{supp}(c^{-1}d)=\lbrace g_{1}M,\dots, g_{r}M\rbrace$, then $p\in (g_{i}M)^{+R}$ for any $1\le i\le r$. Now, using Lemma~\ref{lem:Hausdorffdistancebetweencosets}, one has 
\begin{equation*}
    d_{\text{Haus}}(M,g_{i}M)\le d_{G}(M,g_{i}M)=d_{G}(M,Mg_{i}) \le \ell_{S}(g_{i})
\end{equation*}
for any $1\le i\le r$, and thus $p\in M^{+(R+\ell_{S}(g_{i}))}$ for any $1\le i\le r$. Taking $Q\defeq R+\max_{1\le i\le r}\ell_{S}(g_{i})$, the conclusion follows. 
\end{proof}

\begin{proof}[Proof of Corollary~\ref{cor:QIareleafpreserving}]
Fix $q\colon E\wr_{G/M}G \longrightarrow F\wr_{H/N}H$ an $(A,B)-$quasi-isometry and $\overline{q}$ one of its quasi-inverses such that $q\circ\overline{q}$ and $\overline{q}\circ q$ are within distance $\le B$ from the identities. Notice that, since $M$ has infinite index, $E\wr_{G/M}G$ is infinitely presented by Lemma~\ref{lem:permutationallamplightersareinfinitelypresented}, so $F\wr_{H/N}H$ is also infinitely presented, and thus $N$ also has infinite index in $H$ (otherwise, as $F$ is finite, $F\wr_{H/N}H=F^{H/N}\rtimes H$ would be quasi-isometric to $H$, and thus would be finitely presented). By Theorem~\ref{thm:embeddingthmPWP}, $q(G)$ lies in a neighbourhood of an $H-$coset, say $H$ itself. Conversely, $\overline{q}(H)$ lies in a neighbourhood of a $G-$coset, say $aG$. It follows that $\overline{q}(q(G))$ lies in a neighbourhood of $aG$. Since the Hausdorff distance between $G$ and $\overline{q}(q(G))$ is at most $B$, this in turn implies that $G$ lies in a neighbourhood of $aG$, say $(aG)^{+K}$. 
We now claim that this implies $aG=G$. Indeed, suppose for a contradiction that these two leaves are distinct. Then, by Lemma~\ref{lem:coarseintersectionsofzonesisaneighborhoodof1zone}, there is a constant $Q\ge 0$ such that
\begin{equation*}
    G=G\cap (aG)^{+K}\subset \lbrace (\mathbf{1}, p)\in E\wr_{G/M}G : p\in M^{+Q}\rbrace.
\end{equation*}
In other words, $M$ is quasi-dense in $G$, which implies that $M$ has finite index in $G$, contrary to our assumption. We thus deduce that $aG=G$. The same reasoning shows that any other leaf $cG$ is sent close to a leaf $\alpha(c)H$ by $q$, which is sent back close to $cG$ by $\overline{q}$. It follows that $q$ is leaf-preserving. 
\end{proof}

Here is an example that shows that the corollary does not hold anymore if one drops the assumption on the coarse separation of the base groups.

\begin{example}\label{ex:aQIwhichisnotleafpreserving}
Consider $E=\lbrace -1,1\rbrace$ the cyclic group of order $2$ and $G=\Z^2$, $M=\Z$. Let $S=\lbrace \pm (1,0), \pm(0,1)\rbrace$ denote the standard generating set of $G$. Clearly $G$ is coarsely separated by $M$, so Corollary~\ref{cor:QIareleafpreserving} does not apply. Consider the map 
\begin{align*}
    \varphi\colon E\wr_{G/M}G &\longrightarrow E\wr_{G/M}G \\
    (c,p)&\longmapsto \begin{cases}
        (c,p) &\mbox{if $p\in ([0,+\infty)\cap\Z)\times\Z$} \\
        (c\delta_{0},p) &\mbox{otherwise}
    \end{cases}
\end{align*}
where $\delta_{0}$ denotes the colouring having only coloured the $\Z-$coset containing $(0,0)$. By construction, $\varphi$ is not leaf-preserving (even up to finite distance), and we show it is indeed a quasi-isometry.

Indeed, let $a$ and $b$ be adjacent vertices in $\text{Cay}(E\wr_{G/M}G, E\cup S)$. We consider two cases. Assume first that $a=(c_{1},p)$ and $b=(c_{2},p)$ where $c_{1},c_{2}$ only differ on the coset containing $p$. If $p\in ([0,+\infty)\cap\Z)\times\Z$ then $\varphi(a)$ and $\varphi(b)$ are adjacent in $\text{Cay}(E\wr_{G/M}G, E\cup S)$; and otherwise $\varphi(a)=(c_{1}\delta_{0},p)$, $\varphi(b)=(c_{2}\delta_{0},p)$ are adjacent as well since $c_{1}\delta_{0}$, $c_{2}\delta_{0}$ only differ on the coset containing $p$.

Assume now $a=(c,p)$ and $b=(c,p+s)$ for some $s\in S$. If $p\in \lbrace -1\rbrace\times\Z$ and $p+s\in\lbrace 0\rbrace\times\Z$, then $\varphi(a)=(c\delta_{0},p)$ while $\varphi(b)=(c, p+s)$, so 
\begin{equation*}
    d(\varphi(a),\varphi(b))=2.
\end{equation*}
One gets the same value in the symmetric situation $p\in\lbrace 0\rbrace\times\Z$, $p+s\in \lbrace -1\rbrace\times\Z$. In the other cases, one easily checks that $\varphi$ sends adjacent vertices to adjacent vertices.

Hence we conclude from Lemma~\ref{lem:Lipschitzbetweengraphs} that $\varphi$ is $2-$Lipschitz. Now, observe that $\varphi$ is a bijection, whose inverse is 
\begin{align*}
    \psi\colon E\wr_{G/M}G &\longrightarrow E\wr_{G/M}G \\
    (c,p)&\longmapsto \begin{cases}
        (c,p) &\mbox{if $p\in ([0,+\infty)\cap\Z)\times\Z$} \\
        (c\delta_{0}^{-1},p) &\mbox{otherwise}
    \end{cases}
\end{align*}
and a similar computation shows that $\psi$ is $2-$Lipschitz. Thus it follows that 
\begin{equation*}
    \frac{1}{2}\cdot d(a,b)\le d(\varphi(a),\varphi(b))\le 2\cdot d(a,b)
\end{equation*}
for any $a,b\in E\wr_{G/M}G$, so $\varphi$ is a biLipschitz equivalence. 
\end{example}

We now turn to (non-)aptolic quasi-isometries. Let us first adapt the terminology from~\cite{GT24a, GT24b} in our case. 

\begin{definition}
A quasi-isometry $q\colon E\wr_{G/M}G\longrightarrow F\wr_{H/N}H$ is of~\textit{aptolic form} if there exist two maps $\alpha \colon E^{(G,\mathcal{C}_{M})}\longrightarrow F^{(H,\mathcal{C}_{N})}$, $\beta\colon G\rightarrow H$ such that $q(c,p)=(\alpha(c), \beta(p))$ for any $(c,p)\in E\wr_{G/M}G$. Moreover, $q$ is~\textit{aptolic} if it is of aptolic form and it has a quasi-inverse of aptolic form. 
\end{definition}

The following proposition indicates that a strategy similar to the one followed in~\cite{GT24a} cannot be replicated in the context of permutational wreath products.

\begin{proposition}\label{prop:non-aptolicQI}
Let $F$ be a finite group. Let $H$ be a finitely generated group, $N\leqslant H$ an infinite index subgroup containing an element $z\in N$ of infinite order which is central in $H$. Then the map 
\begin{align*}
    \varphi\colon F\wr_{H/N}H &\longrightarrow F\wr_{H/N}H \\
    (c,p)&\longmapsto (c,pz^{|\text{supp}(c)|})
\end{align*}
is a leaf-preserving quasi-isometry which is not at finite distance from an aptolic quasi-isometry.
\end{proposition}

\begin{proof}
Denote $S$ a finite symmetric generating set of $H$, and let $G\defeq F\wr_{H/N}H$. We first show that $\varphi$ is Lipschitz. Fix $a$ and $b$ two adjacent vertices in $\text{Cay}(G, F\cup S)$. There are two cases to consider:
\begin{itemize}
    \item Assume first that $a=(c,p)$ and $b=(c,ps)$ for some $s\in S$. Then 
    \begin{align*}
        d_{G}(\varphi(a),\varphi(b))&=d_{G}\big((c,pz^{|\text{supp}(c)|}), (c,psz^{|\text{supp}(c)|})\big) \\
        &=d_{H}(pz^{|\text{supp}(c)|}, pz^{|\text{supp}(c)|}s) \\
        &=d_{H}(1_{H},s) \\
        &=1 \\
        &\le 1+\ell_{H}(z)
    \end{align*}
    using Lemma~\ref{lem:leavesareconvex} and the fact that $z$ (or any of its power) is central for the second equality.  
    \item Assume now that $a=(c,p)$ and $b=(c',p)$ where $c,c'$ only differ on the $N-$coset containing $p$. Then 
    \begin{align*}
        d_{G}(\varphi(a),\varphi(b))&=d_{G}\big((c,pz^{|\text{supp}(c)|}), (c',pz^{|\text{supp}(c')|})\big) \\
        &=1+d_{H}(pz^{|\text{supp}(c)|}, pz^{|\text{supp}(c')|}) \\
        &=1+d_{H}(z^{|\text{supp}(c)|}, z^{|\text{supp}(c')|}) \\
        &\le 1+\ell_{H}(z)
    \end{align*}
    using in the last equality that $|\text{supp}(c')|\in\big\lbrace |\text{supp}(c)|-1, |\text{supp}(c)|, |\text{supp}(c)|+1\big\rbrace$.
\end{itemize}
We conclude from Lemma~\ref{lem:Lipschitzbetweengraphs} that $\varphi$ is $\big(1+\ell_{H}(z)\big)-$Lipschitz. Now, notice that $\varphi$ is a bijection, whose inverse is 
\begin{align*}
    \psi\colon G&\longrightarrow G \\
    (c,p)&\longmapsto (c,pz^{-|\text{supp}(c)|}).
\end{align*}
In particular, $\varphi$ is coarsely surjective, and a similar computation as the one above shows that $\psi$ is also $\big(1+\ell_{H}(z)\big)-$Lipschitz. Therefore
\begin{equation*}
    d_{G}(x,y)=d_{G}\big(\psi(\varphi(x)), \psi(\varphi(y))\big) \le \big(1+\ell_{H}(z)\big)\cdot d_{G}(\varphi(x), \varphi(y))
\end{equation*}
for any $x,y\in G$, whence 
\begin{equation*}
    \frac{1}{1+\ell_{H}(z)}\cdot d_{G}(x,y) \le d_{G}(\varphi(x), \varphi(y)) \le \big(1+\ell_{H}(z)\big)\cdot d_{G}(x,y)
\end{equation*}
for all $x,y\in G$. Thus $\varphi\colon G\rightarrow G$ is a biLipschitz equivalence. 

\noindent Next, $\varphi$ is leaf-preserving by construction, so we are only left to prove that $\varphi$ is not aptolic, up to finite distance. Towards a contradiction, assume there is $R\ge 0$ and 
\begin{align*}
    q\colon G&\longrightarrow G \\
    (c,p)&\longmapsto (\alpha(c),\beta(p))
\end{align*}
an aptolic self-quasi-isometry of $G$ such that $d(q,\varphi)\le R$. Then, for any $(c,p)\in G$, denoting $\beta_{c}\colon H\rightarrow H$, $p\longmapsto pz^{|\text{supp}(c)|}$,
one gets
\begin{align*}
    d_{H}(\beta_{c}(p),\beta(p))&\le d_{G}\big((c,\beta_{c}(p)), (\alpha(c),\beta(p))\big) \\
    &=d_{G}\big(q(c,p), \varphi(c,p)\big) \\
    &\le d(q,\varphi) \\
    &\le R
\end{align*}
and it follows that $d(\beta_{c},\beta)\le R$. In particular, all $\beta_{c}$ lie at distance at most $2R$ from $\beta_{\mathbf{1}}=\text{Id}_{H}$. On the other hand, since $z$ has infinite order and balls in $\text{Cay}(H,S)$ are finite, there is an integer $m\in\N$ such that $d_{H}(z^{m},1_{H})>2R$, and thus, fixing an arbitrary colouring $c$ supported on $m$ cosets of $N$, we see that 
\begin{equation*}
    d_{H}(\beta_{c}(p),p)=d_{H}(p, pz^{|\text{supp}(c)|})=d_{H}(1_{H},z^{m})>2R 
\end{equation*}
for any $p\in H$, whence $d(\beta_{c}, \text{Id}_{H})>2R$. This contradiction proves that $\varphi$ is not at finite distance from an aptolic quasi-isometry, concluding the proof.
\end{proof}

\begin{remark}
In this proof, the fact that $F$ is finite plays no role, and one can show the same result for $F$ an arbitrary finitely generated group, as soon as one modifies the definition of lamplighter graphs over partitioned graphs as follows: declare now that the colourings take values in another graph $Y$, and that two pairs $(c,p)$, $(c',p)$ are adjacent if $c,c'$ only differ on the zone $C$ containing $p$ and that $c(C), c'(C)$ are adjacent in $Y$. 
\end{remark}

Based on the same idea, one can exhibit plenty of leaf-preserving non-aptolic quasi-isometries of various permutational halo products, showing that the general rigidity result proved in~\cite[Corollary~6.11]{GT24a} fails in the permutational case.

\subsection{Preserving zones}\label{subsubsection3.4.3} The goal of this subsection is now to deduce from the property of preserving the leaves that in fact any quasi-isometry between two permutational wreath products must also quasi-preserve zones inside a single leaf. More formally:

\begin{proposition}\label{prop:leafpreservingQIareQIofpairs}
Let $E,F$ be non-trivial finite groups, and $G,H$ be finitely generated groups, with normal infinite subgroups $M\lhd G$, $N\lhd H$. Suppose that $M$ has infinite index in $G$. If a quasi-isometry $q\colon E\wr_{G/M}G \longrightarrow F\wr_{H/N}H$ is leaf-preserving, then it is a quasi-isometry of pairs 
\begin{equation*}
    (E\wr_{G/M}G,M)\longrightarrow (F\wr_{H/N}H, N).
\end{equation*}
\end{proposition}

\begin{proof} 
Fix $A\ge 1$, $B\ge 0$ such that $q$ and a quasi-inverse $q'\colon F\wr_{H/N}H\longrightarrow E\wr_{G/M}G$ are both $(A,B)-$quasi-isometries. As $q,q'$ are leaf-preserving we may write 
\begin{equation*}
    q(c,p)=(\alpha(c),\beta_{c}(p)), \; (c,p)\in E\wr_{G/M}G
\end{equation*}
for some maps $\alpha\colon E^{(G,\mathcal{C}_{M})}\longrightarrow F^{(H,\mathcal{C}_{N})}$ and $\beta_{c}\colon G\rightarrow H$, $c\in E^{(G,\mathcal{C}_{M})}$; as well as 
\begin{equation*}
    q'(c,p)=(\alpha'(c), \beta_{c}'(p)),\; (c,p)\in F\wr_{H/N}H
\end{equation*}
for some $\alpha'\colon F^{(H,\mathcal{C}_{N})}\longrightarrow E^{(G,\mathcal{C}_{M})}$ and $\beta_{c}'\colon H\rightarrow G$, $c\in F^{(H,\mathcal{C}_{N})}$. Let us record the following observation for future use:
\begin{claim}\label{claim:claim3.50}
The maps $\alpha$, $\alpha'$ are bijections, inverses of each other, and for any $c\in E^{(G,\mathcal{C}_{M})}$, $\beta_{c}\colon G\rightarrow H$ is an $(A,B)-$quasi-isometry, with $\beta_{\alpha(c)}'\colon H\rightarrow G$ as a quasi-inverse.
\end{claim}

\renewcommand{\qedsymbol}{$\blacksquare$}
\begin{proof}[Proof of Claim~\ref{claim:claim3.50}]
Fix a finitely supported colouring $c\in E^{(G,\mathcal{C}_{M})}$. Then, for any $p\in G$, we have 
\begin{equation}\label{eq3.4}
    d\big((c,p), (\alpha'\circ\alpha(c), \beta'_{\alpha(c)}\circ\beta_{c}(p))\big)=d\big((c,p), (q'\circ q)(c,p)\big)\le B
\end{equation}
which implies that $c$ and $\alpha'\circ\alpha(c)$ may only differ on zones at distance at most $B$ from $p$, for any $p\in G$. As $M$ has infinite index in $G$, it follows that $\alpha'\circ\alpha(c)=c$. Similarly, using that $q\circ q'$ is at distance $\le B$ from $\text{Id}_{F\wr_{H/N}H}$, we deduce also $\alpha\circ\alpha'(c)=c$ for all $c\in F^{(H,\mathcal{C}_{N})}$. Thus $\alpha$ is a bijection, whose inverse is $\alpha'$. 

\noindent Now let $p_{1},p_{2}\in G$. Then 
\begin{align*}
    d_{H}(\beta_{c}(p_{1}), \beta_{c}(p_{2})) &= d\big((\alpha(\mathbf{1}), \beta_{c}(p_{1})),(\alpha(\mathbf{1}), \beta_{c}(p_{2}))\big) \\
    &=d\big(q(\mathbf{1},p_{1}),q(\mathbf{1},p_{2})\big) \\
    &\le A\cdot d\big((\mathbf{1},p_{1}), (\mathbf{1},p_{2})\big)+B \\
    &=A\cdot d_{G}(p_{1},p_{2})+B
\end{align*}
and similarly one gets 
\begin{equation*}
    d_{H}(\beta_{c}(p_{1}), \beta_{c}(p_{2})) \ge \frac{1}{A}\cdot d_{G}(p_{1},p_{2})-B.
\end{equation*}
Thus $\beta_{c}$ is an $(A,B)-$quasi-isometric embedding, and additionally for any $h\in H$ there exists $(c,p)\in E\wr_{G/M}G$ such that $d\big(q(c,p), (0, h)\big)\le B$, whence 
\begin{equation*}
    d_{H}(\beta_{c}(p), h) \le d\big((\alpha(c), \beta_{c}(p)), (\mathbf{1},h)\big)=d\big(q(c,p), (\mathbf{1}, h)\big)\le B.
\end{equation*}
Hence $\beta_{c}$ is coarsely surjective as well, and it is then an $(A,B)-$quasi-isometry  $G\rightarrow H$. Likewise, we prove that $\beta'_{\alpha(c)}\colon H\rightarrow G$ is an $(A,B)-$quasi-isometry, and from (\ref{eq3.4}) it follows that $\beta'_{\alpha(c)}\circ\beta_{c}$ is at distance $\le B$ from $\text{Id}_{G}$. Using similarly that $q\circ q'$ is at distance $\le B$ from the identity, we deduce that $\beta_{c}\circ \beta'_{\alpha(c)}$ is at distance $\le B$ from $\text{Id}_{H}$. We conclude that $\beta_{c}$ and $\beta'_{\alpha(c)}$ are quasi-inverses of each other.
\end{proof}
\renewcommand{\qedsymbol}{$\Box$}

\noindent To prove the statement, it is enough to prove that each $\beta_{c}\colon G\rightarrow H$ is a quasi-isometry of pairs $(G,M)\longrightarrow (H,N)$ (here we identify geometrically the leaf $L(c)$ (resp. $L(\alpha(c))$) in $E\wr_{G/M}G$ (resp. in $F\wr_{H/N}H$) with $G$ (resp. $H$)).

\noindent Fix then a colouring $c\in E^{(G,\mathcal{C}_{M})}$, a zone $gM\subset L(c)$, and let $c'\in E^{(G,\mathcal{C}_{M})}$ be a colouring that differs from $c$ only on the zone $gM$. Let $p\in gM$. One has then 
\begin{equation*}
    d\big((\alpha(c),\beta_{c}(p)), (\alpha(c'),\beta_{c'}(p))\big)=d\big(q(c,p), q(c',p)\big) \le A\cdot d\big((c,p),(c',p)\big)+B=A+B
\end{equation*}
since $d\big((c,p),(c',p)\big)=1$. This inequality implies that 
\begin{equation*}
    \text{supp}\left(\alpha(c)^{-1}\alpha(c')\right) \subset \lbrace \beta_{c}(p)\rbrace_{+(A+B)}
\end{equation*}
so if we write explicitly $\text{supp}\left(\alpha(c)^{-1}\alpha(c')\right)=\lbrace h_{1}N,\dots, h_{r}N\rbrace$ where $r\ge 1$, we get 
\begin{equation}\label{eq3.5}
    d_{H}(h_{i}N, \beta_{c}(p))\le A+B 
\end{equation}
for any $1\le i\le r$. Thus $\beta_{c}(p)\in \dis \bigcap_{i=1}^{r}(h_{i}N)^{+(A+B)}$. Since this holds for any $p\in gM$, we have proved that
\begin{equation*}
    \beta_{c}(gM)\subset \bigcap_{i=1}^{r}(h_{i}N)^{+(A+B)}.
\end{equation*}
Now, for any $1\le i\le r$, it also follows from (\ref{eq3.5}) and from Lemma~\ref{lem:Hausdorffdistancebetweencosets} that 
\begin{equation*}
    d_{\text{Haus}}(h_{i}N, \beta_{c}(g)N) \le A+B
\end{equation*}
for any $1\le i\le r$, so that
\begin{equation*}
    \beta_{c}(gM)\subset \bigcap_{i=1}^{r}(h_{i}N)^{+(A+B)} \subset \bigcap_{i=1}^{r}\left((\beta_{c}(g)N)^{+(A+B)}\right)^{+(A+B)}\subset\big(\beta_{c}(g)N\big)^{+2(A+B)}
\end{equation*}
using Lemma~\ref{lem:neighborhoodsandQI}\textit{(ii)}. Similarly, one shows that given any zone $hN\subset L(\alpha(c))$, we have 
\begin{equation*}
    \beta_{\alpha(c)}'(hN)\subset \big(\beta_{\alpha(c)}'(h)M\big)^{+2(A+B)}.
\end{equation*}
In particular, when $h=\beta_{c}(g)$, this gives
\begin{equation*}
    \beta_{\alpha(c)}'(\beta_{c}(g)N)\subset \left((\beta_{\alpha(c)}'(\beta_{c}(g))M\right)^{+2(A+B)}
\end{equation*}
and applying $\beta_{c}$ to this inclusion and using Lemma~\ref{lem:neighborhoodsandQI}\textit{(iii)}, and that $\beta_{c}\circ\beta_{\alpha(c)}'$ is at distance $\le B$ from $\text{Id}_{H}$, it follows that 
\begin{align*}
    \beta_{c}(g)N &\subset \left(\beta_{c}\circ\beta_{\alpha(c)}'(\beta_{c}(g)N)\right)^{+B} \\
    &\subset \beta_{c}\left(\left(\beta_{\alpha(c)}'(\beta_{c}(g))M\right)^{+2(A+B)}\right)^{+B} \\
    &\subset \left(\beta_{c}\left(\beta_{\alpha(c)}'(\beta_{c}(g))M\right)^{+A(2(A+B))+B}\right)^{+B} \\
    &\subset \beta_{c}\left(\beta_{\alpha(c)}'(\beta_{c}(g))M\right)^{+(2A(A+B)+2B)}.
\end{align*}
It just remains to notice that, since $\beta_{\alpha(c)}'(\beta_{c}(g))$ and $g$ are at distance at most $B$ in $G$, Lemma~\ref{lem:Hausdorffdistancebetweencosets} ensures that the Hausdorff distance between $\beta_{\alpha(c)}'(\beta_{c}(g))M$ and $gM$ is at most $B$ as well, and we finally get 
\begin{align*}
    \beta_{c}(g)N  &\subset \beta_{c}\left(\beta_{\alpha(c)}'(\beta_{c}(g))M\right)^{+(2A(A+B)+2B)} \\
    &\subset \beta_{c}((gM)^{+B})^{+(2A(A+B)+2B)} \\
    &\subset \beta_{c}(gM)^{+((2A(A+B)+2B)+AB+B)}.
\end{align*}
We conclude that 
\begin{equation*}
    d_{\text{Haus}}\big(\beta_{c}(gM), \beta_{c}(g)N\big) \le 2A(A+B)+(A+3)B
\end{equation*}
showing that $\beta_{c}$ is indeed a quasi-isometry of pairs $(G,M)\longrightarrow (H,N)$.
\end{proof}

\begin{corollary}\label{cor:QIareQIofpairs}
Let $E$ and $F$ be non-trivial finite groups. Let $G,H$ be finitely presented groups with normal finitely generated infinite subgroups $M\lhd G$, $N\lhd H$. Suppose that $M$ has infinite index in $G$ and that $G$ (resp. $H$) is not coarsely separable by any collection of subspaces that uniformly quasi-isometrically embed into $M$ (resp. $N$). Then any quasi-isometry $E\wr_{G/M}G \longrightarrow F\wr_{H/N}H$ is a quasi-isometry of pairs
\begin{equation*}
     (E\wr_{G/M}G,M)\longrightarrow (F\wr_{H/N}H, N).
\end{equation*}
\end{corollary}

\begin{proof}
By Corollary~\ref{cor:QIareleafpreserving}, such a quasi-isometry can be assumed to be leaf-preserving, so the conclusion follows from Proposition~\ref{prop:leafpreservingQIareQIofpairs}. 
\end{proof}

\begin{figure}[H]
  \centering
  \includegraphics[width=0.7\linewidth]{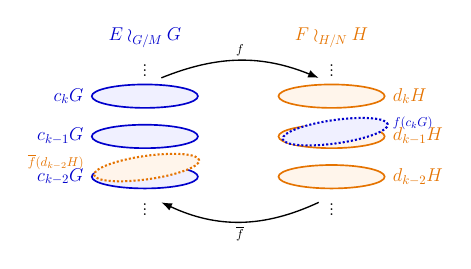}
  \caption{The behaviour of an arbitrary quasi-isometry $E\wr_{G/M}G\longrightarrow F\wr_{H/N}H$.}
\end{figure}

\begin{remark}
In the realm of quasi-isometric rigidity, an interesting and challenging question, already raised in~\cite[Question~2.3]{HMS21}, is the following: given finitely generated groups $G,H$ and a collection $\mathcal{P}$ of subgroups of $G$, under which conditions does there exist a collection $\mathcal{Q}$ of subgroups of $H$ such that any quasi-isometry $G\rightarrow H$ extends to a quasi-isometry of pairs $(G,\mathcal{P})\longrightarrow (H,\mathcal{Q})$? Corollary~\ref{cor:QIareQIofpairs} shows that permutational wreath products of the form $E\wr_{G/M}G$, with our running assumptions, provide a class of groups answering positively to this question. 
\end{remark}

\subsection{Coning-off graphs}\label{subsubsection3.4.4} In this part, we introduce our main tool that will help us to deduce some valuable information from the existence of a quasi-isometry between two permutational wreath products, even if the latter is not aptolic. 

\begin{definition}
Let $X$ be a graph with a partition $\mathcal{C}$ of $V(X)$. The~\textit{cone-off of $X$ over $\mathcal{C}$}, denoted $\text{CO}(X,\mathcal{C})$, is the graph obtained from $X$ by adding an edge between two vertices whenever they belong to a common piece $C\in\mathcal{C}$. 
\end{definition}

\begin{figure}[H]
  \centering
  \includegraphics[width=0.75\linewidth]{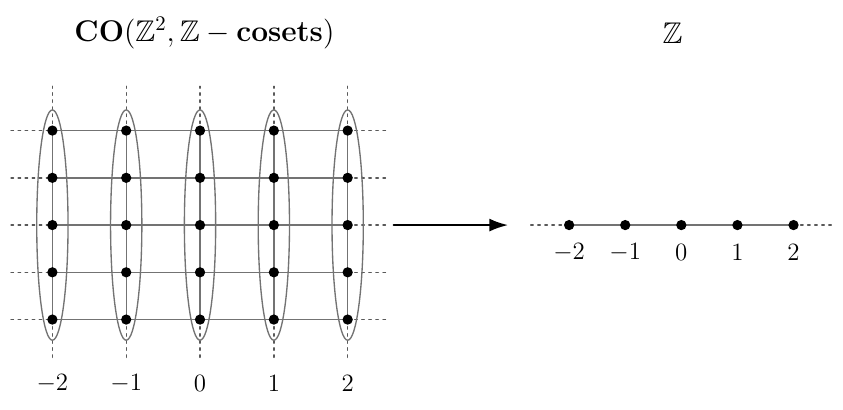}
  \caption{A quasi-isometry between the cone-off of $\Z^2$ with respect to its $\Z-$cosets and $\Z$.}
  \label{fig:chemin}
\end{figure}

As a graph, $\text{CO}(X,\mathcal{C})$ carries a natural metric. Moreover, since now any two vertices in the same piece are adjacent, we have $d_{\text{CO}(X,\mathcal{C})}(x,y) \le d_{X}(x,y)$ for any $x,y\in X$. 

We record some basic facts about cone-offs of quasi-isometric graphs. These are probably well known by the experts, but we include the proofs for completeness.

\begin{lemma}\label{lem:inducedQIbetweenConeOffs}
Let $X,Y$ be two graphs, $f\colon X\rightarrow Y$ be a quasi-isometry, and $\mathcal{C}$ be a partition of $X$. Then $f$ induces a quasi-isometry $\hat{f}\colon \text{CO}(X,\mathcal{C})\longrightarrow \text{CO}(Y, f(\mathcal{C}))$, where $f(\mathcal{C})\defeq\lbrace f(C) : C\in \mathcal{C}\rbrace$. 
\end{lemma}

\begin{proof}
Fix constants $A\ge 1$, $B\ge 0$ such that $f$ and one of its quasi-inverse $f'$ are both $(A,B)-$quasi-isometric embeddings with $f\circ f', f'\circ f$ at distance $\le B$ from the identities.

\noindent As $f(X)$ is quasi-dense in $Y$, so is $f(\text{CO}(X,\mathcal{C}))$ in $\text{CO}(Y, f(\mathcal{C}))$. Now let $x,y\in X$, and let $x=x_{0},x_{1},\dots,x_{n}=y$ be the vertices of a geodesic connecting $x$ and $y$ in $\text{CO}(X,\mathcal{C})$. For any $0\le i\le n-1$, there are then two possibilities: either $x_{i}$ and $x_{i+1}$ belong to a common piece $C\in\mathcal{C}$, so $d_{\text{CO}(Y,f(\mathcal{C}))}(f(x_{i}), f(x_{i+1})) \le 1\le A+B$; or $x_{i}$ and $x_{i+1}$ are adjacent in $X$, in which case 
\begin{equation*}
    d_{\text{CO}(Y,f(\mathcal{C}))}(f(x_{i}), f(x_{i+1})) \le d_{Y}(f(x_{i}), f(x_{i+1})) \le A\cdot d_{X}(x_{i},x_{i+1})+B=A+B.
\end{equation*}
Therefore one deduces that
\begin{align*}
    d_{\text{CO}(Y,f(\mathcal{C}))}(f(x),f(y)) &\le \sum_{i=0}^{n-1}d_{\text{CO}(Y,f(\mathcal{C}))}(f(x_{i}), f(x_{i+1})) \\
    &\le (A+B)\cdot n \\
    &= (A+B)\cdot d_{\text{CO}(X,\mathcal{C})}(x,y).
\end{align*}
For the other inequality, let now $f(x)=x_{0},x_{1},\dots,x_{n}=f(y)$ denote the vertices of a geodesic connecting $f(x)$ and $f(y)$ in $\text{CO}(Y,f(\mathcal{C}))$. Here also there are two cases: if for some $0\le i\le n-1$, $x_{i}$ and $x_{i+1}$ are in the same piece $f(C)$ for some $C\in\mathcal{C}$, then 
\begin{equation*}
    d_{\text{CO}(X,\mathcal{C})}(f'(x_{i}), f'(x_{i+1})) \le d_{\text{CO}(X,\mathcal{C})}(f'(x_{i}), \mathcal{C})+d_{\text{CO}(X,\mathcal{C})}(f'(x_{i+1}),\mathcal{C})+1 \le 2B+1
\end{equation*}
where the second inequality follows from 

\begin{claim}\label{claim:claim3.55}
If $a\in f(C)$ for some $C\in\mathcal{C}$, then $d_{X}(f'(a), C) \le B$.
\end{claim}

\begingroup
\renewcommand{\qedsymbol}{$\blacksquare$}
\begin{proof}[Proof of Claim~\ref{claim:claim3.55}]
Write $a=f(b)$ with $b\in \mathcal{C}$, so that
\begin{equation*}
    d_{X}(f'(a), \mathcal{C}) \le d_{X}(f'(a),b)\le d_{X}\big(f'(a), f'(f(b))\big)+d_{X}(f'(f(b)), b)\le B
\end{equation*}
since $f'\circ f$ is at distance $\le B$ from $\text{Id}_{X}$. 
\end{proof}
\endgroup

\begingroup
\renewcommand{\qedsymbol}{\openbox}
\noindent In the case where $x_{i}$ and $x_{i+1}$ are adjacent in $Y$, then 
\begin{equation*}
   d_{\text{CO}(X,\mathcal{C})}(f'(x_{i}), f'(x_{i+1})) \le d_{X}(f'(x_{i}), f'(x_{i+1})) \le A\cdot d_{Y}(x_{i},x_{i+1})+B=A+B.
\end{equation*}
Hence it follows that 
\begin{align*}
    d_{\text{CO}(X,\mathcal{C})}(x,y) &\le 2B+d_{\text{CO}(X,\mathcal{C})}\big(f'(f(y)), f'(f(x))\big) \\
    &\le 2B+\sum_{i=0}^{n-1}d_{\text{CO}(X,\mathcal{C})}(f'(x_{i}), f'(x_{i+1})) \\
    &\le 2B+\max(A+B, 2B+1)\cdot n \\
    &=2B+\max(A+B, 2B+1)\cdot d_{\text{CO}(Y,f(\mathcal{C}))}(f(x),f(y)).
\end{align*}
We conclude that $f$ induces a quasi-isometry $\text{CO}(X,\mathcal{C})\longrightarrow \text{CO}(Y, f(\mathcal{C}))$, as desired.
\end{proof}
\endgroup

Our second observation ensures that the cone-offs of the same graph over two different partitions are quasi-isometric if the two partitions coarsely coincide.

\begin{lemma}\label{lem:Quasi-isometricConeOffs}
Let $X$ be a graph, and $\mathcal{C}_{1}, \mathcal{C}_{2}$ two partitions of $X$. Assume there exists $Q\ge 0$ such that, for any $C_{1}\in \mathcal{C}_{1}$ (resp. $C_{2}\in\mathcal{C}_{2}$), there exists $C_{2}\in\mathcal{C}_{2}$ (resp. $C_{1}\in\mathcal{C}_{1}$) such that the Hausdorff distance between $C_{1}$ and $C_{2}$ is at most $Q$. Then the identity map $X\rightarrow X$ induces a quasi-isometry $\text{CO}(X,\mathcal{C}_{1})\longrightarrow \text{CO}(X,\mathcal{C}_{2})$.
\end{lemma}

\begin{proof}
Let $x,y\in X$, and fix a geodesic $x=x_{0},x_{1},\dots,x_{n}=y$ connecting $x$ and $y$ in $\text{CO}(X,\mathcal{C}_{1})$. For every $0\le i\le n-1$, either $x_{i}$ and $x_{i+1}$ are adjacent in $X$, so $d_{\text{CO}(X,\mathcal{C}_{2})}(x_{i},x_{i+1})=1\le 2Q+1$, or $x_{i}$ and $x_{i+1}$ belong to a common piece of $\mathcal{C}_{1}$, so that $d_{\text{CO}(X,\mathcal{C}_{2})}(x_{i},x_{i+1}) \le 2Q+1$. Thus we get 
\begin{equation*}
    d_{\text{CO}(X,\mathcal{C}_{2})}(x,y) \le \sum_{i=0}^{n-1}d_{\text{CO}(X,\mathcal{C}_{2})}(x_{i},x_{i+1}) \le (2Q+1)\cdot n = (2Q+1)\cdot d_{\text{CO}(X,\mathcal{C}_{1})}(x,y).
\end{equation*}
By symmetry, one also has $\frac{1}{2Q+1}\cdot d_{\text{CO}(X,\mathcal{C}_{1})} \le d_{\text{CO}(X,\mathcal{C}_{2})}$, which concludes.
\end{proof}

The following claim is then a direct consequence of Lemmas~\ref{lem:inducedQIbetweenConeOffs} and~\ref{lem:Quasi-isometricConeOffs}.

\begin{corollary}\label{cor:QIofpairsinduceQIofConeOffs}
Let $X$ (resp. $Y$) be a graph with a partition $\mathcal{C}$ (resp. $\mathcal{D}$), and let $f\colon (X,\mathcal{C})\longrightarrow (Y,\mathcal{D})$ be an $(A,B,Q)-$quasi-isometry of pairs. Then $f$ induces an $(A',B')-$quasi-isometry $f^{\text{in}}\colon \text{CO}(X,\mathcal{C})\longrightarrow \text{CO}(Y,\mathcal{D})$, where 
\begin{equation*}
    A' \defeq (2Q+1)\cdot \max(A+B,2B+1), \;B'\defeq (2Q+1)\cdot \frac{2B}{\max(A+B,2B+1)}.
\end{equation*}
\end{corollary}

Observe that, given a graph $X$ and a partition $\mathcal{C}$, leaves form a partition of $\mathcal{L}_{n}(X,\mathcal{C})$. Moreover, each leaf is itself a copy of $X$ and thus carries its own partition $\mathcal{C}$. Hence the collection of all pieces constituting all the leaves form a partition $\widehat{\mathcal{C}}$ of $\mathcal{L}_{n}(X,\mathcal{C})$, and one can consider the cone-off $\text{CO}\big(\mathcal{L}_{n}(X,\mathcal{C}), \widehat{\mathcal{C}}\big)$. A crucial observation about this cone-off is the following:

\begin{lemma}\label{lem:isobetweenConeOffandLamplighter}
Let $n\ge 2$. Let $X$ be a graph with a partition $\mathcal{C}$, and let $\widehat{\mathcal{C}}$ be the induced partition of $\mathcal{L}_{n}(X,\mathcal{C})$. Then there is a graph isomorphism 
\begin{equation*}
    \text{CO}\big(\mathcal{L}_{n}(X,\mathcal{C}), \widehat{\mathcal{C}}\big) \cong \mathcal{L}_{n}(\text{CO}(X,\mathcal{C}), \mathcal{C}).
\end{equation*}
\end{lemma}

\begin{proof}
Define the map 
\begin{align*}
    \varphi\colon \text{CO}\big(\mathcal{L}_{n}(X,\mathcal{C}), \widehat{\mathcal{C}}\big) &\longrightarrow \mathcal{L}_{n}(\text{CO}(X,\mathcal{C}), \mathcal{C}) \\
    (c,p)&\longmapsto (c,p)
\end{align*}
where, in the target space, $c$ is seen as a colouring of the partition $\mathcal{C}$ of the graph $\text{CO}(X,\mathcal{C})$. Clearly, $\varphi$ sends bijectively vertices of $\text{CO}\big(\mathcal{L}_{n}(X,\mathcal{C}), \widehat{\mathcal{C}}\big)$ to vertices of $\mathcal{L}_{n}(\text{CO}(X,\mathcal{C}), \mathcal{C})$, so we only need to check that it sends adjacent vertices to adjacent vertices. Thus, let $a=(c_{1},p_{1}), b=(c_{2},p_{2})$ be adjacent vertices in $\text{CO}\big(\mathcal{L}_{n}(X,\mathcal{C}), \widehat{\mathcal{C}}\big)$. There are three cases to consider.
\begin{itemize}
    \item Assume first that $c_{1}=c_{2}$ and that $p_{1}, p_{2}$ are in the same zone. Then $\varphi(a)=(c_{1},p_{1})$, $\varphi(b)=(c_{1},p_{2})$ are in the same leaf and they are adjacent in $\mathcal{L}_{n}(\text{CO}(X,\mathcal{C}), \mathcal{C})$ since $p_{1},p_{2}$ are adjacent in $\text{CO}(X,\mathcal{C})$;
    \item Assume next that $c_{1}=c_{2}$ and that $p_{1},p_{2}$ are in different zones. Then $d_{X}(p_{1},p_{2})=1$, so $\varphi(a),\varphi(b)$ are adjacent in $\mathcal{L}_{n}(\text{CO}(X,\mathcal{C}), \mathcal{C})$;
    \item Lastly, suppose that $c_{1}\neq c_{2}$ and $p_{1}=p_{2}$. Then $c_{1}$ and $c_{2}$ differ only on the zone containing $p_{1}$, and this remains true in $\mathcal{L}_{n}(\text{CO}(X,\mathcal{C}), \mathcal{C})$, so that $\varphi(a)$ and $\varphi(b)$ are adjacent in this case as well.
\end{itemize}

\noindent Conversely, consider the map
\begin{align*}
    \psi\colon \mathcal{L}_{n}(\text{CO}(X,\mathcal{C}), \mathcal{C}) &\longrightarrow \text{CO}\big(\mathcal{L}_{n}(X,\mathcal{C}), \widehat{\mathcal{C}}\big) \\
    (c,p)&\longmapsto (c,p)
\end{align*}
where, in the target space, $c$ is seen as a finitely supported colouring of the vertices of $X$, constant on each piece of $\mathcal{C}$. Let $a=(c_{1},p_{1})$, $b=(c_{2},p_{2})$ be adjacent vertices in $\mathcal{L}_{n}(\text{CO}(X,\mathcal{C}), \mathcal{C})$. 
\begin{itemize}
\item Assume first that $c_{1}=c_{2}$ and $p_{1}$, $p_{2}$ are in the same zone. Then $a=(c,p_{1})$ and $b=(c,p_{2})$ are in the same piece of $\widehat{\mathcal{C}}$, so $\psi(a)$ and $\psi(b)$ are adjacent in $\text{CO}\big(\mathcal{L}_{n}(X,\mathcal{C}), \widehat{\mathcal{C}}\big)$;
\item Assume next that $c_{1}=c_{2}$ but $p_{1}$ and $p_{2}$ are in different zones. Then $d_{X}(p_{1},p_{2})=1$, so $a$ and $b$ are adjacent in $\mathcal{L}_{n}(X,\mathcal{C})$, and thus also in $\text{CO}\big(\mathcal{L}_{n}(X,\mathcal{C}), \widehat{\mathcal{C}}\big)$;
\item Lastly, suppose that $c_{1}\neq c_{2}$ and $p_{1}=p_{2}$. Then $c_{1}$ and $c_{2}$ differ only on the zone containing $p_{1}$, and this stays true when coning-off $\mathcal{L}_{n}(X,\mathcal{C})$ with respect to $\widehat{\mathcal{C}}$. Thus $\psi(a)$ and $\psi(b)$ are adjacent as well. 
\end{itemize}
We conclude that $\varphi$ is a graph isomorphism, with $\psi$ as its inverse.
\end{proof}

\subsection{Proof of Theorem \ref{thm:rigiditypart}}\label{subsubsection3.4.5} We are finally ready to exploit our work in the above parts to deduce the main theorem of this section. 

\begin{proof}[Proof of Theorem~\ref{thm:rigiditypart}]

Let $T$ (resp. $S$) be a finite generating set of $G$ (resp. $H$), and denote $\pi_{G}$ (resp. $\pi_{H}$) the natural projection of $G$ (resp. $H$) onto $G/M$ (resp. $H/N$).

\noindent Fix $q\colon E\wr_{G/M}G\longrightarrow F\wr_{H/N}H$ such an $(A,B)-$quasi-isometry. By Corollary~\ref{cor:QIareleafpreserving}, we know that, up to finite distance, $q$ can be taken to be leaf-preserving, so we can write 
\begin{equation*}
    q(c,p)=(\alpha(c),\beta_{c}(p)), \: (c,p)\in E\wr_{G/M}G
\end{equation*}
for some bijection $\alpha\colon E^{(G,\mathcal{C}_{M})}\longrightarrow F^{(H,\mathcal{C}_{N})}$ and some quasi-isometries $\beta_{c}\colon G\rightarrow H$, $c\in E^{(G,\mathcal{C}_{M})}$. Using Proposition~\ref{prop:leafpreservingQIareQIofpairs}, we also fix $Q\ge 0$ such that 
\begin{equation*}
   q\colon (E\wr_{G/M}G,M)\longrightarrow (F\wr_{H/N}H, N)
\end{equation*}
is an $(A,B,Q)-$quasi-isometry of pairs.

\noindent Now, the image under $q$ of the partition $\widehat{\mathcal{C}_{M}}$ of $E\wr_{G/M}G$ coarsely coincide with $\widehat{\mathcal{C}_{N}}$, according to Corollary~\ref{cor:QIareQIofpairs}. Thus, from Corollary~\ref{cor:QIofpairsinduceQIofConeOffs}, we deduce that $q$ induces a quasi-isometry 
\begin{equation*}
    q^{\text{in}}\colon \text{CO}\big(E\wr_{G/M}G, \widehat{\mathcal{C}_{M}}\big) \longrightarrow \text{CO}\big(F\wr_{H/N}H, \widehat{\mathcal{C}_{N}}\big)
\end{equation*}
or equivalently, by Lemma~\ref{lem:isobetweenConeOffandLamplighter}, a quasi-isometry
\begin{equation*}
    q^{\text{in}}\colon \mathcal{L}_{|E|}\big(\text{CO}(\text{Cay}(G,T),\mathcal{C}_{M}), \mathcal{C}_{M}\big) \longrightarrow \mathcal{L}_{|F|}\big(\text{CO}(\text{Cay}(H,S),\mathcal{C}_{N}), \mathcal{C}_{N}\big).
\end{equation*}
Now, pieces of $\text{CO}(\text{Cay}(G,T),\mathcal{C}_{M})$ (resp. $\text{CO}(\text{Cay}(H,S),\mathcal{C}_{N})$) are uniformly bounded, so we may apply Proposition~\ref{prop:PLovergraphswithboundedpiecesreducetoCL} to deduce that there are graphs $Y_{M}$ and $Y_{N}$, quasi-isometric to $\text{CO}(\text{Cay}(G,T),\mathcal{C}_{M})$ and $\text{CO}(\text{Cay}(H,S),\mathcal{C}_{N})$ respectively, and a quasi-isometry
\begin{equation*}
    \mathcal{L}_{|E|}(Y_{M}) \longrightarrow \mathcal{L}_{|F|}(Y_{N}).
\end{equation*}
Additionally, it follows from the proof of Proposition~\ref{prop:PLovergraphswithboundedpiecesreducetoCL} that $Y_{M}$ (resp. $Y_{N}$) is isomorphic to $\text{Cay}(G/M,\pi_{G}(T))$ (resp. $\text{Cay}(H/N,\pi_{H}(S))$), so that there is a quasi-isometry
\begin{equation*}
    \varphi\colon\mathcal{L}_{|E|}\big(\text{Cay}(G/M,\pi_{G}(T))\big) \longrightarrow \mathcal{L}_{|F|}\big(\text{Cay}(H/N,\pi_{H}(S))\big)
\end{equation*}
given by $(c,pM)\longmapsto \big(\alpha(\Tilde{c})^{\flat}, \overline{\beta_{\Tilde{c}}}(pM)\big)$ (here the notations are the ones introduced in the proof of Proposition~\ref{prop:PLovergraphswithboundedpiecesreducetoCL}). In such a situation, the proof of~\cite[Corollary~6.11]{GT24a} shows that all $\overline{\beta_{c}}\colon G/M\rightarrow H/N$ lie at a uniform bounded distance from $\overline{\beta_{\mathbf{1}}}$: indeed, the source and target spaces of $\varphi$ are in fact quasi-isometric to Cayley graphs of $\Z_{|E|}\wr G/M$ and $\Z_{|F|}\wr H/N$ respectively. These groups are finitely generated halo groups in the sense of~\cite{GT24a}, and they have finite altitude~\cite[Lemma~6.12]{GT24a}. Hence $\varphi$ is at bounded distance from the aptolic quasi-isometry
\begin{align*}
    \mathcal{L}_{|E|}\big(\text{Cay}(G/M,\pi_{G}(T))\big) &\longrightarrow \mathcal{L}_{|F|}\big(\text{Cay}(H/N,\pi_{H}(S))\big) \\
    (c,pM) &\longmapsto (\alpha(\Tilde{c})^{\flat},\overline{\beta_{\mathbf{1}}}(pM)).
\end{align*}

\begin{figure}[H]
  \centering
  \includegraphics[width=0.7\linewidth]{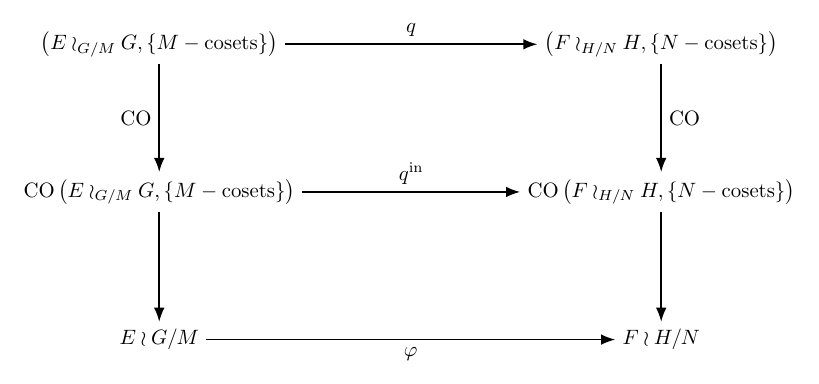}
  \caption{The situation in the proof of Theorem~\ref{thm:rigiditypart}}
  \label{fig:chemin}
\end{figure}

\noindent We deduce from~\cite[Proposition~3.4]{GT24b} that $|E|$ and $|F|$ must have the same prime divisors, which proves~\textit{(i)} and, moreover, if $M$ is co-amenable in $G$, we conclude from~\cite[Theorem~3.9]{GT24b} that $|E|=n^{r}$, $|F|=n^{s}$ for some $n,r,s\ge 1$, and that $\overline{\beta_{\mathbf{1}}}\colon G/M\rightarrow H/N$ is quasi-$\frac{s}{r}$-to-one. 

\noindent To conclude, we have proved in Proposition~\ref{prop:leafpreservingQIareQIofpairs} that $\beta_{\mathbf{1}}$ is a quasi-isometry of pairs, and the existence of such a quasi-isometry implies that $G/M$ is quasi-isometric to $H/N$ by Proposition~\ref{prop:inducedmapbetweenquotients}. Lastly, the restriction of a quasi-isometry of pairs $(G,M)\longrightarrow (H,N)$ to $M$ provides a quasi-isometry $g\colon (M, d_{T})\longrightarrow (N, d_{S})$. If we fix a finite generating set $R$ (resp. $U$) for $M$ (resp. for $N$), we know that $\text{Id}_{M}\colon (M, d_{R})\longrightarrow (M,d_{T})$ and $\text{Id}_{N}\colon (N, d_{U})\longrightarrow (N, d_{S})$ are coarse equivalences, and thus 
\begin{equation*}
    \text{Id}_{N}\circ g\circ \text{Id}_{M}\colon (M,d_{R})\longrightarrow (N,d_{U})
\end{equation*}
is a coarse equivalence. Such a coarse equivalence is in fact a quasi-isometry (Proposition~\ref{prop:boostingcoarsemapstolargescalemaps}(\textit{ii})), concluding the proof. 
\end{proof}

\section{Constructing quasi-isometries of permutational wreath products}\label{subsection3.5}

In this part, we focus on the converse of Theorem~\ref{thm:rigiditypart}, and show how to construct quasi-isometries between permutational wreath products, given various assumptions on the cardinalities of lamp groups or on the geometry of the base groups. 

\subsection{Non-amenable quotients}\label{subsubsection3.5.1} We start with the non-amenable case, by showing that, if the subgroup under consideration is not co-amenable into the ambient group, then two permutational wreath products are quasi-isometric as soon as the cardinalities of the lamp groups have the same prime divisors, namely: 

\begin{proposition}\label{prop:constructingQIbetweenNAPL}
Let $n,m\ge 2$. Let $H$ be a finitely generated group, and let $N\lhd H$ be a normal subgroup of $H$. If $N$ is not co-amenable in $H$ and if $n$ and $m$ have the same prime divisors, then there exists an aptolic quasi-isometry
\begin{equation*}
    \mathcal{L}_{n}(H, \mathcal{C}_{N}) \longrightarrow \mathcal{L}_{m}(H, \mathcal{C}_{N}).
\end{equation*}
\end{proposition}

Recall that $\mathcal{L}_{n}(H, \mathcal{C}_{N})$ stands for $\mathcal{L}_{n}(\text{Cay}(H,S), \mathcal{C}_{N})$ where $S$ is a finite generating set of $H$.

We emphasize that, in this statement as well as in the other ones in this section, the subgroup $N$ need not be finitely generated.

\begin{proof}
Fix a finite generating set $S$ of $H$. Denote $\pi_{H}\colon H\rightarrow H/N$ the canonical projection.

\noindent As a first reduction, we fix integers $m\ge 1$, $n\ge 2$ and a prime number $p$, and we show that there exists an aptolic quasi-isometry
\begin{equation*}
    \mathcal{L}_{mp}(H, \mathcal{C}_{N}) \longrightarrow \mathcal{L}_{mp^{n}}(H, \mathcal{C}_{N}).
\end{equation*}
This is sufficient to deduce the proposition. 

\noindent First of all, since $H/N$ is not amenable, there exists an $n$-to-one map $\overline{f}\colon H/N\rightarrow H/N$ lying at bounded distance from the identity $\text{Id}_{H/N}$. Indeed, as $H/N$ is not amenable, the embedding $\iota \colon H/N \hookrightarrow H/N\oplus \Z_{n}$ lies at finite distance, say $Q\ge 0$, from a bijection $g\colon H/N\longrightarrow H/N\oplus\Z_{n}$. If $p\colon H/N\oplus\Z_{n}\longrightarrow H/N$ denotes the natural projection, then $\overline{f}\defeq p\circ g$ is $n$-to-one and 
\begin{equation*}
    d_{H/N}(x,\overline{f}(x))=d_{H/N}\big(p(\iota(x)),p(g(x))\big) \le d_{H}(\iota(x),g(x))\le Q
\end{equation*}
for any $x\in H/N$, so $\overline{f}$ is at distance $\le Q$ from $\text{Id}_{H/N}$. 

\noindent Fix an enumeration of $H/N$ and identify $\Z_{mp}$ (resp. $\Z_{mp^{n}}$) with $\Z_{m}\oplus\Z_{p}$ (resp. $\Z_{m}\oplus\Z_{p}^{n}$). We denote by $\pi_{1}$ and $\pi_{2}$ the projections on the first and second coordinates in both $\Z_{m}\oplus\Z_{p}$ and $\Z_{m}\oplus\Z_{p}^{n}$.

\noindent Given a colouring $c\colon H/N \longrightarrow \Z_{m}\oplus\Z_{p}$, we define a new colouring $\Tilde{c}\colon H/N\longrightarrow \Z_{m}\oplus\Z_{p}^{n}$ as follows. Given a point $x\in H/N$, set 
\begin{itemize}
    \item $\pi_{1}(\Tilde{c}(x)) \defeq \pi_{1}(c(x))$;
    \item enumerate $\overline{f}^{-1}(\lbrace x\rbrace)$ as $\lbrace x_{1},\dots,x_{n}\rbrace$ following the order of our enumeration, and set
    \begin{equation*}
        \pi_{2}(\Tilde{c}(x))\defeq \big(\pi_{2}(c(x_{1})),\dots,\pi_{2}(c(x_{n}))\big)
    \end{equation*}
\end{itemize}
We can then define 
\begin{align*}
    \varphi\colon \mathcal{L}_{mp}(H, \mathcal{C}_{N}) &\longrightarrow \mathcal{L}_{mp^{n}}(H, \mathcal{C}_{N}) \\
    (c,p)&\longmapsto (\Tilde{c},p)
\end{align*}
and we claim that $\varphi$ is a quasi-isometry. To prove this claim, fix two vertices $a=(c,p), b=(d,q)$ of $\mathcal{L}_{mp}(H,\mathcal{C}_{N})$.

\noindent Notice that if $\Tilde{c}(pN)\neq\Tilde{d}(pN)$ for some $pN\in H/N$, then either $\pi_{1}(\Tilde{c}(pN))\neq \pi_{1}(\Tilde{d}(pN))$, i.e. $\pi_{1}(c(pN))\neq \pi_{1}(d(pN))$, hence $pN\in \text{supp}(c^{-1}d)$; or $\pi_{2}(\Tilde{c}(pN))\neq \pi_{2}(\Tilde{d}(pN))$, i.e. there exists some $p'N\in \overline{f}^{-1}(pN)$ such that $\pi_{2}(c(p'N)) \neq \pi_{2}(d(p'N))$, whence $p'N\in \text{supp}(c^{-1}d)$. We deduce that
\begin{equation}\label{eq1}
    \text{supp}\big(\Tilde{c}^{-1}\Tilde{d}\big) \subset \text{supp}(c^{-1}d)\cup \overline{f}\big(\text{supp}(c^{-1}d)\big)
\end{equation}

\noindent Next, if $c(pN)\neq d(pN)$, then either $\pi_{1}(c(pN))\neq \pi_{1}(d(pN))$ hence $\pi_{1}(\Tilde{c}(pN))\neq\pi_{1}(\Tilde{d}(pN))$; or $\pi_{2}(c(pN)) \neq\pi_{2}(d(pN))$, in which case $\Tilde{c}\big(\overline{f}(pN)\big)\neq \Tilde{d}\big(\overline{f}(pN)\big)$. Thus it follows that
\begin{equation}\label{eq2}
    \text{supp}(c^{-1}d)\subset \text{supp}\big(\Tilde{c}^{-1}\Tilde{d}\big)\cup \overline{f}^{-1}\big(\text{supp}\big(\Tilde{c}^{-1}\Tilde{d}\big)\big).
\end{equation}

\noindent Now, fix a path $\alpha$ in $\text{Cay}(H,S)$ starting from $p$, visiting all zones in $\text{supp}(c^{-1}d)$ and ending at $q$, such that 
\begin{equation*}
    \text{length}(\alpha)+|\text{supp}(c^{-1}d)|=d(a,b).
\end{equation*}
We write explicitly $\text{supp}(c^{-1}d)=\lbrace p_{1}N,\dots, p_{r}N\rbrace$. For any $1\le i\le r$, denote by $x_{i}\in p_{i}N$ the point where the color of the zone is modified. Since 
\begin{equation*}
    d_{H/N}\big(p_{i}N, \overline{f}(p_{i}N)\big)\le Q
\end{equation*}
there is a loop $\overline{\gamma_{i}}$ in $\text{Cay}(H/N,\pi_{H}(S))$ based at $p_{i}N$ and passing through $\overline{f}(p_{i}N)$, of length $\le 2Q$. Thus there is in $\text{Cay}(H,S)$ a loop $\gamma_{i}$ based at $x_{i}$ and passing through the zone $\overline{f}(p_{i}N)$, of length $\le 2Q$. We concatenate this loop to $\alpha$. Doing this for any $1\le i\le r$, we get a path $\alpha'$ in $\text{Cay}(H,S)$ starting at $p$, visiting all zones in $\text{supp}(c^{-1}d)\cup \overline{f}(\text{supp}(c^{-1}d))$ and ending at $q$. A fortiori, $\alpha'$ visits all zones in $\text{supp}\big(\Tilde{c}^{-1}\Tilde{d}\big)$, according to (\ref{eq1}). Observe also that 
\begin{equation*}
    \text{length}(\alpha')\le \text{length}(\alpha)+2Q\cdot |\text{supp}(c^{-1}d)| \le (2Q+1)\cdot \text{length}(\alpha)
\end{equation*}
and thus it follows that
\begin{align*}
    d(\varphi(a),\varphi(b))&=d\big((\Tilde{c},p), (\Tilde{d},q)\big) \le \text{length}(\alpha')+\big|\text{supp}\big(\Tilde{c}^{-1}\Tilde{d}\big)\big| \\
    &\le 2\cdot\text{length}(\alpha') \\
    &\le 2(2Q+1)\cdot \text{length}(\alpha) \\
    &\le 2(2Q+1)\cdot d(a,b).
\end{align*}
Conversely, fix a path $\beta$ in $\text{Cay}(H,S)$ starting from $p$, visiting all zones in $\text{supp}\big(\Tilde{c}^{-1}\Tilde{d}\big)$, ending at $q$, such that 
\begin{equation*}
    \text{length}(\beta)+\big|\text{supp}\big(\Tilde{c}^{-1}\Tilde{d}\big)\big|=d(\varphi(a),\varphi(b)).
\end{equation*}
Here also, we write $\text{supp}\big(\Tilde{c}^{-1}\Tilde{d}\big)=\lbrace q_{1}N,\dots, q_{s}N\rbrace$ and for any $1\le i\le s$, we let $y_{i}\in q_{i}N$ be the point where the color of the zone is changed. For any $1\le i\le s$, and for any $q'N\in \overline{f}^{-1}(\lbrace q_{i}N\rbrace)$, there is a loop $\overline{\eta_{i}}$ in $\text{Cay}(H/N,\pi_{H}(S))$ based at $q_{i}N$ and passing through $q'N$, of length $\le 6Q$, since 
\begin{equation*}
    d_{H/N}(q_{i}N, q'N)\le d_{H/N}\big(q_{i}N, \overline{f}(q_{i}N)\big)+d_{H/N}\big(\overline{f}(q_{i}N), \overline{f}(q'N)\big)+d_{H/N}\big(\overline{f}(q'N), q'N\big) \le 3Q.
\end{equation*}
Thus there is a loop $\eta_{i}$ in $\text{Cay}(H,S)$ based at $y_{i}$ and passing through the zone $q'N$, of length $\le 6Q$. We concatenate this loop to $\beta$. Doing this for any $1\le i\le s$ and for any pre-image of $q_{i}N$ under $\overline{f}$, we get a path $\beta'$ in $\text{Cay}(H,S)$ starting from $p$, visiting all zones in $\text{supp}\big(\Tilde{c}^{-1}\Tilde{d}\big)\cup \overline{f}^{-1}\big(\text{supp}\big(\Tilde{c}^{-1}\Tilde{d}\big)\big)$ and ending at $q$. A fortiori, $\beta'$ visits all zones in $\text{supp}(c^{-1}d)$ according to (\ref{eq2}). Moreover, one has 
\begin{equation*}
    \text{length}(\beta')\le \text{length}(\beta)+6nQ\cdot \big|\text{supp}\big(\Tilde{c}^{-1}\Tilde{d}\big)\big| \le (6nQ+1)\cdot \text{length}(\beta).
\end{equation*}
Hence we deduce that
\begin{align*}
    d(a,b) &=d((c,p),(d,q)) \le \text{length}(\beta')+|\text{supp}(c^{-1}d)| \\
    &\le 2\cdot\text{length}(\beta') \\
    &\le 2(6nQ+1)\cdot \text{length}(\beta) \\
    &\le 2(6nQ+1)\cdot d(\varphi(a),\varphi(b)).
\end{align*}
This concludes the proof that $\varphi$ is a quasi-isometry. 
\end{proof}

\begin{remark}\label{rm5.2}
More generally, and as already noticed in~\cite{BGT24} for the standard case, the same strategy proves that if $A,B$ and $H$ are finitely generated groups and that $N\lhd H$ is not co-amenable in $H$, then there is a quasi-isometry 
\begin{equation*}
    (A\times B)\wr_{H/N}H \longrightarrow (A\times B^{n})\wr_{H/N}H
\end{equation*}
for any $n\ge 1$. 
\end{remark}

Next, we show that a quasi-isometry of pairs inducing a quasi-one-to-one quasi-isometry between the quotients allows us to construct aptolic quasi-isometries between permutational wreath products. 

\begin{proposition}\label{prop:aptolicQIfromscalingQIbetweenquotients}
Let $n\ge 2$. Let $G, H$ be finitely generated groups with normal subgroups $M\lhd G$, $N\lhd H$. Let $f\colon (G,M)\longrightarrow (H,N)$ be a quasi-isometry of pairs such that $\overline{f}\colon G/M\rightarrow H/N$ is quasi-one-to-one. Then there exists an aptolic quasi-isometry
\begin{equation*}
    \mathcal{L}_{n}(G,\mathcal{C}_{M})\longrightarrow \mathcal{L}_{n}(H,\mathcal{C}_{N}).
\end{equation*}
\end{proposition}

Note that, according to Remark~\ref{rm2.12}, the assumption made on $\overline{f}$ does not depend on a specific choice of the induced quasi-isometry. 

\begin{proof} Since $\overline{f}$ is quasi-one-to-one, we know from Theorem~\ref{thm:Whytethm} that it lies at finite distance, say $Q\ge 0$, from a bijection $g\colon G/M\rightarrow H/N$. Given now $p\in G$, define $h(p)\in H$ as being a point of the coset $g(pM)$ that minimizes the distance to $f(p)$. Such a point exists, and as
\begin{equation*}
    d_{H/N}\big(g(pM), \overline{f}(pM)\big) \le Q
\end{equation*}
it is at distance at most $Q$ in $H$ from $f(p)$. The map $h\colon G\rightarrow H$ thus defined lies at distance $\le Q$ from $f$, so it is itself a quasi-isometry, and by construction it sends an $M-$coset into an $N-$coset. Moreover, still by construction, one can take $\overline{h}=g$. We thus have a quasi-isometry of pairs $h\colon (G,M)\longrightarrow (H,N)$ inducing a bijection at the level of quotients. We now claim that 
\begin{align*}
    \varphi\colon \mathcal{L}_{n}(G,\mathcal{C}_{M})&\longrightarrow \mathcal{L}_{n}(H,\mathcal{C}_{N}) \\
    (c,p)&\longmapsto (c\circ \overline{h}^{-1}, h(p))
\end{align*}
is the quasi-isometry we are looking for. 

\noindent Denote $T$ (resp. $S$) a finite generating set for $G$ (resp. $H$). Let $C\ge 1$, $K\ge 0$ be such that $h\colon G\rightarrow H$ and a quasi-inverse $h'\colon H\rightarrow G$ are $(C,K)-$quasi-isometries, with $h\circ h'$ and $h'\circ h$ at distance $\le K$ from the identities. Moreover, in this situation we may assume that $\overline{h'}=\overline{h}^{-1}$.

\noindent Notice first that $\varphi$ is $(C+K)-$Lipschitz. Indeed, if $a,b\in \mathcal{L}_{n}(G,\mathcal{C}_{M})$ are adjacent vertices, then 
\begin{itemize}
    \item either $a=(c,p)$ and $b=(c,ps)$ for some $s\in T$, so 
    \begin{equation*}
        d(\varphi(a),\varphi(b))=d_{H}(h(p),h(ps))\le C\cdot d_{G}(p,ps)+K=C+K;
    \end{equation*}
    \item or $a=(c,p)$ and $b=(c',p)$, where $c,c'$ only differ on $pM$. In this case, we get 
    \begin{equation*}
        c\circ \overline{h}^{-1}(h(p)N)=c\circ \overline{h}^{-1}\big(\overline{h}(pM)\big)=c(pM)
    \end{equation*}
    and likewise $c'\circ\overline{h}^{-1}(h(p)N)=c'(pM)$, so $c\circ \overline{h}^{-1}, c'\circ \overline{h}^{-1}$ differ on the coset containing $h(p)$; and if $yN\neq h(p)N$ is another coset on which these two colourings differ, then we may write $yN=\overline{h}(zM)$ for some $zM\in G/M$, and we get
    \begin{equation*}
         c\circ \overline{h}^{-1}(yN)=c(zM), \; c'\circ \overline{h}^{-1}(yN)=c'(zM)
    \end{equation*}
    so $c,c'$ differ also on $zM\neq pM$, a contradiction. Thus $c\circ \overline{h}^{-1}, c'\circ \overline{h}^{-1}$ differ only on $h(p)N$, and we conclude that $\varphi(a), \varphi(b)$ are adjacent in $\mathcal{L}_{n}(H,\mathcal{C}_{N})$.
\end{itemize}
In both cases, we conclude that $d(\varphi(a),\varphi(b))\le C+K$, and Lemma~\ref{lem:Lipschitzbetweengraphs} now ensures that $\varphi$ is $(C+K)-$Lipschitz. 

\noindent Next, consider the map 
\begin{align*}
    \psi\colon \mathcal{L}_{n}(H,\mathcal{C}_{N}) &\longrightarrow \mathcal{L}_{n}(G,\mathcal{C}_{M}) \\
    (c,p)&\longrightarrow (c\circ \overline{h}, h'(p))
\end{align*}
and note that 
\begin{equation*}
    d\big(\psi\circ\varphi(c,p), (c,p)\big)=d_{G}\big(h'(h(p)), p\big)\le K
\end{equation*}
for any $(c,p)\in \mathcal{L}_{n}(G,\mathcal{C}_{M})$, so $\psi\circ\varphi$ is at distance $\le K$ from $\text{Id}_{\mathcal{L}_{n}(G,\mathcal{C}_{M})}$; and similarly $\varphi\circ\psi$ is at distance $\le K$ from $\text{Id}_{\mathcal{L}_{n}(H,\mathcal{C}_{N})}$. 

\noindent Finally, let us prove that $\psi$ is also $(C+K)-$Lipschitz. Fix $a$ and $b$ two adjacent vertices in $\mathcal{L}_{n}(H,\mathcal{C}_{N})$, and consider two cases:
\begin{itemize}
    \item if $a=(c,p)$ and $b=(c,ps)$ for some $s\in S$, then 
    \begin{equation*}
        d(\psi(a), \psi(b))=d_{G}\big(h'(p),h'(ps)\big) \le C\cdot d_{H}(p,ps)+K=C+K;
    \end{equation*}
    \item and if $a=(c,p)$, $b=(c',p)$ where $c,c'$ differ only on $pN$, then $c\circ\overline{h}$, $c'\circ \overline{h}$ differ on $h'(p)N$, since 
    \begin{equation*}
        c\circ\overline{h}\big(h'(p)M\big)=c\circ\overline{h}\big(\overline{h'}(pN)\big)=c(pN), \; c'\circ\overline{h}\big(h'(p)M\big)=c'\circ\overline{h}\big(\overline{h'}(pN)\big)=c'(pN)
    \end{equation*}
    and since $c,c'$ differ on $pN$ by assumption (here we use also that $\overline{h'}=\overline{h}^{-1}$); and if by contradiction $c\circ\overline{h}$, $c'\circ\overline{h}$ differ on another coset $yM\neq h'(p)M$, then we write this coset as $yM=\overline{h}^{-1}(zN)=\overline{h'}(zN)$ to get that $c,c'$ differ on the coset $zN \neq pN$, a contradiction. Thus $c\circ\overline{h}$ and $c'\circ\overline{h}$ differ only on the coset $h'(p)M$, so $\psi(a)$ and $\psi(b)$ are adjacent in $\mathcal{L}_{n}(G,\mathcal{C}_{M})$.
\end{itemize}
Hence we conclude that $d(\psi(a),\psi(b))\le C+K$, so $\psi$ is $(C+K)-$Lipschitz by Lemma~\ref{lem:Lipschitzbetweengraphs}. Thus we conclude that
\begin{align*}
     d(\varphi(a),\varphi(b)) &\ge \frac{1}{C+K}\cdot d\big(\psi(\varphi(a)),\psi(\varphi(b))\big) \\
     &\ge \frac{1}{C+K}\big(d(\psi(\varphi(a)), a)-d(a,b)-d(b,\psi(\varphi(b)))\big) \\
     &\ge \frac{1}{C+K}\cdot d(a,b)-\frac{2K}{C+K}
\end{align*}
for any $a,b\in\mathcal{L}_{n}(G,\mathcal{C}_{M})$, so $\varphi$ is indeed a quasi-isometry, with $\psi$ as a quasi-inverse.
\end{proof}

We can then deduce that:

\begin{corollary}\label{cor:aptolicQIfromaQIofpairsbetweenbases}
Let $n, m\ge 2$. Let $G, H$ be finitely generated groups with normal subgroups $M\lhd G$, $N\lhd H$. Assume that $M$ is not co-amenable in $G$. If there exists a quasi-isometry of pairs $f\colon (G,M)\longrightarrow (H,N)$ and if $n$ and $m$ have the same prime divisors, then there exists an aptolic quasi-isometry
\begin{equation*}
    \mathcal{L}_{n}(G,\mathcal{C}_{M})\longrightarrow \mathcal{L}_{m}(H,\mathcal{C}_{N}).
\end{equation*}
\end{corollary}

\begin{proof}
As the quotient groups $G/M$, $H/N$ are not amenable, the induced quasi-isometry $\overline{f}\colon G/M\rightarrow H/N$ is quasi-one-to-one, and Proposition~\ref{prop:aptolicQIfromscalingQIbetweenquotients} then shows that there is an aptolic quasi-isometry $\mathcal{L}_{n}(G,\mathcal{C}_{M})\longrightarrow \mathcal{L}_{n}(H,\mathcal{C}_{N})$. Additionally, we know that $n$ and $m$ have the same prime divisors, so we deduce from Proposition~\ref{prop:constructingQIbetweenNAPL} the existence of an aptolic quasi-isometry
\begin{equation*}
    \mathcal{L}_{n}(H, \mathcal{C}_{N}) \longrightarrow \mathcal{L}_{m}(H, \mathcal{C}_{N}).
\end{equation*}
Composing these two quasi-isometries yields the desired conclusion.
\end{proof}

\subsection{Amenable quotients}\label{subsubsection3.5.2}

We now shift our attention to the amenable case. We start with a statement giving sufficient conditions to construct aptolic quasi-isometries in a general situation.

\begin{proposition}\label{prop:2mapstogetanAptolicQIbetweenAPL}
Let $n,m\ge 2$. Let $G,H$ be finitely generated groups with normal subgroups $M\lhd G$, $N\lhd H$. Suppose that $\alpha\colon \Z_{n}^{(G/M)} \longrightarrow \Z_{m}^{(H/N)}$, $\beta\colon G\rightarrow H$ are two maps such that:
\begin{enumerate}[label=(\roman*)]
    \item $\alpha$ is a bijection;
    \item $\beta\colon (G,M)\longrightarrow (H,N)$ is a quasi-isometry of pairs;
    \item there exists a choice of $\overline{\beta}$ for which there exists a constant $Q\ge 0$ such that, for any $c_{1},c_{2}\in \Z_{n}^{(G/M)}$, the Hausdorff distance between $\overline{\beta}\big(\text{supp}(c_{1}^{-1}c_{2})\big)$ and $\text{supp}\big(\alpha(c_{1})^{-1}\alpha(c_{2})\big)$ is at most $Q$.  
\end{enumerate}
Then the map 
\begin{align*}
    q\colon \mathcal{L}_{n}(G,\mathcal{C}_{M}) &\longrightarrow \mathcal{L}_{m}(H,\mathcal{C}_{N})  \\
    (c,p)&\longmapsto (\alpha(c),\beta(p))
\end{align*}
is an aptolic quasi-isometry. 
\end{proposition}

Notice that, as in the non-amenable case above (cf. subsection~\ref{subsubsection3.5.1}), we treat colourings here as colourings of the quotient groups. 

\begin{proof}
Let $T$ (resp. $S$) denote a finite generating set of $G$ (resp. $H$), and let $\pi_{G}$ (resp. $\pi_{H}$) be the canonical projection of $G$ (resp. $H$) onto $G/M$ (resp. $H/N$). Up to finite distance, we assume that $\beta$ sends $M-$cosets into $N-$cosets, and that a quasi-inverse $\beta^{\text{qi}}\colon H\rightarrow G$ sends $N-$cosets into $M-$cosets. Fix also constants $C,K\ge 0$ such that $\beta$, $\beta^{\text{qi}}$ are $(C,K)-$quasi-isometries, with $\beta\circ\beta^{\text{qi}}$, $\beta^{\text{qi}}\circ\beta$ within distance $\le K$ from the identities. Up to increasing those constants, we also assume that $\overline{\beta}$, $\overline{\beta^{\text{qi}}}$ are $(C,K)-$quasi-isometries. Lastly, recall from Lemma~\ref{lem:stabilitypropertiesQIofpairs}\textit{(iii)} that a quasi-inverse of $\overline{\beta}$ is precisely given by $\overline{\beta^{\text{qi}}}$, i.e. $\overline{\beta}^{\text{qi}}=\overline{\beta^{\text{qi}}}$.

\noindent To start, we prove that $q$ is Lipschitz. Fix $a,b\in \mathcal{L}_{n}(G,\mathcal{C}_{M})$ two adjacent vertices. We treat two cases:
\begin{itemize}
    \item Assume first that $a=(c,p)$ and $b=(c,ps)$ for some $s\in T$. Then one gets
    \begin{equation*}
        d(q(a),q(b))=d\big((\alpha(c),\beta(p)), (\alpha(c),\beta(ps))\big)=d_{H}(\beta(p),\beta(ps))\le C+K;
    \end{equation*}
    \item Assume now that $a=(c,p)$ and $b=(c',p)$, where $c,c'$ only differ on $pM$. Then $\overline{\beta}\big(\text{supp}(c^{-1}c')\big)$ is reduced to a point, and it follows from assumption~\textit{(iii)} that 
    \begin{equation*}
        \text{supp}\big(\alpha(c)^{-1}\alpha(c')\big) \subset B_{H/N}\big(\overline{\beta}(pM),Q\big)=B_{H/N}(\beta(p)N,Q).
    \end{equation*}
    the ball centered at $\beta(p)N$ of radius $Q$ in $H/N$. We deduce from this inclusion that
    \begin{align*}
        d(q(a),q(b)) &\le 2Q\cdot \left|\text{supp}(\alpha(c)^{-1}\alpha(c'))\right| + \left|\text{supp}(\alpha(c)^{-1}\alpha(c'))\right|\\
        &\le (2Q+1)\cdot D^{Q}
    \end{align*}
    where $D\ge 1$ is an integer larger than the maximal degree of a vertex in $\text{Cay}(H/N,\pi_{H}(S))$. 
\end{itemize}
We conclude from Lemma~\ref{lem:Lipschitzbetweengraphs} that $q$ is $\max(C+K,(2Q+1)\cdot D^{Q})-$Lipschitz. 

\noindent Next, consider the map 
\begin{align*}
    q'\colon \mathcal{L}_{m}(H,\mathcal{C}_{N})  &\longrightarrow \mathcal{L}_{n}(G,\mathcal{C}_{M})  \\
    (c,p)&\longmapsto (\alpha^{-1}(c), \beta^{\text{qi}}(p)).
\end{align*}
Then we have 
\begin{equation*}
    d\big((c,p), q'\circ q(c,p)\big)=d\big((c,p), (c,\beta^{\text{qi}}\circ\beta(p))\big)=d_{H}\big(p,\beta^{\text{qi}}\circ\beta(p)\big) \le K
\end{equation*}
for any $(c,p)\in \mathcal{L}_{n}(G,\mathcal{C}_{M})$, as well as 
\begin{equation*}
    d\big((c,p), q\circ q'(c,p)\big)=d\big((c,p), (c,\beta\circ\beta^{\text{qi}}(p))\big)=d_{G}\big(p,\beta\circ\beta^{\text{qi}}(p)\big) \le K
\end{equation*}
for any $(c,p)\in \mathcal{L}_{m}(H,\mathcal{C}_{N}) $. 

\noindent Notice that the map $q'$ satisfies also points~\textit{(i)-(iii)} of the statement. Points~\textit{(i)},~\textit{(ii)} are clearly satisfied. To prove~\textit{(iii)}, fix two colourings $c_{1},c_{2}\in \Z_{m}^{(H/N)}$. We know by assumption that 
\begin{equation*}
    d_{\text{Haus}}\big(\overline{\beta}\big(\text{supp}(\alpha^{-1}(c_{1})^{-1}\alpha^{-1}(c_{2}))\big), \text{supp}(c_{1}^{-1}c_{2})\big) \le Q
\end{equation*}
which, by Lemma~\ref{lem:neighborhoodsandQI}\textit{(iv)}, implies
\begin{equation*}
d_{\text{Haus}}\big(\overline{\beta}^{\text{qi}}\circ\overline{\beta}\big(\text{supp}(\alpha^{-1}(c_{1})^{-1}\alpha^{-1}(c_{2}))\big), \overline{\beta}^{\text{qi}}\big(\text{supp}(c_{1}^{-1}c_{2})\big)\big) \le C\cdot Q+K.
\end{equation*}
On the other hand, the Hausdorff distance between $\overline{\beta}^{\text{qi}}\circ\overline{\beta}\big(\text{supp}(\alpha^{-1}(c_{1})^{-1}\alpha^{-1}(c_{2}))\big)$ and $\text{supp}\big(\alpha^{-1}(c_{1})^{-1}\alpha^{-1}(c_{2})\big)$ is at most $K$, so we conclude that the Hausdorff distance between $\text{supp}\big(\alpha^{-1}(c_{1})^{-1}\alpha^{-1}(c_{2})\big)$ and $\overline{\beta}^{\text{qi}}\big(\text{supp}(c_{1}^{-1}c_{2})\big)$ is at most $C\cdot Q+2K$. Recalling that $\overline{\beta}^{\text{qi}}=\overline{\beta^{\text{qi}}}$ shows~\textit{(iii)} for $q'$. Henceforth, we can reproduce the argument we did for $q$ to conclude that $q'$ is $\max\big(C+K, (2Q'+1)\cdot D'^{Q'}\big)-$Lipschitz, where $Q'\defeq C\cdot Q+2K$ and where $D'\ge 1$ is an integer larger than the maximal degree of a vertex of $\text{Cay}(G/M, \pi_{G}(T))$. It follows that 
\begin{align*}
    &d\big(q(c_{1},p_{1}),q(c_{2},p_{2})\big) \ge \frac{1}{\max(C+K, (2Q'+1)\cdot D'^{Q'})}\cdot d\big(q'\circ q(c_{1},p_{1}), q'\circ q(c_{2},p_{2})\big) \\
    &\ge \frac{1}{\max(C+K, (2Q'+1)\cdot D'^{Q'})}\cdot d((c_{1},p_{1}),(c_{2},p_{2}))-\frac{2K}{\max(C+K, (2Q'+1)\cdot D'^{Q'})}
\end{align*}
for any $(c_{1},p_{1}), (c_{2},p_{2})\in \mathcal{L}_{n}(G,\mathcal{C}_{M})$; and that 
\begin{align*}
    &d\big(q'(c_{1},p_{1}),q'(c_{2},p_{2})\big) \ge \frac{1}{\max(C+K, (2Q+1)\cdot D^{Q})}\cdot d\big(q\circ q'(c_{1},p_{1}), q\circ q'(c_{2},p_{2})\big) \\
    &\ge \frac{1}{\max(C+K, (2Q+1)\cdot D^{Q})}\cdot d((c_{1},p_{1}),(c_{2},p_{2}))-\frac{2K}{\max(C+K, (2Q+1)\cdot D^{Q})}
\end{align*}
for any $(c_{1},p_{1}), (c_{2},p_{2})\in \mathcal{L}_{m}(H,\mathcal{C}_{N})$. Thus $q$ is an aptolic quasi-isometry, with $q'$ as a quasi-inverse.
\end{proof}

We can therefore deduce the following consequence.

\begin{proposition}\label{prop:aptolicQIfromscalingQIbetweenquotients-Acase}
Let $n,m\ge 2$. Let $G$ and $H$ be finitely generated groups with normal subgroups $M\lhd G$, $N\lhd H$. Suppose $n=k^{r}$, $m=k^{s}$ for some $k,r,s\ge 1$. If there exists a quasi-isometry of pairs $f\colon (G,M)\longrightarrow (H,N)$ such that $\overline{f}\colon G/M\rightarrow H/N$ is quasi-$\frac{s}{r}$-to-one, then there exists an aptolic quasi-isometry
\begin{equation*}
    \mathcal{L}_{n}(G,\mathcal{C}_{M})\longrightarrow \mathcal{L}_{m}(H,\mathcal{C}_{N}).
\end{equation*}
\end{proposition}

We emphasize that, in this statement as well, the subgroups $N$ and $M$ need not be finitely generated.

\begin{proof}
As usual, fix $T$ (resp. $S$) a finite generating set for $G$ (resp. for $H$) and let $\pi_{G}$ (resp. $\pi_{H}$) be the canonical projection onto $G/M$ (resp. $H/N$). Up to finite distance, we may assume that $f(pM)\subset f(p)N$ for any $p\in G$. 

\noindent As $\overline{f}\colon G/M\rightarrow H/N$ is quasi-$\frac{s}{r}$-to-one, we know from Theorem~\ref{thm:rationalscalingfactor} that there exists a partition $\mathcal{P}$ (resp. $\mathcal{Q}$) of $\text{Cay}(G/M, \pi_{G}(T))$ (resp. $\text{Cay}(H/N, \pi_{H}(S))$) with uniformly bounded pieces of size $s$ (resp. of size $r$), a bijection $\psi\colon \mathcal{P}\rightarrow \mathcal{Q}$ and a map $\beta\colon G/M\rightarrow H/N$ at finite distance, say $C\ge 0$, from $\overline{f}$ such that $\beta(P)\subset \psi(P)$ for any $P\in\mathcal{P}$. 

\noindent Now, given $p\in G$, define $h(p)\in H$ as a point of the coset $\beta(pM)\subset H$ minimizing the distance to $f(p)$. Such a point exists, and it is at distance $\le C$ from $f(p)$: indeed, as 
\begin{equation*}
    d_{H/N}\big(\beta(pM), \overline{f}(pM)\big) \le C
\end{equation*}
any point of $\overline{f}(pM)=f(p)N$ is at distance $\le C$ in $H$ from a point in $\beta(pM)\subset H$. In particular, the map $h\colon G\rightarrow H$ thus defined is itself a quasi-isometry of pairs, and by construction it satisfies $\overline{h}=\beta$. 

\noindent Next, fix a bijection $\sigma\colon \Z_{n}^{s}\longrightarrow \Z_{m}^{r}$ such that $\sigma(0)=0$, and define a bijection $\alpha\colon \Z_{n}^{(G/M)}\longrightarrow \Z_{m}^{(H/N)}$ such that $\alpha$ sends $\mathcal{L}_{n}(P)$ (recall that the latter stands for colourings supported exclusively on elements of $P$) to $\mathcal{L}_{m}(\psi(P))$ through $\sigma$, for any $P\in\mathcal{P}$. Namely, given a finitely supported colouring $c\in \Z_{n}^{(G/M)}$, its values on a piece $P\in\mathcal{P}$ form a vector $v$ with $s$ components and with entries in $\Z_{n}$, i.e. an element of $\Z_{n}^{s}$. Then $\sigma(v)$ is a vector with $r$ components and with entries in $\Z_{m}$. These entries are then the values of the colouring $\alpha(c)$ on the piece $\psi(P)$. This way, $\mathcal{L}_{n}(P)$ is sent to $\mathcal{L}_{m}(\psi(P))$ for any $P\in \mathcal{P}$, and the condition $\sigma(0)=0$ ensures that finitely supported colourings are sent to finitely supported colourings. 

\noindent Define then
\begin{align*}
    q\colon \mathcal{L}_{n}(G,\mathcal{C}_{M}) &\longrightarrow \mathcal{L}_{m}(H,\mathcal{C}_{N}) \\
    (c,p)&\longmapsto (\alpha(c), h(p)).
\end{align*}
We show that $q$ is an aptolic quasi-isometry by checking points~\textit{(i)-(iii)} of Proposition~\ref{prop:2mapstogetanAptolicQIbetweenAPL}. Points~\textit{(i)} and~\textit{(ii)} are satisfied by construction. For~\textit{(iii)}, let $c_{1}, c_{2}\in \Z_{n}^{(G,\mathcal{C}_{M})}$, and denote $P_{1},\dots, P_{n}$ the pieces of $\mathcal{P}$ containing points of $\text{supp}(c_{1}^{-1}c_{2})$. Then the pieces
$\psi(P_{1}),\dots, \psi(P_{n})$ are the pieces of $\mathcal{Q}$ containing the points of $\text{supp}\big(\alpha(c_{1})^{-1}\alpha(c_{2})\big)$. Because the pieces of $\mathcal{Q}$ are uniformly bounded, we deduce that there is $D\ge 0$ (independent of $c_{1},c_{2}$) such that
\begin{equation*}
    d_{\text{Haus}}\big(\text{supp}(\alpha(c_{1})^{-1}\alpha(c_{2})), \psi(P_{1})\cup\dots\cup\psi(P_{n})\big) \le D.
\end{equation*}
Also, since $\beta$ sends each piece $P$ of $\mathcal{P}$ into $\psi(P)$, we see that
\begin{equation*}
    \overline{h}\big(\text{supp}(c_{1}^{-1}c_{2})\big)=\beta\big(\text{supp}(c_{1}^{-1}c_{2})\big)
\end{equation*}
lies in $\psi(P_{1})\cup\dots\cup\psi(P_{n})$ and has a point in each of these pieces. Once again, we deduce that 
\begin{equation*}
    d_{\text{Haus}}\big(\beta(\text{supp}(c_{1}^{-1}c_{2})), \psi(P_{1})\cup\dots\cup\psi(P_{n})\big)\le D'
\end{equation*}
for some $D'\ge 0$ independent of $c_{1}$ and $c_{2}$. We conclude that the Hausdorff distance between $\overline{h}\big(\text{supp}(c_{1}^{-1}c_{2})\big)$ and $\text{supp}\big(\alpha(c_{1})^{-1}\alpha(c_{2})\big)$ is bounded above by a constant, independently of $c_{1},c_{2}$. Thus~\textit{(iii)} of Proposition~\ref{prop:2mapstogetanAptolicQIbetweenAPL} holds as well, and we deduce from the latter that $q\colon \mathcal{L}_{n}(G,\mathcal{C}_{M})\longrightarrow \mathcal{L}_{m}(H,\mathcal{C}_{N})$ is an aptolic quasi-isometry, as desired.
\end{proof}

\subsection{Proof of Theorem~\ref{thm:classificationPWPuptoQI}}\label{subsubsection3.5.3} We can now deduce our main theorem from the introduction.

\begin{proof}[Proof of Theorem~\ref{thm:classificationPWPuptoQI}]
Assume that $E\wr_{G/M}G$ and $F\wr_{H/N}H$ are quasi-isometric. By Theorem~\ref{thm:rigiditypart}, we know that $|E|$ and $|F|$ have the same prime divisors and that there is a quasi-isometry of pairs $\beta\colon (G,M)\longrightarrow (H,N)$, which is the desired conclusion in the case where $M$ is not co-amenable in $G$. If $M$ is co-amenable in $G$, we know from Theorem~\ref{thm:rigiditypart} that $|E|=n^{r}$ and $|F|=n^{s}$ are powers of a common number and that $\overline{\beta}\colon G/M\rightarrow H/N$ must be quasi-$\frac{s}{r}$-to-one.

\noindent Conversely, in the case where $M$ is co-amenable, the conclusion follows from Proposition~\ref{prop:aptolicQIfromscalingQIbetweenquotients-Acase}, while if $M$ is not co-amenable, the conclusion follows from Corollary~\ref{cor:aptolicQIfromaQIofpairsbetweenbases}.
\end{proof}

\section{Consequences of Theorem~\ref{thm:classificationPWPuptoQI}}\label{subsection3.6}

\subsection{Proof of Corollaries~\ref{cor:classificationofPWPoverZ^duptoQI} and~\ref{cor:classificationofPWPoverdirectproducts}}\label{subsubsection3.6.1} We start this part with some concrete applications of our criterion.

\begin{proof}[Proof of Corollary~\ref{cor:classificationofPWPoverZ^duptoQI}]
The subgroup $M=\Z^{m-n}$ (resp. $N=\Z^{m'-n'}$) has growth degree $m-n\le m-2$ (resp. $m'-n'\le m'-2$), so that its coarsely embedded subspaces do not coarsely separate $G=\Z^m$ (resp. $H=\Z^{m'}$) by Theorem~\ref{thm:noncoarseseparationinpolynomialgrowthgroups}. 
Additionally, notice that given any $k>0$, there is a quasi-isometry of pairs $(\Z^m, \Z^{m-n})\longrightarrow (\Z^m, \Z^{m-n})$ inducing a quasi-$k$-to-one quasi-isometry $\Z^n\rightarrow \Z^n$. Indeed, since $\text{Sc}(\Z^n)=\R_{>0}$, fix any quasi-$k$-to-one quasi-isometry $f\colon \Z^{n}\rightarrow \Z^n$, and extend it to $\Z^{m}$ by setting 
\begin{align*}
    g\colon \Z^m=\Z^{m-n}\times\Z^n &\longrightarrow   \Z^{m-n}\times\Z^n=\Z^{m}\\
    (x,y)&\longmapsto (x,f(y)).
\end{align*}
Then $g$ is a quasi-isometry of pairs and $\overline{g}=f$. Hence the conclusion follows from Theorem~\ref{thm:classificationPWPuptoQI}.
\end{proof}

As highlighted in the introduction, since $\text{Sc}(\Z^{n})=\R_{>0}$ for any $n\ge 1$, this result is a particular case of Corollary~\ref{cor:classificationofPWPoverdirectproducts}, that we prove now using the same argument.

\begin{proof}[Proof of Corollary~\ref{cor:classificationofPWPoverdirectproducts}]
The left to right direction is the first point of Theorem~\ref{thm:classificationPWPuptoQI}. 

\noindent Conversely, assume that $|E|=n^{r}$, $|F|=n^{s}$ for some integers $n,r,s\ge 1$, and that $\frac{s}{r}\in \text{Sc}(K)$. Fix then $f\colon K\rightarrow K$ a quasi-$\frac{s}{r}$-to-one quasi-isometry, and extend it to $G=M\times K$ via 
\begin{align*}
    g\colon M\times K &\longrightarrow M\times K \\
    (m,k)&\longmapsto (m,f(k)).
\end{align*}
Then $g$ is a quasi-isometry of pairs $(G,M)\longrightarrow (G,M)$, and by construction $\overline{g}=f$. Thus we deduce from Theorem~\ref{thm:classificationPWPuptoQI} that $E\wr_{K}G$ and $F\wr_{K}G$ are quasi-isometric. 
\end{proof}

\subsection{BiLipschitz equivalences of permutational wreath products}\label{subsubsection3.6.2} Let us now shift our attention to Corollary~\ref{cor:classificationPWPoverZ^duptoBILIP}. It will be a straightforward consequence of the next result:

\begin{proposition}\label{prop:BILIPbetweenPLfromscalingQIbetweenquotients}
Let $E$, $F$ be two non-trivial finite groups. Let $G$, $H$ be finitely generated groups, with normal subgroups $M\lhd G$, $N\lhd H$. Suppose that $M$ is co-amenable in $G$.
If there exist $n,r,s\ge 1$ such that $|E|=n^{r}$, $|F|=n^{s}$ and there exists a quasi-one-to-one quasi-isometry of pairs $f\colon (G,M)\longrightarrow (H,N)$ such that $\overline{f}\colon G/M\rightarrow H/N$ is quasi-$\frac{s}{r}$-to-one, then $E\wr_{G/M}G$ and $F\wr_{H/N}H$ are biLipschitz equivalent.
\end{proposition}

This proposition is proved using the same techniques as Theorem~\ref{thm:classificationPWPuptoQI}, with the following additional observation:

\begin{proposition}\label{prop:scalingQIbetweenPL}
Let $n,m\ge 2$. Let $G,H$ be finitely generated groups with normal subgroups $M\lhd G$, $N\lhd H$ of infinite index. Suppose that $\alpha\colon \Z_{n}^{(G/M)} \longrightarrow \Z_{m}^{(H/N)}$, $\beta\colon G\rightarrow H$ are two maps such that
\begin{align*}
    q\colon \mathcal{L}_{n}(G,\mathcal{C}_{M}) &\longrightarrow \mathcal{L}_{m}(H,\mathcal{C}_{N}) \\
    (c,p)&\longrightarrow (\alpha(c),\beta(p))
\end{align*}
is an aptolic quasi-isometry. If $\beta$ is quasi-$k$-to-one for some $k>0$, then so is $q$. 
\end{proposition}

\begin{proof}
Fix a finite set $A\subset \mathcal{L}_{m}(H,\mathcal{C}_{N})$. Let $\mathcal{S} \subset F^{(H,\mathcal{C}_{N})}$ denote the set of colourings appearing as first coordinates of elements of $A$, and for any $c\in \mathcal{S}$, let $A_{c}\subset A$ denote the subset of elements of $A$ having $c$ as first coordinate. Then 
\begin{equation*}
    |A|=\sum_{c\in\mathcal{S}}|A_{c}|, \; |q^{-1}(A)|=\sum_{c\in\mathcal{S}}|q^{-1}(A_{c})|
\end{equation*}
and letting $B_{c}\defeq \pi_{H}(A_{c})$, where $\pi_{H}\colon \mathcal{L}_{m}(H,\mathcal{C}_{N})\longrightarrow H$ is the canonical projection, we have $|q^{-1}(A_{c})|=|\beta^{-1}(B_{c})|$ for any $c\in\mathcal{S}$. Since $\beta$ is quasi-$k$-to-one, there exists a constant $C>0$ such that 
\begin{align*}
    \left|k|A|-|q^{-1}(A)|\right|&\le \sum_{c\in\mathcal{S}}\left|k|A_{c}|-|q^{-1}(A_{c})|\right| \\
    &=\sum_{c\in\mathcal{S}}\left|k|B_{c}|-|\beta^{-1}(B_{c})|\right| \\
    &\le\sum_{c\in\mathcal{S}} C\cdot|\partial_{H} B_{c}|. 
\end{align*}
It remains to notice that $\bigsqcup_{c\in\mathcal{S}}\lbrace c\rbrace\times \partial_{H}B_{c} \subset \partial A$ to deduce that 
\begin{equation*}
    \left|k|A|-|q^{-1}(A)|\right| \le C\cdot |\partial A|
\end{equation*}
for some constant $C>0$ independent of $A$. Thus $q$ is quasi-$k$-to-one and the proof is complete. 
\end{proof}

\begin{proof}[Proof of Proposition~\ref{prop:BILIPbetweenPLfromscalingQIbetweenquotients}]
Assume that $M$ is co-amenable in $G$, that $|E|=n^{r}$, $|F|=n^{s}$ for some $n,r,s\ge 1$ and that there is a quasi-one-to-one quasi-isometry of pairs 
\begin{equation*}
    f \colon (G,M)\longrightarrow (H,N)
\end{equation*}
such that $\overline{f}\colon G/M\rightarrow H/N$ is quasi-$\frac{s}{r}$-to-one. Up to finite distance, we may assume that $f(pM)\subset f(p)N$ for any $p\in G$. We proceed exactly as in the proof of Proposition~\ref{prop:aptolicQIfromscalingQIbetweenquotients-Acase}: we use Theorem~\ref{thm:rationalscalingfactor} to get partitions $\mathcal{P},\mathcal{Q}$ of $G/M$ and $H/N$ respectively, a bijection $\psi\colon \mathcal{P}\rightarrow \mathcal{Q}$ and a map $\beta\colon G/M\rightarrow H/N$, sending any piece $P\in\mathcal{P}$ into $\psi(\mathcal{P})$, that lies at finite distance from $\overline{f}$. Then, we can construct a map $h\colon G\rightarrow H$ at finite distance from $f$ such that $\overline{h}=\beta$, and this map allows us to define an aptolic quasi-isometry 
\begin{equation*}
    q\colon E\wr_{G/M}G\longrightarrow F\wr_{H/N}H
\end{equation*}
as in the proof of Proposition~\ref{prop:aptolicQIfromscalingQIbetweenquotients-Acase}. Since now $h$ is at finite distance from $f$, it is quasi-one-to-one, so $q$ is also quasi-one-to-one according to Proposition~\ref{prop:scalingQIbetweenPL}. We conclude from Theorem~\ref{thm:Whytethm} that this quasi-isometry lies at finite distance from a bijection 
\begin{equation*}
    E\wr_{G/M}G\longrightarrow F\wr_{H/N}H
\end{equation*}
which is the desired biLipschitz equivalence.
\end{proof}

We are in position to deduce Corollary~\ref{cor:classificationPWPoverZ^duptoBILIP}. We emphasize here that, even though the conclusion is the same as Corollary~\ref{cor:classificationofPWPoverZ^duptoQI}, the same proof does not work: indeed, given $|E|=n^{r}$ and $|F|=n^{s}$, we do not only need a quasi-isometry of pairs inducing a quasi-$\frac{s}{r}$-to-one at the level of quotients, but we must require this quasi-isometry of pairs to be quasi-one-to-one. Proceeding as in the proof of Corollary~\ref{cor:classificationofPWPoverZ^duptoQI} would only give us a quasi-$\frac{s}{r}$-to-one quasi-isometry of the base groups.

We must then explicitly construct such a quasi-isometry of pairs, in such a way that the induced map satisfies the required scaling condition. Let us first explain the proof below on a simple example. Consider $\Z^2$ with its $\Z-$cosets seen as vertical lines. Given a vertical line $\ell=(k,0)\Z \subset \Z^2$, the idea is to “split” it into two vertical lines $(2k,0)\Z$ and $(2k+1,0)\Z$, sending bijectively the “even part” of $\ell$ to $(2k,0)\Z$ and the “odd part” of $\ell$ to $(2k+1,0)\Z$. 

\begin{figure}[H]
  \centering
  \includegraphics[width=0.8\linewidth]{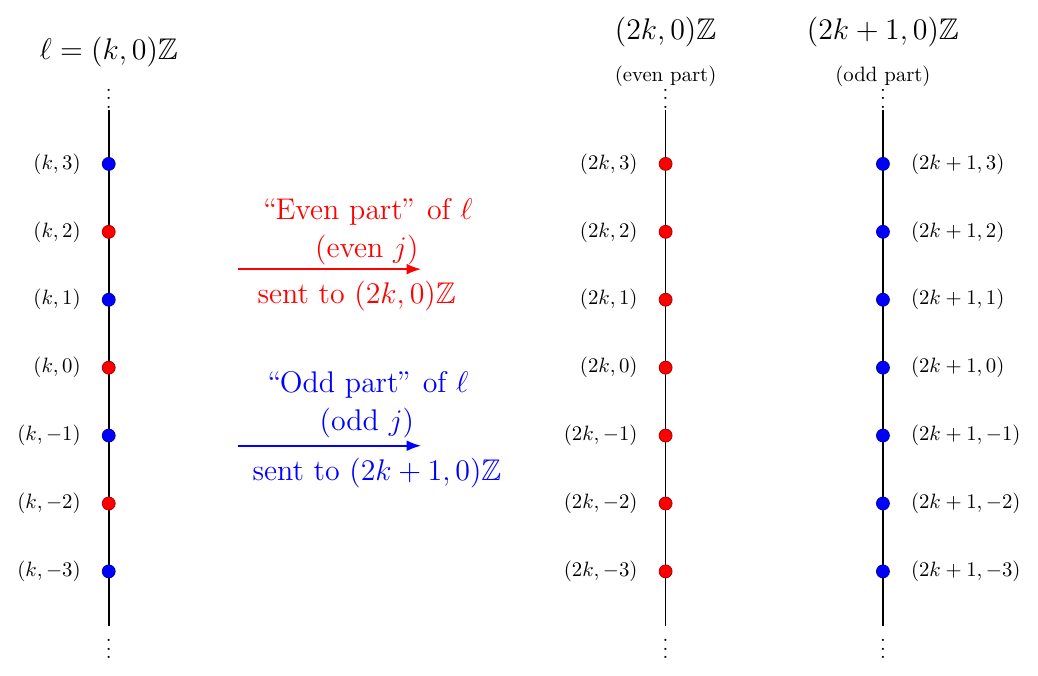}
  \caption{Splitting cosets of $\Z^2$.}
  \label{fig:chemin}
\end{figure}

More formally, we define the bijection 
\begin{align*}
    \gamma\colon \Z^2 &\longrightarrow \Z^2 \\
                 (i,j) &\longmapsto \begin{cases}
                     (2i, \frac{j}{2}) &\mbox{if $j$ is even} \\
                     (2i+1, \frac{j+1}{2}) &\mbox{if $j$ is odd}
                 \end{cases}
\end{align*}
that sends the coset $(k,0)\Z$ at Hausdorff distance at most one from the coset $(2k,0)\Z$, so that the induced map $\overline{\gamma}(k)=2k$ is quasi-$\frac{1}{2}$-to-one. Thus there is a biLipschitz equivalence
\begin{equation*}
    \Z_{2}^{2}\wr_{\Z}\Z^2 \longrightarrow \Z_{2}\wr_{\Z}\Z^2.
\end{equation*}

In higher dimensions, the idea is the same, except that we must split a coset into more cosets if we have a bigger lamp group. 

\begin{proof}[Proof of Corollary~\ref{cor:classificationPWPoverZ^duptoBILIP}]
If $E\wr_{\Z^{n}}\Z^{m}$ and $F\wr_{\Z^{n}}\Z^{m}$ are biLipschitz equivalent, they are quasi-isometric and we know from Corollary~\ref{cor:classificationofPWPoverZ^duptoQI} that $|E|,|F|$ are powers of a common number. 

\noindent Conversely, suppose that $|E|=k^{r}$ and $|F|=k^{s}$. Without loss of generality, we may assume that $E=\Z_{k}^{r}$ and $F=\Z_{k}^{s}$. As a first reduction, notice that it is enough to prove that there is a biLipschitz equivalence
\begin{equation*}
    \Z_{k}^{r}\wr_{\Z^{n}}\Z^{m}\longrightarrow \Z_{k}\wr_{\Z^{n}}\Z^{m}.
\end{equation*}
To prove this claim, following Proposition~\ref{prop:BILIPbetweenPLfromscalingQIbetweenquotients}, it is enough to exhibit a quasi-one-to-one quasi-isometry $\gamma\colon \Z^{m}\rightarrow \Z^{m}$ quasi-preserving the cosets of $\Z^{m-n}$, such that the induced quasi-isometry $\overline{\gamma}\colon \Z^{n}\rightarrow \Z^{n}$ is quasi-$\frac{1}{r}$-to-one. Consider the map 
\begin{align*}
    \gamma\colon \Z^{m}&\longrightarrow \Z^{m} \\
    (x_{1},\dots, x_{m})&\longmapsto 
    \begin{cases}
        \left(rx_{1},x_{2},\dots,x_{m-1},\frac{x_{m}}{r}\right)  &\mbox{if}\; x_{m}\equiv 0 \;[r] \\
        \left(rx_{1}+1,x_{2},\dots,x_{m-1},\frac{x_{m}+(r-1)}{r}\right) &\mbox{if}\; x_{m}\equiv 1 \;[r] \\
        \quad\quad\quad\quad\quad\quad\vdots & \quad\quad\vdots \\
        \left(rx_{1}+(r-1),x_{2},\dots,x_{m-1},\frac{x_{m}+1}{r}\right) &\mbox{if} \;x_{m}\equiv r-1 \;[r]
    \end{cases}
\end{align*}
and its inverse 
\begin{align*}
    \eta\colon \Z^{m}&\longrightarrow \Z^{m} \\
    (p_{1},\dots, p_{m})&\longmapsto \begin{cases}
        \left(\frac{p_{1}}{r},p_{2},\dots,p_{m-1},rp_{m}\right)  &\mbox{if}\; p_{1}\equiv 0 \;[r] \\
        \left(\frac{p_{1}-1}{r},p_{2},\dots,p_{m-1},rp_{m}-(r-1)\right)  &\mbox{if}\; p_{1}\equiv 1 \;[r] \\
        \quad\quad\quad\quad\quad\quad\vdots & \quad\quad\vdots \\
        \left(\frac{p_{1}-(r-1)}{r},p_{2},\dots,p_{m-1},rp_{m}-1\right)  &\mbox{if}\; p_{1}\equiv r-1 \;[r] \\
        \end{cases}.
\end{align*}
One directly checks with Lemma~\ref{lem:Lipschitzbetweengraphs} that $\gamma$ and $\eta$ are both Lipschitz maps, inverses of each other, and thus are biLipschitz equivalences. In particular, they are quasi-one-to-one, and they induce two maps given by
\begin{align*}
    \overline{\gamma}\colon \Z^{n}&\longrightarrow \Z^{n} \\
    (x_{1},\dots, x_{n}) &\longmapsto (rx_{1},x_{2},\dots, x_{n})
\end{align*}
and
\begin{align*}
    \overline{\eta}\colon \Z^{n}&\longrightarrow \Z^{n} \\
    (p_{1},\dots, p_{n})&\longmapsto \begin{cases}
        \left(\frac{p_{1}}{r},p_{2},\dots,p_{n}\right)  &\mbox{if}\; p_{1}\equiv 0 \;[r] \\
        \left(\frac{p_{1}-1}{r},p_{2},\dots,p_{n}\right)  &\mbox{if}\; p_{1}\equiv 1 \;[r] \\
        \quad\quad\quad\quad\vdots  & \quad\quad\vdots \\
        \left(\frac{p_{1}-(r-1)}{r},p_{2},\dots,p_{n}\right)  &\mbox{if}\; p_{1}\equiv r-1 \;[r] \\
        \end{cases}.
\end{align*}
We know from Lemma~\ref{prop:compatibilityoperationsquotients}\textit{(ii)} that $\overline{\eta}$ is a quasi-inverse of $\overline{\gamma}$, and given any $(p_{1},\dots, p_{n})\in\Z^n$, we see that $\overline{\eta}^{-1}\big(\lbrace (p_{1},\dots, p_{n})\rbrace\big)$ equals
\begin{equation*}
    \big\lbrace (rp_{1},\dots,p_{n}), (rp_{1}+1,p_{2},\dots, p_{n}),\dots, (rp_{1}+(r-1),p_{2},\dots, p_{n}) \big\rbrace. 
\end{equation*}
In particular, $\overline{\eta}$ is $r$-to-one, a fortiori quasi-$r$-to-one, and it follows from Proposition~\ref{prop:stabilitypropertiesforscalingQI}\textit{(iii)} that $\overline{\gamma}$ is quasi-$\frac{1}{r}$-to-one, as desired. This concludes the proof. 
\end{proof}

\subsection{More wreath products}\label{subsubsection3.6.3} Let us now focus on Corollary~\ref{cor:classificationNAlamplightersoverPWP}. As indicated in the introduction, we decompose its proof into several intermediate steps, which individually use less assumptions than in the statement. We start with the flexibility part: 

\begin{proposition}\label{prop6.3}
Let $n,m,p,q\ge 2$. Let $H$ be a finitely generated group with a normal subgroup $N\lhd H$ of infinite index. Suppose that $N$ is not co-amenable in $H$. If $n$ and $m$ have the same prime divisors, and $p$ and $q$ have the same prime divisors, then there exists a quasi-isometry 
\begin{equation*}
    \Z_{n}\wr(\Z_{p}\wr_{H/N}H) \longrightarrow \Z_{m}\wr(\Z_{q}\wr_{H/N}H).
\end{equation*}
\end{proposition}

\begin{proof}
Assume that $n$ and $m$ have the same prime divisors, and $p$ and $q$ have the same prime divisors. The latter implies that there is a quasi-isometry 
\begin{equation*}
    \Z_{p}\wr_{H/N}H\longrightarrow \Z_{q}\wr_{H/N}H
\end{equation*}
using the second point of Theorem~\ref{thm:classificationPWPuptoQI}. As $N$ is not co-amenable in $H$, $H$ cannot be amenable (otherwise $H/N$ would be), hence $\Z_{p}\wr_{H/N}H$ and $\Z_{q}\wr_{H/N}H$ are not amenable. We can then use Theorem~\ref{thm1.2} to deduce that there is an aptolic quasi-isometry
\begin{equation*}
    \Z_{n}\wr(\Z_{p}\wr_{H/N}H) \longrightarrow \Z_{m}\wr(\Z_{q}\wr_{H/N}H)
\end{equation*}
as claimed. 
\end{proof}

\begin{proposition}\label{prop:lamplightersoverPWP-rigiditypart}
Let $n,m,p,q\ge 2$. Let $H$ be a one-ended finitely presented group with a finitely generated infinite normal subgroup $N\lhd H$ of infinite index. Suppose that $H$ is not coarsely separable by any collection of subspaces that uniformly quasi-isometrically embed into $N$. Suppose that $\Z_{n}\wr(\Z_{p}\wr_{H/N}H)$ and $\Z_{m}\wr(\Z_{q}\wr_{H/N}H)$ are quasi-isometric. The following claims hold.
\begin{enumerate}[label=(\roman*)]
    \item If $H$ is amenable, then $n$ and $m$ are powers of a common number, $p=k^{r}$ and $q=k^{s}$ are powers of a common number and there exists a quasi-isometry of pairs $(H,N)\longrightarrow (H,N)$ inducing a quasi-$\frac{s}{r}$-to-one quasi-isometry $H/N\rightarrow H/N$;
    \item If $N$ is not co-amenable in $H$, then $n$ and $m$ have the same prime divisors, and $p$ and $q$ have the same prime divisors.
\end{enumerate}
\end{proposition}

\begin{proof}
\textit{(i)} By Theorem~\ref{thm:AlamplightersoverTBPgroups-rigiditypart}, which applies since $\Z_{p}\wr_{H/N}H$, $\Z_{q}\wr_{H/N}H$ have the thick bigon property (see Corollary~\ref{cor:PWPhaveTBP}) and are both amenable, we deduce that $n$ and $m$ are powers of a common number and that there exists a (measure-scaling) quasi-isometry 
\begin{equation*}
    \Z_{p}\wr_{H/N}H\longrightarrow \Z_{q}\wr_{H/N}H.
\end{equation*}
Hence the conclusion follows from the first point of Theorem~\ref{thm:classificationPWPuptoQI}. 

\noindent \textit{(ii)} Assume that there is a quasi-isometry 
\begin{equation}\label{eq6.1}
    \Z_{n}\wr(\Z_{p}\wr_{H/N}H) \longrightarrow \Z_{m}\wr(\Z_{q}\wr_{H/N}H).
\end{equation}
In this case, we start rather with Proposition~\ref{prop:lamplightersoverTBPgroups-rigiditypart}, which applies since $\Z_{p}\wr_{H/N}H$, $\Z_{q}\wr_{H/N}H$ have the thick bigon property (see Corollary~\ref{cor:PWPhaveTBP}), to deduce that $n$ and $m$ have the same prime divisors. Additionally, as explained in the proof of Proposition~\ref{prop:lamplightersoverTBPgroups-rigiditypart}, the quasi-isometry (\ref{eq6.1}) can be assumed to be aptolic (in the sense of Definition~\ref{def:aptolicQIbetweenPWP}), and thus provides a quasi-isometry $\Z_{p}\wr_{H/N}H \longrightarrow \Z_{q}\wr_{H/N}H$. Now Theorem~\ref{thm:classificationPWPuptoQI} implies that $p$ and $q$ have the same prime divisors, and the proof is complete. 
\end{proof}

\begin{proof}[Proof of Corollary~\ref{cor:classificationNAlamplightersoverPWP}]
Combine Proposition~\ref{prop6.3} and point~\textit{(ii)} of Proposition~\ref{prop:lamplightersoverPWP-rigiditypart}.
\end{proof}

In the situation of an amenable permutational wreath product, Theorem~\ref{thm:classificationPWPuptoQI} imposes an additional scaling condition, so we must determine scaling groups of such permutational wreath products in order to prove Proposition~\ref{prop:classificationAlamplightersoverPWPoverZ^d}. We can compute some of them thanks to the next observation:

\begin{lemma}\label{lem:scalinggroupsofdirectproducts}
Let $G$ and $H$ be finitely generated groups. Then $\text{Sc}(G)$ and $\text{Sc}(H)$ are subgroups of $\text{Sc}(G\times H)$.
\end{lemma}

\begin{proof}
The claim is symmetric in $G$ and $H$ so we do the proof for $G$ only. 

\noindent Let $k\in\text{Sc}(G)$ and let $f\colon G\rightarrow G$ be a $(C,K)-$quasi-isometry which is quasi-$k$-to-one. We define then 
\begin{align*}
    \Tilde{f}\colon G\times H&\longrightarrow G\times H \\
    (g,h)&\longmapsto (f(g), h).
\end{align*}
One directly checks that $\Tilde{f}$ is a quasi-isometry, and we claim it is also quasi-$k$-to-one. Indeed, fix a finite subset $A\subset G\times H$. Then $A$ intersects finitely many (say $r\ge 1$) $G-$cosets in $G\times H$, and we can decompose 
\begin{equation*}
    A=\bigsqcup_{i=1}^{r}A_{i}
\end{equation*}
where, for any $1\le i\le r$, $A_{i}$ is contained in a single $G-$coset. Moreover, denoting $\pi_{G}\colon G\times H \longrightarrow G$ the projection on the first factor, one has 
\begin{equation*}
    |A_{i}|=|\pi_{G}(A_{i})| \; \text{and}\; \Tilde{f}^{-1}(A_{i})=f^{-1}(\pi_{G}(A_{i}))
\end{equation*}
for any $1\le i\le r$. Hence it follows that
\begin{align*}
    \left|k|A|-|\Tilde{f}^{-1}(A)|\right| &\le \sum_{i=1}^{r}\left|k|A_{i}|-|\Tilde{f}^{-1}(A_{i})|\right| \\
    &\le \sum_{i=1}^{r}\big|k|\pi_{G}(A_{i})|-|f^{-1}(\pi_{G}(A_{i}))|\big| \\
    &\le C\cdot \sum_{i=1}^{r}|\partial_{G} \pi_{G}(A_{i})| \\
    &\le C\cdot \left|\partial_{G\times H} A\right|
\end{align*}
for some constant $C>0$, using that $f$ is quasi-$k$-to-one and the inclusion
\begin{equation*}
    \bigsqcup_{i=1}^{r}\partial_{G} \pi_{G}(A_{i}) \subset \partial_{G\times H} A
\end{equation*} 
for the last inequality. This proves that $k\in\text{Sc}(G\times H)$ as claimed. 
\end{proof}

\begin{corollary}\label{cor:scalinggroupsofPWP}
Let $F$, $N$ and $K$ be finitely generated groups. Then $\text{Sc}(N)$ is a subgroup of $\text{Sc}\big(F\wr_{K}(N\times K)\big)$.
\end{corollary}

\begin{proof}
This follows from Lemma~\ref{lem:scalinggroupsofdirectproducts} and the fact that $F\wr_{K}(N\times K) \cong N\times (F\wr K)$.
\end{proof}

Thus, for instance, if $F$ is finite and $m>n\ge 1$, $\text{Sc}(F\wr_{\Z^n}\Z^m)=\R_{>0}$ as it contains $\text{Sc}(\Z^{m-n})=\R_{>0}$. This leads us to a proof of the next statement:

\begin{proposition}\label{prop6.7}
Let $n,m,p\ge 2$. Let $N$ and $K$ be finitely presented infinite amenable groups, and let $G\defeq N\times K$. If $\text{Sc}(N)=\R_{>0}$, then $\Z_{n}\wr(\Z_{p}\wr_{K}G)$ and $\Z_{m}\wr(\Z_{p}\wr_{K}G)$ are quasi-isometric if and only if $n$ and $m$ are powers of a common number.
\end{proposition}

\begin{proof}
By Theorem~\ref{thm:AlamplightersoverTBPgroups-rigiditypart}, our two groups are quasi-isometric if and only if there exist $k,r,s\ge 1$ such that $n=k^{r}$, $m=k^{s}$ and $\frac{s}{r}\in\text{Sc}(\Z_{p}\wr_{K}G)$. By Corollary~\ref{cor:scalinggroupsofPWP}, $\text{Sc}(\Z_{p}\wr_{K}G)$ contains $\text{Sc}(N)=\R_{>0}$, so that $\text{Sc}(\Z_{p}\wr_{K}G)=\R_{>0}$. Hence we conclude that $\Z_{n}\wr(\Z_{p}\wr_{K}G)$ and $\Z_{m}\wr(\Z_{p}\wr_{K}G)$ are quasi-isometric if and only if $n$ and $m$ are powers of a common number.
\end{proof}

With the same techniques, we may now deduce a proof of Proposition~\ref{prop:classificationAlamplightersoverPWPoverZ^d}.

\begin{proof}[Proof of Proposition~\ref{prop:classificationAlamplightersoverPWPoverZ^d}]
Suppose first that $\Z_{n}\wr(\Z_{p}\wr_{\Z^{k}}\Z^{d})$ and $\Z_{m}\wr(\Z_{q}\wr_{\Z^{k}}\Z^{d})$ are quasi-isometric. Then Theorem~\ref{thm:AlamplightersoverTBPgroups-rigiditypart} implies that $n=a^{r}$, $m=a^{s}$ for some $a,r,s\ge 1$ and there exists a quasi-$\frac{s}{r}$-to-one quasi-isometry $\Z_{p}\wr_{\Z^{k}}\Z^d \longrightarrow \Z_{q}\wr_{\Z^{k}}\Z^d$. Then Corollary~\ref{cor:classificationofPWPoverZ^duptoQI} shows that $p$ and $q$ must be powers of a common number, as claimed. 

\noindent Conversely, suppose that $n=a^{r}$ and $m=a^{s}$ are powers of a common number, and likewise that $p$ and $q$ are powers of a common number. This assumption implies, according to Corollary~\ref{cor:classificationPWPoverZ^duptoBILIP}, that there is a biLipschitz equivalence 
\begin{equation*}
    \Z_{p}\wr_{\Z^{k}}\Z^{d} \longrightarrow \Z_{q}\wr_{\Z^{k}}\Z^{d}.
\end{equation*}
We also know that $\text{Sc}(\Z_{q}\wr_{\Z^{k}}\Z^{d})=\R_{>0}$, so $\frac{s}{r}\in \text{Sc}(\Z_{q}\wr_{\Z^{k}}\Z^{d})$, and thus there exists a quasi-$\frac{s}{r}$-to-one quasi-isometry $\Z_{q}\wr_{\Z^{k}}\Z^{d} \longrightarrow \Z_{q}\wr_{\Z^{k}}\Z^{d}$. Composing this quasi-isometry with the above biLipschitz equivalence provides a quasi-$\frac{s}{r}$-to-one quasi-isometry
\begin{equation*}
    \Z_{p}\wr_{\Z^{k}}\Z^{d} \longrightarrow \Z_{q}\wr_{\Z^{k}}\Z^{d}.
\end{equation*}
Since we know that $n=a^{r}$ and $m=a^{s}$,~\cite[Proposition~3.12]{GT24b} implies that $\Z_{n}\wr(\Z_{p}\wr_{\Z^{k}}\Z^{d})$ and $\Z_{m}\wr(\Z_{q}\wr_{\Z^{k}}\Z^{d})$ are quasi-isometric. The proof is complete. 
\end{proof}

In~\cite[Corollary~1.16]{GT24b} are also recorded the first examples of amenable groups having all their self-quasi-isometries at a bounded distance from a bijection, namely $F\wr H$ where $F$ is finite and $H$ is amenable finitely presented and one-ended. In particular, the scaling group of such wreath products is always trivial. In contrast, permutational wreath products have many measure-scaling quasi-isometries that do not lie within bounded distance from a bijection. An explicit example is the following: let $G\defeq \Z_{2}\wr_{\Z^2}\Z^3$, and let $\pi\colon G\rightarrow \Z^3$ be the natural surjection. Then $K_{2}\defeq \pi^{-1}(\Z\times \Z\times 2\Z)$ is a proper subgroup of $G$, of index $2$. Thus the natural inclusion $\iota\colon K_{2}\hookrightarrow G$ is quasi-$\frac{1}{2}$-to-one, and since $G\cong K_{2}$, composing this isomorphism with $\iota$ provides a quasi-$\frac{1}{2}$-to-one quasi-isometry $G\rightarrow G$. 

This example also allows us to exclude the natural analog of~\cite[Corollary~1.17]{GT24b} in the permutational case: $K_{2}$ and $K_{3}\defeq \pi^{-1}(\Z\times \Z\times 3\Z)$ are biLipschitz equivalent (even isomorphic) finite index subgroups of $G=\Z_{2}\wr_{\Z^2}\Z^3$, but they do not have the same index in $G$.

\section{Concluding remarks and questions}\label{subsection3.7}

We conclude the article with a list of questions related to the coarse geometry of permutational wreath products. 

First of all, the proof of Theorem~\ref{thm:classificationPWPuptoQI} requires strong assumptions, and it is natural to ask what happens if we drop some of these assumptions. For instance:

\begin{question}
Let $n,m\ge 2$. Let $G$ (resp. $H$) be a finitely presented group with a finitely generated subgroup $M\leqslant G$ (resp. $N\leqslant H$) of infinite index. Assume that $G$ (resp. $H$) is not coarsely separable by any collection of subspaces that uniformly quasi-isometrically embed into $M$ (resp. $N$). When are $\Z_{n}\wr_{G/M}G$ and $\Z_{m}\wr_{H/N}H$ quasi-isometric? biLipschitz equivalent?
\end{question}

In the same vein, notice that in~\cite{GT24a} the authors extended their embedding theorem from~\cite{GT24b} to include some standard lamplighters with infinitely presented base groups. Namely, they introduced a new quasi-isometry invariant, referred to as the~\textit{thick bigon property} (see Chapter~\ref{chap:chapter2}), and proved a version of the embedding theorem for geodesic metric spaces satisfying this property (see Theorem~\ref{thm:embedding2}). Therefore, one might wonder whether there is a version of this embedding theorem for permutational halo products. As an application, this could provide an answer to:

\begin{question}
Let $n,m,p,q\ge 2$, and let $H$ be a finitely generated group with $N$ a normal subgroup of infinite index. Letting $\Z_{p}\wr_{H/N}H$ (resp. $\Z_{q}\wr_{H/N}H$) act on its quotient $H$, when are $\Z_{n}\wr_{H}(\Z_{p}\wr_{H/N}H)$ and $\Z_{m}\wr_{H}(\Z_{q}\wr_{H/N}H)$ quasi-isometric? 
\end{question}

Some partial answers can be deduced from our previous results. For instance, in the case where $N$ is not co-amenable in $H$, if $p$ and $q$ have the same prime divisors, then by Proposition~\ref{prop:constructingQIbetweenNAPL} there exists a quasi-isometry of pairs 
\begin{equation*}
    \left(\Z_{p}\wr_{H/N}H, \;\bigoplus_{H/N}\Z_{p}\right) \longrightarrow \left(\Z_{q}\wr_{H/N}H,\;\bigoplus_{H/N}\Z_{q}\right)
\end{equation*}
and if additionally $n$ and $m$ have the same prime divisors, there exists a quasi-isometry
\begin{equation*}
    \Z_{n}\wr_{H}(\Z_{p}\wr_{H/N}H) \longrightarrow \Z_{m}\wr_{H}(\Z_{q}\wr_{H/N}H).
\end{equation*}
still by Proposition~\ref{prop:constructingQIbetweenNAPL}. A possible refinement of the previous question is then:

\begin{question}
Let $n,m,p,q\ge 2$, and let $H$ be a finitely generated group with $N$ a normal subgroup of infinite index. Assume that $N$ is not co-amenable in $H$. If $\Z_{n}\wr_{H}(\Z_{p}\wr_{H/N}H)$ and $\Z_{m}\wr_{H}(\Z_{q}\wr_{H/N}H)$ are quasi-isometric, must $n$ and $m$ (resp. $p$ and $q$) have the same prime divisors?
\end{question}

In the context of iterated wreath products, it is worth to notice that the assumption of $N$ being infinite in Corollary~\ref{cor:classificationNAlamplightersoverPWP} is crucial to apply results from~\cite{GT24a}. Indeed, as explained in the latter, permutational wreath products of the form $E\wr_{G/M}G$ with $M$ infinite satisfy the thick bigon property, which allows to use the embedding theorem from~\cite{GT24a} to get information on the quasi-isometries between wreath products of the form $\Z_{n}\wr(E\wr_{G/M}G)$. On the other hand, standard wreath products of the form $(\text{finite})\wr G$ usually do not have the thick bigon property, even if $G$ has it. Hence, the following question remains:

\begin{question}
Let $n,m,p,q\ge 2$ be four integers. Let $G$ and $H$ be one-ended finitely presented groups. If $\Z_{n}\wr(\Z_{p}\wr G)$ and $\Z_{m}\wr(\Z_{q}\wr H)$ are quasi-isometric, must $n$ and $m$ (resp. $p$ and $q$) have the same prime divisors? Must $G$ and $H$ be quasi-isometric?
\end{question}

In another direction, it is natural to wonder what happens when the lamp groups are infinite, in the same spirit as~\cite[Theorem~6.33 and Corollary~6.36]{BGT24}. For instance:

\begin{question}
Let $A_{1}$, $A_{2}$ be two nilpotent groups, and $H_{1}$, $H_{2}$ two finitely presented groups, with normal infinite finitely generated subgroups $N_{1}\lhd H_{1}$, $N_{2}\lhd H_{2}$. When are $A_{1}\wr_{H_{1}/N_{1}}H_{1}$ and $A_{2}\wr_{H_{2}/N_{2}}H_{2}$ quasi-isometric?
\end{question}

The case of free abelian groups would be interesting already. For instance, if $\Z^{n}\wr_{\Z^{k}}\Z^{d}$ and $\Z^{n'}\wr_{\Z^{k'}}\Z^{d'}$ are quasi-isometric, does it follow that $d=d'$ and $k=k'$? Or does one have more flexibility? Note that when the subgroup is not co-amenable into the ambient group, Remark~\ref{rm5.2} provides some flexibility in the lamp groups. 

Note that, as explained below in Remark~\ref{rem:remark4.5.2}, a quasi-isometry between $\Z^{n}\wr_{\Z^{k}}\Z^{d}$ and $\Z^{n'}\wr_{\Z^{k'}}\Z^{d'}$ does not necessarily imply that $n=n'$. 

Lastly, Theorem~\ref{thm:classificationPWPuptoQI} requires a strong assumption of non-coarse separation, and it is natural to ask what happens in the remaining cases. For instance:

\begin{question}\label{question7.6}
If $E$ and $F$ are two non-trivial finite groups and $m\ge 2$, and if $E\wr_{\Z}\Z^{m}$ and $F\wr_{\Z}\Z^{m}$ are quasi-isometric, must $|E|$ and $|F|$ be powers of a common number?
\end{question}

We conjecture a positive answer. Note that Example~\ref{ex:aQIwhichisnotleafpreserving} provides an example of a non-leaf-preserving quasi-isometry $\Z_{2}\wr_{\Z}\Z^2\longrightarrow \Z_{2}\wr_{\Z}\Z^2$, hence the need of another strategy. When looking at quasi-isometries between spaces that are direct products, such as
\begin{equation*}
    E\wr_{\Z}\Z^{m} \cong \Z^{m-1}\times (E\wr\Z),
\end{equation*}
another natural form of rigidity that can occur is the one of “product quasi-isometries” (see for instance~\cite{KL97, EFW12, EFW13}), namely those quasi-isometries that are given by a quasi-isometry of the first factor in the first component and a quasi-isometry of the second factor in the second component. Moreover, in our case, such a rigidity would allow us to deduce a positive answer to Question~\ref{question7.6}, thanks to~\cite[Theorem~1.2]{EFW13}. 

However, it is not hard to construct automorphisms (and thus quasi-isometries) of such products that do not lie at bounded distance from a 
product quasi-isometry. One such example is given by $\Z\times (E\wr\Z) \longrightarrow \Z\times (E\wr\Z)$, $(k,(c,p))\longmapsto (k+p, (c,p))$. Hence, for now, the quasi-isometric classification of direct products of $\Z$ with lamplighters over $\Z$ stays out of reach.  

Next, regarding biLipschitz equivalences between permutational wreath products, it is natural to ask whether the sufficient condition given by Proposition~\ref{prop:BILIPbetweenPLfromscalingQIbetweenquotients} is also necessary, at least under the assumptions of Theorem~\ref{thm:classificationPWPuptoQI}.
\begin{question}
Let $n,m\ge 2$ and let $G,H$ be finitely presented groups with finitely generated infinite normal subgroups $M\lhd G$, $N\lhd H$ of infinite index. Suppose that $G$ (resp. $H$) is not coarsely separable by any collection of subspaces that uniformly quasi-isometrically embed into $M$ (resp. $N$). If $\Z_{n}\wr_{G/M}G$ and $\Z_{m}\wr_{H/N}H$ are biLipschitz equivalent, does there exist $k,r,s\ge 1$ such that $n=k^r$, $m=k^s$ and a quasi-one-to-one quasi-isometry of pairs $(G,M)\longrightarrow (H,N)$ inducing a quasi-$\frac{s}{r}$-to-one quasi-isometry $G/M\rightarrow H/N$?
\end{question}

We showed in Corollary~\ref{cor:scalinggroupsofPWP} that $\text{Sc}(F\wr_{G/N}G)$ always contains $\text{Sc}(N)$ if $G$ splits as a direct product $N\times K$. It is therefore natural to ask:

\begin{question}
Let $F$ and $G$ be finitely generated groups, and let $N\leqslant G$ be finitely generated. Is it true that $\text{Sc}(N)$ injects into $\text{Sc}(F\wr_{G/N}G)$?
\end{question}

In fact, we showed in Lemma~\ref{lem:scalinggroupsofdirectproducts} that the scaling group of a direct product always contains scaling groups of the factors. More generally:

\begin{question}
Let $G$ and $H$ be finitely generated groups, and let $\varphi\colon H\longrightarrow \text{Aut}(G)$ be a group morphism. Is it true that $\text{Sc}(G\rtimes_{\varphi}H)$ contains $\text{Sc}(G)$?
\end{question}

On the other hand, the scaling group of $G\rtimes_{\varphi}H$ usually does not contain $\text{Sc}(H)$. For instance, $\text{Sc}(\Z_{2}\wr\Z^2)=\lbrace 1\rbrace$ by~\cite[Corollary~1.16]{GT24b}, while $\text{Sc}(\Z^2)=\R_{>0}$.

\afterpage{\blankpage}
\chapter{Quasi-isometric rigidity of lamplighters with lamps of polynomial growth}\label{chap:chapter4}

This chapter presents the content of the paper~\cite{Dum26}, that explains a rigidity phenomenon for quasi-isometries between wreath products with infinite lamp groups having polynomial growth. 

\vspace{0.3cm}

\minitoc

\section{Introduction}\label{sec:introchapter4}

\setcounter{theoremletter}{0}

In geometric group theory, a challenging question is to understand the collection of all maps between two given metric spaces that are compatible with the large-scale geometry of our spaces. Such collections encompass for instance quasi-isometries, coarse and quasi-isometric embeddings, or regular maps. The motivation behind this program is that, on the one hand finitely generated groups are naturally metric spaces and can thus be studied from a geometric point of view, and on the other hand their geometric properties are closely related to their algebraic structure. 

\sloppy Several milestones have been achieved in the study of the quasi-isometric rigidity of many classes of groups, among which abelian and non-abelian free groups~\cite{Dun85}, mapping class groups~\cite{Beh+12}, non-uniform lattices of $\text{Isom}(\mathbb{H}^{n})$~\cite{Sch96}, $n\ge 3$, lamplighters over $\Z$ and others $\text{SOL}-$like groups~\cite{EFW12, EFW13}, or Baumslag-Solitar groups~\cite{FM98, FM99, Why01}.

In a recent work~\cite{GT24b}, Genevois and Tessera established a complete classification up to quasi-isometry of lamplighters over finitely presented one-ended groups. Precisely, they showed that the existence of a quasi-isometry between $E\wr H$ and $F\wr K$ is parametrized by the existence of a quasi-isometry between $H$ and $K$ which is compatible in a certain way with the lamp groups $E$ and $F$, in the sense that it must be measure-scaling for some scaling factor that depends on the cardinalities of $E$ and $F$. Roughly speaking, a quasi-isometry $X\rightarrow Y$ between two metric spaces is quasi-$k$-to-one if the size of the pre-image of a finite subset $A$ of $Y$ is $k\left|A\right|$, up to an error controlled uniformly over $A\subset Y$ (a precise definition is given in Definition~\ref{def:measurescalingQI}). This quantification is also motivated by a famous result of Whyte (see Theorem~\ref{thm:Whytethm}), which says that being quasi-one-to-one is the same as lying within bounded distance from a bijection. Furthermore, the property of being measure-scaling is well-behaved under composition and quasi-inverses (see Proposition~\ref{prop:stabilitypropertiesforscalingQI}), which yields to define the scaling group of an amenable space $X$ as
\begin{equation*}
    \text{Sc}(X) \defeq \left\lbrace k>0 : \exists \; \text{a quasi-$k$-to-one quasi-isometry $X\rightarrow X$}\right\rbrace. 
\end{equation*}

\sloppy In this way, $\text{Sc}(X)$ is a subgroup of $(\R_{>0},\cdot)$. Given an amenable space $X$, its scaling group is an object of independent interest, that encodes many possible behaviours for its self-quasi-isometries. In the case of $X=G$ a finitely generated group, identifying its scaling group can give access to information on its algebraic structure. For instance, it is not hard to check that, given a finite-index subgroup $H$ of $G$, the natural inclusion $H\hookrightarrow G$ is quasi-$\frac{1}{[G:H]}$-to-one. Thus, if $G$ has trivial scaling group, it cannot contain proper finite index subgroups isomorphic to itself. We refer to Section~\ref{sec:ScalingQI} for more details and explanations on scaling quasi-isometries and scaling groups. 

Scaling groups of several classes of amenable groups, among which Carnot groups and Baumslag-Solitar groups, have been computed in~\cite[Corollary~6.6]{GT22}. Additionally, as an important consequence of the classification established in~\cite[Corollary~1.6]{GT24b}, \mbox{$\text{Sc}(F\wr K)=\lbrace 1\rbrace$} if $F$ is finite and $K$ is amenable, finitely presented and one-ended. As explained in Section~\ref{sec:Scalinggroupsofhaloproducts}, this computation heavily relies on the aptolicity (in the sense of Definition~\ref{def:aptolicity}) rigidity phenomenon of quasi-isometries between lamplighters.

Here we show that, for lamplighters with lamps of polynomial growth, if all quasi-isometries are aptolic (up to bounded distances), then they are also all measure-scaling, with a scaling factor depending on the involved growth degrees. 

\begin{theoremletter}\label{thm:maintheorem}
Let $N$ and $M$ be finitely generated groups of polynomial growth, with growth degrees $n$ and $m$ respectively. Let $G$ and $H$ be finitely generated amenable groups such that any quasi-isometry between $N\wr G$ and $M\wr H$ is (up to bounded distance) aptolic. Then any quasi-isometry $N\wr G\longrightarrow M\wr H$ is quasi-$\frac{m}{n}$-to-one.
\end{theoremletter}

Explicit examples of groups $G$ and $H$ for which any quasi-isometry $N\wr G\longrightarrow M\wr H$ is aptolic, up to a finite distance change, are provided in~\cite{BGT24}. Thus:

\begin{theoremletter}\label{thm:maintheoremfortheclassMexp}
Let $N$ and $M$ be finitely generated groups of polynomial growth, with growth degrees $n$ and $m$ respectively. Let $G$ and $H$ be finitely presented amenable groups from $\mathcal{M}_{\text{exp}}$. Then any quasi-isometry $N\wr G \longrightarrow M\wr H$ is quasi-$\frac{m}{n}$-to-one.  \end{theoremletter}

We refer to Section~\ref{subsection4.2} for the precise definition of the class $\mathcal{M}_{\text{exp}}$. For now, let us simply mention that it contains many amenable groups: solvable Baumslag-Solitar groups $\text{BS}(1,n)$, $n\ge 2$, $\text{SOL}(\Z)$, lamplighter groups, as well as any direct product of one of these groups with an arbitrary finitely generated amenable group.

We also emphasize that, in these theorems, the amenability condition is necessary to avoid an empty statement: between non-amenable spaces, a quasi-isometry is measure-scaling for all possible scaling factors; see Remark~\ref{rem:ScalingQINAspaces}.

\paragraph{Quasi-isometries and biLipschitz equivalences.} Refining our understanding of quasi-isometries, another interesting problem is to understand to what extent being quasi-isometric and being biLipschitz equivalent differ. Whyte's theorem mentioned above (see Theorem~\ref{thm:Whytethm}) ensures that, for non-amenable spaces, a quasi-isometry can always be turned into a biLipschitz equivalence by a finite distance change. For amenable spaces, the question is much more subtle, and no analog of Whyte's theorem can hold: indeed, for a proper finite-index subgroup $H$ of an amenable finitely generated group $G$, the inclusion $H\hookrightarrow G$ is quasi-$\frac{1}{[G:H]}$-to-one and thus does not lie at bounded distance from a bijection. The first example of a pair of amenable quasi-isometric groups that are not biLipschitz equivalent appears in~\cite{Dym10} and involves lamplighters over $\Z$. 

Concerning biLipschitz equivalences between lamplighters with infinite lamp groups, we can deduce from our main result the next consequence:

\begin{corollaryletter}\label{cor:BLE}
Let $N$ and $M$ be finitely generated groups with polynomial growth. Let $G$ and $H$ be finitely presented amenable groups from $\mathcal{M}_{\text{exp}}$. If $N\wr G$ and $M\wr H$ are biLipschitz equivalent, then $N$ and $M$ have the same growth degree and $G$ and $H$ are biLipschitz equivalent. 
\end{corollaryletter}

It is worth noticing that the converse statement is not true: according to~\cite[Corollary~6.36]{BGT24}, if $N$ and $M$ are nilpotent and $N\wr G$ and $M\wr H$ are quasi-isometric, then $N$ and $M$ must have the same nilpotent step. Thus, for instance, $\Z^4\wr\text{BS}(1,n)$ and $\text{Heis}(\Z)\wr\text{BS}(1,n)$ are not quasi-isometric, even though $\Z^4$ and $\text{Heis}(\Z)$ have equal growth degrees. Restricting to virtually abelian groups, we can completely describe the quasi-isometry class of our wreath products. 

\begin{theoremletter}\label{thm:classificationforvirtuallyabeliangroups}
Let $A_{1}$ and $A_{2}$ be infinite virtually abelian finitely generated groups, with growth degrees $d_{1}$ and $d_{2}$ respectively. Let $G$ and $H$ be finitely presented amenable groups from $\mathcal{M}_{\text{exp}}$. Then:
\begin{enumerate}[label=(\roman*)]
    \item $A_{1}\wr G$ and $A_{2}\wr H$ are quasi-isometric if and only if there exists a quasi-$\frac{d_{2}}{d_{1}}$-to-one quasi-isometry $G\rightarrow H$;
    \item $A_{1}\wr G$ and $A_{2}\wr H$ are biLipschitz equivalent if and only if $d_{1}=d_{2}$ and there exists a biLipschitz equivalence $G\rightarrow H$.
\end{enumerate}
\end{theoremletter}

We thus get additional examples of pairs of quasi-isometric groups that are not biLipschitz equivalent: for instance, as $\text{Sc}(\text{BS}(1,k))=\R_{>0}$ (see Proposition~\ref{prop:examplesofscalinggroups}), $\Z^{n}\wr\text{BS}(1,k)$ and $\Z^{m}\wr\text{BS}(1,k)$ are quasi-isometric for all $n,m\ge 1$ and $k\ge 2$, while they are biLipschitz equivalent if and only if $n=m$. 

In the opposite direction, as mentioned above, Genevois and Tessera showed in~\cite{GT24b} that a lamplighter group $F\wr K$, where $F$ is finite and $K$ is amenable, one-ended and finitely presented, has all its self-quasi-isometries at bounded distance from bijections. As a consequence of Theorem~\ref{thm:maintheoremfortheclassMexp}, the same phenomenon occurs for lamplighters with lamps of polynomial growth:

\begin{corollaryletter}\label{cor:scalinggroupsdesired}
Let $N$ and $M$ be finitely generated groups with polynomial growth. Let $G$ and $H$ be finitely presented amenable groups from $\mathcal{M}_{\text{exp}}$. If $N$ and $M$ have equal growth degrees, then any quasi-isometry $N\wr G\longrightarrow M\wr H$ lies within bounded distance from a bijection. 

In particular, any quasi-isometry $N\wr G\longrightarrow N\wr G$ lies at bounded distance from a bijection, and thus $\text{Sc}(N\wr G)=\lbrace 1\rbrace$.
\end{corollaryletter}

As explained in Section~\ref{sec:Scalinggroups}, the fact that $N\wr G$ has a trivial scaling group allows us to deduce some algebraic facts on its subgroups of finite index. For instance, two biLipschitz equivalent finite-index subgroups of $N\wr G$ must have the same index (Corollary~\ref{cor:biLipsubgroupshavesameindex}), and, in particular, $N\wr G$ does not have any proper finite-index subgroup isomorphic to itself (Corollary~\ref{cor:nopropersubgroupisotothegroup}), a fact that is not obvious to prove from a purely algebraic point of view. In turn, this has the consequence that any non-surjective monomorphism $N\wr G \hookrightarrow N\wr G$ has an image of infinite index (such morphisms exist as soon as $N$ or $G$ is not co-Hopfian, e.g. $N=\Z^{n}$ or $G=\text{BS}(1,n)$ for $n\ge 1$; see~\cite[Theorem~6.1]{BFF24}). 

\paragraph{Application to lamplighter-rigidity.} So far, related to the quasi-isometric classification of wreath products of the form $(\text{finite})\wr(\text{finitely generated})$, all results recorded in the literature exhibit a flexibility part, namely given a finitely generated group $H$ and integers $n,m\ge 2$ satisfying some arithmetic condition, it is often possible to construct a quasi-isometry between $\Z_{n}\wr H$ and $\Z_{m}\wr H$. This observation then suggests the next question:
\begin{question}\label{question1.5}
Does there exist a finitely generated group $H$ such that $\Z_{n}\wr H$ and $\Z_{m}\wr H$ are never quasi-isometric if $n\neq m$?
\end{question}

We know the answer when $H$ belongs to various classes of groups: for instance, if $H$ is non-amenable, $\Z_{n}\wr H$ and $\Z_{m}\wr H$ are quasi-isometric as soon as $n$ and $m$ have the same prime divisors, by Theorem~\ref{thm:constructingQIbetweenNAlamplightersGT21}. Hence an example of a group $H$ must necessarily come from the class of amenable groups.

But, further, if $H$ is amenable, it has finitely many ends, and thus it has either zero, one or two ends. The zero end case is trivial, as then both $\Z_{n}\wr H$ and $\Z_{m}\wr H$ are finite and thus quasi-isometric, for any choice of $n$ and $m$. The two-ended case is also completely understood by the work of Eskin-Fisher-Whyte~\cite{EFW12, EFW13}, and in this case it suffices to take $n$ and $m$ to be powers of a common number to construct a quasi-isometry $\Z_{n}\wr H\longrightarrow \Z_{m}\wr H$. Therefore, for a positive answer to Question~\ref{question1.5}, one has to look into amenable one-ended groups. In this class, the case of finitely presented groups is also completely understood, and here as well it suffices to take $n=a^{s}$ and $m=a^{r}$ powers of a common number to deduce that $\Z_{n}\wr H$ and $\Z_{m}\wr H$ are quasi-isometric if $\frac{s}{r}\in\text{Sc}(H)$~\cite[Theorem~3.12]{GT24b}. Hence, within this class, an example of a group $H$ answering positively Question~\ref{question1.5} should have trivial scaling group. However, we do not know a single example of a finitely presented one-ended group with trivial scaling group. 

On the other hand, despite the fact that our wreath products $N\wr G$ are not finitely presented, they do have trivial scaling group, and we show that they indeed provide examples of groups answering positively Question~\ref{question1.5}:

\begin{propositionletter}\label{prop:lamplighterrigidity}
Let $n,m\ge 2$ be two integers. Let $N$ be a finitely generated group with polynomial growth, and let $G$ be a finitely presented amenable group from $\mathcal{M}_{\text{exp}}$. Then the lamplighters $\Z_{n}\wr(N\wr G)$ and $\Z_{m}\wr(N\wr G)$ are quasi-isometric if and only if $n=m$. 
\end{propositionletter}

More generally, we know from~\cite{GT24b} that the existence of a quasi-isometry \mbox{$\Z_{n}\wr G\longrightarrow \Z_{m}\wr H$} provides a measure-scaling quasi-isometry $G\rightarrow H$, where the scaling factor depends on $n$ and $m$. Thus, when $G$ and $H$ are themselves lamplighters with lamps of polynomial growth, the existence of such a quasi-isometry imposes a compatibility condition between $n$, $m$ and the growth degrees of lamp groups. This is the content of the next statement.

\begin{propositionletter}\label{prop:mixingofscalingconditions}
Let $n,m\ge 2$ be two integers. Let $N_{1}$ and $N_{2}$ be finitely generated groups of polynomial growth, with growth degrees $n_{1}$ and $n_{2}$ respectively. Let $G$ and $H$ be finitely presented amenable groups from $\mathcal{M}_{\text{exp}}$. If \;$\Z_{n}\wr(N_{1}\wr G)$ and $\Z_{m}\wr(N_{2}\wr H)$ are quasi-isometric, then there exist $a,r,s\ge 1$ such that $n=a^{r}$, $m=a^{s}$, and $\frac{s}{r}=\frac{n_{2}}{n_{1}}$. 
\end{propositionletter}

\paragraph{Iterated wreath products.} In the light of Proposition~\ref{prop:mixingofscalingconditions}, it is also natural to try to fully classify iterated wreath products. For such groups, most of the technology developed in~\cite{GT24a, GT24b} does not apply, because lamplighter groups $F\wr K$ with $F$ finite do not have the~\textit{thick bigon property}, a key quasi-isometry invariant introduced in~\cite{GT24a}. Roughly speaking, a finitely generated group $G$ has the thick bigon property if given any two points $x,y\in G$ connected by a path $\gamma_{1}$ and any point $p\in \gamma_{1}$ far enough from $x$ and $y$, it is always possible to connect $x$ and $y$ by a path $\gamma_{2}$ which is coarsely homotopic to $\gamma_{1}$ and that avoids arbitrary large balls centered at $p$ (see Definition~\ref{def:TBP} for a more precise definition). It turns out that, as soon as $F$ is infinite, wreath products $F\wr K$ have this property~\cite[Proposition~1.3]{GT24a}, and building on Theorem~\ref{thm:classificationforvirtuallyabeliangroups}, we get a complete classification of wreath products of the form 
\begin{equation*}
    (\text{finite})\wr((\text{infinite virtually abelian})\wr(\text{finitely presented from $\mathcal{M}_{\text{exp}}$})).
\end{equation*}

\begin{corollaryletter}\label{cor:classificationofiteratedwreathproducts}
Let $n,m\ge 2$ be two integers. Let $A_{1}$ and $A_{2}$ be infinite virtually abelian finitely generated groups, with growth degrees $d_{1}$ and $d_{2}$ respectively. Let $G$ and $H$ be finitely presented amenable groups from $\mathcal{M}_{\text{exp}}$. Then:
\begin{enumerate}[label=(\roman*)]
    \item $\Z_{n}\wr(A_{1}\wr G)$ and $\Z_{m}\wr(A_{2}\wr H)$ are quasi-isometric if and only if there exist $a,r,s\ge 1$ such that $n=a^{r}$, $m=a^{s}$, $\frac{s}{r}=\frac{d_{2}}{d_{1}}$, and there exists a quasi-$\frac{d_{2}}{d_{1}}$-to-one quasi-isometry \;$G\rightarrow H$;
    \item $\Z_{n}\wr(A_{1}\wr G)$ and $\Z_{m}\wr(A_{2}\wr H)$ are biLipschitz equivalent if and only if $n=m$, $d_{1}=d_{2}$ and there exists a biLipschitz equivalence \;$G\rightarrow H$.
\end{enumerate}
\end{corollaryletter}

For instance, this corollary rules out the existence of a quasi-isometry 
\begin{equation*}
\Z_{2}\wr(\Z^{2}\wr \text{BS}(1,n)) \longrightarrow \Z_{4}\wr(\Z^{3}\wr \text{BS}(1,n))
\end{equation*}
where $n\ge 2$.

It is worth noticing that, if the amenability assumption on $G$ and $H$ is removed, a classification can also be easily deduced from \cite{BGT24} and Proposition~\ref{prop:constructionsofQI} below: $\Z_{n}\wr(A_{1}\wr G)$ and $\Z_{m}\wr(A_{2}\wr H)$ are quasi-isometric if and only if $n$ and $m$ have the same prime divisors and $G$ and $H$ are quasi-isometric.

\paragraph{More wreath products.} Observe that, when comparing two wreath products over the same bases, many quasi-isometry invariants (such as the volume growth, the number of ends, F\o lner functions, divergence, asymptotic dimension to name a few) fail to distinguish them. In general, much more involved methods, as the ones developed in~\cite{EFW12, EFW13, GT24a, GT24b}, appear necessary to exhibit a significant geometric difference between our groups. Those strategies require additional assumptions, relying on notions such as the thick bigon property in~\cite{GT24a} or coarse separation in~\cite{BGT24}. We conclude the article by noticing that our construction of quasi-isometries, combined with results of Erschler~\cite{Dyu00, Ers03}, also provides classification results in situations where techniques from~\cite{GT24a, BGT24} cannot be applied.  

\begin{corollaryletter}\label{cor:classificationforvirtuallyabeliangroups2}
Let $A_{1},A_{2},B_{1},B_{2}$ be infinite virtually abelian finitely generated groups. Then $A_{1}\wr B_{1}$ and $A_{2}\wr B_{2}$ are quasi-isometric if and only if $B_{1}$ and $B_{2}$ have same growth degrees. 
\end{corollaryletter}

For instance, $\Z^d\wr\Z^{k}$ and $\Z^{d'}\wr\Z^{k'}$ are quasi-isometric if and only if $k=k'$.

\paragraph{Plan of the chapter.} Section~\ref{subsection4.2} introduces notations and the necessary background on the class $\mathcal{M}_{\text{exp}}$, and Section~\ref{subsection4.3} the necessary results on aptolic quasi-isometries between wreath products. Section~\ref{subsection4.4} is dedicated to the proof of our main result. Section~\ref{subsection4.5} presents the various mentioned applications, and Section~\ref{subsection4.6} records several questions related to the article. 

\section{Notations and preliminaries}\label{subsection4.2}

\subsection{Notations}\label{subsubsection4.2.1} 
For an integer $n\ge 1$, $\Z_{n}$ denotes the cyclic group of order $n$. Given a metric space $(X,d_{X})$, we denote 
\begin{equation*}
    B_{X}(x,r) \defeq \lbrace y\in X : d_{X}(x,y)\le r\rbrace
\end{equation*}
the closed ball centered at $x\in X$ of radius $r>0$. 

Any graph in this text is unoriented and simplicial, that is without loops nor multiple edges. For such a graph $\Gamma$, $V(\Gamma)$ refers to its set of vertices, while $E(\Gamma)$ stands for the set of edges. Given any subset $A\subset\Gamma$, its boundary in $\Gamma$ is denoted $\partial_{\Gamma}A$ and is defined as
\begin{equation*}
    \partial_{\Gamma}A \defeq \left\lbrace y\in \Gamma\setminus A : \exists x\in A, (x,y)\in E(\Gamma)\right\rbrace.
\end{equation*}

\subsection{Coarse separation} We now define coarse separation of spaces, as in~\cite{BGT24}. For this, recall that a metric space $(X,d_{X})$ is \textit{$k-$coarsely connected}, with $k>0$, if for any $x,y\in X$, there is a sequence of points $x=x_{0},x_{1},\dots,x_{n-1},x_{n}=y$ such that $d_{X}(x_{i-1},x_{i})\le k$ for any $i=1,\dots, n$. 

\begin{definition}
Let $(X,d_{X})$ be a metric space, and let $\mathcal{W}$ be a family of subsets of $X$. We say that $\mathcal{W}$ \textit{coarsely separates} $X$ if there exist $k>0$ and $L\ge 0$ such that for any $D\ge 0$, there exists $W\in\mathcal{W}$ such that $X\setminus W^{+L}$ has at least two $k-$coarsely connected components with points at distance $\ge D$ from $W$.
\end{definition}

Given a family $\mathcal{W}$ of subsets of $X$, we define its \textit{growth function} as 
\begin{equation*}
    V_{\mathcal{W}}(r)=\sup_{W\in\mathcal{W}, \;w\in W}|W\cap B(w,r)|, \; r\ge 0
\end{equation*}
and we say that $\mathcal{W}$ has \textit{subexponential growth} if $\limsup_{r\rightarrow\infty}\frac{1}{r}\ln(V_{\mathcal{W}}(r))=0$. We denote $\mathcal{M}_{\text{exp}}$ the class of bounded degree graphs such that any coarsely separating family of subsets must have exponential growth. Note that this class is invariant under quasi-isometries.

The following theorem is one of the main results of~\cite{BGT24}. 

\begin{theorem}[{\cite[Theorem~1.4]{BGT24}}]
For $i=1,2$, let $X_{i}$ be either a $(k+1)-$regular tree with $k\ge 2$ or a rank one symmetric space of non-compact type. Then the horocyclic product $X_{1}\bowtie X_{2}$ belongs to $\mathcal{M}_{\text{exp}}$, as well as any direct product of $X_{1}\bowtie X_{2}$ with an arbitrary bounded degree graph. 
\end{theorem}

As already mentioned, this theorem provides many examples of amenable groups that belong to $\mathcal{M}_{\text{exp}}$, among which solvable Baumslag-Solitar groups $\text{BS}(1,n)$, $n\ge 2$, lamplighter groups and $\text{SOL}(\Z)$. 

\section{Aptolic quasi-isometries between wreath products}\label{subsection4.3}

As explained earlier, a quasi-isometry between wreath products is \textit{aptolic} if it preserves the lamplighter structure. In~\cite{GT24b}, Genevois and Tessera proved that when $N$ and $M$ are finite and $G$ and $H$ are both finitely presented and one-ended, then any quasi-isometry $N\wr G\longrightarrow M\wr H$ is at bounded distance from an aptolic quasi-isometry. Later, with Bensaid, they proved that the same phenomenon occurs when $N$ and $M$ are infinite and of subexponential growth, provided that the base groups $G$ and $H$ belong to $\mathcal{M}_{\text{exp}}$: 

\begin{theorem}[{\cite[Theorem~6.30]{BGT24}}]\label{thm:finitedistancefromaptolicmapsBGT24}
Let $A_{1},A_{2}$ be two finitely generated groups of subexponential growth, and let $B_{1},B_{2}$ be two finitely presented groups from $\mathcal{M}_{\text{exp}}$. Then any quasi-isometry $A_{1}\wr B_{1}\longrightarrow A_{2}\wr B_{2}$ is within bounded distance from an aptolic quasi-isometry. 
\end{theorem}

Similarly to the strategy for lamplighters over finitely presented one-ended groups from~\cite{GT24b}, it is proved in~\cite{BGT24} that any quasi-isometry $A_{1}\wr B_{1}\longrightarrow A_{2}\wr B_{2}$, with the assumptions of Theorem~\ref{thm:finitedistancefromaptolicmapsBGT24}, has the property of being \textit{leaf-preserving}, i.e. it sends any $B_{1}-$coset of $A_{1}\wr B_{1}$ at a uniform bounded distance from a $B_{2}-$coset of $A_{2}\wr B_{2}$ and has a quasi-inverse that does the same with roles of $B_{1}$ and $B_{2}$ reversed. This fact is deduced from another version of the embedding theorem proved in~\cite[Theorem~6.1]{BGT24} and that also relies on quasi-median geometry. The second step of the strategy, namely the fact that leaf-preservingness implies aptolicity (up to finite distance), is mentioned 
in~\cite[Theorem~6.28]{BGT24} and actually follows from~\cite[Corollary~6.11 and Lemma~6.12]{GT24a}.

In~\cite[Proposition~3.1]{GT24b} is proved a characterisation of aptolic quasi-isometries between wreath products with~\textit{finite} lamp groups. We will need in our computations a similar characterisation for wreath products with~\textit{arbitrary} lamp groups. 

\begin{proposition}\label{prop:characterisationofaptolicity}
Let $M,N,G,H$ be finitely generated groups, and consider two maps $\alpha\colon N^{(G)}\longrightarrow M^{(H)}$ and $\beta\colon G\longrightarrow H$. Then the map 
\begin{equation*}
    \varphi\colon N\wr G\longrightarrow M\wr H, (c,p)\longmapsto (\alpha(c), \beta(p))
\end{equation*}
is an aptolic quasi-isometry if and only if the following conditions hold:
\begin{enumerate}[label=(\roman*)]
    \item $\alpha\colon N^{(G)}\longrightarrow M^{(H)}$ is a bijection;
    \item $\beta\colon G\rightarrow H$ is a quasi-isometry;
    \item There exists $Q\ge 0$ such that, for any colourings $c_{1},c_{2}\in N^{(G)}$, the Hausdorff distance between $\text{supp}\big(\alpha(c_{1})^{-1}\alpha(c_{2})\big)$ and $\beta\big(\text{supp}(c_{1}^{-1}c_{2})\big)$ is at most $Q$;
    \item There exists $L\ge 0$ such that, for any colourings $c_{1},c_{2}\in N^{(G)}$ that differ on a single point $p\in G$ and such that $c_{1}(p),c_{2}(p)$ are adjacent in $N$, one has
    \begin{equation*}
        d_{M}\big(\alpha(c_{1})(t), \alpha(c_{2})(t)\big) \le L
    \end{equation*} 
    for all $t\in H$;
    \item There exists $L'\ge 0$ such that, for any colourings $c_{1},c_{2}\in M^{(H)}$ that differ on a single point $p\in H$ and such that $c_{1}(p),c_{2}(p)$ are adjacent in $M$, one has
    \begin{equation*}
        d_{N}\big(\alpha^{-1}(c_{1})(t), \alpha^{-1}(c_{2})(t)\big) \le L'
    \end{equation*}
    for all $t\in G$.
\end{enumerate}
In particular, every aptolic quasi-inverse of $\varphi$ is of the form 
\begin{equation*}
    (c,p)\longmapsto (\alpha^{-1}(c), \overline{\beta}(p)),\; (c,p)\in M\wr H
\end{equation*}
where $\overline{\beta}\colon H\rightarrow G$ is a quasi-inverse of $\beta\colon G\rightarrow H$. 
\end{proposition}

\begin{proof}
Suppose first that $\varphi\colon N\wr G\longrightarrow M\wr H$ is a $(C,K)-$quasi-isometry.
The proof of~\textit{(i)},~\textit{(ii)},~\textit{(iii)} and of the last claim of the statement follow the same lines as the proof of Proposition~\ref{prop:characterisationofaptolicQI:finitecase}.  

\noindent Let us then turn to \textit{(iv)}. Let $c_{1},c_{2}\in N^{(G)}$ be two colourings of $G$ that differ only on $p\in G$, and that take adjacent values on that point. This means that $(c_{1},p)$ and $(c_{2},p)$ are neighbours in $N\wr G$, and thus 
\begin{align*}
    d_{M\wr H}\big((\alpha(c_{1}),\beta(p)), (\alpha(c_{2}),\beta(p))\big)&=d_{M\wr H}(\varphi(c_{1},p), \varphi(c_{2},p)) \\
    &\le C\cdot d_{N\wr G}((c_{1},p), (c_{2},p))+K \\
    &=C+K.
\end{align*}
This distance being less than $C+K$ implies that $\alpha(c_{1})$ and $\alpha(c_{2})$ can only differ on points at distance $\le C+K$ from $\beta(p)$ and, at each point $t\in H$ where they effectively differ, one can go in $M$ from $\alpha(c_{1})(t)$ to $\alpha(c_{2})(t)$ in at most $C+K$ steps. Thus
\begin{equation*}
    d_{M}\big(\alpha(c_{1})(t), \alpha(c_{2})(t)\big)\le C+K
\end{equation*}
for any $t\in H$. 

\noindent The proof of~\textit{(v)} is analogous, using a quasi-inverse $\overline{\varphi}$ of $\varphi$ that takes the form 
\begin{equation*}
    (c,p)\longmapsto (\alpha^{-1}(c), \overline{\beta}(p)),\; (c,p)\in M\wr H
\end{equation*}
where $\overline{\beta}\colon H\rightarrow G$ is a quasi-inverse of $\beta$. 

\noindent Conversely, assume that $\textit{(i)}-\textit{(v)}$ hold. Fix constants $C\ge 1,K\ge 0$ such that $\beta\colon G\rightarrow H$ and a quasi-inverse $\overline{\beta}\colon H\rightarrow G$ are $(C,K)-$quasi-isometries, and such that $\overline{\beta}\circ\beta$, $\beta\circ\overline{\beta}$ are at distance $\le K$ from $\text{Id}_{G}$, $\text{Id}_{H}$ respectively. 

\noindent We start by showing that $\varphi$ is Lipschitz. Let $a=(c_{1},p_{1})$, $b=(c_{2},p_{2})$ be two adjacent vertices in $N\wr G$. We distinguish two cases:
\begin{itemize}
    \item First, say that $c_{1}=c_{2}$ and that $p_{1}$ and $p_{2}$ are neighbours in $G$. Then we have 
    \begin{align*}
        d_{M\wr H}(\varphi(a),\varphi(b))&=d_{M\wr H}\big((\alpha(c_{1}),\beta(p_{1})), (\alpha(c_{1}),\beta(p_{2}))\big) \\
        &=d_{H}(\beta(p_{1}), \beta(p_{2})) \\
        &\le C\cdot d_{G}(p_{1},p_{2})+K \\
        &=C+K.
    \end{align*}
    \item Now, assume that $p_{1}=p_{2}$ and that $c_{1},c_{2}$ differ only on $p_{1}$, and that $c_{1}(p_{1})$ and $c_{2}(p_{1})$ are adjacent in $N$. From~\textit{(iii)}, we know that
    \begin{equation*}
        \text{supp}\big(\alpha(c_{1})^{-1}\alpha(c_{2})\big)\subset B_{H}(\beta(p_{1}),Q)
    \end{equation*}
    and, from~\textit{(iv)}, we know that for any $t\in \text{supp}\big(\alpha(c_{1})^{-1}\alpha(c_{2})\big)$, there is a path of length $\le L$ in $M$ from $\alpha(c_{1})(t)$ to $\alpha(c_{2})(t)$. This yields the estimate
    \begin{align*}
        d_{M\wr H}(\varphi(a),\varphi(b))&=d_{M\wr H}\big((\alpha(c_{1}),\beta(p_{1})), (\alpha(c_{2}),\beta(p_{1}))\big) \\
        &\le L\cdot \left|B_{H}(\beta(p_{1}),Q)\right| \\
        &\le L\cdot D^{Q}
    \end{align*}
    where $D\ge 3$ is a fixed integer larger than the degree of a vertex in $\text{Cay}(H,T)$.
\end{itemize}
In all cases, we conclude that 
\begin{equation*}
    d_{M\wr H}(\varphi(a),\varphi(b)) \le \max(C+K, L\cdot D^{Q})
\end{equation*}
so that $\varphi$ is \mbox{$\max(C+K, L\cdot D^{Q})-$}Lipschitz according to Lemma~\ref{lem:Lipschitzbetweengraphs}. 

\noindent Now, consider the map 
\begin{align*}
    \psi \colon M\wr H &\longrightarrow N\wr G \\
    (c,p)&\longmapsto (\alpha^{-1}(c), \overline{\beta}(p)). 
\end{align*}
For any $(c,p)\in N\wr G$, one computes that 
\begin{align*}
    d_{N\wr G}\left(\psi\circ\varphi(c,p), (c,p)\right)&=d_{N\wr G}\left((c, \overline{\beta}(\beta(p))), (c,p)\right) \\
    &=d_{G}\left(\overline{\beta}\circ\beta(p), p\right) \\
    &\le K
\end{align*}
so $\psi\circ\varphi$ is at distance $\le K$ from $\text{Id}_{N\wr G}$, and similarly, using that $\beta\circ\overline{\beta}$ is at distance $\le K$ from $\text{Id}_{H}$, one shows that $\varphi\circ \psi$ is at distance $\le K$ from $\text{Id}_{M\wr H}$. Thus it only remains to prove that $\psi$ is Lipschitz. For this, note that:

\begin{claim}\label{claim3.4}
There exists $Q'\ge 0$ such that, for any $c_{1},c_{2}\in M^{(H)}$, the Hausdorff distance between $\text{supp}\big(\alpha^{-1}(c_{1})^{-1}\alpha^{-1}(c_{2})\big)$ and $\overline{\beta}\big(\text{supp}(c_{1}^{-1}c_{2})\big)$ is at most $Q'$.
\end{claim}

{
\renewcommand{\proofname}{Proof of Claim \ref{claim3.4}.}
\renewcommand{\qedsymbol}{$\blacksquare$}
\begin{proof}
Fix $c_{1},c_{2}\in M^{(H)}$. Applying~\textit{(iii)} to $\alpha^{-1}(c_{1})$ and $\alpha^{-1}(c_{2})$, we have  
\begin{equation*}
    d_{\text{Haus}}\left(\text{supp}(c_{1}^{-1}c_{2}), \beta\left(\text{supp}(\alpha^{-1}(c_{1})^{-1}\alpha^{-1}(c_{2}))\right)\right) \le Q
\end{equation*}
and applying $\overline{\beta}$, it follows from Lemma~\ref{lem:neighborhoodsandQI}\textit{(iv)} that 
\begin{equation*}
    d_{\text{Haus}}\left(\overline{\beta}(\text{supp}(c_{1}^{-1}c_{2})), \overline{\beta}\left(\beta\left(\text{supp}(\alpha^{-1}(c_{1})^{-1}\alpha^{-1}(c_{2}))\right)\right)\right) \le C\cdot Q+K.
\end{equation*}
As the Hausdorff distance between $\overline{\beta}\big(\beta\big(\text{supp}(\alpha^{-1}(c_{1})^{-1}\alpha^{-1}(c_{2}))\big)\big)$ and $\text{supp}\big(\alpha^{-1}(c_{1})^{-1}\alpha^{-1}(c_{2})\big)$ is bounded by $K$, we get the claim with $Q'\defeq C\cdot Q+2K$. 
\end{proof}
\renewcommand{\qedsymbol}{$\square$}
}
\noindent Using Claim~\ref{claim3.4} and~\textit{(v)}, computations as the ones above show that $\psi$ sends any two adjacent vertices of $M\wr H$ to vertices that are at distance at most 
\begin{equation*}
    \max\left(C+K, L'\cdot D'^{Q'}\right)
\end{equation*}
in $N\wr G$, where $D'\ge 3$  is a fixed integer larger than the degree of a vertex in $\text{Cay}(G,S)$. Invoking once again Lemma~\ref{lem:Lipschitzbetweengraphs}, we conclude that $\psi$ is $\max(C+K, L'\cdot D'^{Q'})-$Lipschitz. Finally, we get that 
\begin{align*}
d_{M\wr H}(\varphi(a),\varphi(b))&\ge \frac{1}{\max(C+K, L'\cdot D'^{Q'})}d_{N\wr G}\big(\psi(\varphi(a)), \psi(\varphi(b))\big) \\
&\ge \frac{1}{\max(C+K, L'\cdot D'^{Q'})}\cdot d_{N\wr G}(a,b)-\frac{2K}{\max(C+K, L'\cdot D'^{Q'})} 
\end{align*}
for any $a,b\in N\wr G$, which concludes the proof that $\varphi\colon N\wr G\longrightarrow M\wr H$ is a quasi-isometry, with $\psi\colon M\wr H\longrightarrow N\wr G$ as a quasi-inverse. 
\end{proof}

Here is then the key statement that will provide us with the rigidity we are looking for. 

\begin{proposition}\label{prop:finiteunionofcosets}
Let $N,M,G,H$ be finitely generated groups. Let $\varphi\colon N\wr G\longrightarrow M\wr H$ be an aptolic quasi-isometry, i.e. there are two maps $\alpha\colon N^{(G)}\longrightarrow M^{(H)}$ and $\beta\colon G\rightarrow H$ such that $\varphi(c,p)=(\alpha(c),\beta(p))$ for any $(c,p)\in N\wr G$. 

For any quasi-inverse $\overline{\beta}$ of $\beta$, there exists a constant $Q\ge 0$ such that, for any finite subset $A\subset G$ and any $Q'\ge Q$, $\alpha^{-1}\left(\mathcal{L}(\beta(A)^{+Q'})\right)$ is a union of cosets of $\mathcal{L}(A)$; and conversely, for any finite subset $B\subset H$ and any $Q'\ge Q$, $\alpha\left(\mathcal{L}(\overline{\beta}(B)^{+Q'})\right)$ is a union of cosets of $\mathcal{L}(B)$.
\end{proposition}

Recall that, for a subset $A\subset G$, $\mathcal{L}(A)$ stands for the subgroup of $N^{(G)}$ of colourings supported on $A$. 

The proof of this observation requires two intermediate lemmas.

\begin{lemma}\label{lm:finiteunionofcosets}
Let $N,M,G,H$ be finitely generated groups. Consider two maps $\alpha\colon N^{(G)}\longrightarrow M^{(H)}$ and $\beta\colon G\rightarrow H$. Suppose that there exists $Q\ge 0$ such that, for any $c_{1},c_{2}\in N^{(G)}$ with $\text{supp}(c_{1}^{-1}c_{2})\subset \lbrace p\rbrace$ for some $p\in G$, we have $\text{supp}\big(\alpha(c_{1})^{-1}\alpha(c_{2})\big) \subset B_{H}(\beta(p),Q)$. 

Then, for any finite subset $A\subset G$ and any colouring $c\in N^{(G)}$, we have 
\begin{equation*}
    \alpha\big(c\mathcal{L}(A)\big)\subset \alpha(c)\mathcal{L}\big(\beta(A)^{+Q}\big).
\end{equation*}
\end{lemma}

\begin{proof}
Let $c\in N^{(G)}$. We prove the claim by induction on $|A|$. Suppose to start that $|A|=1$ and $c\in N^{(G)}$, and let $c'\in c\mathcal{L}(A)$. Then $c$ and $c'$ differ on at most one point $p$ of $G$, so our assumption implies that $\alpha(c)$ and $\alpha(c')$ can only differ on points from
\begin{equation*}
    B_{H}(\beta(p),Q)=\beta(\lbrace p\rbrace)^{+Q}=\beta(A)^{+Q}.
\end{equation*}
This means that $\alpha(c')\in \alpha(c)\mathcal{L}\big(\beta(A)^{+Q}\big)$, so the claim holds for subsets reduced to a point. 

\noindent Suppose now that it holds for any colouring and for subsets of cardinality $k\ge 1$, and let $A\subset G$ be a set of cardinality $k+1$. Let $c\in N^{(G)}$ and let $c'\in c\mathcal{L}(A)$. Hence $c^{-1}c'\in\mathcal{L}(A)$, and we can write $c^{-1}c'=c''d$ for some colourings $c''\in \mathcal{L}(A\setminus\lbrace a\rbrace)$, $d\in\mathcal{L}(\lbrace a\rbrace)$, where $a\in A$ is an arbitrary point. Thus $c'=(cc'')d$, and it follows that 
\begin{equation*}
    \alpha(c')\in \alpha(cc'')\mathcal{L}\big(\beta(\lbrace a\rbrace)^{+Q}\big).
\end{equation*}
On the other hand, the inductive assumption applied to $A\setminus\lbrace a\rbrace$ shows that \mbox{$\alpha(cc'')\in\alpha(c)\mathcal{L}\big(\beta(A\setminus\lbrace a\rbrace)^{+Q}\big)$}, whence 
\begin{equation*}
    \alpha(c')\in \alpha(cc'')\mathcal{L}\big(\beta(\lbrace a\rbrace)^{+Q}\big) \subset \alpha(c)\mathcal{L}\big(\beta(A\setminus\lbrace a\rbrace)^{+Q}\big)\mathcal{L}\big(\beta(\lbrace a\rbrace)^{+Q}\big)=\alpha(c)\mathcal{L}\big(\beta(A)^{+Q}\big).
\end{equation*}
Thus we conclude that $\alpha\big(c\mathcal{L}(A)\big)\subset \alpha(c)\mathcal{L}\big(\beta(A)^{+Q}\big)$, and the induction is complete. 
\end{proof}

{
\renewcommand{\proofname}{Proof of Proposition~\ref{prop:finiteunionofcosets}.}

\begin{proof} 
By Proposition~\ref{prop:characterisationofaptolicity}, there exists $Q\ge 0$ such that the Hausdorff distance between $\beta\left(\text{supp}(c_{1}^{-1}c_{2})\right)$ and $\text{supp}\big(\alpha(c_{1})^{-1}\alpha(c_{2})\big)$ is at most $Q$, for any $c_{1},c_{2}\in N^{(G)}$. In particular, the assumption of Lemma~\ref{lm:finiteunionofcosets} is satisfied. 

\noindent Fix now any finite subset $A\subset G$ and a number $Q'\ge Q$. Let $c'\in \alpha^{-1}\left(\mathcal{L}(\beta(A)^{+Q'})\right)$, and let $d\in \mathcal{L}(A)$. Then, applying Lemma~\ref{lm:finiteunionofcosets}, one has 
\begin{equation*}
    \alpha(c'd) \in \alpha(c')\mathcal{L}\big(\beta(A)^{+Q}\big) \subset \alpha(c')\mathcal{L}\big(\beta(A)^{+Q'}\big)
\end{equation*}
and since $\alpha(c')\in \mathcal{L}\big(\beta(A)^{+Q'}\big)$, we deduce that $\alpha(c'd)\in \mathcal{L}(\beta(A)^{+Q'})$, i.e. 
\begin{equation*}
    c'd\in \alpha^{-1}\left(\mathcal{L}(\beta(A)^{+Q'})\right).
\end{equation*}
Thus $\alpha^{-1}\left(\mathcal{L}(\beta(A)^{+Q'})\right)$ is stable under multiplication by elements of $\mathcal{L}(A)$, which implies that it is a union of cosets of $\mathcal{L}(A)$.
\end{proof}}

Let us conclude this section with a sufficient condition for aptolic quasi-isometries to be scaling.

\begin{lemma}\label{lm:scalingfactorsforwreathproducts}
Let $N,M,G,H$ be finitely generated groups. Consider two maps $\alpha\colon N^{(G)}\longrightarrow M^{(H)}$ and $\beta\colon G\rightarrow H$ such that the map
\begin{align*}
    \varphi\colon N\wr G&\longrightarrow M\wr H \\
    (c,p)&\longmapsto (\alpha(c),\beta(p))
\end{align*}
is a quasi-isometry. If $\beta$ is quasi-$k$-to-one for some $k>0$, then $\varphi$ is quasi-$k$-to-one. 
\end{lemma}

\begin{proof}
The proof of Lemma~\ref{lem:scalingQIbetweenhalos} applies word for word. 
\end{proof}

\section{Proof of Theorem~\ref{thm:maintheorem}}\label{subsection4.4}

We can now combine all these intermediate observations to prove our main result. Fix then \mbox{$\varphi\colon N\wr G\longrightarrow M\wr H$} an arbitrary quasi-isometry, where $N$ and $M$ have polynomial growth of degrees $n$ and $m$ respectively. Our goal is to prove that $\varphi$ is quasi-$\frac{m}{n}$-to-one.

\vspace{0.15cm}

\noindent Let $C\ge 1$, $K\ge 0$ be such that $\varphi$ and a quasi-inverse $\overline{\varphi}\colon M\wr H\longrightarrow N\wr G$ are $(C,K)-$quasi-isometries, and such that  $\overline{\varphi}\circ\varphi$, $\varphi\circ\overline{\varphi}$ are within distance $K$ from $\text{Id}_{N\wr G}$, $\text{Id}_{M\wr H}$ respectively. By assumption, up to a bounded distance, $\varphi$ can be chosen aptolic, i.e. there are two maps $\alpha\colon N^{(G)}\longrightarrow M^{(H)}$, $\beta\colon G\rightarrow H$ such that
\begin{equation*}
    \varphi(c,p)=(\alpha(c),\beta(p))
\end{equation*}
for any $(c,p)\in N\wr G$. It follows from Proposition~\ref{prop:characterisationofaptolicity} that $\alpha$ is a bijection, that $\beta$ is a $(C,K)-$quasi-isometry, and that $\overline{\varphi}$ takes the form
\begin{equation*}
    \overline{\varphi}(c,p)=(\alpha^{-1}(c), \overline{\beta}(p))
\end{equation*}
for any $(c,p)\in M\wr H$, where $\overline{\beta}\colon H\rightarrow G$ is a quasi-inverse of $\beta$ (with same parameters $C$ and $K$). Lastly, by increasing $K$ if needed, we may assume that it is larger than the constant $Q\ge 0$ given by Proposition~\ref{prop:characterisationofaptolicity}\textit{(iii)}. 

\begin{claim}\label{claim4.1}
The equality $\left|A\right|=\frac{m}{n}\cdot \left|\beta(A)^{+K}\right|$ holds for any finite subset $A\subset G$. 
\end{claim}

{
\renewcommand{\proofname}{Proof of Claim~\ref{claim4.1}.}
\renewcommand{\qedsymbol}{$\blacksquare$}
\begin{proof}
Let $A\subset G$ be finite. Notice that, by Proposition~\ref{prop:finiteunionofcosets}, $\alpha^{-1}\left(\mathcal{L}(\beta(A)^{+2K})\right)$ is a union of cosets of $\mathcal{L}(A)$, say
\begin{equation*}
    \alpha^{-1}\left(\mathcal{L}(\beta(A)^{+2K})\right)=\bigsqcup_{i=1}^{k}d_{i}\mathcal{L}(A)
\end{equation*}
for some $k\ge 1$ and $d_{1},\dots, d_{k}\in N^{(G)}$. Consider then the subset 
\begin{equation*}
    S \defeq \left\lbrace (d,q) \in M\wr H : d\in \mathcal{L}(\beta(A)^{+K}), q\in \beta(A)^{+K}\right\rbrace. 
\end{equation*}
Notice that $S^{+K}\subset \left\lbrace (d,q)\in M\wr H : d\in \mathcal{L}(\beta(A)^{+2K}), q\in \beta(A)^{+2K}\right\rbrace$, so that 
\begin{align*}
    \varphi^{-1}(S^{+K}) &\subset \varphi^{-1}\left(\left\lbrace (d,q)\in M\wr H : d\in \mathcal{L}(\beta(A)^{+2K}), q\in \beta(A)^{+2K}\right\rbrace\right) \\
    &=\left\lbrace (c,p) \in N\wr G : c\in \alpha^{-1}\left(\mathcal{L}(\beta(A)^{+2K})\right), p\in \beta^{-1}\left(\beta(A)^{+2K}\right)\right\rbrace \\
    &=\left\lbrace (c,p) \in N\wr G : c\in \bigsqcup_{i=1}^{k}d_{i}\mathcal{L}(A), p\in \beta^{-1}\left(\beta(A)^{+2K}\right)\right\rbrace
\end{align*}
and the latter has growth degree $n|A|$. Hence $\varphi^{-1}(S^{+K})$ has growth degree at most $n|A|$. On the other hand, $\varphi^{-1}(S^{+K})$ coarsely coincides with $\overline{\varphi}(S)$ by Lemma~\ref{lm:preimagesandquasiinverses}, so their growth degrees coincide, and as moreover the growth degree is preserved by quasi-isometries (Proposition~\ref{prop:invarianceofgrowth}), it follows that $\varphi^{-1}(S^{+K})$ has the same growth degree as $S$, which is $m\left|\beta(A)^{+K}\right|$. We thus get the inequality
\begin{equation}\label{eq4.1}
    m\left|\beta(A)^{+K}\right| \le n\left|A\right|
\end{equation}
for any finite subset $A\subset G$. The same reasoning with $\overline{\varphi}$ shows that 
\begin{equation}\label{eq4.2}
    n\left|\overline{\beta}(B)^{+K}\right| \le m\left|B\right|
\end{equation}
for any finite subset $B\subset H$. Thus, fixing any finite subset $A\subset G$ and applying (\ref{eq4.2}) with $B=\beta(A)^{+K}$ and then (\ref{eq4.1}), it follows that 
\begin{equation*}
    n\left|\overline{\beta}\left(\beta(A)^{+K}\right)^{+K}\right| \le m\left|\beta(A)^{+K}\right| \le n\left|A\right|
\end{equation*}
Since $A\subset \overline{\beta}\left(\beta(A)^{+K}\right)^{+K}$ (any $x\in A$ is within distance $K$ from $\overline{\beta}(\beta(x))$ with $\beta(x)\in\beta(A)^{+K}$), we get that 
\begin{equation*}
    n\left|A\right| \le n\left|\overline{\beta}\left(\beta(A)^{+K}\right)^{+K}\right| \le m\left|\beta(A)^{+K}\right| \le n\left|A\right|
\end{equation*}
and finally $m\left|\beta(A)^{+K}\right|=n\left|A\right|$, i.e. $\left|A\right|=\frac{m}{n}\cdot \left|\beta(A)^{+K}\right|$ for any finite subset $A\subset G$. This proves Claim~\ref{claim4.1}.
\end{proof}
\renewcommand{\qedsymbol}{$\square$}
}

\begin{claim}\label{claim4.2}
The quasi-isometry $\beta\colon G\rightarrow H$ is quasi-$\frac{m}{n}$-to-one.
\end{claim}

{
\renewcommand{\proofname}{Proof of Claim~\ref{claim4.2}.}
\renewcommand{\qedsymbol}{$\blacksquare$}
\begin{proof}
Fix any finite subset $A\subset H$. Using Claim~\ref{claim4.1}, we estimate 
\begin{align*}
    \left|\frac{m}{n}|A|-|\beta^{-1}(A)|\right| &= \left|\frac{m}{n}|A|-\frac{m}{n}|\beta(\beta^{-1}(A))^{+K}|\right| \\
    &\le \left|\frac{m}{n}|A|-\frac{m}{n}|A^{+K}|\right|+\left|\frac{m}{n}|A^{+K}|-\frac{m}{n}|\beta(\beta^{-1}(A))^{+K}|\right| \\
    &\le \frac{m}{n}\left|A^{+K}\setminus A\right|+\frac{m}{n}\left|A^{+K}\setminus \beta(\beta^{-1}(A))^{+K}\right|.
\end{align*}
Now, we claim that $A^{+K}\setminus \beta(\beta^{-1}(A))^{+K} \subset (\partial_{H}A)^{+(2K-1)}$. Indeed, let $y\in A^{+K}\setminus \beta(\beta^{-1}(A))^{+K}$. We distinguish two cases. If $y\notin A$, then $y\in A^{+K}\setminus A \subset (\partial_{H}A)^{+(K-1)}\subset (\partial_{H}A)^{+(2K-1)}$ as wanted.

\noindent We may therefore assume that $y\in A$. We know that there is $x\in G$ with $d_{H}(y,\beta(x))\le K$. In particular, $\beta(x)\in A^{+K}$. On the other hand, since $y\notin  \beta(\beta^{-1}(A))^{+K}$, it follows that $x\notin \beta^{-1}(A)$, i.e. $\beta(x)\notin A$. Thus $\beta(x)\in A^{+K}\setminus A \subset (\partial_{H}A)^{+(K-1)}$. As $y$ is within distance $K$ from $\beta(x)$, it follows that $y$ is within distance $2K-1$ from $\partial_{H}A$, as wanted. 

\noindent In any case, we get $y\in (\partial_{H}A)^{+(2K-1)}$, which proves that
\begin{equation*}
    A^{+K}\setminus \beta(\beta^{-1}(A))^{+K} \subset (\partial_{H}A)^{+(2K-1)}.
\end{equation*}

\noindent Therefore, it follows from Lemma~\ref{lem:Boundingneighborhoodsingraphs}\textit{(i)} and~\textit{(ii)} that there exist constants $L>0$, $R>0$ (depending only on $K$ and $H$) such that \begin{equation*}
    \left|A^{+K}\setminus A\right| \le R\cdot \left|\partial_{H}A\right|
\end{equation*}
and 
\begin{equation*}
    \left|A^{+K}\setminus \beta(\beta^{-1}(A))^{+K}\right| \le \left|(\partial_{H}A)^{+(2K-1)}\right| \le L\cdot \left|\partial_{H}A\right|.
\end{equation*}
Hence we finish our estimation above as
\begin{equation*}
    \left|\frac{m}{n}|A|-|\beta^{-1}(A)|\right| \le \frac{m}{n}\cdot (L+R)\cdot |\partial_{H}A|
\end{equation*}
which proves Claim~\ref{claim4.2}.
\end{proof}
\renewcommand{\qedsymbol}{$\square$}
}

\noindent Now, the fact that $\varphi$ is quasi-$\frac{m}{n}$-to-one directly follows from Claim~\ref{claim4.2} and Lemma~\ref{lm:scalingfactorsforwreathproducts}, which proves Theorem~\ref{thm:maintheorem}.

\section{Applications}\label{subsection4.5}

We can now focus on the consequences mentioned in the introduction.

{
\renewcommand{\proofname}{Proof of Theorem~\ref{thm:maintheoremfortheclassMexp}}
\begin{proof}
By Theorem~\ref{thm:finitedistancefromaptolicmapsBGT24}, any quasi-isometry $N\wr G\longrightarrow M\wr H$ is at a bounded distance from an aptolic quasi-isometry. Thus Theorem~\ref{thm:maintheorem} applies and gives the conclusion.
\end{proof}}

{
\renewcommand{\proofname}{Proof of Corollary~\ref{cor:BLE}}
\begin{proof}
Assume that $N\wr G$ and $M\wr H$ are biLipschitz equivalent. Such a biLipschitz equivalence is quasi-one-to-one, and up to a finite distance change, it can be taken aptolic. This new quasi-isometry, that we denote 
\begin{align*}
    \varphi\colon N\wr G &\longrightarrow M\wr H \\
    (c,p)&\longmapsto (\alpha(c),\beta(p))
\end{align*}
is still quasi-one-to-one by Proposition~\ref{prop:stabilitypropertiesforscalingQI}(\textit{i}), and in addition it must be quasi-$\frac{m}{n}$-to-one by Theorem~\ref{thm:maintheoremfortheclassMexp}. From Lemma~\ref{lem:uniquenessscalingfactor}, it follows that $m=n$, and we also deduce that $\beta\colon G\rightarrow H$ is quasi-one-to-one. From Theorem~\ref{thm:Whytethm}, $\beta$ lies at finite distance from a bijection, which is the desired biLipschitz equivalence $G\rightarrow H$. 
\end{proof}}

The proofs of Theorem~\ref{thm:classificationforvirtuallyabeliangroups} and Corollary~\ref{cor:classificationforvirtuallyabeliangroups2} require, in the flexibility part, to be able to construct quasi-isometries between wreath products in general situations. This is the goal of the next criterion. 

\begin{proposition}\label{prop:constructionsofQI}
Let $N,M,G,H$ be finitely generated groups. Let $n,m\ge 2$. If there exist a biLipschitz equivalence $N^{m}\rightarrow M^{n}$ and a quasi-$\frac{m}{n}$-to-one quasi-isometry $G\rightarrow H$, then there exists an aptolic quasi-isometry
\begin{equation*}
    N\wr G\longrightarrow M\wr H. 
\end{equation*}
\end{proposition}

\begin{proof}
Let $\sigma\colon N^{m}\rightarrow M^{n}$ be a $C-$biLipschitz equivalence, and let $\beta\colon G\rightarrow H$ be quasi-$\frac{m}{n}$-to-one. From Theorem~\ref{thm:rationalscalingfactor}, there exist a partition $\mathcal{P}$ (resp. $\mathcal{Q}$) of $G$ (resp. of $H$) with uniformly bounded pieces of size $m$ (resp. of size $n$) and a bijection $\psi\colon \mathcal{P}\rightarrow \mathcal{Q}$ such that $\beta(P)\subset \psi(P)$ for all $P\in \mathcal{P}$. Up to postcomposing $\sigma$ with a translation by an element of $M^{n}$, we can assume that $\sigma(1_{N^{m}})=1_{M^{n}}$. We can now define a bijection $\alpha\colon N^{(G)}\longrightarrow M^{(H)}$ in such a way that $\alpha$ sends $\mathcal{L}(P)$ into $\mathcal{L}(\psi(P))$ through $\sigma$ for any $P\in\mathcal{P}$, and $\alpha^{-1}$ sends $\mathcal{L}(Q)$ into $\mathcal{L}(\psi^{-1}(Q))$ through $\sigma^{-1}$ for any $Q\in\mathcal{Q}$.
Set 
\begin{align*}
    f\colon N\wr G &\longrightarrow M\wr H \\
    (c,p) &\longmapsto (\alpha(c), \beta(p)).
\end{align*}
We show that $f$ is a quasi-isometry by checking the five points of Proposition~\ref{prop:characterisationofaptolicity}. Points \textit{(i)} and \textit{(ii)} are satisfied by construction. For \textit{(iv)}, fix two colourings $c_{1},c_{2}\in N^{(G)}$ that differ on a single point $p\in G$ and such that $d_{N}(c_{1}(p), c_{2}(p))=1$. Let $P\in \mathcal{P}$ be the piece containing $p$. Let $u\in N^{m}$ (resp. $v\in N^{m}$) be the vector formed by colours of $c_{1}$ (resp. of $c_{2}$) on the piece $P$. Then $d_{N^{m}}(u,v)=1$, and since $\sigma$ is $C-$Lipschitz, it follows that 
\begin{equation}\label{eq5.2}
    d_{M^{n}}(\sigma(u),\sigma(v))\le C.
\end{equation}
By construction, the components of $\sigma(u)\in M^{n}$ (resp. of $\sigma(v)\in M^{n}$) are the colors of $\alpha(c_{1})$ (resp. of $\alpha(c_{2})$) on the piece $\psi(P)$. Thus, from (\ref{eq5.2}), it follows that 
\begin{equation*}
    d_{M}\big(\alpha(c_{1})(t), \alpha(c_{2})(t)\big) \le C
\end{equation*}
for any $t\in\psi(P)$. Additionally, if $t\in H\setminus \psi(P)$, then it belongs to another piece $Q\in \mathcal{Q}$, that we may write $Q=\psi(P')$ for some $P'\neq P$. By assumption, $c_{1}$ and $c_{2}$ agree on $P'$, so $\alpha(c_{1})$ and $\alpha(c_{2})$ agree on $\psi(P')$, in particular on $t\in \psi(P')$. Thus we have proved 
\begin{equation*}
    d_{M}\big(\alpha(c_{1})(t), \alpha(c_{2})(t)\big) \le C
\end{equation*}
for any $t\in H$, which is exactly \textit{(iv)} of Proposition~\ref{prop:characterisationofaptolicity}. Point \textit{(v)} is checked in a similar manner.

\noindent Finally, let us focus on~\textit{(iii)}. Fix two colourings $c_{1},c_{2}\in N^{(G)}$, and let $P_{1},\dots, P_{k}\in\mathcal{P}$ be the pieces of $\mathcal{P}$ containing points of $\text{supp}(c_{1}^{-1}c_{2})$. Then $\psi(P_{1}),\dots,\psi(P_{k})$ are the pieces of $\mathcal{Q}$ containing points of $\beta\left(\text{supp}(c_{1}^{-1}c_{2})\right)$. Since those pieces are uniformly bounded, we deduce that there is a constant $D\ge 0$ such that
\begin{equation*}
    d_{\text{Haus}}\left(\text{supp}(\alpha(c_{1})^{-1}\alpha(c_{2})), \psi(P_{1})\cup\dots\cup\psi(P_{k})\right) \le D.
\end{equation*}
Additionally, since $\beta\left(\text{supp}(c_{1}^{-1}c_{2})\right) \subset \psi(P_{1})\cup\dots\cup\psi(P_{k})$ and has a point in each of these pieces, there is some $D'\ge 0$ such that 
\begin{equation*}
    d_{\text{Haus}}\left(\beta\left(\text{supp}(c_{1}^{-1}c_{2})\right), \psi(P_{1})\cup\dots\cup\psi(P_{k})\right) \le D'.
\end{equation*}
We conclude that 
\begin{equation*}
    d_{\text{Haus}}\left(\beta\left(\text{supp}(c_{1}^{-1}c_{2})\right), \text{supp}\big(\alpha(c_{1})^{-1}\alpha(c_{2})\big)\right) \le D+D'
\end{equation*}
where $D$ and $D'$ are independent of $c_{1},c_{2}$. This shows \textit{(iii)} of Proposition~\ref{prop:characterisationofaptolicity} and completes the proof that $f$ is a quasi-isometry.
\end{proof}

We already know from~\cite[Fact~6.35]{BGT24} that the converse statement is not true; namely a quasi-isometry between $N\wr G$ and $M\wr H$ does not necessarily provide a quasi-isometry between a power of $N$ and a power of $M$.

\begin{remark}\label{rem:remark4.5.2}
Proposition~\ref{prop:constructionsofQI} can be adapted to construct quasi-isometries between permutational wreath products with infinite lamp groups, thus providing amenable analogs of the flexibility already observed in Remark~\ref{rm5.2}. For instance, combining the above method with the one of Proposition~\ref{prop:aptolicQIfromscalingQIbetweenquotients-Acase} shows that $\Z^{n}\wr_{\Z^k}\Z^d$ and $\Z^{n'}\wr_{\Z^k}\Z^d$ are quasi-isometric for any $n,n'\ge 1$ and $d\ge k\ge 1$. 
\end{remark}

We can now deduce Theorem~\ref{thm:classificationforvirtuallyabeliangroups} from these observations and our main theorem. 

{
\renewcommand{\proofname}{Proof of Theorem~\ref{thm:classificationforvirtuallyabeliangroups}}
\begin{proof}
We start proving~\textit{(i)}. First, if there is a quasi-isometry $\varphi\colon A_{1}\wr G\longrightarrow A_{2}\wr H$, then up to finite distance, $\varphi$ can be chosen aptolic, so we write
\begin{equation*}
    \varphi(c,p)=(\alpha(c),\beta(p))
\end{equation*}
for any $(c,p)\in A_{1}\wr G$. From (the proof of) Theorem~\ref{thm:maintheorem}, we know then that $\beta\colon G\rightarrow H$ is quasi-$\frac{d_{2}}{d_{1}}$-to-one. 

\noindent Conversely, assume that there exists a quasi-$\frac{d_{2}}{d_{1}}$-to-one quasi-isometry $G\rightarrow H$. Since there exists a biLipschitz equivalence $(\Z^{d_{1}})^{d_{2}}\longrightarrow (\Z^{d_{2}})^{d_{1}}$, Proposition~\ref{prop:constructionsofQI} shows that there is a quasi-isometry 
\begin{equation*}
    \varphi\colon\Z^{d_{1}}\wr G\longrightarrow \Z^{d_{2}}\wr H.
\end{equation*}
We therefore have our conclusion once we have proved the following: 
\begin{claim}\label{claim:virtuallyabelianBILIPtoZ^d}
An infinite virtually abelian group of growth degree $d\ge 1$ is biLipschitz equivalent to $\Z^d$. 
\end{claim}

{
\renewcommand{\proofname}{Proof of Claim~\ref{claim:virtuallyabelianBILIPtoZ^d}.}
\renewcommand{\qedsymbol}{$\blacksquare$}
\begin{proof}
If $A$ is virtually abelian and has growth degree $d$, it contains $\Z^d$ as a finite-index subgroup, and thus there is a quasi-$[A:\Z^d]$-to-one quasi-isometry $f\colon A\rightarrow \Z^d$ (any quasi-inverse of the inclusion $\Z^d \hookrightarrow A$). Since $\text{Sc}(\Z^d)=\R_{>0}$, we may fix an arbitrary quasi-$\frac{1}{[A:\Z^d]}$-to-one quasi-isometry $g\colon \Z^{d}\rightarrow \Z^{d}$, and the composition $g\circ f\colon A\rightarrow \Z^{d}$ is therefore quasi-one-to-one according to Proposition~\ref{prop:stabilitypropertiesforscalingQI}\textit{(ii)}. Thus, from Theorem~\ref{thm:Whytethm}, we deduce that it lies at finite distance from a bijection, which is the desired biLipschitz equivalence $A\rightarrow \Z^{d}$.
\end{proof}
\renewcommand{\qedsymbol}{$\square$}
}

\noindent Now, Claim~\ref{claim:virtuallyabelianBILIPtoZ^d} shows that $A_{1}$ is biLipschitz equivalent to $\Z^{d_{1}}$, so Proposition~\ref{prop:BiLipwreathproductscolor} implies that there is a biLipschitz equivalence $f_{1}\colon A_{1}\wr G\longrightarrow \Z^{d_{1}}\wr G$. Likewise, there is a biLipschitz equivalence $f_{2}\colon \Z^{d_{2}}\wr H\longrightarrow A_{2}\wr H$. The composition 
\begin{equation*}
    f_{2}\circ\varphi\circ f_{1}\colon A_{1}\wr G\longrightarrow A_{2}\wr H
\end{equation*}
is the quasi-isometry we are looking for. 

\noindent Let us now focus on~\textit{(ii)}. The left-to-right direction is a particular case of Corollary~\ref{cor:BLE}. Conversely, if $d_{1}=d_{2}$, it follows from Claim~\ref{claim:virtuallyabelianBILIPtoZ^d} that $A_{1}$ and $A_{2}$ are biLipschitz equivalent, so Proposition~\ref{prop:BiLipwreathproductscolor} implies that $A_{1}\wr G$ and $A_{2}\wr G$ are biLipschitz equivalent. Additionally, Proposition~\ref{prop:BiLipwreathproductsbase} ensures that $A_{2}\wr G$ and $A_{2}\wr H$ are biLipschitz equivalent, and composing these two biLipschitz equivalences provides a biLipschitz equivalence between $A_{1}\wr G$ and $A_{2}\wr H$.
\end{proof}}

{
\renewcommand{\proofname}{Proof of Corollary~\ref{cor:classificationforvirtuallyabeliangroups2}}
\begin{proof}
Assume that $A_{1}\wr B_{1}$ is quasi-isometric to $A_{2}\wr B_{2}$. In particular, F\o lner functions of these two groups must coincide up to $\simeq$ (see Theorem~\ref{thm:profilQIinvariant}), which, by Theorem~\ref{thm:formulaforprofilewreathproducts}, implies that 
\begin{equation*}
    \left(n^{\text{deg}(A_{1})}\right)^{n^{\text{deg}(B_{1})}} \simeq  \left(n^{\text{deg}(A_{2})}\right)^{n^{\text{deg}(B_{2})}}.
\end{equation*}
It follows that $n^{\text{deg}(B_{1})}\ln(n) \simeq n^{\text{deg}(B_{2})}\ln(n)$, and thus $n^{\text{deg}(B_{1})} \simeq n^{\text{deg}(B_{2})}$. The latter implies $\text{deg}(B_{1})=\text{deg}(B_{2})$ (Example~\ref{ex:examplesofasymptoticequivalence}\textit{(i)}) as claimed. 

\noindent Conversely, if $B_{1}$ and $B_{2}$ both have growth degree $d$, they are both biLipschitz equivalent to $\Z^{d}$ by Claim~\ref{claim:virtuallyabelianBILIPtoZ^d}, 
and since $\text{Sc}(\Z^d)=\R_{>0}$ (see Proposition~\ref{prop:examplesofscalinggroups}), we may fix an arbitrary quasi-$\frac{\text{deg}(A_{2})}{\text{deg}(A_{1})}$-to-one quasi-isometry $\Z^{d}\rightarrow\Z^{d}$. The latter then provides us a quasi-$\frac{\text{deg}(A_{2})}{\text{deg}(A_{1})}$-to-one quasi-isometry $B_{1}\rightarrow B_{2}$. Additionally, as there is a biLipschitz equivalence $A_{1}^{\text{deg}(A_{2})}\longrightarrow A_{2}^{\text{deg}(A_{1})}$, we conclude with Proposition~\ref{prop:constructionsofQI} that there is a quasi-isometry
\begin{equation*}
    A_{1}\wr B_{1} \longrightarrow A_{2}\wr B_{2}
\end{equation*}
as desired. The proof is complete. 
\end{proof}}

Concerning iterated wreath products, we prove first Proposition~\ref{prop:mixingofscalingconditions}, and Proposition~\ref{prop:lamplighterrigidity} follows immediately.

{
\renewcommand{\proofname}{Proof of Proposition~\ref{prop:mixingofscalingconditions}}
\begin{proof}
Assume that there exists a quasi-isometry 
\begin{equation*}
    \Z_{n}\wr(N_{1}\wr G)\longrightarrow \Z_{m}\wr(N_{2}\wr H).
\end{equation*}
From~\cite[Proposition~1.3]{GT24a}, $N_{1}\wr G$ and $N_{2}\wr H$ have the thick bigon property, and they are both amenable, so we may apply Theorem~\ref{thm:AlamplightersoverTBPgroups-rigiditypart} to deduce that there exist $a, r, s\ge 1$ such that $n=a^{r}$, $m=a^{s}$ and a quasi-$\frac{s}{r}$-to-one quasi-isometry $N_{1}\wr G\longrightarrow N_{2}\wr H$. By Theorem~\ref{thm:maintheorem}, the latter is also quasi-$\frac{n_{2}}{n_{1}}$-to-one, and thus Lemma~\ref{lem:uniquenessscalingfactor} forces $\frac{s}{r}=\frac{n_{2}}{n_{1}}$.
\end{proof}
}

We conclude with the proof of Corollary~\ref{cor:classificationofiteratedwreathproducts}.

{
\renewcommand{\proofname}{Proof of Corollary~\ref{cor:classificationofiteratedwreathproducts}}
\begin{proof}\textit{(i)} The left-to-right direction is a particular case of Proposition~\ref{prop:mixingofscalingconditions}. Conversely, given a quasi-$\frac{d_{2}}{d_{1}}$-to-one quasi-isometry $G\rightarrow H$, we apply Proposition~\ref{prop:constructionsofQI} to construct an aptolic quasi-isometry
\begin{equation*}
    A_{1}\wr G\longrightarrow A_{2}\wr H.
\end{equation*}
From the proof of Proposition~\ref{prop:constructionsofQI} and Lemma~\ref{lm:scalingfactorsforwreathproducts}, this quasi-isometry is quasi-$\frac{d_{2}}{d_{1}}$-to-one, thus also quasi-$\frac{s}{r}$-to-one by assumption. As $n=a^{r}$ and $m=a^{s}$, we apply~\cite[Proposition~3.12]{GT24b} to deduce that there is a quasi-isometry
\begin{equation*}
    \Z_{n}\wr(A_{1}\wr G) \longrightarrow \Z_{m}\wr(A_{2}\wr H)
\end{equation*}
as claimed. 

\noindent Lastly,~\textit{(ii)} follows from a combination of~\textit{(i)} and of Lemma~\ref{lem:uniquenessscalingfactor}.
\end{proof}}

\section{Comments and questions}\label{subsection4.6}

Results from~\cite{BGT24} still apply when lamp groups have superpolynomial but subexponential growth. For such groups, the growth function can exhibit several behaviors, and our understanding is much more limited, but from the proof of Theorem~\ref{thm:maintheorem} we can deduce the following:
\begin{corollary}\label{cor:intermediategrowth}
Let $\Gamma$ and $\Lambda$ be finitely generated groups with subexponential and superpolynomial growth. Let $G$ and $H$ be finitely presented groups from $\mathcal{M}_{\text{exp}}$. 
\begin{enumerate}[label=(\roman*)]
    \item If $\gamma_{\Gamma}(n)\simeq e^{n^{\alpha}}$, $\gamma_{\Lambda}(n)\simeq e^{n^{\beta}}$ for some $\alpha,\beta\in (0,1)$, and if \;$\Gamma\wr G$ and $\Lambda\wr H$ are quasi-isometric, then $\alpha=\beta$.
    \item If $\gamma_{\Gamma}(n)\simeq e^{\ln(n)n^{\alpha}}$, $\gamma_{\Lambda}(n)\simeq e^{\ln(n)n^{\beta}}$ for some $\alpha,\beta\in (0,1)$, and if \;$\Gamma\wr G$ and $\Lambda\wr H$ are quasi-isometric, then $\alpha=\beta$.
    \item If $\gamma_{\Gamma}(n)\simeq e^{\frac{n}{\ln^{\circ k}(n)}}$, $\gamma_{\Lambda}(n)\simeq e^{\frac{n}{\ln^{\circ d}(n)}}$ for some integers $d,k\ge 1$, and if \;$\Gamma\wr G$ and $\Lambda\wr H$ are quasi-isometric, then $k=d$.
\end{enumerate}
\end{corollary}

Here, $\gamma_{G}$ stands for the growth function of a group $G$, and for $k\ge 1$, $\ln^{\circ k}=\ln(\ln(\dots))$ denotes the $k-$th iteration of the log with itself.

\begin{proof}
The proof of Claim~\ref{claim4.2} in Section~\ref{subsection4.4} in fact shows that some powers of $\gamma_{\Gamma}$ are equivalent to some powers of $\gamma_{\Lambda}$. More precisely, given our $(C,K)-$quasi-isometry 
\begin{equation*}
    \Gamma\wr G \longrightarrow \Lambda\wr H, \; (c,p)\longmapsto (\alpha(c),\beta(p)),
\end{equation*}
one has $\gamma_{\Gamma}^{|A|}(n) \simeq \gamma_{\Lambda}^{|\overline{\beta}(\beta(A)^{+K})^{+K}|}(n)$ for any finite subset $A\subset G$. A direct computation then shows that, for the three possible behaviours, such equivalences force equalities of the parameters. 
\end{proof}

Additionally, one also deduces from the equivalences
\begin{equation*}
    \gamma_{\Gamma}^{|A|}(n) \simeq \gamma_{\Lambda}^{|\overline{\beta}(\beta(A)^{+K})^{+K}|}(n),\; A\subset G \;\text{finite}
\end{equation*}
that these three ranges of behaviours are quasi-isometrically distinct when taking wreath products: for instance, if $\gamma_{\Gamma}(n)\simeq e^{n^{\alpha}}$ and $\gamma_{\Lambda}(n)\simeq e^{\ln(n)n^{\beta}}$ for some $\alpha,\beta\in (0,1)$ and if $G,H$ are finitely presented groups from $\mathcal{M}_{\text{exp}}$, then $\Gamma\wr G$ and $\Lambda\wr H$ are not quasi-isometric.

Moreover, we know that these three types of behaviour occur for finitely generated groups. It is shown in~\cite[Theorem~B]{EZ18} that some periodic Grigorchuk groups, among which the first Grigorchuk group, fall into the first class, and the logarithmic growth exponent is computed explicitly. Other examples appear in~\cite[Theorem~1]{BE12}, where the authors construct two families of groups $(K_{k})_{k\in\N}$ and $(H_{k})_{k\in\N}$ such that 
\begin{equation*}
    \gamma_{K_{k}}(n)\simeq e^{n^{1-(1-\alpha)^{k}}}\; \text{and}\; \gamma_{H_{k}}(n)\simeq e^{\ln(n)n^{1-(1-\alpha)^{k}}}
\end{equation*}
where $\alpha \cong 0.7674$ is a fixed constant. Thus, for instance, it follows from Corollary~\ref{cor:intermediategrowth} that, for $n\ge 2$, $K_{r}\wr \text{BS}(1,n)$ and $K_{s}\wr \text{BS}(1,n)$ are quasi-isometric if and only if $r=s$. 

Lastly, the existence of groups with growth functions asymptotically equal to $e^{\frac{n}{\ln^{\circ k}(n)}}$ for integers $k\ge 1$ is proved in~\cite[Theorem~A]{BE14}. 

However, in the case of intermediate growth, the fact that powers of the growth functions must be equivalent does not provide any valuable information on scaling properties for $\beta$. Thus we may ask:
\begin{question}
Let $\Gamma$ be an intermediate growth group, and let $G$ be a finitely presented amenable group from $\mathcal{M}_{\text{exp}}$. Is $\text{Sc}(\Gamma\wr G)$ reduced to $\left\lbrace 1\right\rbrace$?
\end{question}

Note that the proof of Proposition~\ref{prop:lamplighterrigidity} only uses the fact that $N\wr G$ has the thick bigon property and a trivial scaling group. Therefore, any group with these two properties is lamplighter-rigid. 

\begin{corollary}
Let $n,m\ge 2$ be two integers. Let $H$ be a finitely generated amenable group satisfying the thick bigon property. If $\text{Sc}(H)=\lbrace 1\rbrace$, then $\Z_{n}\wr H$ and $\Z_{m}\wr H$ are quasi-isometric if and only if $n=m$. 
\end{corollary}

It would be interesting to find other examples of classes of groups satisfying these two properties. Plausible candidates could be direct products of standard wreath products. Hence:

\begin{question}
Let $F_{1}, F_{2}$ be non-trivial finite groups. Let $G,H$ be finitely generated one-ended amenable groups. Is it true that $\text{Sc}\big((F_{1}\wr G)\times(F_{2}\wr H)\big)=\lbrace 1\rbrace$? At least under some assumptions on $G$ and $H$ to determine?
\end{question}

On the other hand, wreath products $F\wr K$ with $F$ finite and $K$ amenable finitely presented and one-ended, do not satisfy the thick bigon property, but still have trivial scaling group~\cite[Proposition 6.8]{GT22}. Hence:

\begin{question}
Let $F$ be a non-trivial finite group, and let $K$ be a finitely presented one-ended amenable group. Is $F\wr K$ lamplighter-rigid?
\end{question}

Even though we focused on wreath products, it is worth noticing that scaling quasi-isometries also parametrize the existence of quasi-isometries between other halo products, as defined in Chapter~\ref{chap:chapter2}. Among these, lampshufflers are of great interest and, for instance, combining~\cite[Corollary~8.9]{GT24a} with Theorem~\ref{thm:classificationforvirtuallyabeliangroups}, it follows that $\shuf{\Z^n\wr G}$ and $\shuf{\Z^m\wr H}$ are quasi-isometric if and only if $n=m$, when $G$ and $H$ are amenable finitely presented groups from $\mathcal{M}_{\text{exp}}$. 

Finally, regarding quasi-isometric classifications of wreath products, it is natural to wonder whether there is an iterated version of Corollary~\ref{cor:classificationforvirtuallyabeliangroups2}, at least for wreath products of the form 
\begin{equation*}
\Z^{n_{1}}\wr\big(\Z^{n_{2}}\wr(\dots\wr(\Z^{n_{r}}\wr\Z^{k}))\big),\; n_{1},\dots,n_{r},k\ge 1. 
\end{equation*}
The existence of a quasi-isometry 
\begin{equation*}
\Z^{n_{1}}\wr\big(\Z^{n_{2}}\wr(\dots\wr(\Z^{n_{r}}\wr\Z^{k}))\big) \longrightarrow \Z^{m_{1}}\wr\big(\Z^{m_{2}}\wr(\dots\wr(\Z^{m_{r}}\wr\Z^{k'}))\big)
\end{equation*}
already imposes $k=k'$, by computing F\o lner functions, and it seems not unreasonable to believe that it imposes an arithmetic condition of the kind $\frac{m_{1}}{n_{1}}=\frac{m_{2}}{n_{2}}=\dots=\frac{m_{r}}{n_{r}}$, in the spirit of Proposition~\ref{prop:mixingofscalingconditions}. Additionally, note that Proposition~\ref{prop:constructionsofQI} already provides the flexibility part, if all these quotients are equal. Hence:
\begin{question}
Let $n_{1},m_{1},\dots,n_{r},m_{r}\ge 1$ and $k,k'\ge 1$. Is it true that $\Z^{n_{1}}\wr\big(\Z^{n_{2}}\wr(\dots\wr(\Z^{n_{r}}\wr\Z^{k}))\big)$ and $\Z^{m_{1}}\wr\big(\Z^{m_{2}}\wr(\dots\wr(\Z^{m_{r}}\wr\Z^{k'}))\big)$ are quasi-isometric if and only if $\frac{m_{1}}{n_{1}}=\frac{m_{2}}{n_{2}}=\dots=\frac{m_{r}}{n_{r}}$ and $k=k'$?
\end{question}

\afterpage{\blankpage}

\chapter{Isoperimetric profiles of lamplighter-like groups}\label{chap:chapter5}

This chapter presents the paper~\cite{cordum25}, written in collaboration with Corentin Correia.

\vspace{0.3cm}

\minitoc

\section{Introduction}

\setcounter{theoremletter}{0}

It is a recurrent theme in geometric group theory to understand the collection of all maps between two given finitely generated groups that are compatible with their large-scale geometries. Such collections include for instance quasi-isometries, coarse embeddings, and more generally regular maps. The motivation behind this program is that the large-scale geometry of a group is in fact deeply related to its algebraic structure. 

\sloppy Several milestones have been achieved in the study of quasi-isometries of many classes of groups, among which abelian free groups, lamplighters over virtually cyclic groups and other $\text{SOL}$-like groups~\cite{EFW12, EFW13}, Baumslag-Solitar groups~\cite{FM98, FM99, Why01}, or lamplighters over one-ended groups~\cite{GT24b}.

If one rather wants to exclude the existence of such maps between some spaces, an efficient strategy is to use invariants, or even monotonous quantities, that are computable in practice. Numerous invariants have been introduced for quasi-isometries, including 
 the isoperimetric profiles which are monotonous under regular maps~\cite[Theorem~5.5]{DKLMT22}. We refer the reader to Section~\ref{sec:appQIandRegularMaps} for other examples of invariants (asymptotic dimension, volume growth).

The aim of the paper is to compute the isoperimetric profiles of lampshufflers, and more generally of the so-called~\textit{halo products} introduced in~\cite{GT24a}. 

\paragraph{Isoperimetric profiles and F\o lner functions.} For a finitely generated group $G$ with a finite generating set $S_{G}$, and for a $p\ge 1$, its $\ell^{p}-$\textit{isoperimetric profile} is the function $\profp{G}\colon \N\rightarrow\R_{+}$ given by 
\begin{equation*}
    \profp{G}(n) \defeq\sup_{\substack{f\colon G\to\R_+\\|\supp{f}|\leq n}}{\frac{\|f\|_p}{\|\nabla f\|_p}}
\end{equation*}
where the support of $f\colon G\rightarrow\R_{+}$ is $\supp{f}\defeq\lbrace g\in G : f(g)\neq 0\rbrace$ and the $\ell^{p}-$norm of its gradient is defined by 
\begin{equation*}
   \|\nabla f\|_{p}^{p}\defeq\sum_{g\in G,\; s\in S_{G}}\left|f(g)-f(gs)\right|^{p}. 
\end{equation*}

\begin{remark}\label{rem:ConventionProfile}
We warn the reader that many authors introduce the $\ell^{p}-$isoperimetric profile with $\|\cdot\|_{p}^{p}$ instead of $\|\cdot\|_{p}$ in the definition, so their $\ell^{p}-$isoperimetric profile is $\profp{G}(x)^p$ with our conventions.
\end{remark}

The $\ell^{p}-$isoperimetric profile of a group $G$ is the generalized inverse of its $\ell^{p}-$\textit{F\o lner function} $\folp{G}\colon\N\rightarrow \R_{+}$, defined as
 \begin{equation*}
     \folp{G}(n) \defeq \inf\left\lbrace |\supp{f}| : \frac{\|\nabla f\|_p}{\|f\|_p} \le \frac{1}{n} \right\rbrace.
 \end{equation*}

In the case $p=1$, these functions are simply called~\textit{isoperimetric profile} and~\textit{F\o lner function}, and have a simpler definition (up to asymptotic behaviour), namely
\begin{equation*}
    \prof{G}(n) \defeq \sup_{|A|\le n}\frac{|A|}{|\partial_{G} A|}\;\text{ and }\;\text{F\o l}_{G}(n) \defeq \inf\left\lbrace |A| : \frac{|\partial_{G}A|}{|A|} \le \frac{1}{n} \right\rbrace,
\end{equation*}
where $\partial_{G}A \defeq AS_{G}\setminus A=\lbrace g\in G\setminus A : \exists s\in S_{G}, \exists h\in A, g=hs\rbrace$ is the boundary of $A$ in $G$. Note that we only find the $\ell^{1}-$F\o lner function in the literature. In this paper we introduce the more general $\ell^{p}$ versions for $p\ge 1$.

Without loss of generality, we may and do assume that the $\ell^{p}-$isoperimetric profile and the $\ell^{p}-$F\o lner function are real inverses of each other, and not only generalized inverses; see Remark~\ref{rem:InverseEachOther}.

Notice that the $\ell^{p}-$isoperimetric profile of a finitely generated group is bounded if and only if the group is not amenable. Therefore, we will only be interested in $\ell^{p}-$isoperimetric profiles of amenable groups. The asymptotic behaviour of the $\ell^{p}-$isoperimetric profile is, somehow, a measurement of its amenability; the faster it goes to infinity, the “more amenable” the group is.

Among amenable groups, the $\ell^{p}-$isoperimetric profile (or equivalently the $\ell^{p}-$F\o lner function) has been computed for many finitely generated groups. Given $p\ge 1$, we have for instance:
\begin{itemize}
    \item $\profp{G}(n) \simeq n^{\frac{1}{d}}$ if $G$ has polynomial growth of degree $d\ge 1$;
    \item $\profp{G}(n)\simeq \ln(n)$ for $G=\text{BS}(1,k)$, $k\ge 2$, or $G=F\wr\Z$, where $F$ is a non-trivial finite group;
    \item $\profp{G}(n) \simeq \ln(n)$ for any polycyclic group $G$ with exponential growth~\cite{Pit95, Pit00}, or more generally any exponential growth group within the class GES of Tessera~\cite[Corollary 5]{Tes13};
    \item $\profp{F\wr N}(n)\simeq \ln(n)^{\frac{1}{d}}$ with $F$ a non-trivial finite group, and $N$ having polynomial growth of degree $d\ge 1$~\cite{Ers03};
    \item for any non-decreasing function $f\colon \R_{+}\rightarrow \R_{+}$ such that $x\longmapsto \frac{x}{f(x)}$ is non-decreasing, Brieussel and Zheng constructed in~\cite[Theorem~1.1]{BZ21} a finitely generated group $H$ with exponential volume growth having isoperimetric profile 
    \begin{equation*}
        \profp{H}(x)\simeq \frac{\ln(x)}{f(\ln(x))}.
    \end{equation*}
    We will refer to such a group as a~\textit{Brieussel-Zheng's group}. 
\end{itemize}

\begin{figure}[H]
  \centering
  \includegraphics[width=0.75\linewidth]{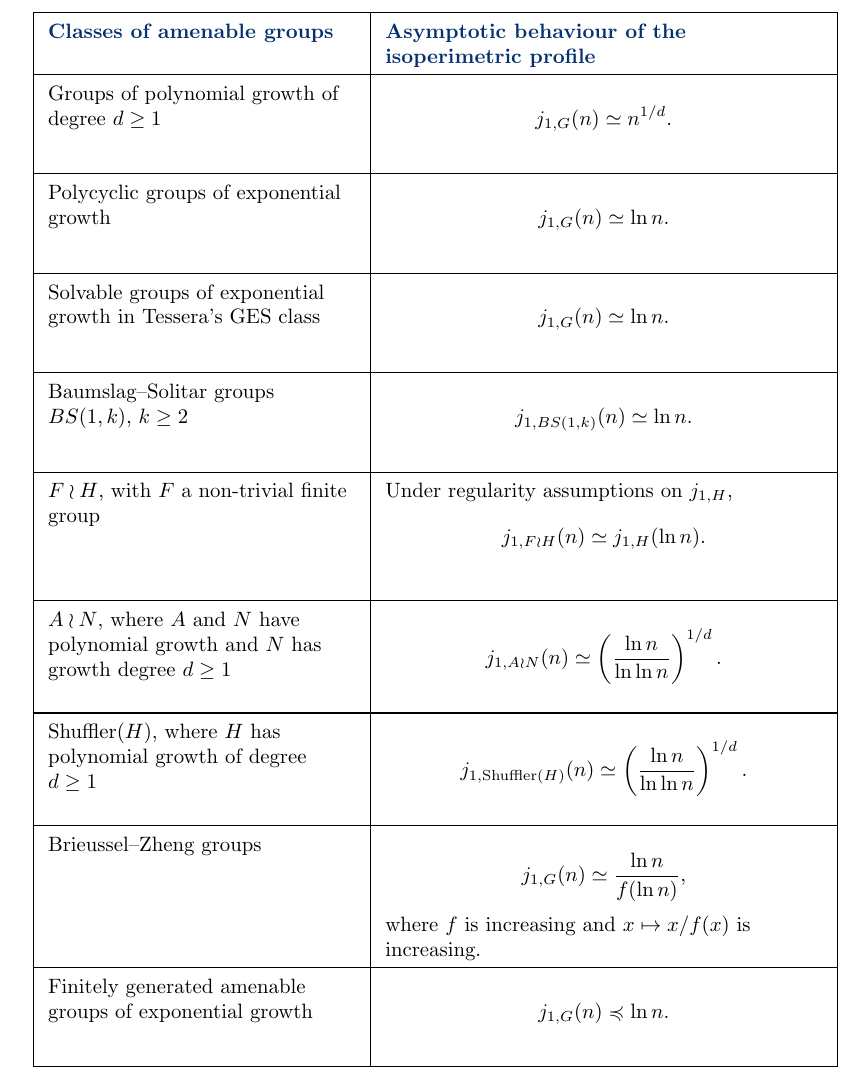}
  \caption{A brief overview of the state of the art on isoperimetric profiles.}
  \label{fig:chemin}
\end{figure}

In fact, isoperimetric profiles and F\o lner functions have been studied in the more general framework of bounded degree graphs; see Section~\ref{sec:preliminaries} for details. For now, let us simply mention that, in this setup, Erschler's estimates for $\ell^{1}-$F\o lner functions of lamplighter graphs~\cite{Ers03,Ers06}  will be a key ingredient in our strategy (see Section~\ref{sec:secondstrategy}).

It is a well-known fact that if $p>q$, then $\profp{G}(x)\preccurlyeq j_{q,G}(x)$, see Lemma~\ref{lem:MonotonuousProfile}. Moreover, a stronger phenomenon is conjectured: $\profp{G}(x)\simeq j_{q,G}(x)$ for all $p,q\ge 1$.

\paragraph{Lampshuffler groups.} As a starting point of our work, let us first focus on~\textit{lampshuffler groups}. Later in the introduction, we will present the more general notion of halo products, constructed in a similar way to wreath products (as lampshufflers). 

Given a group $H$, the~\textit{lampshuffler group} over $H$ is the semi-direct product 
\begin{equation*}
    \shuf{H} \defeq \fsym{H}\rtimes H
\end{equation*}
where $\fsym{H}$ is the group of finitely supported bijections $H\rightarrow H$, and where $H$ acts on the latter as $(h\cdot\sigma)(x)\defeq h\sigma(h^{-1}x)$. These groups already appeared several times in the literature, in relations with many topics of interest in group theory, see for instance~\cite{Yad09, HO16, BZ19, EZ21, SCZ21, GT24a, Sil24}.

Let us first present what is known about profiles of lampshufflers and our results for this class of groups.

\paragraph{Isoperimetric profiles of lampshufflers.} In~\cite[Theorem~6.8]{SCZ21}, Saloff-Coste and Zheng establish a general lower bound on $\profp{\shuf{H}}$ for a finitely generated group $H$, of the form 
\begin{equation*}
    \profp{H}\left(\frac{\ln(x)}{\ln(\ln(x))}\right) \preccurlyeq \profp{\shuf{H}}(x).
\end{equation*}
Their proof will be useful since we will present a natural generalization for halo products, see Corollary~\ref{cor:lowerboundonProf}.

From~\cite[Theorem~6.7]{SCZ21}, we also know upper bounds on $\profp{\shuf{H}}$ for general groups $H$, and~\cite[Corollary~1.4]{EZ21} provides a lower bound on $\text{F\o l}_{1,\shuf{H}}$, equivalently an upper bound on $\prof{\shuf{H}}$.

\begin{theorem}[{\cite[Theorem~6.7]{SCZ21}}]
Let $H$ be a finitely generated group, with a finite generating set $S_{H}$. For $p=1,\; 2$, we have
\begin{equation*}
    \profp{\shuf{H}}(x)\preccurlyeq \beta_{H,S_{H}}^{-1}\left(\frac{\ln(x)}{\ln(\ln(x))}\right).
\end{equation*}
\end{theorem}
Here, $\beta_{H,S_{H}}$ refers to the growth function of $H$ with respect to the finite generating set $S_{H}$. Note that the statement in their paper appears with an exponent $p$, since we do not use the same convention; see Remark~\ref{rem:ConventionProfile}.

\begin{theorem}[{\cite[Corollary~1.4]{EZ21}}]
Let $H$ be a finitely generated group, with a finite generating set $S_H$. Then, we have
\begin{equation*}
    \text{F\o l}_{1,\shuf{H}}(x)\succcurlyeq \beta_{H,S_{H}}(x)^{\beta_{H,S_{H}}(x)}.
\end{equation*}
\end{theorem}

From these theorems, one can deduce for instance the isoperimetric profile of lampshufflers over polynomial growth groups:
\begin{equation*}
    \prof{\shuf{H}}(x)\simeq\left(\frac{\ln(x)}{\ln(\ln(x))}\right)^{\frac{1}{d}},
\end{equation*}
when $H$ has polynomial growth of degree $d\ge 1$.

For many groups, the lower and upper bounds provided by these results are not the same, so they do not provide precise estimates on the isoperimetric profiles of lampshufflers. Theorem~\ref{thm:boundsForProfile intro} below provides finer estimates. It is in fact an application of Corollary~\ref{cor:ProfHalo intro}, that we deduce from Theorem~\ref{th:recallThm intro}, stated in the more general context of halo products. Before defining these groups, let us focus on lampshufflers and the consequences of Theorem~\ref{thm:boundsForProfile intro}.

In our statements, saying that $\profp{H}$ satisfies Assumption~$(\star)$ means that $\profp{H}(Cx)=O(\profp{H}(x))$ for any $C>0$.

\begin{theoremletter}[see Corollary~\ref{cor:encadrementduprofildeshuf}]\label{thm:boundsForProfile intro}
Let $p\ge 1$. Let $H$ be a finitely generated amenable group whose $\ell^{p}-$isoperimetric profile $\profp{H}$ satisfies Assumption~$(\star)$. Then the $\ell^{p}-$isoperimetric profile $\profp{\shuf{H}}$ of $\shuf{H}$ satisfies
\begin{equation*}
    \profp{H}\left(\frac{\ln(x)}{\ln(\ln(x))}\right) \preccurlyeq \profp{\shuf{H}}(x) \preccurlyeq \prof{H}(\ln(x)). 
\end{equation*}
\end{theoremletter}

Assumption $(\star)$ already appeared several times in the literature, see e.g.~\cite{Ers03, Cor25} for the case $p=1$, and does not seem to be restrictive. In fact, to our knowledge, there is no known example of a finitely generated amenable group whose isoperimetric profiles do not satisfy Assumption~$(\star)$. For instance, it is easy to check that Brieussel-Zheng's groups satisfy this assumption (see~\cite[after Corollary 4.1]{Cor25}), as well as all the examples of isoperimetric profiles we mentioned above.

An immediate consequence of Theorem~\ref{thm:boundsForProfile intro} is the next statement.

\begin{corollaryletter}\label{cor:corB}
Let $p\ge 1$. Let $H$ be a finitely generated amenable group whose $\ell^p$-isoperimetric profile $\profp{H}$ satisfies Assumption~$(\star)$. Assume moreover that
\begin{equation*}
   \profp{H}\left(\frac{\ln(x)}{\ln(\ln(x))}\right) \simeq \profp{H}(\ln(x))\;\text{ and }\;\profp{H}(x)\simeq\prof{H}(x).
\end{equation*}
Then one has 
\begin{equation*}
    \profp{\shuf{H}}(x) \simeq \profp{H}(\ln(x)).
\end{equation*}
\end{corollaryletter}

In practice, this result applies for many groups having slow enough profiles, for instance solvable Baumslag-Solitar groups $\text{BS}(1,n)$, lamplighters over $\Z^d$, or polycyclic groups with exponential growth. In particular, for the latter class, we recover~\cite[Corollary~6.9]{SCZ21}.

\begin{remark}
Corollary~\ref{cor:corB} implies that, if the $\ell^{p}-$isoperimetric profiles of $H$ all have the same asymptotic behaviour, then the same holds for $\shuf{H}$, under mild assumptions on $H$.
\end{remark}

One class of groups for which estimates from~\cite{EZ21} and from~\cite{SCZ21} are not optimal is the one of~\textit{iterated lampshufflers}, defined inductively by $\shufn{0}{H}\defeq H$ and $\shufn{n}{H}\defeq \shuf{\shufn{n-1}{H}}$ if $n\ge 1$. It turns out that iterations of Theorem~\ref{thm:boundsForProfile intro} yield finer estimates, that we record in the two following statements. 

\begin{propositionletter}[see Proposition~\ref{prop:profileofshufnofpolynomialgrowthgroups}]\label{prop:profileofshufnofpolynomialgrowthgroupsINTRO}
Let $H$ be a finitely generated group of polynomial growth of degree $d\ge 1$. Then one has 
\begin{equation*}
    \profp{\shufn{n}{H}}(x) \simeq \left(\frac{\ln^{\circ n}(x)}{\ln^{\circ (n+1)}(x)}\right)^{\frac{1}{d}}
\end{equation*}
for any integer $n\ge 1$ and any real number $p\ge 1$.
\end{propositionletter}

\begin{propositionletter}[see Proposition~\ref{prop:profilesofshufn}]\label{prop:profilesofshufnINTRO}
Let $p\geq 1$. Let $H$ be a finitely generated amenable group whose $\ell^{p}-$isoperimetric profile $\profp{H}$ satisfies Assumption~$(\star)$. Suppose that
\begin{equation*}
    \profp{H}\left(\frac{\ln(x)}{\ln(\ln(x))}\right) \simeq \profp{H}(\ln(x))\; \text{ and }\;\profp{H}(x)\simeq\prof{H}(x).
\end{equation*}
Then, we have 
\begin{equation*}
   \profp{\shufn{n}{H}}(x) \simeq \profp{H}(\ln^{\circ n}(x))
\end{equation*}
for all $n\ge 0$.
\end{propositionletter}

\paragraph{Halo products.} In~\cite{GT24a}, Genevois and Tessera introduced a general class of groups, called~\textit{halo products}, as a natural generalization of wreath products. This class encompasses lampshufflers, lampjugglers and lampcloners, and constitutes the suitable framework for our main result Theorem~\ref{th:recallThm intro}.

\begin{definition}\label{def:haloproducts}
Let $X$ be a set. A~\textit{halo of groups $\halo$ over $X$} is the data, for any subset $S\subset X$, of a group $L(S)$ such that:
\begin{itemize}
    \item for all $R,S\subset X$, if $R\subset S$ then $L(R)\leqslant L(S)$;
    \item $L(\emptyset)=\lbrace 1\rbrace$ and $L(X)=\langle L(S) : S\subset X \;\text{finite}\rangle$;
    \item for all $R,S\subset X$, $L(R\cap S)=L(R)\cap L(S)$.
\end{itemize}
\end{definition}

Given an action $H\curvearrowright X$ and a morphism $\alpha\colon H \longrightarrow \text{Aut}(L(X))$ satisfying $\alpha(h)(L(S))=L(hS)$ for any $S\subset X$ and $h\in H$, the~\textit{permutational halo product} $\halo_{X,\alpha}H$ is the semi-direct product 
\begin{equation*}
    \halo_{X,\alpha}H \defeq L(X)\rtimes_{\alpha}H.
\end{equation*}

\sloppy In this paper, we focus on the case $X=H$. As mentioned, examples of halo products include
\begin{itemize}
    \item wreath products $F\wr H=(\bigoplus_H F)\rtimes H$, for which $L(S)=\bigoplus_S F$;
    \item lampshufflers $\shuf{H}=\fsym{H}\rtimes H$, for which $L(S)=\fsym{S}$,
\end{itemize}
and many other examples are introduced and studied in~\cite{GT24a}, such as
\begin{itemize}
    \item lampjugglers $\juggler{s}{H}=\fsym{H\times\lbrace 1,\dots,s\rbrace}\rtimes H$, with an integer $s\geq 1$, for which $L(S)=\fsym{S\times\lbrace 1,\dots,s\rbrace}$;
    \item lampcloners $\cloner{H}=\text{FGL}(H)\rtimes H$, with a field $\field$, for which $L(S)=\text{FGL}(S)$;
    \item lampdesigners $\designer{H}=(F\wr_{H}\fsym{H})\rtimes H$, with a non-trivial finite group $F$, for which $L(S)=F\wr_{S}\fsym{S}$, 
\end{itemize}
where $S$ denotes any subset of $H$. Here, $\text{FGL}(H)$ denotes the group of linear automorphisms of the abstract $\field-$vector space $V_{H}$ admitting $H$ as a basis, fixing all but finitely many basis vectors. We refer the reader to Section~\ref{sec:halo} for more details.

The motivation in~\cite{GT24a} to introduce such a general framework is that the semi-direct product structure provides a foliation of these spaces that must be, if $H$ satisfies additional mild assumptions, “quasi-preserved” by quasi-isometries, allowing the authors to show strong rigidity phenomena for quasi-isometries between such spaces, and thus extending the classification already obtained in~\cite{GT24b}. 

\paragraph{Isoperimetric profiles of halo products.} In this paper, our aim is to show that the halo structure is also particularly well-suited for tracking isoperimetric profiles of these groups. Namely, we prove the following two estimates on their F\o lner functions. The terminologies and notations are explained just after the statement.

\begin{theoremletter}[see Proposition~\ref{prop:upperboundonFol} and Theorem~\ref{th:lowerboundonFol}]\label{th:recallThm intro}
Let $p\ge 1$. Let $H$ be a finitely generated amenable group and let $S_{H}$ be a finite generating set of $H$. Let $\halo H$ be a naturally generated halo product over $H$.
\begin{enumerate}[label=(\roman*)]
    \item\label{item:1} If $\halo H$ is large-scale commutative and has finitely generated blocks, then for any $s_{0}\in S_{H}$, there exists a constant $C>0$ such that
\begin{equation*}
    \left(\text{F\o l}_{L(\lbrace1_{H},s_{0}\rbrace)}(x)\right)^{C\cdot\text{F\o l}_{H}(x)} \preccurlyeq \folp{\halo H}(x).
\end{equation*}
    \item If $\halo H$ has consistent blocks, then there exists a constant $C>0$ such that
\begin{equation*}
    \folp{\halo H}(x) \preccurlyeq \folp{H}(x)\cdot \Lambda_{\halo H}(C\cdot \folp{H}(x)).
\end{equation*}
\end{enumerate}
\end{theoremletter}

A halo product $\halo H$ has~\textit{finite} (resp.~\textit{finitely generated})~\textit{blocks} if $L(S)$ is finite (resp. finitely generated) for any finite subset $S\subset H$, and $\halo H$ has~\textit{consistent} blocks if its blocks are finite and moreover the cardinality of $L(S)$ only depends on $|S|$. This assumption allows to define, as in~\cite{GT24a}, a function $\Lambda_{\halo H}\colon\N\rightarrow\N$ sending any $n\in\N$ to $|L(S)|$ where $|S|=n$, called the~\textit{lamp growth sequence} of $\halo H$.

Moreover, $\halo H$ is~\textit{large-scale commutative} if there is $D\ge 0$ such that for any subsets $R,S\subset H$ that are at least $D$ far apart in $H$, $L(R)$ and $L(S)$ commute in $L(H)$. Such a notion has been introduced in~\cite{GT24a} as a key assumption to understand the general form of quasi-isometries between halo groups. Lastly, $\halo H$ is~\textit{naturally generated} if it admits the natural and simplest generating set that we can imagine for a halo product, in view of the classical finite generating sets for lamplighters and lampshufflers.

For instance, a lamplighter $F\wr H$ and a lampshuffler $\shuf{H}$ are large-scale commutative (with $D=0$ for $F\wr H$, $D=1$ for $\shuf{H}$), are naturally generated, and have consistent blocks, with lamp growth sequences given by
\begin{equation*}
    \Lambda_{F\wr H}(n)=|F|^n \;\; \text{and} \;\; \Lambda_{\shuf{H}}(n)=n!.
\end{equation*}

\paragraph{Towards the proof of Theorem~\ref{th:recallThm intro}.} The lower bound is a direct computation, presented in Section~\ref{sec:UpperBoundFolner}, and inspired from the computations in the proof of~\cite[Theorem~6.8]{SCZ21} in the case of lampshufflers. We exhibit an explicit sequence of almost invariant functions $\halo H\rightarrow \R$ from one such sequence of $H$. Namely, a sequence $(f_{n})_{n\in\N}$ of functions $H\rightarrow\R$, realizing the $\ell^{p}-$isoperimetric profile of $H$ (or equivalently its $\ell^{p}-$F\o lner function), gives rise to a sequence $(g_{n})_{n\in\N}$ for $\halo H$, defined by
\begin{equation*}
    \begin{array}{llcl}
    g_{n}\colon &\halo H &\longrightarrow &\R\\
    &(\sigma,h)&\longmapsto &f_{n}(h)\cdot \mathds{1}_{\sigma\in L(V_{n})}
    \end{array}
\end{equation*}
with $V_{n}\defeq\bigcup_{s\in S_{H}}{(\supp f_{n})s}$. This naturally provides a lower bound for the $\ell^{p}-$isoperimetric profile of $\halo H$.

The technical part is on the upper bound, and Section~\ref{sec:secondstrategy} provides such a bound in a general situation. In the particular case of lampshufflers, a strategy, well-known to the experts, consists in finding a “good” lamplighter subgroup of $\shuf{H}$, in the sense that this lamplighter should be based on a subgroup $K$ of $H$ which is quasi-isometric to $H$, or at least has the same isoperimetric profile. For this,~\cite[Proposition~2.4]{Sil24} is helpful. We may refer the reader to Appendix~\ref{appendixA} where the aforementioned method is presented and is instructive for the sequel. The upper bound then follows from the monotonicity of the $\ell^{p}-$isoperimetric profile when passing to finitely generated subgroups. Finding such lamplighter subgroups requires some algebraic assumptions on the base group, such as being non perfect or non-co-Hopfian. Such classes of groups provide, at first glance, a nice framework (see Remark~\ref{rem:Preservation} and Proposition~\ref{prop:A7}) and encompass already many classical examples (e.g. all solvable groups).

The goal is to find another strategy for the lampshufflers or other halo products which do not contain a “good” lamplighter as a subgroup. For these specific cases, the idea of finding substructures still remains fruitful. In Section~\ref{sec:secondstrategy}, we therefore make use of the more general notion of lamplighter graphs, that turn out to appear naturally in halo products as subgraphs. Large-scale commutativity will be a key ingredient since we need configurations of lamps to commute. In the particular case of lampshufflers, this novelty has the advantage, compared to Appendix~\ref{appendixA}, of requiring no assumptions on the base group $H$. We then conclude by establishing the monotonicity of the F\o lner function when passing to such subgraphs, in a similar manner as in~\cite[Lemma~4]{Ers03}.

\paragraph{Consequences of Theorem~\ref{th:recallThm intro}.} We first deduce from Theorem~\ref{th:recallThm intro} that, if $\halo H$ is large-scale commutative, is naturally generated and has consistent blocks, then for every $p\ge 1$, its $\ell^{p}-$F\o lner function satisfies
\begin{equation*}
    K^{\text{F\o l}_{H}(x)}\preccurlyeq\folp{\halo H}(x) \preccurlyeq \folp{H}(x)\cdot \Lambda_{\halo H}(C\cdot \folp{H}(x))
\end{equation*}
for some positive constants $C,K>0$.

Now, in terms of isoperimetric profiles, the main result is the following.

\begin{corollaryletter}[see Corollaries~\ref{cor:upperboundonProf} and~\ref{cor:lowerboundonProf}]\label{cor:ProfHalo intro}
Let $p\ge 1$. Let $H$ be a finitely generated amenable group whose $\ell^{p}-$isoperimetric profile $\profp{H}$ satisfies Assumption~$(\star)$. Let $\halo H$ be a naturally generated halo product over $H$ having finite blocks.
\begin{enumerate}[label=(\roman*)]
    \item If $\halo H$ is large-scale commutative, then we have
    \begin{equation*}
            \profp{\halo H}(x)\preccurlyeq\prof{H}(\ln(x)).
    \end{equation*}
    \item If $\halo H$ has consistent blocks, then we have
    \begin{equation*}
                \profp{\halo H}(x) \succcurlyeq \profp{H}(\varphi^{-1}(x))
    \end{equation*}
    where $\varphi(x)=x\cdot \Lambda_{\halo H}(x)$ and where $\Lambda_{\halo H}$ is the lamp growth sequence of $\halo H$.
    \end{enumerate}
\end{corollaryletter}

We then deduce Theorem~\ref{thm:boundsForProfile intro} from this corollary. This result also implies that the isoperimetric profiles of lampjugglers and lampdesigners behave as the isoperimetric profiles of lampshufflers.

\begin{corollaryletter}[see Corollaries~\ref{cor:encadrementduprofildeshuf} and~\ref{cor:profileoflampdesigners}]\label{cor:profilesofJugglersandDesignersOverPolyGrowthGroups}
Theorem~\ref{thm:boundsForProfile intro} also holds for lampjugglers and lampdesigners. Moreover, if $H$ has polynomial growth of degree $d\ge 1$, we have
\begin{equation*}
    \profp{\juggler{s}{H}}(x)\simeq\profp{\designer{H}}(x)\simeq\left(\frac{\ln(x)}{\ln(\ln(x))}\right )^{\frac{1}{d}},
\end{equation*}
for any $s\ge 1$, any real number $p\ge 1$ and any non-trivial finite group $F$, similarly to lampshufflers.
\end{corollaryletter}

Finally, let us also illustrate Corollary~\ref{cor:ProfHalo intro} with lampcloners.

\begin{corollaryletter}[see Corollary~\ref{cor:encadrementduprofildeclone}]\label{cor:encadrementduprofildecloneINTRO}
Let $p\ge 1$. Let $H$ be a finitely generated amenable group, whose $\ell^{p}-$isoperimetric profile $\profp{H}$ satisfies Assumption~$(\star)$. Let $\field$ be a finite field. Then one has 
\begin{equation*}
    \profp{H}\left(\sqrt{\ln(x)}\right) \preccurlyeq \profp{\cloner{H}}(x) \preccurlyeq \prof{H}(\ln(x)).
\end{equation*}
\end{corollaryletter}

Hence, in the setting of the above corollary, we have
\begin{equation*}
    \profp{\cloner{H}}(x)\simeq\prof{H}(\ln(x))
\end{equation*}
when $\profp{H}\left(\sqrt{\ln(x)}\right)\simeq\prof{H}(\ln(x))$, and
\begin{equation}\label{eq:ClonerPolynomial}
    \ln(x)^{\frac{1}{2d}}\preccurlyeq\profp{\cloner{H}}(x)\preccurlyeq \ln(x)^{\frac{1}{d}}
\end{equation}
when $H$ has polynomial growth of degree $d\ge 1$.

In the case of a polynomial growth group $H$, we have in fact the following slight improvement of~\eqref{eq:ClonerPolynomial} for the upper bound:
\begin{equation*}
    \profp{\cloner{H}}(x)\preccurlyeq\left(\frac{\ln(x)}{\ln(\ln(x))}\right)^{\frac{1}{d}},
\end{equation*}
since $\shuf{H}$ is a subgroup of $\cloner{H}$ (consider permutation matrices in $\mathrm{FGL}(H)$).

\paragraph{Applications to regular maps.} We now turn to the problem of the existence of quasi-isometries and regular maps between commonly studied spaces, which has been widely investigated in the literature, see e.g.~\cite{BST12, Tes25, HMT20, HMT22, HMT25, Ben26} among others. In these articles, the main guideline is, mostly, to associate to spaces new quantities that are monotonous under regular maps, and that are finer than the most obvious ones, such as volume growth or asymptotic dimension. As a concrete example, the volume growth does not say anything about the existence of a regular map 
\begin{equation*}
    \mathbb{H}_{\R}^{m_{1}}\times\R^{d_{1}} \longrightarrow \mathbb{H}_{\R}^{m_{2}}\times\R^{d_{2}}
\end{equation*}
whereas Poincaré profiles, introduced and studied in~\cite{HMT20}, impose a monotonic behaviour for the dimension of hyperbolic spaces and the growth exponent of the second factors~\cite[Corollary~1.14]{HMT22}.

On the amenable side, isoperimetric profiles remain powerful invariants to distinguish groups of exponential growth up to quasi-isometry. As an illustration:

\begin{theoremletter}[see Corollary~\ref{cor:IteratedShufflersPolynomialQIBiLip}]\label{thm:IteratedShufflersPolynomialQIBiLip intro}
Let $n,m\ge 0$. Let $A$ and $B$ be infinite virtually abelian finitely generated groups, with growth degrees $a$ and $b$ respectively. Then the following are equivalent:
\begin{enumerate}[label=(\roman*)]
    \item $\shufn{n}{A}$ and $\shufn{m}{B}$ are quasi-isometric.
    \item $n=m$ and $a=b$.
    \item $\shufn{n}{A}$ and $\shufn{m}{B}$ are biLipschitz equivalent.
\end{enumerate}
\end{theoremletter}

\sloppy Two comments are in order here. Firstly, the fact that $\shufn{n}{A}$ and $\shufn{m}{B}$ are quasi-isometric implies that $a=b$ can also be detected with the asymptotic dimension. Indeed, $\shuf{A}$ (more generally $\shufn{n}{A}$) and $A$ have same asymptotic dimension. On the other hand, the asymptotic dimension does not detect numbers of iterations we make, whereas isoperimetric profiles do. These invariants are therefore more powerful in this respect. Additionally, regarding other monotonous quantities under regular maps, volume growth is unhelpful, as it is exponential for both groups when $n,m\ge 1$.

Secondly, it is worth noticing that, for (iterated) lampshufflers over virtually abelian groups or groups with slow profiles (see Corollary~\ref{cor:shufflersoverslowprofiles}), being quasi-isometric is the same as being biLipschitz equivalent. This rigidity is in sharp contrast with lamplighters over $\Z$~\cite{Dym10} or over one-ended groups~\cite{GT24b}, classes in which there are pairs of quasi-isometric groups that are not biLipschitz equivalent.

We refer the reader to Corollary~\ref{cor:IteratedShufflersPolynomialRegularMap} and Remark~\ref{rm:extensionsfornilpotentgroups} for asymmetric versions of Theorem~\ref{thm:IteratedShufflersPolynomialQIBiLip intro}, about the existence of a regular map 
\begin{equation*}
    \shufn{n}{A}\longrightarrow \shufn{m}{B}
\end{equation*}
for polynomial growth groups $A$ and $B$ (not necessarily virtually abelian). A nice consequence of these studies is the following.

\begin{corollaryletter}[see Corollary~\ref{cor:QIRegularMap}]\label{cor:corV.J}
Let $n,m\ge 0$. Let $A$ and $B$ be infinite virtually abelian finitely generated groups, with growth degrees $a$ and $b$ respectively. Then the following are equivalent:
\begin{enumerate}[label=(\roman*)]
    \item the three equivalent assertions of Theorem~\ref{thm:IteratedShufflersPolynomialQIBiLip intro} hold;
    \item \sloppy $\shufn{n}{A}$ regularly embeds $\shufn{m}{B}$, and $\shufn{m}{B}$ regularly embeds into $\shufn{n}{A}$.
\end{enumerate}
\end{corollaryletter}

Another interesting consequence of our computations of isoperimetric profiles is the following statement, which cannot be reached with methods from~\cite{GT24a}, even for quasi-isometric or coarse embeddings. Indeed, in the latter is introduced a key property, called the~\textit{thick bigon property} (cf. Definition~\ref{def:TBP}), which is a crucial assumption for the study of quasi-isometries between halo products. Unfortunately, this property is not stable under iterations of lampshufflers, and cannot be used for $\shufn{n}{\Z^d}$ for instance.

\begin{propositionletter}[see Corollary~\ref{cor:iteratedshufintoiteratedlamplighterPolynomial}]\label{prop:IteratedShufflerintoIteratedLL}
Let $d,k,n\ge 1$ be three integers. If there exists a regular map 
\begin{equation*}
\shufn{n}{\Z^d}\longrightarrow \Z/2\Z\wr\big(\Z/2\Z\wr(\dots(\Z/2\Z\wr \Z^k))\big)
\end{equation*}
where the wreath product is iterated $n$ times, then $d<k$. 
\end{propositionletter}

In relation with the results from~\cite{GT24a}, we expect in fact that there is no regular map from $\shufn{n}{\Z^d}$ to $\Z/2\Z\wr\big(\Z/2\Z\wr(\dots(\Z/2\Z\wr \Z^k))\big)$, even when $d<k$. 

Finally, we emphasize here that similar results can be obtained for other halo products, such as lampjugglers, lampdesigners, and lampcloners (except when the base group has polynomial growth since we do not have a precise estimate of $\prof{\cloner{H}}$ in this case).

\paragraph{Plan of the chapter.} After a few preliminaries in Section~\ref{sec:preliminaries}, we introduce halo products in Section~\ref{sec:halo}, with the main assumptions we will need to study them. Section~\ref{sec:computationsFolnerFunction} is devoted to the computation of F\o lner functions for halo products, and we deduce estimates for isoperimetric profiles in Section~\ref{sec:computationsIsoProf}. This finally implies existence and non-existence results of regular maps and quasi-isometries between such groups, see Section~\ref{sec:appQIandRegularMaps}. Lastly, Appendix~\ref{appendixA} presents various minimal algebraic assumptions under which lamplighters appear as subgroups of lampshufflers.

\section{Notations and terminologies}\label{sec:preliminaries}

\subsection{Notations}\label{sec:notations} 

Given non-decreasing functions $f,g\colon \R_{>0}\rightarrow \R_{>0}$, we write $f(x)=O(g(x))$ if there exists $C>0$ such that $f(x)\le Cg(x)$ for all $x$ large enough, and $f(x)=o(g(x))$ if $\frac{f(x)}{g(x)}$ goes to $0$ as $x$ goes to $+\infty$. We write $f\sim g$, and we say that $f$ and $g$ are~\textit{equivalent}, if $\frac{f(x)}{g(x)}$ goes to $1$ as $x$ goes to $+\infty$.  

Recall that a map $f\colon (X,d_{X})\longrightarrow (Y,d_{Y})$ between two metric spaces is called~\textit{regular} if it is $C-$Lipschitz for some $C>0$ and pre-images of points have uniformly bounded cardinality: there is $m\ge 1$ such that 
$|f^{-1}(\lbrace y\rbrace)| \le m$ for any $y\in Y$. 

Note that any quasi-isometry is a quasi-isometric embedding, which is itself a coarse embedding, which is itself a regular map, but none of the reverse implications hold. For instance, the inclusion of a closed compactly generated subgroup in a locally compact compactly generated group is always a coarse embedding, while it is a quasi-isometry only if the subgroup is undistorted, and the map $\Z\rightarrow\Z$, $n\longmapsto |n|$, is a regular map, while it is not a coarse embedding.

\subsection{Isoperimetric profiles} 

\paragraph{Isoperimetric profile for groups.} For a finitely generated group $G$ and a finite generating set $S_{G}$, its $\ell^{p}-$\textit{isoperimetric profile}, for $p\ge 1$, is the function $\profp{G}\colon \N\rightarrow \R_{+}$ given by
\begin{equation*}
    \profp{G}(n)\defeq\sup_{\substack{f\colon G\to\R_+\\|\supp{f}|\leq n}}{\frac{\|f\|_{p}}{\|\nabla f\|_{p}}}
\end{equation*}
where the support of $f\colon G\rightarrow\R_{+}$ is $\supp{f}\defeq\lbrace g\in G : f(g)\not=0\rbrace$ and the $\ell^{p}-$norm of its gradient is defined by 
\begin{equation*}
    \|\nabla f\|_{p}^{p}\defeq\sum_{g\in G,\; s\in S_{G}}{|f(g)-f(gs)|^p}.
\end{equation*}
For $p=1$, the $\ell^{1}-$isoperimetric profile is simply called~\textit{isoperimetric profile} and one has
\begin{equation*}
    \prof{G}(n) \simeq \sup_{|A|\le n}\frac{|A|}{|\partial_{G} A|}
\end{equation*}
where $\partial_{G}A \defeq AS_{G}\setminus A=\lbrace g\in G\setminus A : \exists s\in S_{G}, \exists h\in A, g=hs\rbrace$ is the~\textit{boundary} of $A$ in $G$.

Recall also that the isoperimetric profile of a group $G$ is the generalized inverse of its~\textit{F\o lner function} $\text{F\o l}_{G}\colon \N\rightarrow\R_{+}$, defined as
 \begin{equation*}
     \text{F\o l}_{G}(n) \defeq \inf\left\lbrace |A| : \frac{|\partial_{G}A|}{|A|} \le \frac{1}{n} \right\rbrace.
 \end{equation*}

We more generally define the $\ell^{p}-$F\o lner function $\folp{G}\colon\N\rightarrow \R_{+}$ for every $p\ge 1$, as 
\begin{equation*}
     \folp{G}(n) \defeq \inf\left\lbrace \left|\supp{f}\right| : \frac{\|\nabla f\|_p}{\|f\|_p} \le \frac{1}{n} \right\rbrace.
\end{equation*}
For every $p\geq 1$, $\folp{G}$ and $\profp{G}$ are generalized inverses of each other, and we have $\text{F\o l}_{1,G}(x)\simeq\text{F\o l}_{G}(x)$. Thus, in the sequel, we will always write $\text{F\o l}_{G}$ instead of $\text{F\o l}_{1,G}$.

\begin{remark}\label{rem:InverseEachOther}
Notice that, given the asymptotic behaviour of the $\ell^{p}-$F\o lner function, we can deduce the asymptotic behaviour of the $\ell^{p}-$isoperimetric profile, even though they are not real inverses of each other, but only generalized inverses~\textit{a priori}. Indeed, a non-decreasing function $\R_{+}\rightarrow \R_{+}$ is always asymptotically equivalent to an increasing function $\R_{+}\rightarrow \R_{+}$ (cf.~\cite[Remark~1.2]{Cor25}) and it is not hard to check that $\simeq$ is preserved when passing to generalized inverses. Thus, in the sequel, we can and will assume that the $\ell^{p}-$F\o lner function and the $\ell^{p}-$isoperimetric profile are injective and then real inverses of each other. Hence, studying the asymptotic behaviour of $\profp{G}$ is the same as studying the asymptotic behaviour of $\folp{G}$. 
\end{remark}

Note that a group is amenable if and only if its isoperimetric profile is unbounded. The idea to keep in mind is that the isoperimetric profile is a measurement of how much amenable a group is. The faster the isoperimetric profile tends to infinity, the more the group is amenable. In particular, the isoperimetric profile is a particularly well suited invariant to distinguish amenable groups with exponential growth up to quasi-isometries or regular maps, and it has now been computed for many classes of groups, among which:
\begin{itemize}
    \item $\profp{G}(n) \simeq n^{\frac{1}{d}}$ if $G$ has polynomial growth of degree $d\ge 1$;
    \item $\profp{G}(n)\simeq \ln(n)$ for solvable Baumslag-Solitar groups and lamplighters $F\wr\Z$, where $F$ is a non-trivial finite group;
    \item $\profp{G}(n) \simeq \ln(n)$ for any polycyclic group with exponential growth~\cite{Pit95, Pit00}, or more generally any exponential growth group within the class GES of Tessera~\cite[Corollary 5]{Tes13};
    \item $\profp{F\wr N}(n)\simeq \ln(n)^{\frac{1}{d}}$ for $F$ a non-trivial finite group, and $N$ having polynomial growth of degree $d\ge 1$~\cite{Ers03}.
\end{itemize}

More generally,~\cite[Theorem 1]{Ers03} provides an explicit formula for computing F\o lner functions of wreath products, so for instance the last example can be extended to iterated wreath products. Let us explain with more details Erschler's result, since it will play an important role in the sequel.

Let $G$ and $H$ be finitely generated groups and assume that for every $C>0$, there exists $k>0$ such that $\text{F\o l}_H(kn)>C\cdot \text{F\o l}(n)$ for any large enough integer $n$. Then we have
\begin{equation*}
    \text{F\o l}_{G\wr H}(x)\simeq\text{F\o l}_{G}(x)^{\text{F\o l}_{H}(x)}.
\end{equation*}
This assumption on $\text{F\o l}_{H}$, introduced in~\cite{Ers03}, can also be stated in terms of isoperimetric profile in the following way: for every $C>0$, $\prof{H}(Cx)=O\left(\prof{H}(x)\right)$. This assumption also appeared in~\cite{Cor25} and we do not have any example of a finitely generated group for which it does not hold. In Section~\ref{sec:computationsIsoProf}, we will make use of this mild assumption, and we will refer to it as \textit{Assumption}~($\star$). In~\cite{Ers03}, this assumption is used to get rid of some constant appearing in the lower bound: there exists $C>0$ such that
\begin{equation}\label{eq:ConstantAssumptionStar}
    \text{F\o l}_{G\wr H}(x)\succcurlyeq\left (\text{F\o l}_{G}(x)\right )^{C\text{F\o l}_{H}(x)};
\end{equation}
whereas the upper bound is exactly $\text{F\o l}_{G\wr H}(x)\preccurlyeq \text{F\o l}_{G}(x)^{\text{F\o l}_{H}(x)}$. In the particular case of a non-trivial finite group $G$, we have
\begin{equation*}
    \prof{G\wr H}(x)\simeq\prof{H}(\ln(x))
\end{equation*}
if $\prof{H}$ satisfies Assumption~($\star)$. The lower bound~\eqref{eq:ConstantAssumptionStar} also holds in the context of lamplighter graphs, see Section~\ref{sec:secondstrategy}. To relate the work of Erschler on the $\ell^{1}-$F\o lner function with the $\ell^{p}-$F\o lner function (for $p\ge 1$) that we want to compute for halo products, we will first have to reduce the proof to the case $p=1$, thanks to the well-known fact that if $p\ge 1$, then the $\ell^{p}-$F\o lner function dominates the $\ell^{1}-$F\o lner function (see Lemma~\ref{lem:MonotonuousProfile}). In fact, it is conjectured that the $\ell^{p}-$F\o lner functions, for $p\ge 1$, all have the same asymptotic behaviour.

A fundamental result in geometric group theory is the one of Coulhon and Saloff-Coste, who proved in~\cite{CSC93} that for a finitely generated group $G$ and a finite symmetric generating set $S$ of $G$, one has 
\begin{equation*}
    \frac{|\partial_{G}F|}{|F|} \ge \frac{1}{4|S|}\cdot\frac{1}{\Phi_{G}(2|F|)}
\end{equation*}
for any finite set $F\subset G$, where $\Phi_{G}\colon \R_{>0} \rightarrow \N$, $\Phi_{G,S}(t) \defeq \min\lbrace n \ge 0 : \beta_{G,S}(n)>t\rbrace$ is the~\textit{inverse growth function} of $G$. Since then, it has constantly been improved to finer inequalities, see for instance~\cite{PS22}. Inverting this inequality and taking the sup, one directly gets the upper bound 
\begin{equation*}
    \prof{G}(n) \preccurlyeq \Phi_{G,S}(n)
\end{equation*}
on the $\ell^{1}-$isoperimetric profile of $G$. This upper bound is optimal when $G$ has polynomial growth, while if it has exponential growth, one only gets $\prof{G}(n) \preccurlyeq \ln(n)$. In fact,~\cite[Theorem~1.1]{BZ21} describes a large class of possible asymptotic behaviours for isoperimetric profiles of finitely generated groups with exponential growth, namely for any non-decreasing function $f$ such that $x\longmapsto \frac{x}{f(x)}$ is non-decreasing, there exists a finitely generated group of exponential growth whose $\ell^{p}-$isoperimetric profile is $\simeq \frac{\ln(x)}{f(\ln(x))}$ for every $p\ge 1$.

Isoperimetric profiles are also particularly studied for their relations with return probabilities of random walks on groups, see for instance~\cite{SCZ15, SCZ16, SCZ18, BZ21}.

\paragraph{Isoperimetric profile for graphs.} Isoperimetric profiles can be defined, more generally, in the framework of bounded degree graphs, without any underlying algebraic structure. In this paper, we only focus on the case $p=1$, but similar definitions can be made for $p>1$.

Since there will be no ambiguity, we will abusively use the same notation for a graph and the set of its vertices. Given a graph $Y$, the presence of an edge between two vertices $v$ and $w$ will be denoted by $v\sim_{Y}w$.

Given a graph $Y$, its isoperimetric profile is the function $\prof{Y}\colon\N\rightarrow\R_{+}$ defined by
\begin{equation*}
    \prof{Y}(n)\defeq\sup_{|A|\le n}{\frac{|A|}{|\partial_{Y} A|}}
\end{equation*}
where, given a finite set $A\subset Y$ of vertices, $\partial_{Y}A\defeq\lbrace v\in Y\setminus A : \exists a\in A, v\sim_{Y} a\rbrace$ is the~\textit{boundary} of $A$ in the graph $Y$.

In the case of Cayley graphs of finitely generated groups, we recover the corresponding notion of isoperimetric profile of groups, defined above. Moreover, the invariance of isoperimetric profile under quasi-isometry is still valid in this more general setup. Finally, we analogously define the F\o lner function of a graph.

This setup of graphs will be crucial in our paper. Indeed, in Section~\ref{sec:secondstrategy}, we will define lamplighter graphs and will require a lower bound of their F\o lner functions.

\section{Halo products}\label{sec:halo}

In this section, we define halo products and the main classes we are interested in.

\subsection{Halo products: definition and main examples}\label{sec:defHalo} 

\begin{definition}
Let $X$ be a set. A~\textit{halo of groups $\halo$ over $X$} is the data, for any subset $S\subset X$, of a group $L(S)$ such that:
\begin{itemize}
    \item for all $R,S\subset X$, if $R\subset S$ then $L(R)\leqslant L(S)$;
    \item $L(\emptyset)=\lbrace 1\rbrace$ and $L(X)=\langle L(S) : S\subset X \;\text{finite}\rangle$;
    \item for all $R,S\subset X$, $L(R\cap S)=L(R)\cap L(S)$.
\end{itemize}
\end{definition}

Given an action $H\curvearrowright X$ and a morphism $\alpha\colon H \longrightarrow \text{Aut}(L(X))$ satisfying $\alpha(h)(L(S))=L(hS)$ for any $S\subset X$ and $h\in H$, the~\textit{permutational halo product} $\halo_{X,\alpha}H$ is the semi-direct product 
\begin{equation*}
    \halo_{X,\alpha}H \defeq L(X)\rtimes_{\alpha}H.
\end{equation*}

The definition is motivated by permutational wreath products, which are basic examples of permutational halo products. Indeed, given groups $F,H$ and an action $H\curvearrowright X$, set $L(S)\defeq \bigoplus_{S}F$ for any $S\subset X$. Then $\halo_{X,\alpha}H$ coincides with $F\wr_{X}H$, where $\alpha$ is the action of $H$ on $\bigoplus_{X}F$ obtained by permuting the coordinates through the initial action $H\curvearrowright X$. In particular, for $X=H$ and the left-multiplication action of $H$ on itself, we recover a description of the wreath product $F\wr H$ as a halo product. 

Let us now describe other examples of halo products. From now on, we only focus on halo products with $X=H$, that we simply denote by $\halo H$, for the natural action of $H$ on itself by left-multiplication.

\paragraph{Lampshufflers.} Let $H$ be a group, and let $\fsym{H}$ be the group of~\textit{finitely supported} permutations of $H$, that is the group of bijections $H\rightarrow H$ that are the identity outside a finite subset of $H$. The group $H$ acts naturally on $\fsym{H}$, via 
\begin{equation*}
    (h\cdot\sigma)(x) \defeq h\sigma(h^{-1}x), \; x\in H
\end{equation*}
for any $h\in H$ and $\sigma\in\fsym{H}$. Indeed, if $\sigma\colon H\rightarrow H$ is a finitely supported bijection and $h\in H$, then so is $h\cdot \sigma$ and $\supp{(h\cdot\sigma)}=h\cdot\supp{(\sigma)}$, where $\supp{(\sigma)} \defeq\lbrace x\in H : \sigma(x) \neq x\rbrace$.

The~\textit{lampshuffler group over $H$}, denoted $\shuf{H}$, is then defined as the semidirect product 
\begin{equation*}
    \shuf{H}\defeq \fsym{H} \rtimes H.
\end{equation*}

It coincides with the halo product $\halo H$ where $L(S)\defeq \fsym{S}$ for any $S\subset H$. Additionally, if $H$ is finitely generated and $S_{H}$ denotes a finite generating set, then $\shuf{H}$ is generated by the finite set 
\begin{equation*}
    \Sigma_{H}\defeq\left\lbrace(\tau_{1_{H},s},1_{H}) : s\in S_{H}\right \rbrace\cup\left\lbrace(\text{id},s) : s\in S_{H}\right\rbrace
\end{equation*}
where, given any $x,y\in H$, $\tau_{x,y}\in\fsym{H}$ is the transposition that swaps $x$ and $y$, that is $\tau_{x,y}(x)=y$, $\tau_{x,y}(y)=x$ and $\tau_{x,y}(h)=h$ for any $h\notin\lbrace x,y\rbrace$. 

An element $(\sigma,h)\in\shuf{H}$ can be seen as a labelling of the vertices of the Cayley graph $\text{Cay}(H,S_{H})$ (a vertex $p\in H$ carries the label $\sigma(p)$) together with an arrow pointing at some vertex $h\in H$, and there are two types of moves in $\text{Cay}(\shuf{H},\Sigma_{H})$ to go from $(\sigma, h)$ to a neighbouring vertex:
\begin{itemize}
    \item either the arrow goes from $h$ to a neighbouring vertex in $H$;
    \item or the arrow stands on the vertex $h\in H$, and swaps its label with the label of one of its neighbours in $H$.
\end{itemize} 

\paragraph{Lampjugglers.} Lampshufflers are in fact particular instances of a broader family of groups, called~\textit{lampjugglers}. Given a group $H$ and an integer $r\ge 1$, the~\textit{lampjuggler over $H$} is the semi-direct product
\begin{equation*}
    \juggler{r}{H} \defeq \fsym{H\times\lbrace 1,\dots,r\rbrace} \rtimes H
\end{equation*}
where $H$ acts on $\fsym{H\times\lbrace 1,\dots,r\rbrace}$ through its initial action on $H\times\lbrace 1,\dots,r\rbrace$ given by $h\cdot (x,i) \defeq (hx, i)$. It can be described as the halo group $\halo H$ where 
\begin{equation*}
    L(S)\defeq \fsym{S\times\lbrace 1,\dots, r\rbrace},\;S\subset H.
\end{equation*}
As for lampshufflers, lampjugglers over finitely generated groups are finitely generated, and one can check that if $S_{H}$ is a finite generating set for $H$, then the finite set
\begin{equation*}
    \lbrace (\tau_{(1_{H},i),(s,j)}, 1_{H}) : s\in S_{H}, 1\le i,j\le r\rbrace  \cup \lbrace (\text{id},s) : s\in S_{H}\rbrace
\end{equation*}
generates $\juggler{r}{H}$. Here, an element $(\sigma, h)\in\juggler{r}{H}$ can be seen as a labelling of the vertices of $\text{Cay}(H,S_{H})\times\lbrace 1,\dots,r\rbrace$ together with an arrow pointing at some vertex $h\in H$. Right-multiplying $(\sigma, h)$ by a generator from the above set amounts either to move the arrow from $h$ to a neighbouring vertex $hs$ in $H$, or to keep the arrow on $h\in H$ and switching the labels of two vertices in $h \times\lbrace 1, \dots, r\rbrace$ and $hs \times \lbrace1, \dots, r\rbrace$ for some neighbour $hs$ of $h$.

\paragraph{Lampdesigners.} Let $F$ and $H$ be two groups. The~\textit{lampdesigner over $H$} is the semi-direct product 
\begin{equation*}
    \designer{H} \defeq (F\wr_{H}\fsym{H})\rtimes H
\end{equation*}
where $H$ acts on $\bigoplus_{H}F$ by permuting the coordinates through its action on itself by left-multiplication and acts on $\fsym{H}$ as described above. It is the halo product $\halo H$ for the collection $L(S)\defeq F\wr_{S}\fsym{S}$, $S\subset H$. 

Lampdesigners are close to lampjuggler groups, and in fact if $F$ is finite, $\designer{H}$ is a subgroup of $\juggler{|F|}{H}$, via the map 
\begin{align*}
    \begin{array}{cll}
    \designer{H} &\longrightarrow &\juggler{|F|}{H} \\
    ((f,\sigma),h) &\longmapsto &(\sigma', h)
    \end{array}
\end{align*}
where, given a pair $(f,\sigma)\in F\wr_{H}\fsym{H}$, $\sigma'$ is the permutation of $H\times F$ given by $\sigma'(h,i)=(\sigma(h), f(h)i)$. Note also that $\designer{H}$ contains $\shuf{H}$ as a subgroup. 

\paragraph{Lampcloners.} Let $H$ be a group and let $\field$ be a field. Denote $V_{H}$ the $\field$-vector space admitting $H$ as a basis, and denote $\lbrace e_{u} : u\in H\rbrace$ a formal basis. Let $\text{FGL}(H)$ be the group of linear automorphisms $V_{H}\rightarrow V_{H}$ that fix all but finitely many basis elements. This group can also be seen as the group of finitely supported invertible matrices with coefficients in $\field$ whose entries are indexed by $H\times H$. Once again, the action of $H$ on itself naturally yields an action of $H$ on $\text{FGL}(H)$. The~\textit{lampcloner over $H$} is the semi-direct product 
\begin{equation*}
    \cloner{H} \defeq \text{FGL}(H)\rtimes H.
\end{equation*}
It is a halo product, for the collection $L(S)\defeq \text{FGL}(S)$, for every $S\subset H$, where $\text{FGL}(S)$ is thought of as the subgroup of $\text{FGL}(H)$ of linear automorphisms $V_{H}\rightarrow V_{H}$ that fix $H\setminus S$ and that stabilise the subspace $\langle S\rangle\subset V_{H}$. 

In addition, if $\field$ is finite and if $H=\langle S_{H}\rangle$ is finitely generated, then the finite set
\begin{equation*}
    \lbrace (\delta_{1_{H}}(\lambda), 1_{H}) : \lambda \in \field\setminus\lbrace 0\rbrace\rbrace\cup\lbrace (\tau_{1_{H},s}(\lambda),1_{H}) : s\in S_{H}, \lambda \in \field\setminus\lbrace 0\rbrace\rbrace \cup\lbrace (\text{id}_{V_H}, s) : s\in S_{H}\rbrace
\end{equation*}
generates $\cloner{H}$, where, given $p,q\in H$ and $\lambda\in \field\setminus\lbrace 0\rbrace$, $\delta_{p}(\lambda)$ is the~\textit{diagonal matrix} 
\begin{align*}
    \delta_{p}(\lambda)\colon
    \begin{array}{cll}
    V_{H}&\longrightarrow & V_{H}\\
    \displaystyle \sum_{h\in H}\mu_{h}e_{h} &\longmapsto &\displaystyle \sum_{h\neq p}\mu_{h}e_{h}+\lambda\mu_{p}e_{p}
    \end{array}
\end{align*}
and $\tau_{pq}(\lambda)$ is the \textit{transvection} 
\begin{align*}
    \tau_{pq}(\lambda)\colon
    \begin{array}{cll}
    V_{H}&\longrightarrow &V_{H} \\
    \displaystyle \sum_{h\in H}\mu_{h}e_{h} &\longmapsto &\displaystyle \sum_{h\neq p}\mu_{h}e_{h}+(\mu_{p}+\lambda\mu_{q})e_{p}
    \end{array}.
\end{align*}
Thus, thinking of an element $(\varphi,p)\in\cloner{H}$ as a labelling of $\text{Cay}(H,S_{H})$ (the vertex $h\in H$ has the label $\varphi(e_{h})\in V_{H}$), together with an arrow pointing at $p\in H$, right multiplying $(\varphi,p)$ by a generator from the above set amounts either to move the arrow to an adjacent vertex $q$ of $p$ in $H$; or to keep the arrow where it stands and multiply $\varphi(e_{p})$ by a non-trivial element of $\field$; or to keep the arrow where it stands and to~\textit{clone} the label $\varphi(e_{p})$ and add it to the label of a neighbour of $p$ after multiplication by an element of $\field\setminus\lbrace 0\rbrace$. 

We refer the reader to~\cite[Section~2]{GT24a} for many other possible constructions, such as lampbraiders and verbal halo products, that encompass for instance nilpotent and metabelian wreath products.

\subsection{Important assumptions}\label{sec:ImportantAssumption} 

Our main results deal with halo products satisfying various important assumptions that we introduce in this section.

\paragraph{Large-scale commutativity.} The first one has been introduced in~\cite{GT24a}, under the terminology~\textit{large-scale commutativity}.
\begin{definition}
Let $\halo H$ be a halo product over a finitely generated group $H$, and let $S_{H}$ be a finite generating set of $H$. We say that $\halo H$ is~\textit{large-scale commutative} if there exists a constant $D\ge 0$ such that, for any $R,S\subset H$ with $d_{S_{H}}(R,S)\ge D$, the subgroups $L(R)$ and $L(S)$ commute in $L(H)$.
\end{definition}

This notion plays a key role in the quasi-isometry classification of halo groups established in~\cite{GT24a}, see e.g.~\cite[Theorem~6.3 and Theorem~6.6]{GT24a}. Examples of large-scale commutative halo products include lamplighters ($D=0$), lampshufflers ($D=1$), lampcloners ($D=1$).

\paragraph{Finite generating sets.} Let us now turn to terminologies more specific to halo groups over finitely generated groups. Inspired by lampshufflers, when we are looking for a generating set of a general halo product, there is a natural candidate. The~\textit{natural generation property}, that we now introduce, is by definition satisfied by a halo product having
this natural candidate as generating set.

\begin{definition}
Let $H$ be a finitely generated group, with a finite generating subset $S_{H}$. We say that a halo product $\halo H$ over $H$ is \textit{naturally generated} if it is generated by the set
\begin{equation*}
\lbrace (1_{L(H)},s) : s\in S_{H}\rbrace \cup \bigcup_{s\in S_{H}}\big\lbrace (\sigma, 1_{H})\in\halo H : \sigma\in L(\lbrace 1_{H},s\rbrace)\big\rbrace.
\end{equation*}
\end{definition}

We already know from the previous section that our running examples are naturally generated.

Let us now introduce a terminology relative to the generation for blocks of a halo product. 

\begin{definition}
We say that a halo product $\halo H$ has~\textit{finite} (resp.~\textit{finitely generated}) blocks if, for any finite subset $S\subset H$, $L(S)$ is finite (resp. finitely generated).
\end{definition}

For instance, a wreath product $F\wr H$, where $F$ is finitely generated, has finitely generated blocks. Moreover, lampjugglers, lampdesigners and lampcloners have finite blocks, so they have finitely generated blocks.

If a naturally generated halo product has finitely generated blocks, then it has a natural finite generating set:

\begin{fact}\label{fact:BlocksHaloFG}
Let $H$ be a finitely generated group and let $S_{H}$ be a finite symmetric generating set of $H$. Let $\halo H$ be a halo product over $H$. Suppose that $\halo H$ is naturally generated and has finitely generated blocks. Then the finite set
\begin{equation*}
    S_{\halo H}\defeq\lbrace (1_{L(H)},s) : s\in S_{H}\rbrace \cup \bigcup_{s\in S_{H}}\lbrace (\sigma, 1_{H})\in\halo H : \sigma\in S(s)\rbrace
\end{equation*}
generates $\halo H$, where $S(s)$ is any finite generating subset of $L(\lbrace 1_{H},s\rbrace)$.\qed
\end{fact}

\paragraph{Growth of lamps.} Recall from Definitions~\ref{def:consistencyforhalos} and~\ref{def:Growthoflamps} that if a halo product $\halo H$ is consistent, then we may define its~\textit{lamp growth sequence} as  
\begin{equation*}
    \Lambda_{\halo H}\colon n\longmapsto |L(S)|, \;\text{where}\; |S|=n.
\end{equation*}
It has been computed in~\cite[Facts~7.12-7.16]{GT24a} for many halo products, such as lamplighters, lampshufflers, lampdesigners and lampcloners. For instance, for any group $H$ and finite group $F$, one has $\Lambda_{F\wr H}(n)=|F|^{n}$, $\Lambda_{\juggler{r}{H}}(n)=(rn)!$ and $\Lambda_{\designer{H}}(n)=|F|^{n}n!$. Moreover, given a finite field $\field$, we have
\begin{equation*}
    \Lambda_{\cloner{H}}(n)=\prod_{i=0}^{n-1}{(|\field|^n-|\field|^i)}.
\end{equation*}

In fact, it turns out that the asymptotic behaviour of this sequence is invariant under a special class of quasi-isometries (and more generally coarse embeddings), referred to as~\textit{aptolic quasi-isometries} in~\cite[Section~6]{GT24a}. As proved in~\cite[Corollary~6.11]{GT24a}, under additional assumptions, any quasi-isometry between two halo products is aptolic (up to finite distance). Thus, for these halo products, the asymptotic behaviour of the lamp growth sequence is an invariant of quasi-isometry. 

\section{Estimates of F\o lner functions of halo products}\label{sec:computationsFolnerFunction}

\subsection{A general upper bound on the F\o lner functions: finding almost invariant functions}\label{sec:UpperBoundFolner} 

In this section, we provide an upper bound on the $\ell^{p}-$F\o lner function of many halo products. The $\ell^{p}-$F\o lner function of a finitely generated group $G$ being an infimum over finitely supported functions $G\rightarrow\R$, the strategy is to exhibit “good” such functions, namely almost invariant functions. This constitutes the first step towards proving Theorem~\ref{th:recallThm intro}.

Observe that, if $H$ is an amenable group, then $\halo H$ is amenable if and only if $L(H)$ is amenable, and the latter is often true regardless of $H$. For instance, if blocks are finite, then $L(H)$ is locally finite and thus amenable (cf. Corollary~\ref{cor:amenabilitybyfgsubgroups}). This is the case when $\halo H$ is a lamplighter $F\wr H$ (i.e. $F$ is finite), a lampshuffler $\shuf{H}$, or a lampcloner $\cloner{H}$ over a finite field $\field$.

Therefore, we know that almost invariant functions exist when our halo product has finite blocks. Here our goal is, in particular, to construct a suitable sequence of almost invariant functions for our halo product from such a sequence for the base group $H$. The estimates we can derive enable us to prove the following.

\begin{proposition}\label{prop:upperboundonFol}
Let $H$ be a finitely generated amenable group, and let $\halo H$ be a halo product over $H$. Suppose that $\halo H$ is naturally generated and has consistent blocks. Then, for any $p\ge 1$, there exists a constant $C>0$ such that  
\begin{equation*}
    \folp{\halo H}(x)\preccurlyeq\folp{H}(x)\cdot \Lambda_{\halo H}(C\cdot\folp{H}(x)),
\end{equation*}
where $\Lambda_{\halo H}$ is the lamp growth sequence of $\halo H$. 
\end{proposition}

\begin{proof}
As usual, denote $S_{H}$ a finite generating set for $H$. By assumption, the subset
\begin{equation*}
    S_{\halo H}=\left\lbrace (1_{L(H)},s) : s\in S_H\right\rbrace \cup\bigcup_{s\in S_{H}} \big\lbrace (\sigma_{s},1_{H}) : \sigma_{s}\in L(\lbrace 1_{H},s\rbrace)\big\rbrace
\end{equation*}
generates $\halo H$. Let $(f_{n})_{n\ge 0}$ be a sequence of functions $H\rightarrow \R$ that realizes $\folp{H}$, i.e. $\folp{H}(n)=|\supp f_{n}|$ and $\frac{\|\nabla_{S_H}f_{n}\|_{p}}{\|f_{n}\|_{p}} \le \frac{1}{n}$ for any $n\ge 0$. Given $n\ge 0$, set 
\begin{equation*}
    U_{n}\defeq\supp{f_{n}},\; V_{n}\defeq\bigcup_{s\in S_{H}}{U_{n}s}
\end{equation*}
and
\begin{equation*}
    \begin{array}{llcl}
    g_{n}\colon &\halo H &\longrightarrow &\R\\
    &(\sigma,h)&\longmapsto &f_{n}(h)\mathds{1}_{\sigma\in L(V_n)}
    \end{array}.
\end{equation*}
Let $(\sigma,h)\in\halo H$, $s\in S_{H}$ and $\sigma_{s}\in L(\lbrace 1,s\rbrace)$. The composition law of $\halo H$ directly implies that $(\sigma,h)(1_{L(H)},s)=(\sigma,hs)$ and $(\sigma,h)(\sigma_{s},1_{H})=(\sigma(h\cdot\sigma_{s}),h)$. This implies that 
\begin{equation*}
    g_{n}\left ((\sigma,h)(1_{L(H)},s)\right )-g_{n}(\sigma,h)=(f_{n}(hs)-f_{n}(h))\mathds{1}_{\sigma\in L(V_n)}
\end{equation*}
as well as
\begin{equation*}
    g_{n}\big((\sigma,h)(\sigma_{s},1_{H})\big)-g_{n}(\sigma,h)=0
\end{equation*}
using that $\sigma\in L(V_{n})$ if and only if $ \sigma(h\cdot\sigma_s)\in L(V_{n})$ when $h\in U_{n}$. We thus have
\begin{align*}
    \|\nabla_{S_{\halo H}}g_n\|_{p}^{p}&=\sum_{(\sigma,h)\in\halo H}{\sum_{s\in S_{H}}{\left|g_{n}\left((\sigma,h)(1_{L(H)},s)\right)-g_{n}(\sigma,h)\right|^{p}}}\\
    &=\sum_{(\sigma,h)\in\halo H}{\sum_{s\in S_{H}}{|(f_{n}(hs)-f_{n}(h))\mathds{1}_{\sigma\in L(V_{n})}|^p}}\\
    &=|L(V_{n})|\cdot\|\nabla_{S_H}f_{n}\|_{p}^{p}.
\end{align*}
We also have
\begin{equation*}
    \|g_{n}\|_{p}^{p}=|L(V_{n})|\cdot\|f_{n}\|_{p}^{p}.
\end{equation*}
We finally get
\begin{equation*}
\frac{\|\nabla_{S_{\halo H}}g_{n}\|_{p}}{\|g_{n}\|_{p}}=\frac{\|\nabla_{S_{H}}f_{n}\|_{p}}{\|f_{n}\|_{p}}\le\frac{1}{n}
\end{equation*}
so, by the definition of the $\ell^{p}-$F\o lner function, it follows that
\begin{align*}
    \folp{\halo H}(n)&\le |\supp g_{n}|=|L(V_{n})|\cdot |U_{n}|\\
    &=|U_{n}|\cdot \Lambda_{\halo H}(|V_{n}|)\\
    &\le \folp{H}(n)\cdot\Lambda_{\halo H}\left(|S_{H}|\cdot \folp{H}(n)\right).
\end{align*}
This concludes the proof.
\end{proof}

\subsection{A general lower bound on the F\o lner functions: finding lamplighter subgraphs}\label{sec:secondstrategy}

The goal of this section is to find a lower bound of the $\ell^{p}-$F\o lner function of a halo product. In Section~\ref{sec:computationsIsoProf}, we deduce, for specific cases, an upper bound of the $\ell^{p}-$isoperimetric profile which will often be optimal.

Here is then the general statement. 

\begin{theorem}\label{th:lowerboundonFol}
Let $H$ be a finitely generated amenable group and let $S_{H}$ be a finite generating set. Let $\halo H$ be a naturally generated and large-scale commutative halo product having finitely generated blocks. Then,
for any $p\ge 1$ and any $s_{0}\in S_{H}$, there exists a constant $C>0$ such that
\begin{equation*}
    \folp{\halo H}(x) \succcurlyeq \left (\text{F\o l}_{L(\{1_H,s_0\})}(x)\right )^{C\cdot\text{F\o l}_H(x)}.
\end{equation*}
\end{theorem}

In the particular case of a halo product with finite blocks, we thus get
\begin{equation*}
    \folp{\halo H}(x) \succcurlyeq K^{\text{F\o l}_H(x)}
\end{equation*}
for some positive constant $K>0$ and any $p\ge 1$.

The following lemma will allow us to reduce to the case $p=1$. This is a well-known result on isoperimetric profiles, mentioned in~\cite{Cou00}, stating that for every finitely generated group $G$, $\profp{G}$ is monotonous in the variable $p\ge 1$ for the order given by $\preccurlyeq$. Here we state it in terms of F\o lner functions and we provide a proof for the sake of completeness. Recall that it is conjectured that the asymptotic behaviour of $\profp{G}$ does not depend on $p$.

\begin{lemma}\label{lem:MonotonuousProfile}
Let $G$ be a finitely generated group and let $p,q$ be real numbers such that $p>q\ge 1$. Then we have 
\begin{equation*}
    \folp{G}(x)\succcurlyeq \text{F\o l}_{q,G}(x).
\end{equation*}
\end{lemma}

For Theorem~\ref{th:lowerboundonFol}, we will apply this lemma to $q=1$.

\begin{proof}
Let $f\colon G\rightarrow \R$ be a finitely supported function, and consider the function $h\defeq |f|^{v}$ where $v\defeq \frac{p}{q}>1$. Using the inequality $|a^v-b^v|\le v\max{(a,b)}^{v-1}|a-b|$ that holds for every positive real numbers $a,b\ge 0$, we get
\begin{align*}
        &\|\nabla_{S_{G}} h\|_{q}^{q}=\sum_{g\in G,\; s\in S_{G}}{|h(g)-h(gs)|^q}\\
        &\le v^{q}\cdot \sum_{g\in G,\; s\in S_{G}}|f(g)|^{q(v-1)}\big||f(g)|-|f(gs)|\big|^{q}+v^{q}\cdot\sum_{g\in G,\; s\in S_{G}}|f(gs)|^{q(v-1)}\big| |f(g)|-|f(gs)|\big|^{q}\\
        &\le v^{q}\sum_{g\in G,\; s\in S_{G}}|f(g)|^{q(v-1)}|f(g)-f(gs)|^{q}+v^{q}\sum_{g\in G,\; s\in S_{G}}|f(gs)|^{q(v-1)}|f(g)-f(gs)|^{q}.
\end{align*}
Setting $P=\frac{p}{q}$ and $Q=\frac{p}{p-q}$, we have $\frac{1}{P}+\frac{1}{Q}=1$, and Hölder's inequality provides
\begin{align*}
    \sum_{g\in G,\; s\in S_{G}}{|f(g)|^{q(v-1)}|f(g)-f(gs)|^q} &\le \left(\sum_{g\in G,\; s\in S_{G}}{|f(g)|^{Qq(v-1)}}\right)^{\frac{1}{Q}}\left(\sum_{g\in G,\; s\in S_{G}}{|f(g)-f(gs)|^{Pq}}\right)^{\frac{1}{P}}\\
    &=\left(\sum_{g\in G,\; s\in S_{G}}{|f(g)|^{p}}\right)^{\frac{p-q}{p}}\left(\sum_{g\in G,\; s\in S_{G}}{|f(g)-f(gs)|^{p}}\right)^{\frac{q}{p}}\\
    &=|S_{G}|^{\frac{p-q}{p}}\cdot \|f\|_{p}^{p-q}\cdot \|\nabla_{S_G}f\|_{p}^{q}
\end{align*}
and similarly for $\sum_{g\in G,\; s\in S_{G}}{|f(gs)|^{q(v-1)}|f(g)-f(gs)|^q}$, so that we get
\begin{equation*}
        \|\nabla_{S_{G}} h\|_{q} \le 2^{\frac{1}{q}}\cdot |S_G|^{\frac{p-q}{pq}}\cdot v\cdot \|f\|_{p}^{\frac{p-q}{q}}\cdot \|\nabla_{S_G}f\|_{p}
\end{equation*}
which in turn implies
\begin{equation*}
        \frac{\|\nabla_{S_{G}} h\|_q}{\|h\|_q} \le 2^{\frac{1}{q}}\cdot |S_{G}|^{\frac{p-q}{pq}}\cdot v \cdot \frac{\|\nabla_{S_{G}} f\|_p}{\|f\|_p}
\end{equation*}
since $\|h\|_{q}^{q}=\|f\|_{p}^{p}$. This inequality holds for every finitely supported function $f\colon G\rightarrow \R$, namely we proved that for every such $f$, this inequality holds for some $h\colon G\rightarrow\R$ having the same support, so the statement follows directly from the definition of F\o lner functions.
\end{proof}

We now move on to the proof of Theorem~\ref{th:lowerboundonFol}. At first reading, the reader may look at Appendix~\ref{appendixA}, where the strategy is to find “good” lamplighters as subgroups of a lampshuffler $\shuf{H}$, namely a lamplighter group based on a finitely generated subgroup $K$ of $H$ having the same isoperimetric profile as $H$. This is achieved with some algebraic assumptions on the finitely generated group $H$, covering a large class of groups. We finally conclude using the result analogous to Lemma~\ref{lem:monotonieprofilErschler} for the F\o lner function. Moreover, this first strategy provides an interesting framework since the algebraic assumptions on the base group $H$ are stable in many cases when taking iterations of lampshufflers; see Remark~\ref{rem:Preservation} and Proposition~\ref{prop:A7}.

In this section, we focus on a less restrictive substructure than subgroups, namely subgraphs. The strategy is to find some subgraph $X_0$ of $H$, quasi-isometric to it, playing the role of a subgroup $K$ as described in the above first strategy, and a lamplighter graph on $X_0$ as a subgraph of $\halo H$, in such a manner that we can prove the monotonicity of the F\o lner function in this context, as in Lemma~\ref{lem:monotonieprofilErschler}. We conclude thanks to the lower bounds for the F\o lner functions of lamplighter graphs obtained in~\cite{Ers06}.

\paragraph{Lamplighter graphs.} Let $A$ and $B$ be two graphs, with a base vertex $b_0$ in $B$. Given a map $f\colon A\rightarrow B$, we define its support by $\supp{f}\defeq\lbrace a\in A : f(a)\neq b_{0}\rbrace$. The lamplighter graph of $B$ and $A$, denoted by $B\wr A$, is the graph
\begin{itemize}
    \item whose vertices are pairs $(f,a)$, where $a$ is a vertex of $A$ and $f\colon A\rightarrow B$ has finite support;
    \item whose edges connect $(f,a)$ and $(f',a')$ if either $a=a'$, $f(a)\sim_{B}f'(a)$ and $f(v)=f'(v)$ for every $v\in A\setminus\lbrace a\rbrace$, or if $f=f'$ and $a\sim_{A}a'$.
\end{itemize}

In the case where the graphs $A$ and $B$ are Cayley graphs of finitely generated groups $G$ and $H$ respectively, we recover a Cayley graph of the wreath product $H\wr G$.

Let us now prove Theorem~\ref{th:lowerboundonFol} within this framework, using the following lower bound proved by Erschler~\cite[Theorem~4.3]{Ers06}: there exists $C>0$ such that
\begin{equation*}
    \text{F\o l}_{B\wr A}(x)\succcurlyeq\left (\text{F\o l}_{B}(x)\right)^{C\text{F\o l}_{A}(x)}.
\end{equation*}

\begin{proof}[Proof of Theorem~\ref{th:lowerboundonFol}]
By Lemma~\ref{lem:MonotonuousProfile}, we have 
\begin{equation*}
    \folp{\halo H}(x)\succcurlyeq\text{F\o l}_{\halo H}(x)
\end{equation*}
so it is enough to prove the theorem for $p=1$.

\noindent Given a connected graph $Y$, we denote by $d_{Y}(\cdot,\cdot)$ its path metric. When considering a finitely generated group $H=\langle S\rangle$, we write $d_{H,S}(\cdot,\cdot)$ for the path metric on its Cayley graph $\text{Cay}(H,S)$ (identified with $H$ itself), to specify the choice of a finite generating subset $S$.

\noindent Now, let us fix a finite generating set $S$ of $H$, a constant $D\geq 0$ of large-scale commutativity for $\halo H$, and let us consider $S_{2D+5}\defeq\big\lbrace s_{1}s_2\ldots s_{2D+5} : s_{i}\in S\cup\lbrace 1_{H}\rbrace\big\rbrace$. Note that, for every $x,y\in H$, we have 
\begin{equation*}
    d_{H,S}(x,y)\le 2D+5\Longleftrightarrow d_{H,S_{2D+5}}(x,y)\le 1.
\end{equation*}
Let $X_{0}$ be a maximal $(D+2)-$separated subset of $H$, for the metric $d_{H,S}$, and let us endow $X_{0}$ with the graph structure induced by $d_{H,S_{2D+5}}$, namely $x,y\in X_{0}$ are adjacent if and only if $d_{H,S_{2D+5}}(x,y)=1$. It is straightforward to see that $(X_{0},d_{X_{0},S_{2D+5}})$ is a subgraph of $(H,d_{H,S_{2D+5}})$. In addition, we also have the following.
    
\begin{claim}\label{claim:qiToH}
The graphs $(X_{0},d_{X_{0},S_{2D+5}})$ and $(H,d_{H,S_{2D+5}})$ are quasi-isometric.
\end{claim}

\renewcommand{\qedsymbol}{$\blacksquare$}
\begin{proof}[Proof of Claim~\ref{claim:qiToH}]
Let us prove that the natural inclusion $X_{0}\hookrightarrow H$ is a quasi-isometry. First of all, by maximality, $X_{0}$ is $(D+2)-$dense in $(H,d_{S})$, and this directly implies $d_{H,S_{2D+5}}(h,X_{0})\le 1$ for every $h\in H$.
        
\noindent Let $x,y\in X_{0}$. It is straightforward to show that $d_{X_0,S_{2D+5}}(x,y)\ge d_{H,S_{2D+5}}(x,y)$. The other way around, let $n\defeq d_{H,S_{2D+5}}(x,y)$. By definition, there exist 
points 
\begin{equation*}
    x_0=x,x_1,\ldots,x_{n-1},x_n=y\in H
\end{equation*}
such that $d_{H,S_{2D+5}}(x_{i},x_{i+1})=1$ for every $i\in\lbrace 0,1,\ldots,n-1\rbrace$. Given such an index $i$, the definition of $S_{2D+5}$ implies that there exist points
\begin{equation*}
x_{i,0}=x_{i},x_{i,1},\ldots,x_{i,2D+4},x_{i,2D+5}=x_{i+1}\in H
\end{equation*}
such that $d_{H,S}(x_{i,j},x_{i,j+1})\le 1$ for every $j\in\lbrace 0,1,\ldots, 2D+4\rbrace$. Since we have $x_{i,2D+5}=x_{i+1,0}$ for every $i\in\lbrace 0,1,\ldots,n-1\rbrace$, we have found a path of length $\le (2D+5)n$ in $(H,d_{S})$ that connects $x$ to $y$. Approximating every vertex of this path by an element of $X_{0}$ within $d_{H,S}-$distance less than $D+2$ ($x$ and $y$ being approximated by themselves), we get a sequence 
\begin{equation*}
   w_{0}=x,w_{1},\ldots,w_{(2D+5)n-1},w_{(2D+5)n}=y 
\end{equation*}
of elements in $X_{0}$ satisfying $d_{H,S}(w_{i},w_{i+1})\le 2(D+2)+1=2D+5$, whence $d_{H,S_{2D+5}}(w_{i},w_{i+1})\le 1$. This way, we get a path from $x$ to $y$, of length $\le (2D+5)n$, in $(X_{0},d_{S_{2D+5}})$. Thus 
\begin{equation*}
 d_{X_{0},S_{2D+5}}(x,y)\le (2D+5)\cdot d_{H,S_{2D+5}}(x,y)   
\end{equation*}
and the proof of the claim is complete.
\end{proof}
\renewcommand{\qedsymbol}{$square$}
    
\noindent Let us fix some distinguished generator $s_{0}\in S\setminus\lbrace 1_{H}\rbrace$. By $(D+2)-$separation, for every $x\in X_{0}$, $xs_{0}$ does not lie in $X_{0}$, and large-scale commutativity thus implies that the groups $L(\lbrace x,xs_{0}\rbrace)$, for $x\in X_{0}$, commute. Let us now introduce the subgroup $\mathcal{T}$ of $L(H)$ defined by
\begin{align*}
    \mathcal{T}&\defeq\left\lbrace\prod_{x\in I}{\sigma_x}\ :  I \subset X_{0}\;\text{is finite}, \sigma_{x}\in L(\lbrace x,xs_{0}\rbrace)\right\rbrace \\
    &=\bigoplus_{x\in X_{0}}{L(\lbrace x,xs_{0}\rbrace)}\\
    &=\bigoplus_{x\in X_{0}}{\alpha(x)L(\lbrace 1_{H},s_{0}\rbrace)}.
\end{align*}
Since $\halo H$ has finitely generated blocks, we can fix a finite generating subset $S(s_{0})$ of $L(\lbrace 1_{H},s_{0}\rbrace)$.
Let us now consider the set $Y_{\star}\defeq\mathcal{T}\times X_{0}$ equipped with a graph structure where two vertices $(\rho,x)$ and $(\rho',x')$ of $Y_\star$ are adjacent if 
\begin{itemize}
        \item either $x=x'$ and $\rho^{-1}\rho'=\alpha(x)(\sigma)$ for some $\sigma\in S(s_0)$;
        \item or $\rho=\rho'$ and $d_{X_0,S_{2D+5}}(x,x')= 1$,
\end{itemize}
which can be reformulated as
\begin{itemize}
        \item either $(\rho,x)(\sigma,1_H)=(\rho',x')$, for some $\sigma\in S(s_0)$;
        \item or $(\rho,x)(\mathrm{id}_{H},h)=(\rho',x')$, where $h$ lies in $S_{2D+5}$.
\end{itemize}
Note that the graph $Y_{\star}$ is isomorphic to the lamplighter graph $L(\lbrace 1_{H},s_{0}\rbrace)\wr X_0$. Thus, we endow $\halo H$ with the finite generating set $S_{\halo H}$ given by
\begin{equation*}
    S_{\halo H}\defeq\lbrace (\sigma,1_{H}) : \sigma\in S(s_{0})\rbrace \cup\lbrace (1_{L(H)},h) : h\in S_{2D+5}\rbrace.
\end{equation*}
Now, let us consider the partition of $L(H)$ in $\mathcal{T}-$cosets:
\begin{equation*}
        L(H)=\bigsqcup_{c\in C}{\kappa_{c}\mathcal{T}}
\end{equation*}
with $\kappa_{c_{0}}=1_{L(H)}$ for the index $c_{0}\in C$ of the coset $\mathcal{T}$. For every $c\in C$, let us consider the subset
\begin{equation*}
    Y_{c}\defeq (\kappa_{c}\mathcal{T})\times H,
\end{equation*}
equipped with a graph structure where two vertices $(\kappa_{c}\rho,x)$ and $(\kappa_{c}\rho',x')$ are adjacent if 
\begin{itemize}
\item either $x=x'$, $x$ lies in $X_{0}$ and $\rho^{-1}\rho'=\alpha(x)(\sigma)$ for some $\sigma\in S(s_{0})$, namely $(\kappa_{c}\rho,x)(\sigma,1_{H})=(\kappa_{c}\rho',x')$;
\item or $\rho=\rho'$ and $d_{H,S_{2D+5}}(x,x')= 1$, namely $(\kappa_{c}\rho,x)(1_{L(H)},h)=(\kappa_{c}\rho',x')$ where $h$ lies in $S_{2D+5}$.
\end{itemize}
By definition of $S_{\halo H}$, $(Y_{c})_{c\in C}$ is a family of subgraphs of $(\halo H,d_{S_{\halo H}})$ partitioning the set of its vertices. Moreover the graph $Y_{c}$ is the left translation by $(\kappa_{c},1_{H})$ of the graph $Y_{c_{0}}$. These observations imply
\begin{equation*}
    \text{F\o l}_{\halo H}(n)\succcurlyeq \text{F\o l}_{Y_{c_{0}}}(n),
\end{equation*}
as an immediate adaptation of~\cite[Lemma~4]{Ers03}. The next claim is the final step required for the proof. 

\begin{claim}\label{claim:piecesQI}
The graph $Y_{c_{0}}$ is quasi-isometric to $Y_{\star}$.
\end{claim}

\renewcommand{\qedsymbol}{$\blacksquare$}
\begin{proof}[Proof of Claim~\ref{claim:piecesQI}]
The proof relies on the same technique as in the proof of Claim~\ref{claim:qiToH}. We prove that $Y_{\star}\hookrightarrow Y_{c_{0}}$ is a quasi-isometry. The $1-$density of its image is straightforward, as well as the inequality 
\begin{equation*}
    d_{Y_\star}\big((\rho,x),(\rho',x')\big)\ge d_{Y_{c_{0}}}\big((\rho,x),(\rho',x')\big)
\end{equation*}
for every $(\rho,x),(\rho',x')\in Y_{\star}$.

\noindent The other way around, notice that edges of a path of length $n\defeq d_{Y_{c_{0}}}\big((\rho,x),(\rho',x')\big)$ in $Y_{c_{0}}$ consists in either modifying the permutation on the first coordinate, or moving the arrow pointing at some element of $H$ in the second coordinate. Thus, with the same ideas as in the proof of Claim~\ref{claim:qiToH}, it suffices to approximate elements in the second coordinate by elements of $X_{0}$, so that we get a new path in $Y_{\star}$ of length $\le (2D+5)n$. This concludes the proof.
\end{proof}
\renewcommand{\qedsymbol}{$\Box$}
    
\noindent Combining the above claims, we finally get that 
\begin{equation*}
    \text{F\o l}_{\halo H}(n) \succcurlyeq \text{F\o l}_{Y_{c_{0}}}(n) \simeq \text{F\o l}_{Y_{\star}}(n) \simeq \text{F\o l}_{L(\{1_H,s_0\})\wr X_{0}}(n)
\end{equation*}
and the latter dominates $\left (\text{F\o l}_{L(\lbrace 1_{H},s_{0}\rbrace)}(n)\right )^{C'\text{F\o l}_{X_0}(n)}$, for some constant $C'>0$, using~\cite[Theorem~4.3]{Ers06}. From Claim~\ref{claim:qiToH}, $\text{F\o l}_{X_0}$ is asymptotically equivalent to $\text{F\o l}_{H}$, and thus 
\begin{equation*}
   \text{F\o l}_{\halo H}(n) \succcurlyeq \left(\text{F\o l}_{L(\lbrace 1_{H},s_{0}\rbrace)}(n)\right)^{C\text{F\o l}_{H}(n)} 
\end{equation*}
for some constant $C>0$.
\end{proof}

\section{Estimates of isoperimetric profiles for some examples of halo products}~\label{sec:computationsIsoProf}

The goal of this section is to establish our estimates of $\ell^{p}-$isoperimetric profiles of many halo products and their iterated versions, from our estimates on F\o lner functions. Here we use the fact that the $\ell^{p}-$F\o lner function and the $\ell^{p}-$isoperimetric profile are inverses of each other, using Remark~\ref{rem:InverseEachOther}.

Recall that we proved in the previous section:

\begin{theorem}[see Proposition~\ref{prop:upperboundonFol} and Theorem~\ref{th:lowerboundonFol}]\label{th:recallThm}
Let $p\ge 1$. Let $H$ be a finitely generated amenable group and let $S_{H}$ be a finite generating set. Let $\halo H$ be a naturally generated halo product over $H$.
\begin{enumerate}[label=(\roman*)]
    \item If $\halo H$ is large-scale commutative and has finitely generated blocks, then for any $s_{0}\in S_{H}$, there exists a constant $C>0$ such that
\begin{equation*}
    \left(\text{F\o l}_{L(\lbrace 1_{H},s_{0} \rbrace)}(x)\right)^{C\text{F\o l}_{H}(x)} \preccurlyeq \folp{\halo H}(x).
\end{equation*}
    \item If $\halo H$ has consistent blocks, then there exists a constant $C>0$ such that
\begin{equation*}
    \folp{\halo H}(x) \preccurlyeq \folp{H}(x)\cdot \Lambda_{\halo H}(C\cdot \folp{H}(x)).
\end{equation*}
\end{enumerate}
\end{theorem}

As an easy consequence, if $\halo H$ is large-scale commutative, naturally generated and has consistent blocks, then its $\ell^{p}-$F\o lner function satisfies
\begin{equation*}
    K^{\text{F\o l}_{H}(x)}\preccurlyeq\folp{\halo H}(x) \preccurlyeq \folp{H}(x)\cdot \Lambda_{\halo H}(C\cdot \folp{H}(x))
\end{equation*}
for some positive constants $C,K>0$.

Let us first discuss a terminology that will be useful for our main results.

\subsection{Assumption~$(\star)$}\label{sec:assumptionstar}

\begin{definition}
We say that a non-decreasing map $h\colon\R_{+}\rightarrow \R_{+}$ satisfies \textit{Assumption}~$(\star)$ if one has 
\begin{equation*}
    \forall C>0,\; h(Cx)=O(h(x)). 
\end{equation*}
\end{definition}

Assumption~$(\star)$ already appeared in the literature~\cite{Ers03,Cor25} in the case where $h=\profp{H}$ is the $\ell^{p}-$isoperimetric profile of a finitely generated group $H$, and it seems that $\profp{H}$ satisfies this assumption for many choices of groups $H$. In fact, to our knowledge, there is currently no known example of a finitely generated group whose $\ell^{p}-$isoperimetric profiles do not satisfy Assumption~$(\star)$.

First, let us record in a statement an easy implication of Assumption~$(\star)$.

\begin{lemma}\label{lem:profileandequivalents}
Let $h\colon\R_{+}\rightarrow \R_{+}$ be a non-decreasing map satisfying Assumption~$(\star)$. Let $f,g\colon \R_{+}\rightarrow \R_{+}$ be unbounded non-decreasing maps such that $f(x)\preccurlyeq g(x)$. Then one has 
\begin{equation*}
    h(f(x)) \preccurlyeq h(g(x)).
\end{equation*}
\end{lemma}

\begin{proof}
By assumption, $f(x) \preccurlyeq g(x)$, so there is a constant $C>0$ such that 
\begin{equation*}
    f(x) \le Cg(Cx)
\end{equation*}
for all $x$ large enough. As $h$ is non-decreasing, we thus get 
\begin{equation*}
    h(f(x)) \le h\big(Cg(Cx)\big).
\end{equation*}
Now, we use 
$h(Cy)=O(h(y))$ to deduce that there is a constant $K>0$ such that 
\begin{equation*}
    h(f(x)) \le h\big(Cg(Cx)\big) \le Kh(g(Cx))
\end{equation*}
for all $x$ large enough. Thus $h(f(x)) \preccurlyeq h(g(x))$ and we are done.
\end{proof}

Another useful claim is the following.

\begin{lemma}\label{lem:UsefulLemmaProfile}
Let $f,g,\varphi\colon \left[1,+\infty\right[\longrightarrow \left[1,+\infty\right[$ be three non-decreasing maps, with $\varphi$ injective and satisfying $\varphi(x)\underset{x\to +\infty}{\to}+\infty$. Assume that there exists a positive constant $D>0$ such that
\begin{equation*}
        g(x)\preccurlyeq f(D\varphi(x)).
\end{equation*}
Then we have
\begin{equation*}
        g(K\varphi^{-1}(x))\preccurlyeq f(x)
\end{equation*}
for some positive constant $K>0$. If furthermore $g$ satisfies Assumption~$(\star)$, then we have
\begin{equation*}
        g(\varphi^{-1}(x))\preccurlyeq f(x).
\end{equation*}
\end{lemma}

\begin{proof}
By assumption, there exists a positive constant $C>0$ such that
\begin{equation*}
    g(x)\le Cf\big(D\varphi(Cx)\big)
\end{equation*}
for all $x$ large enough. Let $x$ be a real number greater than $D\varphi(C)$, and let $n\geq 1$ be an integer such that $D\varphi(Cn)\le x\le D\varphi(C(n+1))$. Then we have
\begin{equation*}
    Cf(x)\ge Cf\big(D\varphi(Cn)\big)\ge g(n)\ge g\left(\frac{n+1}{2}\right)\ge g\left (\frac{\varphi^{-1}(\frac{x}{D})}{2C}\right),
\end{equation*}
which can be reformulated as
\begin{equation*}
    g\left(K\varphi^{-1}(y)\right)=O\left(f(Dy)\right),
\end{equation*}
taking $y=\frac{x}{D}$ and $K=\frac{1}{2C}$. This shows the first part of the statement. Additionally, if $g$ satisfies Assumption~$(\star)$, then we have
\begin{equation*}
g\left(\varphi^{-1}\left(y\right)\right)=g\left(\frac{1}{K}\cdot K\varphi^{-1}(y)\right)=O\left(g\left(K\varphi^{-1}(y)\right)\right)=O\left(f(Dy)\right),
\end{equation*}
which concludes the proof.
\end{proof}

With Assumption~$(\star)$, we can deduce two applications for the computation of isoperimetric profiles. The first one reformulates the inequality
\begin{equation*}
    \left(\text{F\o l}_{L(\lbrace1_{H},s_{0}\rbrace)}(x)\right)^{C\text{F\o l}_{H}(x)} \preccurlyeq \folp{\halo H}(x)
\end{equation*}
from Theorem~\ref{th:lowerboundonFol}, in the easier case where the block $L(\lbrace 1_{H},s_{0} \rbrace)$ is finite.

\begin{corollary}\label{cor:upperboundonProf}
Let $H$ be a finitely generated amenable group. Let $\halo H$ be a naturally generated and large-scale commutative halo product having finite blocks. Given $p\ge 1$, if $\prof{H}$ satisfies Assumption~$(\star)$, then the $\ell^{p}-$isoperimetric profile of $\halo H$ satisfies
\begin{equation*}
    \profp{\halo H}(x)\preccurlyeq\prof{H}(\ln(x)).
\end{equation*}
\end{corollary}

\begin{proof}
By Theorem~\ref{th:lowerboundonFol} and the fact that $L(\lbrace 1_{H},s_{0}\rbrace)$ is finite for every generator $s_{0}\in S_{H}$, we have
    \begin{equation*}
    K^{\text{F\o l}_{H}(y)} \preccurlyeq \folp{\halo H}(y)
\end{equation*}
for some positive constant $K>0$. The logarithm satisfies Assumption~$(\star)$, as well as $\prof{H}$, so applying twice Lemma~\ref{lem:profileandequivalents} yields
\begin{equation*}
    y \preccurlyeq \prof{H}\big(\ln\big (\folp{\halo H}(y)\big)\big),
\end{equation*}
meaning that there exists some positive constant $C>0$ such that
\begin{equation*}
    y \le C\cdot \prof{H}\big(\ln\big(\folp{\halo H}(Cy)\big)\big)
\end{equation*}
for all $y$ large enough. If now $x$ is large enough, it suffices to take $y=\frac{\profp{\halo H}(x)}{C}$ in the above inequality to deduce the corollary.
\end{proof}

As a second application, we record in a statement the formulation of the inequality
\begin{equation*}
    \folp{\halo H}(x)\preccurlyeq\folp{H}(x)\cdot \Lambda_{\halo H}(C\cdot\folp{H}(x))
\end{equation*}
from Proposition~\ref{prop:upperboundonFol} in terms of isoperimetric profiles.

\begin{corollary}\label{cor:lowerboundonProf}
Let $H$ be a finitely generated amenable group, and let $\halo H$ be a halo product over $H$. Suppose that $\halo H$ is naturally generated and has consistent blocks. Given $p\ge 1$, if $\profp{H}$ satisfies Assumption~$(\star)$, then one has 
\begin{equation*}
    \profp{\halo H}(x) \succcurlyeq \profp{H}(\varphi^{-1}(x))
\end{equation*}
where $\varphi(x)=x\cdot \Lambda_{\halo H}(x)$ and where $\Lambda_{\halo H}$ is the lamp growth sequence of $\halo H$. 
\end{corollary}

\begin{proof}
Let $p\ge 1$. From Proposition~\ref{prop:upperboundonFol}, we know that there is a constant $D>0$ such that
\begin{equation*}
    \folp{\halo H}(x) \le D\cdot \varphi\big(D\cdot\folp{H}(Dx)\big)
\end{equation*}
for all $x$ large enough. Putting $x=\frac{\profp{H}(y)}{D}$ in this inequality for $y$ large enough, and applying $\profp{\halo H}$ which is increasing, one gets
\begin{equation*}
    \frac{\profp{H}(y)}{D} \le \profp{\halo H}\big(D\cdot \varphi(Dy)\big)
\end{equation*}
for all $y$ large enough, i.e. $\profp{H}(y) \preccurlyeq \profp{\halo H}\big(D\varphi(y)\big)$. Since $\profp{H}$ satisfies Assumption~$(\star)$, we may apply Lemma~\ref{lem:UsefulLemmaProfile}, and we get
\begin{equation*}
    \profp{H}(\varphi^{-1}(x)) \preccurlyeq \profp{\halo H}(x)
\end{equation*}
as claimed. 
\end{proof}

\subsection{Proof of Theorem~\ref{thm:boundsForProfile intro} for lampjugglers}

We now apply our estimates on isoperimetric profiles to concrete examples, using~\cite[Facts~7.12-7.16]{GT24a} that compute the lamp growth sequences of most examples of halo products we are interested in. In this section, we address the case of lampshufflers, lampjugglers and their iterated versions. 

Let us recall that $\ell^{p}-$isoperimetric profiles of lampshufflers over polynomial growth groups are known:
\begin{equation*}
    \profp{\shuf{H}}(x) \simeq \left(\frac{\ln(x)}{\ln(\ln(x))}\right)^{\frac{1}{d}}
\end{equation*}
for any $p\ge 1$, when $H$ has growth degree $d\ge 1$. This can be directly deduced from~\cite{SCZ21} and~\cite{EZ21}.

From our work, we can deduce the following bounds on profiles of lampjugglers. 

\begin{corollary}\label{cor:encadrementduprofildeshuf}
Let $p\ge 1$. Let $H$ be a finitely generated amenable group, whose $\ell^{p}-$ and $\ell^{1}-$isoperimetric profiles satisfy Assumption~$(\star)$. Let $s\ge 1$ be an integer. Then the $\ell^{p}-$isoperimetric profile of $\juggler{s}{H}$ satisfies 
\begin{equation*}
    \profp{H}\left(\frac{\ln(x)}{\ln(\ln(x))}\right) \preccurlyeq \profp{\juggler{s}{H}}(x) \preccurlyeq \prof{H}(\ln(x)).
\end{equation*}
\end{corollary}

\begin{proof}
For the upper bound, it suffices to apply Corollary~\ref{cor:upperboundonProf}, since $\juggler{s}{H}$ is naturally generated and has finite blocks. 

\noindent Let us focus on the lower bound. From Corollary~\ref{cor:lowerboundonProf}, we know that 
\begin{equation}\label{eq:lowerboundwithinverse}
    \profp{H}(\varphi^{-1}(x)) \preccurlyeq \profp{\juggler{s}{H}}(x)
\end{equation}
where $\varphi(x)=x\cdot \Lambda_{\juggler{s}{H}}(x)=x\cdot (sx)!$. It remains to find the asymptotics of $\varphi^{-1}(x)$. We have by definition $\varphi^{-1}(x)(s\varphi^{-1}(x))!=x$, and from Stirling's formula we know that
\begin{equation*}
    \varphi^{-1}(x)(s\varphi^{-1}(x))! \sim \varphi^{-1}(x)\left(\frac{s\varphi^{-1}(x)}{e}\right)^{s\varphi^{-1}(x)}\sqrt{2\pi\cdot s\varphi^{-1}(x)}.
\end{equation*}
Taking the logarithm yields $\ln(x)=\ln\left(\varphi^{-1}(x)(s\varphi^{-1}(x))!\right) \sim s\varphi^{-1}(x)\ln(s\varphi^{-1}(x))$ and taking the logarithm once more, it follows that 
\begin{equation*}
    \ln(\ln(x)) \sim \ln(s\varphi^{-1}(x)).
\end{equation*}
Combining these two equivalences, this gives 
\begin{equation*}
    s\varphi^{-1}(x)\sim\frac{\ln(x)}{\ln(\ln(x))}. 
\end{equation*}
Using Lemma~\ref{lem:profileandequivalents} and inequality~\eqref{eq:lowerboundwithinverse}, we deduce 
\begin{equation*}
    \profp{H}\left(\frac{\ln(x)}{\ln(\ln(x))}\right) \preccurlyeq \profp{\juggler{s}{H}}(x)
\end{equation*}
as claimed. The proof is complete.
\end{proof}

\begin{remark}\label{rm:EZnotoptimal}
In the particular case where $s=1$ and where $H$ is amenable with $\prof{H}(x)\simeq\ln(x)$, the upper bound on $\prof{\shuf{H}}$ provided by this result is $\simeq \ln(\ln(x))$, which is the same that one can get by inverting the lower bound 
\begin{equation*}
    \text{F\o l}_{\shuf{H}}(x) \succcurlyeq V_{H}(x)^{V_{H}(x)}\simeq (e^{x})^{e^{x}}
\end{equation*}
obtained in~\cite[Corollary~1.4]{EZ21} in the case of an exponential growth group $H$. However, as we will see below, this result is not optimal anymore when taking iteration of lampshufflers, since the growth function of such iterations stays exponential, while the isoperimetric profile gets slower at each iteration (see Proposition~\ref{prop:profilesofshufn}). 
\end{remark}

In particular, one deduces that, if $H$ has polynomial growth of degree $d\ge 1$, the estimate on $\profp{\shuf{H}}$ coming from~\cite[Corollary~1.4]{EZ21} is valid for any lampjuggler $\juggler{s}{H}$, $s\ge 2$. Indeed, such a lampjuggler contains a lampshuffler as a subgroup, so that 
\begin{equation*}
    \profp{\juggler{s}{H}}(x) \preccurlyeq \profp{\shuf{H}}(x) \simeq \left(\frac{\ln(x)}{\ln(\ln(x))}\right)^{\frac{1}{d}},
\end{equation*}
by Lemma~\ref{lem:monotonieprofilErschler}, and the lower bound follows from Corollary~\ref{cor:encadrementduprofildeshuf}. In fact, we have more generally the next consequence.

\begin{corollary}\label{cor:profilesoflampjugglers}
Let $p\geq 1$. Let $H$ be a finitely generated amenable group such that $\profp{H}(x)\simeq \prof{H}(x)$. Assume that
\begin{itemize}
    \item either $H$ has polynomial growth;
    \item or its $\ell^{p}-$isoperimetric profile $\profp{H}$ satisfies Assumption~$(\star)$ and $\profp{H}\left(\frac{\ln(x)}{\ln(\ln(x))}\right)\simeq \profp{H}(\ln(x))$.
\end{itemize}
Then one has 
\begin{equation*}\label{eq:profilesoflampjugglers}
    \profp{\juggler{s}{H}}(x) \simeq \profp{\shuf{H}}(x)
\end{equation*}
for all $s\ge 1$.\qed
\end{corollary}

Thus, lampjuggler groups often have the same $\ell^{p}-$isoperimetric profile as lampshufflers, even if the two are not quasi-isometric. For instance, if $H=\Z^{d}\wr\text{BS}(1,n)$, $d\ge 1$, $n\ge 2$, then $\shuf{H}$ and $\juggler{s}{H}$ are not quasi-isometric (by Theorems~\ref{thm:QIrigidityofjugglers} and~\ref{thm:classificationforvirtuallyabeliangroups}) but both have $\ell^{p}-$isoperimetric profile $\simeq \profp{H}(\ln(x)) \simeq \ln(\ln(\ln(x)))$.

From Corollary~\ref{cor:encadrementduprofildeshuf}, we can deduce the $\ell^{p}-$isoperimetric profile of lampshufflers over groups $H$ having isoperimetric profiles that satisfy
\begin{equation*}
   \profp{H}\left(\frac{\ln(x)}{\ln(\ln(x))}\right) \simeq \prof{H}(\ln(x))
\end{equation*}
for instance: 
\begin{itemize}
    \item Solvable Baumslag-Solitar groups $\text{BS}(1,n)$, $n\ge 2$, have $\ell^{p}-$isoperimetric profile $\simeq \ln(x)$. Thus $\shuf{\text{BS}(1,n)}$ has $\ell^{p}-$isoperimetric profile $\simeq \ln(\ln(x))$. The same applies for lamplighters $F\wr\Z$, where $F$ is a non-trivial finite group;
    \item more generally, lamplighters $F\wr \Z^d$, where $d\ge 1$ and $F$ is non-trivial and finite, have $\ell^{p}-$isoperimetric profile $\simeq \ln(x)^{\frac{1}{d}}$. In this case, we get that 
    \begin{equation*}
        \profp{\shuf{F\wr\Z^d}}(x) \simeq \ln(\ln(x))^{\frac{1}{d}}. 
    \end{equation*}
    \item for $d\ge 1$, the group $H=\Z\wr\Z^d$ has $\ell^{p}-$isoperimetric profile $\simeq \left(\frac{\ln(x)}{\ln(\ln(x))}\right)^{\frac{1}{d}}$, so that 
    \begin{equation*}
        \profp{\shuf{\Z\wr\Z^d}}(x) \simeq \left(\frac{\ln(\ln(x))}{\ln(\ln(\ln(x)))}\right)^{\frac{1}{d}}.
    \end{equation*}
\end{itemize}

\begin{example}\label{ex:IsoProfBrieusselZheng}
For any non-decreasing function $f\colon \R_{+}\rightarrow\R_{+}$ such that $x\longmapsto \frac{x}{f(x)}$ is non-decreasing, Brieussel and Zheng constructed in~\cite[Theorem~1.1]{BZ21} a finitely generated group $H$ with exponential volume growth having $\ell^{p}-$isoperimetric profile 
\begin{equation*}
    \profp{H}(x)\simeq \frac{\ln(x)}{f(\ln(x))}. 
\end{equation*}
It is proved in~\cite{Cor25} that $\profp{H}$ satisfies Assumption~$(\star)$ (see the discussion right after Corollary 4.1 in~\cite{Cor25}). Lastly, $H$ also satisfies 
\begin{equation*}
    \profp{H}\left(\frac{\ln(x)}{\ln(\ln(x))}\right) \simeq \prof{H}(\ln(x)).
\end{equation*}
This follows from the fact that $f$ preserves equivalents: if $g,h\colon \R_{+}\rightarrow\R_{+}$ are equivalent, then $(1-\varepsilon)h(x)\le g(x)\le (1+\varepsilon)h(x)$ for some $\varepsilon>0$ and for large enough $x$, so that
\begin{equation*}
    f\left((1-\varepsilon)h(x)\right) \le f(g(x)) \le f\left((1+\varepsilon)h(x)\right)
\end{equation*}
since $f$ is non-decreasing. Since $(1-\varepsilon)h(x) \le h(x) \le (1+\varepsilon)h(x)$, one has also
\begin{equation*}
    \frac{(1-\varepsilon)h(x)}{f((1-\varepsilon)h(x))} \le \frac{h(x)}{f(h(x))} \le \frac{(1+\varepsilon)h(x)}{f((1+\varepsilon)h(x))}.
\end{equation*}
since $x\longmapsto\frac{x}{f(x)}$ is non-decreasing, and thus 
\begin{equation*}
    (1-\varepsilon)f(h(x)) \le f\big((1-\varepsilon)h(x)\big) \le f(g(x)) \le f\big((1+\varepsilon)h(x)\big) \le (1+\varepsilon)f(h(x))
\end{equation*}
for all large enough $x$, whence $f(g(x)) \sim f(h(x))$. Thus, we may apply Theorem~\ref{cor:encadrementduprofildeshuf}, and we obtain that 
\begin{equation*}
    \profp{\shuf{H}}(x) \simeq \frac{\ln(\ln(x))}{f(\ln(\ln(x)))}. 
\end{equation*}
\end{example}

\paragraph{Iterated lampshufflers.} For a group $H$ and an integer $n\ge 0$, let $\shufn{n}{H}$ denote the $n-$th iterated lampshuffler over $H$, defined as  
\begin{equation*}
    \shufn{n}{H} \defeq \shuf{\shufn{n-1}{H}}.
\end{equation*}
if $n\ge 1$, and $\shufn{0}{H} \defeq H$. 

For such groups, we show the following estimates. 

\begin{proposition}\label{prop:profilesofshufn}
Let $p\ge 1$. Let $H$ be a finitely generated amenable group whose $\ell^{p}-$isoperimetric profile $\profp{H}$ satisfies Assumption~$(\star)$. Suppose that 
\begin{equation*}
    \profp{H}\left(\frac{\ln(x)}{\ln(\ln(x))}\right) \simeq \profp{H}(\ln(x)) \; \text{and} \; \profp{H}(x) \simeq \prof{H}(x).
\end{equation*}
Then, we have 
\begin{equation*}
   \profp{\shufn{n}{H}}(x) \simeq \profp{H}(\ln^{\circ n}(x))
\end{equation*}
for all $n\ge 0$. 
\end{proposition}

Here, recall that $\ln^{\circ k}(x) \defeq \ln(\ln(\ln(\dots\ln(x))))$ denotes the $k-$th iteration of the logarithm with itself, with the convention that $\ln^{\circ 0}$ is the identity.  

\begin{proof}
Let $p\ge 1$. We prove the statement by induction over $n$. For $n=0$, it clearly holds, and the case $n=1$ is settled by Corollary~\ref{cor:encadrementduprofildeshuf}. Now, assume that it holds for some $n\ge 0$, so that $\profp{\shufn{n}{H}}(x)\simeq \profp{H}(\ln^{\circ n}(x))$ for any $H$ satisfying the assumptions of the statement. Then, note that $\shuf{H}$ is finitely generated, amenable, and since its $\ell^{p}-$isoperimetric profile is $\simeq \profp{H}(\ln(x))$ and that $\profp{H}$ satisfies Assumption~$(\star)$, $\profp{\shuf{H}}$ satisfies Assumption~$(\star)$ and $\profp{\shuf{H}}\simeq \prof{\shuf{H}}$ as well. Additionally, note that 
\begin{equation*}
    \profp{\shuf{H}}\left(\frac{\ln(x)}{\ln(\ln(x))}\right) \simeq \profp{H}\big(\ln^{\circ 2}(x)-\ln^{\circ 3}(x)\big)\simeq \profp{H}(\ln^{\circ 2}(x)) \simeq \profp{\shuf{H}}(\ln(x))
\end{equation*}
where the second asymptotic equivalence follows from the combination of $\ln^{\circ 2}(x)-\ln^{\circ 3}(x) \simeq \ln^{\circ 2}(x)$ and Lemma~\ref{lem:profileandequivalents}. Thus, applying the inductive assumption and Corollary~\ref{cor:encadrementduprofildeshuf}, we obtain
\begin{align*}
    \profp{\shufn{(n+1)}{H}}(x)&=\profp{\shufn{n}{\shuf{H}}}(x) \\
    &\simeq \profp{\shuf{H}}(\ln^{\circ n}(x)) \\
    &\simeq \profp{H}\big(\ln(\ln^{\circ n}(x))\big) \\
    &=\profp{H}(\ln^{\circ (n+1)}(x))
\end{align*}
and the induction is complete.
\end{proof}

However, note that the assumption that $\profp{H}\left(\frac{\ln(x)}{\ln(\ln(x))}\right) \simeq \profp{H}(\ln(x))$ in Proposition~\ref{prop:profilesofshufn} does not hold for polynomial growth groups. Thus, for this class, we compute isoperimetric profiles of iterated lampshufflers separately. 

\begin{proposition}\label{prop:profileofshufnofpolynomialgrowthgroups}
Let $H$ be a finitely generated group of polynomial growth of degree $d\ge 1$. Then one has 
\begin{equation*}
    \profp{\shufn{n}{H}}(x) \simeq \left(\frac{\ln^{\circ n}(x)}{\ln^{\circ (n+1)}(x)}\right)^{\frac{1}{d}}
\end{equation*}
for any $n\ge 1$ and any real number $p\ge 1$. 
\end{proposition}

\begin{proof}
The case $n=1$ is settled by~\cite[Corollary~1.4]{EZ21}. For $n\ge 2$, it suffices to note that $\shuf{H}$ satisfies the assumptions of Proposition~\ref{prop:profilesofshufn}, and the latter provides
\begin{equation*}
    \profp{\shufn{n}{H}}(x)=\profp{\shufn{n-1}{\shuf{H}}}(x) \simeq \profp{\shuf{H}}\big(\ln^{\circ (n-1)}(x)\big) \simeq \left(\frac{\ln^{\circ n}(x)}{\ln^{\circ (n+1)}(x)}\right)^{\frac{1}{d}}
\end{equation*}
as claimed.
\end{proof}

The assumption $\profp{H}\left(\frac{\ln(x)}{\ln(\ln(x))}\right) \simeq \profp{H}(\ln(x))$ from Proposition~\ref{prop:profilesofshufn} is satisfied for many known behaviours of profiles, and thus motivates the next question. 

\begin{question}\label{q:profile}
Is it true that all finitely generated amenable groups which do not have polynomial growth satisfy $\profp{H}\left(\frac{\ln(x)}{\ln(\ln(x))}\right) \simeq \profp{H}(\ln(x))$?

\end{question}

\begin{remark}\label{rem:IteratedJuggler}
The same strategy shows that if $H$ is a finitely generated amenable group whose isoperimetric profiles satisfy Assumption~$(\star)$ and 
\begin{equation*}
    \profp{H}\left(\frac{\ln(x)}{\ln(\ln(x))}\right) \simeq \profp{H}(\ln(x))\; \text{and} \;\profp{H}(x) \simeq \prof{H}(x),
\end{equation*}
then one has 
\begin{equation*}
    \profp{    \juggler{s_{1}}{\juggler{s_{2}}{\dots\juggler{s_{n}}{H}}}       }(x) \simeq \profp{H}(\ln^{\circ n}(x))
\end{equation*}
for any integers $n\ge 1$ and $s_{1},\dots,s_{n}\ge 1$, and any real number $p\ge 1$. 
\end{remark}

\subsection{Proof of Corollary~\ref{cor:profilesofJugglersandDesignersOverPolyGrowthGroups} for lampdesigners}\label{sec:lampdesigner}

Lampdesigners are close to lampjuggler groups, and in fact if $F$ is finite, $\designer{H}$ is a subgroup of $\juggler{|F|}{H}$, via the map 
\begin{align*}
    \begin{array}{cll}
    \designer{H} &\longrightarrow &\juggler{|F|}{H} \\
    ((f,\sigma),h) &\longmapsto &(\sigma', h)
    \end{array}
\end{align*}
where, given a pair $(f,\sigma)\in F\wr_{H}\fsym{H}$, $\sigma'$ is the permutation of $H\times F$ given by $\sigma'(h,i)=(\sigma(h), f(h)i)$. Hence, from Lemma~\ref{lem:monotonieprofilErschler} and Proposition~\ref{prop:upperboundonFol}, we get directly a lower bound on $\ell^{p}-$isoperimetric profiles of lampdesigners, namely
\begin{equation*}
    \profp{\juggler{|F|}{H}}(x)\preccurlyeq \profp{\designer{H}}(x).
\end{equation*}
Additionally, note that $\designer{H}$ contains $\shuf{H}$ as a subgroup (and also as a quotient), hence
\begin{equation*}
    \profp{\designer{H}}(x) \preccurlyeq \profp{\shuf{H}}(x).
\end{equation*}
Moreover, recall that, when $H$ satisfies $\profp{H}(x)\simeq\prof{H}(x)$ and one of the two hypotheses:
\begin{itemize}
    \item $H$ has polynomial growth;
    \item $\profp{H}$ satisfies Assumption~$(\star)$, \begin{equation*}
    \profp{H}\left(\frac{\ln(x)}{\ln(\ln(x))}\right) \simeq \profp{H}(\ln(x)),
\end{equation*}
\end{itemize}
then Corollary~\ref{cor:profilesoflampjugglers} ensures that $\juggler{|F|}{H}$ and $\shuf{H}$ have same $\ell^{p}-$isoperimetric profile, and thus:
\begin{corollary}\label{cor:profileoflampdesigners}
Let $F$ be a non-trivial finite group. Let $p\ge 1$. Let $H$ be a finitely generated amenable group such that $\profp{H}(x)\simeq\prof{H}(x)$. Assume that one of the following holds:
\begin{itemize}
    \item $H$ has polynomial growth;
    \item $\profp{H}$ has Assumption~$(\star)$ and satisfies 
\begin{equation*}
    \profp{H}\left(\frac{\ln(x)}{\ln(\ln(x))}\right) \simeq \profp{H}(\ln(x)).
\end{equation*}
\end{itemize}
Then one has 
\begin{equation*}
    \profp{\designer{H}}(x)\simeq\profp{\shuf{H}}(x).
\end{equation*}
\end{corollary}

Furthermore, we can also deduce $\ell^{p}-$isoperimetric profiles of iterated lampdesigners.

\subsection{Proof of Corollary~\ref{cor:encadrementduprofildecloneINTRO} for lampcloners}

Let us now turn to lampcloners and lampupcloners over finite fields. 

\begin{corollary}\label{cor:encadrementduprofildeclone}
Let $p\ge 1$. Let $H$ be a finitely generated amenable group whose $\ell^{p}-$isoperimetric profile $\profp{H}$ satisfies Assumption~$(\star)$. Let $\field$ be a finite field. Then one has 
\begin{equation*}
    \profp{H}\left(\sqrt{\ln(x)}\right) \preccurlyeq \profp{\cloner{H}}(x) \preccurlyeq \prof{H}(\ln(x)).
\end{equation*}
\end{corollary}

\begin{proof}
The upper bound directly follows from Corollary~\ref{cor:upperboundonProf}, so we focus on the lower bound. 

\noindent We know from Corollary~\ref{cor:lowerboundonProf} that we must determine the asymptotic behaviour of $\varphi^{-1}$, where $\varphi(x)=x\cdot\Lambda_{\cloner{H}}(x)$. By definition, we have that $\varphi^{-1}(x)\cdot\Lambda_{\cloner{H}}(\varphi^{-1}(x))=x$, and thus 
\begin{equation*}
    \ln(\varphi^{-1}(x))+\ln\big(\Lambda_{\cloner{H}}(\varphi^{-1}(x))\big)=\ln(x).
\end{equation*}
From~\cite[Fact~7.16]{GT24a}, $\ln(\Lambda_{\cloner{H}}(y)) \sim C\cdot y^2$ for some $C>0$, so the above equation tells us that 
\begin{equation*}
    \varphi^{-1}(x)^2 \simeq \ln(x)
\end{equation*}
whence $\varphi^{-1}(x) \simeq \sqrt{\ln(x)}$. From Lemma~\ref{lem:UsefulLemmaProfile}, it follows that
\begin{equation*}
    \profp{H}(\varphi^{-1}(x))\simeq \profp{H}\left(\sqrt{\ln(x)}\right)
\end{equation*}
and we are done. 
\end{proof}

From here, we then directly deduce the following consequence. 

\begin{corollary}\label{cor:profCloner}
Let $p\ge 1$. Let $H$ be a finitely generated amenable group whose $\ell^{p}-$isoperimetric profile $\profp{H}$ satisfies Assumption~$(\star)$. Let $\field$ be a finite field. If $\profp{H}\left(\sqrt{\ln(x)}\right) \simeq \prof{H}(\ln(x))$ and $\profp{H}(x)\simeq \prof{H}(x)$, then 
\begin{equation*}
    \profp{\cloner{H}}(x) \simeq \profp{H}(\ln(x)).
\end{equation*}
\end{corollary}

This corollary applies to many groups that have slow profiles, for instance:
\begin{itemize}
    \item Baumslag-Solitar groups $\text{BS}(1,n)$, $n\ge 2$, whose $\ell^{p}-$isoperimetric profile is $\simeq\ln(x)$. Thus $\cloner{\text{BS}(1,n)}$ has $\ell^{p}-$isoperimetric profile $\simeq \ln(\ln(x))$ for any $n\ge 2$ and finite field $\field$;
    \item The lamplighter $F\wr\Sigma$, where $\Sigma$ has polynomial growth of degree $d\ge 1$, has $\ell^{p}-$isoperimetric profile $\simeq \ln(x)^{\frac{1}{d}}$, whence 
    \begin{equation*}
        \profp{\cloner{F\wr\Sigma}}(x) \simeq \ln(\ln(x))^{\frac{1}{d}}.
    \end{equation*}
\end{itemize}

In the polynomial growth case, we get the following bounds.

\begin{corollary}\label{cor:clonerPolynomial}
Let $H$ be a finitely generated group of polynomial growth of degree $d\ge 1$. Let $\field$ be a finite field. Then we have
\begin{equation*}
    \ln(x)^{\frac{1}{2d}}\preccurlyeq\profp{\cloner{H}}(x) \preccurlyeq \ln(x)^{\frac{1}{d}}
\end{equation*}
for any real number $p\ge 1$. 
\end{corollary}

Inspired by the case of lampshufflers, we would expect that, when $H$ has polynomial growth of degree $d\ge 1$, the $\ell^{p}-$isoperimetric profile of $\cloner{H}$ is the lower bound that we found in the above statement, namely $\ln(x)^{\frac{1}{2d}}$. Recall that for lampshufflers, we applied the upper bound from~\cite[Corollary~1.4]{EZ21} which is still optimal in the polynomial growth case, but its proof seems difficult to generalize for lampcloners.

\begin{remark}
We can slightly improve the upper bound in Corollary~\ref{cor:clonerPolynomial}, since for any group $H$, $\shuf{H}$ is a subgroup of $\cloner{H}$, considering the linear automorphisms permuting the vectors of the canonical basis provided by $H$. Hence, if $H$ has polynomial growth of degree $d\ge 1$, we have
\begin{equation*}
    \profp{\cloner{H}}(x)\preccurlyeq\left(\frac{\ln(x)}{\ln(\ln(x))}\right)^{\frac{1}{d}}.
\end{equation*}
\end{remark}

\paragraph{Iterated lampcloners.} For a group $H$ and an integer $n\ge 0$, let $\clonern{n}{H}$ denote the $n$-th iterated lampcloner over $H$, defined as  
\begin{equation*}
    \clonern{n}{H} \defeq \cloner{\clonern{n-1}{H}}.
\end{equation*}
if $n\ge 1$, and $\clonern{0}{H} \defeq H$. 

A similar strategy as the one above for iterated lampshufflers allows one to prove the next statement. 
\begin{corollary}
Let $p\ge 1$. Let $H$ be a finitely generated amenable group whose $\ell^{p}-$isoperimetric profile $\profp{H}$ satisfies Assumption~$(\star)$. Suppose that 
\begin{equation*}
    \profp{H}\left(\sqrt{\ln(x)}\right) \simeq \profp{H}(\ln(x)) \; \text{and} \; \profp{H}(x) \simeq \prof{H}(x).
\end{equation*}
Then, we have 
\begin{equation*}
   \profp{\clonern{n}{H}}(x) \simeq \profp{H}(\ln^{\circ n}(x))
\end{equation*}
for all $n\ge 1$. \qed
\end{corollary}

Finally, Corollary~\ref{cor:profCloner} motivates a similar question as Question~\ref{q:profile}.

\begin{question}
Is it true that all finitely generated amenable groups which do not have polynomial growth satisfy $\profp{H}\left(\sqrt{\ln(x)}\right) \simeq \profp{H}(\ln(x))$?
\end{question}

\section{Applications to quasi-isometric classifications and regular maps}\label{sec:appQIandRegularMaps} 

This section is dedicated to our applications about the existence of regular maps between halo products and their iterated versions. It relies on computations realized in Section~\ref{sec:computationsIsoProf}. In fact, for simplicity and conciseness, we will be focusing mainly on lampshufflers, but analogous statements can be derived for lampjugglers, lampdesigners and lampcloners.

Let us first distinguish iterated lampshufflers over amenable groups. 

\begin{corollary}\label{cor:shufflersoverslowprofiles}
Let $n,m\ge 1$. Let $H$ be a finitely generated amenable group. Assume that one of the following holds:
\begin{enumerate}[label=(\roman*)]
    \item\label{item:QIpolynomial} $H$ has polynomial growth of degree $d\ge 1$;
    \item\label{item:QInonpolynomial} 
    $\prof{H}$ satisfies Assumption~$(\star)$, $\prof{H}\left(\frac{\ln(x)}{\ln(\ln(x))}\right) \simeq \prof{H}(\ln(x))$ and the following property for any integers $k,\ell\ge 0$:
    \begin{equation*}
        \prof{H}(\ln^{\circ k}(x)) \simeq \prof{H}(\ln^{\circ \ell}(x)) \Longrightarrow k=\ell.
    \end{equation*}
\end{enumerate}
Then $\shufn{n}{H}$ and $\shufn{m}{H}$ are quasi-isometric if and only if $n=m$.
\end{corollary}

\begin{proof}
Suppose that $\shufn{n}{H}$ and $\shufn{m}{H}$ are quasi-isometric. In particular, their isoperimetric profiles are asymptotically equivalent, and if $H$ has polynomial growth, we get
\begin{equation*}
    \left(\frac{\ln^{\circ n}(x)}{\ln^{\circ (n+1)}(x)}\right)^{\frac{1}{d}} \simeq \left(\frac{\ln^{\circ m}(x)}{\ln^{\circ (m+1)}(x)}\right)^{\frac{1}{d}}
\end{equation*}
by Proposition~\ref{prop:profileofshufnofpolynomialgrowthgroups}, which in turn implies $n=m$. If we are in case~\textit{\ref{item:QInonpolynomial}}, then by Proposition~\ref{prop:profilesofshufn}, $\shufn{n}{H}$ has $\ell^{1}-$profile $\simeq \prof{H}(\ln^{\circ n}(x))$ and $\shufn{m}{H}$ has $\ell^{1}-$profile $\simeq \prof{H}(\ln^{\circ m}(x))$. Thus $n=m$ using our assumption, and we are done. 
\end{proof}

In practice, assumptions of~\textit{\ref{item:QInonpolynomial}} are easy to check. It holds for instance for any amenable group whose profile is of the form $\prof{H}(x)\simeq \left (\ln^{\circ k}(x)\right )^{\alpha}$ for $\alpha>0$ and $k\ge 0$, such as solvable Baumslag-Solitar groups or lamplighters over polynomial growth groups.

In fact, the isoperimetric profile being monotonous under regular maps between finitely generated amenable groups (cf. Theorem~\ref{thm:profilesmonotonuousregularmaps}), we get more generally:

\begin{corollary}
Let $n,m\ge 1$. Let $H$ be a finitely generated amenable group whose isoperimetric profile $\prof{H}$ satisfies Assumption~$(\star)$. Suppose that $\prof{H}\left(\frac{\ln(x)}{\ln(\ln(x))}\right) \simeq \prof{H}(\ln(x))$ and the following property holds for any integers $k,\ell\ge 0$:
\begin{equation*}
    \prof{H}(\ln^{\circ \ell}(x)) \preccurlyeq \prof{H}(\ln^{\circ k}(x)) \Longrightarrow k\leq\ell.
\end{equation*}
Then there exists a regular map from $\shufn{n}{H}$ to $\shufn{m}{H}$ if and only if $n\le m$.\qed
\end{corollary}

\subsection{Proofs of Theorem~\ref{thm:IteratedShufflersPolynomialQIBiLip intro} and Corollary~\ref{cor:corV.J}}

We have similar consequences at the other side of the spectrum:

\begin{corollary}\label{cor:IteratedShufflersPolynomialQIBiLip}
Let $n,m\ge 0$. Let $A$ and $B$ be infinite virtually abelian finitely generated groups, with growth degrees $a$ and $b$ respectively. Then the following are equivalent:
\begin{enumerate}[label=(\roman*)]
    \item\label{item:1QIbilip} $\shufn{n}{A}$ and $\shufn{m}{B}$ are quasi-isometric.
    \item\label{item:2QIbilip} $n=m$ and $a=b$.
    \item\label{item:3QIbilip} $\shufn{n}{A}$ and $\shufn{m}{B}$ are biLipschitz equivalent.
\end{enumerate}
\end{corollary}

\begin{proof}
The implication~\textit{\ref{item:3QIbilip}} $\Longrightarrow$~\textit{\ref{item:1QIbilip}} is obvious. 

\noindent We prove~\textit{\ref{item:1QIbilip}} $\Longrightarrow$~\textit{\ref{item:2QIbilip}}. Assume that $\shufn{n}{A}$ and $\shufn{m}{B}$ are quasi-isometric, so that they have asymptotically equivalent isoperimetric profiles. By Proposition~\ref{prop:profileofshufnofpolynomialgrowthgroups}, we then have 
\begin{equation}\label{eq:ComparisonProfile}
    \left(\frac{\ln^{\circ n}(x)}{\ln^{\circ (n+1)}(x)}\right)^{\frac{1}{a}} \simeq \left(\frac{\ln^{\circ m}(x)}{\ln^{\circ (m+1)}(x)}\right)^{\frac{1}{b}}
\end{equation}
and taking the logarithm, it follows that 
\begin{equation*}
    \ln\left(\frac{\ln^{\circ n}(x)}{\ln^{\circ (n+1)}(x)}\right) \simeq \ln\left(\frac{\ln^{\circ m}(x)}{\ln^{\circ (m+1)}(x)}\right).
\end{equation*}
The left-hand side is equivalent to $\ln^{\circ (n+1)}(x)$ and the right-hand side is equivalent to $\ln^{\circ (m+1)}(x)$, so that $n+1=m+1$, i.e. $n=m$. Re-injecting this information in~\eqref{eq:ComparisonProfile} now implies that $a=b$. 

\noindent If $n=m$ and $a=b$, then $A$ and $B$ are both biLipschitz equivalent to $\Z^{a}$ by Claim~\ref{claim:virtuallyabelianBILIPtoZ^d}, and thus $A$ and $B$ are biLipschitz equivalent. Thus, by~\cite[Lemma~8.8]{GT24a}, there is a biLipschitz equivalence from $\shuf{A}$ to $\shuf{B}$. Iterating this, we get a biLipschitz equivalence 
\begin{equation*}
    \shufn{n}{A}\longrightarrow \shufn{n}{B}
\end{equation*}
as claimed. This shows~\textit{\ref{item:2QIbilip}} $\Longrightarrow$~\textit{\ref{item:3QIbilip}} and concludes the proof.
\end{proof}

\begin{remark}
For the broader class of virtually nilpotent groups, some implications still hold and some may fail. For instance,~\textit{\ref{item:1QIbilip}} $\Longrightarrow$~\textit{\ref{item:2QIbilip}} remains true, but the converse is false. For instance, $\Z^4$ and the Heisenberg group $H$ over $\Z$ both have growth degree $4$, but $\shuf{\Z^4}$ and $\shuf{H}$ are not quasi-isometric by Corollary~\ref{cor:BILIPbetweenLampshufflers}, since $\Z^4$ and $H$ are not biLipschitz equivalent (e.g. they have different asymptotic dimensions).  
\end{remark}

For more general maps (e.g. regular maps), the isoperimetric profile is not sufficient to detect a constraint on polynomial growth degrees. However, asymptotic dimension does provide an inequality since, if $A$ has finite asymptotic dimension, then
\begin{equation*}
    \text{asdim}(\shuf{A})=\text{asdim}(A). 
\end{equation*}
Indeed, since $A$ is a subgroup of $\shuf{A}$, one has $\text{asdim}(A) \le \text{asdim}(\shuf{A})$, and on the other hand, since $\shuf{A}$ fits into a short exact sequence with kernel $\fsym{A}$, whose asymptotic dimension is $0$ as it is locally finite, and quotient $A$, one also has $\text{asdim}(\shuf{A}) \le \text{asdim}(A)$ (we refer the reader to the nice survey~\cite{BD08} for all these facts on asymptotic dimension). Iterating, we get 
\begin{equation*}
    \text{asdim}(\shufn{n}{A})=\text{asdim}(A)
\end{equation*}
for all $n\ge 0$. 

Note also that, if $A$ is virtually abelian, then its asymptotic dimension coincides with its growth degree. 

\begin{corollary}\label{cor:IteratedShufflersPolynomialRegularMap}
Let $n$ and $m$ be two natural integers. Let $A$ and $B$ be infinite virtually abelian finitely generated groups, with growth degrees $a$ and $b$ respectively. If there exists a regular map 
\begin{equation*}
    \shufn{n}{A}\longrightarrow \shufn{m}{B}
\end{equation*}
then $n\le m$ and $a\le b$. 
\end{corollary}

\begin{proof}
Assume that such a map exists. If $n=0$ there is nothing to prove, so we assume that $n\ge 1$. In this case, we cannot have $m=0$, because a group of exponential growth cannot regularly embed into a polynomial growth group. Hence $m\ge 1$ as well. Now, by Theorem~\ref{thm:profilesmonotonuousregularmaps} and Proposition~\ref{prop:profileofshufnofpolynomialgrowthgroups}, one has 
\begin{equation}\label{eq:ComparisonProfile2}
    \left(\frac{\ln^{\circ m}(x)}{\ln^{\circ (m+1)}(x)}\right)^{\frac{1}{b}} \preccurlyeq \left(\frac{\ln^{\circ n}(x)}{\ln^{\circ (n+1)}(x)}\right)^{\frac{1}{a}}
\end{equation}
and taking the logarithm implies 
\begin{equation*}
    \ln\left(\frac{\ln^{\circ m}(x)}{\ln^{\circ (m+1)}(x)}\right) \preccurlyeq \ln\left(\frac{\ln^{\circ n}(x)}{\ln^{\circ (n+1)}(x)}\right).
\end{equation*}
The left-hand side is $\simeq \ln^{\circ(m+1)}(x)$, and the right-hand side is $\simeq \ln^{\circ (n+1)}(x)$, so it follows that $n+1\le m+1$, i.e. $n\le m$. Additionally, since asymptotic dimension is monotonous under regular maps one gets 
\begin{equation*}
    a=\text{asdim}(A)=\text{asdim}(\shufn{n}{A})\le \text{asdim}(\shufn{m}{B})=\text{asdim}(B)=b
\end{equation*}
as claimed. 
\end{proof}

\begin{remark}\label{rm:extensionsfornilpotentgroups}
On the other hand, for more general polynomial growth groups $A$ and $B$, we can only conclude that the 
existence of a quasi-isometry between $\shufn{n}{A}$ and $\shufn{m}{B}$ imposes $n=m$, $a=b$ and $\text{asdim}(A)=\text{asdim}(B)$, and the existence of a regular map 
\begin{equation*}
    \shufn{n}{A}\longrightarrow \shufn{m}{B}
\end{equation*}
implies $n\le m$ and $\text{asdim}(A) \le \text{asdim}(B)$. 
\end{remark}

As an immediate consequence of Corollary~\ref{cor:IteratedShufflersPolynomialRegularMap}, we have the following.

\begin{corollary}\label{cor:QIRegularMap}
Let $n,m\ge 0$. Let $A$ and $B$ be infinite virtually abelian finitely generated groups, with growth degrees $a$ and $b$ respectively. Then the following are equivalent:
\begin{enumerate}[label=(\roman*)]
    \item the three equivalent assertions of Corollary~\ref{cor:IteratedShufflersPolynomialQIBiLip} are satisfied;
    \item there exist a regular map from $\shufn{n}{A}$ to $\shufn{m}{B}$, and a regular map from $\shufn{m}{B}$ to $\shufn{n}{A}$.
\end{enumerate}
\end{corollary}

Thus, asymptotic dimension is an obstruction to the existence of a regular map $\shuf{\Z^d}\longrightarrow \shuf{\Z^k}$ when $d>k$. Hence, in the spirit of~\cite[Question~5.4]{BST12}, a natural question arises: can we also rule out the existence of such maps if we increase the asymptotic dimension of the target space, for instance with a polynomial growth factor? It turns out that the answer is positive, and that the isoperimetric profile still gives an obstruction, whereas asymptotic dimension becomes inefficient.

Indeed, thanks to the following lemma, under the assumption that $\prof{H}\succcurlyeq \prof{G}$, we have general estimates on the isoperimetric profile of $G\times H$ in terms of $\prof{G}$ and $\prof{H}$.

\begin{lemma}\label{lem:profileofdirectproducts}
Let $G$ and $H$ be finitely generated amenable groups. If $\prof{H}(n)\succcurlyeq \prof{G}(n)$, then one has 
\begin{equation*}
    \prof{G}\left(\sqrt{n}\right) \preccurlyeq \prof{G\times H}(n) \preccurlyeq \prof{G}(n).
\end{equation*}
\end{lemma}

\begin{proof}
The estimate $\prof{G\times H}(n) \preccurlyeq \prof{G}(n)$ is a consequence of the fact that $G$ is a subgroup of $G\times H$ and Lemma~\ref{lem:monotonieprofilErschler}. Let us focus on the other inequality. Fix $n\in\N$ and subsets $A_{n}\subset G$, $B_{n}\subset H$ that realize $\prof{G}(n)$ and $\prof{H}(n)$ respectively, i.e. $|A_{n}|, |B_{n}|\le n$ and 
\begin{equation*}
    \prof{G}(n)=\frac{\left|A_{n}\right|}{\left|\partial_{G} A_{n}\right|},\; \prof{H}(n)=\frac{\left|B_{n}\right|}{\left|\partial_{H} B_{n}\right|}.
\end{equation*}
Then $A_{n}\times B_{n}\subset G\times H$ has cardinality $\le n^2$, and its boundary is given by 
\begin{equation*}
    \partial_{G\times H}(A_{n}\times B_{n}) = (\partial_{G}A_{n}\times B_{n})\cup(A_{n}\times \partial_{H}B_{n})
\end{equation*}
whence $\left|\partial_{G\times H}(A_{n}\times B_{n})\right| \le |\partial_{G}A_{n}|\cdot|B_{n}|+|A_{n}|\cdot|\partial_{H}B_{n}|$. Thus one gets 
\begin{equation*}
    \frac{\left|\partial_{G\times H}(A_{n}\times B_{n})\right|}{|A_{n}\times B_{n}|} \le \frac{|\partial_{G}A_{n}|\cdot|B_{n}|+|A_{n}|\cdot|\partial_{H}B_{n}|}{|A_{n}||B_{n}|}=\frac{|\partial_{G}A_{n}|}{|A_{n}|}+\frac{|\partial_{H}B_{n}|}{|B_{n}|}
\end{equation*}
and it follows that 
\begin{equation*}
    \prof{G\times H}(n^2) \ge \frac{|A_{n}\times B_{n}|}{\left|\partial_{G\times H}(A_{n}\times B_{n})\right|}\ge \frac{1}{\frac{|\partial_{G}A_{n}|}{|A_{n}|}+\frac{|\partial_{H}B_{n}|}{|B_{n}|}}=\frac{1}{\frac{1}{\prof{G}(n)}+\frac{1}{\prof{H}(n)}}.
\end{equation*}
Now, using Lemma~\ref{lem:UsefulLemmaProfile}, there exists a positive constant $K>0$ such that
\begin{equation*}
    \prof{G\times H}(x)\succcurlyeq\frac{1}{\frac{1}{\prof{G}\left(K\sqrt{x}\right)}+\frac{1}{\prof{H}\left(K\sqrt{x}\right)}}=\frac{1}{\frac{1}{\prof{G}\left(\sqrt{K^2x}\right)}+\frac{1}{\prof{H}\left(\sqrt{K^2x}\right)}}\succcurlyeq \frac{1}{\frac{1}{\prof{G}\left(\sqrt{x}\right)}+\frac{1}{\prof{H}\left(\sqrt{x}\right)}}
\end{equation*}
for all $x$ large enough. By assumption, $\prof{H}(n)\succcurlyeq \prof{G}(n)$, so that 
\begin{equation*}
    \prof{G\times H}(\cdot)\succcurlyeq \frac{1}{\frac{2}{\prof{G}(\sqrt{\cdot})}} \simeq \prof{G}\left(\sqrt{\cdot}\right)
\end{equation*}
as claimed.
\end{proof}

Thus, if additionally the isoperimetric profile of $G$ satisfies $\prof{G}\left(\sqrt{\cdot}\right)\simeq \prof{G}(\cdot)$, then
\begin{equation*}
    \prof{G\times H}(x)\simeq \prof{G}(x). 
\end{equation*}
This happens for many groups $G$ that have slow enough profiles, for instance:
\begin{itemize}
    \item Any polycyclic group with exponential growth, and more generally any GES group with exponential growth~\cite[Corollary~5]{Tes13};
    \item $F\wr\Sigma$, or $\shuf{\Sigma}$, where $F$ is finite and $\Sigma$ has polynomial growth;
    \item $\shufn{n}{H}$, where $H$ has profile $\prof{H}(x)\simeq \left (\ln^{\circ k}(x)\right )^{\alpha}$, for some $\alpha>0$ and integer $k\ge 1$.
\end{itemize}

As a concrete example, we have for instance: 

\begin{corollary}
Let $d,k,p\ge 1$ be three integers. There exists a regular map 
\begin{equation*}
    \shuf{\Z^d} \longrightarrow \Z^{p}\times\shuf{\Z^k}
\end{equation*}
if and only if $d\le k$.\qed
\end{corollary}

\subsection{Proof of Proposition~\ref{prop:IteratedShufflerintoIteratedLL}}

Finally, we want to compare lampshufflers to lamplighters. We already know from Proposition~\ref{prop:A2} that a wreath product over a subgroup of $H$ coarsely embeds into $\shuf{H}$. In~\cite[Corollary~7.17]{GT24a}, given two groups $H$ and $G$ satisfying some mild assumptions, it is proved that $\shuf{H}$ does not quasi-isometrically or coarsely embed into a lamplighter $E\wr G$. Our computations allow us to prove an iterated version of this result for free abelian groups: there is no regular map 
\begin{equation*}
\shufn{n}{\Z^d}\longrightarrow \Z/2\Z\wr(\Z/2\Z\wr(\dots(\Z/2\Z\wr \Z^d)))
\end{equation*}
where the wreath product is iterated $n$ times. In fact, we have more generally the following.

\begin{corollary}\label{cor:iteratedshufintoiteratedlamplighter}
Let $G$ and $H$ be finitely generated amenable groups. Suppose that there is a regular map 
\begin{equation*}
\shufn{n}{H}\longrightarrow \Z/2\Z\wr\big(\Z/2\Z\wr(\dots\wr(\Z/2\Z\wr G))\big)
\end{equation*}
where the wreath product is iterated $n$ times, $n\ge 1$. Then the following holds.
\begin{enumerate}[label=(\roman*)]
      \item\label{item:1shufintolamp} If $H$ has polynomial growth of degree $d\ge 1$, then
      \begin{equation*}
         \prof{G}(x)\preccurlyeq\frac{\prof{H}(x)}{\ln(x)^{\frac{1}{d}}};
      \end{equation*}
      \item\label{item:2shufintolamp} If $\prof{H}$ satisfies Assumption ~$(\star)$ and $\prof{H}\left(\frac{\ln(x)}{\ln(\ln(x))}\right)\simeq\prof{H}(\ln(x))$, then
      \begin{equation*}
        \prof{G}(x)\preccurlyeq\prof{H}(x).
      \end{equation*}
\end{enumerate}
\end{corollary}

\begin{proof}
In case~\textit{\ref{item:1shufintolamp}}, we get
\begin{equation*}
        \prof{G}(\ln^{\circ n}(x))\preccurlyeq\left(\frac{\ln^{\circ n}(x)}{\ln^{\circ (n+1)}(x)}\right)^{\frac{1}{d}}=\frac{\prof{H}(\ln^{\circ n}(x))}{\left(\ln^{\circ (n+1)}(x)\right)^{\frac{1}{d}}}.
\end{equation*}
Using Assumption~$(\star)$ for the logarithm and for $\prof{H}$, we have
\begin{equation*}
        \prof{G}(\ln^{\circ n}(x))=O\left(\frac{\prof{H}(\ln^{\circ n}(x))}{\left(\ln^{\circ (n+1)}(x)\right)^{\frac{1}{d}}}\right)
\end{equation*}
and the change of variable $x'=\ln^{\circ n}(x)$ gives the result. In case~\textit{\ref{item:2shufintolamp}}, we get rather 
\begin{equation*}
    \prof{G}(\ln^{\circ n}(x))=O\big(\prof{H}(\ln^{\circ n}(x))\big)
\end{equation*}
and we conclude similarly.
\end{proof}

We get the following consequence in the case of polynomial growth groups.

\begin{corollary}\label{cor:iteratedshufintoiteratedlamplighterPolynomial}
Let $d,k,n\ge 1$, and let $G$ and $H$ be finitely generated groups of polynomial growth with growth degrees $k$ and $d$ respectively. If there is a regular map 
\begin{equation*}
\shufn{n}{H}\longrightarrow \Z/2\Z\wr\big(\Z/2\Z\wr(\dots\wr(\Z/2\Z\wr G))\big)
\end{equation*}
then $d<k$, where the wreath product is iterated $n$ times. 
\end{corollary}

Note that this statement cannot be reached with methods from~\cite{GT24a}, even for quasi-isometric or coarse embeddings, since the thick bigon property used in~\cite{GT24a} is not stable under iterations of lampshufflers.

\begin{proof}
By Corollary~\ref{cor:iteratedshufintoiteratedlamplighter}, we have $x^{\frac{1}{k}}\preccurlyeq\left(\frac{x}{\ln(x)}\right)^{\frac{1}{d}}$, which immediately implies $d<k$.
\end{proof}

Note that, if $G$ is a proper subgroup of $H$, then an iteration of Proposition~\ref{prop:A2} ensures that $\shufn{n}{H}$ contains $\Z/2\Z\wr\big(\Z/2\Z\wr(\dots(\Z/2\Z\wr G))\big)$ (iterated $n$ times) as a subgroup, and thus we get a regular map 
\begin{equation*}
    \Z/2\Z\wr\big(\Z/2\Z\wr(\dots\wr(\Z/2\Z\wr G))\big) \longrightarrow \shufn{n}{H}. 
\end{equation*}

\section{Appendix: Lamplighter subgroups in lampshu{f}flers}\label{appendixA}

In the article, the strategy for getting optimal upper bounds on the isoperimetric profile of halo products is to find subgraphs that are quasi-isometric to lamplighter graphs. In this appendix, our aim is to prove that, under additional mild algebraic assumptions on the base group, we can find lamplighter~\textit{subgroups} inside lampshufflers.

For other halo products, such as lampjugglers, lampdesigners and lampcloners, we already observed in Section~\ref{sec:halo} that they contain wreath products over the same base group as subgroups. 

Let us recall the following definition: we say that a non-decreasing map $h\colon\R_{+}\rightarrow \R_{+}$ satisfies Assumption~$(\star)$ if 
\begin{equation*}
    \forall C>0,\; h(Cx)=O(h(x)). 
\end{equation*} 

The motivation for this strategy comes from the following result, due to Silva. In the upcoming result (Proposition~\ref{prop:A2}), we will use the main idea of its proof. 

\begin{proposition}[{\cite[Proposition~2.3]{Sil24}}]\label{prop:A1}
Let $H$ be an infinite non-co-Hopfian group. Then, for any finite group $F$, $\shuf{H}$ has a subgroup isomorphic to $F\wr H$. 
\end{proposition}

Recall that a group $H$ is \textit{co-Hopfian} if any injective morphism $H\rightarrow H$ is also surjective. Equivalently, a group is co-Hopfian if it has no proper subgroup isomorphic to itself.

Many amenable groups are known to be non-co-Hopfian, including:
\begin{itemize}
    \item $\Z^d$, $d\ge 1$, and more generally any finitely generated abelian group;
    \item Some torsion-free nilpotent groups, such as the Heisenberg group over the integers~\cite{Corn16};
    \item Solvable Baumslag-Solitar groups $\text{BS}(1,n)$, $n\ge 1$~\cite{NP11};
    \item Wreath products $N\wr G$ where at least one of the two groups is not co-Hopfian~\cite[Proposition~6.1]{BFF24};
    \item Houghton's groups $H_{n}$, $n\ge 2$~\cite[Theorem~7.1]{BCMR16};
    \item The Grigorchuk's group~\cite{Lys85}.
\end{itemize}
In the non-amenable side, they also include for instance non-abelian free groups~\cite[Section III.22]{dlH00} or right-angled Artin groups~\cite{Cas16}.

Thus, using Lemma~\ref{lem:monotonieprofilErschler}, it directly follows that for $H$ an amenable and non-co-Hopfian group, one has already
\begin{equation*}
    \profp{\shuf{H}}(x) \preccurlyeq \prof{H}(\ln(x))
\end{equation*}
when $\prof{H}$ satisfies Assumption~($\star$).

In fact, we can derive from Silva's proof the following more general result.

\begin{proposition}\label{prop:A2}
Let $H$ be a group. If $K$ is a proper subgroup of $H$, then $\shuf{H}$ contains a subgroup isomorphic to $\fsym{[H:K]}\wr K$. 
\end{proposition}

\begin{proof}
Denote $[H:K]\defeq m \in \lbrace 2,3,\dots\rbrace\cup\lbrace \infty\rbrace$, and let $S\subset H$ be a set of representatives of $H/K$, so that $|S|=m$ and we have a partition 
\begin{equation*}
    H=\bigsqcup_{k\in K}kS.
\end{equation*}
Consider then 
\begin{equation*}
    G \defeq \left\lbrace (\sigma, h) \in\shuf{H} : h\in K, \;\sigma(k'S)=k'S\; \text{for all $k'\in K$}\right\rbrace.
\end{equation*}
It is not hard to check that $G$ is a subgroup of $\shuf{H}$, and given $\sigma\in\fsym{H}$ satisfying $\sigma(k'S)=k'S$ for all $k'\in K$, we can define a map 
\begin{align*}
    f_{\sigma}\colon
    \begin{array}{cll}    
    K &\longrightarrow &\fsym{m} \\
    k&\longmapsto &(k^{-1}\cdot\sigma)\big|_{S}
    \end{array}.
\end{align*}
Then, a direct computation shows that $f_{k\cdot \sigma}=k\cdot f_{\sigma}$ and $f_{\sigma\circ\tau}=f_{\sigma}f_{\tau}$ for any $k\in K$ and any $\sigma,\tau\in\fsym{H}$ satisfying $\sigma(k'S)=k'S$ for all $k'\in K$. This implies that the map 
\begin{align*}
    \begin{array}{ccl} G &\longrightarrow &\fsym{m}\wr K \\
    (\sigma, h)&\longmapsto &(f_{\sigma},h)
    \end{array}
\end{align*}
is a group isomorphism. This concludes the proof. 
\end{proof}

In particular, we can characterize immediately when lampshufflers are residually finite:

\begin{corollary}\label{LampshufflersarenotRF}
Let $H$ be a group. Then $\shuf{H}$ is residually finite if and only if $H$ is finite. 
\end{corollary}

Recall here that a group $G$ is~\textit{residually finite} if any non-trivial element of $G$ stays non-trivial in a finite quotient of $G$. 

\begin{proof}
If $H$ is finite, $\shuf{H}$ is finite and therefore residually finite. Conversely, assume that $H$ is infinite. Since $H$ has a proper subgroup $K\leqslant H$ of index $\ge 3$, $\shuf{H}$ contains a subgroup isomorphic to $\fsym{[H:K]}\wr K$, and the latter is not residually finite. As residual finiteness passes to subgroups, $\shuf{H}$ cannot be residually finite either. 
\end{proof}

One can also recover Corollary~\ref{cor:upperboundonProf} for lampshufflers:

\begin{corollary}\label{cor:A3}
Let $p\ge 1$. Let $H$ be an amenable finitely generated group with a finitely generated proper subgroup $K$ such that $\profp{H}(x)\simeq \profp{K}(x)$. If $\prof{H}$ satisfies Assumption~$(\star)$, then one has 
\begin{equation*}
    \profp{\shuf{H}}(x) \preccurlyeq \profp{H}(\ln(x)).
\end{equation*}
\end{corollary}

The corollary applies, of course, to every non-co-Hopfian group, but more generally it also applies to any finitely generated group $H$ having at least one finite-index subgroup $K$. One algebraic criterion to ensure the latter is to be \textit{non-perfect}.

\begin{definition}
Let $H$ be a group. We say that $H$ is \textit{perfect} if $[H,H]=H$, where $[H,H]$ is the commutator subgroup of $H$.
\end{definition}

Examples of perfect groups include finite alternating groups $A_{n}$ for $n\ge 5$ and linear groups $\text{SL}(n,K)$ for $n\ge 3$ and non-commutative fields $K$. Among the three Thompson's groups $F\subset T\subset V$, $T$ and $V$ are simple, thus perfect, as well as the commutator subgroup $[F,F]$ of $F$~\cite{CFP96}. 

As mentioned above, we are in fact interested in groups that are not perfect, due to the next statement.

\begin{lemma}\label{lm:finiteindexsubgroupsinnonperfectgroups}
Let $H$ be a finitely generated group which is not perfect. Then $H$ has a proper finite-index subgroup.
\end{lemma}

\begin{proof}
As $H$ is not perfect, $[H,H]$ is a proper subgroup and the quotient $H/[H,H]$ is a non-trivial abelian finitely generated group. It has therefore a proper finite-index subgroup, and lifting the latter provides a proper finite-index subgroup for $H$, that contains $[H,H]$.  
\end{proof}

On the amenable side, the class of non-perfect groups is huge. It includes for instance
\begin{itemize}
    \item All solvable groups, in particular nilpotent and polycyclic groups;
    \item Lamplighters, lampshufflers and lampjugglers over amenable non-perfect groups;
    \item More generally, any semi-direct product $N\rtimes Q$ where $N$ and $Q$ are amenable and $Q$ is non-perfect, as well as any subgroup of $N\rtimes Q$;
    \item Houghton's groups $H_{n}$, $n\ge 2$.
\end{itemize}

Henceforth, for an amenable non-perfect group $H$, Corollary~\ref{cor:A3} directly provides an upper bound on the $\ell^{p}-$isoperimetric profiles of $\shuf{H}$.

We emphasize here that Corollary~\ref{cor:A3} can also cover cases that are not covered by Proposition~\ref{prop:A1}, since they are indeed examples of finitely generated torsion-free nilpotent groups that are co-Hopfian~\cite{Bel03}. 

It is also worth mentioning that there do exist amenable perfect groups. Such groups have been constructed by Juschenko and Monod in~\cite{JM13}, and are even simple. Thus Corollary~\ref{cor:A3} does not apply to them, and on the other hand it seems to be an open problem whether they are co-Hopfian or not; see the list of open questions in~\cite{Corn14}. Note also that we do not know the asymptotic behaviour of the isoperimetric profiles of these groups. 

\begin{remark}\label{rem:Preservation}
As pointed out above, lampshuffler groups over non-perfect groups are themselves non-perfect. This simply follows from the set inclusions 
\begin{equation*}
    [\shuf{H},\shuf{H}]\subset \fsym{H}\times [H,H] \subsetneq \shuf{H}.
\end{equation*}
In fact, it is true more generally that a lampshuffler $\shuf{H}$ over a non-trivial group $H$ is never perfect, regardless of the perfectness of $H$. Indeed, since the action of $H$ on $\fsym{H}$ preserves the parity of the permutations, the commutator subgroup of $\shuf{H}$ is contained in $\mathcal{A}(H)\times [H,H]$, where $\mathcal{A}(H)\subset\fsym{H}$ is the set of even permutations. Since $\mathcal{A}(H)$ is not equal to $\fsym{H}$ if $|H|\ge 3$, and since $[H,H]$ is trivial if $|H|=2$, we thus see that $\shuf{H}$ is not equal to its commutator subgroup. 
\end{remark}

Lastly, regarding lampshufflers, non co-Hopficity can also be used to deduce the isoperimetric profile of iterated lampshufflers, since it is stable under iteration of lampshufflers:

\begin{proposition}\label{prop:A7}
If $H$ is not co-Hopfian, then $\shuf{H}$ is not co-Hopfian. 
\end{proposition}

\begin{proof}
If $H$ is not co-Hopfian, fix an injective morphism $\psi\colon H\rightarrow H$ which is not surjective, and define the map 
\begin{align*}
\varphi\colon
\begin{array}{cll}
\shuf{H}&\longrightarrow &\shuf{H} \\
(\sigma, g)&\longmapsto &(\overline{\sigma}, \psi(g))
\end{array}
\end{align*}
where, for any $\sigma\in \fsym{H}$, $\overline{\sigma}\in\fsym{H}$ is defined as 
\begin{align*}
    \overline{\sigma}\colon
    \begin{array}{cll}
    H&\longrightarrow &H \\
    g&\longmapsto &\begin{cases}\psi(\sigma(\psi^{-1}(g))) &\mbox{if $g$ is in the image of $\psi$} \\ g &\mbox{otherwise}\end{cases}
    \end{array}.
\end{align*}
Then one checks directly that the correspondence $\sigma\longmapsto\overline{\sigma}$ is well-defined and that the following two properties hold:
\begin{enumerate}[label=(\roman*)]
    \item $\overline{\sigma\circ \tau}=\overline{\sigma}\circ\overline{\tau}$ for any $\sigma,\tau\in\fsym{H}$;
    \item $\overline{p\cdot\sigma}=\psi(p)\cdot\overline{\sigma}$ for any $p\in H$ and $\sigma\in \fsym{H}$.
\end{enumerate}
These two points imply that $\varphi$ is a morphism, which is injective since $\psi$ is, and which is not surjective since $\psi$ is not.
\end{proof}

Going further, it is also natural to ask whether these algebraic assumptions are preserved under iteration of other halo products. First, it is immediate that a halo product over a non-perfect group is not perfect, since the definition of the product law of $\halo H$ directly gives
\begin{equation*}
    [\halo H,\halo H]\subset L(H)\times [H,H] \subsetneq \halo H,
\end{equation*}
as we pointed out in Remark~\ref{rem:Preservation}, in the special case of lampshufflers. In this remark, we also pointed out that lampshufflers were in fact never perfect, using the parity of the permutations, which is invariant by the action of the base group. Thus, for a general halo group $\halo H$, it is enough to find a non-trivial morphism from $L(H)$ to an abelian group, which is invariant by the action of the base group. For instance, for lampcloners, we can use the determinant of linear maps.

Regarding preservation of non-co-hopficity under iterations of halo products, we can only find proofs in concrete examples. For lampcloners, we can adapt the above proof for lampshufflers. Indeed, let us first note that a group morphism $\psi\colon H\rightarrow H$ gives rise to a linear map $\tilde{\psi}\colon V_{H}\rightarrow V_{H}$, by permuting the vector of the canonical basis given by $H$, and if $\psi$ is injective but not bijective, then so is $\tilde{\psi}$. Finally, it remains to
define the map 
\begin{align*}
\varphi\colon
\begin{array}{cll}
\cloner{H}&\longrightarrow &\cloner{H} \\
(\sigma, g)&\longmapsto &(\overline{\sigma}, \psi(g))
\end{array}
\end{align*}
where, for any $\sigma\in \text{FGL}(H)$, $\overline{\sigma}\in \text{FGL}(H)$ is defined as 
\begin{align*}
    \overline{\sigma}\colon
    \begin{array}{cll}
    V_H&\longrightarrow &V_H \\
    v&\longmapsto &\begin{cases}\tilde{\psi}(\sigma(\tilde{\psi}^{-1}(v))) &\mbox{if $v\in V_{H}$ is in the image of $\tilde{\psi}$} \\ v &\mbox{otherwise}\end{cases}
    \end{array},
\end{align*}
and to check that $\varphi$ is an injective, but not bijective, group morphism.

\chapter{On quantitative orbit equivalence for lamplighter-like groups}\label{chap:chapter6}

This chapter presents the paper~\cite{cordum26}, written in collaboration with Corentin Correia. 

\vspace{0.3cm}

\minitoc

\section{Introduction}

\setcounter{theoremletter}{0}

In this paper, we focus on the following class of groups, introduced in~\cite{GT24a}:

\begin{definition}
Let $X$ be a set. A~\textit{halo of groups $\halo$ over $X$} is the data, for any subset $S\subset X$, of a group $L(S)$ such that:
\begin{itemize}
    \item for all $R,S\subset X$, if $R\subset S$ then $L(R)\leqslant L(S)$;
    \item $L(\emptyset)=\lbrace 1\rbrace$ and $L(X)=\langle L(S) : S\subset X \;\text{finite}\rangle$;
    \item for all $R,S\subset X$, $L(R\cap S)=L(R)\cap L(S)$.
\end{itemize}
\end{definition}

Given an action $H\curvearrowright X$ and a morphism $\alpha\colon H \longrightarrow \text{Aut}(L(X))$ satisfying $\alpha(h)(L(S))=L(hS)$ for any $S\subset X$ and $h\in H$, the~\textit{permutational halo product} $\halo_{X,\alpha}H$ is the semi-direct product 
\begin{equation*}
    \halo_{X,\alpha}H \defeq L(X)\rtimes_{\alpha}H.
\end{equation*}
In the case where $X=H$ and $H$ acts on itself by left-multiplication, we simply denote $\halo H$ the corresponding halo product. 

\sloppy Examples of halo products include for instance wreath products, but also~\textit{lampshufflers}, defined as 
\begin{equation*}
    \shuf{H}=\fsym{H}\rtimes H
\end{equation*}
and for which $L(S)=\fsym{S}$. These groups already appeared several times in the literature, in relations with many topics of interest in group theory, see for instance~\cite{Yad09, HO16, BZ19, EZ21, SCZ21, GT24a, Sil24}.

As another instance of this construction, we can also mention~\textit{lampcloners}, of the form
\begin{equation*}
    \cloner{H}=\text{FGL}(H)\rtimes H
\end{equation*}
with a field $\field$, for which $L(S)=\text{FGL}(S)$, where $S$ denotes any subset of $H$. We refer the reader to Section~\ref{sec:halo} for other examples. 

The motivation of~\cite{GT24a} to introduce such a general framework is that the semi-direct product structure provides a foliation of these spaces that must be, if $H$ satisfies additional mild assumptions, “quasi-preserved” by quasi-isometries, allowing the authors to show strong rigidity phenomena for quasi-isometries between such spaces, and thus extending the classification already obtained in~\cite{GT24b}.

The starting point of this article is therefore the following, rather vague question: 

\begin{question}\label{question:startingquestion}
In the case where two halo products $\halo H$ and $\halo K$ are not quasi-isometric, how much do their geometries differ?
\end{question}

As explained in Section~\ref{sec:OE}, quantitative orbit equivalence turns out to be a suitable notion to compare groups, as it captures in a non-trivial way some important geometric information, such as the volume growth or the isoperimetric profile. Moreover, it is an important result of Shalom that, for amenable groups, $\ld^{\infty}-$orbit equivalence coincides with the notion of biLipschitz equivalence, and the weaker notion of mutual cobounded $\ld^{\infty}-$measure equivalence coincides with the corresponding weakening of biLipschitz equivalences, namely quasi-isometries~\cite{Sha04}. Additionally, it is an important observation contained (at least implicitly) in~\cite{GT24a} that, for many halo products $\halo$, a quasi-isometry between $\halo H$ and $\halo K$ can often be upgraded to a biLipschitz equivalence. This interplay can be summed up as:

\vspace{0.15cm}

\begin{figure}[H]
  \centering
  \includegraphics[width=0.9\linewidth, trim=0 0 0 4mm, clip]{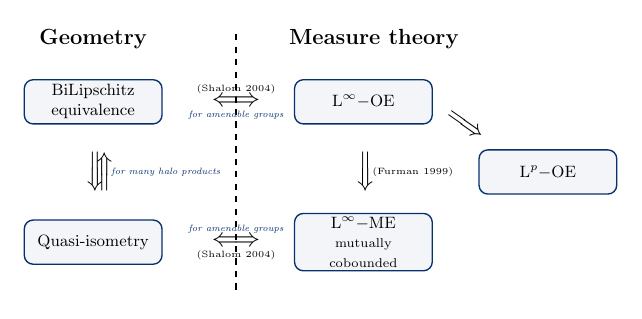}
  \caption{Geometric and measure-theoretic comparison of geometries of groups}
  \label{fig:cayley}
\end{figure}

In particular, if $\halo H$ and $\halo K$ are quasi-isometric, then they are $\ld^p-$orbit equivalent for all $p>0$.

This observation allows us to formulate a more precise version of the question we want to answer in this chapter:

\begin{question}
In the case where two halo products $\halo H$ and $\halo K$ are not quasi-isometric, can they still be $\ld^p-$orbit equivalent for some values of $p>0$?
\end{question}

The goal of this chapter is therefore to explain how to construct optimal orbit equivalence couplings between halo products.

\paragraph{Orbit equivalence couplings for halo products.} In~\cite[Corollary~7.3]{DKLMT22} is established a stability result for the existence of a coupling between two lamplighters $\Lambda\wr H$ and $\Lambda\wr K$ provided a coupling between the base groups $H$ and $K$. See Theorem~\ref{thm:dklmtWreath} for a more precise statement.

The first contributions of our article are similar stability results for other examples of halo products, with a preservation of quantification. 

\begin{theoremletter}[see Theorem~\ref{thm:stabilityofcouplings+quantificationLampjugglers}]\label{thm:stabilityofcouplings+quantificationLampjugglersINTRO}
If two groups $H$ and $K$ are orbit equivalent, then $\shuf{H}$ and $\shuf{K}$ are orbit equivalent. Moreover, given non-decreasing maps $\varphi,\psi\colon\R_{+}\rightarrow\R_{+}$, if $H$ and $K$ are finitely generated and if there exists a $(\varphi,\psi)-$integrable orbit equivalence coupling from $H$ to $K$, then the same holds from $\shuf{H}$ to $\shuf{K}$.
\end{theoremletter}

\sloppy This theorem is proved by constructing an explicit and natural free p.m.p. action of $\shuf{H}$ on a standard probability space, provided such an action for $H$. In fact, we explain in Section~\ref{sec:generalmethodstability} a general method to construct a free p.m.p. action of a halo $\halo H$ on a standard probability space given such an action of $H$. We therefore get similar results for lampjugglers, lampcloners and lampdesigners.

As a consequence of Theorem~\ref{thm:stabilityofcouplings+quantificationLampjugglersINTRO}, given $k>d$, there exists an $\ld^p$ orbit equivalence between $\shuf{\Z^d}$ and $\shuf{\Z^k}$ for every $p<\frac{d}{k}$. By induction, the same holds for $\shufn{n}{\Z^d}$ and $\shufn{n}{\Z^k}$, where $\shufn{n}{H}$ denotes an iterated lampshuffler, defined inductively by 
\begin{equation*}
    \shufn{0}{H} \defeq H, \; \shufn{n}{H} \defeq \shuf{\shufn{(n-1)}{H}}.
\end{equation*}

In~\cite{DKLMT22} is also introduced another technique, that relies on~\textit{F\o lner tiling sequences}, to construct couplings between iterated wreath products but with different numbers of iterations. We therefore also prove a stability property for the existence of F\o lner tiling sequences of lampshufflers. 

\begin{theoremletter}
Let $H$ be a finitely generated amenable group. If $H$ has a F\o lner tiling sequence, then so has $\shuf{H}$.
\end{theoremletter}

See Theorem~\ref{thm:folnertilingshuffler} for a detailed statement and the additional properties satisfied by this tiling. Here also we shall notice that the technique is not specific to lampshufflers, and can be applied in a similar fashion to other instances of halo products.

As an application, we construct and quantify precisely couplings between iterated lampshufflers, with different numbers of iterations:

\begin{theoremletter}[see Corollary~\ref{cor:couplingsbetweeniteratedshufflers}]\label{thm:couplingsbetweeniteratedshufflersINTRO}
Let $n,m\ge 0$ be natural integers such that $m>n$. Let $d,k\ge 1$. Then there exists an orbit equivalence coupling from $\shufn{m}{\Z^k}$ to $\shufn{n}{\Z^d}$, which is $(\mathrm{L}^{<\infty},\varphi_{m-n,k,\varepsilon})-$integrable for every $\varepsilon>0$, where 
\begin{equation*}
    \varphi_{i,k,\varepsilon}(x)\defeq\frac{\left (\ln^{\circ i}(x)\right )^{\frac{1}{k}}}{\left(\ln^{\circ (i+1)}(x)\right)^{1+\frac{1}{k}+\varepsilon}}.
\end{equation*}
\end{theoremletter}

For an integer $i\ge 0$, $\ln^{\circ i}$ denotes the $i-$th iteration of the logarithm, with the convention that $\ln^{\circ 0}$ is the identity.

Observe that Theorem~\ref{thm:ObstructionDKLMT}\textit{(i)} imposes an asymptotic upper bound on the integrability of cocycles of an orbit equivalence between two finitely generated groups. Therefore, to determine if our previous couplings are optimal, we need to know isoperimetric profiles of lampshuffler groups.

\paragraph{Optimality of orbit equivalence couplings.} We now go back to the optimality of the couplings provided by Theorem~\ref{thm:stabilityofcouplings+quantificationLampjugglersINTRO} and Theorem~\ref{thm:couplingsbetweeniteratedshufflersINTRO}. As a first application, using the profile of $\shuf{\Z^n}$ which was already known by~\cite{SCZ21} and~\cite{EZ21}, we obtain that, for integers $k>d$, $\shuf{\Z^k}$ and $\shuf{\Z^d}$ are $\ld^p$ orbit equivalent if and only if $p<\frac{d}{k}$. Our new computations of isoperimetric profiles enable us to get even more:

\begin{theoremletter}[see Theorem~\ref{thm:optimalitypolynomialgrowthJugglers}]\label{thm:optimalitypolynomialgrowthJugglersINTRO}
Let $k,d\ge 1$ be positive integers such that $k> d$. Let $n\ge 0$ be an integer.
Then $\shufn{n}{\Z^k}$ and $\shufn{n}{\Z^d}$ are $\ld^{p}$ orbit equivalent if and only if $p<\frac{d}{k}$.
\end{theoremletter}

\begin{theoremletter}[see Corollary~\ref{cor:optimalityIteratedLampshuffler}]
\sloppy Let $n,m\ge 1$ be two integers such that $m>n$. There exists a $\left(\left (\frac{\ln^{\circ (m-n)}}{\ln^{\circ (m-n+1)}}\right )^{p},\mathrm{L}^0\right)-$integrable orbit equivalence coupling from $\shufn{m}{\Z^k}$ to $\shufn{n}{\Z^d}$ if and only if $p<\frac{1}{k}$.
\end{theoremletter}

We emphasize here that optimality of many other couplings can be deduced from the computations of isoperimetric profiles. We refer to the end of Section~\ref{sec:folnertilingjugglers} for more details, more precisely Corollary~\ref{cor:boundOfIntegrability} and Remark~\ref{rem:FolnerWreathProduct}.

On the other hand, we do not have precise estimates on the isoperimetric profile of a lampcloner over a polynomial growth group, so we do not know if our couplings are quantitatively optimal. The following statement on lampcloners is thus weaker than our results on lampshufflers.

\begin{theoremletter}[see Theorem~\ref{thm:optimalitypolynomialgrowthCloners}]\label{thm:optimalitypolynomialgrowthClonersINTRO}
Let $\field$ be a finite field and let $k,d\ge 1$ be integers such that $k>d$. Then for every integer $n\ge 0$, for every $p<\frac{d}{k}$, $\clonern{n}{\Z^k}$ and $\clonern{n}{\Z^d}$ are $\ld^{p}$ orbit equivalent. Moreover, if $\cloner{\Z^d}$ and $\cloner{\Z^k}$ are $\ld^{p}$ orbit equivalent, then $p\le\frac{2d}{k}$.
\end{theoremletter}

\begin{theoremletter}[see Corollary~\ref{cor:couplingsbetweeniteratedjugglers}]\label{thm:couplingsbetweeniteratedjugglers INTRO}
Let $\field$ be a finite field. Let $n,m\ge 0$ be natural integers such that $m>n$. Let $d,k\ge 1$. Then there exists an orbit equivalence coupling from $\clonern{m}{\Z^k}$ to $\clonern{n}{\Z^d}$, which is $(\mathrm{L}^{<\infty},\varphi_{m-n,k,\varepsilon}(x))-$integrable for every $\varepsilon>0$, where 
\begin{equation*}
    \varphi_{i,k,\varepsilon}(x)\defeq\frac{\ln^{\circ i}(x)^{\frac{1}{2k}}}{\left(\ln^{\circ (i+1)}(x)\right)^{1+\varepsilon}}.
\end{equation*}
\end{theoremletter}

\paragraph{Plan of the chapter.} Section~\ref{sec:stability} is devoted to the construction of orbit equivalence couplings between lampjugglers and lampcloners, by building a natural action of the halo product provided an action of the base group. For lampjugglers, we deduce that they are quantitatively optimal, using precise estimates on the isoperimetric profile. In Section~\ref{sec:folnerTiling}, we will use F\o lner tiling sequences to get other orbit equivalence couplings, for instance between $\Z^d$ and $\juggler{r}{\Z^k}$. Finally, Section~\ref{sec:commentsandquestions} records various questions related to the article.

\section{Stability results for quantitative orbit equivalence}\label{sec:stability}

\subsection{General method}\label{sec:generalmethodstability}

\paragraph{Orbit equivalence coupling between lamplighters, by~\cite{DKLMT22}.} Let us start by explaining a construction of orbit equivalence between the first examples of halo products: lamplighter groups. Delabie, Koivisto, Le Maître and Tessera proved the following:

\begin{theorem}[{\cite[Corollary~7.3]{DKLMT22}}]\label{thm:dklmtWreath}
Let $\Lambda$ be a finite group. Let $H$ and $K$ be finitely generated groups and $\varphi,\psi\colon\R_{+}\rightarrow\R_{+}$ be increasing maps. Assume that there exists a $(\varphi,\psi)-$integrable orbit equivalence coupling from $H$ to $K$. Then there exists a $(\varphi,\psi)-$integrable orbit equivalence coupling from $\Lambda\wr H$ to $\Lambda\wr K$.
\end{theorem}

Their statement is even stronger. In fact, they show that $\Lambda_{1}\wr H$ and $\Lambda_{2}\wr K$ are $(\varphi,\psi)-$integrably orbit equivalent if there exists such a coupling from a finitely generated group $\Lambda_{1}$ to another finitely generated group $\Lambda_{2}$, and such a coupling from $H$ to $K$. Here we only focus on the case where $\Lambda\defeq\Lambda_{1}=\Lambda_{2}$ is a finite group, and we directly explain the construction they use within the appropriate formalism for our purposes. We refer the reader to~\cite[Section~7.1]{DKLMT22} where the more general notion of wreath product of equivalence relations is introduced to prove Theorem~\ref{thm:dklmtWreath}. 

Given a p.m.p. action $H\curvearrowright (X,\mu)$, we can define a p.m.p. action of the group $\Lambda\wr H$ on $(\Lambda^{H}\times X,{m_{\Lambda}}^{\otimes H}\otimes \mu)$, where $m_{\Lambda}$ is the uniform distribution on $\Lambda$, in the following way. Let us denote by $\mathds{1}\colon H\to\Lambda$ the constant map equal to $1_{\Lambda}$. Then, for every $h\in H$, every finitely supported $f\colon H\to\Lambda$ and every $((\ell_{g})_{g\in H},x)\in \Lambda^{H}\times X$, we set
\begin{equation*}
    (\mathds{1},h)\cdot \left((\ell_{g})_{g\in H},x\right)=\left((\ell_{h^{-1}g})_{g\in H},h\cdot x\right)
\end{equation*}
and
\begin{equation*}
    (f,1_{H})\cdot \left((\ell_{g})_{g\in H},x\right)=\left((f(g^{-1})\ell_{g})_{g\in H},x\right).
\end{equation*}
A direct computation shows that this provides a p.m.p. action $\Lambda\wr H\curvearrowright (\Lambda^{H}\times X,{m_{\Lambda}}^{\otimes H}\otimes \mu)$.

Now, given $(X,\mu)$ a $(\varphi,\psi)-$integrable orbit equivalence from $H$ to $K$, if $c_{H,K}$ and $c_{K,H}$ are the associated cocycles, then the standard probability spaces $(\Lambda^H\times X,{m_{\Lambda}}^{\otimes H}\otimes \mu)$ and $(\Lambda^K\times X,{m_{\Lambda}}^{\otimes K}\otimes \mu)$ provide an orbit equivalence coupling between $\Lambda\wr H$ and $\Lambda\wr K$, via the identification
\begin{equation*}
    \begin{array}{ccl}
    \Lambda^{H}\times X &\longrightarrow &\Lambda^K\times X\\
    ((\ell_{h})_{h\in H},x)&\longmapsto & ((\ell_{c_{K,H}(k,x)})_{k\in K},x).
\end{array}
\end{equation*}
One can check directly that the latter is a measured isomorphism.

Additionally, we can describe the associated cocycle $c_{\Lambda\wr H,\Lambda\wr K}\colon (\Lambda\wr H)\times (\Lambda^H\times X)\longrightarrow \Lambda\wr K$ (and similarly $c_{\Lambda\wr K,\Lambda\wr H}\colon(\Lambda\wr K)\times(\Lambda^K\times X)\longrightarrow \Lambda\wr H$) on the generating set 
\begin{equation*}
    S_{\Lambda\wr H}\defeq\lbrace(\delta^{H}_{\lambda},1_{H}) : \lambda\in\Lambda\rbrace\cup\lbrace(\mathds{1},h) : h\in H\rbrace
\end{equation*}
of $\Lambda\wr H$, where $\delta^{H}_{\lambda}\colon H\rightarrow\Lambda$ is the map which maps $1_{H}$ to $\lambda$ and any $h\neq 1_{H}$ to $1_{\Lambda}$. For every $\lambda\in\Lambda$, we have
\begin{equation*}
    c_{\Lambda\wr H,\Lambda\wr K}\left((\delta^{H}_{\lambda},1_{H}),((\ell_{h})_{h\in H},x)\right)=(\delta^{K}_{\lambda},1_{K}),
\end{equation*}
and, for every $h\in H$,
\begin{equation*}
    c_{\Lambda\wr H,\Lambda\wr K}\left((\mathds{1},h),((\ell_h)_{h\in H},x)\right)=(\mathds{1},c_{H,K}(h,x)).
\end{equation*}

Inspired by this construction of orbit equivalence couplings between lamplighters, we can build such couplings between some examples of halo products. Before presenting the applications for lampjugglers and lampcloners, let us explain informally the two main steps of the strategy.

\begin{enumerate}
    \item\label{step1} First, given a free p.m.p. action $H\curvearrowright (X,\mu)$, we want to build a standard probability space $(X_{H},\mu_{H})$ and a free p.m.p. action $\halo H\curvearrowright (X_{H},\mu_{H})$;
    \item\label{step2} If now $(X,\mu)$ is an orbit equivalence coupling between $H$ and $K$, we want a natural identification between $(X_{H},\mu_{H})$ and $(X_{K},\mu_{K})$ in order to get an orbit equivalence coupling between $\halo H$ and $\halo K$.
\end{enumerate}

For step~\eqref{step1}, to build a free p.m.p. action $\halo H\curvearrowright (X_{H},\mu_{H})$ from a free p.m.p. action $H\curvearrowright (X,\mu)$, we first need to find a free p.m.p. action $\triangle$ of $L(H)$ on some probability space $Z_{H}$. Informally, $H$ must “appear” in this set $Z_{H}$. Let us be more explicit on this point:
\begin{itemize}
\item We want the action of $H$ on itself to induce a free p.m.p. $L(H)-$action on $Z_{H}$. For instance, we will choose $Z_{H}=\Lambda^H$ for $\halo H=\Lambda\wr H$; we will choose $Z_{H}=[0,1]^H$ for $\halo H=\shuf{H}$; $Z_{H}=[0,1]^{V_{H}}$ for $\halo H=\cloner{H}$ where the action of $H$ on $[0,1]^{V_H}$ is induced by the action of $H$ on $V_{H}$, which is itself induced by the actions of $H$ on itself by left-multiplication. Combined with an action $\star$ of $H$ on $Z_{H}$, we thus have a $H-$ and a $L(H)-$action on $X_{H}\defeq Z_{H}\times X$, which motivates the following formula to deduce an action $\square$ of $\halo H$ on $X_{H}$:
    \begin{equation*}
        (\sigma,h)\square (z,x)=(\sigma\triangle (h\star z),h\cdot x)
    \end{equation*}
    for every $(\sigma,h)\in\halo H$ and every $(z,x)\in X_{H}$. For this formula to provide a well-defined $\halo H-$action, the $H-$ and $L(H)-$actions on $Z_{H}$ must be compatible in a certain way, namely they must satisfy
    \begin{equation}\label{eq:compatible}
        h\star (\sigma\triangle z)=(h\cdot_{\alpha}\sigma)\triangle (h\star z)
    \end{equation}
    for every $h\in H$, $\sigma\in L(H)$ and $z\in Z_{H}$, where $\cdot_{\alpha}$ denotes the action $H\curvearrowright L(H)$ used to define the halo product $\halo H$. Indeed, if this relation holds, then we have 
    \begin{align*}
        (\sigma_{1}, h_{1})\square \big( (\sigma_{2},h_{2})\square(z,x)\big) &= (\sigma_{1}, h_{1})\square (\sigma_{2}\triangle (h_{2}\star z),h_{2}\cdot x) \\
        &=\big(\sigma_{1}\triangle (h_{1}\star(\sigma_{2}\triangle (h_{2}\star z))), h_{1}\cdot (h_{2}\cdot x)\big)\\
        &=\big(\sigma_{1}\triangle [(h_{1}\cdot_{\alpha}\sigma_{2})\triangle(h_{1}\star(h_{2}\star z))], h_{1}h_{2}\cdot x\big)\\
        &=\big((\sigma_{1}\circ h_{1}\cdot_{\alpha}\sigma_{2})\triangle (h_{1}h_{2}\star z), h_{1}h_{2}\cdot x\big) \\
        &=[(\sigma_{1},h_{1})(\sigma_{2},h_{2})]\square(z,x)
    \end{align*}
    for all $(\sigma_{1},h_{1}),(\sigma_{2},h_{2})\in \halo H$, $(z,x)\in X_{H}$, using (\ref{eq:compatible}) in the fourth equality.
\item The second reason for $H$ to “appear” in $Z_{H}$ is for step~\eqref{step2}. To define a measured spaces isomorphism between $Z_{H}\times X$ and $Z_{K}\times X$, we use the fact that for every $x\in X$, the cocycle $c_{K,H}(\cdot,x)\colon K\rightarrow H$ is a bijection and, in some sense, it allows us to go from $Z_{H}$ to $Z_{K}$, with $x$ being the second coordinate of an element of $Z_{H}\times X$. For instance, for the lamplighters $\Lambda\wr H$ and $\Lambda\wr K$, we define
    \begin{equation*}
        \begin{array}{ccl}
    \Lambda^H\times X &\longrightarrow &\Lambda^K\times X\\
    \left((\ell_{h})_{h\in H},x\right)&\longmapsto & \left((\ell_{c_{K,H}(k,x)})_{k\in K},x\right)
    \end{array}.
    \end{equation*}
    For lampshufflers, where $Z_{H}=[0,1]^H$, we will rather use the bijection $c_{K,H}((\cdot)^{-1},x)^{-1}\colon K\rightarrow H$ and define
    \begin{equation*}
    \begin{array}{ccl}
    [0,1]^H\times X &\longrightarrow &[0,1]^K\times X\\
    \left((\varepsilon_h)_{h\in H},x\right)&\longmapsto & \left((\varepsilon_{c_{K,H}(k^{-1},x)^{-1}})_{k\in K},x\right).
    \end{array}
    \end{equation*}
    More subtle will be the case of lampcloners, where $Z_{H}=[0,1]^{V_{H}}$. In this case, we will set
    \begin{equation*}
    \begin{array}{ccl}
    [0,1]^{V_H}\times X &\longrightarrow &[0,1]^{V_K}\times X\\
    \left((\varepsilon_v)_{v\in V_H},x\right)&\longmapsto & \left((\varepsilon_{\varphi_{K,H}(w,x)})_{w\in V_K},x\right)
    \end{array}
    \end{equation*}
where $\varphi_{K,H}\colon V_{K}\rightarrow V_{H}$ is a bijection mapping $\sum_{k\in K}{\mu_{k}e_{k}}\in V_{K}$ to $\sum_{h\in H}{\mu_{c_{H,K}(h^{-1},x)^{-1}}e_{h}}$ (again we use the bijection $H\rightarrow K$ provided by the cocycle $c_{H,K}$).
\end{itemize}

Finally, steps~\eqref{step1} and~\eqref{step2} being achieved, it will remain to check, case by case, that we indeed get an orbit equivalence coupling between $\halo H$ and $\halo K$, and to quantify the cocycles. This strategy will then provide similar statements as the one of Theorem~\ref{thm:dklmtWreath} for lampjugglers and lampcloners.

\subsection{Proof of Theorem~\ref{thm:stabilityofcouplings+quantificationLampjugglersINTRO}}

Similarly to Theorem~\ref{thm:dklmtWreath}, let us prove the following.

\begin{theorem}\label{thm:stabilityofcouplings+quantificationLampjugglers}
Let $H$ and $K$ be groups and let $r\ge 1$ be an integer. If $H$ and $K$ are orbit equivalent, then $\juggler{r}{H}$ and $\juggler{r}{K}$ are orbit equivalent. Furthermore, given non-decreasing maps $\varphi,\psi\colon\R_{+}\rightarrow\R_{+}$, if $H$ and $K$ are finitely generated and if there exists a $(\varphi,\psi)-$integrable orbit equivalence coupling from $H$ to $K$, then the same holds from $\juggler{r}{H}$ to $\juggler{r}{K}$.
\end{theorem}

For the proof of the second part of this theorem, we will need the following intermediate observation on word lengths of elements of $\juggler{r}{H}$. Recall from Section~\ref{sec:defHalo} that 
\begin{equation*}
    S_{\juggler{r}{H}} \defeq \left\lbrace(\tau_{(1_{H},i),(s,j)},1_{H}) : s\in S_{H}, 1\le i,j\le r\right \rbrace\cup\left\lbrace(\text{id}_{H\times\lbrace1,\dots,r\rbrace},h) : h\in S_H\right\rbrace
\end{equation*}
is a finite generating set for $\juggler{r}{H}$, constructed from a finite generating set $S_{H}$ of $H$.

\begin{lemma}\label{lem:wordlengthinlampshuffler}
Let $H$ be a finitely generated group and let $S_{H}$ be a finite generating set. For every $g\in H$ and all integers $1\le i,j\le r$, we have
\begin{equation*}
    |(\tau_{(1_{H},i),(g,j)},1_{H})|_{S_{\juggler{r}{H}}}\le 4|g|_{S_{H}}.
\end{equation*}
\end{lemma}

\begin{proof}
Let us write $g=h_{1}\ldots h_{n}$ with $n=|g|_{S_H}$ and $h_{1},\ldots,h_{n}\in S_{H}\cup (S_{H})^{-1}$. Let us also set $\tau^{(\ell)}\defeq\tau_{(1_H,i),(h_{\ell}\ldots h_{n},j)}$ and $\sigma^{(\ell)}\defeq\tau_{(1_H,i),(h_{\ell}^{-1},i)}$, for any $1\le \ell\le n$. Then we have
\begin{equation*}
h^{-1}_{\ell}\cdot\tau^{(\ell)}=\tau_{(h_{\ell}^{-1},i),(h_{\ell+1}\ldots h_{n},j)}
\end{equation*}
and
\begin{equation*}
    \sigma^{(\ell)}(h^{-1}_{i}\cdot\tau^{(\ell)})\sigma^{(\ell)}=\tau^{(\ell+1)},
\end{equation*}
which implies 
\begin{equation*}
    (\sigma^{(\ell)},1_{H})(\mathrm{id}_{H\times\lbrace 1,\ldots, r\rbrace},h^{-1}_{\ell})(\tau^{(\ell)},1_{H})(\mathrm{id}_{H\times\lbrace 1,\ldots, r\rbrace},h_{\ell})(\sigma^{(\ell)},1_H)=(\tau^{(\ell+1)},1_{H}).
\end{equation*}
Thus it follows that $|(\tau^{(\ell)},1_{H})|_{S_{\juggler{r}{H}}} \le 4+|(\tau^{(\ell+1)}, 1_{H})|_{S_{\juggler{r}{H}}}$, so by induction we get 
\begin{equation*}
|(\tau_{(1_{H},i),(g,j)},1_{H})|_{S_{\juggler{r}{H}}}=|(\tau^{(1)},1_{H})|_{S_{\juggler{r}{H}}}\le 4\cdot n=4\cdot |g|_{S_{H}}
\end{equation*}
and we are done.
\end{proof}

\begin{proof}[Proof of Theorem~\ref{thm:stabilityofcouplings+quantificationLampjugglers}]
First, we must define free and p.m.p. $\juggler{r}{H}-$ and $\juggler{r}{K}-$actions, provided such $H-$ and $K-$actions. Following the techniques we explained just above, we first define free p.m.p. $H-$ and $\fsym{H\times\lbrace1,\dots,r\rbrace}-$actions on $\left([0,1]^{H\times\lbrace 1,\ldots,r\rbrace},\mathrm{Leb}^{\otimes (H\times\lbrace 1,\ldots,r\rbrace)}\right)$ by
\begin{equation*}
    h\cdot (\varepsilon_{g,i})_{(g,i)\in H\times\lbrace 1,\dots,r\rbrace}\defeq (\varepsilon_{h^{-1}g,i})_{(g,i)\in H\times\lbrace1,\dots,r\rbrace}
\end{equation*}
and
\begin{equation*}
    \sigma\cdot(\varepsilon_{g,i})_{(g,i)\in H\times\lbrace 1,\dots,r\rbrace}\defeq (\varepsilon_{\sigma^{-1}(g,i)})_{(g,i)\in H\times\lbrace1,\dots,r\rbrace}
\end{equation*}
for every $h\in H$, $\sigma\in\fsym{H\times\lbrace1,\dots,r\rbrace}$ and $(\varepsilon_{g,i})_{(g,i)\in H\times\lbrace1,\dots,r\rbrace}\in [0,1]^{H\times\lbrace1,\dots,r\rbrace}$. The compatibility condition~\eqref{eq:compatible} is satisfied, so the formula 
\begin{equation*}
(\sigma,h)\cdot\left((\varepsilon_{g,i})_{(g,i)\in H\times\lbrace 1,\ldots,r\rbrace},x\right)\defeq \left((\varepsilon_{h^{-1}\sigma^{-1}(g,i)})_{(g,i)\in H\times\lbrace1,\dots,r\rbrace},h\cdot x\right),
\end{equation*}
where $h^{-1}\sigma^{-1}(g,i)=(h^{-1}g',i')$ if $\sigma^{-1}(g,i)=(g',i')$ (action of $H$ on $H\times\lbrace1,\dots,r\rbrace$), defines a $\juggler{r}{H}-$action on $X_{H}\defeq [0,1]^H\times X$, endowed with the probability measure $\mu_{H}\defeq\mathrm{Leb}^{\otimes (H\times\lbrace 1,\ldots,r\rbrace)}\otimes\mu$. It is not difficult to prove that this action is p.m.p. and free. We similarly define a free p.m.p. $\juggler{r}{K}-$action on 
\begin{equation*}
    (X_{K},\mu_{K})\defeq ([0,1]^{K\times\lbrace 1,\dots,r\rbrace}\times X,\mathrm{Leb}^{\otimes (K\times\lbrace 1,\ldots,r\rbrace)}\otimes\mu).
\end{equation*}
If now $(X,\mu)$ is an orbit equivalence coupling between $H$ and $K$, let us prove that $(X_{H},\mu_{H})$ and $(X_K,\mu_K)$ are orbit equivalent couplings between $\juggler{r}{H}$ and $\juggler{r}{K}$, via the identification
\begin{equation*}
    \theta\colon\begin{array}{ccl}
    [0,1]^{H\times\{1,\ldots,r\}}\times X &\longrightarrow &[0,1]^{K\times\lbrace1,\dots,r\rbrace}\times X\\
    \left((\varepsilon_{h,i})_{(h,i)\in H\times\lbrace 1,\dots,r\rbrace},x\right)&\longmapsto & \left((\varepsilon_{c_{K,H}(k^{-1},x)^{-1},i})_{(k,i)\in K\times\lbrace1,\dots,r\rbrace},x\right)
\end{array}.
\end{equation*}
We first have to check that $\theta$ is a measured isomorphism between these probability spaces, this can be done by fixing the second coordinate $x$ so that we fix the bijection $K\to H$ given by the cocycle $c_{K,H}(\cdot,x)$. Secondly, let us prove that we have an orbit equivalence. Given generating sets $S_{H}$ and $S_{K}$ of $H$ and $K$, we get generating sets $S_{\juggler{r}{H}}$ and $S_{\juggler{r}{K}}$ of $\juggler{r}{H}$ and $\juggler{r}{K}$ respectively, as described before Lemma~\ref{lem:wordlengthinlampshuffler} (note that these sets are not necessarily finite for now). We prove that we have an orbit equivalence using Lemma~\ref{lem:EqOrbitsGenerating}. Let $\left((\varepsilon_{h,i})_{(h,i)\in H\times\lbrace 1,\dots,r\rbrace},x\right)\in X_{H}$. Given $s\in S_{H}$, we have
\begin{align*}
    &\theta\left((\mathrm{id},s)\cdot \left((\varepsilon_{h,i})_{(h,i)\in H\times\lbrace 1,\dots,r\rbrace},x\right)\right)\\
    &=\theta\left((\varepsilon_{s^{-1}h,i})_{(h,i)\in H\times\lbrace 1,\dots,r\rbrace},s\cdot x\right)\\
    &=\left((\varepsilon_{s^{-1}c_{K,H}(k^{-1},s\cdot x)^{-1},i})_{(k,i)\in K\times\lbrace1,\dots,r\rbrace},s\cdot x\right)\\
    &=\left((\varepsilon_{(c_{K,H}(k^{-1},s\cdot x)s)^{-1},i})_{(k,i)\in K\times\lbrace1,\dots,r\rbrace},s\cdot x\right)\\
    &=\left((\varepsilon_{c_{K,H}(k^{-1}c_{H,K}(s,x),x)^{-1},i})_{(k,i)\in K\times\lbrace1,\dots,r\rbrace},c_{H,K}(s,x)\cdot x\right)\\
    &=\left((\varepsilon_{c_{K,H}((c_{H,K}(s,x)^{-1}k)^{-1},x)^{-1},i})_{(k,i)\in K\times\lbrace1,\ldots,r\rbrace},c_{H,K}(s,x)\cdot x\right)\\
    &=\left(\mathrm{id},c_{H,K}(s,x)\right)\cdot \left((\varepsilon_{c_{K,H}(k^{-1},x)^{-1},i})_{(k,i)\in K\times\lbrace1,\dots,r\rbrace},x\right)\\
    &=\left(\mathrm{id},c_{H,K}(s,x)\right)\cdot\theta \left((\varepsilon_{h,i})_{(h,i)\in H\times\lbrace1,\dots,r\rbrace},x\right)
\end{align*}
where the equality $c_{K,H}(k^{-1},s\cdot x)s=c_{K,H}(k^{-1}c_{H,K}(s,x),x)$ holds since we have 
\begin{equation*}
    c_{K,H}(k^{-1},s\cdot x)s\cdot x=k^{-1}\cdot (s\cdot x)=k^{-1} c_{H,K}(s,x)\cdot x=c_{K,H}(k^{-1} c_{H,K}(s,x),x)\cdot x
\end{equation*}
for every $x\in X$. For the sequel, we set the following notations: $c_{K,H}((k,i),x)\defeq (c_{K,H}(k,x),i)$ and $(k,i)^{-1}\defeq (k^{-1},i)$. Then, for every $s\in S_{H}$ and $j,\ell\in\lbrace1,\ldots,r\rbrace$, we have
\begin{align*}
    &\theta\left((\tau_{(1_{H},j),(s,\ell)},1_{H})\cdot \left((\varepsilon_{h,i})_{(h,i)\in H\times\lbrace1,\dots,r\rbrace},x\right)\right)\\
    &=\theta\left( (\varepsilon_{\tau_{(1_H,j),(s,\ell)}(h,i)})_{(h,i)\in H\times\lbrace1,\dots,r\rbrace},x\right)\\
    &=\left((\varepsilon_{\tau_{(1_{H},j),(s,\ell)}(c_{K,H}(k^{-1},x)^{-1},i)})_{(k,i)\in K\times\lbrace1,\dots,r\rbrace},x\right)\\
    &=\left((\varepsilon_{c_{K,H}([\tau_{(1_K,j),(c_{H,K}(s^{-1},x)^{-1},\ell)}(k,i)]^{-1},x)^{-1}})_{(k,i)\in K\times\lbrace1,\dots,r\rbrace},x\right)\\
    &=\left(\tau_{(1_K,j),(c_{H,K}(s^{-1},x)^{-1},\ell)},1_{K}\right)\cdot \left((\varepsilon_{c_{K,H}(k^{-1},x)^{-1},i})_{(k,i)\in K\times\lbrace1,\dots,r\rbrace},x\right)\\
    &=\left(\tau_{(1_K,j),(c_{H,K}(s^{-1},x)^{-1},\ell)},1_{K}\right)\cdot \theta\left ((\varepsilon_{h,i})_{(h,i)\in H\times\lbrace1,\dots,r\rbrace},x)\right)
\end{align*}
where we prove the equality 
\begin{equation*}
    \tau_{(1_{H},j),(s,\ell)}\left(c_{K,H}(k^{-1},x)^{-1},i\right)=c_{K,H}\left([\tau_{(1_{K},j),(c_{H,K}(s^{-1},x)^{-1},\ell)}(k,i)]^{-1},x\right)^{-1}
\end{equation*}
depending on whether or not $(k,i)$ lies in $\left\lbrace (1_{K},j),(c_{H,K}(s^{-1},x)^{-1},\ell)\right\rbrace$.

\noindent Thus, we have proved
\begin{equation*}
    \theta(S_{\juggler{r}{H}}\cdot z)\subset \juggler{r}{K}\cdot\theta(z)
\end{equation*}
for every $z\in X_{H}$ and we similarly get $S_{\juggler{r}{K}}\cdot \theta(z)\subset \theta(\juggler{r}{H}\cdot z)$, so Lemma~\ref{lem:EqOrbitsGenerating} implies that we have built an orbit equivalence.

\noindent Let us finally assume that $H$ and $K$ are finitely generated (i.e. $S_{H}$ and $S_{K}$ are finite) and that $(X,\mu)$ is a $(\varphi,\psi)-$integrable orbit equivalence from $H$ to $K$. The last computations in fact provide the cocycles on the generators: for every $s\in S_{H}$, one has
\begin{equation*}
    c_{\juggler{r}{H},\juggler{r}{K}}\left((\mathrm{id},s),((\varepsilon_h)_{h\in H},x)\right)=\left(\mathrm{id},c_{H,K}(s,x)\right)
\end{equation*}
and
\begin{equation*}
    c_{\juggler{r}{H},\juggler{r}{K}}\left(\left(\tau_{(1_{H},j),(s,\ell)},1_{H}\right),((\varepsilon_h)_{h\in H},x)\right)=\left(\tau_{(1_{K},j),(c_{H,K}(s^{-1},x)^{-1},\ell)},1_{K}\right)
\end{equation*}
and similarly for the cocycle $c_{\juggler{r}{K},\juggler{r}{H}}$. Hence, using Lemma~\ref{lem:wordlengthinlampshuffler}, we easily get
\begin{equation*}
    \left|c_{\juggler{r}{H},\juggler{r}{K}}\left((\mathrm{id},s),((\varepsilon_h)_{h\in H},x)\right)\right|_{S_{\juggler{r}{K}}}=|c_{H,K}(s,x)|_{S_{K}}
\end{equation*}
for every $s\in S_{H}$, and 
\begin{equation*}
    \left|c_{\juggler{r}{H},\juggler{r}{K}}\left(\left(\tau_{(1_{H},j),(s,\ell)},1_{H}\right),((\varepsilon_h)_{h\in H},x)\right)\right|_{S_{\juggler{r}{K}}} \le  4|c_{H,K}(s^{-1},x)|_{S_{K}},
\end{equation*}
and similarly for $c_{\juggler{r}{K},\juggler{r}{H}}$. Thus the $\varphi-$ and $\psi-$integrabilities of $c_{H,K}$ and $c_{K,H}$ directly imply the same properties for $c_{\juggler{r}{H},\juggler{r}{K}}$ and $c_{\juggler{r}{K},\juggler{r}{H}}$ on generators, and thus on all elements of the groups by Remark~\ref{rem:checkongenerators}, so we get a $(\varphi,\psi)-$integrable orbit equivalence coupling from $\juggler{r}{H}$ to $\juggler{r}{K}$.
\end{proof}

\paragraph{Lampcloners.}

Similarly to Theorems~\ref{thm:dklmtWreath} and~\ref{thm:stabilityofcouplings+quantificationLampjugglers}, let us prove the following.

\begin{theorem}\label{thm:stabilityofcouplings+quantificationLampcloners}
Let $H$ and $K$ be groups and let $\field$ be a finite field. If $H$ and $K$ are orbit equivalent, then $\cloner{H}$ and $\cloner{K}$ are orbit equivalent. Furthermore, given non-decreasing maps $\varphi,\psi\colon\R_{+}\rightarrow\R_{+}$, if $H$ and $K$ are finitely generated and if there exists a $(\varphi,\psi)-$integrable orbit equivalence coupling from $H$ to $K$, then the same holds from $\cloner{H}$ to $\cloner{K}$.
\end{theorem}

For the proof, we will need the following intermediate observation on word lengths of elements of $\cloner{H}$. Recall from Section~\ref{sec:defHalo} that 
\begin{equation*}
    S_{\cloner{H}}\defeq\lbrace (\delta_{1_{H}}(\lambda), 1_{H}) : \lambda \in \field\setminus\lbrace 0\rbrace\rbrace\cup\lbrace (\tau_{1_{H},s}(\lambda),1_{H}) : s\in S_{H}\rbrace \cup\lbrace (\text{id}, s) : s\in S_{H}\rbrace
\end{equation*}
is a finite generating set for $\cloner{H}$, constructed from a finite generating set $S_{H}$ of $H$.

\begin{lemma}\label{lem:wordlengthinlampcloner}
Let $H$ be a finitely generated group. For every $g\in H$ and every $\lambda\in\field$, we have
\begin{equation*}
    \left|(\tau_{1_H,g}(\lambda),1_{H})\right|_{S_{\cloner{H}}}\le 14|g|_{S_{H}}.
\end{equation*}
\end{lemma}

\begin{proof}
Let us write $g=h_{1}\ldots h_{n}$ with $n=|g|_{S_H}$ and $h_{1},\ldots,h_{n}\in S_{H}\cup (S_{H})^{-1}$. Let us also set, for any $1\le i\le n$, $\tau^{(i)}\defeq\tau_{1_H,h_{i}\ldots h_{n}}(\lambda)$ and $\sigma^{(i)}$ the linear automorphism which acts on the canonical basis of $V_{H}$ in the following way: it swaps $e_{1_{H}}$ and $e_{h_{i}^{-1}}$, and it fixes $e_{h}$ for any other $h\in H\setminus\lbrace1_{H},h^{-1}_{i}\rbrace$. Then we have
\begin{equation*}
h^{-1}_{i}\cdot\tau^{(i)}=\tau_{h_{i}^{-1},h_{i+1}\dots h_{n}}(\lambda)
\end{equation*}
and
\begin{equation*}
    \sigma^{(i)}(h^{-1}_{i}\cdot\tau^{(i)})\sigma^{(i)}=\tau^{(i-1)}
\end{equation*}
which implies
\begin{equation*}
    (\sigma^{(i)},1_{H})(\mathrm{id}_{H},h^{-1}_{i})(\tau^{(i)},1_{H})(\mathrm{id}_{H},h_{i})(\sigma^{(i)},1_{H})=(\tau^{(i-1)},1_{H}).
\end{equation*}
Moreover, $h_{i}\cdot \sigma^{(i)}$ can be written as a product of four dilatations and transvections supported in $\lbrace1_{H},h_{i}\rbrace$, so we have $|(\sigma^{(i)},1_{H})|_{\cloner{H}}\le 6$. Thus it follows that $|(\tau^{(i)},1_{H})|_{S_{\cloner{H}}} \le 14+|(\tau^{(i-1)}, 1_{H})|_{S_{\cloner{H}}}$, and by induction we get 
\begin{equation*}
|(\tau_{1_{H},g}(\lambda), 1_{H})|_{S_{\cloner{H}}}=|(\tau^{(n)},1_{H})|_{S_{\cloner{H}}}\le 14n=14|g|_{S_{H}}
\end{equation*}
as claimed. The proof is complete.
\end{proof}

\begin{proof}[Proof of Theorem~\ref{thm:stabilityofcouplings+quantificationLampcloners}]
Given generating sets $S_{H}$ and $S_{K}$ (not necessarily finite) of $H$ and $K$, we get generating subsets $S_{\cloner{H}}$ and $S_{\cloner{K}}$ of $\cloner{H}$ and $\cloner{K}$, as described before Lemma~\ref{lem:wordlengthinlampcloner}.

\noindent We must first define free and p.m.p. $\cloner{H}-$ and $\cloner{K}-$actions provided such $H-$ and $K-$actions on $(X,\mu)$. Let us set 
\begin{equation*}
    (X_{H},\mu_{H})\defeq \left([0,1]^{V_{H}}\times X,\mathrm{Leb}^{\otimes V_{H}}\otimes\mu\right). 
\end{equation*}
Then $H$ and $\mathrm{FGL}(H)$ act on $[0,1]^{V_{H}}$ in the following way:
\begin{equation*}
        \sigma\cdot(\varepsilon_v)_{v\in V_{H}}\defeq (\varepsilon_{\sigma^{-1}(v)})_{v\in V_{H}}
\end{equation*}
and
\begin{equation*}
        h\cdot (\varepsilon_v)_{v\in V_H}\defeq (\varepsilon_{h\cdot v})_{v\in V_H}
\end{equation*}
for every $\sigma\in\mathrm{FGL}(H)$, $h\in H$ and $(\varepsilon_v)_{v\in V_H}\in [0,1]^{V_H}$, using the action of $H$ on $V_{H}$ given by
\begin{equation*}
    h\cdot\left(\sum_{g\in H}{\mu_{g} e_{g}}\right)\defeq\sum_{g\in H}{\mu_{g}e_{h^{-1}g}}.
\end{equation*}
These $H-$ and $\mathrm{FGL}(H)-$actions are compatible in the sense of~\eqref{eq:compatible}, and we get a free p.m.p. $\cloner{H}-$action on $(X_{H},\mu_{H})$ defined by
\begin{equation*}
    (\sigma,h)\cdot \left((\varepsilon_v)_{v\in V_{H}},x\right)\defeq \left((\varepsilon_{h\cdot\sigma^{-1}(v)})_{v\in V_{H}},h\cdot x\right).
\end{equation*}
We similarly define a free p.m.p. $\cloner{K}-$action on 
\begin{equation*}
    (X_{K},\mu_{K})\defeq \left([0,1]^{V_{K}}\times X,\mathrm{Leb}^{\otimes V_{K}}\otimes\mu\right). 
\end{equation*}
Now the map 
\begin{equation*}
    \theta\colon\begin{array}{ccl}
    [0,1]^{V_{H}}\times X &\longrightarrow &[0,1]^{V_K}\times X\\
    \left((\varepsilon_{v})_{v\in V_H},x\right) &\longmapsto & \left((\varepsilon_{\varphi_{K,H}(w,x)})_{w\in V_{K}},x\right)
\end{array}
\end{equation*}
is a measured isomorphism, where $\varphi_{K,H}\colon V_{K}\rightarrow V_{H}$ is a bijection mapping $\sum_{k\in K}{\mu_{k}e_{k}}\in V_{K}$ to $\sum_{h\in H}{\mu_{c_{H,K}(h^{-1},x)^{-1}}e_{h}}$. Moreover, for every $s\in S_{H}$, one has
\begin{equation*}
    \theta\left((\mathrm{id},s)\cdot \left((\varepsilon_{v})_{v\in V_{H}},x\right)\right)=(\mathrm{id},c_{H,K}(s,x))\cdot\theta \left((\varepsilon_{v})_{v\in V_{H}},x\right)
\end{equation*}
as well as
\begin{equation*}
    \theta\left((\tau_{1_{H},s}(\lambda),1_{H})\cdot ((\varepsilon_{v})_{v\in V_H},x)\right)=(\tau_{1_{K},c_{H,K}(s^{-1},x)^{-1}},1_{K})\cdot \theta\left((\varepsilon_{v})_{v\in V_{H}},x\right)
\end{equation*}
and
\begin{equation*}
    \theta\left((\delta_{1_{H}}(\lambda),1_{H})\cdot ((\varepsilon_{v})_{v\in V_{H}},x)\right)=(\delta_{1_{K}}(\lambda),1_{K})\cdot \theta\left((\varepsilon_{v})_{v\in V_{H}},x\right).
\end{equation*}

\noindent We thus have proved $\theta(S_{\cloner{H}}\cdot z)\subset \cloner{K}\cdot\theta(z)$ for every $z\in X_H$ and we similarly get 
\begin{equation*}
    S_{\cloner{K}}\cdot \theta(z)\subset \theta\left(\cloner{H}\cdot z\right)
\end{equation*}
so Lemma~\ref{lem:EqOrbitsGenerating} implies that we have built an orbit equivalence. Finally, in the case of finitely generated groups $H$ and $K$, since we have identified the cocycles $c_{\cloner{H},\cloner{K}}$, $c_{\cloner{K},\cloner{H}}$ on the generators, we conclude, similarly to the proof of Theorem~\ref{thm:stabilityofcouplings+quantificationLampjugglers} (using this time Lemma~\ref{lem:wordlengthinlampcloner}), that a $(\varphi,\psi)-$integrable orbit equivalence coupling from $H$ to $K$ provides a $(\varphi,\psi)-$integrable orbit equivalence coupling from $\cloner{H}$ to $\cloner{K}$.
\end{proof}

\subsection{Proofs of Theorems~\ref{thm:optimalitypolynomialgrowthJugglersINTRO} and~\ref{thm:optimalitypolynomialgrowthClonersINTRO}}

We now establish that the couplings constructed in the previous section are optimal in many cases. As a first application:

\begin{theorem}\label{thm:optimalitypolynomialgrowthJugglers}
Let $k,d,r\ge 1$ be integers such that $k>d$. Let $n\ge 0$ be an integer. Let $p>0$. If $H$ and $K$ are polynomial growth groups of degrees $k$ and $d$ respectively, then $\jugglern{n}{r}{H}$ and $\jugglern{n}{r}{K}$ are $\ld^{p}$ orbit equivalent if and only if $p<\frac{d}{k}$.
\end{theorem}

\begin{proof}
Suppose first that $p<\frac{d}{k}$. From~\cite[Theorem~1.6]{DLIT25}, we know that $H$ and $K$ are $\ld^p$ orbit equivalent, so the conclusion follows from Theorem~\ref{thm:stabilityofcouplings+quantificationLampjugglers}.

\noindent Conversely, assume that $\jugglern{n}{r}{H}$ and $\jugglern{n}{r}{K}$ are $\ld^p$ orbit equivalent. 
Then Theorem~\ref{thm:ObstructionDKLMT}\textit{(i)} implies that \begin{equation*}
    \left(\prof{\jugglern{n}{r}{K}}(x)\right)^{p} \preccurlyeq \prof{\jugglern{n}{r}{H}}(x)
\end{equation*}
and thus, from Proposition~\ref{prop:profileofshufnofpolynomialgrowthgroupsINTRO}, we get 
\begin{equation*}
    \left(\frac{\ln^{\circ n}(x)}{\ln^{\circ (n+1)}(x)}\right)^{\frac{p}{d}} \preccurlyeq \left(\frac{\ln^{\circ n}(x)}{\ln^{\circ (n+1)}(x)}\right)^{\frac{1}{k}}.
\end{equation*}
This domination forces $\frac{p}{d}\le \frac{1}{k}$, i.e. $p\le \frac{d}{k}$. Lastly, from Theorem~\ref{thm:threshold} we deduce that there is no $\ld^{\frac{d}{k}}$ orbit equivalence coupling between $\jugglern{n}{r}{H}$ and $\jugglern{n}{r}{K}$.
\end{proof}

Other examples of quantitatively optimal orbit equivalence couplings have been found in~\cite{DKLMT22}, for instance between $\Z$ and $F\wr\Z$, where $F$ is a non-trivial finite group. Combined with our stability result, we get the following.

\begin{theorem}
Let $F$ be a non-trivial finite group. Let $n,r\ge 1$ be integers. Let $p>0$. Then there is an $(\exp,\ln^p)-$integrable orbit equivalence coupling from $\jugglern{n}{r}{\Z}$ to $\jugglern{n}{r}{F\wr\Z}$ if and only if $p<1$.
\end{theorem}

\begin{proof}
Given $\varepsilon>0$, let $\psi_{\varepsilon}\colon \R_{+}\rightarrow \R_{+}$ be the map defined by
\begin{equation*}
        \psi_{\varepsilon}(x)=\frac{\ln(x)}{\ln(\ln(x))^{1+\varepsilon}}.
\end{equation*}
From~\cite[Proposition~6.20]{DKLMT22}, we know that, for any $\varepsilon>0$, there exists an $(\exp,\psi_{\varepsilon})-$integrable orbit equivalence coupling from $\Z$ to $F\wr\Z$. In particular, this coupling is $(\exp,\ln^p)-$integrable for every $p<1$. Hence Theorem~\ref{thm:stabilityofcouplings+quantificationLampjugglers} provides an $(\exp,\ln^p)-$integrable orbit equivalence coupling from $\jugglern{n}{r}{\Z}$ to $\jugglern{n}{r}{F\wr\Z}$ for every $p<1$.

\noindent Now let us assume the existence of an $(\exp,\ln^p)-$integrable orbit equivalence coupling from $\jugglern{n}{r}{\Z}$ to $\jugglern{n}{r}{F\wr\Z}$ and let us prove that $p<1$. By Theorem~\ref{thm:ObstructionDKLMT}\textit{(i)}, we get
\begin{equation*}
        \left(\ln(\prof{\jugglern{n}{r}{\Z}}(x))\right)^{p}\preccurlyeq\prof{\jugglern{n}{r}{F\wr\Z}}(x)
\end{equation*}
namely $\left (\ln^{\circ (n+1)}(x)\right )^p\preccurlyeq\ln^{\circ (n+1)}(x)$ using Propositions~\ref{prop:profileofshufnofpolynomialgrowthgroupsINTRO} and~\ref{prop:profilesofshufnINTRO}, which forces $p\le 1$. Lastly, an orbit equivalence in the case $p=1$ is excluded by Theorem~\ref{thm:threshold}.
\end{proof}

For lampcloners, we get similar results, except that we lose the “if and only if” in the analogue of Theorem~\ref{thm:optimalitypolynomialgrowthJugglers} since we do not have precise estimates for the isoperimetric profile of iterated lampcloners over polynomial growth groups.

\begin{theorem}\label{thm:optimalitypolynomialgrowthCloners}
Let $\field$ be a finite field and let $k,d\ge 1$ be integers such that $k>d$. Then, for every integer $n\ge 0$, for every $p<\frac{d}{k}$, $\clonern{n}{\Z^k}$ and $\clonern{n}{\Z^d}$ are $\ld^{p}$ orbit equivalent. Conversely, if $\cloner{\Z^d}$ and $\cloner{\Z^k}$ are $\ld^{p}$ orbit equivalent, then $p\le\frac{2d}{k}$.
\end{theorem}

\begin{proof}
From~\cite[Theorem~6.12]{DKLMT22}, we know that $\Z^k$ and $\Z^d$ are $\ld^p$ orbit equivalent when $p<\frac{d}{k}$, so Theorem~\ref{thm:stabilityofcouplings+quantificationLampcloners} provides such a coupling between $\clonern{n}{\Z^k}$ and $\clonern{n}{\Z^d}$. Conversely, if $\cloner{\Z^k}$ and $\cloner{\Z^d}$ are $\ld^{p}$ orbit equivalent, then Theorem~\ref{thm:ObstructionDKLMT}\textit{(i)} implies
\begin{equation}\label{eq:3.2}
    \left(\prof{\cloner{\Z^d}}(x)\right)^p\preccurlyeq\prof{\cloner{\Z^k}}(x).
\end{equation}
We know from Corollary~\ref{cor:encadrementduprofildecloneINTRO} that
\begin{equation*}
    \prof{\cloner{\Z^k}}(x)\preccurlyeq \ln(x)^{\frac{1}{k}}
\end{equation*}
and
\begin{equation*}
    \prof{\cloner{\Z^d}}(x)\succcurlyeq \ln(x)^{\frac{1}{2d}}
\end{equation*}
and thus (\ref{eq:3.2}) forces $\frac{p}{2d}\leq\frac{1}{k}$, namely $p\leq\frac{2d}{k}$.
\end{proof}

\section{Construction of orbit equivalence couplings using F\o lner tiling sequences}\label{sec:folnerTiling}

In this section, given a halo product $\halo$ (for instance a lampshuffler or a lampjuggler), we construct orbit equivalence couplings between $\Z^d$ and $\halo \Z^k$. A powerful tool to construct such equivalences is provided by \textit{F\o lner tiling sequences}, introduced first in~\cite{DKLMT22} and that we present now. Combined with a composition result of couplings from~\cite{DKLMT22} and our stability results (Theorems~\ref{thm:stabilityofcouplings+quantificationLampjugglers} and~\ref{thm:stabilityofcouplings+quantificationLampcloners}), we also construct couplings between iterated halo products. These applications are presented in the subsequent subsections.

\subsection{Preliminaries}\label{sec:PreliminariesFolnerTiling}

Given an amenable group $ G$, a (right) F\o lner tiling sequence $(F_{n})_{n\ge 0}$ is a (right) F\o lner sequence of $G$ satisfying the \textit{tiling condition}: for every $n\ge 0$, there exists a finite subset $\Sigma_{n}$ of $G$ such that
\begin{itemize}
    \item the translates $g F_{n}$, for $g\in\Sigma_{n}$, are pairwise disjoint;
    \item $F_{n+1}=\Sigma_{n} F_{n}$.
\end{itemize}
The set $F_{n}$ is a \textit{tile}, $\Sigma_{n}$ is a set of \textit{shifts}, and the sequence $(\Sigma_{n})_{n\ge 0}$ is the \textit{F\o lner tiling shift} associated to the F\o lner tiling sequence $(F_{n})_{n\ge 0}$. We can define analogously left F\o lner tiling sequences, but we will only work with right F\o lner tiling sequences in this paper. For convenience, we also assume $F_{0}=\lbrace1_{G}\rbrace$, so that any element of $F_{n+1}$ can be uniquely written as $x_{n}x_{n-1}\dots x_{0}$ with $x_{i}\in\Sigma_{i}$.

\begin{example}
Given integer $m\ge 2$ and $d\ge 1$, the sequence $(F_{n})_{n\ge 0}$ defined by
\begin{equation*}
        F_{n}\defeq \lbrace 0,1,\ldots,m^{n}-1\rbrace^{d}
\end{equation*}
is a F\o lner tiling sequence of $\Z^{d}$, with F\o lner tiling shifts $(\Sigma_{n})_{n\ge 0}$ given by
\begin{equation*}
        \Sigma_{n}\defeq \lbrace 0,m^n,2m^n,\dots,(m-1)m^n\rbrace^{d}.
\end{equation*}
\end{example}

A F\o lner tiling sequence gives rise to a probability measure-preserving $G-$action on the Cantor set
\begin{equation*}
    X_{F}\defeq\prod_{n\ge 0}{\Sigma_{n}}
\end{equation*}
endowed with the product measure of uniform distributions on each $\Sigma_{n}$. It is defined in the following way: the F\o lner condition implies that, for almost every $(x_{n})_{n\ge 0}$, for every $g\in G$, we have
\begin{equation}\label{eq:defactionfolnertiling}
    x_{N}x_{N-1}\dots x_{0} g\in F_{N+1}
\end{equation}
for a large enough integer $N$, so it can uniquely be written as $x'_{N}\dots x'_{0}$ with $x'_{i}\in\Sigma_{i}$, and we define
\begin{equation*}
g \cdot (x_{n})_{n\ge 0}\defeq (x'_{n})_{n\ge 0}
\end{equation*}
with $x'_{i}\defeq x_{i}$ for every $i\ge N+1$. This definition does not depend on the integer $N$ for which Equation~\eqref{eq:defactionfolnertiling} occurs.

The main interest of this action is that the equivalence relation it generates is (up to a null set) the cofinite equivalence relation: for almost every $\mathbf{x}=(x_{n})_{n\ge 0}\in X_F$, for every $g\in G$, we have $x_{n}=(g\mathbf{x})_{n}$ for large enough integers $n$. Moreover, if $(F'_{n})_{n\ge 0}$ is a F\o lner tiling sequence of another finitely generated group $ H$ satisfying $|F_{n}|=|F'_{n}|$ for every $n\ge 0$, then using a bijection $\Sigma_{n}\rightarrow\Sigma_{n}'$ for each $n\ge 0$ , we get a measure isomorphism between $X_{F}$ and $X_{F'}$. This is actually an orbit equivalence between the underlying actions since they generate the cofinite equivalence relation. We now present a criterion, proved in~\cite{DKLMT22}, to quantify this orbit equivalence.

\begin{definition}
Let $ G$ be a finitely generated group with a finite generating set $S_{G}$. Given sequences $R=(R_{n})_{n\ge 0}$ and $\varepsilon=(\varepsilon_{n})_{n\ge 0}$ of positive real numbers, we say that the sequence $(F_{n})_{n\ge 0}$ is an $(R,\varepsilon)-$F\o lner tiling sequence of $G$ if
\begin{itemize}
    \item $(F_{n})_{n\ge 0}$ is a F\o lner tiling sequence of $ G$;
    \item for every $n\ge 0$, $\diam{F_{n}}\le R_{n}$;
    \item for every $n\ge 0$ and every $s\in S_{G}$, $\frac{|F_{n}s\setminus F_{n}|}{|F_{n}|}\le\varepsilon_{n}$,
\end{itemize}
where $\diam{F}$ denotes the diameter of a finite subset $F$ of $G$, defined as
\begin{equation*}
    \diam{F}\defeq\sup_{f,g\in F}{|f^{-1}g|_{S_{G}}}.
\end{equation*}
\end{definition}

Here is then a sufficient condition to check on the tilings to ensure that the couplings constructed above have the desired level of integrability.

\begin{theorem}[{\cite[Proposition~6.9]{DKLMT22}}]\label{thm:TilingSufficientConditionQuantitative}
Let $ G$ and $ H$ be two finitely generated groups, and $\varphi,\psi\colon\R_{+}\rightarrow\R_{+}$ be non-decreasing maps. Assume that there exist sequences $R=(R_{n})_{n\ge 0}$, $R'=(R'_{n})_{n\ge 0}$, $\varepsilon=(\varepsilon_{n})_{n\ge 0}$ and $\varepsilon'=(\varepsilon'_{n})_{n\ge 0}$ of positive real numbers, an $(R,\varepsilon)-$F\o lner tiling sequence $(F_{n})_{n\in\N}$ of $G$ and an $(R',\varepsilon')-$F\o lner tiling sequence $(F_{n}')_{n\in\N}$ of $H$ such that $|F_{n}|=|F_{n}'|$ for any $n\in\N$, and assume that the series
\begin{equation*}
    \sum_{n\ge 0} \varphi(R_{n+1}')\varepsilon_{n}\; \text{and}\; \sum_{n\ge 0} \psi(R_{n+1})\varepsilon_{n}'
\end{equation*}
converge. Then the orbit equivalence coupling built above from $G$ to $H$ is $(\varphi,\psi)-$integrable.
\end{theorem}

Note that, when one wants to get quantitative orbit equivalences between two groups $ G$ and $ H$, F\o lner tiling sequences are particularly useful if one of the two groups is $\Z$. Indeed, if we find a F\o lner tiling sequence of $ G$, then a natural F\o lner tiling sequence $(F'_{n})_{n\ge 0}$ for $ H=\Z$, satisfying $|F_{n}|=|F_{n}'|$, is
\begin{equation*}
    F_{n}'\defeq \lbrace 0,1,\ldots, |F_{n}|-1\rbrace.
\end{equation*}
However the tiling condition and the condition on the cardinalities make the problem harder when $ H$ is bigger than $\Z$, even for $ H=\Z^2$. A clever choice of F\o lner tiling sequence of free abelian groups yields the following.

\begin{theorem}[{\cite[Theorem~6.12]{DKLMT22}}]\label{thm:ExplicitCouplingZd}
Let $k>d$ be two positive integers. Then there exists an orbit equivalence coupling from $\Z^{k}$ to $\Z^{d}$ which is $(\varphi_{\varepsilon},\psi_{\varepsilon})-$integrable for every $\varepsilon>0$, where
\begin{equation*}
        \varphi_{\varepsilon}(x)=\frac{x^{\frac{d}{k}}}{\ln{(x)}^{1+\varepsilon}}\; \text{and}\; \psi_{\varepsilon}(x)=\frac{x^{\frac{k}{d}}}{\ln{(x)}^{1+\varepsilon}}.
\end{equation*}
\end{theorem}

Note that, if $k>d$ and $p\ge\frac{d}{k}$, Theorems~\ref{thm:ObstructionDKLMT} and~\ref{thm:threshold} imply that we cannot have an $(\ld^p,\ld^0)$ orbit equivalence coupling from 
$\Z^{k}$ to $\Z^{d}$. Therefore the orbit equivalence coupling built in the last statement is optimal.

It is also possible to find a nice F\o lner tiling sequence for the lamplighter group $\Z/2\Z\wr\Z$, producing a quantitatively optimal coupling with $\Z$ (see~\cite[Proposition~6.20]{DKLMT22}). Using these techniques, Escalier also built in~\cite{Esc24} optimal couplings between $\Z$ and Brieussel-Zheng's groups defined in~\cite{BZ21}, which have prescribed isoperimetric profiles.

\subsection{F\o lner tiling sequences of halo products}

\paragraph{A general observation.}\label{sec:generalobservationtilinghalo}

By~\cite[Proposition~6.19]{DKLMT22}, we know a F\o lner tiling $(L_n)_{n\ge 1}$ of the lamplighter $F\wr\Z$, where $F$ is a finite group. It is defined by
\begin{equation}\label{eq:FolnerTilingLamplighter}
    L_{n}\defeq \left\lbrace (f,h)\in F\wr\Z : \supp{f}\subset\lbrace 0,1,\ldots,n-1\rbrace, h\in\lbrace 0,1,\ldots,n-1\rbrace\right\rbrace .
\end{equation}

In this section, we describe a general phenomenon that will be useful to construct F\o lner tiling sequences for our favorite halo products.

\begin{proposition}\label{prop:HowToBuildFolnerTiling}
Let $H$ be an amenable group and $\halo H$ be a halo product over $H$. Assume that $(F_{n})_{n\in\N}$ is a F\o lner tiling sequence of $H$, with an associated shift sequence $(\Sigma_{n})_{n\in\N}$, and that $L(F_{n})$ is finite for any $n\in\N$. Assume furthermore that for every $h\in\Sigma_{n}$, there exists a finite subset $\Sigma_{h,n}$ of $L(F_{n+1})$ such that we have
\begin{equation*}
        L(F_{n+1})=\bigsqcup_{\sigma\in\Sigma_{h,n}}{\sigma L(hF_{n})}.
\end{equation*}
Then $(L(F_{n})\times F_{n})_{n\in\N}$ satisfies the tiling condition (as defined in the beginning of Section~\ref{sec:PreliminariesFolnerTiling}), with shift sequence $\left(\bigcup_{h\in \Sigma_{n}}{(\Sigma_{h,n}\times\lbrace h\rbrace)}\right)_{n\in\N}$.
\end{proposition}

For instance, for the F\o lner tiling sequence described in~\eqref{eq:FolnerTilingLamplighter}, we start from the F\o lner tiling sequence $(\lbrace 0,1,\dots,2^n-1\rbrace)_{n\ge 0}$ of $\Z$, with shift sequence $(\lbrace 0,2^n\rbrace)_{n\ge 1}$, and we get $L_{n}=L(\lbrace 0,1,\dots,2^n-1\rbrace)$. Given $h\in\lbrace 0,2^n\rbrace$, $L(\lbrace h,h+1,\dots,h+2^n-1\rbrace)$ is the set of functions $\Z\rightarrow F$ of support in $\lbrace h,h+1,\dots,h+2^n-1\rbrace$. Thus, to get $L(\lbrace 0,1,\dots,2^{n+1}-1\rbrace)$ from $L(\lbrace h,h+1,\dots,h+2^n-1\rbrace)$, we have to translate it by all functions supported in $\lbrace 0,1,\dots,2^{n+1}\rbrace\setminus\lbrace h,h+1,\dots,h+2^{n}-1\rbrace$, this is exactly $\Sigma_{h,n}$.

\begin{proof}
Using the composition law, we see that
\begin{align*}
    \bigcup_{h\in\Sigma_{n}}{\bigcup_{\sigma\in\Sigma_{h,n}}{(\sigma,h)(L(F_{n})\times F_{n})}}&=\bigcup_{h\in\Sigma_{n}}{\bigcup_{\sigma\in\Sigma_{h,n}}{(\sigma L(hF_{n})\times hF_{n})}}\\
    &=\bigcup_{h\in\Sigma_n}{\left (\left (\bigcup_{\sigma\in\Sigma_{h,n}}{\sigma L(hF_{n})}\right)\times hF_{n}\right)}\\
    &=\bigcup_{h\in\Sigma_{n}}{(L(F_{n+1})\times hF_{n})}\\
    &=L(F_{n+1})\times F_{n+1},
\end{align*}
and it is straightforward to prove that the unions are disjoint.
\end{proof}

For the sequel, we will need to check that the sequence $(L(F_{n})\times F_{n})_{n\in\N}$ provided by this proposition is in fact a F\o lner sequence, and, in view of Theorem~\ref{thm:TilingSufficientConditionQuantitative}, we will actually need to quantify this F\o lner sequence, namely finding sequences $R=(R_{n})_{n\ge 0}$ and $\varepsilon=(\varepsilon_{n})_{n\ge 0}$ such that $(L(F_{n})\times F_{n})_{n\geq 0}$ is an $(R,\varepsilon)-$F\o lner tiling sequence. The following result goes in this direction.

\begin{proposition}\label{prop:sequencesinhaloproducts}
Let $H$ be a finitely generated amenable group and $\halo H$ be a halo product over $H$ with finite blocks. Let $(F_{n})_{n\in \N}$ be a sequence of non empty subsets of $H$ and $S_{H}$ be a finite generating set of $H$. For any $n\ge 0$, let 
\begin{equation*}
     L_{n}\defeq\lbrace (\sigma,h)\in\halo H :  \sigma\in L(F_n), h\in F_{n}\rbrace.
\end{equation*}
Then, for any $s\in S_H$, one has
\begin{equation*}
     \frac{|L_{n}\cdot (1_{L(H)},s)\setminus L_{n}|}{|L_{n}|}=\frac{|F_{n}s\setminus F_{n}|}{|F_{n}|}
\end{equation*}
as well as 
\begin{equation*}
         \frac{|L_{n}\cdot (\sigma_{s},1_{H})\setminus L_{n}|}{|L_{n}|}\leq\frac{|F_{n}s\setminus F_{n}|}{|F_{n}|}
\end{equation*}
for any $\sigma_{s}\in L(\lbrace 1_{H},s\rbrace)$. 
\end{proposition}

\begin{proof}
Let $(\sigma,h)\in L_{n}$ and let $s\in S_{H}$. The composition law of $\halo H$ directly implies that $(\sigma,h)(1_{L(H)},s)=(\sigma,hs)$ and $(\sigma,h)(\sigma_s,1_H)=(\sigma(h\cdot\sigma_s),h)$. Therefore $(\sigma,h)(1_{L(H)},s)$ lies in $L_n$ if and only if $hs$ lies in $F_n$, and since $h\cdot\sigma_s$ lies in $L(\{h,hs\})$, we know that $(\sigma,h)(\sigma_s,1_H)$ lies in $L_n$ if $hs$ lies in $F_n$. It thus follows that
\begin{equation*}
     \frac{|L_{n}\cdot(1_{L(H)},s)\setminus L_{n}|}{|L_{n}|}=\frac{|L_{n}|\cdot|F_{n}s\setminus F_{n}|}{|L_{n}|\cdot |F_n|}=\frac{|F_{n}s\setminus F_{n}|}{|F_{n}|}
\end{equation*}
and that
\begin{equation*}
     \frac{|L_{n}\cdot(\sigma_{s},1_{H})\setminus L_{n}|}{|L_{n}|}\leq\frac{|L_{n}|\cdot|F_{n}s\setminus F_{n}|}{|L_{n}|\cdot |F_{n}|}=\frac{|F_{n}s\setminus F_{n}|}{|F_{n}|}
\end{equation*}
and we are done.
\end{proof}

\begin{corollary}\label{cor:quantitativeFolner}
Let $H$ be a finitely generated amenable group and $\halo H$ be a halo product over $H$. Suppose that $\halo H$ is naturally generated and has finite blocks. Let $(F_{n})_{n\in \N}$ be a F\o lner sequence of $H$ and $S_{H}$ be a finite generating set of $H$. For any $n\ge 0$, let $L_{n} \defeq L(F_{n})\times F_n$. Then the set
\begin{equation*}
    S_{\halo H}\defeq\lbrace (1_{L(H)},s) : s\in S_{H}\rbrace \cup \bigcup_{s\in S_{H}}\left\lbrace (\sigma, 1_{H})\in\halo H : \sigma\in L(\lbrace 1_{H},s\rbrace)\right\rbrace
\end{equation*}
generates $\halo H$. Furthermore, for any $n\geq 1$, one has
\begin{equation*}
        \sup_{\varphi\in S_{\halo H}}{\frac{|L_{n}\varphi\setminus L_{n}|}{|L_{n}|}}\leq\sup_{s\in S_H}{\frac{|F_{n}s\setminus F_{n}|}{|F_{n}|}}.
\end{equation*}
\end{corollary}

This applies for instance to lamplighters, as was done in~\cite[Proposition~6.19]{DKLMT22}, and also to lampjugglers and lampcloners, as we will see in the next subsection. Thus, we find $\varepsilon_{n}$ from the quantitative information of the F\o lner tiling sequence $(F_n)_{n\ge 0}$ of $H$, and we compute the diameter $R_n$ with the generating set $S_{\halo H}$ provided by this last result. Furthermore, if the halo has consistent blocks, then the sets of the F\o lner tiling sequence have cardinality $\Lambda_{\halo H}(|F_{n}|)\cdot |F_{n}|$, where $\Lambda_{\halo H}$ is the lamp growth sequence of $\halo H$.

In the sequel, we illustrate these results for lampjugglers and lampcloners. Note that, for simplicity, we will only consider F\o lner tiling sequences $(F_{n})_{n\in\N}$ for the base group $H$ satisfying $1_{H}\in F_{n}$ for every $n\in\N$. This will be the case in our applications, when picking $H=\Z^{d}$.

\paragraph{F\o lner tiling sequences of lampshufflers and lampjugglers.}

Let us recall that if $S_{H}$ is a finite generating subset of $H$, then a finite generating set of $\juggler{r}{H}$ is 
\begin{equation*}
    S_{\halo H}=\left\lbrace(\tau_{(1_H,i),(s,j)},1_H) : s\in S_{H}, 1\le i,j\le r\right\rbrace\cup\left\lbrace(\mathrm{id},s) : s\in S_{H}\right\rbrace
\end{equation*}
where $\tau_{(g,i),(h,j)}\colon H\times\{1,\ldots,r\}\longrightarrow H\times\{1,\ldots,r\}$ denotes the transposition that swaps $(g,i)$ and $(h,j)$.

Let us now apply Proposition~\ref{prop:HowToBuildFolnerTiling} and Corollary~\ref{cor:quantitativeFolner} to get a F\o lner tiling sequence of $\juggler{r}{H}$ from a tiling sequence of $H$, with the associated quantitative information.

\begin{theorem}\label{thm:folnertilingshuffler}
Let $H$ be a finitely generated amenable group. If $(F_{n})_{n\in\N}$ is a right F\o lner tiling sequence of $H$ satisfying
\begin{equation*}
        \frac{|F_{n}s\setminus F_{n}|}{|F_{n}|}\le\varepsilon_n
\end{equation*}
for all $s\in S_{H}$ and some sequence $(\varepsilon_{n})_{n\in\N}$ tending to $0$,
then there exists a right F\o lner tiling sequence $(L_{n})_{n\in\N}$ of $\juggler{r}{H}$ with the following properties:
\begin{enumerate}[label=(\roman*)]
        \item\label{item:1Folnertilingjuggler} for any $n\in\N$, $|L_{n}|=(r|F_{n}|)!\cdot|F_{n}|$;
        \item\label{item:2Folnertilingjuggler} for any $n\in\N$, $\diam{L_{n}}\le 5\cdot|F_{n}|\cdot\diam{F_{n}}$;
        \item\label{item:3Folnertilingjuggler} for all $\varphi\in S_{\juggler{r}{H}}$, $\frac{|L_{n}\varphi\setminus L_{n}|}{|L_{n}|}\le\varepsilon_{n}$.
\end{enumerate}
\end{theorem}

To prove the upper bound on the diameter, we will need the following claim.

\begin{lemma}\label{lem:diamsymjuggler}
Let $F$ be a finite subset of $H$. Then, for every $\varphi\in\sym{F\times\lbrace 1,\dots,r\rbrace}$, we have
\begin{equation*}
     \left|(\varphi,1_{H})\right|_{\juggler{r}{H}}\le 2K+(5|F|-6)\diam{F},
\end{equation*}
where $K\defeq\displaystyle\max_{h\in F}{|h|_{S_H}}$.
\end{lemma}

\begin{proof}
Let us decompose $\varphi$ as a product of cycles $\varphi=c_{1}c_{2}\dots c_{k}$ on $F\times\lbrace 1,\dots,r\rbrace$ with disjoint supports. For $1\le i \le k$, $c_{i}$ is an $\ell_{i}-$cycle
\begin{equation*}
    (x_{i,1},\ldots,x_{i,\ell_{i}})=\left((h_{i,1},j_{i,1}),\dots,(h_{i,\ell_{i}},j_{i,\ell_{i}})\right)
\end{equation*}
which can itself be written as a product of transpositions
\begin{equation*}
c_{i}=\tau_{x_{i,1},x_{i,2}}\tau_{x_{i,2},x_{i,3}}\dots\tau_{x_{i,\ell_i-1},x_{i,\ell_i}}.
\end{equation*}
For any $1\le m\le \ell_{i}-1$, the transposition $\tau_{x_{i,m},x_{i,m+1}}$ can be written as
\begin{equation*}
\tau_{x_{i,m},x_{i,m+1}}=h_{i,m}\cdot\tau_{(1_H,j_{i,m}),(h_{i,m}^{-1}h_{i,m+1},j_{i,m+1})}
\end{equation*}
with the convention that $h_{i,\ell_{i}+1}\defeq h_{i,1}$ and $j_{i,\ell_{i}+1}\defeq j_{i,1}$, so that we have 
\begin{equation*}
(\tau_{x_{i,m},x_{i,m+1}},1_H)=(\mathrm{id},h_{i,m})(\tau_{(1_H,j_{i,m}),(h_{i,m}^{-1}h_{i,m+1},j_{i,m+1})},1_H)(\mathrm{id},h_{i,m}^{-1}).
\end{equation*}
Setting $y_{i,m}\defeq (1_H,j_{i,m})$ and $z_{i,m}\defeq (h_{i,m}^{-1}h_{i,m+1},j_{i,m+1})$, the element $(c_{i},1_{H})$ is thus equal to
\begin{equation*}
    (\mathrm{id},h_{i,1})(\tau_{y_{i,1},z_{i,1}},1_{H})(\mathrm{id},h_{i,1}^{-1}h_{i,2})(\tau_{y_{i,2},z_{i,2}},1_{H})\dots (\mathrm{id},h_{i,\ell_i-2}^{-1}h_{i,\ell_i-1}) (\tau_{y_{i,\ell_i-1},z_{i,\ell_i-1}},1_H)(\mathrm{id},h_{i,\ell_i-1}^{-1}),
\end{equation*}
namely $(c_{i},1_{H})=(\mathrm{id},h_{i,1}) \kappa_{i} (\mathrm{id},h_{i,\ell_i-1}^{-1})$, where 
\begin{equation*}
    \left|\kappa_{i}\right|_{\juggler{r}{H}}\le \left(4(\ell_{i}-1)+(\ell_{i}-2)\right)\diam{F}=(5\ell_{i}-6)\diam{F}
\end{equation*}
using Lemma~\ref{lem:wordlengthinlampshuffler}. Since we have
\begin{equation*}
        (\varphi,1_{H})=(\mathrm{id},h_{1,1})\kappa_1(\mathrm{id},h_{1,\ell_1-1}^{-1}h_{2,1})\kappa_2(\mathrm{id},h_{2,\ell_2-1}^{-1}h_{3,1})\ldots (\mathrm{id},h_{k-1,\ell_{k-1}-1}^{-1}h_{k,1})\kappa_k(\mathrm{id},h_{k,\ell_k-1}^{-1}),
\end{equation*}
we finally get
\begin{align*}
    \left|(\varphi,1_{H})\right|_{\juggler{r}{H}}&\le 2K+(k-1)\cdot \diam{F}+\sum_{i=1}^{k}{\left(5\ell_i-6\right)\cdot \diam{F}}\\
    &\le 2K+\left(k-1+(5|F|-6k)\right)\cdot \diam{F}\\
    &\le 2K+(5|F|-6)\cdot \diam{F}
\end{align*}
which concludes the proof.
\end{proof}

\begin{proof}[Proof of Theorem~\ref{thm:folnertilingshuffler}]
For $n\in\N$, let us set
\begin{equation*}
L_{n}\defeq\left\lbrace(\sigma,h)\in \juggler{r}{H} : \sigma\in \sym{F_{n}\times\lbrace 1,\dots ,r\rbrace}, h\in F_{n}\right\rbrace
\end{equation*}
where $\sym{F\times\lbrace 1,\dots ,r\rbrace}$ denotes the set of permutations supported in a finite subset $F\times\lbrace 1,\dots ,r\rbrace$ of $H$. Clearly,~\textit{\ref{item:1Folnertilingjuggler}} holds, as well as~\textit{\ref{item:3Folnertilingjuggler}} by Corollary~\ref{cor:quantitativeFolner}.

\noindent Given $(\sigma,h)$ and $(\sigma',h')$ in $L_{n}$, let us find an upper bound for $\left|(\sigma,h)^{-1}(\sigma',h')\right|_{S_{\juggler{r}{H}}}$. We first have
\begin{align*}
    (\sigma,h)^{-1}(\sigma',h')&=(h^{-1}\cdot\sigma^{-1},h^{-1})(\sigma',h')  \\
    &=\left([h^{-1}\cdot\sigma^{-1}]\circ [h^{-1}\cdot\sigma'],h^{-1}h'\right)  \\
    &=\left([h^{-1}\cdot\sigma^{-1}]\circ [h^{-1}\cdot\sigma'],1_H\right)\left(\mathrm{id},h^{-1}h'\right)
\end{align*}
so that
\begin{align*}
    \left|(\sigma,h)^{-1}(\sigma',h')\right|_{S_{\juggler{r}{H}}} &\le \left|([h^{-1}\cdot\sigma^{-1}]\circ [h^{-1}\cdot\sigma'],1_H)\right|_{S_{\juggler{r}{H}}}+\left|(\mathrm{id},h^{-1}h')\right|_{S_{\juggler{r}{H}}}  \\
    &\le \left|([h^{-1}\cdot\sigma^{-1}]\circ [h^{-1}\cdot\sigma'],1_{H})\right|_{S_{\juggler{r}{H}}}+\diam{F_{n}}
\end{align*}
and it remains to find an upper bound for $\left|([h^{-1}\cdot\sigma^{-1}]\circ [h^{-1}\cdot\sigma'],1_H)\right|_{S_{\juggler{r}{H}}}$. But the permutation $[h^{-1}\cdot\sigma^{-1}]\circ [h^{-1}\cdot\sigma']$ lies in $\sym{(h^{-1}F_{n})\times\lbrace 1,\dots,r\rbrace}$, so using Lemma~\ref{lem:diamsymjuggler}, it follows that
\begin{align*}
    \left|([h^{-1}\cdot\sigma^{-1}]\circ [h^{-1}\cdot\sigma'],1_H)\right|_{S_{\juggler{r}{H}}}&\le 2K+(5|h^{-1}F_n|-6)\cdot\diam{h^{-1}F_{n}},
\end{align*}
with $K\defeq\displaystyle\max_{k\in F_{n}}{|h^{-1}k|_{S_{H}}}\leq\diam{F_{n}}$, so that we have
\begin{align*}
    \left|([h^{-1}\cdot\sigma^{-1}]\circ [h^{-1}\cdot\sigma'],1_H)\right|_{S_{\juggler{r}{H}}}&\le (5|F_{n}|-4)\cdot\diam{F_{n}}.
\end{align*}
Hence, we have 
\begin{equation*}
    \left|(\sigma,h)^{-1}(\sigma',h')\right|_{S_{\juggler{r}{H}}}\le 5|F_{n}|\cdot\diam{F_{n}}
\end{equation*}
and we are done for the proof of~\textit{\ref{item:2Folnertilingjuggler}}.

\noindent Lastly, we show the tiling condition, using Proposition~\ref{prop:HowToBuildFolnerTiling}. Let $(\Sigma_{n})_{n\in\N}$ be a sequence of F\o lner shifts for $(F_{n})_{n\in\N}$, satisfying $F_{n+1}=\Sigma_{n} F_{n}$. Let us fix $h\in\Sigma_n$ and let us find a finite subset $\Sigma_{h,n}$ of $L(F_{n+1})$ such that we have the following disjoint union:
\begin{equation}\label{eq:DesiredDisjointUnionjuggler}
    \sym{F_{n+1}}=\bigsqcup_{\sigma\in\Sigma_{h,n}}{\sigma\  \sym{hF_n}}.
\end{equation}
For every $A\subset F_{n+1}\times\lbrace 1,\dots,r\rbrace$ of cardinality $r|F_{n}|$ and every family $\textbf{$x$}=(x_{g,i})_{(g,i)\in (F_{n+1}\setminus hF_n)\times\{1,\ldots,r\}}$ satisfying 
\begin{equation}\label{eq:CompatibleJuggler}
    A\sqcup\lbrace x_{g,i} : (g,i)\in (F_{n+1}\setminus hF_{n})\times\lbrace 1,\dots,r\rbrace\rbrace=F_{n+1}\times\lbrace 1,\dots,r\rbrace,
\end{equation}
let us choose one permutation $\tau^{A,\textbf{$x$}}\in\sym{F_{n+1}}$ satisfying
\begin{equation*}
        \tau^{A,\textbf{$x$}}\left((hF_{n})\times\lbrace 1,\dots,r\rbrace\right)=A \;\text{and} \; \forall (g,i)\in (F_{n+1}\setminus hF_{n})\times\lbrace 1,\dots,r\rbrace, \tau^{A,\textbf{$x$}}(g,i)=x_{g,i},
\end{equation*}
then the set
\begin{equation*}
    \left\lbrace\tau^{A,\textbf{$x$}}\circ\sigma : \sigma\in\sym{hF_{n}\times\lbrace 1,\dots,r\rbrace}\right\rbrace
\end{equation*}
describes all permutations $\rho\in\sym{F_{n+1}\times\lbrace 1,\dots,r\rbrace}$ satisfying
\begin{equation*}
\rho(hF_{n}\times\lbrace 1,\dots,r\rbrace)=A\;\text{and}\; \forall (g,i)\in (F_{n+1}\setminus hF_{n})\times\lbrace 1,\dots,r\rbrace, \rho(g,i)=x_{g,i}.
\end{equation*}
It remains to define $\Sigma_{h,n}$ as the set of all chosen permutations $\tau^{A,\textbf{$x$}}$ for every subset $A\subset F_{n+1}\times\lbrace 1,\dots,r\rbrace$ of cardinality $r|F_{n}|$ and every family $\textbf{$x$}=(x_g)_{g\in (F_{n+1}\setminus hF_{n})\times\lbrace 1,\dots r\rbrace}$ such that~\eqref{eq:CompatibleJuggler} holds, and we get~\eqref{eq:DesiredDisjointUnionjuggler}. So the F\o lner sequence $(L_{n})_{n\in\N}$ satisfies the tiling condition by Proposition~\ref{prop:HowToBuildFolnerTiling}.
\end{proof}

\paragraph{F\o lner tiling sequences of lampcloners.}

Let us recall that if $S_{H}$ is a finite generating subset of $H$, then a finite generating set of $\cloner{H}$ is 
\begin{equation*}
    S_{\cloner{H}}\defeq\lbrace (\delta_{1_{H}}(\lambda), 1_{H}) : \lambda \in \field\setminus\lbrace 0\rbrace\rbrace\cup\lbrace (\tau_{1_{H},s}(\lambda),1_{H}) : s\in S_{H}\rbrace \cup\lbrace (\text{id}, s) : s\in S_{H}\rbrace
\end{equation*}
where the elements $\delta_{p}(\lambda)$ are diagonal matrices and elements $\tau_{pq}(\lambda)$ are transvections.

Let us now apply Propositions~\ref{prop:HowToBuildFolnerTiling} and Corollary~\ref{cor:quantitativeFolner} to get a F\o lner tiling sequence of $\cloner{H}$ from a tiling sequence of $H$, with the associated quantitative information on it.

\begin{theorem}\label{thm:folnertilingcloner}
Let $H$ be a finitely generated amenable group. If $(F_{n})_{n\in\N}$ is a right F\o lner tiling sequence of $H$ satisfying
\begin{equation*}
        \frac{|F_{n}s\setminus F_{n}|}{|F_{n}|}\le\varepsilon_n
\end{equation*}
for all $s\in S_{H}$ and some sequence $(\varepsilon_{n})_{n\in\N}$ tending to $0$,
then there exists a right F\o lner tiling sequence $(L_{n})_{n\in\N}$ of $\cloner{H}$ with the following properties:
\begin{enumerate}[label=(\roman*)]
        \item\label{item:1Folnertilingcloner} for any $n\in\N$, $|L_{n}|=\Lambda_{\cloner{H}}(|F_n|)\cdot|F_{n}|$;
        \item\label{item:2Folnertilingcloner} there exists a constant $C>0$ such that for any $n\in\N$, $\diam{L_{n}}\le C|F_{n}|^3\diam{F_{n}}$;
        \item\label{item:3Folnertilingcloner} for all $\varphi\in S_{\cloner{H}}$, $\frac{|L_{n}\varphi\setminus L_{n}|}{|L_{n}|}\le\varepsilon_{n}$.
\end{enumerate}
\end{theorem}

To prove the upper bound on the diameter, we will need the following claim.

\begin{lemma}\label{lem:diamsymcloner}
Let $F$ be a finite subset of $H$. Then there exists a constant $C>0$ such that, for every $\varphi\in\mathrm{FGL(F)}$, we have
\begin{equation*}
     \left|(\varphi,1_{H})\right|_{S_{\cloner{H}}}\le 2K+C|F|^{3}\cdot \diam{F}
\end{equation*}
with $K\defeq\displaystyle\max_{h\in F}{|h|_{S_{H}}}$.
\end{lemma}

\begin{proof}
By Gaussian elimination, we need $N$ basic operations to get $\varphi$ from the identity in $\mathrm{FGL(F)}$, with $N=O(|F|^3)$. By operations, we mean transvections and dilatations of support in $F$. Let us notice that we have 
\begin{equation*}
    \tau_{g,h}(\lambda)=g\cdot\tau_{1_{H},g^{-1}h}(\lambda)
\end{equation*}
for transvections, so that we get 
\begin{equation*}
    (\tau_{g,h}(\lambda),1_{H})=(\mathrm{id},g)(\tau_{1_H,g^{-1}h}(\lambda),1_{H})(\mathrm{id},g^{-1})
\end{equation*}
in $\cloner{H}$. Also, for dilatations, we have
\begin{equation*}
    \delta_g(\lambda)=g\cdot\delta_{1_{H}}(\lambda)
\end{equation*}
and this can be written in $\cloner{H}$ as 
\begin{equation*}
    (\delta_g(\lambda),1_{H})=(\mathrm{id},g)(\delta_{1_{H}}(\lambda),1_{H})(\mathrm{id},g^{-1}).
\end{equation*}
We thus get $(\varphi,1_{H})$ as a product of the following form
\begin{equation*}
    (\varphi,1_{H})=(\mathrm{id},h_{1})(\varphi_{1},1_{H})(\mathrm{id},h_{1}^{-1}h_{2})(\varphi_{2},1_{H})(\mathrm{id},h_{2}^{-1}h_{3})\dots (\mathrm{id},h_{N-1}^{-1}h_{N})(\varphi_{N},1_{H})(\mathrm{id},h_{N}^{-1})
\end{equation*}
with basic operations $\varphi_{i}$ of support in $\lbrace 1_{H},F^{-1}F\rbrace$ and elements $h_{i}\in F$. Now using Lemma~\ref{lem:wordlengthinlampcloner}, we get
\begin{equation*}
    |(\varphi,1_{H})|_{S_{\cloner{H}}}\le 2K+N\cdot\diam{F}+14N\cdot\diam{F}
\end{equation*}
and we are done.
\end{proof}

\begin{proof}[Proof of Theorem~\ref{thm:folnertilingcloner}]
For $n\in\N$, let us set
\begin{equation*}
L_{n}\defeq\left\lbrace(\sigma,h)\in \cloner{H} : \sigma\in \mathrm{FGL}(F_n), h\in F_{n}\right\rbrace.
\end{equation*}
Clearly,~\textit{\ref{item:1Folnertilingcloner}} holds by definition of the lamp growth sequence $\Lambda_{\cloner{H}}$, as well as~\textit{\ref{item:3Folnertilingcloner}} by Corollary~\ref{cor:quantitativeFolner}.

\noindent Given $(\sigma,h)$ and $(\sigma',h')$ in $L_n$, let us find an upper bound for $\left|(\sigma,h)^{-1}(\sigma',h')\right|_{S_{\cloner{H}}}$. We first have
\begin{align*}
    (\sigma,h)^{-1}(\sigma',h')&=(h^{-1}\cdot\sigma^{-1},h^{-1})(\sigma',h')  \\
    &=\left([h^{-1}\cdot\sigma^{-1}]\circ [h^{-1}\cdot\sigma'],h^{-1}h'\right)  \\
    &=\left([h^{-1}\cdot\sigma^{-1}]\circ [h^{-1}\cdot\sigma'],1_{H}\right)\left(\mathrm{id},h^{-1}h'\right)
\end{align*}
so that
\begin{align*}
    \left|(\sigma,h)^{-1}(\sigma',h')\right|_{S_{\cloner{H}}} &\le \left|([h^{-1}\cdot\sigma^{-1}]\circ [h^{-1}\cdot\sigma'],1_{H})\right|_{S_{\cloner{H}}}+\left|(\mathrm{id},h^{-1}h')\right|_{S_{\cloner{H}}}  \\
    &\le \left|([h^{-1}\cdot\sigma^{-1}]\circ [h^{-1}\cdot\sigma'],1_{H})\right|_{S_{\cloner{H}}}+\diam{F_{n}}
\end{align*}
and it remains to find an upper bound for $\left|([h^{-1}\cdot\sigma^{-1}]\circ [h^{-1}\cdot\sigma'],1_{H})\right|_{S_{\cloner{H}}}$. But the linear automorphism $[h^{-1}\cdot\sigma^{-1}]\circ [h^{-1}\cdot\sigma']$ lies in $\mathrm{FGL}(h^{-1}F_{n})$, so using Lemma~\ref{lem:diamsymcloner}, it follows that
\begin{align*}
    \left|([h^{-1}\cdot\sigma^{-1}]\circ [h^{-1}\cdot\sigma'],1_{H})\right|_{S_{\cloner{H}}}&\le 2K+C|h^{-1}F_{n}|^3\cdot\diam{h^{-1}F_{n}},
\end{align*}
with $K\defeq\max_{k\in F_{n}}{|h^{-1}k|_{S_H}}\le\diam{F_{n}}$ and for some constant $C>0$, so that we have
\begin{align*}
    \left|([h^{-1}\cdot\sigma^{-1}]\circ [h^{-1}\cdot\sigma'],1_{H})\right|_{S_{\cloner{H}}}&\le (C+2)|F_{n}|^3\diam{F_{n}}.
\end{align*}
Hence, we have 
\begin{equation*}
    \left|(\sigma,h)^{-1}(\sigma',h')\right|_{S_{\cloner{H}}}\le (C+3)|F_{n}|^3\cdot\diam{F_{n}}
\end{equation*}
and we are done for the proof of~\textit{\ref{item:2Folnertilingcloner}}.

\noindent Lastly, we show the tiling condition, using Proposition~\ref{prop:HowToBuildFolnerTiling}. Let $(\Sigma_{n})_{n\in\N}$ be a sequence of F\o lner shifts for $(F_{n})_{n\in\N}$, satisfying $F_{n+1}=\Sigma_{n} F_{n}$. Let us fix $h\in\Sigma_n$ and let us find a finite subset $\Sigma_{h,n}$ of $L(F_{n+1})$ such that we have the following disjoint union:
\begin{equation}\label{eq:DesiredDisjointUnioncloner}
    \mathrm{FGL}(F_{n+1})=\bigsqcup_{\sigma\in\Sigma_{h,n}}{\sigma\  \mathrm{FGL}(hF_n)}.
\end{equation}
For every subspace $A$ of $V_{F_{n+1}}$, of dimension $|F_{n}|$ and every family $\textbf{$x$}=(x_{v})_{v\in V_{F_{n+1}\setminus hF_{n}}}$ satisfying 
\begin{equation}\label{eq:CompatibleCloner}
    A\sqcup\lbrace x_{v} : v\in V_{F_{n+1}\setminus hF_n}\rbrace=V_{F_{n+1}},
\end{equation}
let us choose, if it exists, one linear automorphism $\tau^{A,\textbf{$x$}}\in\mathrm{FGL}(F_{n+1})$ satisfying
\begin{equation*}
        \tau^{A,\textbf{$x$}}(V_{hF_{n}})=A \;\text{and} \; \forall v\in V_{F_{n+1}\setminus hF_{n}}, \tau^{A,\textbf{$x$}}(v)=x_{v},
\end{equation*}
then the set
\begin{equation*}
    \left\lbrace\tau^{A,\textbf{$x$}}\circ\sigma : \sigma\in\mathrm{FGL}(hF_{n})\right\rbrace
\end{equation*}
describes all the linear automorphisms $\rho\in\mathrm{FGL}(F_{n+1})$ satisfying
\begin{equation*}
\rho(V_{hF_{n}})=A\;\text{and}\; \forall v\in V_{F_{n+1}\setminus hF_{n}}, \rho(v)=x_{v}.
\end{equation*}
It remains to define $\Sigma_{h,n}$ as the set of all the chosen linear automorphisms $\tau^{A,\textbf{$x$}}$ for every subspace $A$ of $V_{F_{n+1}}$, of dimension $|F_{n}|$ and every family $\textbf{$x$}=(x_{v})_{v\in V_{F_{n+1}\setminus hF_n}}$ such that~\eqref{eq:CompatibleCloner} holds and such that such a $\tau^{A,\textbf{$x$}}$ exists, and we get~\eqref{eq:DesiredDisjointUnioncloner}. Hence the F\o lner sequence $(L_{n})_{n\in\N}$ satisfies the tiling condition by Proposition~\ref{prop:HowToBuildFolnerTiling}.
\end{proof}

\subsection{Applications to quantitatively optimal orbit equivalence couplings}

\paragraph{Lampshufflers and lampjugglers.}\label{sec:folnertilingjugglers}

As a first application of our construction of F\o lner tiling sequences, we give an explicit orbit equivalence coupling between $\juggler{r}{\Z^k}$ and $\Z^d$.

\begin{theorem}\label{thm:couplingsZdandshufZk}
Let $d,k,r\ge 1$ be three integers. There exists an orbit equivalence coupling from $\juggler{r}{\Z^k}$ to $\Z^d$, which is $(\varphi_{k,\varepsilon},\psi)-$integrable for every $\varepsilon>0$, where 
\begin{equation*}
    \varphi_{k,\varepsilon}(x)\defeq\frac{\ln(x)^{\frac{1}{k}}}{\ln(\ln(x))^{1+\frac{1}{k}+\varepsilon}}\; \text{and} \; \psi(x)=x^{x^{\frac{k}{k+1}}}.
\end{equation*}
\end{theorem}

\begin{remark}
Note that $\left(\frac{\ln(x)}{\ln(\ln(x))}\right)^{\frac{1}{k}}$ is (asymptotically) the isoperimetric profile of $\juggler{r}{\Z^k}$. Theorem~\ref{thm:couplingsZdandshufZk} implies in particular that, for every $\varepsilon>0$, there is a $\left(\frac{\prof{\juggler{r}{\Z^k}}(x)}{\ln(\ln(x))^{1+\varepsilon}},\ld^{<\infty}\right)-$integrable orbit equivalence coupling from $\juggler{r}{\Z^{k}}$ to $\Z^{d}$.
\end{remark}

\begin{proof}[Proof of Theorem~\ref{thm:couplingsZdandshufZk}]
Endow $\Z^{k}$ (resp. $\Z^d$) with its standard generating set $S_{\Z^k}$ (resp. $S_{\Z^d}$). Using Theorem~\ref{thm:folnertilingshuffler}, the F\o lner tiling sequence 
\begin{equation*}
    (F_{n})_{n\in\N}\defeq \left(\lbrace 0,\dots,d^{dn}-1\rbrace^{k}\right)_{n\in\N}
\end{equation*}
of $\Z^k$ provides a F\o lner tiling sequence $(L_{n})_{n\in\N}$ of $\juggler{r}{\Z^k}$ satisfying
\begin{itemize}
        \item $|L_{n}|=(rd^{dkn})!\cdot d^{dkn}$, for all $n\in\N$;
        \item $\diam{L_{n}}\le 5 d^{d(k+1)n+1}$ for all $n\in\N$;
        \item $\frac{|L_{n}s\setminus L_{n}|}{|L_{n}|}\le\frac{1}{d^{dn}}$, for all $s\in S_{\juggler{r}{\Z^k}}$.
\end{itemize}
Given a positive integer $n$, let us define, for every $i\in\lbrace 1,2,\ldots d\rbrace$, the integer
\begin{equation*}
    \ell_{i}(n)\defeq\prod_{\substack{j\in\N \\ dj+i\le n}}{(dj+i)}.
\end{equation*}
These quantities satisfy
\begin{equation*}
    \ell_{1}(n)\ell_{2}(n)\dots \ell_{d}(n)=n!,
\end{equation*}
\begin{equation*}
    \ell_{1}(dn)\le \ell_{2}(dn)\le\dots\le \ell_{d}(dn)\le \ell_{1}(d(n+1))
\end{equation*}
and
\begin{equation*}
    \ell_{d}(dn)=d^n\cdot n!\sim \left(\frac{dn}{e}\right)^{n}\sqrt{2\pi n}
\end{equation*}
using Stirling's formula. Now, we find a F\o lner tiling sequence $(G_{n})_{n\in\N}$ of $\Z^{d}$ satisfying $|L_{n}|=|G_{n}|$, setting
\begin{equation*}
    G_{n}\defeq\prod_{1\le i\le d}{\lbrace 0,1,\dots,d^{kn}\ell_i(rd^{dkn})\rbrace}.
\end{equation*}
It satisfies
\begin{itemize}
    \item for all $n\in\N$, $\diam{G_{n}}\le d^{kn}(\ell_1(rd^{dkn}) +\ldots +\ell_d(rd^{dkn}))\le d^{kn+1}\ell_d(rd^{dkn})$;
    \item for all $s\in S_{\Z^d}$, $\frac{|(s+G_{n})\setminus G_{n}|}{|G_{n}|}\le\frac{1}{d^{kn}\min{(\ell_1(rd^{dkn}),\ldots ,\ell_d(rd^{dkn}))}}\le\frac{1}{d^{kn}\ell_d(rd^{dkn-1})}$.
\end{itemize}

\noindent Observe that we have
\begin{equation*}
    \ln(d^{kn+1}\ell_{d}(rd^{dkn}))\sim rd^{dkn-1}\ln(rd^{dkn})=rd^{dkn-1}dkn\ln(d)
\end{equation*}
and
\begin{equation*}
    \ln(\ln(d^{kn+1}\ell_d(rd^{dkn})))\sim dkn\ln(d),
\end{equation*}
which in turn implies
\begin{equation*}
    \frac{\varphi_{k,\varepsilon}(d^{kn+1}\ell_{d}(d^{dkn}))}{d^{dn}}\sim\frac{K}{n^{1+\varepsilon}}
\end{equation*}
for some constant $K>0$. Hence the series $\displaystyle\sum_{n\ge 0}{\frac{\varphi_{k,\varepsilon}(d^{kn+1}\ell_d(d^{dkn}))}{d^{dn}}}$ converges. Secondly, we get 
\begin{equation*}
    \frac{\psi(5Cd^{d(k+1)n+1})}{d^{kn}\ell_{d}(rd^{dkn-1})}=O\left(\left ((5Cd)^{(5Cd)^{\frac{k}{k+1}}}\left (\frac{e}{r}\right )^{\frac{r}{d^2}}d^{(5Cd)^{\frac{k}{k+1}}d(k+1)n-\frac{rkn}{d}+\frac{r}{d^2}}\right )^{d^{dkn}}\right)
\end{equation*}
for every $C>0$. Taking $C$ small enough, the quantity $(5Cd)^{\frac{k}{k+1}}d(k+1)-\frac{rk}{d}$ is negative and the series $\sum{\frac{\psi(5C d^{d(k+1)n+1})}{d^{kn}\ell_d(rd^{dkn-1})}}$ converges. By Theorem~\ref{thm:TilingSufficientConditionQuantitative}, we get a $(\varphi_{k,\varepsilon},\psi)-$integrable orbit equivalence coupling from $\juggler{r}{\Z^{k}}$ to $\Z^{d}$.
\end{proof}

Next, we consider the case of iterated lampjugglers over free abelian groups. Using the notion of composition of couplings (see Theorem~\ref{thm:CompositionCouplings}), we deduce the following.

\begin{corollary}\label{cor:couplingsbetweeniteratedshufflers}
Let $n,m\ge 0$ be natural integers such that $m>n$. Let $d,k,r\ge 1$. Then there exists an orbit equivalence coupling from $\jugglern{m}{r}{\Z^k}$ to $\jugglern{n}{r}{\Z^d}$, which is $(\varphi_{m-n,k,\varepsilon}(x), \mathrm{L}^{<\infty})-$integrable for every $\varepsilon>0$, where 
\begin{equation*}
    \varphi_{i,k,\varepsilon}(x)\defeq\frac{\ln^{\circ i}(x)^{\frac{1}{k}}}{\left(\ln^{\circ (i+1)}(x)\right)^{1+\frac{1}{k}+\varepsilon}}.
\end{equation*}
\end{corollary}

One remark is in order before the proof.

\begin{remark}\label{rem:couplingsbetweeniteratedshufflers}
Note that we have
\begin{equation}\label{eq:profileComposedInverseProfile}
    \prof{\jugglern{m}{r}{\Z^k}}\circ\prof{\jugglern{n}{r}{\Z^d}}^{-1}(x)\simeq\left(\frac{\ln^{\circ (m-n)}(x)}{\ln^{\circ (m-n+1)}(x)}\right)^{\frac{1}{k}}
\end{equation}
for all natural integers $m,n,k,d$ satisfying $m>n$ and $d,k\ge 1$. Indeed, if we denote by $f$ the inverse of the isoperimetric profile of $\jugglern{n}{r}{\Z^d}$, we get
\begin{equation*}
    \ln^{\circ (n+1)}(f(x))\sim d\ln{x}
\end{equation*}
using the equality $\prof{\jugglern{n}{r}{\Z^d}}(f(x))=x$. This implies
\begin{align*}
    \left(\prof{\jugglern{m}{r}{\Z^k}}\circ\prof{\jugglern{n}{r}{\Z^d}}^{-1}(x)\right)^{k}&=\frac{\ln^{\circ m}{(f(x))}}{\ln^{\circ (m+1)}{(f(x))}} \\
    &=\frac{\ln^{\circ (m-n-1)}({\ln^{\circ (n+1)}(f(x)))}}{\ln^{\circ (m-n)}(\ln^{\circ (n+1)}(f(x)))}\\
    &\simeq\frac{\ln^{\circ (m-n)}(x)}{\ln^{\circ (m-n+1)}(x)}
\end{align*}
and~\eqref{eq:profileComposedInverseProfile} follows.

In particular, Corollary~\ref{cor:couplingsbetweeniteratedshufflers} implies that, for every $\varepsilon>0$, there is a $\left(\varphi,\ld^{<\infty}\right)$-integrable orbit equivalence coupling from $\jugglern{m}{r}{\Z^k}$ to $\jugglern{n}{r}{\Z^d}$, with
\begin{equation*}
    \varphi(x)=\frac{\prof{\jugglern{m}{r}{\Z^k}}\circ\prof{\jugglern{n}{r}{\Z^d}}^{-1}(x)}{\left (\ln^{\circ (m-n+1)}{(x)}\right )^{1+\varepsilon}}.
\end{equation*}
\end{remark}

\begin{proof}[Proof of Corollary~\ref{cor:couplingsbetweeniteratedshufflers}]
\sloppy From Theorem~\ref{thm:couplingsZdandshufZk}, there is an orbit equivalence from $\juggler{r}{\Z^d}$ to $\Z^d$ which is $(\varphi_{1,d,\varepsilon},\ld^{q})-$integrable for all $q>0$. Applying $i\ge 1$ times Theorem~\ref{thm:stabilityofcouplings+quantificationLampjugglers}, it follows that there is an orbit equivalence coupling from $\jugglern{(i+1)}{r}{\Z^d}$ to $\jugglern{i}{r}{\Z^d}$, which is $(\varphi_{1,d,\varepsilon},\ld^{q})-$integrable for all $q>0$. In particular, we have a coupling from $\jugglern{(n+1)}{r}{\Z^d}$ to $\jugglern{n}{r}{\Z^d}$, a coupling from $\jugglern{(n+2)}{r}{\Z^d}$ to $\jugglern{(n+1)}{r}{\Z^d}$,$\dots$, and a coupling from $\jugglern{(m-1)}{r}{\Z^d}$ to $\jugglern{(m-2)}{r}{\Z^d}$, which are all $(\varphi_{1,d,\varepsilon},\ld^{q})-$integrable for all $q>0$. By Theorem~\ref{thm:CompositionCouplings}, composing these successive couplings yields a coupling from $\jugglern{(m-1)}{r}{\Z^d}$ to $\jugglern{n}{r}{\Z^d}$, which is $(\varphi_{1,d,\varepsilon}^{\circ (m-n-1)},\ld^{q})-$integrable for all $q>0$. Using the asymptotic equivalence 
\begin{equation*}
    \varphi_{1,d,\varepsilon}^{\circ (m-n-1)}(x)\simeq\varphi_{m-n-1,d,\varepsilon}(x)
\end{equation*}
and the second item of Remark~\ref{rem:checkongenerators}, this coupling is $(\varphi_{m-n-1,d,\varepsilon},\ld^{q})-$integrable for all $q>0$.

\noindent By the same techniques, there is an orbit equivalence coupling from $\jugglern{m}{r}{\Z^k}$ to $\jugglern{(m-1)}{r}{\Z^d}$, which is $(\varphi_{1,k,\varepsilon},\ld^{q})-$integrable for all $q>0$, and composing this coupling with the one of the previous paragraph, we get a $(\varphi_{m-n,k,\varepsilon},\ld^{<\infty})-$integrable orbit equivalence from $\jugglern{m}{r}{\Z^k}$ to $\jugglern{n}{r}{\Z^d}$. This proves the corollary.
\end{proof}

Hence, we can also prove that:

\begin{corollary}\label{cor:optimalityIteratedLampshuffler}
\sloppy Let $n,m\ge 1$ be two integers such that $m>n$. There exists a $\left(\left(\frac{\ln^{\circ (m-n)}(x)}{\ln^{\circ (m-n+1)}(x)}\right)^{p},\mathrm{L}^0\right)-$integrable orbit equivalence coupling from $\jugglern{m}{r}{\Z^k}$ to $\jugglern{n}{r}{\Z^d}$ if and only if $p<\frac{1}{k}$. 
\end{corollary}

Once again, before the proof, let us recall that $\prof{\jugglern{m}{r}{\Z^k}}\circ\prof{\jugglern{n}{r}{\Z^d}}^{-1}(x)$ is asymptotically equivalent to $\left (\frac{\ln^{\circ (m-n)}(x)}{\ln^{\circ (m-n+1)}(x)}\right )^{\frac{1}{k}}$ (see Remark~\ref{rem:couplingsbetweeniteratedshufflers}).

\begin{proof}
Assume first that such an orbit equivalence exists. Then Theorem~\ref{thm:ObstructionDKLMT}\textit{(i)} implies that 
\begin{equation}\label{eq:ComparisonProfileIterated}
    \left (\frac{\ln^{\circ (m-n)}(\prof{\jugglern{n}{r}{\Z^d}}(x))}{\ln^{\circ (m-n+1)}(\prof{\jugglern{n}{r}{\Z^d}}(x))}\right )^p\preccurlyeq \prof{\jugglern{m}{r}{\Z^k}}(x).
\end{equation}
From Proposition~\ref{prop:profileofshufnofpolynomialgrowthgroupsINTRO}, we know that, for any integers $i,j\ge 1$, the isoperimetric profile of $\jugglern{j}{r}{\Z^{i}}$ is $\simeq \left(\frac{\ln^{\circ j}(x)}{\ln^{\circ (j+1)}(x)}\right)^{\frac{1}{i}}$, so that~\eqref{eq:ComparisonProfileIterated} implies 
\begin{equation*}
    \left(\frac{\ln^{\circ m}(x)}{\ln^{\circ (m+1)}(x)}\right)^{p} \preccurlyeq \left (\frac{\ln^{\circ m}(x)}{\ln^{\circ (m+1)}(x)}\right)^{\frac{1}{k}}.
\end{equation*}
This inequality now forces $p\le \frac{1}{k}$, and once again the fact that there is no $\left(\left(\frac{\ln^{\circ (m-n)}(x)}{\ln^{\circ (m-n+1)}(x)}\right)^{\frac{1}{k}}, \mathrm{L}^{0}\right)-$orbit equivalence coupling follows from Remark~\ref{rem:couplingsbetweeniteratedshufflers} and Theorem~\ref{thm:threshold}.

\noindent Conversely, assume that $p<\frac{1}{k}$. From Corollary~\ref{cor:couplingsbetweeniteratedshufflers}, we have an orbit equivalence coupling from $\jugglern{m}{r}{\Z^k}$ to $\jugglern{n}{r}{\Z^d}$ which is $(\varphi_{m-n,k,\varepsilon},\mathrm{L}^{<\infty})-$integrable for all $\varepsilon>0$. Since $p<\frac{1}{k}$, we have
\begin{equation*}
    \left(\frac{\ln^{\circ (m-n)}(x)}{\ln^{\circ (m-n+1)}(x)}\right )^{p}=O(\varphi_{m-n,k,\varepsilon}(x))
\end{equation*}
and it follows that our coupling is also $\left (\left(\frac{\ln^{\circ (m-n)}(x)}{\ln^{\circ (m-n+1)}(x)}\right )^{p},\mathrm{L}^{0}\right)-$integrable, as claimed. 
\end{proof}

About iterated lampjugglers with different numbers of iterations, let us point out the following result when the base groups are finitely generated and have slow profiles, which does not use F\o lner tiling sequences but our work on isoperimetric profiles presented in Chapter 5.

\begin{corollary}\label{cor:boundOfIntegrability}
\sloppy Let $n,m\ge 1$ be two integers such that $m>n$. Let $H$ be a finitely generated amenable group with isoperimetric profile $\prof{H}(x)\simeq \left (\ln^{\circ k}(x)\right )^{\alpha}$, for some integer $k>0$ and $\alpha>0$. If there exists a $\left(\left(\ln^{\circ (m-n)}(x)\right )^{p},\mathrm{L}^0\right)-$integrable orbit equivalence coupling from $\jugglern{m}{r}{H}$ to $\jugglern{n}{r}{H}$, then $p<\alpha$. 
\end{corollary}

\begin{proof}
We immediately get $p\leq\alpha$ by Theorem~\ref{thm:ObstructionDKLMT}. With similar techniques as in Remark~\ref{rem:couplingsbetweeniteratedshufflers}, we show that $\prof{\jugglern{m}{r}{H}}\circ\prof{\jugglern{n}{r}{H}}^{-1}$ is asymptotically equivalent to $\left (\ln^{\circ (m-n)}(x)\right)^{\alpha}$, so that Theorem~\ref{thm:threshold} implies $p\neq\alpha$.  
\end{proof}

\begin{remark}\label{rem:FolnerWreathProduct}
Following our methods, we can deduce much more examples of quantitatively optimal couplings. For instance, we get an optimal coupling between $\Z$ and $\juggler{r}{F\wr\Z}$ (with a non-trivial finite group $F$) by composing two couplings:
\begin{itemize}
    \item an optimal coupling between $\Z$ and $\juggler{r}{\Z}$;
    \item an optimal coupling between $\juggler{r}{\Z}$ and $\juggler{r}{F\wr\Z}$, coming from our stability result (Theorem~\ref{thm:stabilityofcouplings+quantificationLampjugglers}) and a coupling between $\Z$ and $F\wr\Z$, coming from~\cite[Proposition~6.20]{DKLMT22}.
\end{itemize}
\end{remark}

\paragraph{Lampcloners.} Let us finally apply our construction of F\o lner tiling sequences for couplings between $\cloner{\Z^k}$ and $\Z^d$.

\begin{theorem}\label{thm:couplingsZdandclonerZk}
Let $d,k\ge 1$ be integers. There exists an orbit equivalence coupling from $\cloner{\Z^k}$ to $\Z^d$, which is $(\varphi_{k,\varepsilon},\psi)-$integrable for every $\varepsilon>0$, where
\begin{equation*}
    \varphi_{k,\varepsilon}(x)\defeq\frac{\ln(x)^{\frac{1}{2k}}}{\ln(\ln(x))^{1+\varepsilon}}\; \text{and} \; \psi(x)=x^{x^{\frac{k}{3k+1}}}.
\end{equation*}
\end{theorem}

\begin{remark}
Theorem~\ref{thm:couplingsZdandshufZk} tells us that, for every $\varepsilon>0$, there is a $\left(\frac{j(x)}{\ln(\ln(x))^{1+\varepsilon}},\ld^{<\infty}\right)-$integrable orbit equivalence coupling from $\cloner{\Z^{k}}$ to $\Z^{d}$, where $j(x)\defeq \ln(x)^{\frac{1}{2k}}$ is a lower bound for $\prof{\cloner{\Z^k}}(x)$ that we find in Corollary~\ref{cor:encadrementduprofildecloneINTRO} (recall that we do not have precise estimates of the isoperimetric profile in this case).
\end{remark}

\begin{proof}[Proof of Theorem~\ref{thm:couplingsZdandclonerZk}]
Endow $\Z^{k}$ (resp. $\Z^d$) with its standard generating set $S_{\Z^k}$ (resp. $S_{\Z^d}$). Using Theorem~\ref{thm:folnertilingcloner}, the F\o lner tiling sequence 
\begin{equation*}
    (F_{n})_{n\in\N}\defeq \left(\lbrace 0,\dots,d^{dn}-1\rbrace^{k}\right)_{n\in\N}
\end{equation*}
of $\Z^k$ provides a F\o lner tiling sequence $(L_{n})_{n\in\N}$ of $\cloner{\Z^k}$ satisfying
\begin{itemize}
        \item $|L_{n}|=\Lambda_{\cloner{\Z^k}}(d^{dnk})\cdot d^{dkn}$, for all $n\in\N$;
        \item there exists a constant $C>0$ such that $\diam{L_{n}}\le C d^{d(3k+1)n+1}$ for all $n\in\N$;
        \item $\frac{|L_{n}s\setminus L_{n}|}{|L_{n}|}\le\frac{1}{d^{dn}}$, for all $s\in S_{\cloner{\Z^k}}$.
\end{itemize}
Let us recall that the lamp growth sequence of $\cloner{\Z^k}$ is given by
\begin{equation*}
    \forall n\ge 0,\; \Lambda_{\cloner{\Z^k}}(n)=\prod_{i=1}^{n}{(q^n-q^{n-i})}
\end{equation*}
where $q\defeq |\field|$. Given a positive integer $n$, let us define, for every $i\in\lbrace 1,2,\ldots d\rbrace$, the integer
\begin{equation*}
    \ell_{i}(n)\defeq\prod_{\substack{j\in\N \\ dj+i\le n}}{(q^n-q^{n-(dj+1)})}.
\end{equation*}
These quantities satisfy
\begin{equation*}
    \ell_{1}(n)\ell_{2}(n)\ldots \ell_{d}(n)=\Lambda_{\cloner{\Z^k}}(n),
\end{equation*}
\begin{equation*}
    \ell_{1}(dn)\le \ell_{2}(dn)\le\ldots\le \ell_{d}(dn)\le \ell_{1}(d(n+1))
\end{equation*}
and
\begin{equation*}
    \ell_{d}(dn)=\prod_{j=0}^{n-1}{(q^{dn}-q^{dn-d(j+1)})}=q^{dn^2}\prod_{j=0}^{n-1}{(1-q^{-d(j+1)})}\sim C'q^{dn^2}
\end{equation*}
where $C'$ is the limit of the convergent product $\prod_{j=0}^{\infty}{(1-q^{-d(j+1)})}$. Now, we find a F\o lner tiling sequence $(G_{n})_{n\in\N}$ of $\Z^{d}$ satisfying $|L_{n}|=|G_{n}|$, setting
\begin{equation*}
    G_{n}\defeq\prod_{1\le i\le d}{\lbrace 0,1,\dots,d^{kn}\ell_{i}(d^{dkn})\rbrace}.
\end{equation*}
It satisfies
\begin{itemize}
    \item for all $n\in\N$, $\diam{G_{n}}\le d^{kn}(\ell_{1}(d^{dkn}) +\dots +\ell_{d}(d^{dkn}))\le d^{kn+1}\ell_{d}(d^{dkn})$;
    \item for all $s\in S_{\Z^d}$, $\frac{|(s+G_{n})\setminus G_{n}|}{|G_{n}|}\le\frac{1}{d^{kn}\min{(\ell_{1}(d^{dkn}),\ldots ,\ell_{d}(d^{dkn}))}}\le\frac{1}{d^{kn}\ell_{d}(d^{dkn-1})}$.
\end{itemize}

\noindent Observe that we have
\begin{equation*}
    \ln{(d^{kn+1}\ell_d(d^{dkn}))}\sim d^{2dkn-1}\ln{q}
\end{equation*}
and
\begin{equation*}
    \ln(\ln(d^{kn+1}\ell_d(rd^{dkn})))\sim 2dkn\ln(d),
\end{equation*}
which in turn implies
\begin{equation*}
    \frac{\varphi_{k,\varepsilon}(d^{kn+1}\ell_{d}(d^{dkn}))}{d^{dn}}\sim\frac{K}{n^{1+\varepsilon}}
\end{equation*}
for some constant $K>0$. Hence the series $\displaystyle\sum_{n\ge 0}{\frac{\varphi_{k,\varepsilon}(d^{kn+1}\ell_{d}(d^{dkn}))}{d^{dn}}}$ converges. Secondly, we get 
\begin{equation*}
    \frac{\psi(C d^{d(3k+1)n+1})}{d^{kn}\ell_{d}(rd^{dkn-1})}=O\left(\left((Cd)^{(Cd)^{\frac{k}{3k+1}}}d^{(Cd)^{\frac{k}{3k+1}}d(3k+1)n}q^{-d^{dkn}}\right)^{d^{dkn}}\right)
\end{equation*}
so the series $\sum{\frac{\psi(C d^{d(3k+1)n+1})}{d^{kn}\ell_{d}(d^{dkn-1})}}$ converges. By Theorem~\ref{thm:TilingSufficientConditionQuantitative}, we get a $(\varphi_{k,\varepsilon},\psi)-$integrable orbit equivalence coupling from $\cloner{\Z^k}$ to $\Z^d$.
\end{proof}

Hence, as in Corollary~\ref{cor:couplingsbetweeniteratedshufflers}, we can deduce the following.

\begin{corollary}\label{cor:couplingsbetweeniteratedjugglers}
Let $\field$ be a finite field. Let $n,m\ge 0$ be natural integers such that $m>n$. Let $d,k,r\ge 1$. Then there exists an orbit equivalence coupling from $\clonern{n}{\Z^d}$ to $\clonern{m}{\Z^k}$, which is $(\mathrm{L}^{<\infty},\varphi_{m-n,k,\varepsilon}(x))-$integrable for every $\varepsilon>0$, where 
\begin{equation*}
    \varphi_{i,k,\varepsilon}(x)\defeq\frac{\ln^{\circ i}(x)^{\frac{1}{2k}}}{\left(\ln^{\circ (i+1)}(x)\right)^{1+\varepsilon}}.
\end{equation*}
\end{corollary}

\begin{proof}
    This is the same proof as Corollary~\ref{cor:couplingsbetweeniteratedshufflers}.
\end{proof}

\section{Comments and questions}\label{sec:commentsandquestions}

In this final section, we record some questions that naturally arise from what has been done in the article.

One of our main results is the construction of an explicit orbit equivalence coupling between the lampjugglers $\juggler{r}{H}$ and $\juggler{r}{K}$, provided a coupling between $H$ and $K$. Additionally, our construction preserves quantification. Hence:

\begin{question}
Let $r,s\ge 1$, $r\neq s$. How can one construct an orbit equivalence coupling between $\juggler{s}{H}$ and $\juggler{r}{K}$ from a coupling between $H$ and $K$? Moreover, if $H$ and $K$ are $(\varphi,\psi)-$integrably orbit equivalent, does the same hold for $\juggler{s}{H}$ and $\juggler{r}{K}$?
\end{question}

Beyond the quasi-isometric classification of lampshufflers established in~\cite{GT24a}, another important direction of the latter is to compare geometrically different classes of halo products, to rule out or to establish the existence of a quasi-isometry (resp. quasi-isometric/coarse/regular embedding) between, for instance, a lampshuffler and a lamplighter, or a lampdesigner and a lamplighter. Analogous questions are also relevant from the measured point of view. For example, note that the isoperimetric profiles of $\shuf{\Z}$ and $\Z/2\Z \wr\Z$ are asymptotically equivalent up to a logarithmic factor, so we may wonder:

\begin{question}
Are $\shuf{\Z}$ and $\Z/2\Z \wr\Z$ $\ld^{p}$-orbit equivalent for $p\in \left]0,1\right[$?
\end{question}

Moreover, there are pairs of amenable groups with asymptotically equivalent isoperimetric profiles, such as
\begin{itemize}
    \item $\shuf{\Z}$ and $\Z\wr\Z$;
    \item $\shuf{H}$ and $\Z/2\Z\wr H$, with $H$ such that Assumption~$(\star)$ holds and $\prof{H}\left(\frac{\ln(x)}{\ln(\ln(x))}\right) \simeq \prof{H}(\ln(x))$,
\end{itemize}
and thus, we may wonder whether these pairs of groups are $\ld^{p}$ orbit equivalent for every $p>0$.

More generally, even without quantifications, it would be interesting to have an explicit description of actions from a wreath product and a lampshuffler on a common probability space sharing the same orbits.

\begin{question}
What could be an explicit orbit equivalence coupling between a lampshuffler and a wreath product?
\end{question}

Alternatively, another useful technique to produce orbit equivalences is that of F\o lner tiling sequences, recalled and implemented in Section~\ref{sec:folnerTiling} in the case of lampshufflers. More precisely, Theorem~\ref{thm:folnertilingshuffler} provides an explicit description of a F\o lner tiling sequence for $\juggler{r}{H}$ from such a sequence for $H$. However, this technique remains hard to apply for the comparison between lamplighters and lampjugglers, and we did not manage to find F\o lner tiling sequences of both groups whose tiles have the same cardinalities. Also for lampcloners, we have an explicit description of F\o lner tiling sequences (Theorem~\ref{thm:folnertilingcloner}) which does not enable us to get a coupling between a lampcloner and a lamplighter, or between a lampcloner and a lampshuffler.

\afterpage{\blankpage}

\fancyhead{} % clear all header fields

%----------------------------------------------------------------------------------------
%	BIBLIOGRAPHY
%----------------------------------------------------------------------------------------
%\printbibliography %Prints bibliography
\fancyhead{} % clear all header fields
\fancyhead[OL]{\textsc{References}}
\printbibliography[
heading=bibintoc,
title={References}]

@article{AM24,
    author = {Abbott, Carolyn and Martinez-Pedroza, Eduardo},
    title = {The quasi-isometry invariance of the coset intersection complex},
    journal = {Algebraic \& Geometric Topology},
    volume = {26},
    number = {2},
    pages = {659--698},
    year = {2026},
}

@article{AT19,
    author = {Arzhantseva, Goulnara and Tessera, Romain},
    title = {Admitting a coarse embedding is not preserved under group extensions},
    journal = {International Mathematics Research Notices},
    volume = {2019},
    number = {20},
    pages = {6480--6498},
    year = {2019},
}

@article{Aus16,
    author = {Austin, Tim},
    title = {Integrable measure equivalence for groups of polynomial growth},
    journal = {Groups, Geometry and Dynamics},
    volume = {10},
    number = {1},
    year = {2016},
    pages = {117--154},
    shorthand = {Aus16},
}

@article{BCMR16,
  title = {Commensurations and metric properties {{of Houghton}}'s groups},
  author = {Burillo, José and Cleary, Sean and Martino, Armando and R{\"o}ver, Claas},
  year = {2016},
  journal = {Pacific Journal of Mathematics},
  volume = {285},
  number = {2},
  pages = {289--301},
  shorthand = {BCMR16},
}

@article{BD08,
    author = {Bell, Greg and Dranishnikov, Alexander},
    title = {Asymptotic dimension},
    journal = {Topology and its Applications},
    year = {2008},
    volume = {155},
    pages = {1265--1296},
}

@article{BE12,
    author = {Bartholdi, Laurent and Erschler, Anna},
    title = {Growth of permutational extensions},
    journal = {Inventiones Mathematicae},
    volume = {189},
    pages = {431--455},
    year = {2012},
}

@article{BE14,
    author = {Bartholdi, Laurent and Erschler, Anna},
    title = {Groups of given intermediate word growth},
    journal = {Annales de l'Institut Fourier},
    volume = {64},
    number = {5},
    pages = {2003--2036},
    year = {2014},
}

@article{BE17,
    author = {Bartholdi, Laurent and Erschler, Anna},
    title = {Poisson-Furstenberg boundary and growth of groups},
    journal = {Probability Theory and Related Fields},
    volume = {168},
    number = {1-2},
    pages = {347--372},
    year = {2017},
}

@article{Beh+12,
  title = {Geometry and rigidity of mapping class groups},
  author = {Behrstock, Jason and Kleiner, Bruce and Minsky, Yair and Mosher, Lee},
  journal = {Geometry and Topology},
  volume = {16},
  number = {2},
  pages = {781--888},
  year = {2012},
}

@article{Bel03,
    author = {Belegradek, Igor},
    title = {On co-Hopfian nilpotent groups},
    journal = {Bulletin of the London Mathematical Society},
    year = {2003},
    volume = {35},
    number = {6},
    pages = {805--811}, 
}

@article{Ben26,
     author = {Bensaid, Oussama},
     title = {Coarse embeddings of symmetric spaces and Euclidean buildings},
     journal = {Mathematische Zeitschrift},
     year = {2026},
     volume = {312},
     number = {62},
}

@article{BFF24,
    author = {Bradford, Henry and Fournier-Facio, Francesco},
    title = {Hopfian wreath products and the stable finiteness conjecture},
    journal ={Mathematische Zeitschrift},
    volume = {308},
    number = {58},
    year = {2024},
    shorthand = {BFF24},
}

@unpublished{BGT24,
    author = {Bensaid, Oussama and Genevois, Anthony and Tessera, Romain},
    title = {Coarse separation and the large-scale geometry of wreath products},
    note = {arXiv:2401.18025},
    year = {2024},
}

@unpublished{BGT26a,
    author = {Bensaid, Oussama and Genevois, Anthony and Tessera, Romain},
    title = {Coarse separation and splittings in hyperbolic groups},
    note = {arXiv:2603.17852},
    year = {2026},
    shorthand = {BGT26a},
}

@unpublished{BGT26b,
    author = {Bensaid, Oussama and Genevois, Anthony and Tessera, Romain},
    title = {Coarse separation and splittings in right-angled Artin groups},
    note = {arXiv:2603.24706},
    year = {2026},
    shorthand = {BGT26b},
}

@unpublished{BLP15,
    author = {Burillo, José and Lopez-Platon, Eric},
    title = {Metric properties and distortion in wreath products},
    note = {arXiv:1506.06935},
    year = {2015},
}

@article{BST12,
    author = {Benjamini, Itai and Schramm, Oded and Timar, Adam},
    title = {On the separation profile of infinite graphs},
    journal = {Groups, Geometry and Dynamics},
    year = {2012},
    volume = {6},
    number = {4},
    pages = {639--658},
}

@article{BZ19,
  title = {Shalom's property {{$H_{FD}$}} and extensions by {{$\mathbb{Z}$}} of locally finite groups},
  author = {Brieussel, J{\'e}r{\'e}mie and Zheng, Tianyi},
  year = {2019},
  journal = {Israel Journal of Mathematics},
  volume = {230},
  number = {1},
  pages = {45--70},
}

@article{BZ21,
  title = {Speed of random walks, isoperimetry and compression of finitely generated groups},
  author = {Brieussel, J{\'e}r{\'e}mie and Zheng, Tianyi},
  journal = {Annals of Mathematics},
  volume = {193},
  number = {1},
  pages = {1--105},
  year = {2021},
}

@article{Cas16,
  title = {Embeddability and quasi-isometric classification of partially commutative groups},
  author = {{Casals-Ruiz}, Montserrat},
  year = {2016},
  journal = {Algebraic \& Geometric Topology},
  volume = {16},
  number = {1},
  pages = {597--620},
}

@book{CecAdd21,
    author = {Ceccherini-Silberstein, Tullio and D'Adderio, Michele},
    title = {Topics in Groups and Geometry},
    publisher = {Springer Monographs in Mathematics},
    year = {2021},
    shorthand = {CecAdd21},
}

@misc{cordum25,
  title = {Isoperimetric profiles of lamplighter-like groups},
  author = {Correia, Corentin and Dumoncel, Vincent},
  note = {arXiv:2506.13235},
  year = {2025},
  shorthand = {CD25},
}

@misc{cordum26,
  title = {On quantitative orbit equivalence between lamplighter-like groups},
  author = {Correia, Corentin and Dumoncel, Vincent},
  note = {arXiv:2604.14945},
  year = {2026},
  shorthand = {CD26},
}

@article{CFP96,
    author = {Cannon, James and Floyd, William and Parry, Walter},
    title = {Introductory notes on
    Richard Thompson’s groups},
    journal = {L'Enseignement Math\'ematique},
    year = {1996},
    volume = {42},
    number = {3},
    pages = {215--256},
}

@book{CH16,
  title = {Metric geometry of locally compact groups},
  author = {Cornulier, Yves and de la Harpe, Pierre},
  publisher = {European Mathematical Society},
  year = {2016},
  shorthand = {CH16},
}

@article{CK11,
    author = {Cornulier, Yves and Kar, Aditi},
    title = {On property (FA) for wreath products},
    journal = {Journal of Group Theory},
    volume = {14},
    number = {1},
    pages = {165-174},
    year = {2011},
}

@article{Cor06,
  title = {Finitely presented wreath products and double coset decompositions},
  author = {Cornulier, Yves},
  journal = {Geometriae Dedicata},
  volume = {122},
  pages = {89--108},
  year = {2006},
}

@article{Cor25,
    author = {Correia, Corentin},
    title = {On the absence of quantitatively critical measure equivalence couplings},
    journal = {Proceedings of the American Mathematical Society},
    year = {2025},
}

@article{Corn14,
    author = {Cornulier, Yves},
    title = {Groupes pleins-topologiques d'après Matui, Juschenko, Monod, ...},
    journal = {Astérisque},
    pages = {177--217},
    year = {2014},
    shorthand = {Corn14},
}

@article{Corn16,
  title = {Gradings on {{Lie}} algebras, systolic growth, and cohopfian properties of nilpotent groups},
  author = {Cornulier, Yves},
  year = {2016},
  journal = {Bulletin de la Société Mathématique de France},
  volume = {144},
  number = {4},
  pages = {693--744},
  shorthand = {Corn16},
}

@incollection{Cou00,
  title = {Random walks and geometry on infinite graphs},
  booktitle = {Lecture {{Notes}} on {{Analysis}} in {{Metric Spaces}}},
  author = {Coulhon, Thierry},
  year = {2000},
  pages = {5--36},
  publisher = {Pisa: Scuola Normale Superiore},
  address = {Trento, Italy}
}

@article{CSC93,
    author = {Coulhon, Thierry and Saloff-Coste, Laurent},
    title = {Isopérimétrie pour les groupes et les variétés},
    journal = {Revista Matem{\'a}tica Iberoamericana},
    volume = {9},
    number = {2},
    pages = {293--314},
    year = {1993},
}

@article{CSV08,
    author = {Cornulier, Yves and Stalder, Yves and Valette, Alain},
    title = {Proper actions of lamplighter groups associated with free groups},
    journal = {C. R. Math. Acad. Sci. Paris},
    volume = {346},
    number = {3},
    pages = {173--176},
    year = {2008},
}

@article{CSV12,
    author = {Cornulier, Yves and Stalder, Yves and Valette, Alain},
    title = {Proper actions of wreath products and generalizations},
    journal = {Transactions of the American Mathematical Society},
    volume = {364},
    number = {6},
    pages = {3159--3184},
    year = {2012},
}

@article{DFX23,
   author = {Dymarz, Tullia and Fisher, David and Xiangdong, Xie},   
   title = {A fibered Tukia theorem for nilpotent Lie groups},
   journal = {Annales Fennici Mathematici},
   volume = {48},
   number = {2},
   year = {2023},
   shorthand = {DFX23},
}

@article{DFX25,
   author = {Dymarz, Tullia and Fisher, David and Xiangdong, Xie},   
   title = {A Tukia-type theorem for nilpotent Lie groups and quasi-isometric rigidity of solvable groups},
   journal = {Advances in Mathematics},
   volume = {468},
   year = {2025},
   shorthand = {DFX25},
}

@book{DK18,
    author = {Drutu, Cornelia and Kapovich, Michael},
    title = {Geometric group theory},
    publisher = {American Mathematical Society},
    volume = {63},
    year = {2018},
}

@article{DKLMT22,
  title = {Quantitative measure equivalence between amenable groups},
  author = {Delabie, Thiebout and Koivisto, Juhani and Le Maître, François and Tessera, Romain},
  journal = {Annales Henri Lebesgue},
  volume = {5},
  pages = {1417--1487},
  year = {2022},
  shorthand = {DKLMT22},
}

@book{dlH00,
   author = {de la Harpe, Pierre},
   title = {Topics in Geometric group theory},
   publisher = {University of Chicago Press},
   year = {2000},
}

@misc{DLIT25,
  title = {$L^p$ measure equivalence of nilpotent groups},
  author = {Delabie, Thiebout and Isenrich, Claudio Llosa and Tessera, Romain},
  note = {arXiv:2505.17865},
  year = {2025},
  shorthand = {DLIT25},
}

@article{DO11,
    author = {Davis, Tara and Olshanskii, Alexander},
    title = {Subgroup distortion in wreath products of cyclic groups},
    journal = {Journal of Pure and Applied Algebra},
    volume = {215},
    number = {12},
    pages = {2987--3004},
    year = {2011},
}

@article{DPT15,
    author = {Dymarz, Tullia and Peng, Irine and Taback, Jennifer},
    title = {Bilipschitz versus quasi-isometric equivalence for higher rank lamplighter groups},
    journal = {New York Journal of Mathematics},
    volume = {21},
    pages = {129--150},
    year = {2015},
}

@unpublished{Dum24,
    author = {Dumoncel, Vincent},
    title = {On the quasi-isometric classification of permutational wreath products},
    note = {arXiv:2409.20159},
    year = {2024},
    shorthand = {Dum24},
}

@article{Dum26,
  author = {Dumoncel, Vincent},
  title = {Quasi-isometric rigidity for lamplighters with lamps of polynomial growth},
  journal = {Annales de l'Institut Fourier},
  year = {2026},
  shorthand = {Dum26},
}

@article{Dun85,
    author = {Dunwoody, Martin John},
    title = {The accessibility of finitely presented groups},
    journal = {Inventiones Mathematicae},
    volume = {81},
    number = {3},
    pages = {449--457},
    year = {1985},
}

@article{Dym10,
    author = {Dymarz, Tullia},
    title = {Bilipschitz equivalence is not equivalent to quasi-isometric equivalence for finitely generated groups},
    journal = {Duke Mathematical Journal},
    volume = {154},
    number = {3},
    pages = {509--526},
    year = {2010},
}

@article{Dyu00,
    author = {Dyubina, Anna},
    title = {Instability of the virtual solvability and the property of being virtually torsion-free for quasi-isometric groups},
    journal = {International Mathematical Research Notices},
    volume = {2000},
    number = {21},
    pages = {1097--1101},
    year = {2000},
}

@article{EFW12,
  title = {Coarse differentiation of quasi-isometries I: Spaces not quasi-isometric to Cayley graphs},
  author = {Eskin, Alex and Fisher, David and Whyte, Kevin},
  journal = {Annals of Mathematics},
  number = {2},
  volume = {176},
  number = {1},
  pages = {221--260},
  year = {2012},
}

@article{EFW13,
  title = {Coarse differentiation of quasi-isometries II: Rigidity for Sol and lamplighter groups},
  author = {Eskin, Alex and Fisher, David and Whyte, Kevin},
  journal = {Annals of Mathematics},
  number = {2},
  volume = {177},
  number = {3},
  pages = {869--910},
  year = {2013},
}

@article{Ers03,
    author = {Erschler, Anna},
    title = {On isoperimetric profiles of finitely generated groups},
    journal = {Geometriae Dedicata},
    volume = {100},
    pages = {157--171},
    year = {2003},
}

@article{Ers06,
  title = {Generalized wreath products},
  author = {Erschler, Anna},
  journal = {International Mathematics Research Notices},
  pages = {57835},
  year = {2006},
}

@article{Esc24,
  title = {Building prescribed quantitative orbit equivalence with the integers},
  author = {Escalier, Amandine},
  journal = {Groups, Geometry, and Dynamics},
  year = {2024},
  volume = {18},
  number = {3},
  pages = {1007-1035},
}

@article{EZ18,
    author = {Erschler, Anna and Zheng, Tianyi},
    title = {Growth of periodic Grigorchuk groups},
    journal = {Inventiones Mathematicae},
    volume = {219},
    pages = {1069--1155},
    year = {2020},
}

@article{EZ21,
  title = {Isoperimetric inequalities, shapes of {{F{\o}lner}} sets and groups with {{Shalom}}'s property {{$H_{FD}$}}},
  author = {Erschler, Anna and Zheng, Tianyi},
  year = {2021},
  journal = {Annales de l'Institut Fourier},
  volume = {70},
  number = {4},
  pages = {1363--1402},
}

@article{FM98,
    author = {Farb, Benson and Mosher, Lee},
    title = {A rigidity theorem
    for the solvable Baumslag-Solitar groups},
    journal = {Inventiones Mathematicae},
    volume = {131},
    pages = {419--451},
    year = {1998},
}

@article{FM99,
    author = {Farb, Benson and Mosher, Lee},
    title = {Quasi-isometric rigidity for the solvable Baumslag-Solitar groups, II},
    journal = {Inventiones Mathematicae},
    volume = {137},
    pages = {613--649},
    year = {1999},
}

@article{FS96,
    author = {Farb, Benson and Schwartz, Richard},
    title = {The large-scale geometry of Hilbert modular groups.},
    journal = {Journal of Differential Geometry},
    volume = {44},
    number = {3},
    pages = {435--478},
    year = {1996},
}

@phdthesis{Gen17a,
  title = {Cubical-like geometry of quasi-median graphs and applications to geometric group theory},
  author = {Genevois, Anthony},
  note = {arXiv:1712.01618},
  year = {2017},
}

@article{Gen22,
    author = {Genevois, Anthony},
    title = {Lamplighter groups, median spaces, and a-T-menability},
    journal = {Proceedings of the Edinburgh Mathematical Society},
    volume = {65},
    number = {2},
    pages = {500--529},
    year = {2022},
}

@article{Gen24,
    author = {Genevois, Anthony},
    title = {Automorphisms of graph products of groups and acylindrical hyperbolicity},
    journal = {Memoirs of the American Mathematical Society},
    volume = {301},
    year = {2024},
}

@article{GM19,
    author = {Genevois, Anthony and Martin, Alexandre},
    title = {Automorphisms of graph products of groups from a geometric perspective},
    journal = {Proceedings of the London Mathematical Society},
    volume = {119},
    number = {6},
    pages = {1745-1779},
    year = {2019},
}

@article{Gri84,
   author = {Grigorchuk, Rostislav},
   title = {Degrees of growth of finitely generated groups and the theory of invariant means},
   journal = {Izv. Akad. Nauk SSSR Ser. Mat.},
   volume = {48},
   pages = {939--985},
   year = {1984},
}

@article{Gro81,
    author = {Gromov, Misha},
    title = {Groups of polynomial growth and expanding maps},
    journal = {Publications Math\'ematiques de l'Institut des Hautes \'Etudes Scientifiques},
    volume = {53},
    pages = {53--73},
    year = {1981},
}

@incollection{Gro93,
  title = {Asymptotic invariants of infinite groups},
  author = {Gromov, Misha},
  booktitle = {Geometric group theory, Vol. 2 (Sussex, 1991)},
  publisher = {Cambridge Univ. Press},
  address = {Cambridge},
  pages = {1--295},
  year = {1993},
}

@article{GT22,
  title = {Measure-scaling quasi-isometries},
  author = {Genevois, Anthony and Tessera, Romain},
  journal = {Geometriae Dedicata},
  volume = {216},
  number = {3},
  pages = {34},
  year = {2022},
}

@misc{GT24a,
  title = {Lamplighter-like geometry of groups},
  author = {Genevois, Anthony and Tessera, Romain},
  note = {arXiv:2401.13520},
  year = {2024},
  shorthand = {GT24a},
}

@article{GT24b,
  title = {Asymptotic geometry of lamplighters over one-ended groups},
  author = {Genevois, Anthony and Tessera, Romain},
  year = {2024},
  journal = {Inventiones Mathematicae},
  volume = {238},
  number = {1},
  pages = {1--67},
  shorthand = {GT24b},
}

@article{GT25,
  title = {A note on morphisms to wreath products},
  author = {Genevois, Anthony and Tessera, Romain},
  journal = {Mathematical Proceedings of the Cambridge Philosophical Society},
  year = {2025},
}

@article{GV20,
    author = {Genevois, Anthony and Varghese, Olga},
    title = {Conjugating automorphisms of graph products: Kazhdan's property (T) and SQ-universality},
    journal = {Bulletin of the Australian Mathematical Society},
    volume = {101},
    number = {2},
    pages = {272--282},
    year = {2020},
}

@article{HM23,
    author = {Hughes, Sam and Martinez-Pedroza, Eduardo},
    title = {Hyperbolically embedded subgroups
and quasi-isometries of pairs},
    journal = {Canadian Mathematical Bulletin},
    volume = {66},
    number = {3},
    pages ={827--843},
    year = {2023},
}

@article{HMS21,
    author = {Hughes, Sam and Martinez-Pedroza, Eduardo and S\'anchez Saldana, Luis Jorge},
    title = {A survey on quasi-isometries of pairs: invariants and rigidity},
    journal = {Geometrical methods in group theory: papers dedicated to Ruth Charney, Séminaires et Congrès},
    volume = {34},
    pages = {137--153},
    year = {2025},
}

@article{HMT20,
  title = {Poincar{\'e} profiles of groups and spaces},
  author = {Hume, David and Mackay, John and Tessera, Romain},
  year = {2020},
  journal = {Revista Matem{\'a}tica Iberoamericana},
  volume = {36},
  number = {6},
  pages = {1835--1886},
}

@article{HMT22,
  title = {Poincar{\'e} profiles of {{Lie}} groups and a coarse geometric dichotomy},
  author = {Hume, David and Mackay, John and Tessera, Romain},
  year = {2022},
  journal = {Geometric and Functional Analysis},
  volume = {32},
  number = {5},
  pages = {1063--1133},
}

@article{HMT25,
    author = {Hume, David and Mackay, John and Tessera, Romain},
    title = {Asymptotic dimension for covers with controlled growth},
    journal = {Journal of the London Mathematical Society},
    volume = {111},
    number = {2}, 
    year = {2025},
}

@article{HO16,
  title = {Transitivity degrees of countable groups and acylindrical hyperbolicity},
  author = {Hull, Michael and Osin, Denis},
  year = {2016},
  journal = {Israel Journal of Mathematics},
  volume = {216},
  number = {1},
  pages = {307--353},
}

@article{JM13,
    author = {Juschenko, Kate and Monod, Nicolas},
    title = {Cantor systems, piecewise translations and simple amenable groups},
    journal = {Annals of Mathematics},
    year = {2013},
    volume = {178},
    number = {2},
    pages = {775--787},
}

@article{KL97,
    author = {Kleiner, Bruce and Leeb, Bernhard},
    title = {Rigidity of quasi-isometries for symmetric spaces and Euclidean buildings},
    journal = {Publications math\'ematiques de l'IHES},
    volume = {86},
    pages = {115--197},
    year = {1997},
}

@misc{LB25,
    title = {Rigidity and flexibility results for groups with a common cocompact envelope},
     author = {Le Boudec, Adrien},
     note = {arXiv:2510.24581},
     year = {2025},
}

@article{Lev24,
   author = {Levitin, Daniel},
   title = {Metric spaces of arbitrary finitely generated scaling group},
   journal = {Indiana University Mathematics Journal},
   volume = {73},
   number = {4},
   pages = {1357-1399},
   year = {2024},
}

@article{Li10,
    author = {Li, Sean},
    title = {Compression bounds for wreath products},
    journal = {Proceedings of the American Mathematical Society},
    volume = {138},
    number = {8},
    pages = {2701--2714},
    year = {2010}, 
}

@book{Loh17,
   author = {Löh, Clara},
   title = {Geometric group theory: An introduction},
   publisher = {Springer},
   year = {2017},
   shorthand={Loh17},
}

@article{LPP96,
    author = {Lyons, Russell and Pemantle, Robin and Peres, Yuval},
    title = {Random walks on the lamplighter group},
    journal = {Annals of Probability},
    volume = {24},
    number = {4},
    pages = {1993--2006},
    year = {1996},
}

@article{LS22,
    author = {Leemann, Paul-Henri and Schneeberger, Gr\'egoire},
    title = {Property (FW) and wreath products of groups: a simple approach using Schreier graphs},
    journal = {Expositiones Mathematicae},
    volume = {40},
    number = {4},
    pages = {1261--1270},
    year = {2022},
}

@article{LS24,
    author = {Leemann, Paul-Henri and Schneeberger, Gr\'egoire},
    title = {Wreath products of groups acting with bounded orbits},
    journal = {L'Enseignement Math\'ematique},
    volume = {70},
    number = {2},
    pages = {121--149},
    year = {2024},
}

@article{Lys85,
  title = {A system of defining relations for a {{Grigorchuk}} group},
  author = {Lysenok, Igor},
  year = {1985},
  journal = {Mathematical Notes of the Academy of Sciences of the USSR},
  volume = {38},
  number = {4},
  pages = {784--792},
}

@article{Mar21,
    author = {Margolis, Alex},
    title = {The geometry of groups containing almost normal subgroups},
    journal = {Geometry and Topology},
    volume = {25},
    number = {5},
    pages = {2405--2468},
    year = {2021},
}

@article{Mil68,
    author = {Milnor, John},
    title = {Growth of finitely generated solvable groups},
    journal = {Journal of Differential Geometry},
    volume = {2},
    number = {4},
    pages = {447--449},
    year = {1968},
    shorthand = {Mil68},
}

@article{MO10,
    author = {Monod, Nicolas and Ozawa, Narutaka},
    title = {The Dixmier problem, lamplighters and Burnside groups},
    journal = {Journal of Functional Analysis},
    volume = {258},
    number = {1},
    pages = {255--259},
    year = {2010},
}

@article{Mon22,
    author = {Monod, Nicolas},
    title = {Lamplighters and the bounded cohomology of Thompson's group},
    journal = {Geometric and Functional Analysis},
    volume = {32},
    pages = {662--675},
    year = {2022},
}

@book{NG23,
    author = {Nowak, Piotr and Yu, Guoliang},
    title = {Large-scale Geometry},
    publisher = {European Mathematical Society},
    year = {2023},
    shorthand = {NG23},
}

@article{NP11,
  title = {Scale-invariant groups},
  author = {Nekrashevych, Volodymyr V. and Pete, G{\'a}bor},
  year = {2011},
  journal = {Groups, Geometry, and Dynamics},
  volume = {5},
  number = {1},
  pages = {139--167},
}

@article{OW80,
    author = {Ornstein, Donald and Weiss, Benjamin},
    title = {Ergodic theory of amenable group actions. I: The Rohlin lemma.},
    journal = {Bulletin (New Series) of the American Mathematical Society},
    year = {1980}
}

@article{Pit95,
  title = {F{\o}lner sequences in polycyclic groups},
  author = {Pittet, Christophe},
  year = {1995},
  journal = {Revista Matem{\'a}tica Iberoamericana},
  volume = {11},
  number = {3},
  pages = {675--685},
}

@article{Pit00,
  title = {The isoperimetric profile of homogeneous {{Riemannian}} manifolds},
  author = {Pittet, Christophe},
  year = {2000},
  journal = {Journal of Differential Geometry},
  volume = {54},
  number = {2},
}

@misc{PS22,
    title = {Coulhon$-$Saloff-Coste isoperimetric inequalities for finitely generated groups},
     author = {Pittet, Christophe and Stankov, Bogdan},
     note = {arXiv:2211.03227},
     year = {2022},
}

@article{PSC02,
    author = {Pittet, Christophe and Saloff-Coste, Laurent},
    title = {On random walks on wreath products},
    journal = {Annals of Probability},
    volume = {30},
    number = {2},
    pages = {948--977},
    year = {2002},
}

@article{Ril22,
    author = {Riley, Tim},
    title = {Exponentially distorted subgroups in wreath products},
    journal = {The Quarterly Journal of Mathematics},
    volume = {75},
    number = {4},
    pages = {1355--1361},
    year = {2024},
}

@article{Sch96,
    author = {Schwartz, Richard Evan},
    title = {The quasi-isometry classification of rank one lattices},
    journal = {Publications math\'ematiques de l'IHES},
    volume = {82},
    pages = {133--168},
    year = {1996},
}

@article{SCZ15,
  title = {Random walks on free solvable groups},
  author = {{Saloff-Coste}, Laurent and Zheng, Tianyi},
  year = {2015},
  journal = {Mathematische Zeitschrift},
  volume = {279},
  number = {3-4},
  pages = {811--848},
}

@article{SCZ16,
  title = {Random walks and isoperimetric profiles under moment conditions},
  author = {{Saloff-Coste}, Laurent and Zheng, Tianyi},
  year = {2016},
  journal = {The Annals of Probability},
  volume = {44},
  number = {6},
}

@article{SCZ18,
  title = {Isoperimetric profiles and random walks on some permutation wreath products},
  author = {{Saloff-Coste}, Laurent and Zheng, Tianyi},
  year = {2018},
  journal = {Revista Matem{\'a}tica Iberoamericana},
  volume = {34},
  number = {2},
  pages = {481--540},
}

@article{SCZ21,
  title = {Isoperimetric profiles and random walks on some groups defined by piecewise actions},
  author = {{Saloff-Coste}, Laurent and Zheng, Tianyi},
  journal = {Probability Theory and Related Fields},
  volume = {181},
  pages = {711--756},
  year = {2021},
}

@article{Sha04,
  title = {Harmonic analysis, cohomology, and the large-scale geometry of amenable groups},
  author = {Shalom, Yehuda},
  year = {2004},
  journal = {Acta Mathematica},
  volume = {192},
  number = {2},
  pages = {119--185},
}

@article{Sil24,
  title = {The {{Poisson}} boundary of lampshuffler groups},
  author = {Silva, Eduardo},
  year = {2024},
  journal = {Mathematische Zeitschrift},
  volume = {308},
  number = {3},
  pages = {46},
}

@article{Sta68,
  title = {On torsion-free groups with infinitely many ends},
  author = {Stallings, John},
  journal = {Annals of Mathematics},
  volume = {88},
  pages = {312--334},
  year = {1968},
}

@article{Tes08,
  title = {Large scale Sobolev inequalities on
  metric measure spaces and applications},
  author = {Tessera, Romain},
  journal = {Revista Matem{\'a}tica Iberoamericana},
  volume = {24},
  number = {3},
  pages = {825--864},
  year = {2008},
}

@article{Tes13,
    title = {Isoperimetric profiles and random walks on locally compact solvable groups},
    author = {Tessera, Romain},
    year = {2013},
    journal = {Revista Matem{\'a}tica Iberoamericana},
    volume = {29},
    number = {2},
    pages = {715--737},
}

@article{Tes25,
    title = {On coarse embeddings of amenable groups into hyperbolic graphs},
    author = {Tessera, Romain},
    journal = {Proceedings of the American Mathematical Society},
    volume = {153},
    pages = {4545--4552},
    year = {2025},
}

@article{Tuk86,
    title = {On quasiconformal groups},
    author = {Tukia, Pekka},
    journal = {Journal d'Analyse Mathématique},
    volume = {46},
    pages = {318--346},
    year = {1986},
    shorthand = {Tuk86},
}

@article{Tuk94,
    title = {Convergence groups and Gromov’s metric hyperbolic spaces},
    author = {Tukia, Pekka},
    journal = {New Zealand Journal of Mathematics},
    volume = {23},
    pages = {157--187},
    year = {1994},
    shorthand = {Tuk94},
}

@article{Why99,
    author = {Whyte, Kevin},
    title = {Amenability, bi-Lipschitz equivalence, and the Von Neumann conjecture},
    journal = {Duke Mathematical Journal},
    volume = {99},
    number = {1},
    pages = {93--112},
    year = {1999},
}

@article{Why01,
    author = {Whyte, Kevin},
    title = {The large scale geometry of the higher Baumslag-Solitar groups},
    journal = {Geometric and Functional Analysis},
    volume = {11},
    number = {6},
    pages = {1327--1343},
    year = {2001},
}

@article{Woe05,
    author = {Woess, Wolfgang},
    title = {Lamplighters, Diestel-Leader graphs, random walks, and harmonic functions},
    journal = {Combinatorics, Probability and Computing},
    volume = {14},
    number = {3},
    pages = {415--433},
    year = {2005},
}

@article{Wol68,
    author = {Wolf, Joseph},
    title = {Growth of finitely generated solvable groups and curvature of Riemannian manifolds},
    journal = {Journal of Differential Geometry},
    volume = {2},
    number = {4},
    pages = {421--446},
    year = {1968},
    shorthand = {Wol68},
}

@article{Yad09,
  title = {Rate of {{escape}} of the {{mixer chain}}},
  author = {Yadin, Ariel},
  journal = {Electronic Communications in Probability},
  volume = {14},
  pages = {347--357},
  year = {2009},
 }

%----------------------------------------------------------------------------------------

%: ----------------------- glossary ------------------------
%\printindex
%\fancyhead{} % clear all header fields
%\fancyhead[OL]{\textsc{Glossaire}}
%\printglossary[title=Glossaire,toctitle=Glossaire]

%\input{auxilliaires/glossaire} 

\end{document}